\documentclass[12pt]{amsbook}
\usepackage{amssymb}
\usepackage[all]{xy}

\usepackage[colorlinks=true,linkcolor=magenta,citecolor=magenta]{hyperref}
\makeindex

\newtheorem{theorem}{Theorem}[chapter]
\newtheorem{answer}[theorem]{Answer}
\newtheorem{cat}[theorem]{Cat}
\newtheorem{claim}[theorem]{Claim}
\newtheorem{comment}[theorem]{Comment}
\newtheorem{conclusion}[theorem]{Conclusion}

\newtheorem{convention}[theorem]{Convention}
\newtheorem{definition}[theorem]{Definition}
\newtheorem{exercise}[theorem]{Exercise}
\newtheorem{fact}[theorem]{Fact}
\newtheorem{goal}[theorem]{Goal}
\newtheorem{idea}[theorem]{Idea}
\newtheorem{principle}[theorem]{Principle}
\newtheorem{problem}[theorem]{Problem}
\newtheorem{proposition}[theorem]{Proposition}
\newtheorem{question}[theorem]{Question}
\newtheorem{questions}[theorem]{Questions}
\newtheorem{speculation}[theorem]{Speculation}

\begin{document}

\title{Basic quantum algebra}

\author{Teo Banica}
\address{Department of Mathematics, University of Cergy-Pontoise, F-95000 Cergy-Pontoise, France. {\tt teo.banica@gmail.com}}


\begin{abstract}
This is an introduction to quantum algebra, from a geometric perspective. The classical spaces $X$, such as the Lie groups, homogeneous spaces, or more general manifolds, are described by various algebras $A$, defined over various fields $F$. These algebras $A$ satisfy a commutativity type condition, and the general idea is that of lifting this condition, and calling quantum spaces the underlying space-like objects $X$. One problem comes from the fact that different fields $F$ lead, via different algebras $A$, to different classes of quantum spaces $X$. Our aim here is to identify and put at the center of the presentation those quantum spaces $X$ which do not depend on the choice of $F$.
\end{abstract}

\maketitle

\chapter*{Preface}

The classical spaces $X$, such as the Lie groups, homogeneous spaces, or more general manifolds, are described by various algebras $A$, defined over various fields $F$. These algebras $A$ typically satisfy a commutativity type condition, such as $fg=gf$ when $A$ is a usual algebra of functions, and the idea of quantum algebra is that of lifting this commutativity condition, and calling quantum spaces the underlying space-like objects $X$. With the hope that these quantum spaces $X$ can be useful in physics.

\bigskip

When performed without further axioms, and precise physics motivations in mind, this procedure can turn into something dangerous, leading to a considerable amount of quantum spaces $X$. Which can be bewildering even for professionals, and with the first victims being of course the beginners, who quickly get lost in the quantum jungle.

\bigskip

In short, quantum algebra is a wide discipline. There are many ways of getting introduced to it, and books dedicated to this, with each book being based of course on some sort of philosophy, usually coming from authors' own research work and interests. The present book is one such introduction, our main philosophy being as follows:

\bigskip

(1) Quantum groups come before quantum spaces. The point indeed is that quantum groups $G$ are far easier to construct and study than quantum spaces $X$, via all sorts of algebras associated to them, such as function algebras, Lie algebras, Brauer algebras, planar algebras, and so on. As for the quantum spaces $X$, these can safely come afterwards, first in the form of quotient spaces $X=G/H$, and then in more general forms.

\bigskip

(2) Quantum groups are quantum versions of classical groups. This is something more subtle, the idea being that a quantum group, when constructed in the most general framework, that of the Hopf algebras, can be related or not to a classical group $G$. Our policy here will be that of calling quantum groups the beasts of type $G^\times$, with $G$ being a classical group, and with $\times$ standing for a liberation and twisting operation. 

\bigskip

(3) The main quantum groups and spaces are the absolute ones. Again, this is something subtle, coming from the fact that, in the general context of quantum algebra, different fields $F$ lead, via different algebras $A$, to different classes of quantum spaces $X$. In view of this, we will center our presentation on the quantum spaces $X$ which are absolute, in the sense that they do not depend on the choice of the ground field $F$.

\bigskip

In addition to this, although we will be certainly committed to algebraic aspects, we will be interested in analysis too, by doing some whenever the conditions permit. Finally, regarding physics, we will be mainly mathematical here, but with the hope that ``the best quantum spaces should lead to the best quantum physics''. And with best for us meaning absolute in the above sense. Time will tell if this is right or wrong.

\bigskip

The book is organized in four parts, as follows:

\bigskip

Part I - General introduction to Hopf algebras, which in quantum space terms means general introduction to quantum group type objects $G$, without further axioms.

\bigskip

Part II - Preliminary look at what can be the quantum groups, with a study, in Hopf algebra terms, of the quantum symmetry groups of various basic objects.

\bigskip

Part III - Introduction to quantum groups, with a quantum group being for us, as explained above, something of type $G^\times$, related to a classical group $G$. 

\bigskip

Part IV - Introduction to quantum spaces, starting with basic quotient spaces of type $X=G/H$, followed by more general types of quantum algebraic manifolds. 

\bigskip

This book, while being independent and readable as such, completes my series of quantum group books \cite{ba1}, \cite{ba2}, \cite{ba3}, in a somewhat circular way. To be more precise, in case you like the material at the end, which is certainly quite specialized, and want to know more about that, my previous books \cite{ba1}, \cite{ba2}, \cite{ba3} are here for that.

\bigskip

The book contains, besides the basics of quantum algebra, some recent contributions as well, and I would like to thank Julien Bichon, Beno\^ it Collins, Steve Curran and the others, for our joint work. Many thanks go as well to Nicol\'as Andruskiewitsch and Sonia Natale, for a number of useful algebraic discussions. Finally, it is a pleasure to thank my cats (meow) for precious help in the preparation of the present book.

\bigskip

\

{\em Cergy, July 2025}

\smallskip

{\em Teo Banica}

\baselineskip=15.95pt
\tableofcontents
\baselineskip=14pt

\part{Hopf algebras}

\ \vskip50mm

\begin{center}
{\em Nu-mi lua iubirea daca pleci

Mai lasa-ma sa cred in ea

In noptile pustii si reci

Sa-mi mai inchipui ca-i a mea}
\end{center}

\chapter{Hopf algebras}

\section*{1a. Hopf algebras}

The classical spaces $X$, such as the Lie groups, homogeneous spaces, or more general manifolds, are described by various algebras $A$, defined over various fields $F$. These algebras $A$ typically satisfy a commutativity type condition, such as $fg=gf$ when $A$ is a usual algebra of functions, and the idea of quantum algebra is that of lifting this commutativity condition, and calling quantum spaces the underlying space-like objects $X$. With the hope that these quantum spaces $X$ can be useful in physics.

\bigskip

In this chapter we start developing quantum algebra, with some inspiration from classical group theory. We would like to develop a theory of suitable algebras $A$, which are not necessarily commutative, corresponding to quantum groups $G$. 

\bigskip 

In what regards the classical constructions $G\to A$ that we have in mind, in view of a noncommutative extension, there are three of them, all well-known and widely used in group theory, which can be informally described, modulo several details, as follows:

\index{function algebra}
\index{convolution product}
\index{group algebra}
\index{Lie algebra}
\index{enveloping Lie algebra}

\begin{fact}
A group $G$ can be typically described by several algebras:
\begin{enumerate}
\item We have the algebra $F(G)$ of functions $f:G\to F$, with the usual, pointwise product of functions. This algebra is commutative, $fg=gf$.

\item We have the algebra $F[G]$ of functions $f:G\to F$, with the convolution product of functions. This algebra is commutative when $G$ is abelian.

\item We have the algebra $U\mathfrak g$ generated by the functions $f:G\to F$ infinitesimally defined around $1$. This algebra is commutative when $G$ is abelian.
\end{enumerate}
\end{fact}

As already mentioned, this is something quite informal, just meant to help us in order to start this book, and do not worry, we will come back to this later, with details. At the present stage of things, the comments to be made on this are as follows:

\bigskip

(1) The construction there is something quite simple and solid, and makes sense as stated, with the remark however that when $G$ is a topological group, things get more complicated, because we can further ask for the functions $f:G\to F$ that we use to be continuous, or even smooth, or vanishing at $\infty$, or be measurable, and so on, leading to several interesting versions of $F(G)$. Also, the recovery of $G$ from the algebra $F(G)$, or one of its versions, when $G$ is topological, is usually a non-trivial question.

\bigskip

(2) In what regards the construction there, pretty much the same comments apply, the point being that when $G$ is a topological group, things get more complicated, again leading to several interesting versions of the algebra $F[G]$, and with the recovery of the group $G$ itself, out of these algebras, being usually a non-trivial question. We will be back to all this later, with details, and in the meantime, you can simply consider, by using Dirac masses, that $F[G]$ is the formal span of the group elements $g\in G$.

\bigskip

(3) What we said there is definitely informal, the idea being that, when $G$ is a Lie group, we can consider its tangent space at the origin, or Lie algebra $\mathfrak g=T_1(G)$, consisting of the functions $f:G\to F$ infinitesimally defined around $1$, and then the corresponding enveloping Lie algebra $U\mathfrak g$, with product such that the Lie bracket is given by $[x,y]=xy-yx$. All this is quite non-trivial, notably with a discussion in relation with the field $F$ being needed, but do not worry, we will come back to it, with details.

\bigskip

Looking now at our list (1,2,3) above, it looks like (1) is the simplest construction, and the most adapted to our noncommutative goals, followed by (2), followed by (3). So, let us formulate the following goal, for the theory that we want to develop:

\begin{goal}
We want to develop a theory of associative algebras $A$ over a given field $F$, with some extra structure, as follows:
\begin{enumerate}
\item As main and motivating examples, we want to have the algebras $F(G)$.

\item We also want our theory to include, later, the algebras $F[G]$ and $U\mathfrak g$.

\item And we also want, later, to discuss what happens when $G$ is topological.
\end{enumerate}
\end{goal}

Needless to say, this goal is formulated quite informally, but this is just a goal, and if at this point you can see right away a complete and rigorous theory doing the job, that would be of course very welcome, and I will look myself for something else to do.

\bigskip

Getting started now, we would first like to have a look at the algebras $F(G)$ that we want to generalize. But before that, let us have a closer look at the groups $G$ themselves, with algebraic motivations in mind, in relation with the algebras $F(G)$. As our first result in this book, we have the following frightening reformulation of the group axioms:

\index{group}
\index{group axioms}

\begin{proposition}
A group is a set $G$ with operations as follows,
$$m:G\times G\to G\quad,\quad
u:\{.\}\to G\quad,\quad
i:G\to G$$
which are subject to the following axioms, with $\delta(g)=(g,g)$:
$$m(m\times id)=m(id\times m)$$
$$m(u\times id)=m(id\times u)=id$$
$$m(i\times id)\delta=m(id\times i)\delta=1$$
In addition, the inverse map $i$ satisfies $i^2=id$.
\end{proposition}

\begin{proof}
Our claim is that the formulae in the statement correspond to the axioms satisfied by the multiplication, unit and inverse map of $G$, given by:
$$m(g,h)=gh\quad,\quad u(.)=1\quad,\quad i(g)=g^{-1}$$

Indeed, let us start with the group axioms for $G$, which are as follows:
$$(gh)k=g(hk)$$
$$1g=g1=g$$
$$g^{-1}g=gg^{-1}=1$$ 

With $\delta(g)=(g,g)$ being as in the statement, these group axioms read:
$$m(m\times id)(g,h,k)=m(id\times m)(g,h,k)$$
$$m(u\times id)(g)=m(id\times u)(g)=g$$
$$m(i\times id)\delta(g)=m(id\times i)\delta(g)=1$$

Now since these must hold for any $g,h,k$, they are equivalent, as claimed, to:
$$m(m\times id)=m(id\times m)$$
$$m(u\times id)=m(id\times u)=id$$
$$m(i\times id)\delta=m(id\times i)\delta=1$$

As for $i^2=id$, this is something which holds too, coming from $(g^{-1})^{-1}=1$:
$$(g^{-1})^{-1}=1\iff i^2(g)=g\iff i^2=id$$

Thus, we are led to the various conclusions in the statement.
\end{proof}

The above result does not look very healthy, and might make Sophus Lie, Felix Klein and the others turn in their graves, but for our purposes here, this is exactly what we need. Indeed, turning now to the algebra $F(G)$, we have the following result:

\index{functional transpose}
\index{comultiplication}
\index{counit}
\index{antipode}
\index{square of antipode}
\index{algebra of functions}

\begin{theorem}
Given a finite group $G$, the functional transposes of the structural maps $m,u,i$, called comultiplication, counit and antipode, are as follows, 
$$\Delta:A\to A\otimes A\quad,\quad 
\varepsilon:A\to F\quad,\quad 
S:A\to A$$
with $A=F(G)$ being the algebra of functions $\varphi:G\to F$. The group axioms read:
$$(\Delta\otimes id)\Delta=(id\otimes \Delta)\Delta$$
$$(\varepsilon\otimes id)\Delta=(id\otimes\varepsilon)\Delta=id$$
$$m(S\otimes id)\Delta=m(id\otimes S)\Delta=\varepsilon(.)1$$
In addition, the square of the antipode is the identity, $S^2=id$.
\end{theorem}

\begin{proof}
This is something which is clear from Proposition 1.3, and from the properties of the functional transpose, with no computations needed. However, since the formalism of the functional transpose might be new to you, here is a detailed proof:

\medskip

(1) Let us first recall that, given a map between two sets $f:X\to Y$, its functional transpose is the morphism of algebras $f^t:F(Y)\to F(X)$ given by:
$$f^t(\varphi)=[x\to\varphi(f(x))]$$

To be more precise, this map $f^t$ is indeed well-defined, and the fact that we obtain in this way a morphism of algebras is clear. Indeed, for the addition, we have:
\begin{eqnarray*}
f^t(\varphi+\psi)
&=&[x\to(\varphi+\psi)(f(x))]\\
&=&[x\to(\varphi(f(x))+\psi(f(x)))]\\
&=&[x\to(\varphi(f(x))]+[x\to\psi(f(x))]\\
&=&f^t(\varphi)+f^t(\psi)
\end{eqnarray*}

As for the multiplication, the verification here is similar, as follows:
\begin{eqnarray*}
f^t(\varphi\psi)
&=&[x\to(\varphi\psi)(f(x))]\\
&=&[x\to(\varphi(f(x))\psi(f(x)))]\\
&=&[x\to(\varphi(f(x))]\cdot [x\to\psi(f(x))]\\
&=&f^t(\varphi)f^t(\psi)
\end{eqnarray*}

(2) Observe now that the operation $f\to f^t$ is by definition contravariant, in the sense that it reverses the direction of the arrows. Also, we have the following formula:
$$(fg)^t=g^tf^t$$

In order to check this, consider indeed two composable maps, as follows:
$$g:X\to Y\quad,\quad f:Y\to Z$$

The transpose of the composed map $fg:X\to Z$ is then given by:
\begin{eqnarray*}
(fg)^t(\varphi)
&=&[x\to\varphi(fg(x))]\\
&=&[x\to(\varphi f)(g(x))]\\
&=&[x\to(f^t(\varphi))(g(x))]\\
&=&g^tf^t(\varphi)
\end{eqnarray*}

Thus, we have indeed $(fg)^t=g^tf^t$, as claimed above. It is of course possible to use this iteratively, with the general formula, that we will often use, being as follows:
$$(f_1\ldots f_n)^t=f_n^t\ldots f_1^t$$

(3) As a second piece of preliminaries, we will need a bit of tensor product calculus too. In case you are familiar with this, say from physics classes, that is good news, and if not, here is a crash course on this. Let us start with something familiar, namely:
$$F^{M+N}=F^M\oplus F^N$$

As a consequence of this, that my students quite often tend to forget, we have:
$$F^{MN}\neq F^M\oplus F^N$$

And so, question for us now, what can be the new, mysterious operation $\otimes$, which is definitely not the direct sum $\oplus$, making the following formula work:
$$F^{MN}=F^M\otimes F^N$$

(4) In answer, let us look at the standard bases of these three vector spaces. We can denote by $\{f_1,\ldots,f_M\}$ the standard basis of $F^M$, and by $\{g_1,\ldots,g_N\}$ the standard basis of $F^N$. As for the space $F^{MN}$ on the left, here you would probably say to use the notation $\{e_1,\ldots,e_{MN}\}$, but I have here something better, namely $\{e_{11},\ldots,e_{MN}\}$, by using double indices. And, with this trick, the solution of our problem becomes clear, namely:
$$e_{ia}=f_i\otimes g_a$$

Thus, as a conclusion, given two vector spaces with bases $\{f_i\}$ and $\{g_a\}$, we can talk about their tensor product, as being the vector space with basis $\{f_i\otimes g_a\}$. And with this, we have the following formula, answering the question raised above: 
$$F^{MN}=F^M\otimes F^N$$

(5) As a continuation of this, in the case where our vector spaces are algebras $A,B$, their tensor product $A\otimes B$ is an algebra too, with multiplication as follows:
$$(a\otimes b)(a'\otimes b')=aa'\otimes bb'$$

Indeed, the algebra axioms are easily seen to be satisfied. And as a verification here, the above identification $F^{MN}=F^M\otimes F^N$ is indeed an algebra morphism, due to:
\begin{eqnarray*}
e_{ia}e_{jb}
&=&(f_i\otimes g_a)(f_j\otimes g_b)\\
&=&f_if_j\otimes g_ag_b\\
&=&\delta_{ij}f_i\otimes\delta_{ab}g_a\\
&=&\delta_{ia,jb}f_i\otimes g_a\\
&=&\delta_{ia,jb}e_{ia}
\end{eqnarray*}

(6) Moving ahead, still in relation with tensor products, we can say, more abstractly, that given two finite sets $X,Y$, we have an isomorphism of algebras as follows:
$$F(X\times Y)=F(X)\otimes F(Y)$$

To be more precise, we have a morphism from right to left, constructed as follows:
$$\varphi\otimes\psi\to[(x,y)\to\varphi(x)\psi(y)]$$

Now since this morphism is injective, and since the dimensions of the domain and range match, both being equal to $|X|\cdot|Y|$, this morphism is an isomorphism.

\medskip

(7) In what follows, the above formula from (6) will be what we will mostly need, in relation with the tensor products. As a last comment here, observe that this is in fact the formula $F^{MN}=F^M\otimes F^N$ that we started with. Indeed, when writing $X=\{1,\ldots,M\}$ and $Y=\{1,\ldots,N\}$, the above formula from (6) takes the following form:
$$F\left(\{1,\ldots,M\}\times\{1,\ldots,N\}\right)=F\left(\{1,\ldots,M\}\right)\otimes F\left(\{1,\ldots,N\}\right)$$

But this obviously translates into the following formula, as claimed:
$$F^{MN}=F^M\otimes F^N$$

And we will end with our preliminaries here. In case all this was not fully clear, may the Gods of Algebra forgive us, and we recommend on one hand reading about tensor products from a solid algebra book, such as Lang \cite{lng}, and on the other hand, doing some physics computations with tensor products, nothing can replace those either.

\medskip

(8) Still with me, I hope, and time now to prove our theorem? With the above ingredients in hand, let us go back indeed to our group theory problems. To start with, the structural maps $m,u,i$ of our group $G$ are maps as follows:
$$m:G\times G\to G\quad,\quad
u:\{.\}\to G\quad,\quad
i:G\to G$$

Thus, with $A=F(G)$ being the algebra of functions $\varphi:G\to F$, their functional transposes are morphisms of algebras $\Delta,\varepsilon,S$ as follows:
$$\Delta:A\to A\otimes A\quad,\quad 
\varepsilon:A\to F\quad,\quad 
S:A\to A$$

(9) Regarding now the formulae of these transposed maps, we know that the structural maps $m,u,i$ of our group $G$ are given by the following formulae:
$$m(g,h)=gh\quad,\quad 
u(.)=1\quad,\quad 
i(g)=g^{-1}$$

Thus, the functional transposes $\Delta,\varepsilon,S$ are given by the following formulae:
$$\Delta(\varphi)=[(g,h)\to \varphi(gh)]\quad,\quad
\varepsilon(\varphi)=\varphi(1)\quad,\quad
S(\varphi)=[g\to\varphi(g^{-1})]$$

(10) Regarding the group axioms, we know from Proposition 1.3 that in terms of the structural maps $m,u,i$, these are as follows, with $\delta(g)=(g,g)$:
$$m(m\times id)=m(id\times m)$$
$$m(u\times id)=m(id\times u)=id$$
$$m(i\times id)\delta=m(id\times i)\delta=1$$

In terms of the functional transposes $\Delta,\varepsilon,S$, these axioms read: 
$$(\Delta\otimes id)\Delta=(id\otimes \Delta)\Delta$$
$$(\varepsilon\otimes id)\Delta=(id\otimes\varepsilon)\Delta=id$$
$$m(S\otimes id)\Delta=m(id\otimes S)\Delta=\varepsilon(.)1$$

Finally, the formula $S^2=id$ comes by transposing $i^2=id$.
\end{proof}

Good news, Theorem 1.4 is all we need, or almost, in order to fulfill Goal 1.2 (1). Indeed, based on what we found above, we can formulate the following definition:

\index{comultiplication}
\index{counit}
\index{antipode}
\index{opposite algebra}
\index{Hopf algebra}
\index{finite quantum group}
\index{quantum group}
\index{deformed Hopf algebra}

\begin{definition}
A Hopf algebra is an algebra $A$, with morphisms of algebras
$$\Delta:A\to A\otimes A\quad,\quad 
\varepsilon:A\to F\quad,\quad 
S:A\to A^{opp}$$
called comultiplication, counit and antipode, satisfying the following conditions:
$$(\Delta\otimes id)\Delta=(id\otimes \Delta)\Delta$$
$$(\varepsilon\otimes id)\Delta=(id\otimes\varepsilon)\Delta=id$$
$$m(S\otimes id)\Delta=m(id\otimes S)\Delta=\varepsilon(.)1$$
If the square of the antipode is the identity, $S^2=id$, we say that $A$ is undeformed. Otherwise, in the case $S^2\neq id$, we say that $A$ is deformed.
\end{definition}

Here everything is standard, based on what we found in Theorem 1.4, we just copied the formulae there, with a different banner for them, except for two points, namely:

\bigskip

(1) We chose to have the antipode $S$ as being a morphism of algebras $S:A\to A^{opp}$, instead of being a morphism $S:A\to A$, as Theorem 1.4 might suggest. Indeed, since the algebra $A=F(G)$ in Theorem 1.4 is commutative, we have $A=A^{opp}$ in that case, so we can make this choice. And, we will see in a moment that $S:A\to A^{opp}$ is indeed the good choice, with this coming from some further examples, that we want our formalism to cover, and more specifically, coming from the algebras $F[G]$, from Goal 1.2 (2).

\bigskip

(2) We also chose the antipode $S$ not to be subject to the condition $S^2=id$. However, this is something debatable, because in the usual group setting $i^2=id$, while formally not being a group axiom, is something so trivial and familiar, that it is ``almost'' a group axiom. We will be back to this issue, on several occasions. In fact, clarifying the relation between Hopf algebras axiomatized with $S^2=id$, and Hopf algebras axiomatized without $S^2=id$, will be a main theme of discussion, throughout this book.

\bigskip

Finally, observe also that we chose not to impose any finite dimensionality condition on our Hopf algebra $A$, and this in contrast with Theorem 1.4, where the group $G$ there is finite. Again, this is something subtle, to be discussed more in detail later on.

\section*{1b. Basic examples}

There are several basic examples of Hopf algebras, which are all undeformed. We first have the following result, which provides a good motivation for our theory:

\index{commutative Hopf algebra}
\index{function algebra}

\begin{theorem}
Given a finite group $G$, the algebra of the $F$-valued functions on it, $F(G)=\{\varphi:G\to F\}$, with the usual pointwise product of functions,
$$(\varphi\psi)(g)=\varphi(g)\psi(g)$$
is a Hopf algebra, with comultiplication, counit and antipode as follows:
$$\Delta(\varphi)=[(g,h)\to \varphi(gh)]$$
$$\varepsilon(\varphi)=\varphi(1)$$
$$S(\varphi)=[g\to\varphi(g^{-1})]$$
This Hopf algebra is finite dimensional, commutative, and undeformed.
\end{theorem}

\begin{proof}
This is a reformulation of Theorem 1.4, by taking into account Definition 1.5, with the remark, already made above, that we have $A=A^{opp}$ in this case, due to the commutativity of $A=F(G)$, and with the last assertion being something clear.
\end{proof}

In view of the above result, we can make the following speculation:

\index{quantum group}
\index{finite quantum group}

\begin{speculation}
We can think of any finite dimensional Hopf algebra $A$ as being of the following form, with $G$ being a finite quantum group: 
$$A=F(G)$$
That is, we can define the category of finite quantum groups $G$ to be the category of finite dimensional Hopf algebras $A$, with the arrows reversed.
\end{speculation} 

Observe that, from the perspective of pure mathematics, all this is not that speculatory, because what we said in the end is something categorical and rigorous, perfectly making sense, and with the category of the usual finite groups $G$ embedding covariantly into the category of the finite quantum groups $G$, thanks to Theorem 1.6.

\bigskip

However, still mathematically speaking, there are some bugs with this. One problem is whether we want to include or not $S^2=id$ in our axioms, and in the lack of $S^2\neq id$ examples, at this stage of things, we are in the dark. Another problem is that, even when assuming $S^2=id$, nothing guarantees that a finite dimensional commutative Hopf algebra must be of the form $A=F(G)$, which would be something desirable to have.

\bigskip

As for the perspective brought by applied mathematics, here things are harsher, because the use of the word ``quantum'' would normally assume that our notion of Hopf algebra has something to do with quantum physics, and this is certainly not the case, now that we are into chapter 1 of the present book. Long way to go here, trust me.

\bigskip

In short, Speculation 1.7 remains a speculation, with our comments on it being:

\begin{comment}
The above $A=F(G)$ picture is something very useful, definitely worth to be kept in mind, but we will have to work some more on our axioms for Hopf algebras $A$, as for the corresponding objects $G$ to deserve the name ``quantum groups''.
\end{comment}

Finally, still in connection with all this, axiomatics, we would like if possible the construction in Theorem 1.6 to cover other groups as well, infinite this time, such as the discrete ones, or the compact ones, or, ideally, the locally compact ones. 

\bigskip

The problem with this, however, is that in the framework of Definition 1.5 this is not exactly possible, due to the fact that the comultiplication $\Delta$ would have to land in the algebra $F(G\times G)$, and for infinite groups $G,H$, we have:
$$F(G\times H)\neq F(G)\otimes F(H)$$

However, there are several tricks in order to overcome this, either by allowing $\otimes$ to be a topological tensor product, or by using Lie algebras. We will be back to this question, which is not trivial to solve, on several occassions, in what follows.

\bigskip

Moving ahead now, let us say that a Hopf algebra $A$ as axiomatized above is cocommutative if, with $\Sigma(a\otimes b)=b\otimes a$ being the flip map, we have the following formula:
$$\Sigma\Delta=\Delta$$

With this notion in hand, we have the following result, providing us with more examples, and that we will soon see to be ``dual'' to Theorem 1.6, in a suitable sense:

\index{group algebra}

\begin{theorem}
Given a group $H$, which can be finite or not, its group algebra
$$F[H]=span(H)$$
is a Hopf algebra, with structural maps given on group elements as follows:
$$\Delta(g)=g\otimes g\quad,\quad
\varepsilon(g)=1\quad,\quad
S(g)=g^{-1}$$
This Hopf algebra is cocommutative, and undeformed.
\end{theorem}

\begin{proof}
This is something elementary, the idea being as follows:

\medskip

(1) As a first observation, we can define indeed linear maps $\Delta,\varepsilon,S$ as in the statement, by linearity, and the maps $\Delta,\varepsilon$ are obviously morphisms of algebras. As for the antipode $S:A\to A^{opp}$, this is a morphism of algebras too, due to the following computation, crucially using the opposite multiplication $(a,b)\to a\cdot b$ on the target algebra:
\begin{eqnarray*}
S(gh)
&=&(gh)^{-1}\\
&=&h^{-1}g^{-1}\\
&=&g^{-1}\cdot h^{-1}\\
&=&S(g)\cdot S(h)
\end{eqnarray*}

(2) We have to verify now that $\Delta,\varepsilon,S$ satisfy the axioms in Definition 1.5, and the verifications here, performed on generators, are as follows:
$$(\Delta\otimes id)\Delta(g)=(id\otimes \Delta)\Delta(g)=g\otimes g\otimes g$$
$$(\varepsilon\otimes id)\Delta(g)=(id\otimes\varepsilon)\Delta(g)=g$$
$$m(S\otimes id)\Delta(g)=m(id\otimes S)\Delta(g)=\varepsilon(g)1=1$$

(3) Finally, it is clear from definitions that our Hopf algebra satisfies indeed the cocommutativity condition $\Sigma\Delta=\Delta$, as well as the condition $S^2=id$.
\end{proof}

The fact that the group $H$ in the above can be infinite comes as good news, and it is tempting to jump on this, and formulate, in analogy with Speculation 1.7:

\begin{speculation}
We can think of any Hopf algebra $A$ as being of the following form, with $H$ being a quantum group: 
$$A=F[H]$$
That is, we can define the category of quantum groups $H$ to be the category of Hopf algebras $A$.
\end{speculation}

However, as before with Speculation 1.7, while this being something useful, providing us with some intuition on what a Hopf algebra is, when looking more in detail at this, there are countless problems with it, which are both purely mathematical, of algebraic and analytic nature, and applied mathematical, in relation with quantum physics, which is certainly something more complicated than what we did in the above. 

\bigskip

Be said in passing, observe that, while both Speculation 1.7 and Speculation 1.10 formally make sense, from a pure mathematics perspective, their joint presence does not make much sense, mathematically, at least with the results that we have so far, because nothing guarantees that the category of finite quantum groups from Speculation 1.7 is indeed a subcategory of the category of quantum groups from Speculation 1.10.

\bigskip

Thinking a bit more at all this, we are led into the following question:

\begin{question}
What is the precise relation between Theorems 1.6 and 1.9, in the finite group case, and can this make peace between Speculations 1.7 and 1.10?
\end{question}

And the point now is that, despite its informal look, this question appears to be well-defined, and quite interesting, and answering it will be our next objective.

\bigskip

For this purpose, we first need to see what happens to Theorem 1.9, when assuming that the group $H$ there is finite. And here, we have the following statement:

\index{cocommutative Hopf algebra}
\index{cocommutative}
\index{convolution product}
\index{Dirac masses}

\begin{theorem}
Given a finite group $H$, the algebra of the $F$-valued functions on it $F[H]=\{\varphi:H\to F\}$, with the convolution product of functions,
$$(\varphi*\psi)(g)=\sum_{g=hk}\varphi(h)\psi(k)$$
is a Hopf algebra, with structural maps given on Dirac masses as follows:
$$\Delta(\delta_g)=\delta_g\otimes \delta_g\quad,\quad 
\varepsilon(\delta_g)=1\quad,\quad 
S(\delta_g)=\delta_{g^{-1}}$$
This Hopf algebra, which coincides with the previous $F[H]$, in the finite group case, is finite dimensional, cocommutative, and undeformed.
\end{theorem}

\begin{proof}
This is what comes from Theorem 1.9, when the group $H$ there is finite. Indeed, in this case the vector space $F[H]=span(H)$ from Theorem 1.9 coincides with the vector space $F[H]=\{\varphi:H\to F\}$ in the statement, with the correspondence being given on the standard generators $g\in span(H)$ by the following formula:
$$g\to\delta_g$$

Regarding now the product operation, the product on $F[H]=span(H)$ from Theorem 1.9 corresponds to the above convolution product on $F[H]=\{\varphi:H\to F\}$, because:
\begin{eqnarray*}
(\delta_r*\delta_s)(g)
&=&\sum_{g=hk}\delta_r(h)\delta_s(k)\\
&=&\delta_{g,rs}\\
&=&\delta_{rs}(g)
\end{eqnarray*}

Thus $\delta_r*\delta_s=\delta_{rs}$, as desired. We conclude that the algebra $F[H]$ from Theorem 1.9 coincides with the algebra $F[H]$ constructed here, and this gives the result.
\end{proof}

In practice now, the above statement has a weakness, coming from the fact that our formulae for $\Delta,\varepsilon,S$ are in terms of the Dirac masses. Here is a better version of it:

\begin{theorem}
Given a finite group $H$, the algebra of the $F$-valued functions on it $F[H]=\{\varphi:H\to F\}$, with the convolution product of functions,
$$(\varphi*\psi)(g)=\sum_{g=hk}\varphi(h)\psi(k)$$
is a Hopf algebra, with structural maps constructed as follows:
$$\Delta(\varphi)=[(g,h)\to\delta_{gh}\varphi(g)]$$
$$\varepsilon(\varphi)=\sum_{g\in H}\varphi(g)$$ 
$$S(\varphi)=[g\to\varphi(g^{-1})]$$
This Hopf algebra, which coincides with the previous $F[H]$, in the finite group case, is finite dimensional, cocommutative, and undeformed.
\end{theorem}

\begin{proof}
This is what comes from Theorem 1.12, by linearity. Indeed, according to our formula of $\Delta$ there, on the Dirac masses, we have:
\begin{eqnarray*}
\Delta(\varphi)(g,h)
&=&\Delta\left(\sum_{k\in H}\varphi(k)\delta_k\right)(g,h)\\
&=&\sum_{k\in H}\varphi(k)\Delta(\delta_k)(g,h)\\
&=&\sum_{k\in H}\varphi(k)(\delta_k\otimes\delta_k)(g,h)\\
&=&\sum_{k\in H}\varphi(k)\delta_{kg}\delta_{kh}\\
&=&\delta_{gh}\varphi(g)
\end{eqnarray*}

Also, according to our formula of $\varepsilon$ there, on the Dirac masses, we have:
\begin{eqnarray*}
\varepsilon(\varphi)
&=&\varepsilon\left(\sum_{g\in H}\varphi(g)\delta_g\right)\\
&=&\sum_{g\in H}\varphi(g)\varepsilon(\delta_g)\\
&=&\sum_{g\in H}\varphi(g)
\end{eqnarray*}

Finally, according to our formula of $S$ there, on the Dirac masses, we have:
\begin{eqnarray*}
S(\varphi)(g)
&=&S\left(\sum_{h\in H}\varphi(h)\delta_h\right)(g)\\
&=&\sum_{h\in H}\varphi(h)S(\delta_h)(g)\\
&=&\sum_{h\in H}\varphi(h)\delta_{h^{-1}}(g)\\
&=&\varphi(g^{-1})
\end{eqnarray*}

Thus, we are led to the conclusions in the statement.
\end{proof}

As a comment here, the proof of the above result relies on Theorem 1.12, which itself relies on Theorem 1.9, and in this type of situation, when things pile up, it is better to work out a new, direct proof, matter of doublechecking everything, and also matter of better understanding what is going on. So, let us do this, as an instructive exercise.

\bigskip

The problem is that of checking the Hopf algebra axioms directly, starting from the formulae of $\Delta,\varepsilon,S$ from Theorem 1.13, and we have here the following result:

\begin{proposition}
We have the following formulae over the algebra $F[H]$,
$$(\Delta\otimes id)\Delta(\varphi)=(id\otimes \Delta)\Delta(\varphi)$$
$$(\varepsilon\otimes id)\Delta(\varphi)=(id\otimes\varepsilon)\Delta(\varphi)=\varphi$$
$$m(S\otimes id)\Delta(\varphi)=m(id\otimes S)\Delta(\varphi)=\varepsilon(\varphi)1$$
guaranteeing that $F[H]$ is indeed a Hopf algebra, for any finite group $H$.
\end{proposition}

\begin{proof}
As mentioned, we already know this, as a consequence of Theorem 1.13, coming as consequence of Theorem 1.12, coming itself as consequence of Theorem 1.9, but time now to prove this directly as well, by using the formulae of $\Delta,\varepsilon,S$, namely:
$$\Delta(\varphi)=[(g,h)\to\delta_{gh}\varphi(g)]$$
$$\varepsilon(\varphi)=\sum_{g\in H}\varphi(g)$$ 
$$S(\varphi)=[g\to\varphi(g^{-1})]$$

We have to prove the following formulae, for any group elements $g,h,k$:
$$(\Delta\otimes id)\Delta(\varphi)(g,h,k)=(id\otimes \Delta)\Delta(\varphi)(g,h,k)$$
$$(\varepsilon\otimes id)\Delta(\varphi)(g)=(id\otimes\varepsilon)\Delta(\varphi)(g)=\varphi(g)$$
$$m(S\otimes id)\Delta(\varphi)(g)=m(id\otimes S)\Delta(\varphi)(g)=\varepsilon(\varphi)1$$

In what regards the first formula, this is clear, because the second iteration $\Delta^{(2)}$ of the comultiplication, no matter how computed, will by given the following formula:
$$\Delta^{(2)}(\varphi)(g,h,k)=\delta_{ghk}\varphi(g)$$

Regarding the second formula, this is again clear, because when applying either of the maps $E_1=\varepsilon\otimes id$ and $E_2=id\otimes\varepsilon$ to the quantity $\Delta(\varphi)$, what we get is:
$$E_i\Delta(\varphi)(g)=\sum_{g=h}\varphi(h)=\varphi(g)$$

As for the third formula, this is similar, because when applying either of the maps $T_1=m(S\otimes id)$ and $T_2=m(id\otimes S)$ to the quantity $\Delta(\varphi)$, what we get is:
$$T_i\Delta(\varphi)(g)=\sum_{g\in H}\varphi(g)1=\varepsilon(\varphi)1$$

Thus, we are led to the conclusions in the statement.
\end{proof}

\section*{1c. Abelian groups}

With the above done, let us try now to understand the relation between the algebras $F(G)$ from Theorem 1.6, and the algebras $F[H]$ from Theorem 1.12, or Theorem 1.13. 

\bigskip

For this purpose, we must first talk about abelian groups, and their duals. And as a starting point here, we first have the following elementary result:

\index{character}
\index{dual group}
\index{group of characters}

\begin{theorem}
Given a finite group $G$, the multiplicative characters 
$$\chi:G\to F^*$$
form a group $\widehat{G}$, called dual group, having the following properties:
\begin{enumerate}
\item $\widehat{G}$ is finite and abelian, and depends on both $G$ and $F$.

\item $\widehat{G}=\widehat{G}_{ab}$, where $G_{ab}=G/[G,G]$ is the abelianization of $G$.

\item We have a morphism $G\to\widehat{\widehat{G}}$, producing a morphism $G_{ab}\to\widehat{\widehat{G}}$.

\item The dual of a product is the product of duals, $\widehat{G\times H}=\widehat{G}\times\widehat{H}$.
\end{enumerate}
\end{theorem}

\begin{proof}
Our first claim is that $\widehat{G}$ as constructed above is indeed a group, with the pointwise multiplication of the characters, given by the following formula:
$$(\chi\rho)(g)=\chi(g)\rho(g)$$

Indeed, if $\chi,\rho$ are characters, so is $\chi\rho$, and so the multiplication is well-defined on $\widehat{G}$. Regarding the unit, this is the trivial character, constructed as follows:
$$1:G\to F^*\quad,\quad 
g\to1$$ 

Finally, we have inverses, with the inverse of $\chi:G\to F^*$ being as follows:
$$\chi^{-1}:G\to F^*\quad,\quad 
g\to\chi(g)^{-1}$$

Thus the dual group $\widehat{G}$ is indeed a group, and regarding now the other assertions:

\medskip

(1) We have several things to be proved here, the idea being as follows:

\medskip

-- Our first claim is that $\widehat{G}$ is finite. Indeed, given a group element $g\in G$, we can talk about its order, which is smallest integer $N\in\mathbb N$ such that $g^N=1$. Now assuming that we have a character $\chi:G\to F^*$, we must have the following formula:
$$\chi(g)^N=1$$

Thus $\chi(g)$ must be one of the $N$-th roots of unity inside $F$, that is, must be one of the roots over $F$ of the polynomial $X^N-1$, and in particular, there are finitely many choices for $\chi(g)$. Thus, there are finitely many choices for $\chi$, and so $\widehat{G}$ is finite, as claimed.

\medskip

-- Next, the fact that $\widehat{G}$ is abelian follows from definitions, because the pointwise multiplication of functions, and in particular of characters, is commutative.

\medskip

-- Finally, the group dual $\widehat{G}$ as constructed above certainly depends on $G$, but the point is that it can depend on the ground field $F$ too. Indeed, for an illustration here, consider the cyclic group $G=\mathbb Z_N$. A character $\chi:\mathbb Z_N\to F^*$ is then uniquely determined by its value $z=\chi(g)$ on the standard generator $g\in\mathbb Z_N$, and this value must satisfy:
$$z^N=1$$

Now over the complex numbers, $F=\mathbb C$, the solutions here are the usual $N$-th roots of unity, so we have $|\widehat{\mathbb Z}_N|=N$. Moreover, by thinking a bit, we have in fact:
$$\widehat{\mathbb Z}_N=\mathbb Z_N$$

In contrast, over the real numbers, $F=\mathbb R$, the possible solutions of $z^N=1$ must be among $z=\pm1$, and we conclude that in this case, the dual is given by:
$$\widehat{\mathbb Z}_N=\begin{cases}
\{1\}&(N\ {\rm odd})\\
\mathbb Z_2&(N\ {\rm even})
\end{cases}$$

There are of course far more things that can be said here, with all this being related to the structure of the group of $N$-th roots of unity inside $F$, but for our present purposes, what we have so far, namely the above illustration using $F=\mathbb R,\mathbb C$, will do.

\medskip

(2) Let us prove now the second assertion, $\widehat{G}=\widehat{G}_{ab}$. We recall that given a group $G$, its commutator subgroup $[G,G]\subset G$ is constructed as follows:
$$[G,G]=\left\{ghg^{-1}h^{-1}\Big|g,h\in G\right\}$$

This subgroup is then normal, and we can define the abelianization of $G$ as being:
$$G_{ab}=G/[G,G]$$

Now given a character $\chi:G\to F^*$, we have the following formula, for any $g,h\in G$, based of course on the fact that the multiplicative group $F^*$ is abelian:
$$\chi(ghg^{-1}h^{-1})=1$$

Thus, our character factorizes as follows, into a character of the group $G_{ab}$:
$$\chi:G\to G_{ab}\to F^*$$

Summarizing, we have constructed an identification $\widehat{G}=\widehat{G}_{ab}$, as claimed.

\medskip

(3) Regarding now the third assertion, as a first remark, given a finite group $G$ we have indeed a morphism $I:G\to\widehat{\widehat{G}}$, appearing via evaluation maps, as follows:
$$I:G\to\widehat{\widehat{G}}\quad,\quad g\to[\chi\to\chi(g)]$$

Since the group duals, and so the group double duals too, are always abelian, we cannot expect $I$ to be injective, in general. In fact, due to $\chi(ghg^{-1}h^{-1})=1$, we have:
$$[G,G]\subset\ker I$$

Thus, the morphism $I:G\to\widehat{\widehat{G}}$ constructed above factorizes as follows:
$$I:G\to G_{ab}\to \widehat{\widehat{G}}$$

And with this, third assertion proved. There are of course some further things that can be said here, in relation with the factorized morphism $G_{ab}\to \widehat{\widehat{G}}$, but the examples in (1) above, where $G$ itself was abelian, show that, when using an arbitrary field $F$, as we are currently doing, this factorized morphism is not an isomorphism, in general.

\medskip

(4) Finally, in what regards the fourth assertion, $\widehat{G\times H}=\widehat{G}\times\widehat{H}$, observe that a character of a product of groups $\chi:G\times H\to F^*$ must satisfy:
\begin{eqnarray*}
\chi(g,h)
&=&\chi\left[(g,1)(1,h)\right]\\
&=&\chi(g,1)\chi(1,h)\\
&=&\chi_{|G}(g)\chi_{|H}(h)
\end{eqnarray*}

Thus $\chi$ must appear as the product of its restrictions $\chi_{|G},\chi_{|H}$, which must be both characters, and this gives the identification in the statement.
\end{proof}

As a conclusion to what we have above, certainly some interesting things, but with the overall situation being not very good, due to the fact that the group dual $\widehat{G}$ depends on the ground field $F$, and with this preventing many interesting things to happen.

\bigskip

So, good time to temporarily break with our policy so far of using an arbitrary field $F$. By choosing our ground field to be the smartest one around, namely $F=\mathbb C$, and by assuming as well that $G$ is abelian, in view of what we found above, we are led to:

\index{cyclic group}
\index{product of cyclic groups}
\index{self-dual group}

\begin{theorem}
Given a finite abelian group $G$, its complex characters 
$$\chi:G\to\mathbb T$$
form a finite abelian group $\widehat{G}$, called Pontrjagin dual, and the following happen:
\begin{enumerate}
\item The dual of a cyclic group is the group itself, $\widehat{\mathbb Z}_N=\mathbb Z_N$.

\item The dual of a product is the product of duals, $\widehat{G\times H}=\widehat{G}\times\widehat{H}$.

\item Any product of cyclic groups $G=\mathbb Z_{N_1}\times\ldots\times\mathbb Z_{N_k}$ is self-dual, $G=\widehat{G}$.
\end{enumerate}
\end{theorem}

\begin{proof}
We already know some of these things from Theorem 1.15, but since we are here simplifying that, the best is to start all over again. We have indeed a finite abelian group, as stated, and regarding the other assertions, their proof goes as follows:

\medskip

(1) A character $\chi:\mathbb Z_N\to\mathbb T$ is uniquely determined by its value $z=\chi(g)$ on the standard generator $g\in\mathbb Z_N$. But this value must satisfy:
$$z^N=1$$

Thus we must have $z\in\mathbb Z_N$, with the cyclic group $\mathbb Z_N$ being regarded this time as being the group of $N$-th roots of unity. Now conversely, any $N$-th root of unity $z\in\mathbb Z_N$ defines a character $\chi:\mathbb Z_N\to\mathbb T$, by setting, for any $r\in\mathbb N$:
$$\chi(g^r)=z^r$$

Thus we have an identification $\widehat{\mathbb Z}_N=\mathbb Z_N$, as claimed.

\medskip

(2) A character of a product of groups $\chi:G\times H\to\mathbb T$ must satisfy:
$$\chi(g,h)=\chi\left[(g,1)(1,h)\right]=\chi(g,1)\chi(1,h)$$

Thus $\chi$ must appear as the product of its restrictions $\chi_{|G},\chi_{|H}$, which must be both characters, and this gives the identification in the statement.

\medskip

(3) This follows from (1) and (2). Alternatively, any character $\chi:G\to\mathbb T$ is uniquely determined by its values $\chi(g_1),\ldots,\chi(g_k)$ on the standard generators of $\mathbb Z_{N_1},\ldots,\mathbb Z_{N_k}$, which must belong to $\mathbb Z_{N_1},\ldots,\mathbb Z_{N_k}\subset\mathbb T$, and this gives $\widehat{G}=G$, as claimed.
\end{proof}

Let us discuss now the structure result for the finite abelian groups. This is something which is more advanced, requiring good knowledge of group theory, as follows:

\index{finite abelian group}

\begin{theorem}
The finite abelian groups are the following groups,
$$G=\mathbb Z_{N_1}\times\ldots\times\mathbb Z_{N_k}$$
and these groups are all self-dual, $G=\widehat{G}$.
\end{theorem}

\begin{proof}
This is something quite tricky, the idea being as follows:

\medskip

(1) In order to prove our result, assume that $G$ is finite and abelian. For any prime number $p\in\mathbb N$, let us define $G_p\subset G$ to be the subset of elements having as order a power of $p$. Equivalently, this subset $G_p\subset G$ can be defined as follows:
$$G_p=\left\{g\in G\Big|\exists k\in\mathbb N,g^{p^k}=1\right\}$$

(2) It is then routine to check, based on definitions, that each $G_p$ is a subgroup. Our claim now is that we have a direct product decomposition as follows:
$$G=\prod_pG_p$$

(3) Indeed, by using the fact that our group $G$ is abelian, we have a morphism as follows, with the order of the factors when computing $\prod_pg_p$ being irrelevant:
$$\prod_pG_p\to G\quad,\quad (g_p)\to\prod_pg_p$$

Moreover, it is routine to check that this morphism is both injective and surjective, via some simple manipulations, so we have our group decomposition, as in (2).

\medskip

(4) Thus, we are left with proving that each component $G_p$ decomposes as a product of cyclic groups, having as orders powers of $p$, as follows:
$$G_p=\mathbb Z_{p^{r_1}}\times\ldots\times\mathbb Z_{p^{r_s}}$$

But this is something that can be checked by recurrence on $|G_p|$, via some routine computations, and so we are led to the conclusions in the statement.

\medskip

(5) So, this was for the quick story of the present result, structure theorem for the finite abelian groups, and for more on all this, technical details, and some useful generalizations too, we recommend learning this from a solid algebra book, such as Lang \cite{lng}.
\end{proof}

Getting back now to the Hopf algebras, we have the following result:

\index{dual group}
\index{finite abelian group}
\index{Pontrjagin duality}
\index{Fourier transform}
\index{discrete Fourier transform}

\begin{theorem}
If $G,H$ are finite abelian groups, dual to each other via Pontrjagin duality, in the sense that each of them is the character group of the other,
$$G=\big\{\chi:H\to\mathbb T\big\}\quad,\quad H=\big\{\rho:G\to\mathbb T\big\}$$
we have an identification of Hopf algebras as follows:
$$F(G)=F[H]$$
In the case $G=H=\mathbb Z_N$, this identification is the usual discrete Fourier transform isomorphism. In general, we obtain a tensor product of such Fourier transforms.
\end{theorem}

\begin{proof}
All this is standard Fourier analysis, the idea being as follows:

\medskip

(1) In the simplest case, where our groups are $G=H=\mathbb Z_N$, we have indeed an identification of algebras as above, which is a Hopf algebra isomorphism, given by the usual discrete Fourier transform isomorphism, whose matrix with respect to the standard basis on each side is the following matrix, with $w=e^{2\pi i/N}$, called Fourier matrix:
$$F_N=
\begin{pmatrix}
1&1&1&\ldots&1\\
1&w&w^2&\ldots&w^{N-1}\\
1&w^2&w^4&\ldots&w^{2(N-1)}\\
\vdots&\vdots&\vdots&&\vdots\\
1&w^{N-1}&w^{2(N-1)}&\ldots&w^{(N-1)^2}
\end{pmatrix}$$

(2) In the general case now, we can invoke the stucture theorem for the abelian groups, which tells us that $G$ must appear as a product of cyclic groups, as follows:
$$G=\mathbb Z_{N_1}\times\ldots\times\mathbb Z_{N_k}$$

Indeed, due to the functorial properties of the Pontrjagin duality, we have as well:
$$H=\mathbb Z_{N_1}\times\ldots\times\mathbb Z_{N_k}$$

Thus, we are led to the conclusions in the statement, with the corresponding isomorphism of Hopf algebras being obtained by tensoring the isomorphisms from (1), and corresponding to the following matrix, called generalized Fourier matrix:
$$F_G=F_{N_1}\otimes\ldots\otimes F_{N_k}$$

(3) Alternatively, it is possible to be more direct on all this, and short-circuiting the heavy results, simply by viewing the identification $F(G)=F[H]$ as appearing by complexifying the characters from the definition of the dual group, namely:
$$H=\big\{\chi:G\to\mathbb T\big\}\quad {\rm or}\quad 
G=\big\{\chi:H\to\mathbb T\big\}$$

Indeed, with this approach, which relies only on the definition of the Pontrjagin dual, there is no need for the computations for $\mathbb Z_N$, or for the structure theorem for the finite abelian groups. Again, we will leave the details here as an instructive exercise.
\end{proof}

As a comment here, we can feel that Theorem 1.18 is related to Fourier analysis, and this is indeed the case. The point is that we have 3 types of Fourier analysis in life:

\medskip

(1) We first have the ``standard'' one, corresponding to $G=\mathbb R$, that you probably know well, and which can be learned from any advanced analysis book.

\medskip

(2) Then we have another one, called the ``Fourier series'' one, which is also something popular and useful, corresponding to $G=\mathbb Z,\mathbb T$, that you probably know well too.

\medskip

(3) And finally we have the ``discrete'' one, as above, over $G=\mathbb Z_N$ and other finite abelian groups. We will be back to this, on several occasions, in this book.

\section*{1d. Duality theory} 

Quite remarkably, the Pontrjagin duality for finite abelian groups can be extended to the general finite group case, in the context of the Hopf algebras. To be more precise, we have the following result, which is something truly remarkable, solving many questions, and which will be our first general theorem on the Hopf algebras:

\index{dual Hopf algebra}
\index{duality of Hopf algebras}
\index{cocommutative algebra}

\begin{theorem}
Given a finite dimensional Hopf algebra $A$, its dual space
$$A^*=\Big\{\varphi:A\to F\ {\rm linear}\Big\}$$
is also a finite dimensional Hopf algebra, with multiplication and unit as follows,
$$\Delta^t:A^*\otimes A^*\to A^*\quad,\quad \varepsilon^t:F\to A^*$$
and with comultiplication, counit and antipode as follows:
$$m^t:A^*\to A^*\otimes A^*\quad,\quad
u^t:A^*\to F\quad,\quad 
S^t:A^*\to A^*$$
This duality makes correspond the commutative algebras to the cocommutative algebras. Also, this duality makes correspond $F(G)$ to $F[G]$, for any finite group $G$.
\end{theorem}

\begin{proof}
At the first glance, we can only expect here something more complicated than for Theorem 1.18, that our result generalizes. However, by the power of abstract algebra, where precise formulations matter a lot, things are in fact quite simple:

\medskip

(1) To start with, we know that $A$ is a Hopf algebra. Thus, as an associative algebra, $A$ has a multiplication map $m$, and a unit map $u$, which are as follows:
$$m:A\otimes A\to A\quad,\quad u:F\to A$$

Also, $A$ has a comultiplication $\Delta$, counit $\varepsilon$ and antipode $S$, which are as follows:
$$\Delta:A\to A\otimes A\quad,\quad\varepsilon:A\to F\quad,\quad S:A\to A^{opp}$$

(2) By taking now the functional transposes of these 5 maps, we obtain 5 other maps, whose domains and ranges are as in the statement. Moreover, it is routine to check that these latter 5 maps are all morphisms of algebras, with this being actually clear for all the maps involved, except for $S^t$, which requires some thinking at opposite algebras.

\medskip

(3) Regarding now the axioms, since $A$ is, before anything, an associative algebra, its multiplication and unit maps $m,u$ are subject to the following axioms:
$$m(m\otimes id)=m(id\otimes m)$$
$$m(u\otimes id)=m(id\otimes u)=id$$

We also know that $A$ is a Hopf algebra, so the following are satisfied too:
$$(\Delta\otimes id)\Delta=(id\otimes \Delta)\Delta$$
$$(\varepsilon\otimes id)\Delta=(id\otimes\varepsilon)\Delta=id$$
$$m(S\otimes id)\Delta=m(id\otimes S)\Delta=\varepsilon(.)1$$

(4) The point now is that the collection of the above 8 formulae is ``self-dual'', in the sense that when transposing, we obtain exactly the same 8 formulae. Indeed, the transposes of the first two formulae are as follows:
$$(m^t\otimes id)m^t=(id\otimes m^t)m^t$$
$$(u^t\otimes id)m^t=(id\otimes u^t)m^t=id$$

As for the transposes of the last three formulae, these are as follows:
$$\Delta^t(\Delta^t\otimes id)=\Delta^t(id\otimes\Delta^t)$$
$$\Delta^t(\varepsilon^t\otimes id)=\Delta^t(id\otimes\varepsilon^t)=id$$
$$\Delta^t(S^t\otimes id)m^t=\Delta^t(id\otimes S^t)m^t=u^t(.)1$$

But, we recognize here the full axioms for Hopf algebras, including those for associative algebras. Thus $A^*$, as constructed in the statement, is indeed a Hopf algebra. 

\medskip

(5) Observe now, as a complement to what is said in the statement, and which is something that is useful to know, in practice, that the operation $A\to A^*$ is indeed a duality, because if we dualize one more time, we obtain $A$ itself:
$$A^{**}=A$$

(6) Regarding the assertion about commutative and cocommutative algebras, this is clear from definitions, because we have the following equivalences:
\begin{eqnarray*}
\Sigma m^t=m^t
&\iff&\Sigma m^t(\varphi)=m^t(\varphi),\forall\varphi\\
&\iff&\Sigma m^t(\varphi)(a,b)=m^t(\varphi)(a,b),\forall\varphi,a,b\\
&\iff&\varphi(ba)=\varphi(ab),\forall\varphi,a,b\\
&\iff&ba=ab,\forall a,b
\end{eqnarray*}

Indeed, this computation gives the result in one sense, and in the other sense, this follows either via a similar computation, or just by dualizing, and using (5).

\medskip

(7) Finally, the last assertion, regarding the group algebras, is clear from definitions too, after a quick comparison with Theorem 1.6 and Theorem 1.13. 

\medskip

(8) Indeed, the point is that we have dual vector spaces, and the Hopf algebra maps from Theorem 1.6 are given by the following formulae:
$$(\varphi\psi)(g)=\varphi(g)\psi(g)$$
$$1=[g\to1]$$
$$\Delta(\varphi)=[(g,h)\to \varphi(gh)]$$
$$\varepsilon(\varphi)=\varphi(1)$$
$$S(\varphi)=[g\to\varphi(g^{-1})]$$

(9) As for the Hopf algebra maps from Theorem 1.13, these are as follows:
$$(\varphi*\psi)(g)=\sum_{g=hk}\varphi(h)\psi(k)$$
$$1=[g\to\delta_{g1}]$$
$$\Delta(\varphi)=[(g,h)\to\delta_{gh}\varphi(g)]$$
$$\varepsilon(\varphi)=\sum_{g\in H}\varphi(g)$$ 
$$S(\varphi)=[g\to\varphi(g^{-1})]$$

Thus, we have indeed a pair of dual Hopf algebras, as stated.

\medskip

(10) As a last observation, in the case where the finite group $G$ is abelian, we recover in this way what we know from Theorem 1.18, proved there the hard way. But, as mentioned there, there are in fact several proofs for that result, with going via the heavy theorems for the finite abelian groups being something practical, but not really necessary.
\end{proof}

As a conclusion to this, we can now answer Question 1.11 in the affirmative, by merging Speculation 1.7 and Speculation 1.10 in the finite case, in the following way:

\begin{speculation}
We can think of any finite dimensional Hopf algebra $A$, not necessarily commutative or cocommutative, as being of the form
$$A=F(G)=F[H]$$
with $G,H$ being finite quantum groups, related by a generalized Pontrjagin duality. And with this generalizing what we know about the abelian groups.
\end{speculation}

All this is very nice, and we will leave some further categorical thinking and clarification of all this, if needed, depending on taste, as an exercise.

\bigskip

As a second comment, despite the above being something nice, and quite deep, the criticisms formulated on the occasion of Speculation 1.7 and Speculation 1.10 remain. One problem, as usual, is whether we want to further include the condition $S^2=id$ in our axioms for the quantum groups. Another problem is that, even when assuming $S^2=id$, nothing guarantees that a finite dimensional commutative Hopf algebra must be of the form $A=F(G)$. And also, once again when assuming $S^2=id$, nothing guarantees that a finite dimensional cocommutative Hopf algebra must be of the form $A=F[H]$. And finally, we have the big question regarding the relation of all this, mathematics developed in a few pages, with quantum physics, which is certainly something more complicated. All these are good problems, and we will be back to them, in what follows.

\bigskip

As a third comment, all the above concerns the finite dimensional case, and when trying to do such things in the infinite dimensional setting, which must be topological, as per the usual Pontrjagin duality for infinite groups, there are plenty of difficulties. Again, these are all good questions, and we will be back to them, later in this book.

\section*{1e. Exercises}

We had a lot of interesting algebra in this chapter, sometimes going towards basic functional analysis, or differential geometry, and as exercises, we have:

\begin{exercise}
Learn everything about tensor products, from a good algebra book.
\end{exercise}

\begin{exercise}
Do as well some physics computations with $\otimes$, matter of loving it.
\end{exercise}

\begin{exercise}
Clarify all details of the functional transpose operation $f\to f^t$.
\end{exercise}

\begin{exercise}
Learn everything about abelian groups, and Pontrjagin duality.
\end{exercise}

\begin{exercise}
Learn about the various types of Fourier transforms and series.
\end{exercise}

\begin{exercise}
Work out all details for the duality theorem for Hopf algebras.
\end{exercise}

As bonus exercise, learn some functional analysis, which is obviously related to all this. The more you will know here, in advance, the better that will be.

\chapter{Basic theory}

\section*{2a. The antipode}

With the Hopf algebras axiomatized, the basic examples $F(G)$ and $F[H]$ discussed, and the duality theory discussed too, in the finite dimensional case, what is next? Many things, and as a list of pressing topics to be discussed, for this chapter, we have:

\bigskip

(1) We would first like to have a more detailed look at the Hopf algebra axioms, and what can be done with them. And notably, know more about the antipode $S$.

\bigskip

(2) Next, at the level of the basic examples, we have unfinished business, or rather unstarted business, with the enveloping Lie algebras $U\mathfrak g$. On our to-do list, too.

\bigskip

(3) Then, with the knowledge of $F(G)$, $F[H]$, $U\mathfrak g$, we can try to emulate, as a continuation of (1), some advanced group technology, inside the arbitrary Hopf algebras $A$.

\bigskip

(4) And finally, again with some inspiration from $F(G)$, $F[H]$, $U\mathfrak g$, which will be our main input here, we can have a discussion about Haar integration.

\bigskip

Which sounds quite good, and with none of the above questions (1,2,3,4) being easy to solve, via some instant thinking, expect lots of new and interesting things, to come.

\bigskip

Getting started, as a first topic for this chapter, let us go back to the arbitrary Hopf algebras, as axiomatized in chapter 1, and have a more detailed look at their antipode $S$. The definition and basic properties of the antipode can be summarized as follows:

\index{antipode}

\begin{theorem}
Given a Hopf algebra $A$, its antipode is the morphism of algebras
$$S:A\to A^{opp}$$
$A^{opp}$ being the opposite algebra, with product $a\cdot b=ba$, subject to the following axiom:
$$m(S\otimes id)\Delta=m(id\otimes S)\Delta=\varepsilon(.)1$$
For $F(G)$ the antipode is the transpose of the inversion map $i:G\to G$. For $F[H]$, the antipode is given by $S(g)=g^{-1}$. In both these cases, the above axiom corresponds to
$$g^{-1}g=gg^{-1}=1$$
and the extra condition $S^2=id$, coming from $(g^{-1})^{-1}=g$, is satisfied.
\end{theorem}

\begin{proof}
This is something that we know well from chapter 1, and for full details about all this, along with slightly more about $S$, we refer to the material there.
\end{proof}

In relation with this, as a first question that you might have, are there Hopf algebras with $S^2\neq id$? Here is a key example, due to Sweedler, which is the simplest one:

\index{Sweedler algebra}
\index{square of antipode}

\begin{theorem}
The Sweedler algebra, $A=span(1,c,x,cx)$ with the relations
$$c^2=1\quad,\quad x^2=0\quad,\quad cx=-xc$$
and with Hopf algebra structure given by
$$\Delta(c)=c\otimes c\quad,\quad\Delta(x)=1\otimes x+x\otimes c\quad,\quad\Delta(cx)=c\otimes cx+cx\otimes1$$
$$\varepsilon(c)=1\quad,\quad\varepsilon(x)=0\quad,\quad\varepsilon(cx)=0$$
$$S(c)=c\quad,\quad S(x)=cx\quad,\quad S(cx)=-x$$
is not commutative, nor cocommutative, and has the property $S^4=id$, but $S^2\neq id$.
\end{theorem}

\begin{proof}
This is something quite tricky, the idea being as follows:

\medskip

(1) Consider the 4-dimensional vector space $A=span(1,c,x,cx)$, with $1,c,x,cx$ being some abstract variables. We can make then $A$ into an associative algebra, with unit 1, by declaring that we have indeed $c\cdot x=cx$, and by imposing the following rules:
$$c^2=1\quad,\quad x^2=0\quad,\quad cx=-xc$$

(2) Next, by using the universal property of $A$, we can define a morphism of algebras $\Delta:A\to A\otimes A$, according to the following formulae, on the algebra generators $c,x$:
$$\Delta(c)=c\otimes c\quad,\quad\Delta(x)=1\otimes x+x\otimes c$$
 
Observe that by multiplying we have as well $\Delta(cx)=c\otimes cx+cx\otimes1$.

\medskip

(3) Now let us try to prove that $A$ is a Hopf algebra. By using the Hopf algebra axioms, we conclude that $\varepsilon,S$ can only be given on $c,x$ by the following formulae: 
$$\varepsilon(c)=1\quad,\quad\varepsilon(x)=0$$
$$S(c)=c\quad,\quad S(x)=cx$$

But, and here comes the point, we can define indeed such morphisms, $\varepsilon:A\to F$ and $S:A\to A^{opp}$, via the above formulae, by using the universality property of $A$. Observe that by multiplying, we obtain as well $\varepsilon(cx)=0$ and $S(cx)=-x$.

\medskip

(4) Summarizing, we have our Hopf algebra, which is clearly not commutative, not cocommutative either, and whose antipode satisfies $S^4=id$, but $S^2\neq id$.
\end{proof}

Getting back now to the general case, in order to further build on Theorem 2.1, our main source of inspiration will be what happens for $A=F(G)$, where the antipode appears as the functional analytic transpose, $S=i^t$, of the inversion map $i(g)=g^{-1}$. 

\bigskip

In view of this, we have the following natural question, which appears:

\begin{question}
In group theory we have many elementary formulae involving products, units and inverses, all coming from the group axioms, such as 
$$1^{-1}=1\quad,\quad (g^{-1})^{-1}=g\quad,\quad (gh)^{-1}=h^{-1}g^{-1}\quad,\quad\ldots$$
which can be reformulated in terms of $m,u,i$, and then transposed, leading to formulae as follows, involving $\Delta,\varepsilon,S$, inside the algebras $A=F(G)$:
$$\varepsilon S=\varepsilon\quad,\quad S^2=id\quad,\quad \Delta S=\Sigma(S\otimes S)\Delta\quad,\quad\ldots$$
But, which of these latter formulae hold in general, inside any Hopf algebra $A$?
\end{question}

And isn't this a good question. Indeed, as we know well from chapter 1, the answer to the above question is definitely ``yes'' for the group axioms themselves, which are as follows, and which reformulate, by definition, into the Hopf algebra axioms:
$$(gh)k=g(hk)\quad,\quad 1g=g1=g\quad,\quad g^{-1}g=gg^{-1}=1$$

However, in what regards for instance the group theory formula $(g^{-1})^{-1}=g$, this reformulates in Hopf algebra terms as $S^2=id$, and we have seen in Theorem 2.2 that we have counterexamples to this. Thus, Question 2.3 definitely makes sense.

\bigskip

In order to get familiar with this, let us first study $1^{-1}=1$. We have here:

\begin{theorem}
We have the following formula, valid over any Hopf algebra,
$$\varepsilon S=\varepsilon$$
and with this coming from $1^{-1}=1$ for $A=F(G)$, and being trivial for $A=F[H]$.
\end{theorem}

\begin{proof}
This is something elementary, the idea being as follows:

\medskip

(1) In order to establish the above formula, we can use the Hopf algebra axiom for $S$. Indeed, by applying the counit to this axiom, we obtain the following formula:
$$\varepsilon m(S\otimes id)\Delta=\varepsilon m(id\otimes S)\Delta=\varepsilon$$

Let us compute the map on the left. By using the counit axiom, we have:
\begin{eqnarray*}
\varepsilon m(S\otimes id)\Delta
&=&(\varepsilon\otimes\varepsilon)(S\otimes id)\Delta\\
&=&(\varepsilon S\otimes\varepsilon)\Delta\\
&=&\varepsilon S(id\otimes\varepsilon)\Delta\\
&=&\varepsilon S\circ id\\
&=&\varepsilon S
\end{eqnarray*}

Similarly, and although not needed here, the map appearing in the middle in the above formula is $\varepsilon S$ too. Thus, our above equality of maps reads $\varepsilon S=\varepsilon S=\varepsilon$, as desired.

\medskip

(2) Regarding now the case $A=F(G)$, here our condition is as follows, as claimed:
\begin{eqnarray*}
\varepsilon S=\varepsilon
&\iff&\varepsilon S(\varphi)=\varepsilon(\varphi)\\
&\iff&\varphi(1^{-1})=\varphi(1)\\
&\iff&1^{-1}=1
\end{eqnarray*}

(3) On the opposite, in the case $A=F[H]$ our condition is trivial, coming from:
\begin{eqnarray*}
\varepsilon S=\varepsilon
&\iff&\varepsilon S(g)=\varepsilon(g)\\
&\iff&\varepsilon(g^{-1})=\varepsilon(g)\\
&\iff&1=1
\end{eqnarray*}

(4) Thus, result proved, with all this being quite trivial, but with the remark that what we have at the end, that $F(G)$ vs $F[H]$ dissymmetry, is however quite interesting.
\end{proof}

Next on our list, still coming from Question 2.3, let us have a look at $(gh)^{-1}=h^{-1}g^{-1}$. Things here are more tricky, and as a first result on the subject, we have:

\begin{proposition}
The antipode of both the algebras $F(G)$ and $F[H]$ satisfies
$$\Delta S=\Sigma(S\otimes S)\Delta$$
with this coming from $(gh)^{-1}=h^{-1}g^{-1}$ for $A=F(G)$, and being trivial for $A=F[H]$.
\end{proposition}

\begin{proof}
For $A=F(G)$ the proof goes as follows, with $\sigma(g,h)=(h,g)$:
\begin{eqnarray*}
(gh)^{-1}=h^{-1}g^{-1}
&\iff&im(g,h)=m(i\times i)\sigma(g,h)\\
&\iff&im=m(i\times i)\sigma\\
&\iff&m^ti^t=\sigma^t(i^t\otimes i^t)m^t\\
&\iff&\Delta S=\Sigma(S\otimes S)\Delta
\end{eqnarray*}

As for the algebra $A=F[H]$, here the verification is trivial, as follows:
\begin{eqnarray*}
\Delta S=\Sigma(S\otimes S)\Delta
&\iff&\Delta S(g)=\Sigma(S\otimes S)\Delta(g)\\
&\iff&\Delta(g^{-1})=\Sigma(S\otimes S)(g\otimes g)\\
&\iff&\Delta(g^{-1})=\Sigma(g^{-1}\otimes g^{-1})\\
&\iff&\Delta(g^{-1})=g^{-1}\otimes g^{-1}
\end{eqnarray*}

Observe that, again, we have an interesting $F(G)$ vs $F[H]$ dissymmetry here.
\end{proof}

Before getting further with our study of $(gh)^{-1}=h^{-1}g^{-1}$, that we will eventually show to hold over any Hopf algebra, let us go back to Question 2.3. We have only recorded 3 possible relations there, but there are infinitely many more, and in relation with this, we have the following result, related to $(gh)^{-1}=h^{-1}g^{-1}$, making a bit of cleanup:

\begin{proposition}
We have the following implication, over any Hopf algebra,
$$\Delta S=\Sigma(S\otimes S)\Delta\quad\implies\quad\Delta^{(2)}S=\Sigma^{(2)}(S\otimes S\otimes S)\Delta^{(2)}$$
corresponding in the case $A=F(G)$ to the following implication, at the group level:
$$(gh)^{-1}=h^{-1}g^{-1}\quad\implies\quad (ghk)^{-1}=k^{-1}h^{-1}g^{-1}$$
Moreover, we can iterate this observation, as many times as we want to.
\end{proposition}

\begin{proof}
We are presently into uncharted territory, so take the above statement as it is, as something a bit informal. The point with all this is that, in the group theory setting, it is quite obvious that by iterating, we have the following series of implications:
\begin{eqnarray*}
(gh)^{-1}=h^{-1}g^{-1}
&\implies&(ghk)^{-1}=k^{-1}h^{-1}g^{-1}\\
&\implies&(ghkl)^{-1}=l^{-1}k^{-1}h^{-1}g^{-1}\\
&\implies&\ldots
\end{eqnarray*}

And the question is that, can we have these implications going, I mean these implications only, not the validity of the formulae themselves, in the Hopf algebra setting. With the answer to this latter question being definitely yes, for the first implication, with the computation being as follows, using twice the condition $\Delta S=\Sigma(S\otimes S)\Delta$:
\begin{eqnarray*}
\Delta^{(2)}S
&=&(\Delta\otimes id)\Delta S\\
&=&(\Delta\otimes id)\Sigma(S\otimes S)\Delta\\
&=&\Sigma_\leftarrow(id\otimes\Delta)(S\otimes S)\Delta\\
&=&\Sigma_\leftarrow(S\otimes\Delta S)\Delta\\
&=&\Sigma_\leftarrow(S\otimes\Sigma(S\otimes S)\Delta)\Delta\\
&=&\Sigma_\leftarrow(id\otimes\Sigma)(S\otimes S\otimes S)(id\otimes \Delta)\Delta\\
&=&\Sigma^{(2)}(S\otimes S\otimes S)\Delta^{(2)}
\end{eqnarray*}

Thus, result proved, and I will leave it to you to figure out what the various versions of $\Sigma,\Delta$ used above exactly mean. In what regards the last assertion, exercise as well.
\end{proof}

Summarizing, things are quite tricky. In order to further discuss this, we will need some abstract algebraic preliminaries. Let us start with something standard, namely:

\index{convolution}

\begin{theorem}
If we define the convolution of linear maps $\varphi,\psi:A\to A$ by
$$\varphi*\psi=m(\varphi\otimes\psi)\Delta$$
then the Hopf algebra axiom for the antipode reads
$$S*id=id*S=\varepsilon(.)1$$
with the map on the right, $\varepsilon(.)1$, being the unit for the operation $*$.
\end{theorem}

\begin{proof}
This is something which comes from the axioms, as follows:

\medskip

(1) The first assertion is clear from the Hopf algebra axiom for the antipode, as formulated in chapter 1, or in Theorem 2.1, which was as follows:
$$m(S\otimes id)\Delta=m(id\otimes S)\Delta=\varepsilon(.)1$$

Indeed, in terms of the convolution operation from the statement, $\varphi*\psi=m(\varphi\otimes\psi)\Delta$, this axiom takes the following more conceptual form, as indicated above:
$$S*id=id*S=\varepsilon(.)1$$

(2) Regarding the second assertion, this follows from the counit axiom, namely:
$$(\varepsilon\otimes id)\Delta=(id\otimes\varepsilon)\Delta=id$$

Indeed, given a linear map $\varphi:A\to A$, we have the following computation:
\begin{eqnarray*}
\varphi*\varepsilon(.)1
&=&m(\varphi\otimes\varepsilon(.)1)\Delta\\
&=&\varphi(id\otimes\varepsilon)\Delta\\
&=&\varphi\circ id\\
&=&\varphi
\end{eqnarray*}

Similarly, again for any linear map $\varphi:A\to A$, we have the following computation:
\begin{eqnarray*}
\varepsilon(.)1*\varphi
&=&m(\varepsilon(.)1\otimes\varphi)\Delta\\
&=&\varphi(\varepsilon\otimes id)\Delta\\
&=&\varphi\circ id\\
&=&\varphi
\end{eqnarray*}

Thus the linear map $\varepsilon(.)1$ is indeed the unit for the operation $*$, as claimed.
\end{proof}

In order to do our next antipode computations, which will be sometimes quite tough, we will need as well the following useful convention, due to Sweedler:

\index{Sweedler notation}

\begin{definition}
We use the Sweedler notation for the comultiplication $\Delta$,
$$\Delta(x)=\sum x_1\otimes x_2$$
with the sum on the right being understood to correspond to the tensor expansion of $\Delta(x)$.
\end{definition}

And in the hope that this will sound quite nice and clever to you, when seeing it for the first time. As illustrations for this Sweedler notation, or rather as a first piece of advertisement for it, let us have a look at the Hopf algebra axioms, namely:
$$(\Delta\otimes id)\Delta=(id\otimes \Delta)\Delta$$
$$(\varepsilon\otimes id)\Delta=(id\otimes\varepsilon)\Delta=id$$
$$m(S\otimes id)\Delta=m(id\otimes S)\Delta=\varepsilon(.)1$$

We know these axioms since the beginning of chapter 1, and we certainly have some knowledge in dealing with them. However, the point is that by using the Sweedler notation above, these axioms take the following form, which is even more digest:
$$\sum\Delta(x_1)\otimes x_2=\sum x_1\otimes \Delta(x_2)$$ 
$$\sum\varepsilon(x_1)x_2=\sum x_1\varepsilon(x_2)=x$$ 
$$\sum S(x_1)x_2=\sum x_1S(x_2)=\varepsilon(x)$$

With this discussed, let us go back now to Theorem 2.7, which is our main theoretical result, so far. With a bit more work, we can further improve this result, as follows:

\begin{theorem}
Given an algebra $A$, with morphisms of algebras
$$\Delta:A\to A\otimes A\quad,\quad 
\varepsilon:A\to F$$
satisfying the usual axioms for a comultiplication and an antipode, namely
$$(\Delta\otimes id)\Delta=(id\otimes \Delta)\Delta$$
$$(\varepsilon\otimes id)\Delta=(id\otimes\varepsilon)\Delta=id$$
this is a Hopf algebra precisely when $id:A\to A$ is invertible with respect to 
$$\varphi*\psi=m(\varphi\otimes\psi)\Delta$$
and in this case, the convolution inverse $S=id^{-1}$ is the antipode.
\end{theorem}

\begin{proof}
This follows indeed from what we have in Theorem 2.7, and from a few extra computations, best done by using the Sweedler notation, the idea being as follows:

\medskip

(1) Assume first that $A$ is a Hopf algebra. According to Theorem 2.7 we have indeed $S=id^{-1}$, inverse with respect to convolution, so the only thing that we have to prove is that this inverse is unique. But this is something purely algebraic, which is valid under very general circumstances, because for any associative multiplication $\cdot$ we have:
$$ab=ba=1,\ ac=ca=1\quad\implies\quad b=bac=c$$ 

(2) Conversely now, assume that the identity map $id:A\to A$ is invertible with respect to the convolution operation $*$, with inverse $S=id^{-1}$. As explained in Theorem 2.7, the condition $S=id^{-1}$ tells us that $S$ satisfies the usual antipode axiom, namely:
$$m(S\otimes id)\Delta=m(id\otimes S)\Delta=\varepsilon(.)1$$

However, we are not done yet, because our map $S:A\to A$ is just a linear map, that we still have to prove to be a morphism, when regarded as map $S:A\to A^{opp}$. Thus, with $S=id^{-1}$ regarded as it comes, as linear map $S:A\to A$, we must prove that we have:
$$S(ab)=S(b)S(a)$$

(3) In order to prove this formula, consider the following three linear maps:
$$m(a\otimes b)=ab\quad,\quad p(a\otimes b)=S(ab)\quad,\quad q(a\otimes b)=S(b)S(a)$$

We have then the following computation, involving the maps $m,p$:
\begin{eqnarray*}
(p*m)(a\otimes b)
&=&\sum p\left((a\otimes b)_1\right)m\left((a\otimes b)_2\right)\\
&=&\sum p(a_1\otimes b_1)m(a_2\otimes b_2)\\
&=&\sum S(a_1b_1)a_2b_2\\
&=&\sum S\left((ab)_1\right)(ab)_2\\
&=&(S*1)(ab)\\
&=&\varepsilon(ab)\\
&=&\varepsilon(a)\varepsilon(b)\\
&=&(\varepsilon\otimes\varepsilon)(a\otimes b)
\end{eqnarray*}

On the other hand, we have as well the following computation, involving $m,q$:
\begin{eqnarray*}
(m*q)(a\otimes b)
&=&\sum m(a_1\otimes b_1)q(a_2\otimes b_2)\\
&=&\sum a_1b_1S(b_2)S(a_2)\\
&=&\sum a_1\varepsilon(b)S(a_2)\\
&=&(1*S)(a)\cdot\varepsilon(b)\\
&=&\varepsilon(a)\varepsilon(b)\\
&=&(\varepsilon\otimes\varepsilon)(a\otimes b)
\end{eqnarray*}

Summarizing, we have proved that we have the following formulae:
$$p*m=m*q=(\varepsilon\otimes\varepsilon)(.)1\otimes1$$

(4) But with this, we can finish our proof, in the following way:
\begin{eqnarray*}
p
&=&p*((\varepsilon\otimes\varepsilon)(.)1\otimes1)\\
&=&p*(m*q)\\
&=&(p*m)*q\\
&=&((\varepsilon\otimes\varepsilon)(.)1\otimes1)*q\\
&=&q
\end{eqnarray*}

Thus we have $p=q$, which means $S(ab)=S(b)S(a)$, as desired.
\end{proof}

With all these preliminaries discussed, time now for our first true theorem. We can indeed formulate something non-trivial regarding the antipode, as follows:

\begin{theorem}
The antipode of a Hopf algebra $S:A\to A^{opp}$ satisfies:
\begin{enumerate}
\item $\varepsilon S=\varepsilon$.

\item $\Delta S=\Sigma(S\otimes S)\Delta$.

\item $S^2=id$, when $A$ is commutative or cocommutative.
\end{enumerate}
\end{theorem}

\begin{proof}
This is something quite tricky, the idea being as follows:

\medskip

(1) This is something that we already know, from Theorem 2.4.

\medskip

(2) We have, by using the Sweedler notation from Definition 2.8:
\begin{eqnarray*}
(\Delta S*\Delta)(x)
&=&\sum\Delta S(x_1)\Delta(x_2)\\
&=&\Delta\left(\sum S(x_1)x_2\right)\\
&=&\Delta(\varepsilon(x)1)\\
&=&\varepsilon(x)\cdot1\otimes1
\end{eqnarray*}

On the other hand, we have as well, by using the Sweedler notation, iterated:
\begin{eqnarray*}
(\Delta*\Sigma(S\otimes S)\Delta)(x)
&=&\sum(x_1\otimes x_2)(S(x_4)\otimes S(x_3))\\
&=&\sum x_1S(x_4)\otimes x_2S(x_3)\\
&=&\sum x_1S(x_3)\otimes\varepsilon(x_2)1\\
&=&\sum x_1S(x_3)\varepsilon(x_2)\otimes1\\
&=&\sum x_1S(x_2)\otimes1\\
&=&\varepsilon(x)\cdot1\otimes1
\end{eqnarray*}

As a conclusion to this, we have proved the following equalities:
$$\Delta S*\Delta=\Delta*\Sigma(S\otimes S)\Delta=\varepsilon(.)1\otimes1$$

Now by using Theorem 2.7, we obtain from this, as desired:
\begin{eqnarray*}
\Delta S
&=&\Delta S*(\varepsilon(.)1\otimes1)\\
&=&\Delta S*(\Delta*\Sigma(S\otimes S)\Delta)\\
&=&(\Delta S*\Delta)*\Sigma(S\otimes S)\Delta\\
&=&(\varepsilon(.)1\otimes1)*\Sigma(S\otimes S)\Delta\\
&=&\Sigma(S\otimes S)\Delta
\end{eqnarray*}

(3) Our first claim is that when $A$ is commutative or cocommutative, we have:
$$\sum S(x_2)x_1=\varepsilon(x)$$

Indeed, in the commutative case this follows from the Hopf algebra axiom for $S$, in Sweedler notation, which reads, as explained before:
$$\sum x_1S(x_2)=\varepsilon(x)$$

As for the cocommutative case, here we can use again the axiom for $S$, namely:
$$m(S\otimes id)\Delta(x)=\varepsilon(x)1$$

Indeed, by replacing $\Delta(x)$ with $\Sigma\Delta(x)$, for flipping the tensors, we obtain the formula claimed in the above. Thus, claim proved, and with this in hand, we have:
\begin{eqnarray*}
m(S\otimes S^2)\Delta(x)
&=&m(S\otimes S^2)\left(\sum x_1\otimes x_2\right)\\
&=&\sum S(x_1)S^2(x_2)\\
&=&S\left(\sum S(x_2)x_1\right)\\
&=&S(\varepsilon(x)1)\\
&=&\varepsilon(x)1
\end{eqnarray*}

Now by using Theorem 2.9, we obtain from this $S^2=id$, as claimed.
\end{proof}

As a first comment on the above result, (2) there can be used in conjunction with Proposition 2.6, and shows that much more is true. However, Question 2.3 as formulated is still there, and knowing what group theory type relations hold, and what don't, inside an arbitrary Hopf algebra, remains something which requires some experience and skill.

\bigskip

Many other things can be said, as a continuation of the above, notably with some general theory for the square of the antipode $S^2:A\to A$, which is more specialized, again based on the general interpretation of the antipode coming from Theorem 2.9. For more on all this, we refer to the specialized Hopf algebra literature.

\section*{2b. Lie algebras}

As a second task for this chapter, again coming as a continuation of what we did in chapter 1, let us try now to better understand what happens, at the fine level, beyond the finite dimensional case, where we already have some good examples and results. 

\bigskip

As already mentioned in chapter 1, when trying to cover various infinite groups, such as the compact, discrete, or more generally locally compact ones, the standard trick is that of modifying a bit the Hopf algebra axioms, by using topological tensor products. But this is something quite technical, and we will discuss this later in this book. 

\bigskip

For the moment, let us show that we can cover the compact Lie groups, without changing our axioms, by using a Lie algebra trick. Our claim is as follows:

\begin{claim}
Given a compact Lie group $G$, with Lie algebra $\mathfrak g$, the corresponding enveloping Lie algebra $U\mathfrak g$ is a Hopf algebra, which is cocommutative.
\end{claim}

Obviously, many non-trivial notions involved here, with this being not your routine abstract algebra statement, that you can understand right away, armed with Love for algebra only. So, we will explain in what follows the various notions involved, namely the Lie groups $G$, the Lie algebras $\mathfrak g$, and the enveloping Lie algebras $U\mathfrak g$, and then we will get of course to what the above claim says, and even provide a proof for it.

\bigskip

Getting started now, a Lie group is by definition a group which is a smooth manifold. So, let us start our discussion with this, smooth manifolds. Here is their definition:

\begin{definition}
A smooth manifold is a space $X$ which is locally isomorphic to $\mathbb R^N$. To be more precise, this space $X$ must be covered by charts, bijectively mapping open pieces of it to open pieces of $\mathbb R^N$, with the changes of charts being $C^\infty$ functions.
\end{definition}

As basic examples of smooth manifolds, we have of course $\mathbb R^N$ itself, or any open subset $X\subset\mathbb R^N$, with only 1 chart being needed here. Other basic examples, in the plane, at $N=2$, include the circle, or various curves like ellipses and so on, somehow for obvious reasons. Here is a more precise statement in this sense, covering the conics:

\index{deformation of manifold}

\begin{proposition}
The following are smooth manifolds, in the plane:
\begin{enumerate}
\item The circles.

\item The ellipses.

\item The non-degenerate conics.

\item Smooth deformations of these.
\end{enumerate}
\end{proposition}

\begin{proof}
All this is quite intuitive, the idea being as follows:

\medskip

(1) Consider the unit circle, $x^2+y^2=1$. We can write then $x=\cos t$, $y=\sin t$, with $t\in[0,2\pi)$, and we seem to have here the solution to our problem, just using 1 chart. But this is of course wrong, because $[0,2\pi)$ is not open, and we have a problem at $0$. In practice we need to use 2 such charts, say with the first one being with $t\in(0,3\pi/2)$, and the second one being with $t\in(\pi,5\pi/2)$. As for the fact that the change of charts is indeed smooth, this comes by writing down the formulae, or just thinking a bit, and arguing that this change of chart being actually a translation, it is automatically linear.

\medskip

(2) This follows from (1), by pulling the circle in both the $Ox$ and $Oy$ directions, and the  formulae here, based on the standard formulae for ellipses, are left to you reader.

\medskip

(3) We already have the ellipses, and the case of the parabolas and hyperbolas is elementary as well, and in fact simpler than the case of the ellipses. Indeed, a parablola is clearly homeomorphic to $\mathbb R$, and a hyperbola, to two copies of $\mathbb R$.

\medskip

(4) This is something which is clear too, depending of course on what exactly we mean by ``smooth deformation'', and by using a bit of multivariable calculus if needed.
\end{proof}

In higher dimensions we have as basic examples the spheres, and I will leave it to you to find a proof, using spherical coordinates, or the stereographic projection. Exercise as well, to find higher dimensional analogues of the other assertions in Proposition 2.13.

\bigskip

Getting now to what we wanted to do here, Lie groups, let us start with:

\index{Lie group}
\index{Lie algebra}

\begin{definition}
A Lie group is a group $G$ which is a smooth manifold, with the corresponding multiplication and inverse maps
$$m:G\times G\to G\quad,\quad i:G\to G$$
being assumed to be smooth. The tangent space at the origin $1\in G$ is denoted
$$\mathfrak g=T_1G$$
and is called Lie algebra of $G$.
\end{definition}

So, this is our definition, and as a first observation, the examples of Lie groups abound, with the circle $\mathbb T$ and with the higher dimensional tori $\mathbb T^N$ being the standard examples. For these, the Lie algebra is obviously equal to $\mathbb R$ and $\mathbb R^N$, respectively. There are of course many other examples, all very interesting, and more on this in a moment.

\bigskip

Before getting into examples, let us discuss a basic question, that you surely have in mind, namely why calling the tangent space $\mathfrak g=T_1G$ an algebra. In answer, since $G$ is a group, with a certain multiplication map $m:G\times G\to G$, we can normally expect this map $m$ to produce some sort of ``algebra structure'' on the tangent space $\mathfrak g=T_1G$.

\bigskip

This was for the general idea, but in practice, things are more complicated than this, because even for very simple examples of Lie groups, what we get in this way is not an associative algebra, but rather a new type of beast, called Lie algebra. 

\bigskip

So, coming as a continuation and complement to Definition 2.14, we have:

\index{Lie algebra}
\index{Lie bracket}
\index{Jacobi identity}

\begin{definition}
A Lie algebra is a vector space $\mathfrak g$ with an operation $(x,y)\to[x,y]$, called Lie bracket, subject to the following conditions:
\begin{enumerate}
\item $[x+y,z]=[x,z]+[y,z]$, $[x,y+z]=[x,y]+[x,z]$.

\item $[\lambda x,y]=[x,\lambda y]=\lambda[x,y]$.

\item $[x,x]=0$.

\item $[[x,y],z]+[[y,z],x]+[[z,x],y]=0$.
\end{enumerate}
\end{definition}

As a basic example, consider a usual, associative algebra $A$. We can define then the Lie bracket on it as being the usual commutator, namely:
$$[x,y]=xy-yx$$

The above axioms (1,2,3) are then clearly satisfied, and in what regards axiom (4), called Jacobi identity, this is satisfied too, the verification being as follows:
\begin{eqnarray*}
&&[[x,y],z]+[[y,z],x]+[[z,x],y]\\
&=&[xy-yx,z]+[yz-zy,x]+[zx-xz,y]\\
&=&xyz-yxz-zxy+zyx+yzx-zyx-xyz+xzy+zxy-xzy-yzx+yxz\\
&=&0
\end{eqnarray*}

We will see in a moment that up to a certain abstract operation $\mathfrak g\to U\mathfrak g$, called enveloping Lie algebra construction, and which is something quite elementary, any Lie algebra appears in this way, with its Lie bracket being formally given by:
$$[x,y]=xy-yx$$

Before that, however, you might wonder where that Gothic letter $\mathfrak g$ in Definition 2.15 comes from. That comes from the following fundamental result, making the connection with the theory of Lie groups from Definition 2.14, denoted as usual by $G$:

\begin{theorem}
Given a Lie group $G$, that is, a group which is a smooth manifold, with the group operations being smooth, the tangent space at the identity
$$\mathfrak g=T_1(G)$$
is a Lie algebra, with its Lie bracket being basically a usual commutator.
\end{theorem}

\begin{proof}
This is something non-trivial, the idea being as follows:

\medskip

(1) Let us first have a look at the orthogonal and unitary groups $O_N,N_N$. These are both Lie groups, and the corresponding Lie algebras $\mathfrak o_N,\mathfrak u_N$ can be computed by differentiating the equations defining $O_N,U_N$, with the conclusion being as follows:
$$\mathfrak o_N=\left\{ A\in M_N(\mathbb R)\Big|A^t=-A\right\}$$
$$\mathfrak u_N=\left\{ B\in M_N(\mathbb C)\Big|B^*=-B\right\}$$

This was for the correspondences $O_N\to\mathfrak o_N$ and $U_N\to\mathfrak u_N$. In the other sense, the correspondences $\mathfrak o_N\to O_N$ and $\mathfrak u_N\to U_N$ appear by exponentiation, the result here stating that, around 1, the orthogonal matrices can be written as $U=e^A$, with $A\in\mathfrak o_N$, and the unitary matrices can be written as $U=e^B$, with $B\in\mathfrak u_N$.

\medskip

(2) Getting now to the Lie bracket, the first observation is that both $\mathfrak o_N,\mathfrak u_N$ are stable under the usual commutator of the $N\times N$ matrices. Indeed, assuming that $A,B\in M_N(\mathbb R)$ satisfy $A^t=-A$, $B^t=-B$, their commutator satisfies $[A,B]\in M_N(\mathbb R)$, and:
\begin{eqnarray*}
[A,B]^t
&=&(AB-BA)^t\\
&=&B^tA^t-A^tB^t\\
&=&BA-AB\\
&=&-[A,B]
\end{eqnarray*}

Similarly, assuming that $A,B\in M_N(\mathbb C)$ satisfy $A^*=-A$, $B^*=-B$, their commutator $[A,B]\in M_N(\mathbb C)$ satisfies the condition $[A,B]^*=-[A,B]$.

\medskip

(3) We conclude from this discussion that both the tangent spaces $\mathfrak o_N,\mathfrak u_N$ are Lie algebras, with the Lie bracket being the usual commutator of the $N\times N$ matrices. 

\medskip

(4) It remains now to understand how the Lie bracket $[A,B]=AB-BA$ is related to the group commutator $[U,V]=UVU^{-1}V^{-1}$ via the exponentiation map $U=e^A$, and this can be indeed done, by making use of the differential geometry of $O_N,U_N$, and the situation is quite similar when dealing with an arbitrary Lie group $G$.
\end{proof}

With this understood, let us go back to the arbitrary Lie algebras, as axiomatized in Definition 2.15. There is an obvious analogy there with the axioms for the usual, associative algebras, and based on this analogy, we can build some abstract algebra theory for the Lie algebras. Let us record some basic results, along these lines:

\begin{proposition}
Let $\mathfrak g$ be a Lie algebra. If we define its ideals as being the vector spaces $\mathfrak i\subset\mathfrak g$ satisfying the condition
$$x\in\mathfrak i,y\in\mathfrak g\implies[x,y]\in\mathfrak i$$
then the quotients $\mathfrak g/\mathfrak i$ are Lie algebras. Also, given a morphism of Lie algebras $f:\mathfrak g\to\mathfrak h$, its kernel $ker(f)\subset\mathfrak g$ is an ideal, and we have $\mathfrak g/\ker(f)=Im(f)$.
\end{proposition}

\begin{proof}
All this is very standard, exactly as in the case of the associative algebras, and we will leave the various verifications here as an instructive exercise.
\end{proof}

Getting now to the point, remember our claim from the discussion after Definition 2.15, stating that up to a certain abstract operation $\mathfrak g\to U\mathfrak g$, called enveloping Lie algebra construction, any Lie algebra appears in fact from the ``trivial'' associative algebra construction, that is, with its Lie bracket being formally a usual commutator:
$$[x,y]=xy-yx$$

Time now to clarify this. The result here, making as well to the link with the various Lie group considerations from Theorem 2.16 and its proof, is as follows:

\index{enveloping Lie algebra}
\index{tensor algebra}

\begin{theorem}
Given a Lie algebra $\mathfrak g$, define its enveloping Lie algebra $U\mathfrak g$ as being the quotient of the tensor algebra of $\mathfrak g$, namely
$$T(\mathfrak g)=\bigoplus_{k=0}^\infty\mathfrak g^{\otimes k}$$
by the following associative algebra ideal, with $x,y$ ranging over the elements of $\mathfrak g$:
$$I=<x\otimes y-y\otimes x-[x,y]>$$
Then $U\mathfrak g$ is an associative algebra, so it is a Lie algebra too, with bracket 
$$[x,y]=xy-yx$$
and the standard embedding $\mathfrak g\subset U\mathfrak g$ is a Lie algebra embedding.\end{theorem}

\begin{proof}
This is something which is quite self-explanatory, and in what regards the examples, illustrations, and other things that can be said, for instance in relation with the Lie groups, we will leave some further reading here as an instructive exercise.
\end{proof}

Importantly, the above enveloping Lie algebra construction makes the link with our Hopf algebra considerations, from the present book, via the following result:

\begin{theorem}
Given a Lie algebra $\mathfrak g$, its enveloping Lie algebra $U\mathfrak g$ is a cocommutative Hopf algebra, with comultiplication, counit and antipode given by
$$\Delta:U\mathfrak g\to U(\mathfrak g\oplus\mathfrak g)=U\mathfrak g\otimes U\mathfrak g\quad,\quad x\to x+x$$
$$\varepsilon:U\mathfrak g\to F\quad,\quad x\to 1$$
$$S:U\mathfrak g\to U\mathfrak g^{opp}=(U\mathfrak g)^{opp}\quad,\quad x\to-x$$
via various standard identifications, for the various associative algebras involved.
\end{theorem}

\begin{proof}
Again, this is something quite self-explanatory, and in what regards the examples, illustrations, and other things that can be said, for instance in relation with the Lie groups, we will leave some further reading here as an instructive exercise.
\end{proof}

We will be back to this, and to Lie algebras in general, on several occasions, in what follows. Among others, we will see later in this book how to reconstruct the Lie group $G$ from the knowledge of the enveloping Lie algebra $U\mathfrak g$, using representation theory.

\section*{2c. Special elements}

In view of the above results regarding the enveloping Lie algebras $U\mathfrak g$, which are cocommutative, and of the results from chapter 1 too, regarding the group algebras $F[H]$, which are cocommutative too, it makes sense to have a more systematic look at the Hopf algebras $A$ which are cocommutative, in our usual sense, namely:
$$\Sigma\Delta=\Delta$$

We already know a bit about such algebras, in the finite dimensional case, and as a complement to that material, we first have the following result:

\index{group-like element}

\begin{theorem}
Let $A$ be a Hopf algebra. The elements satisfying the condition
$$\Delta(a)=a\otimes a$$
which are called group-like, have the following properties:
\begin{enumerate}
\item They form a group $G_A$.

\item They satisfy $\Sigma\Delta(a)=\Delta(a)$.

\item We have a Hopf algebra embedding $F[G_A]\subset A$.

\item For a group algebra $A=F[H]$, this embedding is an isomorphism.
\end{enumerate}
\end{theorem}

\begin{proof}
This is something elementary, the idea being as follows:

\medskip

(1) Let us call indeed group-like the elements $a\in A$ satisfying $\Delta(a)=a\otimes a$, in analogy with the formula $\Delta(g)=g\otimes g$ when $A=F[H]$, and more on this in a moment. The group-like elements are then stable under the product operation, as shown by:
$$\Delta(ab)
=\Delta(a)\Delta(b)
=(a\otimes a)(b\otimes b)
=ab\otimes ab$$

We have as well the stability under taking inverses, with this coming from:
$$\Delta(a^{-1})
=\Delta(a)^{-1}
=(a\otimes a)^{-1}
=a^{-1}\otimes a^{-1}$$

Finally, the formula $\Delta(1)=1\otimes1$ shows that $1\in A$ is group-like. Thus, the set of group-like elements $G_A\subset A$ is indeed a multiplicative subgroup, as claimed.

\medskip

(2) Assuming that $a\in A$ is group-like, we have indeed, as claimed:
$$\Sigma\Delta(a)
=\Sigma(a\otimes a)
=a\otimes a
=\Delta(a)$$

Observe that the converse of what we just proved here does not hold, for instance because in the case of the group algebras $A=F[H]$, which are cocommutative, $\Sigma\Delta=\Delta$, there are many elements $a\in A$ which are not group-like. More on this in a moment.

\medskip

(3) There are several checks here, all being trivial or routine, the only point being that of proving that the group-like elements are linearly independent. So, let us prove this. Assume that we have a linear combination of group-like elements, as follows:
$$a=\sum_i\lambda_ia_i$$

By applying $\Delta$ to this element we obtain, by using the condition $a\in G_A$:
$$\Delta(a)
=a\otimes a
=\sum_{ij}\lambda_i\lambda_ja_i\otimes a_j$$

On the other hand, also by applying $\Delta$, but by using $a_i\in G_A$, we obtain:
$$\Delta(a)
=\sum_i\lambda_i\Delta(a_i)
=\sum_i\lambda_ia_i\otimes a_i$$

We conclude that the following equation must be satisfied, when $a\in G_A$:
$$\sum_{ij}\lambda_i\lambda_ja_i\otimes a_j=\sum_i\lambda_ia_i\otimes a_i$$

But this equation shows that we must have $\#\{i\}=\#\{j\}=1$, as desired. That is, we have proved that the group-like elements are linearly independent, and this gives a Hopf algebra embedding $F[G_A]\subset A$ as in the statement, appearing in the obvious way.

\medskip

(4) This is indeed clear from (3), because in the case of a group algebra $A=F[H]$ we have $G_A=H$, with $G_A\supset H$ being clear, and with $G_A\subset H$ coming from linear independence. Thus, in this case, the embedding $F[G_A]\subset A$ is an isomorphism.
\end{proof}

Many other things can be said about the group-like elements, and we will leave their study in the case algebras $A=F(G)$, and of the algebras $A=\mathfrak g$ too, as an instructive exercise. Moving on, here is another key construction, this time Lie algebra-inspired:

\index{primitive element}

\begin{theorem}
Let $A$ be a Hopf algebra. The elements satisfying the condition
$$\Delta(a)=a\otimes 1+1\otimes a$$
which are called primitive, have the following properties:
\begin{enumerate}
\item They form a Lie algebra $P_A$, with bracket $[a,b]=ab-ba$.

\item They automatically satisfy $\Sigma\Delta(a)=\Delta(a)$.

\item We have a Hopf algebra embedding $UP_A\subset A$.

\item For an enveloping Lie algebra $A=U\mathfrak g$, this embedding is an isomorphism.
\end{enumerate}
\end{theorem}

\begin{proof}
Observe the similarity with Theorem 2.20, and more on this later. Regarding now the proof of the various assertions, this is straightforward, as follows:

\medskip

(1) There are several things to be checked here, all being trivial or routine, the only point being that of proving that the primitive elements are stable under taking commutators. So, let us prove this. Assuming $a,b\in P_A$, we have the following computation:
\begin{eqnarray*}
\Delta([a,b])
&=&\Delta a\cdot\Delta b-\Delta b\cdot\Delta a\\
&=&(a\otimes 1+1\otimes a)(b\otimes 1+1\otimes b)-(b\otimes 1+1\otimes b)(a\otimes 1+1\otimes a)\\
&=&ab\otimes 1-ba\otimes 1+1\otimes ab-1\otimes ba\\
&=&[a,b]\otimes 1+1\otimes[a,b]
\end{eqnarray*}

Thus we have $[a,b]\in P_A$, as desired, and this gives the assertion.

\medskip

(2) Assuming that $a\in A$ is primitive, we have indeed, as claimed:
\begin{eqnarray*}
\Sigma\Delta(a)
&=&\Sigma(a\otimes 1+1\otimes a)\\
&=&1\otimes a+a\otimes 1\\
&=&a\otimes 1+1\otimes a\\
&=&\Delta(a)
\end{eqnarray*}

Observe that the converse of what we just proved does not hold, for instance because in the case of the envelopoing Lie algebras $A=U\mathfrak g$, which are cocommutative, $\Sigma\Delta=\Delta$, there are many elements $a\in A$ which are not primitive. More on this in a moment.

\medskip

(3) This indeed something quite routine, a bit as before, for the group-likes.

\medskip

(4) This follows indeed from (3), again a bit as before, for the group-likes.
\end{proof}

Many other things can be said about the primitive elements, and we will leave their study in the case of the function algebras $A=F(G)$, and of the group algebras $A=F[H]$ too, which is something quite routine, as an instructive exercise.

\bigskip

Along the same lines, as a third and last construction now, motivated this time by the function algebras $A=F(G)$, which are commutative, we have:

\index{central element}
\index{infinite conjugacy class}
\index{ICC property}
\index{central functions}

\begin{theorem}
Given a Hopf algebra $A$, we can talk about its center
$$Z(A)\subset A$$
which is an associative subalgebra, having the following properties:
\begin{enumerate}
\item For $A=F(G)$, and more generally when $A$ is commutative, $Z(A)=A$.

\item For $A=F[H]$, this algebra $Z(A)$ is the algebra of central functions on $H$.

\item In particular, when all conjugacy classes of $H$ are infinite, $Z(A)=F$.
\end{enumerate}
\end{theorem}

\begin{proof}
This is something quite self-explanatory, and a bit in analogy with Theorem 2.20, and Theorem 2.21. Consider indeed the central elements of a Hopf algebra $A$:
$$Z(A)=\left\{a\in A\Big|ab=ba,\forall b\in A\right\}$$

It is then clear that these central elements form an associative subalgebra, and:

\medskip

(1) For $A=F(G)$, and more generally when $A$ is commutative, $Z(A)=A$.

\medskip

(2) For $A=F[H]$, consider a linear combination of group elements, as follows:
$$a=\sum_g\lambda_gg$$

By linearity, this element $a\in F[H]$ belongs to the center of $F[H]$ precisely when it commutes with all the group elements $h\in H$, and this gives:
\begin{eqnarray*}
a\in Z(A)
&\iff&ah=ha\\
&\iff&\sum_g\lambda_ggh=\sum_g\lambda_ghg\\
&\iff&\sum_k\lambda_{kh^{-1}}k=\sum_k\lambda_{h^{-1}k}k\\
&\iff&\lambda_{kh^{-1}}=\lambda_{h^{-1}k}
\end{eqnarray*}

We conclude that $\lambda$ must satisfy $\lambda_{gh}=\lambda_{hg}$, and so must be a central function on $H$, as claimed, with the precise conclusion being that the center is given by:
$$Z(A)=\left\{\sum_g\lambda_gg\Big|\lambda_{gh}=\lambda_{hg}\right\}$$

(3) This is a consequence of what we found in (2), and of the fact that the elements $a\in F[H]$ have by definition finite support. Indeed, when the group $H$ is infinite, having the infinite conjugacy class (ICC) property, there is no central function having finite support, except for scalar multiples of the unit, so we have $Z(A)=F$, as stated.
\end{proof}

Many other things can be said about the central elements, and we will leave their study for the enveloping Lie algebras $A=U\mathfrak g$ as an instructive exercise. 

\bigskip

So long for special elements, inside an arbitrary Hopf algebra. The above results are in fact just the tip of the iceberg, and we will be back to this on several occasions, in what follows, and notably in chapter 3 below, when doing representation theory.

\bigskip

Finally, for our discussion to be complete, many things can be said about the group-like, primitive and central elements, in relation with the various possible operations for the Hopf algebras. But here, again, we will leave all this material for later.

\section*{2d. Haar measure} 

As a last topic for this chapter, which is something which is of key importance too, let us discuss now Haar integration. Let us formulate, indeed:

\index{left integral}
\index{right integral}
\index{integral}
\index{Haar integral}

\begin{definition}
Given a Hopf algebra $A$, a linear form $\int:A\to F$ satisfying
$$\left(\int\otimes\,id\right)\Delta=\int(.)1$$
is called a left integral. Similarly,  a linear form $\int:A\to F$ satisfying
$$\left(id\otimes\int\right)\Delta=\int(.)1$$
is called a right integral. If both conditions are satisfied, we call $\int:A\to F$ an integral.
\end{definition}

These notions are motivated by the Haar integration theory on the various types of groups, such as finite, compact, or locally compact. Among others, and in answer to a question that you might have right now, we have to make the distinction between left and right integrals, because for the generally locally compact groups, these two integrals might differ. But more on such topics, which can be quite technical, later. 

\bigskip

As a first illustration, in the case of the function algebras, $A=F(G)$, with $G$ finite group, these notions are all equivalent, and lead to the uniform integration over $G$:

\begin{theorem}
For a function algebra, $A=F(G)$ with $G$ finite group, with the notation $\int\varphi=\int_G\varphi(g)dg$, the left integral condition takes the following form,
$$\int_G\varphi(gh)dg=\int_G\varphi(g)dg$$
and the right integral condition takes the following form,
$$\int_G\varphi(hg)dg=\int_G\varphi(g)dg$$
and in both cases the unique solution is the uniform integration over $G$,
$$\int_G\varphi(g)dg=\frac{1}{|G|}\sum_{g\in G}\varphi(g)$$
under the normalization assumption $\int 1=1$.
\end{theorem}

\begin{proof}
This is something quite self-explanatory, the idea being as follows:

\medskip

(1) With the convention $\int\varphi=\int_g\varphi(g)dg$ from the statement, we have:
\begin{eqnarray*}
\left(\int\otimes\,id\right)\Delta\varphi
&=&\left(\int\otimes\,id\right)\left[(g,h)\to\varphi(gh)\right]\\
&=&\int_G\varphi(gh)dg
\end{eqnarray*}

Thus, the left integral condition reformulates as in the statement.

\medskip

(2) Again with the convention $\int\varphi=\int_g\varphi(g)dg$ from the statement, we have:
\begin{eqnarray*}
\left(id\otimes\int\right)\Delta
&=&\left(id\otimes\int\right)\Delta\left[(h,g)\to\varphi(hg)\right]\\
&=&\int_G\varphi(hg)dg
\end{eqnarray*}

Thus, the right integral condition reformulates as in the statement.

\medskip

(3) When looking now for solutions, be that either for left invariant forms, or right invariant forms, by taking as input Dirac masses, $\varphi=\delta_g$ with $g\in G$, we are led to the conclusion that our invariant linear form must satisfy the following condition:
$$\int_G\delta_g=\int_G\delta_h\quad,\quad\forall g,h$$

Thus the correspoding density function must be constant over $G$, and under the extra assumption $\int 1=1$, the only solution is the uniform, mass 1 integration, as stated.
\end{proof}

As a second illustration, for the group algebras, $A=F[H]$, with $H$ arbitrary group, the notions in Definition 2.23 are again equivalent, with unique solution, as follows:

\begin{theorem}
For a group algebra, $A=F[H]$, with $H$ arbitrary group, the various invariance notions for a linear form $\int:A\to F$ are equivalent, with the solution being
$$\int g=\delta_{g,1}$$
under the normalization assumption $\int 1=1$. When $H$ is finite and abelian we have
$$\int=\int_{\widehat{H}}$$
with $\widehat{H}$ being the dual finite abelian group.
\end{theorem}

\begin{proof}
This is again something quite self-explanatory, the idea being as follows:

\medskip

(1) In what regards the left invariance condition, we have the following computation, using the fact that the group elements $g\in H$ span the group algebra $F[H]$:
\begin{eqnarray*}
\left(\int\otimes\,id\right)\Delta=\int(.)1
&\iff&\left(\int\otimes\,id\right)\Delta(g)=\int g\cdot1\\
&\iff&\left(\int\otimes\,id\right)(g\otimes g)=\int g\cdot1\\
&\iff&\int g\cdot g=\int g\cdot1\\
&\iff&\left[g\neq1\implies\int g=0\right]
\end{eqnarray*}

Thus with the normalization $\int 1=1$, the solution is unique, $\int g=\delta_{g,1}$, as stated.

\medskip

(2) In what regards the right invariance condition, we have the following computation, using again the fact that the group elements $g\in H$ span the group algebra $F[H]$:
\begin{eqnarray*}
\left(id\otimes\int\right)\Delta=\int(.)1
&\iff&\left(id\otimes\int\right)\Delta(g)=\int g\cdot1\\
&\iff&\left(id\otimes\int\right)(g\otimes g)=\int g\cdot1\\
&\iff&\int g\cdot g=\int g\cdot1\\
&\iff&\left[g\neq1\implies\int g=0\right]
\end{eqnarray*}

Thus, with the normalization $\int 1=1$, the solution is unique, $\int g=\delta_{g,1}$, as stated.

\medskip

(3) This is something which follows from the uniqueness of the integral, both from Theorem 2.24 and from here, and which is clear as well from definitions.
\end{proof}

Inspired by the above, a number of things can be said about integrals in the finite dimensional algebra case, by using the duality $A\leftrightarrow A^*$ from chapter 1, as follows:

\begin{theorem}
For a finite dimensional Hopf algebra $A$, in the context of the duality $A\leftrightarrow A^*$, the left and right integrals, regarded as elements $\int^t\in A^*$, must satisfy
$$\int^t\!\!\cdot\,(\varphi-\varepsilon(\varphi))=0\quad,\quad (\varphi-\varepsilon(\varphi))\cdot\!\!\int^t=0$$
for any $\varphi\in A^*$. At the level of the main examples, of mass $1$, these are as follows:
\begin{enumerate}
\item For $A=F(G)$ the integral is $\int^t=\frac{1}{|G|}\sum_{g\in G}g$, as element of $A^*=F[G]$.

\item For $A=F[H]$ the integral is $\int^t=\delta_1$, as element of $A^*=F(H)$.
\end{enumerate}
\end{theorem}

\begin{proof}
In what regards the left integrals, we have the following equivalences:
\begin{eqnarray*}
\left(\int\otimes\,id\right)\Delta=\int(.)1
&\iff&\left(\int\otimes\,id\right)\Delta=\int\otimes\,u\\
&\iff&\Delta^t\left(\int^t\otimes\,id\right)=\int^t\otimes\,u^t\\
&\iff&\Delta^t\left(\int^t\otimes\,id\right)\varphi=\left(\int^t\otimes\,u^t\right)\varphi\\
&\iff&\int^t\cdot\ \varphi=\int^t\cdot\ \varepsilon(\varphi)
\end{eqnarray*}

Similarly, in what regards the right integrals, we have the following equivalences:
\begin{eqnarray*}
\left(id\otimes\int\right)\Delta=\int(.)1
&\iff&\left(id\otimes\int\right)\Delta=\int\otimes\,u\\
&\iff&\Delta^t\left(id\otimes\int^t\right)=\int^t\otimes\,u^t\\
&\iff&\Delta^t\left(id\otimes\int^t\right)\varphi=\left(\int^t\otimes\,u^t\right)\varphi\\
&\iff&\varphi\cdot\int^t=\int^t\cdot\ \varepsilon(\varphi)
\end{eqnarray*}

Thus, main assertion proved, and in what regards now the illustrations:

\medskip

(1) For $A=F(G)$ we have $A^*=F[G]$, with the isomorphism $F[G]\simeq F(G)^*$ coming via $g\to\delta_g$. But this isomorphism maps $\frac{1}{|G|}\sum_{g\in G}g\to\frac{1}{|G|}\sum_{g\in G}\delta_g$, which is exactly the normalized, mass $1$ integral of $F(G)$, as computed in Theorem 2.24.

\medskip

(2) Similarly, for the algebra $A=F[H]$ we have $A^*=F(H)$, with the isomorphism $F(H)\simeq F[H]^*$ coming via $\delta_g\to\delta_g$. But this isomorphism maps $\delta_1\to\delta_1$, which is exactly the normalized, mass $1$ integral of $F[H]$, as computed in Theorem 2.25.
\end{proof}

As a comment here, with a bit more algebraic work, there are many other things that can be said, in the finite dimensional case, as a continuation of the above. For more on all this, further theory and examples, we refer to any specialized Hopf algebra book.

\bigskip

In the general case now, observe that the invariance conditions in Definition 2.23 can be written as follows, in terms of the usual convolution operation $\varphi*\psi=(\varphi*\psi)\Delta$:
$$\int*\ id=id*\!\int=\int(.)1$$

There is a bit of analogy here with what we did in the beginning of this chapter, in relation with the antipode $S$, and many things can be said here. We have indeed:

\index{idempotent state}

\begin{theorem}
Both the left and right integrals $\int:A\to F$, when normalized as to have $\int 1=1$, satisfy the following idempotent linear form condition:
$$\int*\int=\int$$
At the level of the main examples, for this latter condition, these are as follows:
\begin{enumerate}
\item For $A=F(G)$ this condition is satisfied in fact by the normalized uniform integration form over any subgroup $H\subset G$.

\item For $A=F[H]$ this condition is satisfied in fact by the normalized uniform integration form over any quotient group $H\to K$.
\end{enumerate}
\end{theorem}

\begin{proof}
We have several asssertions here, the idea being as follows:

\medskip

(1) In what regards the first assertion, for a left integral $\int:A\to F$, normalized as to have $\int 1=1$, we have indeed the following computation:
\begin{eqnarray*}
\int*\int
&=&\left(\int\otimes\int\right)\Delta\\
&=&\int\circ\left[\left(\int\otimes\,id\right)\Delta\right]\\
&=&\int\circ\left[\int(.)1\right]\\
&=&\int(.)\int(1)\\
&=&\int(.)
\end{eqnarray*}

(2) Also in what regards the first assertion,  for a right integral $\int:A\to F$, again normalized as to have $\int 1=1$, the computation is similar, as follows:
\begin{eqnarray*}
\int*\int
&=&\left(\int\otimes\int\right)\Delta\\
&=&\int\circ\left[\left(id\otimes\int\right)\Delta\right]\\
&=&\int\circ\left[\int(.)1\right]\\
&=&\int(.)\int(1)\\
&=&\int(.)
\end{eqnarray*}

(3) Finally, regarding the various generalizations of the above computations, in the special cases $A=F(G)$ and $A=F[H]$, as indicated in the statement, we will leave these as an instructive exercise. We will be back to this in the next chapter, when talking about quantum subgroups in general, and subgroups of group duals in particular.
\end{proof}

Summarizing, the theory of integrals for the Hopf algebras brings us right away into a number of interesting topics, featuring duality, subgroups, quotients, and more. We will be back to this later, and discuss as well later the relation with representation theory.

\section*{2e. Exercises}

We had a lot of interesting algebra in this chapter, sometimes going towards basic functional analysis, or differential geometry, and as exercises, we have:

\begin{exercise}
Learn more about the Hopf algebra antipode $S$.
\end{exercise}

\begin{exercise}
Learn more about the square of the antipode $S^2$.
\end{exercise}

\begin{exercise}
Clarify the missing details for the group-like elements.
\end{exercise}

\begin{exercise}
Clarify the missing details for the primitive elements.
\end{exercise}

\begin{exercise}
Work out some further examples for the central elements.
\end{exercise}

\begin{exercise}
Compute the Haar integral, for some algebras of your choice.
\end{exercise}

As bonus exercise, reiterated, learn some functional analysis, which is obviously related to all this. The more you will know here, in advance, the better that will be.

\chapter{Product operations}

\section*{3a. Representations}

We have seen so far that some interesting general theory can be developed for the Hopf algebras, in analogy with the basic theory of groups, by using the Hopf algebra maps $\Delta,\varepsilon,S$, and the axioms satisfied by them. However, when doing group theory, you won't get very far just by playing with $m,u,i$, and the situation is pretty much the same with the Hopf algebras, where you won't get very far just by playing with $\Delta,\varepsilon,S$. 

\bigskip

In order to reach to a more advanced theory, we must talk about actions and coactions, and about representations and corepresentations. Many things can be said here, and in what follows we will present the basics, mostly definitions, that we will use right after for talking about various product operations, and keep for later a more detailed study of this, notably in relation with the notion of semisimplicity, and cosemisimplicity.

\bigskip

Let us begin with something straightforward, namely:

\index{action}
\index{representation}

\begin{definition}
An action, or representation, of a Hopf algebra $A$ on a finite dimensional vector space $V$ is a morphism of associative algebras, as follows:
$$\pi:A\to\mathcal L(V)$$
Equivalently, by using a basis of $V$, this is the same as having a morphism as follows:
$$\pi:A\to M_N(F)$$
In this latter situation, we write $\pi=(\pi_{ij})$, with $\pi_{ij}:A\to F$ given by $\pi_{ij}(a)=\pi(a)_{ij}$.
\end{definition}

Observe that the above notion has nothing to do with the Hopf algebra maps $\Delta,\varepsilon,S$, with only the associative algebra structure of $A$ being involved. However, when $A$ is a Hopf algebra, as above, several interesting things can be said, as we will soon discover. 

\bigskip

To start with, in the context of Definition 3.1, the number $N=\dim V$ is called dimension of the representation $\pi$. The simplest situation, namely $N=1$, corresponds by definition to a representation as follows, also called character of $A$:
$$\pi:A\to F$$

So, let us first study these characters, under our assumption that $A$ is a Hopf algebra, as in Definition 3.1. We can say several things here, as follows:

\index{character}

\begin{theorem}
The characters of a Hopf algebra $\pi:A\to F$ are as follows:
\begin{enumerate}
\item When $A$ is finite dimensional, $\pi\in A^*$ must be a group-like element.

\item When $A=F(G)$ with $|G|<\infty$, we must have $\pi(f)=f(g)$, for some $g\in G$.

\item When $A=F[H]$, our character must come from a group morphism $\rho:H\to F^*$.
\end{enumerate}
\end{theorem}

\begin{proof}
This follows from the general Hopf algebra theory that we developed in chapter 1, the details of the proof, and of the statement too, being as follows:

\medskip

(1) Assuming $\dim A<\infty$, we know from chapter 1 how to construct the dual Hopf algebra $A^*$, consisting of the linear forms $\pi:A\to F$. So, let us pick such a linear form, and see when this form is a character. But this happens precisely when $\pi$ is multiplicative, $\pi(ab)=\pi(a)\pi(b)$, and we can process this latter condition as follows:
\begin{eqnarray*}
\pi(ab)=\pi(a)\pi(b)
&\iff&\pi m(a\otimes b)=m(\pi\otimes\pi)(a\otimes b)\\
&\iff&\pi m=m(\pi\otimes\pi)\\
&\iff&m^t\pi^t=(\pi^t\otimes\pi^t)m^t\\
&\iff&m^t\pi^t(1)=(\pi^t\otimes\pi^t)m^t(1)\\
&\iff&m^t\pi^t(1)=(\pi^t\otimes\pi^t)(1\otimes1)\\
&\iff&m^t\pi^t(1)=\pi^t(1)\otimes\pi^t(1)
\end{eqnarray*}

Now forgetting about $A$, and using the notation $\Delta=m^t$ for the comultiplication of $A^*$, and also by identifying $\pi^t(1)\in A^*$ with $\pi\in A^*$, this condition reads:
$$\Delta(\pi)=\pi\otimes\pi$$

Thus, we are led to the conclusion in the statement.

\medskip

(2) Assume now $A=F(G)$, with $|G|<\infty$. Our Hopf algebra $A$ being finite dimensional, what we found in (1) above applies, and we conclude that the characters of $A$ correspond to the group-like elements of the following Hopf algebra:
$$F(G)^*=F[G]$$

But the group-like elements of $F[G]$ are very easy to compute, due to:
$$\Delta\left(\sum_{g\in G}\lambda_g\,g\right)=\sum_{g\in G}\lambda_g\,g\otimes g$$

Indeed, this formula shows that the group-like elements of $F[G]$ are precisely the group elements $g\in G$. Thus, as a conclusion, the characters $\pi:A\to F$ must come from the group elements $g\in G$, and now by carefully looking at what we did in the above, we can also say that the connecting formula is the one in the statement, namely:
$$\pi(f)=f(g)$$

(3) Again, this is something which comes from the general theory from chapter 1. Indeed, assuming $A=F[H]$, to any character $\pi:A\to F$ we can associate a group morphism $\rho:H\to F^*$, simply as being the following composition:
$$\rho:H\subset F[H]\to F$$

Conversely now, given a group morphism $\rho:H\to F^*$, we can associate to it a Hopf algebra character $\pi:F[H]\to F$, simply by linearizing, as follows:
$$\pi\left(\sum_{h\in H}\lambda_h\,h\right)=\sum_{h\in H}\lambda_h\,\rho(h)$$

Thus, we have our bijection, as claimed in the statement.
\end{proof}

Many other things can be said, as a continuation of the above. Recall for instance from chapter 1 that to any finite abelian group $G$ we can associate its dual $\widehat{G}$ with respect to a field $F$, as being the finite abelian group formed by the group characters of $G$:
$$\widehat{G}=\Big\{\rho:G\to F^*\Big\}$$

In fact, again following chapter 1, we can perform this construction for any group $G$, not necessarily finite, or abelian, and we obtain in this way a certain group $\widehat{G}$. Of course, this construction is not always very interesting, for instance because there are non-trivial, and even infinite groups $G$, for which $\widehat{G}=\{1\}$. However, our construction makes sense, as something rather theoretical, and with this in hand, what we found in Theorem 3.2 (3) says that the characters of $A=F[H]$ come from the following group elements:
$$\rho\in\widehat{H}$$

Moving forward now, in the general context of Definition 3.1, we have so far some good understanding on what happens in the case $\dim V=1$, coming from Theorem 3.2. In general things can be quite complicated, and as a first result here, regarding the main examples of Hopf algebras, namely $F(G)$, with $G$ finite group, and $F[H]$, with $H$ arbitrary group, we can say a few things about representations, as follows: 

\begin{theorem}
The following happen:
\begin{enumerate}
\item Given a finite group $G$, and elements $g_1,\ldots,g_N\in G$, we have a representation $\pi:F(G)\to M_N(F)$, given by $\pi(f)=diag(f(g_1),\ldots,f(g_N))$.

\item Given a group $H$, the representations $\pi:F[H]\to M_N(F)$ come by linearization from the group representations $\rho:H\to GL_N(F)$.
\end{enumerate}
\end{theorem}

\begin{proof}
These assertions come as a continuation of Theorem 3.2 (2) and (3), with their proof, along with a bit more on the subject, being as follows:

\medskip

(1) Given a finite group $G$, and a field $F$, let us try to find representations, in the sense of Definition 3.1, of the corresponding function algebra of $G$:
$$\pi:F(G)\to M_N(F)$$

Now observe that, since the algebra $F(G)$ is commutative, so must be its image $Im(\pi)\subset M_N(F)$. Thus, as a first question, we must look for commutative subalgebras $A\subset M_N(F)$. But the standard choice here is the algebra of diagonal matrices $\Delta\subset M_N(F)$, and its various subalgebras $A\subset\Delta$, which are all commutative.

\medskip

(2) With this idea in mind, in order to find basic examples, let us look for representations as follows, with $\Delta\subset M_N(F)$ being the algebra of diagonal matrices:
$$\pi:F(G)\to\Delta$$

But such a representation must be of the special form $\pi=diag(\pi_1,\ldots,\pi_N)$, with $\pi_i:F(G)\to F$ being certain 1-dimensional representations, or characters. Now since such characters must come from group elements, we must have $\pi_i(f)=f(g_i)$, for certain elements $g_i\in G$, and we are led to the formula in the statement, namely:
$$\pi(f)=\begin{pmatrix}
f(g_1)\\
&\ddots\\
&&f(g_N)
\end{pmatrix}$$

(3) Summarizing, we have proved the result for $F(G)$, along with a bit more, namely the fact that any representation of type $\pi:F(G)\to\Delta\subset M_N(F)$ appears as in the statement. It is possible to say more about this, for instance by spinning our representations with the help of a matrix $U\in GL_N(F)$, but no hurry with this, and we will leave this material for later, when systematically doing representation theory.

\medskip

(4) With this done, let us discuss now the second situation in the statement, where we have a group $H$, which can be finite or not, and a field $F$, and we are looking for representations of the corresponding group algebra, as follows:
$$\pi:F[H]\to M_N(F)$$

By restriction to $H\subset F[H]$, we obtain a certain map, as follows:
$$\rho:H\to M_N(F)$$

(5) Now observe that, since $\pi$ is a morphism of algebras, this map $\rho$ must be multiplicative with respect to the group structure of $H$, in the sense that we must have:
$$\rho(g)\rho(h)=\rho(gh)$$

In particular with $h=g^{-1}$ we can see that each $\rho(g)$ must be invertible, and so our map $\rho$ must be in fact a group morphism, as follows:
$$\rho:H\to GL_N(F)$$

(6) But this gives the result. Indeed, by linearity, the representation $\pi$ is uniquely determined by $\rho$, and conversely, given a group morphism $\rho:H\to GL_N(F)$ as above, by linearizing we obtain an algebra representation $\pi:F[H]\to M_N(F)$, as desired.

\medskip

(7) Thus, done with the second assertion too, and as before with the first assertion, several things can be added to this. For instance when $H$ is finite and abelian, we know from chapter 1 that we have an isomorphism as follows, with $G=\widehat{H}$:
$$F(G)\simeq F[H]$$

In view of this, the question is, how do the examples of representations of $F(G)$ constructed in (1) fit with the arbitrary representations of $F[H]$ classified in (2).

\medskip

(8) In answer, the representations classified in (2) correspond, via Pontrjagin duality, to those constructed in (1), further spinned by a matrix $U\in GL_N(F)$, along the lines suggested in step (3) above. We will leave the clarification of this as an instructive exercise, and we will come back to this subject, with full details, in due time.
\end{proof}

Many other things can be said, as a continuation of the above. We will be back to this, once we will have more general theory to be applied, and examples to be studied.

\section*{3b. Corepresentations} 

Moving forward, in order to reach to a continuation of the above, let us recall that a Hopf algebra $A=(A,m,u,\Delta,\varepsilon,S)$ is a special type of bialgebra $A=(A,m,u,\Delta,\varepsilon)$, which itself is a certain mix of an algebra $(A,m,u)$, and a coalgebra $(A,\Delta,\varepsilon)$. We have not talked about such things so far, but right now is the good time to do it.

\bigskip

Indeed, with Definition 3.1 being something only related to the algebra structure of $A$, the question is, what is the ``dual'' definition, related to the coalgebra structure of $A$. In answer, such a dual definition exists indeed, as follows:

\index{coaction}
\index{corepresentation}

\begin{definition}
A coaction, or corepresentation, of a Hopf algebra $A$ on a finite dimensional vector space $V$ is a linear map $\alpha:V\to V\otimes A$ satisfying the condition
$$(\alpha\otimes id)\alpha=(id\otimes\Delta)\alpha$$
called coassociativity. Equivalently, by using a basis of $V$, and writing
$$\alpha(e_i)=\sum_je_j\otimes u_{ji}$$ 
with $u_{ij}\in A$, the square matrix $u=(u_{ij})$ must satisfy the condition
$$\Delta(u_{ij})=\sum_ku_{ik}\otimes u_{kj}$$
also called coassociativity.
\end{definition}

As a first observation, it looks like we forgot here to say something in relation with $\varepsilon$, but that condition is automatic, and more on this in a moment. Also, the fact that the condition at the end is indeed equivalent to that in the beginning is something that must be checked, and this can be done by comparing the following two formulae:
$$(id\otimes\Delta)\alpha(e_j)
=(id\otimes\Delta)\sum_ie_i\otimes u_{ij}
=\sum_ie_i\otimes\Delta(u_{ij})$$
$$(\alpha\otimes id)\alpha(e_j)
=(\alpha\otimes id)\sum_ke_k\otimes u_{kj}
=\sum_{ik}e_i\otimes u_{ik}\otimes u_{kj}$$

Getting back to Definition 3.4 as formulated, as already said, this appears as a ``dual'' of Definition 3.1, and similar comments can be made about it. Let us start with:

\begin{proposition}
The $1$-dimensional corepresentations of $A$ correspond to the elements $a\in A$ satisfying 
$$\Delta(a)=a\otimes a$$
which are the group-like elements of $A$.
\end{proposition}

\begin{proof}
This is indeed clear from definitions, because at $N=1$ the corepesentations of $A$ are the $1\times1$ matrices $u=(a)$, satisfying $\Delta(a)=a\otimes a$.
\end{proof}

Observe the similarity with what we know from Theorem 3.2 (1). However, at the level of the proofs, Theorem 3.2 (1) was something rather complicated, while the above is something trivial. Quite surprising all this, hope you agree with me. In short, we have something interesting here, philosophically speaking, suggesting that Definition 3.4 is something quite magical, when compared to Definition 3.1. Good to know, and because of this, we will be often prefer Definition 3.4 over Definition 3.1, in what follows.

\bigskip

As another comment, Definition 3.4 only involves the comultiplication $\Delta$, and you might wonder about the role of the counit $\varepsilon$ and the antipode $S$, in relation with corepresentations. In answer, at $N=1$ the situation is very simple, because, as we know well from chapter 2, the group-like elements $a\in A$ are subject to the following formulae:
$$\varepsilon(a)=1\quad,\quad S(a)=a^{-1}$$

A similar phenomenon happens in general, with the result here, which can be regarded as being a useful complement to Definition 3.4, being as follows:

\begin{theorem}
Given a corepresentation $u\in M_N(A)$, we have:
$$(id\otimes\varepsilon)u=1\quad,\quad (id\otimes S)u=u^{-1}$$
Also, the associated coaction $\alpha:F^N\to F^N\otimes A$ is counital, $(id\otimes\varepsilon)\alpha=id$.
\end{theorem}

\begin{proof}
There are several things going on here, the idea being as follows:

\medskip

(1) Let us first prove the second formula, the one involving the antipode $S$. For this purpose, we can use the Hopf algebra axiom for the antipode, namely: 
$$m(S\otimes id)\Delta=m(id\otimes S)\Delta=\varepsilon(.)1$$

Indeed, by applying this to $u_{ij}$, and setting $v=(id\otimes S)u$, we have, as desired:
\begin{eqnarray*}
&&m(S\otimes id)\Delta(u_{ij})=m(id\otimes S)\Delta(u_{ij})=\varepsilon(u_{ij})\\
&\implies&\sum_kS(u_{ik})u_{kj}=\sum_ku_{ik}S(u_{kj})=\delta_{ij}\\
&\implies&\sum_kv_{ik}u_{kj}=\sum_ku_{ik}v_{kj}=\delta_{ij}\\
&\implies&(vu)_{ij}=(uv)_{ij}=\delta_{ij}\\
&\implies&v=u^{-1}
\end{eqnarray*}

(2) Let us prove now the first formula, the one involving the counit $\varepsilon$. For this purpose, we can use the Hopf algebra axiom for the counit, namely: 
$$(\varepsilon\otimes id)\Delta=(id\otimes\varepsilon)\Delta=id$$

Indeed, by applying this to $u_{ij}$, and setting $E=(id\otimes\varepsilon)u$, we have:
\begin{eqnarray*}
&&(\varepsilon\otimes id)\Delta(u_{ij})=(id\otimes\varepsilon)\Delta(u_{ij})=u_{ij}\\
&\implies&\sum_k\varepsilon(u_{ik})u_{kj}=\sum_ku_{ik}\varepsilon(u_{kj})=u_{ij}\\
&\implies&\sum_kE_{ik}u_{kj}=\sum_ku_{ik}E_{kj}=u_{ij}\\
&\implies&(Eu)_{ij}=(uE)_{ij}=u_{ij}\\
&\implies&Eu=uE=u
\end{eqnarray*}

Now since $u$ is invertible by (1), we obtain from this, as desired:
$$E=1$$

(3) Regarding now the last assertion, our claim here is that, in the general context of Definition 3.4, the following two counitality conditions are equivalent:
$$(id\otimes\varepsilon)\alpha=id\iff (id\otimes\varepsilon)u=1$$

But this is something which is clear, coming from the following computation:
$$(id\otimes\varepsilon)\alpha(e_i)
=(id\otimes\varepsilon)\sum_je_j\otimes u_{ji}
=\sum_je_j\varepsilon(u_{ji})$$

Indeed, we obtain from this the following equivalences:
\begin{eqnarray*}
(id\otimes\varepsilon)\alpha
&\iff&(id\otimes\varepsilon)\alpha(e_i)=e_i\\
&\iff&\varepsilon(u_{ji})=\delta_{ij}\\
&\iff&(id\otimes\varepsilon)u=1
\end{eqnarray*}

Thus, we are led to the conclusion in the statement.
\end{proof}

Still talking generalities, in the finite dimensional case we have the following result, making it clear that Definition 3.1 and Definition 3.4 are indeed dual to each other:

\begin{theorem}
Given a finite dimensional Hopf algebra $A$:
\begin{enumerate}
\item The representations of $A$ correspond to the corepresentations of $A^*$.

\item The corepresentations of $A^*$ correspond to the representations of $A$.
\end{enumerate}
\end{theorem}

\begin{proof}
In view of the duality result from chapter 1, it is enough to prove one of the assertions, and we will prove the first one. So, consider a linear map, as follows:
$$\pi:A\to M_N(F)$$

As in Definition 3.1, let us construct the coefficients $\pi_{ij}:A\to F$ of this map by the following formula, which must hold for any $a\in A$:
$$\pi_{ij}(a)=\pi(a)_{ij}$$

Now observe that each of these coefficients $\pi_{ij}:A\to F$ can be regarded as an element of the dual algebra $A^*$. As in Definition 3.4, we denote by $u_{ij}$ these elements:
$$\pi_{ij}=u_{ij}\in A^*$$

With these conventions made, we must prove that $\pi$ is a representation of $A$ precisely when $u$ is a corepresentation of $A^*$. But this can be done as follows:

\medskip

(1) Our first claim is that $\pi$ is associative precisely when $u$ is coassociative. But this is something straightforward, which can be established as follows:
\begin{eqnarray*}
\pi(ab)=\pi(a)\pi(b)
&\iff&\pi_{ij}(ab)=\sum_k\pi_{ik}(a)\pi_{kj}(b)\\
&\iff&u_{ij}(ab)=\sum_ku_{ik}(a)u_{kj}(b)\\
&\iff&\Delta(u_{ij})(a\otimes b)=\left(\sum_ku_{ik}\otimes u_{kj}\right)(a\otimes b)\\
&\iff&\Delta(u_{ij})=\sum_ku_{ik}\otimes u_{kj}
\end{eqnarray*}

(2) Our second claim is that $\pi$ is unital precisely when $u$ is counital. But this is again something straightforward, which can be established as follows:
\begin{eqnarray*}
\pi(1)=1
&\iff&\pi_{ij}(1)=\delta_{ij}\\
&\iff&u_{ij}(1)=\delta_{ij}\\
&\iff&\varepsilon(u_{ij})=\delta_{ij}
\end{eqnarray*}

Thus, we are led to the conclusion in the statement.
\end{proof}

Finally, at the level of basic examples, we have the following result, which is in analogy with what we know about representations, from Theorem 3.3:

\begin{theorem}
The following happen:
\begin{enumerate}
\item Given a finite group $G$, the corepresentations $u\in M_N(F(G))$ come, via $u_{ij}(g)=\rho(g)_{ij}$ from the group representations $\rho:G\to GL_N(F)$.

\item Given an arbitrary group $H$, and elements $g_1,\ldots,g_N\in H$, we have a corepresentation $u\in M_N(F[H])$, given by $u=diag(g_1,\ldots,g_N)$.
\end{enumerate}
\end{theorem}

\begin{proof}
As before with Theorem 3.3, many things can be said here, and we will come back to this, on the several occasions, the idea for now being as follows:

\medskip

(1) To start with, in the case of the finite groups, which produce finite dimensional Hopf algebras, the result formally follows from Theorem 3.3, via the duality from Theorem 3.7. Thus, done with (1), and with (2) being trivial anyway, done.

\medskip

(2) However, all this is a bit abstract, so let us check as well (1) directly. Given a finite group $G$, the question is when $u_{ij}\in F(G)$ satisfy the following conditions:
$$\Delta(u_{ij})=\sum_ku_{ik}\otimes u_{kj}\quad,\quad\varepsilon(u_{ij})=\delta_{ij}$$

But, by doing exactly the same computations as those in the proof of Theorem 3.7, the answer here is that this happens precisely when the following is a representation:
$$g\to\begin{pmatrix}
u_{11}(g)&\ldots&u_{1N}(g)\\
\vdots&&\vdots\\
u_{N1}(g)&\ldots&u_{NN}(g)
\end{pmatrix}$$

Thus, we have a bijection with the representations $\rho:G\to GL_N(F)$, as stated.

\medskip

(3) Now let us look at the second assertion in the statement. Given a group $H$, we know from Proposition 3.5 that all the group elements $g\in H\in F[H]$ are 1-dimensional corepresentations. Now let us perform a diagonal sum of such corepresentations:
$$u=\begin{pmatrix}
g_1\\
&\ddots\\
&&g_N
\end{pmatrix}$$

This matrix is then a corepresentation of $F[H]$, as stated, with both the coassociativity and counitaly axioms being clear. And with the remark that we actually proved a bit more, namely that the diagonal corepresentations of $F[H]$ are those in the statement. 

\medskip

(4) Finally, as before in the proof of Theorem 3.3, some further things can be said here, by conjugating such diagonal corepresentations with matrices $U\in GL_N(F)$, and by making the link with the first assertion, in the abelian group case. But, also as before, on several occasions, we will rather leave this for now as an instructive exercise, and come back to all this later, when systematically doing representation theory.
\end{proof}

As a conclusion to all this, we have now a nice representation and corepresentation theory, up and working, for the Hopf algebras, which is related in a nice way to the duality considerations from chapter 1, and with the main examples, which are all quite illustrating, coming from the group algebras of type $F(G)$ and $F[H]$. However, this remains just the tip of the iceberg, and as questions still to be solved, we have:

\bigskip

(1) Fully clarify the classification of the representations and corepresentations of the  group algebras of type $F(G)$ and $F[H]$, in the missing cases.

\bigskip

(2) Clarify as well what happens to the duality between representations of $A$ and corepresentations of $A^*$, when $A$ is no longer finite dimensional.

\bigskip

(3) And finally, discuss what happens for the enveloping Lie algebras $U\mathfrak g$, in relation with the representations of the associated Lie groups $G$.

\bigskip

These questions are all fundamental, but none being trivial, we will leave them for later, when discussing more systematically representation theory. 

\section*{3c. Product operations}

Let us discuss now some natural operations on the Hopf algebras, inspired by those for the groups. We will heavily rely here on the following speculation, from chapter 1:

\begin{speculation}
We can think of any finite dimensional Hopf algebra $A$, not necessarily commutative or cocommutative, as being of the form
$$A=F(G)=F[H]$$
with $G,H$ being finite quantum groups, related by a generalized Pontrjagin duality. And with this generalizing what we know about the finite abelian groups.
\end{speculation}

As explained in chapter 1, this speculation is here for what it is worth, on one hand encapsulating some non-trivial results regarding the finite abelian groups, and the finite dimensional Hopf algebras, and the duality theory for them, but on the other hand, missing some important aspects of the same theory of finite dimensional Hopf algebras. Thus, interesting speculation that we have here, but to be taken with care.

\bigskip

In relation with various product operations, what we would like to have are product operations for the Hopf algebras $\star$, subject to formulae of the following type:
$$F(G)\star F(H)=F(G\circ H)$$

Equivalently, at the dual level, what we would like to have are product operations for the Hopf algebras $\star$, subject to formulae of the following type:
$$F[G]\star F[H]=F(G\bullet H)$$

But probably too much talking, let us get to work. We first have the tensor products of Hopf algebras, whose construction and main properties are as follows:

\index{product of quantum groups}
\index{tensor product}

\begin{theorem}
Given two Hopf algebras $A,B$, so is their tensor product 
$$C=A\otimes B$$
and as main illustrations for this operation, we have the following formulae:
\begin{enumerate}
\item $F(G\times H)=F(G)\otimes F(H)$.

\item $F[G\times H]=F[G]\otimes F[H]$.
\end{enumerate}
\end{theorem}

\begin{proof}
This is something quite self-explanatory, relying on the general theory of the tensor products $\otimes$ explained in chapter 1, the details being as follows:

\medskip

(1) To start with, given two associative algebras $A,B$, so is their tensor product as vector spaces $A\otimes B$, with multiplicative structure as follows:
$$(a\otimes b)(a'\otimes b')=aa'\otimes bb'\quad,\quad 1=1\otimes1$$

Now assume in addition that $A,B$ are Hopf algebras, each coming with its own $\Delta,\varepsilon,S$ operations. In this case we can define $\Delta,\varepsilon,S$ operations on $A\otimes B$, as follows:
$$\Delta(a\otimes b)=\Delta(a)_{13}\Delta(b)_{24}$$
$$\varepsilon(a\otimes b)=\varepsilon(a)\varepsilon(b)$$
$$S(a\otimes b)=S(a)\otimes S(b)$$

(2) But with the above formulae in hand, the verification of the Hopf algebra axioms is straightforward. Indeed, in what regards the comultiplication axiom, we have:
\begin{eqnarray*}
(\Delta\otimes id)\Delta(a\otimes b)
&=&(\Delta\otimes id\otimes id)(\Delta(a)_{13}\Delta(b)_{24})\\
&=&[(\Delta\otimes id)\Delta(a)]_{135}[(\Delta\otimes id)\Delta(b)]_{246}\\
&=&[(id\otimes\Delta)\Delta(a)]_{135}[(id\otimes\Delta)\Delta(b)]_{246}\\
&=&(id\otimes id\otimes\Delta)(\Delta(a)_{13}\Delta(b)_{24})\\
&=&(id\otimes\Delta)\Delta(a\otimes b)
\end{eqnarray*}

As for the counit and antipode axioms, their verification is similar, with no tricks of any kind involved. To be more precise, in what regards the counit axiom, we have:
\begin{eqnarray*}
(\varepsilon\otimes id)\Delta(a\otimes b)
&=&(\varepsilon\otimes\varepsilon\otimes id\otimes id)(\Delta(a)_{13}\Delta(b)_{24})\\
&=&[(\varepsilon\otimes id)\Delta(a)]_1[(\varepsilon\otimes id)\Delta(b)]_2\\
&=&a_1b_2\\
&=&a\otimes b
\end{eqnarray*}

Similarly, we have the following computation, for the other counit axiom:
\begin{eqnarray*}
(id\otimes\varepsilon)\Delta(a\otimes b)
&=&( id\otimes id\otimes\varepsilon\otimes\varepsilon)(\Delta(a)_{13}\Delta(b)_{24})\\
&=&[(id\otimes\varepsilon)\Delta(a)]_1[(id\otimes\varepsilon)\Delta(b)]_2\\
&=&a_1b_2\\
&=&a\otimes b
\end{eqnarray*}

Finally, for the antipode axiom, we have the following computation:
\begin{eqnarray*}
m(S\otimes id)\Delta(a\otimes b)
&=&(m_{13}m_{24})(S\otimes S\otimes id\otimes id)(\Delta(a)_{13}\Delta(b)_{24})\\
&=&(m(S\otimes id)\Delta(a))_1(m(S\otimes id)\Delta(b))_2\\
&=&(\varepsilon(a)1)_1(\varepsilon(b)1)_2\\
&=&\varepsilon(a)1\otimes\varepsilon(b)1\\
&=&\varepsilon(a\otimes b)1\otimes1
\end{eqnarray*}

Similarly, we have the following computation, for the other antipode axiom:
\begin{eqnarray*}
m(id\otimes S)\Delta(a\otimes b)
&=&(m_{13}m_{24})(id\otimes id\otimes S\otimes S)(\Delta(a)_{13}\Delta(b)_{24})\\
&=&(m(id\otimes S)\Delta(a))_1(m(id\otimes S)\Delta(b))_2\\
&=&(\varepsilon(a)1)_1(\varepsilon(b)1)_2\\
&=&\varepsilon(a)1\otimes\varepsilon(b)1\\
&=&\varepsilon(a\otimes b)1\otimes1
\end{eqnarray*}

We conclude that $C=A\otimes B$ is indeed a Hopf algebra, as stated.

\medskip

(3) In what regards now the formula $F(G\times H)=F(G)\otimes F(H)$, when $G,H$ are finite groups, as well as the formula $F[G\times H]=F[G]\otimes F[H]$, when $G,H$ are arbitrary groups, these are both clear from the definition of the tensor product operation.
\end{proof}

As a continuation of the above, in what regards the special elements, we have:

\begin{proposition}
The special elements of $A\otimes B$ are as follows:
\begin{enumerate}
\item $G_{A\otimes B}$ contains $G_A\times G_B$.

\item $P_{A\otimes B}$ contains $P_A,P_B$.

\item $Z(A\otimes B)=Z(A)\otimes Z(B)$.
\end{enumerate}
\end{proposition}

\begin{proof}
This is something quite self-explanatory, the idea being as follows:

\medskip

(1) In what regards the group-like elements, $\Delta(c)=c\otimes c$, assuming $a\in G_A,b\in G_B$ we have the following computation, showing that we have $a\otimes b\in G_{A\otimes B}$:
\begin{eqnarray*}
\Delta(a\otimes b)
&=&\Delta(a)_{13}\Delta(b)_{24}\\
&=&(a\otimes a)_{13}(b\otimes b)_{24}\\
&=&a\otimes b\otimes a\otimes b
\end{eqnarray*}

But this gives the inclusion in the statement, $G_A\times G_B\subset G_{A\otimes B}$.

\medskip

(2) In what regards the primitive elements, $\Delta(c)=c\otimes 1+1\otimes c$, assuming $a\in P_A$ we have the following computation, showing that we have $a\otimes 1\in P_{A\otimes B}$:
\begin{eqnarray*}
\Delta(a\otimes 1)
&=&\Delta(a)_{13}\\
&=&(a\otimes 1+1\otimes a)_{13}\\
&=&a\otimes 1\otimes 1\otimes 1+1\otimes 1\otimes a\otimes 1
\end{eqnarray*}

Similarly, assuming $b\in P_B$ we have $1\otimes b\in P_{A\otimes B}$. We therefore conclude that the Lie algebra $P_{A\otimes B}$ contains the Lie algebras $P_A,P_B$, as stated.

\medskip

(3) In what regards the center, $Z(A)\otimes Z(B)\subset Z(A\otimes B)$ is clear. Conversely, we have the following computation, assuming that the elements $b_i$ are linearly independent:
\begin{eqnarray*}
\sum_ia_i\otimes b_i\in Z(A\otimes B)
&\implies&\left[\sum_ia_i\otimes b_i,a\otimes 1\right]=0\\
&\implies&\sum_ia_ia\otimes b_i=\sum_iaa_i\otimes b_i\\
&\implies&a_ia=aa_i\\
&\implies&a_i\in Z(A)
\end{eqnarray*}

Similarly, assuming that the elements $a_i$ are linearly independent, we have:
$$\sum_ia_i\otimes b_i\in Z(A\otimes B)\implies b_i\in Z(B)$$

Thus we have the reverse inclusion too, so $Z(A\otimes B)=Z(A)\otimes Z(B)$, as stated.
\end{proof}

We have as well a result regarding the Haar integration, as follows:

\begin{theorem}
The Haar integral of a tensor product $A\otimes B$ appears as
$$\int_{A\otimes B}=\int_A\otimes\int_B$$
with this happening for left integrals, right integrals, and integrals.
\end{theorem}

\begin{proof}
This is again something self-explanatory, the idea being as follows:

\medskip

(1) In what regards the left integrals, the verification goes as follows:
\begin{eqnarray*}
\left(\int_{A\otimes B}\otimes\  id\right)\Delta(a\otimes b)
&=&\left(\int_A\otimes\int_B\otimes\ id\otimes id\right)(\Delta(a)_{13}\Delta(b)_{24})\\
&=&\left[\left(\int_A\otimes\ id\right)\Delta(a)\right]_1\left[\left(\int_B\otimes\ id\right)\Delta(b)\right]_2\\
&=&\left(\int_Aa\cdot 1\right)_1\left(\int_Bb\cdot 1\right)_2\\
&=&\int_Aa\cdot\int_Bb\cdot 1\otimes 1\\
&=&\left(\int_A\otimes\int_B\right)(a\otimes b)\cdot 1\otimes 1\\
&=&\int_{A\otimes B}(a\otimes b)\cdot 1\otimes 1
\end{eqnarray*}

(2) In what regards the right integrals, the verification is similar, as follows:
\begin{eqnarray*}
\left(id\otimes\int_{A\otimes B}\right)\Delta(a\otimes b)
&=&\left(id\otimes id\otimes\int_A\otimes\int_B\right)(\Delta(a)_{13}\Delta(b)_{24})\\
&=&\left[\left(id\otimes\int_A\right)\Delta(a)\right]_1\left[\left(id\otimes\int_B\right)\Delta(b)\right]_2\\
&=&\left(\int_Aa\cdot 1\right)_1\left(\int_Bb\cdot 1\right)_2\\
&=&\int_Aa\cdot\int_Bb\cdot 1\otimes 1\\
&=&\left(\int_A\otimes\int_B\right)(a\otimes b)\cdot 1\otimes 1\\
&=&\int_{A\otimes B}(a\otimes b)\cdot 1\otimes 1
\end{eqnarray*}

(3) Finally, in relation with all this, there is a uniqueness discussion to be made too, which is quite standard, and that we will leave here as an instructive exercise.
\end{proof}

Summarizing, we have now a good understanding of the tensor product operation, with good results all around the spectrum, with respect to the general theory developed in chapters 1-2. Many other things can be said, for instance with some straightforward results regarding the representations and corepresentations, introduced earlier in this chapter. We will be back to tensor products on a regular basis, in what follows.

\bigskip

Moving forward, now that we know about tensor products, we can do exactly the same thing with free products, and we are led in this way to the following result:

\index{free product}
\index{dual free product}

\begin{theorem}
Given two Hopf algebras $A,B$, so is their free product 
$$C=A*B$$
and as main illustrations for this operation, we have the following formulae:
\begin{enumerate}
\item $F(G\,\hat{*}\,H)=F(G)*F(H)$, standing as definition for $G\,\hat{*}\,H$, as quantum group.

\item $F[G*H]=F[G]*F[H]$.
\end{enumerate}
\end{theorem}

\begin{proof}
This is again something self-explanatory, save for the abstract meaning of the object $G\,\hat{*}\,H$ appearing in (1), that we will explain below, the details being as follows:

\medskip

(1) To start with, given two associative algebras $A,B$, so is their free product $A*B$, with the multiplicative structure being by definition as follows:
$$(\ldots a_ib_i\ldots)(\ldots a_i'b_i'\ldots)=\ldots a_ib_i\ldots a_i'b_i'\ldots$$
$$1=1_A=1_B$$

Now assume in addition that $A,B$ are Hopf algebras, each coming with its own $\Delta,\varepsilon,S$ operations. In this case we can define $\Delta,\varepsilon,S$ operations on $A*B$, as follows:
$$\Delta(\ldots a_ib_i\ldots)=\ldots\Delta(a_i)\Delta(b_i)\ldots$$
$$\varepsilon(\ldots a_ib_i\ldots)=\ldots\varepsilon(a_i)\varepsilon(b_i)\ldots$$
$$S(\ldots a_ib_i\ldots)=\ldots S(b_i)S(a_i)\ldots$$

(2) But with the above formulae in hand, the verification of the Hopf algebra axioms is straightforward. Indeed, in what regards the comultiplication axiom, we have:
\begin{eqnarray*}
(\Delta\otimes id)\Delta(\ldots a_ib_i\ldots)
&=&(\Delta\otimes id)(\ldots\Delta(a_i)\Delta(b_i)\ldots)\\
&=&\ldots[(\Delta\otimes id)\Delta(a_i)][(\Delta\otimes id)\Delta(b_i)]\ldots\\
&=&\ldots[(id\otimes\Delta)\Delta(a_i)][(id\otimes\Delta)\Delta(b_i)]\ldots\\
&=&(id\otimes\Delta)(\ldots\Delta(a_i)\Delta(b_i)\ldots)\\
&=&(id\otimes\Delta)\Delta(\ldots a_ib_i\ldots)
\end{eqnarray*}

As for the counit and antipode axioms, their verification is similar, with no tricks of any kind involved. To be more precise, in what regards the counit axiom, we have:
\begin{eqnarray*}
(\varepsilon\otimes id)\Delta(\ldots a_ib_i\ldots)
&=&(\varepsilon\otimes id)(\ldots\Delta(a_i)\Delta(b_i)\ldots)\\
&=&\ldots(\varepsilon\otimes id)\Delta(a_i)(\varepsilon\otimes id)\Delta(b_i)\ldots\\
&=&\ldots a_ib_i\ldots
\end{eqnarray*}

Similarly, we have the following computation, for the other counit axiom:
\begin{eqnarray*}
(id\otimes\varepsilon)\Delta(\ldots a_ib_i\ldots)
&=&(id\otimes\varepsilon)(\ldots\Delta(a_i)\Delta(b_i)\ldots)\\
&=&\ldots(id\otimes\varepsilon)\Delta(a_i)(id\otimes\varepsilon)\Delta(b_i)\ldots\\
&=&\ldots a_ib_i\ldots
\end{eqnarray*}

Finally, for the antipode axiom, we have the following computation:
\begin{eqnarray*}
m(S\otimes id)\Delta(\ldots a_ib_i\ldots)
&=&m(S\otimes id)(\ldots\Delta(a_i)\Delta(b_i)\ldots)\\
&=&m(\ldots(S\otimes id)\Delta(b_i)(S\otimes id)\Delta(a_i)\ldots)\\
&=&\varepsilon(\ldots a_ib_i\ldots)1
\end{eqnarray*}

Similarly, we have the following computation, for the other antipode axiom:
\begin{eqnarray*}
m(id\otimes S)\Delta(\ldots a_ib_i\ldots)
&=&m(id\otimes S)(\ldots\Delta(a_i)\Delta(b_i)\ldots)\\
&=&m(\ldots(id\otimes S)\Delta(b_i)(id\otimes S)\Delta(a_i)\ldots)\\
&=&\varepsilon(\ldots a_ib_i\ldots)1
\end{eqnarray*}

We conclude that $C=A*B$ is indeed a Hopf algebra, as stated.

\medskip

(3) In what regards now the formula $F(G\,\hat{*}\,H)=F(G)*F(H)$, when $G,H$ are finite groups, this stands as a definition for $G\,\hat{*}\,H$, as a quantum group, the point being that, unless $G$ or $H$ is trivial, the Hopf algebra $F(G)*F(H)$ is not commutative, and so cannot be understood as being an algebra of functions. Welcome to noncommutativity.

\medskip

(4) As for the formula $F[G*H]=F[G]*F[H]$, with here $G,H$ being arbitrary groups, possibly infinite, this is something which is clear from definitions.
\end{proof}

As a continuation of the above, in what regards the special elements, we have:

\begin{proposition}
The special elements of $A*B$ are as follows:
\begin{enumerate}
\item $G_{A*B}$ contains $G_A*G_B$.

\item $P_{A*B}$ contains $P_A,P_B$.

\item $Z(A*B)=F$, unless $A=F$, or $B=F$.
\end{enumerate}
\end{proposition}

\begin{proof}
This is something quite similar to Proposition 3.11, with the various computations being very similar to those there, the idea being as follows:

\medskip

(1) As before with tensor products, we have $G_A,G_B\subset G_{A*B}$, which gives the result.

\medskip

(2) Also as before with tensor products, we have $P_A,P_B\subset P_{A*B}$, as claimed.

\medskip

(3) Finally, in what regards the center, things are different with respect to Proposition 3.11. Indeed, since the elements $a\in A-F$ cannot commute with the elements $b\in B-F$, by definition of the free product $A*B$, we generically have $Z(A*B)=F$, as claimed.
\end{proof}

We have as well a result regarding the Haar integration, as follows:

\begin{theorem}
The Haar integral of a free product $A*B$ appears as
$$\int_{A*B}=\int_A*\int_B$$
with this happening for left integrals, right integrals, and integrals.
\end{theorem}

\begin{proof}
This is again something self-explanatory, the idea being as follows:

\medskip

(1) In what regards the left integrals, the verification goes as follows:
\begin{eqnarray*}
\left(\int_{A*B}\otimes\  id\right)\Delta(\ldots a_ib_i\ldots)
&=&\left(\int_A*\int_B\otimes\ id\right)(\ldots\Delta(a_i)\Delta(b_i)\ldots)\\
&=&\ldots\left(\int_A\otimes\ id\right)\Delta(a_i)\left(\int_B\otimes\ id\right)\Delta(b_i)\ldots\\
&=&\ldots\int_Aa_i\cdot 1\int_Bb_i\cdot1\ldots\\
&=&\int_{A*B}(\ldots a_ib_i\ldots)\cdot 1
\end{eqnarray*}

(2) In what regards the right integrals, the verification is similar, as follows:
\begin{eqnarray*}
\left(id\otimes\int_{A*B}\right)\Delta(\ldots a_ib_i\ldots)
&=&\left(id\otimes\int_A*\int_B\right)(\ldots\Delta(a_i)\Delta(b_i)\ldots)\\
&=&\ldots\left(id\otimes\int_A\right)\Delta(a_i)\left(id\otimes\int_B\right)\Delta(b_i)\ldots\\
&=&\ldots\int_Aa_i\cdot 1\int_Bb_i\cdot1\ldots\\
&=&\int_{A*B}(\ldots a_ib_i\ldots)\cdot 1
\end{eqnarray*}

(3) Finally, in relation with all this, there is a uniqueness discussion to be made too, which is quite standard, and that we will leave here as an instructive exercise.
\end{proof}

As before with the tensor products, many other things can be said about free products, for instance with some straightforward results regarding their representations and corepresentations. We will be back to free products on a regular basis, in what follows.

\bigskip

As a further comment here, algebrically speaking, there are several other possible products, which are quite natural, between the tensor products and the free products. But things here are quite technical, and we will discuss them later in this book.

\section*{3d. Quotients, subalgebras} 

Another standard operation, that we would like to discuss now, is that of taking quantum subgroups, with the result here, at the algebraic level, being as follows:

\index{quantum subgroup}
\index{quotient quantum group}
\index{Hopf ideal}

\begin{theorem}
Given a Hopf algebra $A$, so is its quotient $B=A/I$, provided that $I\subset A$ is an ideal satisfying the following conditions, called Hopf ideal conditions,
$$\Delta(I)\subset A\otimes I+I\otimes A\quad,\quad\varepsilon(I)=0\quad,\quad S(I)\subset I$$
and as main illustrations for this operation, we have the following formulae:
\begin{enumerate}
\item $F(G)/I=F(H)$, with $H\subset G$ being a certain subgroup.

\item $F[G]/I=F[H]$, with $G\to H$ being a certain quotient.
\end{enumerate}
\end{theorem}

\begin{proof}
As before, this is something self-explanatory, the idea being as follows:

\medskip

(1) Given an associative algebra $A$ and an ideal $I\subset A$, we can certainly construct the quotient $B=A/I$, which is an associative algebra. Thus, we must just see when the Hopf algebra operations $\Delta,\varepsilon,S$ correctly factorize, from $A$ to $B$.

\medskip

(2) Let us first see when $\Delta$ factorizes. If we denote by $\pi:A\to B$ the canonical projection, the factorization diagram that we are looking for is as follows:
$$\xymatrix@R=50pt@C=70pt{
A\ar[r]^\Delta\ar[d]_\pi&A\otimes A\ar[d]^{\pi\otimes\pi}\\
B\ar@.[r]&B\otimes B
}$$

We can see that the factorization condition is as follows:
$$\pi(a)=0\implies(\pi\otimes\pi)\Delta(a)=0$$

Thus, in terms of the ideal $I\subset A$, the following condition must be satisfied:
$$a\in I\implies(\pi\otimes\pi)\Delta(a)=0$$

But, for an element $b\in A\otimes A$, we have the following equivalence:
$$(\pi\otimes\pi)b=0\iff b\in A\otimes I+I\otimes A$$

We conclude that the factorization condition for $\Delta$ is as follows:
$$a\in I\implies\Delta(a)\in A\otimes I+I\otimes A$$

But this is precisely the first condition on $I$ in the statement, namely:
$$\Delta(I)\subset A\otimes I+I\otimes A$$

(3) Similarly, the counit $\varepsilon$ factorizes precisely when the condition $\varepsilon(I)=0$ in the statement is satisfied, with the factorization diagram here being as follows:
$$\xymatrix@R=40pt@C=70pt{
A\ar[r]^\varepsilon\ar[d]_\pi&F\\
B\ar@.[ur]
}$$

(4) Also, the antipode $S$ factorizes precisely when the condition $S(I)\subset I$ in the statement is satisfied, with the factorization diagram here being as follows: 
$$\xymatrix@R=50pt@C=70pt{
A\ar[r]^S\ar[d]_\pi&A^{opp}\ar[d]^{\pi^{opp}}\\
B\ar@.[r]&B^{opp}
}$$

(5) Together with the remark that the maps $\Delta,\varepsilon,S$, once factorized, will keep satisfying the Hopf algebra axioms, automatically, we are led to the first assertion.

\medskip

(6) In the group setting now, the formula $F(G)/I=F(H)$, with $H\subset G$ being a certain subgroup, is something which clear from definitions.

\medskip

(7) As for the last formula, namely $F[G]/I=F[H]$, with $G\to H$ being a certain quotient, this is something which is clear from definitions too.
\end{proof}

As a continuation of the above, in what regards the special elements, we have:

\begin{proposition}
In what regards the special elements of a quotient $B=A/I$, the quotient map $A\to B$ induces quotient maps as follows:
\begin{enumerate}
\item $G_A\to G_B$.

\item $P_A\to P_B$.

\item $Z(A)\to Z(B)$.
\end{enumerate}
\end{proposition}

\begin{proof}
This is something quite trivial, because in what regards the group-like elements, the primitive elements, and the central elements too, their defining formulae pass to the quotient, in the obvious way. Thus, we are led to the above conclusions.
\end{proof}

We have as well a statement regarding the Haar integration, as follows:

\begin{fact}
The Haar integral of a quotient $B=A/I$ does not appear as for the other Hopf algebra operations, simply by factorizing the following diagram,
$$\xymatrix@R=40pt@C=70pt{
A\ar[r]^{\int_A}\ar[d]_\pi&F\\
B\ar@.[ur]
}$$
but is related however to the Haar integral of $A$, via a number of more technical formulae, and with this happening for left integrals, right integrals, and integrals.
\end{fact}

\begin{proof}
This is obviously something quite informal, that we included here for the sake of symmetry, with respect to the other operations, the idea being as follows:

\medskip

(1) The first assertion is certainly something which happens, coming from the intuitive fact that, by taking for instance $B=A/A$, we have certainly not constructed here the integral of $A$, which is well-known to require some hard work, harder than this.

\medskip

(2) As for the second assertion, this is something more technical, say involving representations and corepresentations. We will leave working out the details here, based on what happens in the cases $A=F(G)$ and $A=F[H]$, as an instructive exercise, and we will come back to such questions later, under a number of suitable extra assumptions.
\end{proof}

As before with the tensor and free products, many other things can be said, for instance with some straightforward results regarding the representations and corepresentations of quotients. We will be back to quotients on a regular basis, in what follows.

\bigskip

Regarding now taking quotient quantum groups, the result here is quite similar, somehow dual to Theorem 3.16, but technically very straightforward, as follows:

\index{quantum subgroup}
\index{quotient quantum group}

\begin{theorem}
Given a Hopf algebra $A$, any subalgebra $B\subset A$ satisfying
$$\Delta(B)\subset B\otimes B\quad,\quad S(B)\subset B$$
is a Hopf algebra, and as main illustrations, we have the following subalgebras:
\begin{enumerate}
\item $F(H)\subset F(G)$, for any quotient group $G\to H$.

\item $F[H]\subset F[G]$, for any subgroup $H\subset G$.
\end{enumerate}
\end{theorem}

\begin{proof}
The main assertion is clear from definitions, because the Hopf algebra axioms being satisfied over $A$, they are satisfied as well over the subalgebra $B\subset A$. As for the main illustrations, in the group case, these are as well both clear from definitions.
\end{proof}

The above result is a bit abstract, and as a useful version of it, providing examples, let us record as well the following statement, that will play an important role later:

\begin{theorem}
Given a Hopf algebra $A$, and a finite dimensional corepresentation $u=(u_{ij})$, the subalgebra generated by the coefficients of $u$,
$$B=<u_{ij}>\subset A$$
is a Hopf algebra. As main illustrations for this operation, we obtain subalgebras:
\begin{enumerate}
\item $F(H)\subset F(G)$, with $G\to H$ being certain quotients.

\item $F[H]\subset F[G]$, with $H\subset G$ being certain subgroups.
\end{enumerate}
\end{theorem}

\begin{proof}
This is something which follows from Theorem 3.19, and that we will actually fully clarify later in this book, the idea being as follows:

\medskip

(1) Given a coalgebra $A$ and a corepresentation $u=(u_{ij})$, we can certainly construct the space of coefficients $C_u=<u_{ij}>\subset A$, which is automatically a coalgebra. Indeed, recall that the corepresentations are subject to the following condition:
$$\Delta(u_{ij})=\sum_ku_{ik}\otimes u_{kj}$$

But, by using this condition, we see that $\Delta$ leaves indeed invariant $C_u$:
$$\Delta(C_u)\subset C_u\otimes C_u$$

(2) In the case now where $A$ is a Hopf algebra, by setting $B=<C>\subset A$, our claim is that we obtain both an algebra and a coalgebra, so that we have a Hopf algebra. Indeed, in what regards $\Delta$, the inclusion found in (1) gives, by multiplicativity:
$$\Delta(B)\subset B\otimes B$$

In what regards now the counit $\varepsilon$, there is no verification needed, because we can simply take the restriction of the counit of $A$, to the subalgebra $B\subset A$. 

\medskip

(3) Finally, in what regards the antipode $S$, things here are more tricky. We can use here the following formula, that we know from Theorem 3.6:
$$(id\otimes S)u=u^{-1}$$

Thus, with the convention that our subalgebra $B=<C>\subset A$ contains the inverses of all invertible elements $c\in C$, we see that $S$ satisfies the following condition:
$$S(B)\subset B$$

Alternatively, the fact that our subalgebra $B=<C>\subset A$ contains the inverses of all invertible elements $c\in C$ can come by theorem, under suitable assumptions on the class of algebras involved. We will be back to this point, later on in this book.

\medskip

(4) As a conclusion, under the assumptions in the statement, the Hopf algebra maps $\Delta,\varepsilon,S$ restrict to the subalgebra $B=<C>\subset A$. But the Hopf algebras axioms being automatic for these restrictions, we are led to the first assertion.

\medskip

(5) In the group setting now, the formula $F(H)\subset F(G)$, with $G\to H$ being a certain quotient, is something which is clear from definitions.

\medskip

(6) As for the last formula, namely $F[H]\subset F[G]$, with $H\subset G$ being a certain subgroup, this is something which is clear from definitions too.
\end{proof}

As a continuation of the above, in what regards the special elements, we have:

\begin{proposition}
The special elements of a subalgebra $B\subset A$ are as follows,
\begin{enumerate}
\item $G_B=G_A\cap B$.

\item $P_B=P_A\cap B$.
\end{enumerate}
and in what regards the center, nothing in particular can be said.
\end{proposition}

\begin{proof}
This is something trivial, the idea being as follows:

\medskip

(1) The first two assertions are clear, since the comultiplication of $B$ appears by definition as the restriction of the comultiplication of $A$. Thus, when it comes to group-like elements, or to primitive elements, we obtain the formulae in the statement.

\medskip

(2) As for the last assertion, this is something informal, and we will leave some thinking here, with various examples and counterexamples, as an instructive exercise, the idea being that the center of $B$ can be substantially smaller, or larger, than the center of $A$.
\end{proof}

We have as well a result regarding the Haar integration, as follows:

\begin{theorem}
The Haar integral of a subalgebra $B\subset A$ appears as a restriction
$$\int_B=\left(\int_A\right)_B$$
and with this happening for left integrals, right integrals, and integrals.
\end{theorem}

\begin{proof}
This is again something quite self-explanatory, and clear from definitions, with the corresponding commuting diagram here being as follows:
$$\xymatrix@R=40pt@C=70pt{
A\ar[r]^{\int_A}&F\\
B\ar[u]^i\ar[ur]_{\int_B}
}$$

Thus, we are led to the conclusion in the statement.
\end{proof}

As before with the products and quotients, many other things can be said, for instance with some straightforward results regarding the representations and corepresentations of subalgebras. We will be back to subalgebras on a regular basis, in what follows.

\bigskip

Finally, in relation with the above quotient and subalgebra operations, and with their main illustrations too, there is some discussion to be made, in the finite dimensional case, in the context of the Hopf algebra duality for such Hopf algebras, from chapter 1. To be more precise, we have here the following result, which clarifies the situation:

\begin{theorem}
In the context of the duality for finite dimensional Hopf algebras
$$A\leftrightarrow A^*$$
the operations of taking quotients and subalgebras are dual to each other.
\end{theorem}

\begin{proof}
This is something straightforward, the idea being as follows:

\medskip

(1) Consider a finite dimensional Hopf algebra $A$, with structural maps as follows:
$$m:A\otimes A\to A$$
$$u:F\to A$$
$$\Delta:A\to A\otimes A$$
$$\varepsilon:A\to F$$
$$S:A\to A^{opp}$$

(2) As explained in chapter 1, the dual vector space $A^*$, consisting of the linear forms $\varphi:A\to F$, is then a Hopf algebra too, with structural maps as follows:
$$\Delta^t:A^*\otimes A^*\to A^*$$
$$\varepsilon^t:F\to A^*$$
$$m^t:A^*\to A^*\otimes A^*$$
$$u^t:A^*\to F$$
$$S^t:A^*\to (A^*)^{opp}$$

(3) But, with these formulae in hand, it is straightforward to check that the Hopf algebra quotients of $A$ correspond to the Hopf subalgebras of $A^*$, and vice versa:
$$B\subset A\quad\longleftrightarrow\quad C=A^*/I$$
$$B=A/I\quad\longleftrightarrow\quad C\subset A^*$$

(4) Moreover, in the context of the duality between the function algebras $A=F(G)$ and the group algebras $A=F[G]$, with $G$ being a finite group, we obtain in this way generalizations of the well-known fact that the quotients of a finite abelian group $G$ correspond to the subgroups of the dual finite abelian group $\widehat{G}$, and vice versa.

\medskip

(5) So, this was for the general idea, that everything comes in the end from what we know about the finite abelian groups, and we will leave the proof of the various assertions formulated above, in the precise order that you prefer, as an instructive exercise. 
\end{proof}

\section*{3e. Exercises}

We had a lot of interesting algebra in this chapter, and as exercises, we have:

\begin{exercise}
Further study the representation theory of $F(G)$.
\end{exercise}

\begin{exercise}
Further study the representation theory of $F[H]$.
\end{exercise}

\begin{exercise}
Fill in the missing details for the $\otimes$ operation.
\end{exercise}

\begin{exercise}
Fill in the missing details for the $*$ operation.
\end{exercise}

\begin{exercise}
Clarify the missing details for the quotient operation.
\end{exercise}

\begin{exercise}
Clarify the missing details for the subalgebra operation.
\end{exercise}

As bonus exercise, figure out what happens to all the above when $F=\mathbb C$.

\chapter{Affine algebras}

\section*{4a. Affine algebras}

We have seen so far the Hopf algebra basics, including the theory of the basic operations for the Hopf algebras. We would like to discuss now a number of more tricky operations on the Hopf algebras, appearing as variations of the above. 

\bigskip

For this purpose, let us introduce the following notion, inspired as usual from group theory, that will play a key role in this book, starting from now, and until the end:

\index{affine algebra}
\index{affine Hopf algebra}
\index{algebraic group}
\index{algebraic quantum group}

\begin{definition}
We call a Hopf algebra $A$ affine when it is of the form
$$A=<u_{ij}>$$
with $u\in M_N(A)$ being a corepresentation, called fundamental corepresentation.
\end{definition}

As already mentioned, this notion is inspired from group theory, and more specifically, from advanced group theory, our motivation coming from the following facts:

\bigskip

(1) In group theory, at a reasonably advanced level, a natural assumption on a group $G$ is that this appears as an algebraic group, $G\subset GL_N(F)$. This is indeed how most of the examples of groups $G$ appear, in practice, as groups of invertible matrices.

\bigskip

(2) Observe that any finite group $G$ appears as above, thanks to the Cayley embedding theorem, $G\subset S_N$ with $N=|G|$, coupled with the standard embedding $S_N\subset GL_N(F)$ given by the permutation matrices, which give an embedding as follows:
$$G\subset S_N\subset GL_N(F)$$

(3) As another key example, which is more advanced, it is known that any compact Lie group $G$ appears as a group of unitary matrices, $G\subset U_N$, so that we have:
$$G\subset U_N\subset GL_N(\mathbb C)$$

(4) Now the point is that, save for a few topological issues, the fact that a group $G$ is algebraic is equivalent to the fact that the Hopf algebra $A=F(G)$ is affine, in the sense of Definition 4.1. Thus, we have here a good motivation for Definition 4.1.

\bigskip

In addition to this, getting now to the group dual level, we have some extra motivations for Definition 4.1, again coming from advanced group theory, which are as follows:

\bigskip

(1) Again in group theory, at a reasonably advanced level, but this time with the discrete group theory in mind, instead of the continuous group theory used above, a natural assumption on a group $H$ is that this group is finitely generated:
$$H=<g_1,\ldots,g_N>$$

(2) Indeed, this is how nearly all the interesting discrete groups appear. In fact, many of these groups appear by definition via generators and relations, as follows:
$$H=\left<g_1,\ldots,g_N\,\Big|\,\mathcal R\right>$$

(3) As another comment here, passed the wealth of examples, the fact of being finitely generated is needed for developing the theory, because talking about Cayley graphs and their metric aspects, random walks and so on, always requires the use of generators.

\bigskip

(4) Now the point is that, save for a few topological issues, the fact that a group $H$ is finitely generated is equivalent to the fact that the Hopf algebra $A=F[H]$ is affine, in the sense of Definition 4.1. Thus, we have again a good motivation for Definition 4.1.

\bigskip

Summarizing, in order to reach to a more advanced Hopf algebra theory, we have all good reasons in this world to assume that our algebras are affine, as in Definition 4.1. And for ending this discussion, let us formulate our conclusions as follows:

\begin{conclusion}
In relation with group advanced theory:
\begin{enumerate}
\item It makes sense to assume that the groups are algebraic, $G\subset GL_N(F)$,

\item Or to assume that the groups are finitely generated, $H=<g_1,\ldots,g_N>$,

\item With this basically corresponding to the fact that $F(G)$ and $F[H]$ are affine,

\item So, we will assume in what follows that our Hopf algebras $A$ are affine.
\end{enumerate}
\end{conclusion}

And more on this later. Also, we will see many examples of affine Hopf algebras in what follows, appearing via various operations, the idea being that the affine Hopf algebras are subject to far more operations than those discussed before. More later.

\bigskip

In order to get started now, we need to know more about corepresentations, and in particular, about the fundamental one. As a first result, that we will need in what follows, the corepresentations of an arbitrary Hopf algebra $A$ are subject to a number of operations, exactly as the group representations in the usual group case, as follows:

\index{sum of corepresentations}
\index{product of corepresentations}
\index{tensor product of corepresentations}
\index{conjugate of corepresentation}
\index{spinned corepresentation}
\index{spinning by unitaries}

\begin{proposition}
The corepresentations of a Hopf algebra $A$ are subject to the following operations, in analogy with what happens for the group representations:
\begin{enumerate}
\item Making sums, $u+v=diag(u,v)$.

\item Making tensor products, $(u\otimes v)_{ia,jb}=u_{ij}v_{ab}$.

\item Spinning by invertible scalar matrices, $u\to VuV^{-1}$.
\end{enumerate}
\end{proposition}

\begin{proof}
Observe first that the result holds indeed for $A=F(G)$, where we obtain the usual operations on the representations of $G$. In general, the proof goes as follows:

\medskip

(1) Everything here is clear from definitions.

\medskip

(2) The comultiplicativity condition follows indeed from the following computation: 
\begin{eqnarray*}
\Delta((u\otimes v)_{ia,jb})
&=&\Delta(u_{ij}v_{ab})\\
&=&\Delta(u_{ij})\Delta(v_{ab})\\
&=&\sum_ku_{ik}\otimes u_{kj}\sum_cv_{ac}\otimes v_{cb}\\
&=&\sum_{kc}u_{ik}v_{ac}\otimes u_{kj}v_{cb}\\
&=&\sum_{kc}(u\otimes v)_{ia,kc}\otimes(u\otimes v)_{kc,jb}
\end{eqnarray*}

(3) The comultiplicativity property of the matrix $v=VuV^{-1}$ in the statement can be checked by doing a straightforward computation. Alternatively, if we write $u\in M_n(F)\otimes A$, the usual comultiplicativity axiom, as formulated in chapter 3, reads:
$$(id\otimes\Delta)u=u_{12}u_{13}$$

Here we use standard tensor calculus conventions. Now when spinning by a scalar matrix, the matrix that we obtain is $v=V_1uV_1^{-1}$, and we have:
\begin{eqnarray*}
(id\otimes\Delta)v
&=&V_1u_{12}u_{13}V_1^{-1}\\
&=&V_1u_{12}V_1^{-1}\cdot V_1u_{13}V_1^{-1}\\
&=&v_{12}v_{13}
\end{eqnarray*}

Thus, with usual notations, $v=VuV^{-1}$ is a corepresentation, as claimed.
\end{proof}

As a comment now, the various operations in Proposition 4.3 can be viewed as operations on the class of affine Hopf algebras, the result here being as follows:

\begin{proposition}
We have the following operations on the affine Hopf algebras, with the convention $A_u=<u_{ij}>\subset A$, for a corepresentation $u\in M_N(A)$:
\begin{enumerate}
\item $(A_u,A_v)\to A_{u+v}$.

\item $(A_u,A_v)\to A_{u\otimes v}$.

\item $(A,u)\to(A,VuV^{-1})$.
\end{enumerate}
\end{proposition}

\begin{proof}
This is indeed something clear, coming from the various operations on the corepresentations constructed in Proposition 4.3, and with the remark that, in the context of the last assertion, we have indeed $<u_{ij}>=<(VuV^{-1})_{ij}>$, as needed there.
\end{proof}

As a further operation, also inspired from usual group theory, we have:

\index{conjugate corepresentation}
\index{contragradient corepresentation}

\begin{theorem}
Given a Hopf algebra corepresentation $u\in M_n(A)$,
$$\bar{u}=(t\otimes id)u^{-1}$$
is a corepresentation too, called contragradient, or conjugate to $u$.
\end{theorem}

\begin{proof}
This is something very standard, the idea being as follows:

\medskip

(1) To start with, we know from chapter 3 that $u$ is indeed invertible, with inverse given by $u^{-1}=(id\otimes S)u$. Thus, $\bar{u}$ is well-defined, and in addition, we have:
$$\bar{u}=(t\otimes S)u$$

(2) But with this latter formula in hand, the proof of the corepresentation property of $\bar{u}$ goes as follows, by using the formula $\Delta S=\Sigma(S\otimes S)\Delta$ from chapter 2:
\begin{eqnarray*}
\Delta(\bar{u}_{ij})
&=&\Delta S(u_{ji})\\
&=&\Sigma(S\otimes S)\Delta(u_{ji})\\
&=&\Sigma(S\otimes S)\sum_ku_{jk}\otimes u_{ki}\\
&=&\sum_kS(u_{ki})\otimes S(u_{jk})\\
&=&\sum_k\bar{u}_{ik}\otimes\bar{u}_{kj}
\end{eqnarray*}

(3) Alternatively, we can check the corepresentation property of $\bar{u}$ as follows, simply by inverting the corepresentation property of $u$, written in compact form:
\begin{eqnarray*}
(id\otimes\Delta)u=u_{12}u_{13}
&\implies&(id\otimes\Delta)u^{-1}=u_{13}^{-1}u_{12}^{-1}\\
&\implies&(t\otimes\Delta)u^{-1}=(t\otimes id)u_{12}^{-1}\cdot(t\otimes id)u_{13}^{-1}\\
&\implies&(id\otimes\Delta)\bar{u}=\bar{u}_{12}\bar{u}_{13}
\end{eqnarray*}

(4) Thus, one way or another, we get the result. As further comments, observe first that $S^2=id$ implies $\bar{\bar{u}}=u$. Observe also that when $F=\mathbb C$ and $u$ is unitary, $u^*=u^{-1}$, the matrix $\bar{u}$ is the usual conjugate, given by $\bar{u}_{ij}=u_{ij}^*$. We will be back to this.
\end{proof}

In analogy now with Proposition 4.4, the above result suggests talking about the operation $A_u\to A_{\bar{u}}$ for the affine Hopf algebras, and also investigating the relation between the conditions $A=<u_{ij}>$ and $A=<\bar{u}_{ij}>$. In relation with this, which actually leads to a bit of rethinking of Definition 4.1, let us make the following convention:

\index{affine algebra}
\index{fully affine algebra}

\begin{convention}
Given a Hopf algebra $A$, with a corepresentation $u\in M_N(A)$,
\begin{enumerate}
\item We keep calling $A$ affine when $A=<u_{ij}>$,

\item We call $A$ fully affine when $A=<u_{ij},\bar{u}_{ij}>$,
\end{enumerate}
with the remark that, when $u\sim\bar{u}$, meaning $u=V\bar{u}V^{-1}$, these notions coincide.
\end{convention}

Many further things can be said here, and we will see later in this book that the complex algebra framework, $F=\mathbb C$, brings some clarification in relation with all this, by using $*$-algebras and unitary corepresentations, which are subject to $\bar{u}_{ij}=u_{ij}^*$.

\bigskip

However, importantly, we will also see later in this book, when talking $F=\mathbb C$, that most of the interesting examples satisfy $u\sim\bar{u}$, and that in fact, up to passing to projective versions, which is something quite natural, we can always assume $u\sim\bar{u}$. Thus, and for closing this discussion, Definition 4.1 as stated is basically the good one. 

\bigskip

Getting now to what can be done with an affine Hopf algebra, we will use as source of inspiration what happens for $A=F(G)$. Given an algebraic group $G\subset GL_N(F)$, a natural construction is that of considering its diagonal torus $T\subset G$, which is given by the following formula, $(F^*)^N\subset GL_N(F)$ being the subgroup of diagonal matrices:
$$T=G\cap(F^*)^N$$

The point now is that we can perform in fact this construction in the general affine Hopf algebra context, that of Definition 4.1, with the result being as follows:

\index{diagonal algebra}
\index{diagonal torus}

\begin{theorem}
Given an undeformed affine Hopf algebra $(A,u)$, the quotient
$$A^\delta=A\Big/\left<u_{ij}=\bar{u}_{ij}=0\Big|\forall i\neq j\right>$$
is an affine Hopf algebra too, called diagonal algebra. Its standard generators 
$$u_{ii}\in A^\delta$$
are group-like, and the algebra $A^\delta$ itself is cocommutative.
\end{theorem}

\begin{proof}
This is something very standard, the idea being as follows:

\medskip

(1) We know from chapter 3 that given a Hopf algebra $A$, so is its quotient $B=A/I$, provided that $I\subset A$ is an ideal satisfying the Hopf ideal conditions, namely:
$$\Delta(I)\subset A\otimes I+I\otimes A\quad,\quad\varepsilon(I)=0\quad,\quad S(I)\subset I$$

In our case, the ideal that we are dividing by is given by the following formula:
$$I=\left<u_{ij},\bar{u}_{ij}\Big|i\neq j\right>$$

(2) So, let us check that this is indeed a Hopf ideal. Regarding the condition involving the comultiplication, on the main generators of $I$ we have, as desired:
\begin{eqnarray*}
\Delta(u_{ij})
&=&\sum_ku_{ik}\otimes u_{kj}\\
&=&u_{ii}\otimes u_{ij}+\sum_{k\neq i}u_{ik}\otimes u_{kj}\\
&\in&A\otimes I+\sum_{k\neq i}I\otimes A\\
&=&A\otimes I+I\otimes A
\end{eqnarray*}

Similarly, in what regards the coefficients of $\bar{u}=(t\otimes id)u^{-1}$, we have:
$$\Delta(\bar{u}_{ij})\in A\otimes I+I\otimes A$$

(3) Next, the condition involving the counit is clear as well, because we have:
$$i\neq j\implies\varepsilon(u_{ij})=\varepsilon(\bar{u}_{ij})=0$$

In what regards now the last condition, involving the antipode, here we can use the following formulae, with the first one being something that we know well, and with the second one coming from it, by using our undeformability assumption $S^2=id$: 
$$\bar{u}=(t\otimes S)u\quad,\quad u=(t\otimes S)\bar{u}$$

Indeed, we obtain from these formulae that for any $i\neq j$ we have, as desired:
$$S(u_{ij})=\bar{u}_{ji}\in I\quad,\quad 
S(\bar{u}_{ij})=u_{ji}\in I$$

(4) We conclude from all this that $A^\delta=A/I$ is indeed an affine Hopf algebra, with its fundamental corepresentation being as follows, with the convention $u_{ii}\in A^\delta$:
$$u^\delta=\begin{pmatrix}
u_{11}\\
&\ddots\\
&&u_{NN}
\end{pmatrix}$$

(5) Regarding now the last assertion, this is clear from the above diagonal form of the fundamental corepresentation, but we can check this directly too. We have, inside $A^\delta$:
\begin{eqnarray*}
\Delta(u_{ii})
&=&\sum_ku_{ik}\otimes u_{ki}\\
&=&u_{ii}\otimes u_{ii}+\sum_{k\neq i}u_{ik}\otimes u_{ki}\\
&=&u_{ii}\otimes u_{ii}+\sum_{k\neq i}0\otimes0\\
&=&u_{ii}\otimes u_{ii}
\end{eqnarray*}

Thus the standard generators are indeed group-like, and this implies of course that the diagonal algebra $A^\delta$ itself is indeed cocommutative, as stated.

\medskip

(6) Finally, as a complement to this, in relation with what was said before the statement, for $A=F(G)$ with $G\subset GL_N(F)$ we obtain $A^\delta=F(T)$, with $T\subset G$ being the diagonal torus. Also, for $A=F[H]$, with $H=<g_1,\ldots,g_N>$ being a finitely generated group, with fundamental corepresentation $u=diag(g_1,\ldots,g_N)$, we obtain $A^\delta=A$.
\end{proof}

As an interesting variation of the above construction, generalizing it, we have:

\index{spinned diagonal algebra}
\index{spinned tori}

\begin{theorem}
Given an undeformed affine Hopf algebra $(A,u)$ and $Q\in GL_N(F)$,
$$A^\delta_Q=A\Big/\left<(QuQ^{-1})_{ij}=0\Big|\forall i\neq j\right>$$
is an affine Hopf algebra too, called spinned diagonal algebra. Its standard generators 
$$(QuQ^{-1})_{ii}\in A^\delta_Q$$
are group-like, and the algebra $A^\delta_Q$ itself is cocommutative.
\end{theorem}

\begin{proof}
This follows indeed from Theorem 4.7 applied to the following affine Hopf algebra, with this latter Hopf algebra coming from Proposition 4.4 (3):
$$(A',u')=(A,QuQ^{-1})$$

Alternatively, we can redo the proof of Theorem 4.7, in the present more general setting, by adding the parameter matrix $Q\in GL_N(F)$, to all the computations there.
\end{proof}

All the above is quite interesting, making a connection with advanced group theory, and more specifically, with the notion of maximal torus, from the Lie group theory. Indeed, at the level of basic examples, we first have the following result:

\begin{theorem}
For a function algebra $A=F(G)$, with $G\subset GL_N(F)$, the diagonal algebra $A\to A^\delta$ is the algebra of functions on the diagonal torus $T=G\cap(F^*)^N$:
$$A^\delta=F(T)$$
More generally, the spinned diagonal algebras $A\to A^\delta_Q$,  with $Q\in GL_N(F)$, are the algebras of functions on the spinned diagonal tori, $T_Q=G\cap Q^{-1}(F^*)^NQ$:
$$A^\delta_Q=F(T_Q)$$
Also, when $F=\mathbb C$, any abelian subgroup $H\subset G$ appears inside such a torus, $H\subset T_Q$.
\end{theorem}

\begin{proof}
This is something quite self-explanatory, the idea being as follows:

\medskip

(1) The first assertion, regarding the diagonal torus, is something routine, coming from definitions, as explained at the end of the proof of Theorem 4.7.

\medskip

(2) The second assertion, regarding the spinned tori, is clear from definitions too, via the same argument, with the above definition for the spinned tori.

\medskip

(3) Finally, in what regards the third assertion, given an abelian subgroup $H\subset G$, we can compose this group embedding with the embedding $G\subset GL_N(F)$:
$$H\subset G\subset GL_N(F)$$

Now the group $H$ being assumed to be abelian, when $F=\mathbb C$, by representation theory, this embedding must be similar to a diagonal embedding. Thus, if we denote by $Q\in GL_N(F)$ the matrix producing this similarity, we have an embedding as follows:
$$H\subset Q^{-1}(F^N)^*Q$$

But this proves our claim, because by intersecting with $G$, we obtain:
$$H\subset G\cap Q^{-1}(F^N)^*Q=T_Q$$

So, this was for the general idea, and in practice, we will be back to this later, with full details and explanations, and with some generalizations too.
\end{proof}

Let us discuss as well the case of the group algebras. Here we have:

\begin{theorem}
For a group algebra $A=F[H]$, with $H=<g_1,\ldots,g_N>$, the diagonal algebra $A\to A^\delta$ is the group algebra of the following quotient group,
$$H_Q=H\Big/\left<g_i=g_j\Big|\exists k,Q_{ki}\neq0,Q_{kj}\neq0\right>$$
with the embedding $T_Q=\widehat{H}_Q\subset G=\widehat{H}$ coming from the quotient map $H\to H_Q$. Also, a similar result holds for $A=F[H]$, with spinned fundamental corepresentation.
\end{theorem}

\begin{proof}
This is something elementary, the idea being as follows:

\medskip

(1) Assume first that $H=<g_1,\ldots,g_N>$ is a discrete group, with dual $\widehat{H}$ diagonally embedded, that is, with fundamental corepresentation of $F[H]$ as follows:
$$u=\begin{pmatrix}
g_1\\
&\ldots\\
&&g_N
\end{pmatrix}$$

With $v=QuQ^{-1}$, we have then the following computation:
\begin{eqnarray*}
\sum_s(Q^{-1})_{si}v_{sk}
&=&\sum_s(Q^{-1})_{si}(QuQ^{-1})_{sk}\\
&=&\sum_{st}(Q^{-1})_{si}Q_{st}u_{tt}(Q^{-1})_{kt}\\
&=&\sum_{st}(Q^{-1})_{si}Q_{st}(Q^{-1})_{kt}g_t\\
&=&\sum_t\delta_{it}(Q^{-1})_{kt}g_t\\
&=&(Q^{-1})_{ki}g_i
\end{eqnarray*}

Thus the condition $v_{ij}=0$ for $i\neq j$, used in Theorem 4.8, gives:
$$(Q^{-1})_{ki}v_{kk}=(Q^{-1})_{ki}g_i$$

(2) But this latter condition tells us that we must have:
$$Q_{ki}\neq0\implies g_i=v_{kk}$$

We conclude from this that we have, as desired:
$$Q_{ki}\neq0,Q_{kj}\neq0\implies g_i=g_j$$

(3) In order to finish now, consider the group in the statement. We must prove that the off-diagonal coefficients of $QuQ^{-1}$ vanish. So, let us look at these coefficients:
$$(QuQ^{-1})_{ij}=\sum_kQ_{ik}u_{kk}(Q^{-1})_{kj}=\sum_kQ_{ik}(Q^{-1})_{kj}g_k$$

In this sum $k$ ranges over the set $S=\{1,\ldots,N\}$, but we can restrict the attention to the subset $S'$ of indices having the property $Q_{ik}(Q^{-1})_{kj}\neq 0$.  But for these latter indices the elements $g_k$ are all equal, say to an element $g$, and we obtain, as desired:
\begin{eqnarray*}
(QuQ^{-1})_{ij}
&=&\left(\sum_{k\in S'}Q_{ik}(Q^{-1})_{kj}\right)g\\
&=&\left(\sum_{k\in S}Q_{ik}(Q^{-1})_{kj}\right)g\\
&=&(QQ^{-1})_{ij}g\\
&=&\delta_{ij}g
\end{eqnarray*}

(4) Finally, in what regards the last assertion, this is again elementary, obtained by adding an extra matrix parameter to the above computations, spinning the fundamental corepresentation of $F[H]$, and we will leave the computations here as an exercise.
\end{proof}

Summarizing, we have here some beginning of Lie theory, for the affine Hopf algebras, going beyond what we know from chapter 2. According to the above results, we can expect the collection of tori $\{T_Q|Q\in GL_N(F)\}$ to encode various algebraic and analytic properties of $G$. We will discuss this later, with a number of results and conjectures.

\section*{4b. Projective versions}

As already mentioned in the beginning of this chapter, the affine Hopf algebras are subject to far more operations than those discussed in chapter 3, for the arbitrary Hopf algebras, and this makes the world of affine Hopf algebras quite interesting. 

\bigskip

In fact, we have already seen a number of such new operations, in Theorem 4.7 and Theorem 4.8, making an interesting link with the advanced Lie group theory. So, let us further explore this subject, what operations can be applied to the affine Hopf algebras. We can first talk about complexifications, in an abstract sense, as follows:

\index{complexification}

\begin{theorem}
Given an affine Hopf algebra $(A,u)$, we can construct its complexification $(A^c,u^c)$ as follows,
$$A^c=<u^c>\subset F[\mathbb Z]\otimes A\quad,\quad u^c=zu$$
with $z=1\in F[\mathbb Z]$. As main illustrations for this operation, we have:
\begin{enumerate}
\item $F(G)^c=F(G^c)$, for a certain group quotient $\mathbb T\times G\to G^c$.

\item $F[H]^c=F[H^c]$, with $H^c\subset\mathbb Z\times H$ being constructed similarly.
\end{enumerate}
\end{theorem}

\begin{proof}
This is something quite routine, the idea being as follows:

\medskip

(1) Regarding the Hopf algebra assertion, this follows from our general results regarding the tensor products and subalgebras, from chapter 3. Indeed, we have:
\begin{eqnarray*}
\Delta(u^c_{ij})
&=&\Delta(z)\Delta(u_{ij})\\
&=&(z\otimes z)\sum_ku_{ik}\otimes u_{kj}\\
&=&\sum_k zu_{ik}\otimes zu_{kj}\\
&=&\sum_ku^c_{ik}\otimes u^c_{kj}
\end{eqnarray*}

Thus $\tilde{u}=zu$ is indeed a corepresentation of the tensor product algebra $F[\mathbb Z]\otimes F(G)$, as constructed in chapter 3, so the results there apply, and gives the result.

\medskip

(2) As a comment here, by Fourier transform, we can define alternatively the complexification $A^c$ as follows, with $z=id\in F[\mathbb T]$ being the standard generator:
$$A^c=<u^c>\subset F(\mathbb T)\otimes A\quad,\quad u^c=zu$$

(3) In what regards now the formula $F(G)^c=F(G^c)$, with $\mathbb T\times G\to G^c$, many things can be said here, and we will leave some study here as an exercise.

\medskip

(4) As for the formula $F[H]^c=F[H^c]$, with the subgroup $H^c\subset\mathbb Z\times H$ being constructed similarly, by multiplying the generators by $z$, this is something self-explanatory too. Again, many things can be said here, and we will leave this as an exercise.
\end{proof}

As a comment here, the above construction is particularly relevant when the ground field is $F=\mathbb C$. We will be back to this later, with comments and illustrations. 

\bigskip

Now still at the general level, we have a result regarding the Haar integration of the complexifications, which is something quite straightforward, as follows:

\begin{theorem}
The Haar integral of a complexification
$$A^c=<u^c>\subset F[\mathbb Z]\otimes A\quad,\quad u^c=zu$$
appears as a restriction of a tensor product, as follows,
$$\int_{A^c}=\left(\int_\mathbb T\otimes\int_A\right)_{A^c}$$
and with this happening for left integrals, right integrals, and integrals.
\end{theorem}

\begin{proof}
This is something quite self-explanatory, and clear from our various results from chapter 3, with the corresponding commuting diagram here being as follows:
$$\xymatrix@R=40pt@C=70pt{
F[\mathbb Z]\otimes A\ar[r]^{\int_\mathbb T\otimes\int_A}&F\\
A^c\ar[u]^i\ar[ur]_{\int_{A^c}}
}$$

Thus, we are led to the conclusion in the statement.
\end{proof}

Moving on, along the same lines, we have as well the following construction:

\index{free complexification}

\begin{theorem}
Given an affine Hopf algebra $(A,u)$, we can construct its free complexification $(\widetilde{A},\tilde{u})$ as follows,
$$\widetilde{A}=<\tilde{u}>\subset F[\mathbb Z]*A\quad,\quad\tilde{u}=zu$$
with $z=1\in F[\mathbb Z]$. As main illustrations for this operation, we have:
\begin{enumerate}
\item $\widetilde{F}(G)=F(\widetilde{G})$, standing as definition for $\widetilde{G}$, as quantum group.

\item $\widetilde{F}[H]=F[\widetilde{H}]$, with $\widetilde{H}\subset\mathbb Z*H$ being constructed similarly.
\end{enumerate}
\end{theorem}

\begin{proof}
This is something quite similar, the idea being as follows:

\medskip

(1) Regarding the Hopf algebra assertion, this follows from our general results regarding free products and subalgebras, from chapter 3. Indeed, we have:
\begin{eqnarray*}
\Delta(\tilde{u}_{ij})
&=&\Delta(z)\Delta(u_{ij})\\
&=&(z\otimes z)\sum_ku_{ik}\otimes u_{kj}\\
&=&\sum_k zu_{ik}\otimes zu_{kj}\\
&=&\sum_k\tilde{u}_{ik}\otimes\tilde{u}_{kj}
\end{eqnarray*}

Thus $\tilde{u}=zu$ is indeed a corepresentation of the free product algebra $F[\mathbb Z]*F(G)$, as constructed in chapter 3, so the results there apply, and gives the result.

\medskip

(2) As a comment here, by Fourier transform, we can define alternatively the free complexification $\widetilde{A}$ as follows, with $z=id\in F[\mathbb T]$ being the standard generator:
$$\widetilde{A}=<\tilde{u}>\subset F(\mathbb T)*A\quad,\quad\tilde{u}=zu$$

(3) In what regards now $\widetilde{F}(G)=F(\widetilde{G})$, this stands as a definition for $\widetilde{G}$, as a quantum group, the point being that, unless $G$ is trivial, the algebra $F(\widetilde{G})$ is not commutative. We will be back to this, with some examples and illustrations, in what follows.

\medskip

(4) As for the formula $\widetilde{F}[H]=F[\widetilde{H}]$, with the subgroup $\widetilde{H}\subset\mathbb Z*H$ being constructed similarly, by multiplying the generators by $z$, this is something self-explanatory. Again, we will be back to this, with some examples and illustrations, in what follows.
\end{proof}

As a comment here, the above construction is particularly relevant when the ground field is $F=\mathbb C$. We will be back to this later in this book, with comments and illustrations. Now still at the general level, we have a result regarding the Haar integration, which is quite similar to what we had before for the usual complexifications, as follows:

\begin{theorem}
The Haar integral of a free complexification
$$\widetilde{A}=<\tilde{u}>\subset F[\mathbb Z]*A\quad,\quad\tilde{u}=zu$$
appears as a restriction of a free product, as follows,
$$\int_{\widetilde{A}}=\left(\int_\mathbb T*\int_A\right)_{\widetilde{A}}$$
and with this happening for left integrals, right integrals, and integrals.
\end{theorem}

\begin{proof}
This is something quite self-explanatory, and clear from our various results from chapter 3, with the corresponding commuting diagram here being as follows:
$$\xymatrix@R=40pt@C=70pt{
F[\mathbb Z]*A\ar[r]^{\int_\mathbb T*\int_A}&F\\
\widetilde{A}\ar[u]^i\ar[ur]_{\int_{\widetilde{A}}}
}$$

Thus, exactly as before, when talking about the usual complexifications, and their Haar functionals, we are led to the conclusion in the statement.
\end{proof}

The above results raise the question of understanding what are the ``intermediate'' complexifications, lying between the usual one, and the free one. More on this later.

\bigskip

Moving ahead, we can introduce now our next basic operation on the affine Hopf algebras, which is something quite fundamental, as follows:

\index{projective version}
\index{projective quantum group}
\index{left projective version}
\index{right projective version}

\begin{theorem}
Given an affine Hopf algebra $(A,u)$, we can construct its projective version $(PA,v)$ by setting
$$PA=<v_{ia,jb}>\subset A\quad,\quad v_{ia,jb}=u_{ij}\bar{u}_{ab}$$
and as main illustrations for this construction, we have the following formulae:
\begin{enumerate}
\item $PF(G)=F(PG)$, with $PG=G/(G\cap\mathbb T^N)$, when $F=\mathbb C$.

\item $PF[H]=F[PH]$, with $PH=<g_ig_j^{-1}>$, assuming $H=<g_i>$.
\end{enumerate}
\end{theorem}

\begin{proof}
As before, this is something self-explanatory, the idea being as follows:

\medskip

(1) Our first claim is that the matrix $v$ in the statement is a corepresentation. But this is something standard, the computation being as follows:
\begin{eqnarray*}
\Delta(v_{ia,jb})
&=&\Delta(u_{ij}\bar{u}_{ab})\\
&=&\Delta(u_{ij})\Delta(\bar{u}_{ab})\\
&=&\sum_ku_{ik}\otimes u_{kj}\sum_c\bar{u}_{ac}\otimes\bar{u}_{cb}\\
&=&\sum_{kc}u_{ik}\bar{u}_{ac}\otimes u_{kj}\bar{u}_{cb}\\
&=&\sum_{kc}v_{ia,kc}\otimes v_{kc,jb}
\end{eqnarray*}

(2) Thus, $v$ is indeed a corepresentation, and so, by the results in chapter 3, the projective version $PA$, as constructed in the statement, is indeed a Hopf subalgebra.

\medskip

(3) Before going further with our study, let us mention that the construction in the statement is that of the standard, left projective version. It is possible to talk as well about right projective versions, constructed by using the following corepresentation:
$$w_{ia,jb}=\bar{u}_{ij}u_{ab}$$

In general, the left and right projective versions do not coincide, as one can see with examples coming from algebras of type $A=F[H]$, discussed below. Many other things can be said here, but with the subject being quite technical, we will basically restrict the attention in what follows to the left projective versions, from the statement.

\medskip

(4) Regarding now the formula $PF(G)=F(PG)$, with $PG=G/(G\cap\mathbb T^N)$, this follows from the elementary fact that, via Gelfand duality, the matrix $v$ in the statement is the matrix of coefficients of the adjoint representation of $G$, whose kernel is the subgroup $G\cap\mathbb T^N$, where $\mathbb T^N\subset U_N$ denotes the subgroup formed by the diagonal matrices.

\medskip

(5) So, this was for the idea with $PF(G)=F(PG)$, and in practice, we will leave some study here, both over $F=\mathbb C$ and in general, as an instructive exercise.

\medskip

(6) As for the last formula, namely $PF[H]=F[PH]$, with $PH=<g_ig_j^{-1}>$, assuming $H=<g_i>$, this is something trivial, which comes from definitions.
\end{proof}

We have as well a result regarding the Haar integration, as follows:

\begin{theorem}
The Haar integral of a projective version
$$PA=<v_{ia,jb}>\subset A\quad,\quad v_{ia,jb}=u_{ij}\bar{u}_{ab}$$
appears as a restriction, as follows,
$$\int_{PA}=\left(\int_A\right)_{PA}$$
and with this happening for left integrals, right integrals, and integrals.
\end{theorem}

\begin{proof}
This is again something quite self-explanatory, and clear from our results from chapter 3, with the corresponding commuting diagram here being as follows:
$$\xymatrix@R=40pt@C=70pt{
A\ar[r]^{\int_A}&F\\
PA\ar[u]^i\ar[ur]_{\int_{PA}}
}$$

Thus, we are led to the conclusion in the statement.
\end{proof}

Importantly now, let us mention that there is an interesting relation here with the notion of free complexification. Indeed, in the context of Theorem 4.13, we have:

\begin{theorem}
Given an affine Hopf algebra $(A,u)$, construct its free complexification $(\widetilde{A},\tilde{u})$, with $\tilde{u}=zu$. We have then an identification
$$P\widetilde{A}=PA$$
and the same happens at the level of right projective versions.
\end{theorem}

\begin{proof}
This is something coming from definitions, the idea being as follows:

\medskip

(1) Let us construct indeed the free complexification $(\widetilde{A},\tilde{u})$, with $\tilde{u}=zu$, as in Theorem 4.13. The conjugate of $\bar{u}$ is then given by the following formula:
$$\bar{\tilde{u}}=\bar{u}z^{-1}$$

Thus, the adjoint of $\bar{u}$ is given by the following formula:
$$\tilde{v}=zvz^{-1}$$

But this gives an identification $P\widetilde{A}=PA$ as in the statement, by conjugating by $z$.

\medskip

(2) Regarding now the right projective versions, as constructed in the proof of Theorem 4.15, things here are in fact even simpler. Indeed, the right adjoint of $\bar{u}$ is given by:
$$\tilde{w}=w$$

Thus, we have a plain equality of right projective versions $P\widetilde{A}=PA$.

\medskip

(3) Finally, let us mention that the story is not over here, quite the opposite. Indeed, an interesting subsequent question, of categorical and algebraic geometry flavor, is that of understanding if the free complexification $\widetilde{A}$ is the biggest affine Hopf algebra having the same projective version as $A$. We will discuss this question later in this book, when systematically discussing the relation between affine and projective geometry.
\end{proof}

So long for the basic operations on the Hopf algebras. Many more things can be said, and we will be back to this later, when discussing the representation theory and Haar integration, for the various products of Hopf algebras constructed above.

\bigskip

Finally, let us mention that it is possible to talk as well about wreath products, and free wreath products, in the Hopf algebra setting. However, this is something more tricky, requiring talking about quantum permutations first, and we will do this later in this book, after developing some more theory, in order to talk about quantum permutations.

\section*{4c. Intersection, generation}

We would like to discuss now a number of further basic operations on the Hopf algebras, again coming from group theory, but which are more subtle, and related to each other, namely the intersection operation, generation operation, and Hopf image operation. 

\bigskip

As before with the usual operations, it is convenient in what follows to have in mind the following informal formula, for an arbitrary Hopf algebra $A$, in terms of certain underlying quantum groups $G$ and $H$, related by some sort of generalized Pontrjagin duality:
$$A=F(G)=F[H]$$

In terms of such quantum groups, and more specifically of the first ones, $A=F(G)$, the questions that we would like to solve, which are quite natural, are as follows: 

\begin{problem}
Given two quantum subgroups $H,K\subset G$ of a quantum group:
\begin{enumerate}
\item How to define their intersection, $H\cap K$?

\item What about the subgroup that they generate, $<H,K>$?
\end{enumerate}
\end{problem}

In order to answer the first question, the best is to start by drawing some diagrams. In the classical case, given a group $G$ and two subgroups $H,K\subset G$, we can indeed intersect these subgroups, and the relevant diagram, which gives $H\cap K$, is as follows:
$$\xymatrix@R=50pt@C=50pt{
G&H\ar[l]\\
K\ar[u]&H\cap K\ar[l]\ar[u]}$$

Now at the level of the corresponding algebras of functions, the diagram, having as arrows surjections which are dual to the above inclusions, is as follows:
$$\xymatrix@R=50pt@C=50pt{
F(G)\ar[r]\ar[d]&F(H)\ar[d]\\
F(K)\ar[r]&F(H\cap K)}$$

But this is exactly what we need, in order to solve our intersection problem formulated above. Indeed, based on this, we can come up with the following solution:

\index{intersection of quantum groups}

\begin{theorem}
Given two quotient Hopf algebras $A\to B,C$, we can construct the universal Hopf algebra quotient $A\to B\sqcap C$ producing the following diagram:
$$\xymatrix@R=50pt@C=50pt{
A\ar[r]\ar[d]&B\ar[d]\\
C\ar[r]&B\sqcap C}$$
As an illustration for this, in the group algebra case we have the formula
$$F(H)\sqcap F(K)=F(H\cap K)$$
for any two subgroups of a given group, $H,K\subset G$.
\end{theorem}

\begin{proof}
We must prove that the universal Hopf algebra in the statement exists indeed. For this purpose, let us pick writings as follows, with $I,J$ being Hopf ideals:
$$B=A/I\quad,\quad 
C=A/J$$

We can then construct our universal Hopf algebra, as follows:
$$B\sqcap C=A/<I,J>$$

Thus, we are led to the conclusions in the statement.
\end{proof}

In the affine Hopf algebra setting now, that of Definition 4.1, the operation $\sqcap$ constructed above can be usually computed by using the following simple fact:

\begin{proposition}
Assuming $A\to B,C$, the intersection $B\sqcap C$ is given by
$$B\sqcap C=A/\{\mathcal R,\mathcal P\}$$
whenever we have writings as follows,
$$B=A/\mathcal R\quad,\quad 
C=A/\mathcal P$$
with $\mathcal R,\mathcal P$ being certain sets of polynomial relations between the coordinates $u_{ij}$.
\end{proposition}

\begin{proof}
This follows from Theorem 4.19, or rather from its proof, and from the following trivial fact, regarding relations and ideals:
$$I=<\mathcal R>\ ,\ J=<\mathcal P>
\quad\implies\quad <I,J>=<\mathcal R,\mathcal P>$$

Thus, we are led to the conclusion in the statement.
\end{proof}

Finally, let us record as well what happens in the group algebra case:

\begin{theorem}
Given two quotient groups $G\to H,K$, we have the formula
$$F[H]\sqcap F[K]=F[H\sqcap K]$$
with the quotient $G\to[H\sqcap K]$ being the one producing the following diagram:
$$\xymatrix@R=50pt@C=50pt{
G\ar[r]\ar[d]&H\ar[d]\\
K\ar[r]&H\sqcap K}$$
Alternatively, we have $[H\sqcap K]=G/<\ker(G\to H),\ker(G\to K)>$.
\end{theorem}

\begin{proof}
This is indeed something self-explanatory, with the first assertion coming from Theorem 4.19, and with the second assertion coming from Proposition 4.20.
\end{proof}

Moving on, regarding now the generation operation question, from Problem 4.18 (2), the theory here is quite similar. In the classical case, given a group $G$ and two subgroups $H,K\subset G$, the relevant diagram, which gives the subgroup $<H,K>$, is as follows:
$$\xymatrix@R=50pt@C=50pt{
G&H\ar[l]\ar[d]\\
K\ar[u]\ar[r]&<H,K>}$$

Now at the level of the corresponding algebras of functions, the diagram, having as arrows surjections which are dual to the above inclusions, is as follows:
$$\xymatrix@R=50pt@C=50pt{
F(G)\ar[r]\ar[d]&F(H)\\
F(K)&F(<H,K>)\ar[l]\ar[u]}$$

But this is exactly what we need, in order to solve our generation problem formulated above. Indeed, based on this, we can come up with the following solution:

\index{generation operation}
\index{topological generation}

\begin{theorem}
Given two quotient Hopf algebras $A\to B,C$, we can construct the universal Hopf algebra quotient $A\to[B,C]$ producing the following diagram:
$$\xymatrix@R=50pt@C=50pt{
A\ar[r]\ar[d]&B\\
C&[B,C]\ar[l]\ar[u]}$$
As an illustration for this, in the group algebra case we have the formula
$$[F(H),F(K)]=F(<H,K>)$$
for any two subgroups of a given group, $H,K\subset G$.
\end{theorem}

\begin{proof}
We must prove that the universal Hopf algebra in the statement exists indeed. For this purpose, let us pick writings as follows, with $I,J$ being Hopf ideals:
$$B=A/I\quad,\quad 
C=A/J$$

We can then construct our universal Hopf algebra, as follows:
$$B\sqcap C=A/(I\cap J)$$

Thus, we are led to the conclusions in the statement.
\end{proof}

As a complement to this, let us record as well what happens for group algebras:

\begin{theorem}
Given two quotient groups $G\to H,K$, we have the formula
$$\big[F[H],F[K]\big]=F\big[[H,K]\big]$$
with the quotient $G\to[H,K]$ being the one producing the following diagram:
$$\xymatrix@R=50pt@C=50pt{
G\ar[r]\ar[d]&H\\
K&[H,K]\ar[l]\ar[u]}$$
Alternatively, we have $[H,K]=G/(\ker(G\to H)\cap\ker(G\to K))$.
\end{theorem}

\begin{proof}
This is indeed something self-explanatory, and elementary, with both the assertions coming from Theorem 4.22, and its proof.
\end{proof}

So long for the basics of the intersection and generation operations. The story is of course not over with the above results, because we still have as job, as apprentice or confirmed algebraists, to explore the obvious dual nature of these operations. And skipping some details here, that we will leave as an exercise, the situation is as follows:

\begin{fact}
In the finite dimensional case, where the quotients $A\to B,C$ correspond to certain subalgebras $B^*,C^*\subset A^*$:
\begin{enumerate}
\item The intersection and generation operations on the algebras $B,C$ can be understood in terms of the dual algebras $B^*,C^*$.

\item For the algebras of functions, and for group algebras too, all this is compatible with known formulae from group theory.
\end{enumerate}
\end{fact}

As a last topic for this section, let us discuss now the connection with the notion of diagonal algebras, introduced earlier in this chapter. As explained there, associated to an affine Hopf algebra is a family of cocommutative diagonal algebras, as follows:
$$A^\Delta=\left\{A^\delta_Q\Big|Q\in GL_N(F)\right\}$$

Generally speaking, we can expect this collection of diagonal algebras to encode the various algebraic and analytic properties of $A$. Here is a basic result on this subject:

\index{generation conjecture}

\begin{theorem}
The following hold, over the complex numbers $F=\mathbb C$, both for the algebras $A=F(G)$ with $G\subset U_N$, and for the group algebras $A=F[H]$:
\begin{enumerate}
\item Injectivity: the construction $A\to A^\Delta$ is injective, in the sense that $A\neq B$ implies $A^\delta_Q\neq B^\delta_Q$, for some $Q\in GL_N(F)$.

\item Monotony: the construction $A\to A^\Delta$ is increasing, in the sense that passing to a quotient $A\to B$ decreases one of the diagonal algebras, $A^\delta_Q\neq B^\delta_Q$.

\item Generation: any quantum group is generated by its tori, or, equivalently, any affine Hopf algebra $A$ has the property $A=[A^\delta_Q|Q\in GL_N(F)]$.
\end{enumerate}
\end{theorem}

\begin{proof}
We have two cases to be investigated, as follows:

\medskip

(1) Assume first that we are in the group algebra case, $A=F(G)$. In order to prove the generation property we use the following formula, established before: 
$$T_Q=G\cap Q^{-1}(F^*)^NQ$$

Now since any group element $U\in G$ is unitary, and so diagonalizable by basic linear algebra, we can write, for certain matrices $Q\in U_N$ and $D\in(F^*)^N$:
$$U=Q^{-1}DQ$$

But this shows that we have $U\in T_Q$, for this precise value of the spinning matrix $Q\in U_N$, used in the construction of the standard torus $T_Q$. Thus we have proved the generation property, and the injectivity and monotony properties follow from this.

\medskip

(2) Regarding now the group algebras, here everything is trivial. Indeed, when these algebras are diagonally embedded we can take $Q=1$, and when they are embedded by using a spinning matrix $Q\in GL_N(F)$, we can use precisely this matrix $Q$.
\end{proof}

Many other things can be said, as a continuation of the above, notably with some further verifications of the above conjectures, say in relation with product operations, and with some more specialized conjectures too. We will be back to this later in this book, when systematically discussing what happens over the complex numbers.

\section*{4d. Images and models}

In order to further discuss now the generation operation $[\,,]$ introduced above, we will need the following construction, which is something of independent interest:

\index{Hopf image}

\begin{theorem}
Given a representation $\pi:A\to C$, with $C$ being an associative algebra, there is a smallest Hopf algebra quotient $A\to B$ producing a factorization
$$\pi:A\to B\to C$$
called Hopf image of $\pi$. More generally, given representations $\pi:A\to C_i$, with $C_i$ being algebras, there is a smallest Hopf algebra quotient $A\to B$ producing factorizations
$$\pi_i:A\to B\to C_i$$
called  joint Hopf image of the family of representations $\pi_i$.
\end{theorem}

\begin{proof}
This is something quite obvious, obtained by dividing by a suitable ideal:

\medskip

(1) Let $I_\pi$ be sum of all Hopf ideals contained in $\ker(\pi)$. It is clear then that $I_\pi$ is a Hopf ideal, and to be more precise, is  the largest Hopf ideal contained in $\ker(\pi)$. Thus, we have the solution to our factorization problem, obtained as follows:
$$B=A/I_\pi$$

(2) As for the second assertion, regarding the Hopf image of an arbitrary family of representations $\pi_i:A\to C_i$, the proof here is similar, again by dividing by a suitable ideal, obtained as the sum of all Hopf ideals contained in all the kernels $\ker(\pi_i)$.
\end{proof}

The above construction might look quite trivial, but under some suitable extra assumptions, such as having the complex numbers as scalars, $F=\mathbb C$, a number of more subtle things can be said about it, and this even at the very general level, as follows:

\bigskip

(1) To start with, the Hopf image construction has a very simple description at the Tannakian level, namely ``the Hom spaces are those in the model'', and this can be taken as a definition for it. But more on this later, when talking Tannakian duality.

\bigskip

(2) As yet another approach, we can talk about idempotent states, and again, we have a simple description of the Hopf image construction, in such terms. But more on this later in this book, when talking Haar integration, and idempotent states. 

\bigskip

Before going further, with some applications of the above construction to the computation of the $[\,,]$ operation, let us formulate the following definition, which is something of theoretical interest, that will appear on a regular basis, in what follows:

\index{inner faithful}

\begin{definition}
We say that a representation $\pi:A\to C$ is inner faithful when there is no proper factorization of type
$$\pi:A\to B\to C$$
that is, when its Hopf image is $A$ itself.
\end{definition}

As before with the notion of Hopf image, this might look like a trivial notion, but under a number of  some suitable extra assumptions, such as having $F=\mathbb C$, some non-trivial things can be said too, including a simple Tannakian description of inner faithfulness, and an indempotent state formula as well. More on this, later in this book.

\bigskip

Getting now to the examples and illustrations, for the notions introduced above, these will come, as usual, from the case of the function algebras $A=F(G)$, and that of the group algebras $A=F[H]$. It is convenient to start with these latter algebras, $A=F[H]$, which are the most illustrating. Regarding them, we have the following result:

\begin{theorem}
Given a matrix representation of a group algebra
$$\pi:F[H]\to M_N(F)$$
coming by linearizing from a usual group representation
$$\rho:H\to GL_N(F)$$
the Hopf image factorization of $\pi$ is obtained by taking the group image
$$\pi:F[H]\to F[\rho(H)]\to M_N(F)$$
and $\pi$ is inner faithful when $\rho$ is faithful, that is, when $H\subset GL_N(F)$. 
\end{theorem}

\begin{proof}
This is something elementary and self-explanatory, with the first assertion coming from definitions, and with the second assertion coming from it.
\end{proof}

The above result is quite interesting, providing us with the key for understanding the Hopf image construction, and the related notion of inner faithfulness. For making things clear here, let us formulate our final conclusion a bit informally, as follows:

\begin{conclusion}
When regarding the Hopf algebras as being of the form $A=F[H]$, with $H$ being a quantum group:
\begin{enumerate}
\item The Hopf image construction produces the algebra $A'=F[H']$, with the quantum group $H'$ being the image of $H$. 

\item A representation of $A=F[H]$ is inner faithful precisely when the corresponding representation of $H$ is faithful.
\end{enumerate}
\end{conclusion}

With this discussed, let us get now to our other class of basic examples, the function algebras $A=F(G)$. Here the result is something quite simple as well, as follows:

\index{diagonal representation}

\begin{theorem}
Given a diagonal representation of a function algebra
$$\pi:F(G)\to M_N(F)\quad,\quad f\to
\begin{pmatrix}
f(g_1)\\
&\ddots\\
&&f(g_K)
\end{pmatrix}$$
coming from an arbitrary family of group elements, as follows,
$$g_1,\ldots,g_K\in G$$
the Hopf image factorization of $\pi$ is obtained by taking the generated subgroup
$$\pi:F(G)\to F(<g_1,\ldots,g_K>)\to M_N(F)$$
and $\pi$ is inner faithful when we have $G=<g_1,\ldots,g_K>$.
\end{theorem}

\begin{proof}
This is again something elementary and self-explanatory, with the first assertion coming from definitions, and with the second assertion coming from it.
\end{proof}

So long for the notions of Hopf image, and inner faithfulness. There are of course some other examples too, and we will be back to this later. Also, as already mentioned, there is some further general theory to be developed too, under suitable extra assumptions, in relation with Tannakian duals and Haar integration, and we will come back to this later too. For the moment, what we have in the above, which is quite illustrating, will do.

\bigskip

Now by getting back to the generation operation $[\,,]$, as introduced before, we have the following result about it, which is something very useful, in practice:

\index{joint Hopf image}

\begin{theorem}
Assuming $A\to B,C$, the Hopf algebra $[B,C]$ is such that
$$A\to [B,C]\to B,C$$
is the joint Hopf image of the following quotient maps:
$$A\to B,C$$
A similar result holds for an arbitrary family of quotients $A\to B_i$.
\end{theorem}

\begin{proof}
In the particular case from the statement, the joint Hopf image appears as the smallest Hopf algebra quotient $D$ producing factorizations as follows:
$$A\to D\to B,C$$

We conclude from this that we have $D=[B,C]$, as desired. As for the extension to the case of an arbitrary family of quotients $A\to B_i$, this is straightforward.
\end{proof}

As an application of the above Hopf image technology, let us discuss now matrix modeling questions. Let us start with something very basic, as follows:

\index{matrix model}
\index{random matrix algebra}

\begin{definition}
A matrix model for an affine Hopf algebra $(A,u)$ is a morphism of associative algebras as follows, with $T$ being a certain space:
$$\pi:A\to M_K(F(T))\quad,\quad u_{ij}\to U_{ij}$$
When this morphism $\pi$ is an inclusion, we say that our model is faithful.
\end{definition}

Obviously, this notion is potentially something quite useful. In practice now, we would like of course our matrix models to be faithful, and this in order for our computations inside the random matrix algebra $M_K(F(T))$ to be relevant, to our questions regarding $A$. In fact, this is why we chose above to use random matrix algebras $M_K(F(T))$ instead of plain matrix algebras $M_K(F)$, as to have more chances to have faithfulness.

\bigskip

But, the problem is that this situation is not always possible, due to a number of analytic reasons, the idea here being that the random matrix algebras $M_K(F(T))$ are quite ``thin'', from a certain functional analytic viewpoint, while the algebra $A$ to be modeled might be ``thick'', from the same functional analytic viewpoint.

\bigskip

We will discuss more in detail such things later in this book, when talking about $F=\mathbb C$, and various analytic aspects. In the meantime, however, we certainly do have a problem, that is quite clear, and the point is that the notion of inner faithfulness from Definition 4.27 provides us with a potential solution to this problem, as follows:

\begin{definition}
We say that a matrix model as above,
$$\pi:A\to M_K(F(T))\quad,\quad u_{ij}\to U_{ij}$$
is inner faithful when there is no Hopf algebra factorization as follows:
$$\pi:A\to B\to M_K(F(T))\quad,\quad u_{ij}\to v_{ij}\to U_{ij}$$
That is, we can use our inner faithfulness notion, for the matrix models.
\end{definition} 

And the point now is that, with this notion in hand, we can model far more affine Hopf algebras than before, with the above-mentioned analytic obstructions dissapearing. In fact, there is no known obstruction on the algebras $A$ than can be modeled as above.

\bigskip

Still in relation with the faithfulness problematics for models, let us record as well:

\index{inner faithful model}

\begin{theorem}
Given an arbitrary matrix model, as before,
$$\pi:A\to M_K(F(T))\quad,\quad u_{ij}\to U_{ij}$$
we can always factorize it via a smallest Hopf algebra, as follows,
$$\pi:A\to B\to M_K(F(T))\quad,\quad u_{ij}\to v_{ij}\to U_{ij}$$
and the resulting factorized model $B\to M_K(F(T))$ is then inner faithful.
\end{theorem} 

\begin{proof}
This is indeed something self-explanatory, coming by dividing by a suitable ideal, with the result itself being a particular case of Theorem 4.26.
\end{proof}

Summarizing, we have some interesting theory going on, for the matrix models. In practice now, in order to reach to something concrete, out of this, far more work is needed, and we will discuss this later in this book. In the meantime, let us formulate:

\begin{conclusion}
In relation with the matrix models $\pi:A\to M_K(F(T))$:
\begin{enumerate}
\item The notion of inner faithful model is the good one, perfectly remembering $A$, and allowing us to model far more algebras $A$, than with usual faithfulness.

\item In general, we can also use the Hopf image technology in order to construct new Hopf algebras, by taking the Hopf image of an arbitrary model.
\end{enumerate}
\end{conclusion}

And we will end this chapter, and the present opening Part I, with this. Good conclusions that we have here, waiting to be explored. We will be back to this, later.

\section*{4e. Exercises}

We had a lot of interesting algebra in this chapter, and as exercises, we have:

\begin{exercise}
Further study the projective versions.
\end{exercise}

\begin{exercise}
Further study the free complexifications.
\end{exercise}

\begin{exercise}
Work out some basic examples of intersections.
\end{exercise}

\begin{exercise}
Work out some basic examples of generations.
\end{exercise}

\begin{exercise}
Clarify the missing details for the Hopf image illustrations.
\end{exercise}

\begin{exercise}
Work out some basic examples of matrix models.
\end{exercise}

As bonus exercise, figure out what happens to all the above when $F=\mathbb C$.

\part{Symmetry groups}

\ \vskip50mm

\begin{center}
{\em Casa mea e un cantec cu acorduri ample

Melancolic si fierbinte si curat

Si ma simt ca intr-un palat de paduri inconjurat

Viata mea in cantec s-a mutat}
\end{center}

\chapter{Rotation groups}

\section*{5a. Free algebras}

We have seen that many things can be said about Hopf algebras, and the underlying quantum group type objects. In this second part of the present book we discuss some concrete examples, obtained via two basic operations, called liberation and twisting. 

\bigskip

The idea with liberation and twisting, inspired by the results of Drinfeld-Jimbo and Woronowicz, and going back to work of Wang from the early 1990s in the case $F=\mathbb C$, and to work of Bichon from the early 2000s in general, is very simple, as follows:

\bigskip

(1) Liberation. Given an algebraic group $G\subset GL_N(F)$, with a bit of patience and know-how, we can normally liberate the commutative algebra $A=F(G)$ into a highly noncommutative algebra $A^+=F(G^+)$, by removing the commutation relations $ab=ba$ between the standard coordinates $u_{ij}:G\to F$. Thus, we get a liberation $G\subset G^+$.

\bigskip

(2) Half-liberation. Once the liberation constructed, or even before that, a natural question is that of understanding the intermediate liberations $G\subset G^\times\subset G^+$. As a basic example here, we usually have a canonical such liberation, denoted $G^*$, obtained by replacing the relations $ab=ba$ with the half-commutation relations $abc=cba$.

\bigskip

(3) Twisting. Along the same lines, we can try to look for twists $G'\subset G^+$, typically obtained from $G$ by replacing the relations $ab=ba$ by commutation/anticommutation relations $ab=\pm ba$. This is not the same as half-liberation, because in this case $G,G'$ are no longer related by an inclusion, but rather by some sort of mirror symmetry.

\bigskip

(4) Half-liberation and twisting. With the above done, and for the picture to be complete, we can further look for intermediate objects $G'\subset G^\times\subset G^+$. As a basic example here, we usually have a canonical such liberation, denoted $G^{*\prime}$, obtained by replacing the relations $ab=\pm ba$ with relations of type $abc=\pm cba$.

\bigskip

This was for the general story, and in practice now, things are quite tricky, depending on the chosen algebraic group $G\subset GL_N(F)$, the general principle being as follows:

\begin{principle}
The liberation and twisting can be a quite tricky business, depending on the chosen group $G\subset GL_N(F)$, with the problem coming from the ``remove the commutation relations $ab=ba$'' procedure required, which remains something ad-hoc.
\end{principle}

And for the story to be complete, as an illustration for the difficulties here, although the basics of liberation and twisting are now quite well understood in the regular Lie cases, ABCD, the final classification results in the ABCD cases are still lacking. As for the exceptional cases EFG, here nothing much is known, at least for the moment.

\bigskip

Excuse me, but cat is here, meowing something. What is it, cat?

\begin{cat}
The commutation relations $ab=ba$ are best removed by using claws.
\end{cat}

Okay, thanks cat, but is this a joke or something, can't you see that I am a peaceful big animal, having no claws, and quite often eating vegetables. So, I cannot really follow your advice, but one day some sort of generalized Lie theory, or algebraic geometry, will provide us humans with the needed claws, in order to deal easily with such questions.

\bigskip

Getting to work now, as a first concrete goal, in order to have our liberation and twisting program started, we would like to talk about quantum rotation groups over arbitrary fields $F$. Let us start with something that we know since chapter 3, namely:

\begin{definition}
A corepresentation of a Hopf algebra $A$ is a matrix $u\in M_N(A)$ satisfying the following condition, called coassociativity:
$$\Delta(u_{ij})=\sum_ku_{ik}\otimes u_{kj}$$
We say that a Hopf algebra $A$ is affine when it is of the form
$$A=<u_{ij}>$$
with $u\in M_N(A)$ being a certain corepresentation, called fundamental corepresentation.
\end{definition}

We refer to chapter 3 for some basic theory and observations, regarding the corepresentations, and their relation with the equivalent notion of coaction. Among others, it follows from there that we have the following formulae, involving the maps $\Delta,\varepsilon,S$:
$$(id\otimes\Delta)u=u_{12}u_{13}\quad,\quad (id\otimes \varepsilon)u=1\quad,\quad (id\otimes S)u=u^{-1}$$

To be more precise, the first formula is equivalent to the coassociativity condition from Definition 5.3, the third formula follows from this, via the antipode axiom, and the second condition follows from the counit axiom, by using the fact that $u$ is invertible.

\bigskip

As a continuation of this, we have the following key definition:

\index{intertwiner}
\index{Hom space}
\index{End space}
\index{Fix space}
\index{fixed points}

\begin{definition}
Given two corepresentations $u\in M_n(A),v\in M_m(A)$, we set: 
$$Hom(u,v)=\left\{T\in M_{m\times n}(F)\Big|Tu=vT\right\}$$
Also, we use the notations $Fix(v)=Hom(1,v)$, and $End(v)=Hom(v,v)$.
\end{definition}

The linear maps, or rather rectangular matrices $T\in M_{m\times n}(F)$ in the above definition, are called intertwiners between $u$ and $v$. In the particular case $u=1$, what we obtain are certain vectors $T\in F^m$, called fixed points of $v$. Also, in the particular case $u=v$, what we obtain are certain square matrices $T\in M_m(F)$, called intertwiners of $v$.

\bigskip

The above notions are something quite natural and intuitive, and in the function algebra case $A=F(G)$ we obtain the usual notions concerning the representations of $G$, that you might already know about, if you have some previous knowledge of representation theory. In case you don't, no worries, we will be doing things in a self-contained way.

\bigskip

Here is now a fundamental result, regarding the above Hom spaces:

\index{tensor category}

\begin{theorem}
We have the following results:
\begin{enumerate}
\item $T\in Hom(u,v),S\in Hom(v,w)\implies ST\in Hom(u,w)$.

\item $S\in Hom(p,q),T\in Hom(v,w)\implies S\otimes T\in Hom(p\otimes v,q\otimes w)$.
\end{enumerate}
In other words, the Hom spaces between corepresentations form a tensor category.
\end{theorem}

\begin{proof}
These assertions are all elementary, as follows:

\medskip

(1) By using our assumptions $Tu=vT$ and $Sv=Ws$ we obtain, as desired:
$$STu
=SvT
=wST$$

(2) Assume indeed that we have $Sp=qS$ and $Tv=wT$. With standard tensor product notations, we have the following computation:
\begin{eqnarray*}
(S\otimes T)(p\otimes v)
&=&S_1T_2p_{13}v_{23}\\
&=&S_1p_{13}T_2v_{23}\\
&=&(Sp)_{13}(Tv)_{23}
\end{eqnarray*}

On the other hand, we have as well the following computation:
\begin{eqnarray*}
(q\otimes w)(S\otimes T)
&=&q_{13}w_{23}S_1T_2\\
&=&q_{13}S_1w_{23}T_2\\
&=&(qS)_{13}(wT)_{23}
\end{eqnarray*}

The quantities on the right being equal, this gives the result.

\medskip

(3) Finally, the last assertion follows from definitions, and from the obvious fact that, in addition to the above properties (1,2), the Hom spaces are linear spaces, and contain the units. In short, this is just a theoretical remark, that we will need later on.
\end{proof}

As a continuation of Definition 5.4, we have the following key notions:

\begin{definition}
We use the following conventions, for corepresentations:
\begin{enumerate}
\item We write $u\sim v$ when $Hom(u,v)$ contains an invertible element.

\item We say that $u$ is irreducible, and write $u\in Irr(A)$, when $End(u)=F1$.
\end{enumerate}
\end{definition}

To be more precise, here it is understood in (1) that in order to have $u\sim v$, we must have $\dim u=\dim v$, prior to the condition on $Hom(u,v)$ there. Also, in (2), the $F1$ at the end stands for the algebra of constant matrices, $F1\subset M_n(F)$, with $n=\dim u$. 

\bigskip

As before with other notions, in the function algebra case $A=F(G)$ we obtain the usual notions concerning the representations of $G$, that you might already know about. Observe also, as another illustration, that in the group dual case, $A=F[H]$, we have:
$$g\sim h\iff g=h$$

As yet another remark, and in the general case now, observe that the condition $u\sim v$ means by definition that $u,v$ are conjugated by an invertible matrix:
$$u\sim v\iff v=TuT^{-1}\quad,\quad T\in GL_n(F)$$ 

Finally, as a last piece of basic representation theory, that we will need in what follows, let us recall from chapter 3 that we have the following key construction:

\begin{definition}
Given a corepresentation $u\in M_n(A)$, we set
$$\bar{u}=(t\otimes id)u^{-1}=(t\otimes S)u$$
and call it contragradient, or conjugate corepresentation to $u$.
\end{definition}

We refer to chapter 3 for more on this notion, which is something quite subtle. As a key observation here, when $F=\mathbb C$, and assuming that $A$ is a $*$-algebra, and that $u$ is unitary, $u^*=u^{-1}$, the conjugate is given by $\bar{u}_{ij}=u_{ij}^*$. In general, we do not have such an interpretation, but we will still agree to call $\bar{u}$ conjugate corepresentation.

\bigskip

Good news, with the above notions in hand, we can now discuss, following Bichon and others, the liberation question for the most basic matrix group of them all, namely:
$$GL_N=GL_N(F)$$

Let us start with some basic facts about the affine Hopf algebras, as follows:

\begin{proposition}
Given an affine Hopf algebra $(A,u)$, with $u\in M_N(A)$, let
$$\bar{u}=(t\otimes S)u\quad,\quad \bar{\bar{u}}=(id\otimes S^2)u$$
be the conjugate, and double conjugate, of the fundamental corepresentation.
\begin{enumerate}
\item $u=\bar{\bar{u}}$ is equivalent to $S^2=id$.

\item $u\sim\bar{\bar{u}}$ means that $v=\bar{u}$ satisfies $v^{-1}=Qu^tQ^{-1}$, with $Q\in GL_N$.
\end{enumerate} 
\end{proposition}

\begin{proof}
Both assertions are elementary, the idea being as follows:

\medskip

(1) This is something trivial, which is clear from definitions.

\medskip

(2) To start with, the condition $u\sim\bar{u}$ can be interpreted as follows:
$$u\sim\bar{u}\iff \bar{u}=RuR^{-1}\quad,\quad R\in GL_N$$

Now in view of the formula $\bar{u}=(t\otimes id)u^{-1}$, with the notation $v=\bar{u}$, we have:
$$u^{-1}=(t\otimes id)\bar{u}=(t\otimes id)v=v^t$$

But this shows that, still with $v=\bar{u}$ as above, the condition that we found reads:
$$(t\otimes id)v^{-1}
=(t\otimes id)\bar{u}^{-1}
=\bar{\bar{u}}
=RuR^{-1}$$

Now by transposing, we obtain the following formula, with $Q=(R^t)^{-1}$:
$$v^{-1}=Qu^tQ^{-1}$$

Thus, we are led to the conclusion in the statement.

\medskip

(3) Let us mention too, as a complement to what is said in (2), that the condition there is equivalent to having a formula of type $S^2=\Phi*id*\Phi^{-1}$, for the square of the antipode. Many other things can be said here, and we will be back to this later.
\end{proof}

The above result is quite interesting, and we are led in this way to considering the universal Hopf algebra generated by variables $u_{ij},v_{ij}$, with the relations that we found:

\begin{theorem}
Given a matrix $Q\in GL_N$, the universal algebra $F(GL_Q^+)$ generated by variables $u_{ij},v_{ij}$ with $i,j=1,\ldots,N$ with the relations
$$u^{-1}=v^t\quad,\quad v^{-1}=Qu^tQ^{-1}$$
is a Hopf algebra, with comultiplication, counit and antipode making $u,v$ corepresentations. We call the quantum groups $GL_Q^+$ free general linear groups.
\end{theorem}

\begin{proof}
Consider indeed the universal algebra in the statement. We must construct Hopf algebra maps $\Delta,\varepsilon,S$, and this can be done as follows:

\medskip

(1) Our first claim is that we can construct a comultiplication map $\Delta$ by using the universal property of our algebra, via the following formulae:
$$(id\otimes\Delta)u=u_{12}u_{13}\quad,\quad (id\otimes\Delta)v=v_{12}v_{13}$$

In order to prove this, consider the matrices appearing on the right, namely:
$$U=u_{12}u_{13}\quad,\quad V=v_{12}v_{13}$$

By using the formula $u^{-1}=v^t$ in the statement, we have the following computation:
\begin{eqnarray*}
U^{-1}
&=&(u_{12}u_{13})^{-1}\\
&=&u_{13}^{-1}u_{12}^{-1}\\
&=&v_{13}^tv_{12}^t\\
&=&(v_{12}v_{13})^t\\
&=&V^t
\end{eqnarray*}

Similarly, by using $v^{-1}=Qu^tQ^{-1}$, we have as well the following computation:
\begin{eqnarray*}
V^{-1}
&=&(v_{12}v_{13})^{-1}\\
&=&v_{13}^{-1}v_{12}^{-1}\\
&=&(Qu^tQ^{-1})_{13}(Qu^tQ^{-1})_{12}\\
&=&Qu_{13}^tQ^{-1}Qu_{12}^tQ^{-1}\\
&=&Qu_{13}^tu_{12}^tQ^{-1}\\
&=&Q(u_{12}u_{13})^tQ^{-1}\\
&=&QU^tQ^{-1}
\end{eqnarray*}

We conclude that the matrices $U,V$ satisfy the same conditions as $u,v$, so we can define $\Delta$ as above, by using the universal property of our algebra.

\medskip

(2) Regarding now the counit map $\varepsilon$, our claim is this can be defined in a similar way, by using the universal property of our algebra, according to:
$$(id\otimes\varepsilon)u=1\quad,\quad (id\otimes\varepsilon)v=1$$

But this is clear, because the matrices $U=1$, $V=1$ obviously satisfy the required relations $U^{-1}=V^t$ and $V^{-1}=QU^tQ^{-1}$. Thus, we have our counit map $\varepsilon$.

\medskip

(3) Finally, regarding the antipode $S$, our claim is that this can be defined in a similar way, via the universal property of our algebra, according to the following formulae:
$$(id\otimes S)u=v^t\quad,\quad (id\otimes S)v=Qu^tQ^{-1}$$

In order to prove this, consider the matrices $U=v^t$ and $V=Qu^tQ^{-1}$, over the opposite algebra, with product $a\cdot b=ba$. We have the following computation:
\begin{eqnarray*}
(U\cdot V^t)_{ij}
&=&\sum_kU_{ik}\cdot V_{jk}\\
&=&\sum_kV_{jk}U_{ik}\\
&=&(VU^t)_{ji}\\
&=&(Qu^tQ^{-1}v)_{ji}\\
&=&(v^{-1}v)_{ji}\\
&=&\delta_{ji}
\end{eqnarray*}

Similarly, we have as well the following computation:
\begin{eqnarray*}
(V^t\cdot U)_{ij}
&=&\sum_kV_{ki}\cdot U_{kj}\\
&=&\sum_kU_{kj}V_{ki}\\
&=&(U^tV)_{ji}\\
&=&(vQu^tQ^{-1})_{ji}\\
&=&(vv^{-1})_{ji}\\
&=&\delta_{ji}
\end{eqnarray*}

We conclude that the first required relation, $U^{-1}=V^t$, is satisfied indeed. 

\medskip

(4) Next, still in relation with the antipode, we have the following computation:
\begin{eqnarray*}
(V\cdot QU^tQ^{-1})_{ij}
&=&\sum_kV_{ik}\cdot(QU^tQ^{-1})_{kj}\\
&=&\sum_k(QU^tQ^{-1})_{kj}V_{ik}\\
&=&[(QU^tQ^{-1})^tV^t]_{ji}\\
&=&[(QvQ^{-1})^t(Qu^tQ^{-1})^t]_{ji}\\
&=&[(Q^{-1})^tv^tQ^t(Q^{-1})^tuQ^t]_{ji}\\
&=&[(Q^{-1})^tv^tuQ^t]_{ji}\\
&=&[(Q^{-1})^tu^{-1}uQ^t]_{ji}\\
&=&[(Q^{-1})^tQ^t]_{ji}\\
&=&\delta_{ji}
\end{eqnarray*}

Similarly, we have as well the following computation:
\begin{eqnarray*}
(QU^tQ^{-1}\cdot V)_{ij}
&=&\sum_k(QU^tQ^{-1})_{ik}\cdot V_{kj}\\
&=&\sum_kV_{kj}(QU^tQ^{-1})_{ik}\\
&=&[V^t(QU^tQ^{-1})^t]_{ji}\\
&=&[(Qu^tQ^{-1})^t(QvQ^{-1})^t]_{ji}\\
&=&[(Q^{-1})^tuQ^t(Q^{-1})^tv^tQ^t]_{ji}\\
&=&[(Q^{-1})^tuv^tQ^t]_{ji}\\
&=&[(Q^{-1})^tuu^{-1}Q^t]_{ji}\\
&=&[(Q^{-1})^tQ^t]_{ji}\\
&=&\delta_{ji}
\end{eqnarray*}

We conclude that the second required relation, namely $V^{-1}=QU^tQ^{-1}$, is satisfied too. Thus that the matrices $U,V$ satisfy the same conditions as $u,v$, so we can define the antipode map $S$ as above, by using the universal property of our algebra.
\end{proof}

Many things can be said about quantum groups $GL_Q^+$ constructed above. As a first observation, the quantum group $GL_Q^+$ depends on the matrix $Q\in GL_N$ up to $Q\to\lambda Q$, and in view of this, it is convenient to make the following normalization, provided that $Tr(Q)\neq0$, and that the number $Tr(Q^{-1})/Tr(Q)$ has a square root in $F$:
$$Tr(Q)=Tr(Q^{-1})$$

Along the same lines, but at a more advanced level, the quantum groups $GL_Q^+$ are invariant under certain transformations of type $Q\to AQB$, at the level of the parameter matrices $Q\in GL_N$, and working out how all this works is an interesting question. We will be back to such things, and to the quantum groups $GL_Q^+$ in general, with further results about them, later in this book, when systematically discussing representation theory.

\bigskip

Let us restrict now the attention to the affine Hopf algebras satisfying $S^2=id$. According to Proposition 5.8 the condition $S^2=id$ is equivalent to $u=\bar{\bar{u}}$, and so to $Q=1$. Thus, we only need the $Q=1$ particular case of Theorem 5.9, which is as follows:

\begin{theorem}
Given an integer $N\in\mathbb N$, the universal algebra $F(GL_N^+)$ generated by variables $u_{ij},v_{ij}$ with $i,j=1,\ldots,N$ with the relations
$$u^{-1}=v^t\quad,\quad v^{-1}=u^t$$
is a Hopf algebra, with comultiplication, counit and antipode making $u,v$ corepresentations. We call the quantum groups $GL_N^+$ free general linear groups.
\end{theorem}

\begin{proof}
As already mentioned, this is simply the $Q=1$ particular case of Theorem 5.9, with the above quantum group $GL_N^+$ being the ``true'' and unique free general linear group, in the framework of the affine Hopf algebras satisfying $S^2=id$.
\end{proof}

Many other things can be said, as a continuation of the above, in particular with the following observations regarding the diagonal algebras, which get liberated too:

\begin{fact}
In the context of the liberation operation $F(GL_N)\to F(GL_N^+)$:
\begin{enumerate}
\item The corresponding diagonal algebras get liberated too, as $F[\mathbb Z^N]\to F[\mathbb Z^{*N}]$.

\item This latter liberation corresponds to the usual group liberation $\mathbb Z^N\to\mathbb Z^{*N}$.
\end{enumerate}
\end{fact} 

In addition to this, let us mention that $F(GL_N^+)$ is known to be generated by its subalgebras $F(GL_N)$ and $F[\mathbb Z^{*N}]$, so in the end, the liberation operation $GL_N\to GL_N^+$ simply consists in liberating the dual of the diagonal torus, $\mathbb Z^N\to\mathbb Z^{*N}$. However, such things are not exactly well understood and axiomatized, and in the lack of something concrete here, we refer to Cat 5.2 for the general principles of liberation.

\bigskip

We refer to Bichon \cite{bi3} for more on all this. We will be back to the quantum groups $GL_Q^+$ and $GL_N^+$ later in this book, with more results about them. 

\section*{5b. Orthogonal groups}

As a continuation of the above, still following Bichon \cite{bi3} and related work, let us discuss now the free rotation groups. In order to keep things simple, we will restrict the attention to the case where the square of the antipode is the identity:
$$S^2=id$$

Our aim will be that of constructing quantum groups as follows, depending on a parameter matrix $J\in GL_N$, whose meaning will become clear in a moment:
$$O_J^+\subset GL_N^+$$

In order to get started, we first have the following version of Proposition 5.8:

\begin{proposition}
Given an affine Hopf algebra $(A,u)$, with $u\in M_N(A)$, let
$$\bar{u}=(t\otimes S)u$$
be the conjugate of the fundamental corepresentation.
\begin{enumerate}
\item $u=\bar{u}$ is equivalent to $u^{-1}=u^t$.

\item $u\sim\bar{u}$ is equivalent to $u^{-1}=Ju^tJ^{-1}$, with $J\in GL_N$.
\end{enumerate} 
\end{proposition}

\begin{proof}
Both the assertions are clear from $\bar{u}=(t\otimes id)u^{-1}$, as follows:

\medskip

(1) This follows indeed from the following equivalences:
\begin{eqnarray*}
u=\bar{u}
&\iff&u=(t\otimes id)u^{-1}\\
&\iff&(t\otimes id)u=u^{-1}\\
&\iff&u^t=u^{-1}
\end{eqnarray*}

(2) This follows from the following equivalences, with $J=(K^t)^{-1}$ at the end:
\begin{eqnarray*}
u\sim\bar{u}
&\iff&\bar{u}=KuK^{-1}\\
&\iff&(t\otimes id)u^{-1}=KuK^{-1}\\
&\iff&u^{-1}=(K^t)^{-1}u^tK^t\\
&\iff&u^{-1}=Ju^tJ^{-1}
\end{eqnarray*}

Thus, we are led to the conclusions in the statement.
\end{proof}

Thus, we are led to the following ``orthogonal'' version of Theorem 5.10:

\begin{theorem}
Given a matrix $J\in GL_N$, the universal algebra $F(O_J^+)$ generated by variables $u_{ij}$ with $i,j=1,\ldots,N$ with the relations
$$u^{-1}=Ju^tJ^{-1}$$
is a Hopf algebra, with comultiplication, counit and antipode making $u$ a corepresentation. We call the quantum groups $O_J^+$ free rotation groups.
\end{theorem}

\begin{proof}
Consider indeed the universal algebra in the statement. We must construct Hopf algebra maps $\Delta,\varepsilon,S$, and this can be done as follows:

\medskip

(1) Our first claim is that we can construct a comultiplication map $\Delta$ by using the universal property of our algebra, via the following formula:
$$(id\otimes\Delta)u=u_{12}u_{13}$$

In order to prove this, consider the matrix appearing on the right, namely:
$$U=u_{12}u_{13}$$

By using the formula $u^{-1}=Ju^tJ^{-1}$, we have the following computation:
\begin{eqnarray*}
U^{-1}
&=&(u_{12}u_{13})^{-1}\\
&=&u_{13}^{-1}u_{12}^{-1}\\
&=&(Ju^tJ^{-1})_{13}(Ju^tJ^{-1})_{12}\\
&=&Ju_{13}^tJ^{-1}Ju_{12}^tJ^{-1}\\
&=&Ju_{13}^tu_{12}^tJ^{-1}\\
&=&J(u_{12}u_{13})^tJ^{-1}\\
&=&JU^tJ^{-1}
\end{eqnarray*}

We conclude that the matrix $U$ satisfies the same conditions as $u$, so we can define indeed $\Delta$ as above, by using the universal property of our algebra.

\medskip

(2) Regarding now the counit map $\varepsilon$, our claim is this can be defined in a similar way, by using the universal property of our algebra, according to:
$$(id\otimes\varepsilon)u=1$$

But this is clear, because the unit matrix $U=1$ obviously satisfies the required relations, namely $U^{-1}=JU^tJ^{-1}$. Thus, we have our counit map $\varepsilon$.

\medskip

(3) Finally, regarding the antipode $S$, our claim is this can be defined in a similar way, by using the universal property of our algebra, according to the following formula:
$$(id\otimes S)u=Ju^tJ^{-1}$$

In order to prove this, consider the matrix $U=Ju^tJ^{-1}$, over the opposite algebra, with product $a\cdot b=ba$. We have the following computation:
\begin{eqnarray*}
(U\cdot JU^tJ^{-1})_{ij}
&=&\sum_kU_{ik}\cdot(JU^tJ^{-1})_{kj}\\
&=&\sum_k(JU^tJ^{-1})_{kj}U_{ik}\\
&=&[(JU^tJ^{-1})^tU^t]_{ji}\\
&=&[(J^{-1})^tUJ^tU^t]_{ji}\\
&=&[(J^{-1})^tJu^tJ^{-1}J^t(J^{-1})^tuJ^t]_{ji}\\
&=&[(J^{-1})^tJu^tJ^{-1}uJ^t]_{ji}\\
&=&[(J^{-1})^tu^{-1}uJ^t]_{ji}\\
&=&[(J^{-1})^tJ^t]_{ji}\\
&=&\delta_{ij}
\end{eqnarray*}

On the other hand, we have as well the following computation:
\begin{eqnarray*}
(JU^tJ^{-1}\cdot U)_{ij}
&=&\sum_k(JU^tJ^{-1})_{ik}\cdot U_{kj}\\
&=&\sum_kU_{kj}(JU^tJ^{-1})_{ik}\\
&=&[U^t(JU^tJ^{-1})^t]_{ji}\\
&=&[U^t(J^{-1})^tUJ^t]_{ji}\\
&=&[(J^{-1})^tuJ^t(J^{-1})^tJu^tJ^{-1}J^t]_{ji}\\
&=&[(J^{-1})^tuJu^tJ^{-1}J^t]_{ji}\\
&=&[(J^{-1})^tuu^{-1}J^t]_{ji}\\
&=&[(J^{-1})^tJ^t]_{ji}\\
&=&\delta_{ij}
\end{eqnarray*}

We conclude that the matrix $U$ satisfies the condition $U^{-1}=JU^tJ^{-1}$, so we can define the antipode map $S$ as above, by using the universal property of our algebra.
\end{proof}

The story is not over here, because we still have to see what conditions must be imposed on $J\in GL_N$, as to avoid degeneracy of the corresponding quantum group $O_J^+$, a natural requirement here being the fact that $u$ must be irreducible. We will need:

\begin{proposition}
In the context of the quantum group $O_J^+$ we have
$$J^tJ^{-1}\in End(u)$$
and so, in order for $u$ to be irreducible, we must assume $J^t=\pm J$.
\end{proposition}

\begin{proof}
This is something which is quite routine, as follows:

\medskip

(1) We have the following computation, using $u^{-1}=Ju^tJ^{-1}$:
\begin{eqnarray*}
(id\otimes S^2)u
&=&(id\otimes S)u^{-1}\\
&=&(id\otimes S)(Ju^tJ^{-1})\\
&=&J[(id\otimes S)u^t]J^{-1}\\
&=&J[(id\otimes S)u]^tJ^{-1}\\
&=&J(u^{-1})^tJ^{-1}\\
&=&J(Ju^tJ^{-1})^tJ^{-1}\\
&=&J(J^t)^{-1}uJ^tJ^{-1}\\
&=&(J^tJ^{-1})^{-1}uJ^tJ^{-1}
\end{eqnarray*}

Now since $S^2=id$, this equality gives, as claimed in the statement:
$$J^tJ^{-1}\in End(u)$$

(2) Regarding the second assertion, by using what we just found, we have:
\begin{eqnarray*}
u=\,{\rm irreducible}
&\implies&J^tJ^{-1}\in F1\\
&\implies&J^tJ^{-1}=c1,\ c\in F\\
&\implies&J^t=cJ,\ c\in F
\end{eqnarray*}

Now observe that by transposing the formula $J^t=cJ$ we obtain $J=cJ^t$, and by putting these two formulae together, we obtain the following relation:
$$J=c^2J$$

But, $J$ being invertible, we conclude that we must have $c^2=1$, and so $c=\pm1$. Thus, the formula $J^t=cJ$ that we found reads $J^t=\pm J$, as claimed in the statement.
\end{proof}

Now with Proposition 5.14 in mind, let us record the following useful result:

\begin{theorem}
Given a matrix $J\in GL_N$ satisfying $J^t=\pm J$, the corresponding algebra $F(O_J^+)$ is generated by variables $u_{ij}$ with $i,j=1,\ldots,N$ with the relations
$$u=J\bar{u}J^{-1}$$
where $\bar{u}=(t\otimes id)u^{-1}$. We call these quantum groups $O_J^+$ the main free rotation groups, and we use the notation $O_N^+$, for the quantum group obtained at $J=1$.
\end{theorem}

\begin{proof}
This is a remake of Theorem 5.13, with the condition $J^t=\pm J$ found in Proposition 5.14 added, and with the following manipulation included too:
\begin{eqnarray*}
u^{-1}=Ju^tJ^{-1}
&\iff&(t\otimes id)u^{-1}=(J^t)^{-1}uJ^t\\
&\iff&(t\otimes id)u^{-1}=J^{-1}uJ\\
&\iff&\bar{u}=J^{-1}uJ\\
&\iff&u=J\bar{u}J^{-1}
\end{eqnarray*}

To be more precise, we have used here, in the beginning, the formula $J^t=\pm J$. Observe also that our statement makes indeed sense, because $u=J\bar{u}J^{-1}$, with the convention $\bar{u}=(t\otimes id)u^{-1}$, is a collection of relations between the variables $u_{ij}$. Finally, in what regards the ``main'' at the end, this comes from what we found in Proposition 5.14.
\end{proof}

As before in the general linear case, many other things can be said, as a continuation of the above, in particular with the following fact regarding the diagonal algebras:

\begin{fact}
In the context of the liberation operation $F(O_N)\to F(O_N^+)$:
\begin{enumerate}
\item The corresponding diagonal algebras get liberated too, as $F[\mathbb Z_2^N]\to F[\mathbb Z_2^{*N}]$.

\item This latter liberation corresponds to the usual group liberation $\mathbb Z_2^N\to\mathbb Z_2^{*N}$.
\end{enumerate}
\end{fact} 

Observe that this is something very similar to Fact 5.11, and for more about all this, we refer to the comments there, which apply to the present situation as well.

\bigskip

As another comment, importantly, as explained by Bichon in \cite{bi3} and related papers, with the condition $J^t=\pm J$ imposed, $u$ becomes indeed irreducible. We will be back to the quantum groups $O_J^+$ and $O_N^+$ later in this book, with more results about them.

\section*{5c. Bistochastic groups} 

As an interesting version of the orthogonal groups and quantum groups, let us discuss now the bistochastic groups and quantum groups. This is more applied linear algebra, and the very basic definition there, that you might already know, is as follows:

\index{row-stochastic}
\index{column-stochastic}
\index{bistochastic}

\begin{definition}
A square matrix $M\in M_N(F)$ is called bistochastic if each row and each column sum up to the same number:
$$\begin{matrix}
M_{11}&\ldots&M_{1N}&\to&\lambda\\
\vdots&&\vdots\\
M_{N1}&\ldots&M_{NN}&\to&\lambda\\
\downarrow&&\downarrow\\
\lambda&&\lambda
\end{matrix}$$
If this happens only for the rows, or only for the columns, the matrix is called row-stochastic, respectively column-stochastic.
\end{definition}

As the name indicates, these matrices are useful in statistics, and perhaps in other fields like graph theory, computer science and so on, with the case of the real matrices having positive entries, which sum up to $\lambda=1$, being the important one. As a basic example of a bistochastic matrix, we have the flat matrix, which is as follows:
$$\mathbb I_N=\begin{pmatrix}
1&\ldots&1\\
\vdots&&\vdots\\
1&\ldots&1
\end{pmatrix}$$

Observe that the rescaling $P_N=\mathbb I_N/N$ has the property mentioned above, namely positive entries, summing up to $1$. In fact, this matrix $P_N=\mathbb I_N/N$ is a very familiar object in linear algebra, being the projection on the all-one vector, namely:
$$\xi=\begin{pmatrix}
1\\
\vdots\\
1
\end{pmatrix}$$ 

Getting back now to the general case, the various notions of stochasticity in Definition 5.17 are closely related to this vector $\xi$, due to the following simple fact:

\begin{proposition}
Let $M\in M_N(F)$ be a square matrix.
\begin{enumerate}
\item $M$ is row stochastic, with sums $\lambda$, when $M\xi=\lambda\xi$.

\item $M$ is column stochastic, with sums $\lambda$, when $M^t\xi=\lambda\xi$.

\item $M$ is bistochastic, with sums $\lambda$, when $M\xi=M^t\xi=\lambda\xi$.
\end{enumerate}
\end{proposition}

\begin{proof}
The first assertion is clear from definitions, because when multiplying a matrix by $\xi$, we obtain the vector formed by the row sums:
$$\begin{pmatrix}
M_{11}&\ldots&M_{1N}\\
\vdots&&\vdots\\
M_{N1}&\ldots&M_{NN}
\end{pmatrix}
\begin{pmatrix}
1\\
\vdots\\
1
\end{pmatrix}
=\begin{pmatrix}
M_{11}+\ldots+M_{1N}\\
\vdots\\
M_{N1}+\ldots+M_{NN}
\end{pmatrix}$$

As for the second, and then third assertion, these are both clear from this.
\end{proof}

As an observation here, we can reformulate if we want the above statement in a purely matrix-theoretic form, by using the flat matrix $\mathbb I_N$, as follows:

\index{flat matrix}

\begin{proposition}
Let $M\in M_N(F)$ be a square matrix.
\begin{enumerate}
\item $M$ is row stochastic, with sums $\lambda$, when $M\mathbb I_N=\lambda\mathbb I_N$.

\item $M$ is column stochastic, with sums $\lambda$, when $\mathbb I_NM=\lambda\mathbb I_N$.

\item $M$ is bistochastic, with sums $\lambda$, when $M\mathbb I_N=\mathbb I_NM=\lambda\mathbb I_N$.
\end{enumerate}
\end{proposition}

\begin{proof}
This follows from Proposition 5.18, and from the fact that both the rows and columns of the flat matrix $\mathbb I_N$ are copies of the all-one vector $\xi$. Alternatively, we have the following formula, $S_1,\ldots,S_N$ being the row sums of $M$, which gives (1):
$$\begin{pmatrix}
M_{11}&\ldots&M_{1N}\\
\vdots&&\vdots\\
M_{N1}&\ldots&M_{NN}
\end{pmatrix}
\begin{pmatrix}
1&\ldots&1\\
\vdots&&\vdots\\
1&\ldots&1
\end{pmatrix}
=\begin{pmatrix}
S_1&\ldots&S_1\\
\vdots&&\vdots\\
S_N&\ldots&S_N
\end{pmatrix}$$

As for the second, and then third assertion, these are both clear from this.
\end{proof}

Getting now to our quantum group business, we first have the following result:

\index{bistochastic group}
\index{complex bistochastic group}

\begin{proposition}
We have a subgroup as follows,
$$B_J\subset O_J$$
consisting of the orthogonal matrices which are bistochastic.
\end{proposition}

\begin{proof}
It is clear indeed that $B_J$ is stable under the multiplication, contains the unit, and is stable by inversion. Thus, we have indeed a subgroup, as stated.
\end{proof}

In fact, we have the following result, dealing with quantum groups too, and which is formulated a bit informally, by using our general $A=F(G)$ philosophy:

\index{bistochsastic group}
\index{free bistochastic group}
\index{complex bistochastic group}
\index{bistochastic quantum group}

\begin{theorem}
We have the following groups and quantum groups:
\begin{enumerate}
\item $B_J\subset O_J$, consisting of the orthogonal matrices which are bistochastic.

\item $B_J^+\subset O_J^+$, coming via $u\xi=\xi$, where $\xi$ is the all-one vector.
\end{enumerate}
Also, we have an inclusion $B_J\subset B_J^+$, which is a liberation.
\end{theorem}

\begin{proof}
There are several things to be proved, the idea being as follows:

\medskip

(1) We already know from Proposition 5.20 that $B_J$ is indeed a group, with this coming from the following formula, with $\xi$ being the all-one vector:
$$B_J=\left\{U\in O_J\Big|U\xi=\xi\right\}$$

(2) In what regards now $B_J^+$, this appears by definition as follows:
$$C(B_J^+)=F(O_J^+)\Big/\Big<\xi\in Fix(u)\Big>$$

But since the relation $\xi\in Fix(u)$ is categorical, producing a Hopf ideal, and we will leave the verifications here as an exercise, we have indeed quantum groups.

\medskip

(3) Finally, in what regards the last assertion, since we already know that $O_J\subset O_J^+$ is a  liberation, we must prove that we have an isomorphism as follows:
$$F(B_J)=F(O_J)\Big/\Big<\xi\in Fix(u)\Big>$$

But this isomorphism is clear from the formula of $B_J$ in (1).
\end{proof}

Quite interesting all this, but the bistochastic groups and quantum groups are in fact not really ``new'', because at $J=1$ for instance, we have the following result:

\begin{theorem}
We have isomorphisms as follows:
\begin{enumerate}
\item $B_N\simeq O_{N-1}$.

\item $B_N^+\simeq O_{N-1}^+$.
\end{enumerate}
\end{theorem}

\begin{proof}
Let us pick indeed a matrix $K\in GL_N$ satisfying the following condition, where $e_0=(1,0,\ldots,0)$ is the first vector of the standard basis of $F^N$:
$$Ke_0=\xi$$

We have then the following computation, for any corepresentation $u$:
\begin{eqnarray*}
u\xi=\xi
&\iff&uKe_0=Ke_0\\
&\iff&K^{-1}uKe_0=e_0\\
&\iff&K^{-1}uK=diag(1,w)
\end{eqnarray*}

Thus we have an isomorphism given by $w_{ij}\to(K^{-1}uK)_{ij}$, as desired.
\end{proof}

With this understood, we would like to discuss now certain versions of the bistochastic groups and quantum groups $B_N,B_N^+$, which are quite interesting, for various reasons. To be more precise, it is possible to talk about bistochastic groups and quantum groups $B_N^I,B_N^{I+}$ defined by using other vectors $\xi=\xi_I$, the result being as follows:

\begin{theorem}
We have a closed subgroup $B_N^{I+}\subset O_N^+$, defined via the formula
$$C(B_N^{I+})=C(O_N^+)\Big/\left<u\xi_I=\xi_I\right>$$
with the vector $\xi_I$ being as follows, with $I\subset\{1,\ldots,N\}$ being a subset:
$$\xi_I=\frac{1}{\sqrt{|I|}}(\delta_{i\in I})_i$$
Moreover, we can talk about classical groups $B_N^I$ too, in a similar way.
\end{theorem}

\begin{proof}
We must check the Hopf algebra axioms, and the proof goes as follows:

\medskip

(1) Let us set $U_{ij}=\sum_ku_{ik}\otimes u_{kj}$. We have then the following computation:
\begin{eqnarray*}
(U\xi_I)_i
&=&\frac{1}{\sqrt{|I|}}\sum_{j\in I}U_{ij}\\
&=&\frac{1}{\sqrt{|I|}}\sum_{j\in I}\sum_ku_{ik}\otimes u_{kj}\\
&=&\sum_ku_{ik}\otimes(u\xi_I)_k\\
&=&\sum_ku_{ik}\otimes(\xi_I)_k\\
&=&\frac{1}{\sqrt{|I|}}\sum_{k\in I}u_{ik}\otimes1\\
&=&(u\xi_I)_i\otimes1\\
&=&(\xi_I)_i\otimes1
\end{eqnarray*}

Thus we can define indeed a comultiplication map, by $\Delta(u_{ij})=U_{ij}$.

\medskip

(2) In order to construct the counit map, $\varepsilon(u_{ij})=\delta_{ij}$, we must prove that the identity matrix $1=(\delta_{ij})_{ij}$ satisfies $1\xi_I=\xi_I$. But this is clear.

\medskip

(3) In order to construct the antipode, $S(u_{ij})=(u^{-1})_{ij}$, we must prove that the inverse matrix $u^{-1}=(u^{-1})_{ij}$ satisfies $u^{-1}\xi_I=\xi_I$. But this is clear from $u\xi_I=\xi_I$.

\medskip

(4) Finally, we can talk about classical groups $B_N^I$ in a similar way, with these appearing as subgroups of the quantum groups $B_N^{I+}$, in the obvious way. 
\end{proof}

We should mention that the above quantum groups are of key importance in connection with the notion of ``affine homogeneous space'', which is something generalizing the usual spheres, and requires an index set $I\subset\{1,\ldots,N\}$, as above. We will be back to this later in this book, in Part IV, with a brief introduction to the affine homogeneous spaces.

\section*{5d. Lower liberations}

We would like to discuss now some ``lower liberations'' of the group $GL_N$ and its subgroups of type $O_N$, and why not of some other subgroups too, which do not have free versions. With this meaning looking for intermediate quantum groups $G\subset G^\times\subset G^+$, or just looking for partial liberations $G\subset G^\times$, when $G^+$ does not really exist. 

\bigskip

All this is something a bit philosophical, beyond Hopf algebras, just in relation with the notion of commutativity. And, we are led right away to the following question:

\begin{question}
What is simple and natural as relation between commutativity
$$ab=ba$$
and no relation at all, or freeness, $\emptyset$?
\end{question}

In order to answer this question, we can use a diagrammatic approach to it. Obviously, the notion of commutativity, $ab=ba$, comes from the following diagram:
$$
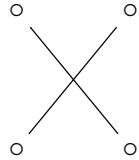
\xymatrix@R=40pt@C=30pt{
\circ\ar@{-}[dr]&\circ\ar@{-}[dl]\\
\circ&\circ}$$

With this idea in mind, we are led right way into the following diagram:
$$
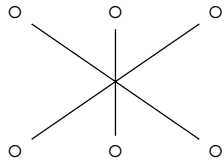
\xymatrix@R=40pt@C=25pt{
\circ\ar@{-}[drr]&\circ\ar@{-}[d]&\circ\ar@{-}[dll]\\
\circ&\circ&\circ}$$

But, by following this path, why not using as well the following diagram:
$$
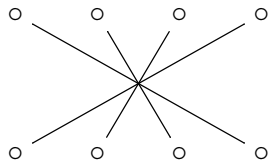
\xymatrix@R=40pt@C=18pt{
\circ\ar@{-}[drrr]&\circ\ar@{-}[dr]&\circ\ar@{-}[dl]&\circ\ar@{-}[dlll]\\
\circ&\circ&\circ&\circ}$$

And so on. A clarification is certainly needed here, and we have:

\begin{claim}
Under suitable combinatorial assumptions, only half-commutativity,
$$abc=cba$$
has its place between commutativity, $ab=ba$, and freeness, $\emptyset$.
\end{claim}

In order to discuss this claim, which is something quite subtle, the idea will be that of regarding the above diagrams as being set-theoretic pairings, and then talk about categories of pairings, rather of pairings themselves. Indeed, let us formulate:

\begin{definition}
A category of pairings is a collection of sets of pairings,
$$D=\bigsqcup_{k,l}D(k,l)$$
with $D(k,l)\subset P_2(k,l)$, having the following properties:
\begin{enumerate}
\item Stability under the horizontal concatenation, $(\pi,\sigma)\to[\pi\sigma]$.

\item Stability under vertical concatenation $(\pi,\sigma)\to[^\sigma_\pi]$, with matching middle symbols.

\item Stability under the upside-down turning $*$, with switching of colors, $\circ\leftrightarrow\bullet$.

\item Each set $P(k,k)$ contains the identity partition $||\ldots||$.

\item The sets $P(\emptyset,\circ\bullet)$ and $P(\emptyset,\bullet\circ)$ both contain the semicircle $\cap$.
\end{enumerate}
\end{definition} 

This might look a bit complicated, at the first glance, but if you are a bit familiar with Tannakian duality, this definition will certainly bring you joy. If not, you can simply take the above definition as such, sort of an abstract notion introduced for our theory below to work, and we will be back to this in chapter 11, with more details about all this.

\bigskip

The point now is that, with the above definition in hand, we can formulate the following statement, making our previous Claim 5.25 into something very precise:

\begin{theorem}
There are only $3$ categories of pairings, namely
$$NC_2\subset P_2^*\subset P_2$$
with $P_2$, all pairings, corresponding to commutativity, $NC_2$, noncrossing pairings, corresponding to freeness, and $P_2^*=<\slash\hskip-2.0mm\backslash\hskip-1.7mm|\hskip0.5mm>$ corresponding to half-commutativity.
\end{theorem}

\begin{proof}
As already mentioned, this is something quite technical and precise, and with the last part being of course a bit informal, making the link with Claim 5.25.

\medskip

(1) In order to prove the uniqueness result, consider a crossing pairing, $\pi\in P_2-NC_2$, having $s\geq 4$ strings. Our claim is that the following happen:

\medskip

-- If $\pi\in P_2-P_2^*$, there exists a semicircle capping $\pi'\in P_2-P_2^*$.

\smallskip

-- If $\pi\in P_2^*-NC_2$, there exists a semicircle capping $\pi'\in P_2^*-NC_2$.

\medskip

Indeed, both these assertions can be easily proved, by drawing pictures.

\medskip

(2) Consider now a partition $\pi\in P_2(k,l)-NC_2(k,l)$. Our claim is that:

\medskip

-- If $\pi\in P_2(k, l)-P_2^*(k,l)$ then $<\pi>=P_2$.

\smallskip

-- If $\pi\in P_2^*(k,l)-NC_2(k,l)$ then $<\pi>=P_2^*$.

\medskip

But this can be indeed proved by recurrence on the number of strings, $s=(k+l)/2$, by using (1), which provides us with a descent procedure $s\to s-1$, at any $s\geq4$.

\medskip

(3) Finally, given a category of pairings $NC_2\subset D\subset P_2$, we have three cases:

\medskip

-- If $D\not\subset P_2^*$, we obtain $D=P_2$.

\smallskip

-- If $D\subset P_2,D\not\subset NC_2$, we obtain $D=P_2^*$.

\smallskip

-- If $D\subset NC_2$, we obtain $D=NC_2$.

\medskip

Thus, we have proved the uniquess result, as desired.
\end{proof}

We will be back to all this later in this book, the idea being that what we just proved is that, under suitable combinatorial assumptions, of Brauer theory type, there are exactly 3 orthogonal quantum groups, namely $O_N\subset O_N^*\subset O_N^+$. More on this later. 

\bigskip

With the above discussed, let us go back now to quantum groups, and discuss the half-liberation operation, obtained by replacing the commutation relations $ab=ba$ between standard coordinates with half-commutation relations, $abc=cba$. We have here:

\begin{theorem}
We have intermediate quantum groups as follows,
$$F(G^*)=F(G^+)\Big/\left<abc=cba\Big|\forall a,b,c\in\{u_{ij},\bar{u}_{ij}\}\right>$$
which have the following properties:
\begin{enumerate}
\item The projective versions $PF(G^*)$ are commutative.

\item The diagonal algebras $F(G^*)^\delta$ are not commutative, for $G=GL_N,O_N$.

\item Thus, the algebras $F(G^*)$ are not commutative, for $G=GL_N,O_N$.

\item However, collapsing $F(G^*)=F(G)$ can happen, for instance for $G=B_N$.
\end{enumerate}
\end{theorem}

\begin{proof}
There are many things going on here, and with this being just the tip of the iceberg, on what can be said, with the idea with all this being as follows:

\medskip

(1) To start with, the fact that $G^*$ as constructed above is indeed a quantum group, lying as a proper intermediate subgroup $G\subset G^*\subset G^+$, can be checked via a routine computation, but the best is to view this via representation theory. Indeed, intuitively, the half-commutation relations $abc=cba$ come from the half-classical crossing:
$$\slash\hskip-2.0mm\backslash\hskip-1.7mm|\hskip0.5mm\in P_2(3,3)$$

Thus, we are led to the conclusion in the statement, by imposing the following conditions, with the convention that $\slash\hskip-2.0mm\backslash\hskip-1.7mm|\hskip0.5mm$ stands for all 8 possible colorings of $\slash\hskip-2.0mm\backslash\hskip-1.7mm|\hskip0.5mm$:
$$\slash\hskip-2.0mm\backslash\hskip-1.7mm|\hskip0.5mm\in End(u^\alpha\otimes u^\beta\otimes u^\gamma)$$

And we will be back to this, correspondence between intertwiners and partitions, later in this book, in chapter 11, when talking Tannakian duality and Brauer theorems.

\medskip

(2) Regarding now the claim that the projective version $PF(G^*)$ is commutative, this is something immediate, coming from the following computation:
\begin{eqnarray*}
ab\cdot cd
&=&abc\cdot d\\
&=&cba\cdot d\\
&=&c\cdot bad\\
&=&c\cdot dab\\
&=&cd\cdot ab
\end{eqnarray*}

(3) Next, in the cases $G=GL_N,O_N$ the diagonal algebras $F(G^*)^\delta$ are respectively equal to the group algebras of $H=\mathbb Z^{\circ N},\mathbb Z_2^{\circ N}$, with the symbol $\circ$ standing for the half-classical product of groups, lying between the classical product $\times$ and the free product $*$, by imposing the commutation relations $ghk=khg$. Now since these latter groups $H$ are not abelian, nor are the corresponding group algebras commutative, as claimed.

\medskip

(4) Still in the cases $G=GL_N,O_N$, since the diagonal algebras $F(G^*)^\delta$ are not commutative, nor are the algebras $F(G^*)$ themselves commutative, as claimed.

\medskip

(5) Finally, in the case $G=B_N$, as a first observation, the above argument involving diagonal algebras will not work, because these diagonal algebras collapse to $F$, which half-liberate and liberate into $F$ too. In fact, as said in the statement, the half-liberation procedure fails for the bishochastic groups, and this because by summing the matrix entries over a row, the relations $abc=cba$ get converted into the relations $ab=ba$.
\end{proof}

Many other things can be said about half-liberation, following Bichon-Dubois-Violette \cite{bd1} and related work, with a key observation here being the fact that the $2\times2$ antidiagonal scalar matrices half-commute. Indeed, we have the following computation:
$$\begin{pmatrix}0&a\\b&0\end{pmatrix}
\begin{pmatrix}0&x\\y&0\end{pmatrix}
\begin{pmatrix}0&c\\d&0\end{pmatrix}
=\begin{pmatrix}ay&0\\0&bx\end{pmatrix}
\begin{pmatrix}0&c\\d&0\end{pmatrix}
=\begin{pmatrix}0&ayc\\bxd&0\end{pmatrix}$$

On the other hand, we have as well the following computation:
$$\begin{pmatrix}0&c\\d&0\end{pmatrix}
\begin{pmatrix}0&x\\y&0\end{pmatrix}
\begin{pmatrix}0&a\\b&0\end{pmatrix}
=\begin{pmatrix}cy&0\\0&dx\end{pmatrix}
\begin{pmatrix}0&a\\b&0\end{pmatrix}
=\begin{pmatrix}0&cya\\dxb&0\end{pmatrix}$$

Thus, we have indeed half-commutation, and based on this, it is possible to work out explicit $2\times2$ matrix models for our quantum groups. More on this, later in this book.

\bigskip

As a last topic, back to the bistochastic groups, it is possible to construct intermediate objects for $B_N\subset B_N^+$. Let us recall from Theorem 5.22 that we have an isomorphism as follows, whenever $K\in GL_N$ satisfies $Ke_0=\frac{1}{\sqrt{N}}\xi$, where $\xi$ is the all-one vector:
$$C(O_{N-1}^+)\to C(B_N^+)\quad,\quad w_{ij}\to(K^{-1}uK)_{ij}$$

Here, and in what follows, we use indices $i,j=0,1,\ldots,N-1$ for the $N\times N$ quantum groups, and indices $i,j=1,\ldots,N-1$ for their $(N-1)\times(N-1)$ subgroups. But with this, we can construct intermediate objects for $B_N\subset B_N^+$, as follows:

\begin{theorem}
Assuming that $K\in GL_N$ satisfies $Ke_0=\frac{1}{\sqrt{N}}\xi$, we have inclusions
$$B_N\subset B_K^\circ\subset B_N^+$$
obtained by taking the image of the standard inclusions 
$$O_{N-1}\subset O_{N-1}^*\subset O_{N-1}^+$$
via the above isomorphism $O_{N-1}^+\simeq B_N^+$ induced by $K$.
\end{theorem}

\begin{proof}
The fact that we have inclusions as in the statement follows from the above isomorphism, which produces a diagram as follows:
$$\xymatrix@R=15mm@C=20mm{
O_{N-1}\ar[r]&O_{N-1}^*\ar[r]&O_{N-1}^+\\
B_N\ar[r]\ar@{=}[u]&B_K^\circ\ar[r]\ar@{=}[u]&B_N^+\ar@{=}[u]
}$$

To be more precise, the quantum group $B_K^\circ$ from the bottom is by definition the image of the quantum group $O_{N-1}^*$ from the top, and this gives the result.
\end{proof}

Many other things can be said, as a continuation of the above. We will back to this, and with some further examples of bistochastic quantum groups too, later on, when systematically discussing the half-liberation operation, over $F=\mathbb C$.

\section*{5e. Exercises}

We had a lot of interesting algebra in this chapter, and as exercises, we have:

\begin{exercise}
Clarify the isomorphisms between the quantum groups $GL_Q^+$.
\end{exercise}

\begin{exercise}
Clarify the deformation of the quantum groups $O_N^+$.
\end{exercise}

\begin{exercise}
Learn more about the bistochastic matrices.
\end{exercise}

\begin{exercise}
Clarify the details, in relation with the half-liberation operation.
\end{exercise}

\begin{exercise}
Try constructing more general half-liberation operations.
\end{exercise}

\begin{exercise}
Further study the bistochastic quantum groups $B_K^\circ$.
\end{exercise}

As bonus exercise, study the representations of our various quantum groups.

\chapter{Symmetric groups}

\section*{6a. Symmetric groups}

Welcome to quantum permutations. Our next goal, coming as a continuation of the general liberation material from chapter 5, will be to talk about quantum permutations, and quantum reflections too, over arbitrary fields $F$. To be more precise, we would first like to talk about quantum permutation groups $S_N^+$, and then about more general quantum reflection groups, appearing as certain products, $H_N^s=\mathbb Z_s\wr_*S_N^+$ with $s\in\mathbb N\cup\{\infty\}$.

\bigskip

As before with the rotation groups, the story with quantum permutations involves some key work of Wang from the 1990s, in the case $F=\mathbb C$, and then some subsequent work of Bichon, from the 2000s, in the case where $F$ is arbitrary. We will be following here the paper of Bichon \cite{bi4}, which includes adaptations of most of Wang's results.

\bigskip

Let us start with the following standard result, regarding the group $S_N$ itself:

\index{magic matrix}
\index{magic unitary}
\index{symmetric group}

\begin{theorem}
Consider the symmetric group $S_N$, viewed as permutation group of the $N$ coordinate axes of $F^N$. The coordinate functions on $S_N\subset O_N$ are given by
$$u_{ij}=\chi\left(\sigma\in G\Big|\sigma(j)=i\right)$$
and the matrix $u=(u_{ij})$ that these functions form is magic, in the sense that:
\begin{enumerate}
\item Its entries are projections, $p^2=p$,

\item Which satisfy $pq=0$, on each row and column,

\item And which sum up to $1$ on each row and column.
\end{enumerate}
\end{theorem}

\begin{proof}
The action of $S_N$ on the standard basis $e_1,\ldots,e_N\in F^N$ being given by $\sigma:e_j\to e_{\sigma(j)}$, this gives the formula of $u_{ij}$ in the statement. As for the fact that the matrix $u=(u_{ij})$ that these functions form is magic, this is clear.
\end{proof}

As a comment here, in the case $F=\mathbb C$ the above projections $u_{ij}$ are orthogonal projections, $p^2=p^*=p$, and by using this, the condition (2) automatically follows from (3), via some standard positivity tricks, that we will leave as an exercise. In the general case, however, as explained in \cite{bi4}, the good magic conditions are (1,2,3) above.

\bigskip

Now back to the general case, where our ground field $F$ is arbitrary, with a bit more effort, we obtain the following nice algebraic characterization of $S_N$:

\index{Gelfand theorem}
\index{commutative algebra}

\begin{theorem}
The algebra of functions on $S_N$ has the following presentation,
$$F(S_N)=F_{comm}\left((u_{ij})_{i,j=1,\ldots,N}\Big|u={\rm magic}\right)$$
and the multiplication, unit and inversion map of $S_N$ appear from the maps
$$\Delta(u_{ij})=\sum_ku_{ik}\otimes u_{kj}\quad,\quad 
\varepsilon(u_{ij})=\delta_{ij}\quad,\quad 
S(u_{ij})=u_{ji}$$
defined at the algebraic level, of functions on $S_N$, by transposing.
\end{theorem}

\begin{proof}
This follows indeed from Theorem 6.1, via some basic algebra, that we will leave here as an exercise, and with the comment that there are no topological issues of any kind involved, both the algebras in the statement being finite dimensional. We will actually do this basic algebra in a moment, when talking about liberations of $S_N$.
\end{proof}

Following now Wang and Bichon, we can liberate $S_N$, as follows:

\index{quantum permutation group}
\index{magic unitary}
\index{free symmetric group}

\begin{theorem}
The following universal algebra, with magic meaning as usual formed by projections, $p^2=p$, satisfying $pq=0$ and $\sum p=1$ on each row and column,
$$F(S_N^+)=F\left((u_{ij})_{i,j=1,\ldots,N}\Big|u={\rm magic}\right)$$
is a Hopf algebra, with comultiplication, counit and antipode given by:
$$\Delta(u_{ij})=\sum_ku_{ik}\otimes u_{kj}\quad,\quad 
\varepsilon(u_{ij})=\delta_{ij}\quad,\quad 
S(u_{ij})=u_{ji}$$
We call the quantum group $S_N^+$ quantum permutation group.
\end{theorem}

\begin{proof}
We must construct morphisms of algebras $\Delta,\varepsilon,S$ given by the formulae in the statement. Consider the following matrices:
$$u^\Delta_{ij}=\sum_ku_{ik}\otimes u_{kj}\quad,\quad 
u^\varepsilon_{ij}=\delta_{ij}\quad,\quad 
u^S_{ij}=u_{ji}$$

Our claim is that, since $u$ is magic, so are these three matrices:

\medskip

(1) Indeed, regarding the matrix $u^\Delta$, its entries are projections, as shown by:
\begin{eqnarray*}
(u_{ij}^\Delta)^2
&=&\sum_{kl}u_{ik}u_{il}\otimes u_{kj}u_{lj}\\
&=&\sum_{kl}\delta_{kl}u_{ik}\otimes\delta_{kl}u_{kj}\\
&=&\sum_ku_{ik}\otimes u_{kj}\\
&=&u_{ij}^\Delta
\end{eqnarray*}

(2) Next, the condition $pq=0$ is satisfied on the rows of $u^\Delta$, as shown by:
\begin{eqnarray*}
u_{ij}u_{ik}
&=&\sum_{lm}u_{il}u_{im}\otimes u_{lj}u_{mk}\\
&=&\sum_{lm}\delta_{lm}u_{il}\otimes u_{lj}u_{mk}\\
&=&\sum_lu_{il}\otimes u_{lj}u_{lk}\\
&=&\sum_lu_{il}\otimes\delta_{jk}u_{lj}\\
&=&\delta_{jk}u_{ij}
\end{eqnarray*}

(3) The condition $pq=0$ is satisfied on the columns of $u^\Delta$ too, as shown by:
\begin{eqnarray*}
u_{ij}u_{kj}
&=&\sum_{lm}u_{il}u_{km}\otimes u_{lj}u_{mj}\\
&=&\sum_{lm}u_{il}u_{km}\otimes\delta_{lm}u_{lj}\\
&=&\sum_lu_{il}u_{kl}\otimes u_{lj}\\
&=&\sum_l\delta_{ik}u_{il}\otimes u_{lj}\\
&=&\delta_{ik}u_{ij}
\end{eqnarray*}

(4) The row sums for the matrix $u^\Delta$ can be computed as follows:
\begin{eqnarray*}
\sum_ju_{ij}^\Delta
&=&\sum_{jk}u_{ik}\otimes u_{kj}\\
&=&\sum_ku_{ik}\otimes 1\\
&=&1
\end{eqnarray*}

(5) As for the column sums, the computation here is similar, as follows:
\begin{eqnarray*}
\sum_iu_{ij}^\Delta
&=&\sum_{ik}u_{ik}\otimes u_{kj}\\
&=&\sum_k1\otimes u_{kj}\\
&=&1
\end{eqnarray*}

(6) We conclude from these computations that the matrix $u^\Delta$ is magic, as claimed. As for $u^\varepsilon,u^S$, these matrices are magic too, and this for obvious reasons. 

\medskip

(7) Thus, all our three matrices $u^\Delta,u^\varepsilon,u^S$ are magic, so we can define $\Delta,\varepsilon,S$ by the formulae in the statement, by using the universality property of $F(S_N^+)$.
\end{proof}

Very nice, but our first task now is to make sure that Theorem 6.3 produces indeed a new quantum group, which does not collapse to $S_N$. Following Wang, we have:

\index{liberation}

\begin{theorem}
We have an embedding $S_N\subset S_N^+$, given at the algebra level by: 
$$u_{ij}\to\chi\left(\sigma\in S_N\Big|\sigma(j)=i\right)$$
This is an isomorphism at $N\leq3$, but not at $N\geq4$, where $S_N^+$ is not classical, nor finite.
\end{theorem} 

\begin{proof}
The fact that we have indeed an embedding as above follows from Theorem 6.2. Observe that in fact more is true, because Theorems 6.2 and 6.3 give:
$$F(S_N)=F(S_N^+)\Big/\Big<ab=ba\Big>$$

Thus, the inclusion $S_N\subset S_N^+$ is a ``liberation'', in the sense that $S_N$ is the classical version of $S_N^+$. We will often use this basic fact, in what follows. Regarding now the second assertion, we can prove this in four steps, as follows:

\medskip

\underline{Case $N=2$}. The fact that $S_2^+$ is indeed classical, and hence collapses to $S_2$, is trivial, because the $2\times2$ magic matrices are as follows, with $p$ being a projection:
$$U=\begin{pmatrix}p&1-p\\1-p&p\end{pmatrix}$$

Indeed, this shows that the entries of $U$ commute. Thus $F(S_2^+)$ is commutative, and so equals its biggest commutative quotient, which is $F(S_2)$. Thus, $S_2^+=S_2$.

\medskip

\underline{Case $N=3$}. By using the same argument as in the $N=2$ case, and the symmetries of the problem, it is enough to check that $u_{11},u_{22}$ commute. We have:
\begin{eqnarray*}
u_{11}u_{22}
&=&u_{11}u_{22}(u_{11}+u_{12}+u_{13})\\
&=&u_{11}u_{22}u_{11}+u_{11}u_{22}u_{13}\\
&=&u_{11}u_{22}u_{11}+u_{11}(1-u_{21}-u_{23})u_{13}\\
&=&u_{11}u_{22}u_{11}
\end{eqnarray*}

On the other hand, we have as well the following computation:
\begin{eqnarray*}
u_{22}u_{11}
&=&(u_{11}+u_{12}+u_{13})u_{22}u_{11}\\
&=&u_{11}u_{22}u_{11}+u_{13}u_{22}u_{11}\\
&=&u_{11}u_{22}u_{11}+u_{13}(1-u_{21}-u_{23})u_{11}\\
&=&u_{11}u_{22}u_{11}
\end{eqnarray*}

By comparing the above formulae we obtain $u_{11}u_{22}=u_{22}u_{11}$, as desired.

\medskip

\underline{Case $N=4$}. Consider the following matrix, with $p,q$ being projections:
$$U=\begin{pmatrix}
p&1-p&0&0\\
1-p&p&0&0\\
0&0&q&1-q\\
0&0&1-q&q
\end{pmatrix}$$ 

This matrix is magic, and we can choose $p,q$ as for the algebra $<p,q>$ to be noncommutative and infinite dimensional. We conclude that $F(S_4^+)$ is noncommutative and infinite dimensional as well, and so $S_4^+$ is non-classical and infinite, as claimed.

\medskip

\underline{Case $N\geq5$}. Here we can use the standard embedding $S_4^+\subset S_N^+$, obtained at the level of the corresponding magic matrices in the following way:
$$u\to\begin{pmatrix}u&0\\ 0&1_{N-4}\end{pmatrix}$$

Indeed, with this in hand, the fact that $S_4^+$ is a non-classical, infinite quantum group implies that $S_N^+$ with $N\geq5$ has these two properties as well.
\end{proof}

As a first observation, as a matter of doublechecking our findings, we are not wrong with our formalism, because as explained once again in \cite{bi4}, \cite{wan}, we have as well:

\index{quantum permutation}
\index{counting measure}
\index{coaction}
\index{standard coaction}

\begin{theorem}
The quantum permutation group $S_N^+$ acts on the set $X=\{1,\ldots,N\}$, the corresponding coaction map $\alpha:F(X)\to F(X)\otimes F(S_N^+)$ being given by:
$$\alpha(e_i)=\sum_je_j\otimes u_{ji}$$
In fact, $S_N^+$ is the biggest quantum group acting on $X$, by leaving the counting measure invariant, in the sense that $(tr\otimes id)\alpha=tr(.)1$, where $tr(e_i)=\frac{1}{N},\forall i$.
\end{theorem}

\begin{proof}
Our claim, which will prove the result, is that given a Hopf algebra $A$, the following formula defines a morphism of algebras, which is a coaction map, leaving the trace invariant, precisely when $u=(u_{ij})$ is a magic corepresentation of $A$: 
$$\alpha(e_i)=\sum_je_j\otimes u_{ji}$$

(1) To start with, we know from chapter 3 that the following happen:
$$(\alpha\otimes id)\alpha=(id\otimes\Delta)\alpha\iff\Delta(u_{ij})=\sum_ku_{ik}\otimes u_{kj}$$
$$(id\otimes\varepsilon)\alpha=id\iff (id\otimes\varepsilon)u=1$$

In other words, $\alpha$ is a coaction precisely when $u$ is a corepresentation.

\medskip

(2) Next, let us determine when the above map $\alpha$ is multiplicative. We have:
$$\alpha(e_i)\alpha(e_k)
=\sum_{jl}e_je_l\otimes u_{ji}u_{lk}
=\sum_je_j\otimes u_{ji}u_{jk}$$
$$\alpha(e_ie_k)
=\delta_{ik}\alpha(e_i)
=\delta_{ik}\sum_je_j\otimes u_{ji}$$

Thus the multiplicativity of $\alpha$ is equivalent to the following conditions:
$$u_{ji}u_{jk}=\delta_{ik}u_{ji}\quad,\quad\forall i,j,k$$

In other words, the entries $u_{ij}$ must be projections, satisfying $pq=0$ on rows.

\medskip

(3) Regarding now the unitality property of $\alpha$, we have the following formulae:
$$\alpha(1)=\alpha\left(\sum_ie_i\right)=\sum_{ij}e_j\otimes u_{ji}$$
$$1=\sum_je_j$$

We conclude that our map $\alpha$ is unital when the following condition is satisfied:
$$\sum_iu_{ji}=1\quad,\quad\forall j$$

(4) Time for a partial conclusion, to all this. We just found that, in order for our linear map $\alpha(e_i)=\sum_je_j\otimes u_{ji}$ to be a morphism of algebras, the elements $u_{ij}$ must be projections, satisfying $pq=0$ on each row, and summing up to 1 on each row. 

\medskip

(5) Regarding now the preservation of the trace condition, observe that we have:
$$(tr\otimes id)\alpha(e_i)=\frac{1}{N}\sum_ju_{ji}$$

Thus the trace is preserved precisely when the following condition is satisfied:
$$\sum_ju_{ij}=1\quad,\quad\forall i$$

(6) Time for an updated conclusion. We found that, in order for $\alpha(e_i)=\sum_je_j\otimes u_{ji}$ to be a morphism of algebras, preserving the trace, the elements $u_{ij}$ must be projections, satisfying $pq=0$ on each row, and summing up to 1 on each row and column. 

\medskip

(7) And with this, we are almost there, but we still have to prove, based on what we have, that the elements $u_{ij}$ satisfy the condition $pq=0$, on any column. Now for this purpose, observe that, by using the conditions found in (6), we have:
$$uu^t=1$$

Now remember from (1) that the matrix $u=(u_{ij})$ must be a corepresentation. Thus $(id\otimes S)u=u^t$, and so $S(u_{ij})=u_{ji}$, and by using the formulae in (2) we get:
\begin{eqnarray*}
u_{ji}u_{jk}=\delta_{ik}u_{ji}
&\implies&S(u_{ji}u_{jk})=\delta_{ik}S(u_{ji})\\
&\implies&u_{ij}u_{kj}=\delta_{ik}u_{ij}
\end{eqnarray*}

Thus the condition $pq=0$ on the columns is indeed satisfied, as desired.

\medskip

(8) Summarizing, we have proved that $\alpha(e_i)=\sum_je_j\otimes u_{ji}$ is a coaction, preserving the trace, precisely when $u=(u_{ij})$ is a magic corepresentation, and this gives the result.
\end{proof}

And with this, end of our preliminary discussion regarding the quantum permutations, these definitely exist, in the real life, we even have multiple proofs for their existence, and so, welcome to them, we will have to live with this, for the rest of our lives.

\bigskip

Retrospectively thinking, the origin of these phenomena was in chapter 5, where we opened the Pandora box of liberations. At that time, the first beasts emerging from that box were not the scary, but now that we got into permutations, we obviously face some truly wild beasts. But do not worry, in the long run we will manage to have some understanding of the free world, and in the end, we will even get to love it.

\section*{6b. Orbits, orbitals} 

Getting now to quantum subgroups $G\subset S_N^+$, there are many of them. In order to have some theory going on for them, again following Bichon \cite{bi4}, let us formulate:

\index{orbits}

\begin{theorem}
Given a quantum subgroup $G\subset S_N^+$, with standard coordinates denoted $u_{ij}\in C(G)$, the following defines an equivalence relation on $X=\{1,\ldots,N\}$,
$$i\sim j\iff u_{ij}\neq0$$
that we call orbit decomposition associated to the action $G\curvearrowright X$. In the classical case, $G\subset S_N$, this is the usual orbit equivalence.
\end{theorem}

\begin{proof}
We first check the fact that we have indeed an equivalence relation. The reflexivity axiom $i\sim i$ follows by using the counit, as follows:
$$\varepsilon(u_{ii})=1
\implies u_{ii}\neq0$$

The symmetry axiom $i\sim j\implies j\sim i$ follows by using the antipode:
$$S(u_{ji})=u_{ij}\implies[u_{ij}\neq0\implies u_{ji}\neq0]$$

As for the transitivity axiom $i\sim k,k\sim j\implies i\sim j$, this is something more tricky. In order to prove this, we can use the following formula:
\begin{eqnarray*}
\Delta(u_{ij})(1\otimes u_{kj})
&=&\left(\sum_lu_{il}\otimes u_{lj}\right)(1\otimes u_{kj})\\
&=&\sum_lu_{il}\otimes u_{lj}u_{kj}\\
&=&u_{ik}\otimes u_{kj}
\end{eqnarray*}

We deduce from this that we have, as required by the transitivity condition:
$$u_{ik}\neq0,u_{kj}\neq0\ \implies\ u_{ij}\neq0$$

Finally, in the classical case, where $G\subset S_N$, the standard coordinates are:
$$u_{ij}=\chi\left(\sigma\in G\Big|\sigma(j)=i\right)$$

Thus $u_{ij}\neq0$ means that $i,j$ must be in the same orbit, as claimed.
\end{proof}

Generally speaking, the theory from the classical case extends well to the quantum group setting, and we have in particular the following result, also from \cite{bi4}:

\index{fixed points}

\begin{theorem}
Given a subgroup $G\subset S_N^+$, the following algebras are equal,
$$Fix(u)=\left\{\xi\in F(X)\Big|u\xi=\xi\right\}$$
$$Fix(\alpha)=\left\{\xi\in F(X)\Big|\alpha(\xi)=\xi\otimes1\right\}$$
$$Fix(\sim)=\left\{\xi\in F(X)\Big|i\sim j\implies \xi_i=\xi_j\right\}$$
where $\sim$ is the orbit equivalence relation constructed in Theorem 6.6.
\end{theorem}

\begin{proof}
The fact that we have $Fix(u)=Fix(\alpha)$ is standard, with this being valid for any corepresentation $u=(u_{ij})$. Indeed, we first have the following computation:
\begin{eqnarray*}
\xi\in Fix(u)
&\iff&u\xi=\xi\\
&\iff&(u\xi)_j=\xi_j,\forall j\\
&\iff&\sum_iu_{ji}\xi_i=\xi_j,\forall j
\end{eqnarray*}

On the other hand, we have as well the following computation:
\begin{eqnarray*}
\xi\in Fix(\alpha)
&\iff&\alpha(\xi)=\xi\otimes1\\
&\iff&\sum_i\alpha(e_i)\xi_i=\xi\otimes1\\
&\iff&\sum_{ij}e_j\otimes u_{ji}\xi_i=\sum_je_j\otimes\xi_j\\
&\iff&\sum_iu_{ji}\xi_i=\xi_j,\forall j
\end{eqnarray*}

Thus we have $Fix(u)=Fix(\alpha)$, as claimed. Regarding now the equality of this algebra with $Fix(\sim)$, observe first that given a vector $\xi\in Fix(\sim)$, we have:
\begin{eqnarray*}
\sum_iu_{ji}\xi_i
&=&\sum_{i\sim j}u_{ji}\xi_i\\
&=&\sum_{i\sim j}u_{ji}\xi_j\\
&=&\sum_iu_{ji}\xi_j\\
&=&\xi_j
\end{eqnarray*}

Thus $\xi\in Fix(u)=Fix(\alpha)$. Finally, for the reverse inclusion, we know from Theorem 6.6 that the magic unitary $u=(u_{ij})$ is block-diagonal, with respect to the orbit decomposition there. But this shows that the algebra $Fix(u)=Fix(\alpha)$ decomposes as well with respect to the orbit decomposition, and this gives the result. See \cite{bi4}.
\end{proof}

As an application of the above, let us discuss now the notion of transitivity. We have here the following result, once again coming from \cite{bi4}:

\begin{theorem}
For a subgroup $G\subset S_N^+$, the following are equivalent:
\begin{enumerate}
\item $G$ is transitive, in the sense that $i\sim j$, for any $i,j$.

\item $Fix(u)=F\xi$, where $\xi$ is the all-one vector.
\end{enumerate}
In the classical case, $G\subset S_N$, this is the usual notion of transitivity.
\end{theorem}

\begin{proof}
We use the standard fact that the fixed point space of a corepresentation coincides with the fixed point space of the associated coaction:
$$Fix(u)=Fix(\alpha)$$

As we explained in Theorem 6.7, the fixed point space of the magic corepresentation $u=(u_{ij})$ has the following interpretation, in terms of orbits:
$$Fix(u)=\left\{\xi\in F(X)\Big|i\sim j\implies \xi(i)=\xi(j)\right\}$$

In particular, the transitivity condition corresponds to $Fix(u)=F\xi$, as stated.
\end{proof}

Following Lupini, Man\v cinska and Roberson, let us discuss now the higher orbitals. Things are quite tricky here, and we have the following result, to start with:

\index{higher orbitals}

\begin{theorem}
For a subgroup $G\subset S_N^+$, with magic corepresentation $u=(u_{ij})$, and $k\in\mathbb N$, the relation 
$$(i_1,\ldots,i_k)\sim(j_1,\ldots,j_k)\iff u_{i_1j_1}\ldots u_{i_kj_k}\neq0$$
is reflexive and symmetric, and is transitive at $k=1,2$. In the classical case, $G\subset S_N$, this relation is transitive at any $k\in\mathbb N$, and is the usual $k$-orbital equivalence.
\end{theorem}

\begin{proof}
This is something quite tricky, the proof being as follows:

\medskip

(1) The reflexivity of $\sim$ follows by using the counit, as follows:
\begin{eqnarray*}
\varepsilon(u_{i_ri_r})=1,\forall r
&\implies&\varepsilon(u_{i_1i_1}\ldots u_{i_ki_k})=1\\
&\implies&u_{i_1i_1}\ldots u_{i_ki_k}\neq0\\
&\implies&(i_1,\ldots,i_k)\sim(i_1,\ldots,i_k)
\end{eqnarray*}

(2) The symmetry follows by applying the antipode, and then the involution:
\begin{eqnarray*}
(i_1,\ldots,i_k)\sim(j_1,\ldots,j_k)
&\implies&u_{i_1j_1}\ldots u_{i_kj_k}\neq0\\
&\implies&u_{j_ki_k}\ldots u_{j_1i_1}\neq0\\
&\implies&u_{j_1i_1}\ldots u_{j_ki_k}\neq0\\
&\implies&(j_1,\ldots,j_k)\sim(i_1,\ldots,i_k)
\end{eqnarray*}

(3) The transitivity at $k=1,2$ is more tricky. Here we need to prove that:
$$u_{i_1j_1}\ldots u_{i_kj_k}\neq0\ ,\ u_{j_1l_1}\ldots u_{j_kl_k}\neq0\implies u_{i_1l_1}\ldots u_{i_kl_k}\neq0$$

In order to do so, we can use the following formula:
$$\Delta(u_{i_1l_1}\ldots u_{i_kl_k})=\sum_{s_1\ldots s_k}u_{i_1s_1}\ldots u_{i_ks_k}\otimes u_{s_1l_1}\ldots u_{s_kl_k}$$

At $k=1$, we already know this. At $k=2$ now, we can use the following trick:
\begin{eqnarray*}
(u_{i_1j_1}\otimes u_{j_1l_1})\Delta(u_{i_1l_1}u_{i_2l_2})(u_{i_2j_2}\otimes u_{j_2l_2})
&=&\sum_{s_1s_2}u_{i_1j_1}u_{i_1s_1}u_{i_2s_2}u_{i_2j_2}\otimes u_{j_1l_1}u_{s_1l_1}u_{s_2l_2}u_{j_2l_2}\\
&=&u_{i_1j_1}u_{i_2j_2}\otimes u_{j_1l_1}u_{j_2l_2}
\end{eqnarray*}

Indeed, we obtain from this the following implication, as desired:
$$u_{i_1j_1}u_{i_2j_2}\neq0,u_{j_1l_1}u_{j_2l_2}\neq0\implies u_{i_1l_1}u_{i_2l_2}\neq0$$

(4) Finally, assume that we are in the classical case, $G\subset S_N$. We have:
$$u_{ij}=\chi\left(\sigma\in G\Big|\sigma(j)=i\right)$$

But this formula shows that we have the following equivalence:
$$u_{i_1j_1}\ldots u_{i_kj_k}\neq0\iff\exists \sigma\in G,\ \sigma(i_1)=j_1,\ldots,\sigma(i_k)=j_k$$

In other words, $(i_1,\ldots,i_k)\sim(j_1,\ldots,j_k)$ happens precisely when $(i_1,\ldots,i_k)$ and $(j_1,\ldots,j_k)$ are in the same $k$-orbital of $G$, and this gives the last assertion.
\end{proof}

The above result raises the question about what exactly happens at $k=3$, in relation with the transitivity property of $\sim$, and the answer here is negative in general. For more on all this, in the case $F=\mathbb C$, we refer to the recent work of McCarthy.

\bigskip

Summarizing, we can talk about orbits and orbitals, but not about higher orbitals. In what regards now the orbitals, we have the following key analogue of Theorem 6.7:

\begin{theorem}
Given a subgroup $G\subset S_N^+$, with magic corepresentation $u=(u_{ij})$, consider the following vector space coaction map, where $X=\{1,\ldots,N\}$:
$$\alpha:F(X\times X)\to F(X\times X)\otimes F(G)\quad,\quad e_{ik}\to\sum_{jl}e_{jl}\otimes u_{ji}u_{lk}$$
The following three algebras are then isomorphic,
$$End(u)=\left\{d\in M_N(F)\Big|du=ud\right\}$$
$$Fix(\alpha)=\left\{\xi\in F(X\times X)\Big|\alpha(\xi)=\xi\otimes1\right\}$$
$$Fix(\sim)=\left\{\xi\in F(X\times X)\Big|(i,k)\sim(j,l)\implies \xi_{ik}=\xi_{jl}\right\}$$
where $\sim$ is the orbital equivalence relation coming from Theorem 6.9.
\end{theorem}

\begin{proof}
This follows by doing some computations which are quite similar to those from the proof of Theorem 6.7, and we refer here to the literature, for the details.
\end{proof}

Finally, as a theoretical application of the theory of orbitals, as developed above, let us discuss now the notion of double transitivity. We have here:

\begin{definition}
Let $G\subset S_N^+$ be a subgroup, with magic corepresentation $u=(u_{ij})$, and consider as before the equivalence relation on $\{1,\ldots,N\}^2$ given by:
$$(i,k)\sim (j,l)\iff u_{ij}u_{kl}\neq0$$ 
\begin{enumerate}
\item The equivalence classes under $\sim$ are called orbitals of $G$.

\item $G$ is called doubly transitive when the action has two orbitals. 
\end{enumerate}
In other words, we call $G\subset S_N^+$ doubly transitive when $u_{ij}u_{kl}\neq0$, for any $i\neq k,j\neq l$.
\end{definition}

To be more precise, it is clear from definitions that the diagonal $D\subset\{1,\ldots,N\}^2$ is an orbital, and that its complement $D^c$ must be a union of orbitals. With this remark in hand, the meaning of (2) is that the orbitals must be $D,D^c$. 

\bigskip

We will be back to this, orbitals, double transitivity and applications, later in this chapter. Moving ahead now, still following Bichon \cite{bi4}, we have:

\begin{theorem}
Given a quotient group $\mathbb Z_{N_1}*\ldots*\mathbb Z_{N_k}\to\Gamma$, we have an embedding $\widehat{\Gamma}\subset S_N^+$, with $N=N_1+\ldots+N_k$, having the following properties:
\begin{enumerate}
\item This embedding appears by diagonally joining the embeddings $\widehat{\mathbb Z_{N_k}}\subset S_{N_k}^+$, and the corresponding magic matrix has blocks of sizes $N_1,\ldots,N_k$.

\item The equivalence relation on $X=\{1,\ldots,N\}$ coming from the orbits of the action $\widehat{\Gamma}\curvearrowright X$ appears by refining the partition $N=N_1+\ldots+N_k$.
\end{enumerate}
\end{theorem}

\begin{proof}
This is something elementary, the idea being as follows:

\medskip

(1) Given a quotient group $\mathbb Z_{N_1}*\ldots*\mathbb Z_{N_k}\to\Gamma$, we have indeed an embedding as follows, with $N=N_1+\ldots+N_k$, with all inclusions being the standard ones:
\begin{eqnarray*}
\widehat{\Gamma}
&\subset&\widehat{\mathbb Z_{N_1}*\ldots*\mathbb Z_{N_k}}
=\widehat{\mathbb Z_{N_1}}\,\hat{*}\,\ldots\,\hat{*}\,\widehat{\mathbb Z_{N_k}}\\
&\simeq&\mathbb Z_{N_1}\,\hat{*}\,\ldots\,\hat{*}\,\mathbb Z_{N_k}
\subset S_{N_1}\,\hat{*}\,\ldots\,\hat{*}\,S_{N_k}\\
&\subset&S_{N_1}^+\,\hat{*}\,\ldots\,\hat{*}\,S_{N_k}^+
\subset S_N^+
\end{eqnarray*}

(2) Regarding now the second assertion in the statement, this is clear from the fact that $u$ is block-diagonal, with blocks corresponding to the partition $N=N_1+\ldots+N_k$. For more on all this, and related topics, we refer as before to the paper of Bichon \cite{bi4}.
\end{proof}

Many other things can be said, as a continuation of the above. We will be back to this later in this book, with a more detailed study in the case $F=\mathbb C$.

\section*{6c. Graph symmetries} 

Switching topics now, we can construct many interesting quantum permutation groups as quantum symmetry groups of graphs, in the following way:

\index{graph symmetry}
\index{quantum symmetry}
\index{quantum automorphism group}

\begin{theorem}
Given a finite graph $X$, with adjacency matrix $d\in M_N(0,1)$, the following construction produces a quantum permutation group, 
$$F(G^+(X))=F(S_N^+)\Big/\Big<du=ud\Big>$$
whose classical version $G(X)$ is the usual  automorphism group of $X$.
\end{theorem}

\begin{proof}
The fact that we have indeed a quantum group comes from the fact that $du=ud$ reformulates as $d\in End(u)$, which makes it clear that $\Delta,\varepsilon,S$ factorize. Regarding now the second assertion, we must establish here the following equality:
$$F(G(X))=F(S_N)\Big/\Big<du=ud\Big>$$

For this purpose, recall that we have $u_{ij}(\sigma)=\delta_{\sigma(j)i}$. We therefore obtain:
$$(du)_{ij}(\sigma)
=\sum_kd_{ik}u_{kj}(\sigma)
=\sum_kd_{ik}\delta_{\sigma(j)k}
=d_{i\sigma(j)}$$

On the other hand, we have as well the following formula:
$$(ud)_{ij}(\sigma)
=\sum_ku_{ik}(\sigma)d_{kj}
=\sum_k\delta_{\sigma(k)i}d_{kj}
=d_{\sigma^{-1}(i)j}$$

Thus $du=ud$ reformulates as $d_{ij}=d_{\sigma(i)\sigma(j)}$, which gives the result.
\end{proof}

Let us work out some examples. With the convention that $\hat{*}$ is the dual free product, obtained by diagonally concatenating the magic unitaries, we have:

\begin{theorem}
The construction $X\to G^+(X)$ has the following properties:
\begin{enumerate}
\item For the $N$-point graph, having no edges at all, we obtain $S_N^+$.

\item For the $N$-simplex, having edges everywhere, we obtain as well $S_N^+$.

\item We have $G^+(X)=G^+(X^c)$, where $X^c$ is the complementary graph.

\item For a disconnected union, we have $G^+(X)\,\hat{*}\,G^+(Y)\subset G^+(X\sqcup Y)$.

\item For the square, we obtain a non-classical, proper subgroup of $S_4^+$.
\end{enumerate}
\end{theorem}

\begin{proof}
All these results are elementary, the proofs being as follows:

\medskip

(1) This follows from definitions, because here we have $d=0$.

\medskip

(2) Here $d=\mathbb I-1$, where $\mathbb I$ is the all-one matrix, and the magic condition gives $u\mathbb I=\mathbb Iu=N\mathbb I$. We conclude that $du=ud$ is automatic, and so $G^+(X)=S_N^+$.

\medskip

(3) The adjacency matrices of $X,X^c$ being related by the following formula: $$d_X+d_{X^c}=\mathbb I-1$$

By using now the above formula $u\mathbb I=\mathbb Iu=N\mathbb I$, we conclude that $d_Xu=ud_X$ is equivalent to $d_{X^c}u=ud_{X^c}$. Thus, we obtain, as claimed, $G^+(X)=G^+(X^c)$.

\medskip

(4) The adjacency matrix of a disconnected union is given by:
$$d_{X\sqcup Y}=\begin{pmatrix}
d_X&0\\
0&d_Y
\end{pmatrix}$$

Now let $w=diag(u,v)$ be the fundamental corepresentation of $G^+(X)\,\hat{*}\,G^+(Y)$. We have then $d_Xu=ud_X$ and $d_Yv=vd_Y$, and we obtain from this, as desired:
$$d_{X\sqcup Y}w=wd_{X\sqcup Y}$$

(5) We know from (3) that we have $G^+(\square)=G^+(|\ |)$. We know as well from (4) that we have $\mathbb Z_2\,\hat{*}\,\mathbb Z_2\subset G^+(|\ |)$. It follows that $G^+(\square)$ is non-classical. Finally, the inclusion $G^+(\square)\subset S_4^+$ is indeed proper, because $S_4\subset S_4^+$ does not act on the square.
\end{proof}

In order to further advance, and explicitly compute various quantum automorphism groups, we need some sharper techniques, based on the spectral decomposition of the adjacency matrix $d\in M_N(0,1)$. However, this decomposition is not something which is available over any field $F$, so we will make the convention, for the material that follows next, that such a decomposition exists indeed. With this convention, we first have:

\index{spectral decomposition}

\begin{theorem}
A subgroup $G\subset S_N^+$ acts on a graph $X$ precisely when
$$P_\lambda u=uP_\lambda\quad,\quad\forall\lambda\in F$$
where $d=\sum_\lambda\lambda\cdot P_\lambda$ is the spectral decomposition of the adjacency matrix of $X$.
\end{theorem}

\begin{proof}
With $d=\sum_\lambda\lambda\cdot P_\lambda$ as above, we have the following formula:
$$<d>=span\left\{P_\lambda\Big|\lambda\in F\right\}$$

Thus $d\in End(u)$ when $P_\lambda\in End(u)$ for all $\lambda\in F$, which gives the result.
\end{proof}

In practice, in order to exploit the invariance condition in Theorem 6.15, we can combine it with the following elementary fact, coming from definitions:

\begin{proposition}
Given a subgroup $G\subset S_N^+$, with associated coaction
$$\alpha:F^N\to F^N\otimes F(G)\quad,\quad e_i\to\sum_je_j\otimes u_{ji}$$
and a linear subspace $V\subset F^N$, the following are equivalent:
\begin{enumerate}
\item The magic matrix $u=(u_{ij})$ commutes with $P_V$.

\item $V$ is invariant, in the sense that $\alpha(V)\subset V\otimes F(G)$.
\end{enumerate}
\end{proposition}

\begin{proof}
Let $P=P_V$, with this being by definition allowed to be any projection having image $V$. For any $i\in\{1,\ldots,N\}$ we have the following formula:
\begin{eqnarray*}
\alpha(P(e_i))
&=&\alpha\left(\sum_kP_{ki}e_k\right)\\
&=&\sum_{jk}P_{ki}e_j\otimes u_{jk}\\
&=&\sum_je_j\otimes (uP)_{ji}
\end{eqnarray*}

On the other hand the linear map $(P\otimes id)\alpha$ is given by a similar formula:
\begin{eqnarray*}
(P\otimes id)(\alpha(e_i))
&=&\sum_kP(e_k)\otimes u_{ki}\\
&=&\sum_{jk}P_{jk}e_j\otimes u_{ki}\\
&=&\sum_je_j\otimes (Pu)_{ji}
\end{eqnarray*}

Thus $uP=Pu$ is equivalent to $\alpha P=(P\otimes id)\alpha$, and the conclusion follows.
\end{proof}

Next, we have the following result, complementary to Theorem 6.15, which again basically assumes that the ground field $F$ is the field of the complex numbers:

\index{color decomposition}
\index{color components}

\begin{theorem}
Given a matrix $p\in M_N(F)$, consider its color decomposition
$$p=\sum_{c\in F}c\cdot p_c$$
with the color components $p_c\in M_N(0,1)$ with $c\in F$ being constructed as follows:
$$(p_c)_{ij}=\begin{cases}
1&{\rm if}\ p_{ij}=c\\
0&{\rm otherwise}
\end{cases}$$
Then a magic matrix $u=(u_{ij})$ commutes with $p$ iff it commutes with all matrices $p_c$.
\end{theorem}

\begin{proof}
Consider the multiplication and counit maps of the algebra $F^N$:
$$M:e_i\otimes e_j\to e_ie_j\quad,\quad 
C:e_i\to e_i\otimes e_i$$

Since $M,C$ intertwine $u,u^{\otimes 2}$, their iterations $M^{(k)},C^{(k)}$ intertwine $u,u^{\otimes k}$, and so:
$$M^{(k)}p^{\otimes k}C^{(k)}
=\sum_{c\in F}c^kp_c
\in End(u)$$

Now since this formula holds for any $k\in\mathbb N$, we obtain the result.
\end{proof}

The above results can be combined, and we are led to the following statement:

\index{spectral-color decomposition}
\index{color-spectral decomposition}
\index{planar algebra}

\begin{theorem}
A subgroup $G\subset S_N^+$ acts on a graph $X$ precisely when 
$$u=(u_{ij})$$
commutes with all the matrices coming from the color-spectral decomposition of $d$.
\end{theorem}

\begin{proof}
This follows by combining Theorem 6.15 and Theorem 6.17, with the ``color-spectral'' decomposition in the statement referring to what comes out by succesively doing the color and spectral decomposition, until the process stabilizes.
\end{proof}

This latter statement is quite interesting, with the color-spectral decomposition there being something quite intriguing. We will be back to this later in this book, when discussing planar algebras, which is the good framework for discussing such things.

\bigskip

As something more concrete now, we have the following result:

\begin{theorem}
In order for a quantum permutation group $G\subset S_N^+$ to act on a graph $X$, having $N$ vertices, the adjacency matrix $d\in M_N(0,1)$ of the graph must be, when regarded as function on the set $\{1,\ldots,N\}^2$, constant on the orbitals of $G$. 
\end{theorem}

\begin{proof}
This follows indeed from the following isomorphism, from Theorem 6.10:
$$End(u)\simeq Fix(\sim)$$

For more on all this, details, examples, applications, we refer to the literature.
\end{proof}

We would like to understand now how the construction $X\to G^+(X)$ behaves, with respect to the various products of graphs. Let us start with the following definition:

\begin{definition}
Let $X,Y$ be two finite graphs.
\begin{enumerate}
\item The direct product $X\times Y$ has vertex set $X\times Y$, and edges:
$$(i,\alpha)-(j,\beta)\Longleftrightarrow i-j,\, \alpha-\beta$$

\item The Cartesian product $X\,\square\,Y$ has vertex set $X\times Y$, and edges:
$$(i,\alpha)-(j,\beta)\Longleftrightarrow i=j,\, \alpha-\beta\mbox{ \rm{or} }i-j,\alpha=\beta$$
\end{enumerate}
We call these operations the standard products of graphs.
\end{definition}

The above products are indeed well-known and standard in graph theory. In relation now with symmetry groups, we first have the following elementary result:

\begin{proposition}
We have embeddings as follows,
$$G^+(X)\times G^+(Y)\subset G^+(X \times Y)$$
$$G^+(X)\times G^+(Y)\subset G^+(X\,\square\,Y)$$
valid for any two graphs $X,Y$.
\end{proposition}

\begin{proof}
We use the following identification, given by $\delta_{(i,\alpha)}=\delta_i\otimes\delta_\alpha$:
$$F(X\times Y)=F(X)\otimes F(Y)$$

(1) The adjacency matrix of the direct product is given by:
$$d_{X\times Y}=d_X\otimes d_Y$$

Thus if $u$ commutes with $d_X$ and $v$ commutes with $d_Y$, then $u\otimes v=(u_{ij}v_{\alpha\beta})_{(i\alpha,j\beta)}$ is a magic matrix that commutes with $d_{X\times Y}$. But this gives a morphism as follows:
$$F(G^+(X\times Y))\to F(G^+(X)\times G^+(Y))$$ 

Finally, the surjectivity of this morphism follows by summing over $i$ and $\beta$.

\medskip

(2) The adjacency matrix of the Cartesian product is given by:
$$d_{X\,\square\,Y}=d_X\otimes1+1\otimes d_Y$$

Thus if $u$ commutes with $d_X$ and $v$ commutes with $d_Y$, then $u\otimes v=(u_{ij}v_{\alpha\beta})_{(i\alpha,j\beta)}$ is a magic unitary that commutes with $d_{X\,\square\,Y}$, and this gives the result.
\end{proof}

The problem now is that of deciding when the embeddings in Proposition 6.21 are isomorphisms. This is something quite subtle, and we have here the following result: 

\index{connected graph}
\index{regular graph}

\begin{theorem}
Let $X$ and $Y$ be finite connected regular graphs. If their spectra $\{\lambda\}$ and $\{\mu\}$ do not contain $0$ and satisfy
$$\big\{\lambda_i/\lambda_j\big\}\cap\big\{\mu_k/\mu_l\big\}=\{1\}$$
then $G^+(X\times Y)=G^+(X)\times G^+(Y)$. Also, if their spectra satisfy
$$\big\{\lambda_i-\lambda_j\big\}\cap\big\{\mu_k-\mu_l\big\}=\{0\}$$
then $G^+(X\,\square\,Y)=G^+(X)\times G^+(Y)$.
\end{theorem}

\begin{proof}
Let $\lambda_1$ be the valence of $X$. Since $X$ is regular we have $\lambda_1\in Sp(X)$, with $1$ as eigenvector, and since $X$ is connected $\lambda_1$ has multiplicity 1. Hence if $P_1$ is the orthogonal projection onto $F1$, the spectral decomposition of $d_X$ is of the following form:
$$d_X=\lambda_1P_1+\sum_{i\neq1}\lambda_iP_i$$

We have a similar formula for the adjacency matrix $d_Y$, namely:
$$d_Y=\mu_1Q_1+\sum_{j\neq1}\mu_jQ_j$$

But this gives the following formulae for the products:
$$d_{X\times Y}=\sum_{ij}(\lambda_i\mu_j)P_{i}\otimes Q_{j}$$
$$d_{X\,\square\,Y}=\sum_{ij}(\lambda_i+\mu_i)P_i\otimes Q_j$$

Here the projections form partitions of unity, and the scalars are distinct, so these are spectral decompositions. The coactions will commute with any of the spectral projections, and hence with both $P_1\otimes1$, $1\otimes Q_1$. In both cases the universal coaction $v$ is the tensor product of its restrictions to the images of $P_1\otimes1$, $1\otimes Q_1$, which gives the result.
\end{proof}

\section*{6d. Wreath products} 

Let us talk now about free wreath products, as a joint continuation of the above, and of the product operation material from Part I. Following Bichon \cite{bi2}, we have the following result, with our usual convention of formulating things in terms of quantum groups, according to the informal formula $A=F(G)$, for the Hopf algebras $A$:

\index{free wreath product}

\begin{theorem}
Given a quantum group $G$, and a subgroup $H\subset S_k^+$, with fundamental corepresentations $u,v$, the following construction produces a quantum group:
$$F(G\wr_*H)=(F(G)^{*k}*F(H))/<[u_{ij}^{(a)},v_{ab}]=0>$$
In the case where $G,H$ are classical, the classical version of $G\wr_*H$ is the usual wreath product $G\wr H$. Also, when $G$ is a quantum permutation group, so is $G\wr_*H$.
\end{theorem}

\begin{proof}
There are many things going on here, the idea being as follows:

\medskip

(1) Consider indeed the matrix $w_{ia,jb}=u_{ij}^{(a)}v_{ab}$, over the quotient algebra in the statement. This matrix is then a corepresentation, with respect to the standard Hopf algebra structure of the quotient algebra in the statement, and this gives the first assertion.

\medskip

(2) Next, we must prove that under the supplementary assumption $G\subset S_N^+$, this matrix $w_{ia,jb}=u_{ij}^{(a)}v_{ab}$ is magic. For this purpose, observe first that the entries of this matrix are projections, because they appear as products of commuting projections.

\medskip

(3) In what regards now the verification of $pq=0$ on rows, this goes as follows:
\begin{eqnarray*}
w_{ia,jb}w_{ia,kc}
&=&u_{ij}^{(a)}v_{ab}u_{ik}^{(a)}v_{ac}\\
&=&u_{ij}^{(a)}u_{ik}^{(a)}v_{ab}v_{ac}\\
&=&\delta_{jk}u_{ij}^{(a)}\delta_{bc}v_{ab}\\
&=&\delta_{jb,kc}w_{ia,jb}
\end{eqnarray*}

As for the verification of $pq=0$ on columns, this is quite similar, as follows:
\begin{eqnarray*}
w_{ia,jb}w_{kc,jb}
&=&u_{ij}^{(a)}v_{ab}u_{kj}^{(c)}v_{cb}\\
&=&u_{ij}^{(a)}u_{kj}^{(c)}v_{ab}v_{cb}\\
&=&u_{ij}^{(a)}u_{kj}^{(c)}\delta_{ac}v_{ab}\\
&=&\delta_{ac}u_{ij}^{(a)}u_{kj}^{(a)}v_{ab}\\
&=&\delta_{ac}\delta_{ik}u_{ij}^{(a)}v_{ab}\\
&=&\delta_{ia,kc}w_{ia,jb}
\end{eqnarray*}

(4) Next, the row sums for the matrix $w$ are computed as follows:
\begin{eqnarray*}
\sum_{jb}w_{ia,jb}
&=&\sum_{jb}u_{ij}^{(a)}v_{ab}\\
&=&\sum_bv_{ab}\sum_ju_{ij}^{(a)}\\
&=&1
\end{eqnarray*}

As for the column sums of $w$, these are computed similarly, as follows:
\begin{eqnarray*}
\sum_{ia}w_{ia,jb}
&=&\sum_{ia}u_{ij}^{(a)}v_{ab}\\
&=&\sum_av_{ab}\sum_iu_{ij}^{(a)}\\
&=&1
\end{eqnarray*}

(5) We conclude that in the case $G\subset S_N^+$, the matrix $w_{ia,jb}=u_{ij}^{(a)}v_{ab}$ is indeed magic, so we obtain a certain quantum permutation group $G\wr_*H\subset S_{Nk}^+$, as claimed. Finally, the assertion regarding the classical version is standard as well. See \cite{bi2}.
\end{proof}

In relation now with the quantum symmetry groups of graphs, let us recall that, as a complement to the product operations from Definition 6.20, we have as well:

\begin{definition}
Given two graphs $X,Y$, we can construct a graph having vertex set $X\times Y$, and edges as follows,
$$(i,\alpha)-(j,\beta)\Longleftrightarrow \alpha-\beta\mbox{ \rm{or} }\alpha=\beta,\, i-j$$
which is denoted $X\circ Y$, and called lexicographic product of $X,Y$.
\end{definition}

As before with the other products of graphs, from Definition 6.20, this is something quite standard in graph theory, and the interesting examples abound. The point now is that, somehow in analogy with Theorem 6.22, we have as well the following result, under our usual present assumptions, that the adjacency matrices are diagonalizable:

\begin{theorem}
Let $X,Y$ be regular graphs, with $X$ connected. If their spectra $\{\lambda_i\}$ and $\{\mu_j\}$ satisfy the condition
$$\big\{\lambda_1-\lambda_i\big|i\neq 1\big\}\cap\big\{-n\mu_j\big\}=\emptyset$$
where $n$ and $\lambda_1$ are the order and valence of $X$, then $G^+(X\circ Y)=G^+(X)\wr_*G^+(Y)$.   
\end{theorem}

\begin{proof}
This is something quite tricky, the idea being as follows:

\medskip

(1) First, it is clear from definitions that we have an embedding as follows:
$$G^+(X)\wr_*G^+(Y)\subset G^+(X\circ Y)$$

(2) Let us denote now by $P_i,Q_j$ the spectral projections corresponding to $\lambda_i,\mu_j$. Since $X$ was assumed to be connected we have $P_1=\mathbb I/n$, and we obtain:
\begin{eqnarray*}
d_{X\circ Y}
&=&d_X\otimes 1+{\mathbb I}\otimes d_Y\\
&=&\left(\sum_i\lambda_iP_i\right)\otimes\left(\sum_jQ_j\right)+\left(nP_1\right)\otimes \left(\sum_i\mu_jQ_j\right)\\
&=&\sum_j(\lambda_1+n\mu_j)(P_1\otimes Q_j) + \sum_{i\not=1}\lambda_i(P_i\otimes 1)
\end{eqnarray*} 

Since in this formula the projections form a partition of unity, and the scalars are distinct, we conclude that this is the spectral decomposition of $d_{X\circ Y}$. 

\medskip

(3) Now let $W$ be the universal magic matrix for $X\circ Y$. Then $W$ must commute with all spectral projections, and in particular, we have:
$$[W,P_1\otimes Q_j]=0$$

Summing over $j$ gives $[W,P_1\otimes 1]=0$, so $1\otimes F(Y)$ is invariant under the coaction. So, consider the restriction of $W$, which gives a coaction of $G^+(X\circ Y)$ on $1\otimes F(Y)$, that we can denote as follows, with $y$ being a certain magic unitary:
$$W(1\otimes e_a)=\sum_b1\otimes e_b\otimes y_{ba}$$

(4) On the other hand, according to our definition of $W$, we can write:
$$W(e_i\otimes 1)=\sum_{jb}e_j\otimes e_b\otimes x_{ji}^b$$  

By multiplying by the previous relation, found in (3), we obtain:
$$W(e_i\otimes e_a)
=\sum_{jb}e_j\otimes e_b\otimes y_{ba}x_{ji}^b
=\sum_{jb}e_j \otimes e_b\otimes x_{ji}^by_{ba}$$

But this shows that the coefficients of $W$ are of the following form:
$$W_{jb,ia}=y_{ba}x_{ji}^b=x_{ji}^b y_{ba}$$

(5) Consider now the matrix $x^b=(x_{ij}^b)$. Since $W$ is a morphism of algebras, each row of $x^b$ is a partition of unity. Also, by using the antipode, we have:
$$S\left(\sum_jx_{ji}^{b}\right)
=S\left(\sum_{ja}W_{jb,ia}\right)
=\sum_{ja}W_{ia,jb}
=1$$

As a conclusion to this, the matrix $x^b$ constructed above is magic. 

\medskip

(6) We check now that both $x^a,y$ commute with $d_X,d_Y$. We have:
$$(d_{X\circ Y})_{ia,jb} = (d_X)_{ij}\delta_{ab} + (d_Y)_{ab}$$

Thus the two products between $W$ and $d_{X\circ Y}$ are given by:
$$(Wd_{X\circ Y})_{ia,kc}=\sum_j W_{ia,jc} (d_X)_{jk} + \sum_{jb}W_{ia,jb}(d_Y)_{bc}$$
$$(d_{X\circ Y}W)_{ia,kc}=\sum_j (d_X)_{ij} W_{ja,kc} + \sum_{jb}(d_Y)_{ab}W_{jb,kc}$$

(7) Now since the magic matrix $W$ commutes by definition with $d_{X\circ Y}$, the terms on the right in the above equations are equal, and by summing over $c$ we get:
$$\sum_j x_{ij}^a(d_X)_{jk} + \sum_{cb} y_{ab}(d_Y)_{bc}
= \sum_{j} (d_X)_{ij}x_{jk}^a + \sum_{cb} (d_Y)_{ab}y_{bc}$$

The second sums in both terms are equal to the valence of $Y$, so we get:
$$[x^a,d_X]=0$$

Now once again from the formula coming from $[W,d_{X\circ Y}]=0$, we get:
$$[y,d_Y] =0$$

(8) Summing up, the coefficients of $W$ are of the following form, where $x^b$ are magic matrices commuting with $d_X$, and $y$ is a magic matrix commuting with $d_Y$: 
$$W_{jb,ia}=x_{ji}^by_{ba}$$

But this gives a morphism $F(G^+(X)\wr_*G^+(Y))\to F(G^+(X\circ Y))$ mapping $u_{ji}^{(b)}\to x_{ji}^b$ and $v_{ba}\to y_{ba}$, which is inverse to the morphism in (1), as desired.
\end{proof}

As a main application of the above result, we have:

\index{connected graph}

\begin{theorem}
Given a connected graph $X$, and $k\in\mathbb N$, we have the formulae
$$G(kX)=G(X)\wr S_k\quad,\quad 
G^+(kX)=G^+(X)\wr_*S_k^+$$
where $kX=X\sqcup\ldots\sqcup X$ is the $k$-fold disjoint union of $X$ with itself.
\end{theorem}

\begin{proof}
There are several things going on here, the idea being as follows:

\medskip

(1) The first formula is something very intuitive and well-known, and technically follows as well from the second formula, by taking the classical version. 

\medskip

(2) Regarding the second formula, this follows from Theorem 6.25, but let us discuss this in detail. Our first claim is that we have an inclusion as follows, for any graph $X$:
$$G^+(X)\wr_*S_k^+\subset G^+(kX)$$

Indeed, we want to construct an action $G^+(X)\wr_*S_k^+\curvearrowright kX$, and this amounts in proving that we have $[w,d]=0$. But, the matrices $w,d$ are given by:
$$w_{ia,jb}=u_{ij}^{(a)}v_{ab}\quad,\quad d_{ia,jb}=\delta_{ij}d_{ab}$$

With these formulae in hand, we have the following computation:
\begin{eqnarray*}
(dw)_{ia,jb}
&=&\sum_kd_{ik}w_{ka,jb}\\
&=&\sum_kd_{ik}u_{kj}^{(a)}v_{ab}\\
&=&(du^{(a)})_{ij}v_{ab}
\end{eqnarray*}

On the other hand, we have as well the following computation:
\begin{eqnarray*}
(wd)_{ia,jb}
&=&\sum_kw_{ia,kb}d_{kj}\\
&=&\sum_ku_{ik}^{(a)}v_{ab}d_{kj}\\
&=&(u^{(a)}d)_{ij}v_{ab}
\end{eqnarray*}

Thus we have commutation, $[w,d]=0$, and from this we obtain, as claimed:
$$G^+(X)\wr_*S_k^+\subset G^+(kX)$$

(3) Regarding now the reverse inclusion, which requires $X$ to be connected, this follows by doing some matrix analysis, by using the commutation with $u$. To be more precise, let us denote by $w$ the fundamental corepresentation of $G^+(kX)$, and set:
$$u_{ij}^{(a)}=\sum_bw_{ia,jb}\quad,\quad 
v_{ab}=\sum_iv_{ab}$$

It is then routine to check that we have here magic matrices. Thus we obtain the reverse inclusion, that we were looking for, and this gives the result, as desired.

\medskip

(4) Alternatively, the reverse inclusion follows as well from Theorem 6.25, because the spectral condition there is trivially satisfied, in the present case.
\end{proof}

Good news, we can now construct quantum reflection groups. Let us start with:

\begin{definition}
A $(s,N)$-sudoku matrix is a magic matrix of size $sN$, of the form
$$m=\begin{pmatrix}
a^0&a^1&\ldots&a^{s-1}\\
a^{s-1}&a^0&\ldots&a^{s-2}\\
\vdots&\vdots&&\vdots\\
a^1&a^2&\ldots&a^0
\end{pmatrix}$$
where $a^0,\ldots,a^{s-1}$ are $N\times N$ matrices.
\end{definition}

The basic examples of such matrices come from the group $H_N^s=\mathbb Z_s\wr S_N$. Indeed, with $w=e^{2\pi i/s}$, each of the $N^2$ matrix coordinates $u_{ij}:H_N^s\to\mathbb C$ decomposes as follows:
$$u_{ij}=\sum_{r=0}^{s-1}w^ra^r_{ij}$$

Here each $a^r_{ij}$ is a function taking values in $\{0,1\}$, and so a projection in the algebraic sense, and it follows from definitions that these projections form a sudoku matrix. 

\bigskip

In fact, regarding the group $H_N^s=\mathbb Z_s\wr S_N$, we have the following result:

\begin{theorem}
The algebra of functions $F(H_N^s)$ on the complex reflection group
$$H_N^s=\mathbb Z_s\wr S_N$$
is the universal commutative algebra generated by the entries of a $(s,N)$-sudoku matrix.
\end{theorem}

\begin{proof}
This is indeed something which is clear from definitions, and from the above observations, and we will leave working out the details here, as an instructive exercise.
\end{proof}

Now based on the above result, we can formulate, a bit as we did before, in the beginning of this chapter, for the permutation groups:

\begin{theorem}
The universal algebra $F(H_N^{s+})$ generated by the entries of a $(s,N)$-sudoku matrix is a Hopf algebra, which appears as a liberation of $F(H_N^s)$.
\end{theorem}

\begin{proof}
As before with Theorem 6.28, this is again something which is clear from definitions, and we will leave working out the details here, as an instructive exercise.
\end{proof}

We call the quantum group type objects $H_N^{s+}$ appearing above quantum reflection groups. Regarding now the study of these quantum groups, we first have:

\index{reflection group}
\index{quantum reflection group}

\begin{theorem}
We have the following results:
\begin{enumerate}
\item $H_N^s=\mathbb Z_s\wr S_N$.

\item $H_N^{s+}=\mathbb Z_s\wr_*S_N^+$.
\end{enumerate}
\end{theorem}

\begin{proof}
Here the first assertion is something that we already know, coming from Theorem 6.28, and the proof of the second assertion is standard, coming from definitions, that we will leave here as an exercise. We will be actually back in a moment to this second assertion, with an alternative proof as well, using our graph technology.
\end{proof}

In order to further advance, we will need the following simple fact: 

\begin{proposition}
A $sN\times sN$ magic unitary commutes with the matrix
$$\Sigma=
\begin{pmatrix}
0&I_N&0&\ldots&0\\
0&0&I_N&\ldots&0\\
\vdots&\vdots&&\ddots&\\
0&0&0&\ldots&I_N\\
I_N&0&0&\ldots&0
\end{pmatrix}$$
precisely when it is a sudoku matrix in the sense of Definition 6.27.
\end{proposition}

\begin{proof}
This follows from the fact that commutation with $\Sigma$ means that the matrix is circulant. Thus, we obtain the sudoku relations from Definition 6.27.
\end{proof}

Now let $Z_s$ be the oriented cycle with $s$ vertices, and consider the graph $NZ_s$ consisting of $N$ disjoint copies of it. Observe that, with a suitable labeling of the vertices, the adjacency matrix of this graph is precisely the above matrix $\Sigma$. We have:

\begin{theorem}
We have the following results:
\begin{enumerate}
\item $H_N^s$ is the symmetry group of $NZ_s$.

\item $H_N^{s+}$ is the quantum symmetry group of $NZ_s$.
\end{enumerate}
\end{theorem}

\begin{proof}
This is again something elementary, the idea being as follows:

\medskip

(1) This follows indeed from definitions.

\medskip

(2) This follows by combining the various results above, because the algebra $F(H_N^{s+})$ follows to be the quotient of the algebra $F(S_{sN}^+)$ by the relations making the fundamental corepresentation commute with the adjacency matrix of $NZ_s$.
\end{proof}

Many other things can be said about the complex reflection groups $H_N^{s+}$, and notably about the hyperoctahedral quantum group $H_N^+=H_N^{2+}$. We will be back to this.

\section*{6e. Exercises}

We had a lot of interesting algebra in this chapter, and as exercises, we have:

\begin{exercise}
Find the simplest proof ever for the isomorphism $S_3=S_3^+$.
\end{exercise}

\begin{exercise}
Try computing the quantum symmetry group of the $N$-cycle.
\end{exercise}

\begin{exercise}
Further investigate the spectral-color decomposition.
\end{exercise}

\begin{exercise}
Read about the classification of complex reflection groups.
\end{exercise}

\begin{exercise}
Try to compute the quantum symmetry group of the hypercube.
\end{exercise}

\begin{exercise}
Study more the quantum reflection groups constructed above.
\end{exercise}

As bonus exercise, study the representations of our various quantum groups.

\chapter{Standard twists}

\section*{7a. Standard twists}

We have seen in chapters 5-6 that many interesting examples of Hopf algebras and quantum groups can be constructed via liberation, and half-liberation:
$$ab=ba\quad\to\quad\emptyset$$
$$ab=ba\quad\to\quad abc=cba$$

We discuss in this chapter, and in the next one too, another operation of the same type, called twisting. We will see that this leads to many new examples of quantum groups, which are sometimes related, in a quite subtle way, to liberation.

\bigskip

Getting started, there are many ways of twisting, and we will first discuss the simplest twists, namely the anticommutation ones. The idea is very simple, consisting of suitably replacing the commutation relations with anticommutation relations:
$$ab=ba\quad\to\quad ab=\pm ba$$
$$abc=cba\quad\to\quad abc=\pm cba$$

However, in practice this is something quite tricky, which will take some time, to be developed. The problem comes from the choice of commutation vs anticommutation, for each pair of relevant variables $\{a,b\}$, typically coefficients of $u,\bar{u}$, which, as we will soon discover, leads us into some combinatorics, which is not exactly of trivial type. 

\bigskip

Let us formulate, to start with, the following definition:

\begin{definition}
We say that an invertible matrix $u\in M_N(A)$ is twisted when
$$ab=\begin{cases}
-ba&{\rm for}\ a,b\ {\rm at\ distinct\ places,\ on\ the\ same\ row\ or\ column\ of\ }u^\times\\
ab&{\rm otherwise}
\end{cases}$$
with $u^\times$ standing for the $N\times N$ array obtained by superposing $u$ and $v=(id\otimes t)u^{-1}$.
\end{definition}

This notion looks like something quite tricky, so let us see how this works. Taking $N=2$, consider an invertible matrix $u\in M_N(A)$, set $v=(id\otimes t)u^{-1}$, and write:
$$u=\begin{pmatrix}u_{11}&u_{12}\\ u_{21}&u_{22}\end{pmatrix}
\quad,\quad 
v=\begin{pmatrix}v_{11}&v_{12}\\ v_{21}&v_{22}\end{pmatrix}$$

The $2\times2$ array obtained by superposing these two matrices is then as follows:
$$u^\times=\begin{pmatrix}u_{11},v_{11}&&u_{12},v_{12}\\ \\u_{21},v_{21}&&u_{22},v_{22}\end{pmatrix}$$

With $[[a,b]]=0$ standing for anticommutation, the anticommutation relations to be imposed, coming from the two rows, then the two columns, are then as follows:
$$[[u_{11},u_{12}]]=[[u_{11},v_{12}]]=[[v_{11},u_{12}]]=[[v_{11},v_{12}]]=0$$
$$[[u_{21},u_{22}]]=[[u_{21},v_{22}]]=[[v_{21},u_{22}]]=[[v_{21},v_{22}]]=0$$
$$[[u_{11},u_{21}]]=[[u_{11},v_{21}]]=[[v_{11},u_{21}]]=[[v_{11},v_{21}]]=0$$
$$[[u_{12},u_{22}]]=[[u_{12},v_{22}]]=[[v_{12},u_{22}]]=[[v_{12},v_{22}]]=0$$

As for the remaining pairs of variables, these are all subject to commutation, $[a,b]=0$. That is, the commutation relations to be imposed are as follows:
$$u_{11},v_{11},u_{22},v_{22}\ {\rm commute}$$
$$u_{12},v_{12},u_{21},v_{21}\ {\rm commute}$$

In general, at arbitrary $N\in\mathbb N$, the picture is quite similar, and with the remark here that at $N>>0$ the variables $u_{ij},v_{ij}$ generically commute, with the anticommutation being reserved for what happens on rows and columns, as indicated in Definition 7.1.

\bigskip

We can now twist the general linear groups, as follows:

\index{twisting}
\index{ad-hoc twisting}
\index{anticommutation}

\begin{theorem}
We have twisted general linear groups
$$GL_N'\subset GL_N^+$$
obtained by assuming that the fundamental corepresentation $u$ is twisted.
\end{theorem}

\begin{proof}
This is somehing quite routine, the idea being as follows:

\medskip

(1) To start with, let us first recall from chapter 5 that the free general linear groups $GL_N^+$ are constructed according to the following formula:
$$F(GL_N^+)=\left<(u_{ij},v_{ij})_{i,j=1,\ldots,N}\,\Big|\,u^{-1}=v^t,\ v^{-1}=u^t\right>$$

To be more precise, as explained in chapter 5, what we have here is a Hopf algebra, with comultiplication, counit and antipode making $u,v$ corepresentations, as follows:
$$\Delta(u_{ij})=\sum_ku_{ij}\otimes u_{kj}\quad,\quad\Delta(v_{ij})=\sum_kv_{ij}\otimes v_{kj}$$
$$\varepsilon(u_{ij})=\varepsilon(v_{ij})=\delta_{ij}$$
$$S(u_{ij})=v_{ji}\quad,\quad S(v_{ij})=u_{ji}$$

(2) Now let us construct $GL_N'\subset GL_N^+$ as in the statement, by performing the following quotient operation, with the relations being those in Definition 7.1:
$$F(GL_N')=F(GL_N^+)\Big/\Big<u={\rm twisted}\Big>$$

In order to prove that what we get is indeed a quantum group, we must prove that we have here a Hopf algebra. That is, we must construct factorizations as follows:
$$\Delta:F(GL_N')\to F(GL_N')\otimes F(GL_N')$$
$$\varepsilon:F(GL_N')\to F$$
$$S:F(GL_N')\to F(GL_N')^{opp}$$

(3) Let us first discuss the construction of the factorized comultiplication map $\Delta$. For this purpose, consider the following elements:
$$U_{ij}=\sum_ku_{ik}\otimes u_{kj}\quad, \quad V_{ij}=\sum_kv_{ik}\otimes v_{kj}$$

We know that the original comultiplication map $\Delta$, that of the Hopf algebra $F(GL_N^+)$, comes via the following formulae, on the standard generators:
$$\Delta(u_{ij})=U_{ij}\quad,\quad\Delta(v_{ij})=V_{ij}$$

Thus, in order to prove that $\Delta$ factorizes indeed, we must prove that we have:
$$u={\rm twisted}\implies U={\rm twisted}$$

(4) So, let us prove this. For any $i$ and $j\neq k$ we have, with $w=u,v$:
\begin{eqnarray*}
U_{ij}W_{ik}
&=&\sum_{st}u_{is}w_{it}\otimes u_{sj}w_{tk}\\
&=&\sum_{s\neq t}u_{is}w_{it}\otimes u_{sj}w_{tk}+\sum_su_{is}w_{is}\otimes u_{sj}w_{sk}\\
&=&\sum_{s\neq t}(-w_{it}u_{is})\otimes w_{tk}u_{sj}+\sum_sw_{is}u_{is}\otimes(-w_{sk}u_{sj})\\
&=&-\sum_{st}w_{it}u_{is}\otimes w_{tk}u_{sj}\\
&=&-W_{ik}U_{ij}
\end{eqnarray*}

Similarly, for any $i$ and $j\neq k$ we have, again with $w=u,v$:
$$U_{ji}W_{ki}=-W_{ki}U_{ji}$$

Thus, the anticommutation relations from Definition 7.1 are satisfied for $U$. 

\medskip

(5) Regarding now the commutation relations to be verified, for $i\neq k,j\neq l$ we have the following computation, as before with $w=u,v$:
\begin{eqnarray*}
U_{ij}W_{kl}
&=&\sum_{st}u_{is}w_{kt}\otimes u_{sj}w_{tl}\\
&=&\sum_{s\neq t}u_{is}w_{kt}\otimes u_{sj}w_{tl}+\sum_su_{is}w_{ks}\otimes u_{sj}w_{sl}\\
&=&\sum_{s\neq t}w_{kt}u_{is}\otimes w_{tl}u_{sj}+\sum_s(-w_{ks}u_{is})\otimes(-w_{sl}u_{sj})\\
&=&\sum_{st}w_{kt}u_{is}\otimes w_{tl}u_{sj}\\
&=&W_{kl}U_{ij}
\end{eqnarray*}

Finally, for any $i,j$ we have the following computation:
\begin{eqnarray*}
U_{ij}V_{ij}
&=&\sum_{st}u_{is}v_{it}\otimes u_{sj}v_{tj}\\
&=&\sum_{s\neq t}u_{is}v_{it}\otimes u_{sj}v_{tj}+\sum_su_{is}v_{is}\otimes u_{sj}v_{sj}\\
&=&\sum_{s\neq t}(-v_{it}u_{is})\otimes(-v_{tj}u_{sj})+\sum_sv_{is}u_{is}\otimes v_{sj}u_{sj}\\
&=&\sum_{s\neq t}v_{it}u_{is}\otimes v_{tj}u_{sj}+\sum_sv_{is}u_{is}\otimes v_{sj}u_{sj}\\
&=&\sum_{st}v_{it}u_{is}\otimes v_{tj}u_{sj}\\
&=&V_{ij}U_{ij}
\end{eqnarray*}

Thus, the commutation relations from Definition 7.1 are satisfied for $U$. 

\medskip

(6) Summarizing, we can define a comultiplication map for $F(GL_N')$, by using the universal property of the algebra $F(GL_N')$, by setting:
$$\Delta(u_{ij})=U_{ij}\quad,\quad\Delta(v_{ij})=V_{ij}$$

(7) Regarding now the counit $\varepsilon$ and the antipode $S$, things are clear here, by using the same method, and with no computations needed, the formulae to be satisfied being trivially satisfied. We conclude that $GL_N'$ is indeed a quantum group, as desired.
\end{proof}

The above result is quite nice, as stated, but all in all, there are many ingredients there, and relations that are imposed. To illustrate this, at $N=2$ we have:

\begin{theorem}
The quantum group $GL_2'\subset GL_2^+$ appears from
$$F(GL_2^+)=\left<(u_{ij},v_{ij})_{i,j=1,2}\,\Big|\,u^{-1}=v^t,\ v^{-1}=u^t\right>$$
by further imposing the following commutation and anticommutation relations:
\begin{enumerate}
\item $u_{11},v_{11},u_{22},v_{22}$ commute.

\item $u_{12},v_{12},u_{21},v_{21}$ commute.

\item Everything else anticommutes.
\end{enumerate}
\end{theorem}

\begin{proof}
This is indeed a reformulation of Theorem 7.2 at $N=2$, by taking into account the $N=2$ discussion made after Definition 7.1.
\end{proof}

At $N=3$ things get obviously more complicated, and formulating a simple statement as Theorem 7.3 is not an easy task. Theorem 7.2 remains our main result there.

\bigskip

At at opposite end now, $N>>0$, let us make the following observation:

\begin{proposition}
The standard coordinates $u_{ij},\bar{u}_{ij}$ of the quantum group 
$$GL_N'\subset GL_N^+$$
with $N>>0$ generically commute.
\end{proposition}

\begin{proof}
This is just a theoretical remark, of rather statistical nature, coming from Definition 7.1. Indeed, we have a total of $2N^2$ variables of type $u_{ij},\bar{u}_{ij}$, which are subject to a total of $\binom{2N^2}{2}$ commutation and anticommutation relations, as follows:

\medskip

(1) First, we have the anticommutation relations on the rows and columns of the double matrix $u^\times=u,\bar{u}$, and since there are $2N$ such rows and columns, with each producing $4\binom{N}{2}$ anticommutation relations, the total number of anticommutation relations is:
$$A=8N\binom{N}{2}=4N^2(N-1)$$

(2) And then, everything else must commute, so the number of commutation relations is given by the following formula, which asympotically means $C\simeq NA/2$:
$$C=\binom{2N^2}{2}-A=N^2(2N^2-4N+3)$$

(3) Thus, we are led to the conclusion in the statement. Before leaving, however, let us take $N=2$, as a matter to see if we have not messed up our computations, and if our figures above are correct. And here, our formula is $C=4\times 3=12$, which is in agreement indeed with Theorem 7.3, and the $6+6=12$ commutation relations there.
\end{proof}

Many other things can be said, as a continuation of this. For instance, we can try to impose twisting conditions to the fundamental corepresentation of $GL_Q^+$, given by:
$$F(GL_Q^+)=\left<(u_{ij},v_{ij})_{i,j=1,\ldots,N}\,\Big|\,u^{-1}=v^t,\ v^{-1}=Qu^tQ^{-1}\right>$$

However, things are quite technical here, and a bit ill-posed too, because it does not make much sense, conceptually speaking, to deform and twist at the same time. We will discuss such questions in the next section, in the orthogonal group setting, where things get more interesting, due to the different nature of the deformation operation there.

\section*{7b. Half-liberation} 

Another interesting question, that we will get into now, is that of twisting the half-liberations introduced in chapter 5, by replacing the half-commutation relations there $abc=cba$ with twisted half-commutation relations, of the following type:
$$abc=\pm cba$$

Let us begin with a summary of what we know about half-liberation, since chapter 5, along with a bit more, in relation with the material from chapter 6, as follows:

\begin{theorem}
We have quantum groups $G\subset G^*\subset G^+$, constructed via
$$F(G^*)=F(G^+)\Big/\left<abc=cba\Big|\forall a,b,c\in\{u_{ij},\bar{u}_{ij}\}\right>$$
which have the following properties:
\begin{enumerate}
\item The projective versions $PG^*$ are classical.

\item The quantum groups $G^*$ are not classical, for $G=GL_N,O_N,H_N^s$ with $s\geq2$.

\item However, collapsing $G=G^*$ can happen, for instance for $G=B_N,S_N$.
\end{enumerate}
\end{theorem}

\begin{proof}
This is something from chapter 5, along with a bit more, as follows:

\medskip

(1) In what regards the first assertion, this is something immediate, coming from:
$$abcd=cbad=cdab$$

(2) This comes by looking at the diagonal algebras $F(G^*)^\delta$, which are respectively equal to the group algebras of $H=\mathbb Z^{\circ N},\mathbb Z_2^{\circ N},\mathbb Z_s^{\circ N}$, with the symbol $\circ$ standing for the half-classical product of groups, obtained via $ghk=khg$. Now since these latter groups $H$ are not abelian, nor are the corresponding group algebras commutative, as desired.

\medskip

(3) Finally, for $G=S_N,B_N$ the above argument involving diagonal algebras will not work, because these diagonal algebras collapse to $F$. In fact, as said in the statement, the half-liberation procedure fails here, because by summing the matrix entries over a row, the relations $abc=cba$ get converted into the commutation relations $ab=ba$.
\end{proof}

As explained in chapter 5, the above result is in fact just the tip of the iceberg, and there are far more things that can be said about half-liberation, including:

\bigskip

(1) A uniqueness result for $G\subset G^*\subset G^+$, and in fact about half-liberation itself, viewed as an abstract operation, under strong combinatorial assumptions.

\bigskip

(2) The fact that the antidiagonal $2\times2$ scalar matrices automatically half-commute, which potentially produces interesting matrix models for the algebras $F(G^*)$.

\bigskip

(3) The related fact that, with the projective version $PG^*$ being classical, the half-liberation $G^*$ can be potentially understood as being the maximal affine lift of $PG^*$.

\bigskip

(4) And finally, some tricks for half-liberating however the bistochastic group $B_N$, via a related procedure, but with these tricks still not applying, however, to $S_N$.

\bigskip

Quite interesting all this, obviously the half-classical world is quite rich, and worth a more detailed study, and as main questions now, coming from this, let us formulate:

\begin{questions}
In relation with the half-liberation operation:
\begin{enumerate}
\item How to further develop its abstract aspects mentioned above?

\item How to develop algebraic and differential geometry, in this setting?

\item Is there any intermediate quantum group $S_N\subset S_N^\times\subset S_N^+$ at all?

\item How to adapt the half-liberation procedure, to the twisted setting?
\end{enumerate}
\end{questions}

We will discuss these questions in what follows, with a look right next at (4), in the remainder of this section, then with a look at (3) too, which is in fact related to twisting, afterwards, and then with a look at (1,2) too, later in this book, directly over the complex numbers, where the theory can be fully developed, along the lines mentioned above.

\bigskip

Getting to work now, in order to solve (4), we are facing the same problem as in the beginning of this chapter, namely lots of choices commutation vs anticommutation. However, things get considerably more complicated in the half-classical setting, because the formulae $abc=\pm cba$ involve far more variables than the formulae $ab=\pm ba$.

\bigskip

So, we will trick here, by using the following noncommutative geometry idea:

\begin{idea}
In order to adapt the half-liberation procedure, to the twisted setting,
\begin{enumerate}
\item We will first construct quantum spheres $S^{N-1}\subset S^{N-1}_*\subset S^{N-1}_+$,

\item Then, we will twist these spheres, into $S^{N-1,\prime}\subset S^{N-1,\prime}_*\subset S^{N-1}_+$,

\item And then, we will talk about $O_N'\subset O_N^{*\prime}\subset O_N^+$, and other quantum groups,
\end{enumerate}
the point being that the sphere case is easier, the sphere coordinates having single indices.
\end{idea}

So, this was for the general idea, and if you have some troubles in understanding, please be sure that the same happened to me, when I first learned this, from my cats. Getting started now, we can formally define a free sphere, in the following way:
$$F(S^{N-1}_+)=F\left(z_1,\ldots,z_N\Big|\sum_iz_i^2=1\right)$$

By imposing commutation relations between the coordinates, we obtain:
$$F(S^{N-1})=F\left(z_1,\ldots,z_N\Big|\sum_iz_i^2=1,z_iz_j=z_jz_i\right)$$

We can define as well a half-classical sphere, in the following way:
$$F(S^{N-1}_*)=F\left(z_1,\ldots,z_N\Big|\sum_iz_i^2=1,z_iz_jz_k=z_kz_jz_i\right)$$

Summarizing, we have three spheres, with inclusions between them, as follows:
$$S^{N-1}\subset S^{N-1}_*\subset S^{N-1}_+$$

Now let us try to twist these spheres. The free sphere is obviously not twistable. Regarding now the classical sphere, this is easy to twist, as follows:

\begin{definition}
We can construct a twisted sphere as follows,
$$F(S^{N-1,\prime})=F\left(z_1,\ldots,z_N\Big|\sum_iz_i^2=1,z_iz_j=(-1)^{\varepsilon_{ij}}z_jz_i\right)$$
where $\varepsilon_{ij}=1-\delta_{ij}$, with $\delta_{ij}$ being Kronecker symbols.
\end{definition}

Observe that this construction is compatible with the quantum group twisting method developed before, via standard actions of type $O_N\curvearrowright S^{N-1}$, and we will leave some computations here as an instructive exercise. In fact, retrospectively thinking, the signs in Definition 7.1 can be thought of as coming from those in Definition 7.8.

\bigskip

Getting now to half-commutation, we can use the same idea, and do some reverse engineering, by first twisting the spheres, and then the quantum groups themselves. Indeed, we can first twist the half-classical sphere, in a quite standard way, as follows:

\begin{definition}
We can construct a twisted half-classical sphere as follows,
$$F(S^{N-1,\prime}_*)=F\left(z_1,\ldots,z_N\Big|\sum_iz_i^2=1,z_iz_jz_k=(-1)^{\varepsilon_{ij}+\varepsilon_{jk}+\varepsilon_{ik}}z_kz_jz_i\right)$$
where $\varepsilon_{ij}=1-\delta_{ij}$ as before, with $\delta_{ij}$ being Kronecker symbols.
\end{definition}

To be more precise here, our claim is that the commutation relations in the above are uniquely determined, in order to have inclusions of twisted spheres, as follows:
$$S^{N-1,\prime}\subset S^{N-1,\prime}_*\subset S^{N-1}_+$$

Indeed, this claim basically comes from the following computation:
\begin{eqnarray*}
z_iz_jz_k
&=&(-1)^{\varepsilon_{ij}}z_jz_iz_k\\
&=&(-1)^{\varepsilon_{ij}+\varepsilon_{ik}}z_jz_kz_i\\
&=&(-1)^{\varepsilon_{ij}+\varepsilon_{jk}+\varepsilon_{ik}}z_kz_jz_i
\end{eqnarray*}

Summarizing, we have solved the twisting question at the level of the spheres, in a quite standard way, with our conclusion being as follows:

\begin{conclusion}
We have twisted quantum spheres, as follows,
$$S^{N-1,\prime}\subset S^{N-1,\prime}_*\subset S^{N-1}_+$$
in analogy with the usual spheres $S^{N-1}\subset S^{N-1}_*\subset S^{N-1}_+$.
\end{conclusion}

We must do now some further computations, in order to deduce from this the solution of the twisting question for the quantum groups themselves:

\bigskip

(1) In what regards the twisting of the classical groups, there is nothing to be done here, because as already mentioned after Definition 7.8, this method leads to the signs in Definition 7.1. In fact, this is how the signs in Definition 7.1 were originally found.

\bigskip

(2) In what regards now the twisting of the half-classical groups, things here are considerably more complicated, and our sphere technology above is more or less the only way of reaching to the correct signs. We will explain this, in what follows.

\bigskip

So, let us explain the solution. Given three coordinates $a,b,c\in\{u_{ij}\}$, let us set, where $r,c\in\{1,2,3\}$ are the number of rows and columns spanned by $a,b,c$:
$$span(a,b,c)=(r,c)$$

In other words, if we write $a=u_{ij},b=u_{kl},c=u_{pq}$ then $r,c\in\{1,2,3\}$ are given by:
$$r=\#\{i,k,p\}\quad,\quad l=\#\{j,l,q\}$$

With this convention, we can now formulate, in analogy with Definition 7.1:

\begin{definition}
We say that an invertible matrix $u\in M_N(A)$ is half-twisted when
$$abc=\begin{cases}
-cba&{\rm for}\ a,b,c\ {\rm with}\ span(a,b,c)=(\leq 2,3)\ {\rm or}\ (3,\leq 2)\\
cba&{\rm otherwise}
\end{cases}$$
for $a,b,c$ belonging to the $N\times N$ array obtained by superposing $u$ and $v=(id\otimes t)u^{-1}$.
\end{definition}

As before with Definition 7.1, this notion is something quite tricky, so let us first see how this works. As a first observation, the rules for the various commutation and anticommutation signs in the statement can be summarized as follows:
$$\begin{matrix}
r\backslash c&1&2&3\\
1&+&+&-\\
2&+&+&-\\
3&-&-&+
\end{matrix}$$

At the level of examples now, things are a bit more complicated than before, because at $N=2$ we cannot have triples $a,b,c$ spanning 3 rows or columns, and so all signs in Definition 7.11 are $+$. Thus, the first non-trivial examples appear at $N=3$. So, let us take $N=3$, consider an invertible matrix $u\in M_N(A)$, set $v=(id\otimes t)u^{-1}$, and write:
$$u=\begin{pmatrix}u_{11}&u_{12}&u_{13}\\ u_{21}&u_{22}&u_{23}\\ u_{31}&u_{32}&u_{33}\end{pmatrix}
\quad,\quad 
v=\begin{pmatrix}v_{11}&v_{12}&v_{13}\\ v_{21}&v_{22}&v_{23}\\ v_{31}&v_{32}&v_{33}\end{pmatrix}$$

The $3\times3$ array obtained by superposing these two matrices is then as follows:
$$u^\times=\begin{pmatrix}u_{11},v_{11}&&u_{12},v_{12}&&u_{13},v_{13}\\ 
\\
u_{21},v_{21}&&u_{22},v_{22}&&u_{23},v_{23}\\
\\
u_{31},v_{31}&&u_{32},v_{32}&&u_{33},v_{33}
\end{pmatrix}$$

This array has 18 entries, and so there are $18^3=5832$ relations to be imposed, and as an exercise for you, find out how many $abc=cba$ and $abc=-cba$ relations there are.

\bigskip

In general, at arbitrary $N\in\mathbb N$, the picture is quite similar, and with the remark here that at $N>>0$ the variables $u_{ij},v_{ij}$ must generically half-commute, and this because the number of rows, and of columns as well, spanned by a triple $a,b,c$ is generically 3.

\bigskip

We can now twist the half-classical general linear groups, as follows:

\index{twisting}
\index{ad-hoc twisting}
\index{anticommutation}

\begin{theorem}
We have twisted half-classical general linear groups
$$GL_N^{*\prime}\subset GL_N^+$$
obtained by assuming that the fundamental corepresentation $u$ is half-twisted.
\end{theorem}

\begin{proof}
This is somehing quite routine, the idea being as follows:

\medskip

(1) To start with, let us first recall from chapter 5 that the free general linear groups $GL_N^+$ are constructed according to the following formula:
$$F(GL_N^+)=\left<(u_{ij},v_{ij})_{i,j=1,\ldots,N}\,\Big|\,u^{-1}=v^t,\ v^{-1}=u^t\right>$$

Now let us construct $GL_N^{*\prime}\subset GL_N^+$ as in the statement, by performing the following quotient operation, with the relations being those in Definition 7.11:
$$F(GL_N^{*\prime})=F(GL_N^+)\Big/\Big<u={\rm half\!-\!twisted}\Big>$$

(2) Let us first discuss the construction of the factorized comultiplication map $\Delta$. For this purpose, consider the following elements:
$$U_{ij}=\sum_ku_{ik}\otimes u_{kj}\quad, \quad V_{ij}=\sum_kv_{ik}\otimes v_{kj}$$

We know that the original comultiplication map $\Delta$, that of the Hopf algebra $F(GL_N^+)$, comes via the following formulae, on the standard generators:
$$\Delta(u_{ij})=U_{ij}\quad,\quad\Delta(v_{ij})=V_{ij}$$

Thus, in order to prove that $\Delta$ factorizes indeed, we must prove that we have:
$$u={\rm half\!-\!twisted}\implies U={\rm half\!-\!twisted}$$

(3) So, let us prove this. We have the following computation, with $w=u,v$ and $W=U,V$, and with the signs being those coming from Definition 7.11:
\begin{eqnarray*}
W_{ia}W_{jb}W_{kc}
&=&\sum_{xyz}w_{ix}w_{jy}w_{kz}\otimes w_{xa}w_{yb}w_{zc}\\
&=&\sum_{xyz}\pm w_{kz}w_{jy}w_{ix}\otimes\pm w_{zc}w_{yb}w_{xa}\\
&=&\pm W_{kc}W_{jb}W_{ia}
\end{eqnarray*}

We must show that, when examining the precise two $\pm$ signs in the middle formula, their product produces the correct $\pm$ sign at the end. But the point is that both these signs depend only on $s=span(x,y,z)$, and for $s=1,2,3$ respectively, we have:

\medskip

-- For a $(3,1)$ span we obtain $+-$, $+-$, $-+$, so a product $-$ as needed.

\smallskip

-- For a $(2,1)$ span we obtain $++$, $++$, $--$, so a product $+$ as needed.

\smallskip

-- For a $(3,3)$ span we obtain $--$, $--$, $++$, so a product $+$ as needed.

\smallskip

-- For a $(3,2)$ span we obtain $+-$, $+-$, $-+$, so a product $-$ as needed.

\smallskip

-- For a $(2,2)$ span we obtain $++$, $++$, $--$, so a product $+$ as needed.

\medskip

Together with the fact that our problem is invariant under $(r,c)\to(c,r)$, and with the fact that for a $(1,1)$ span there is nothing to prove, this finishes the proof for $\Delta$.

\medskip

(4) Regarding now the counit $\varepsilon$ and the antipode $S$, things are clear here, by using the same method, and with no computations needed, the formulae to be satisfied being trivially satisfied. We conclude that $GL_N^{*\prime}$ is indeed a quantum group, as desired.
\end{proof}

As yet another question, we can now talk about quantum groups related to the special linear group $SL_N$. Indeed, recall that $SL_N\subset GL_N$ is defined by imposing the condition $\det U=1$, which amounts in imposing the following condition to the coordinates:
$$\sum_{\sigma\in S_N}\varepsilon(\sigma)u_{1\sigma(1)}\ldots u_{N\sigma(N)}=1$$

As explained in chapter 5, this relation does not liberate, or even half-liberate. However, and here comes our point, this relation has a twisted analogue, namely:
$$\sum_{\sigma\in S_N}u_{1\sigma(1)}\ldots u_{N\sigma(N)}=1$$

To be more precise, at the quantum group level we have the folowing construction, making appear, eventually, $SL_N$ in our quantum group discussion:

\begin{theorem}
We have a quantum subgroup, as follows,
$$SL_N'\subset GL_N'$$
obtained by imposing to $u,v$ the twisted determinant one condition.
\end{theorem}

\begin{proof}
This is something quite tricky, as follows:

\medskip

(1) Let us first see how this construction goes at $N=2$. Consider the fundamental corepresentation of the algebra $F(GL_N')$, denoted as follows:
$$u=\begin{pmatrix}u_{11}&u_{12}\\ u_{21}&u_{22}\end{pmatrix}$$

The twisted determinant one condition, on this matrix, reads:
$$u_{11}u_{22}+u_{21}u_{12}=1$$

In order to construct the comultiplication of the algebra $F(SL_N')$, let us set:
$$U_{ij}=\sum_ku_{ik}\otimes u_{kj}=u_{i1}\otimes u_{1j}+u_{i2}\otimes u_{2j}$$

We must prove that these elements satisfy the twisted determinant one condition, and the verification here, using the commutation and anticommutation relations for the variables $u_{ij}$, and then the twisted determinant one condition for $u$, goes as follows:
\begin{eqnarray*}
&&U_{11}U_{22}+U_{21}U_{12}\\
&=&(u_{11}\otimes u_{11}+u_{12}\otimes u_{21})(u_{21}\otimes u_{12}+u_{22}\otimes u_{22})\\
&&+(u_{21}\otimes u_{11}+u_{22}\otimes u_{21})(u_{11}\otimes u_{12}+u_{12}\otimes u_{22})\\
&=&u_{11}u_{21}\otimes u_{11}u_{12}+u_{11}u_{22}\otimes u_{11}u_{22}
+u_{12}u_{21}\otimes u_{21}u_{12}+u_{12}u_{22}\otimes u_{21}u_{22}\\
&&+u_{21}u_{11}\otimes u_{11}u_{12}+u_{21}u_{12}\otimes u_{11}u_{22}
+u_{22}u_{11}\otimes u_{21}u_{12}+u_{22}u_{12}\otimes u_{21}u_{22}\\
&=&u_{11}u_{21}\otimes u_{11}u_{12}+u_{11}u_{22}\otimes u_{11}u_{22}
+u_{21}u_{12}\otimes u_{21}u_{12}+u_{12}u_{22}\otimes u_{21}u_{22}\\
&&-u_{11}u_{21}\otimes u_{11}u_{12}+u_{21}u_{12}\otimes u_{11}u_{22}
+u_{11}u_{22}\otimes u_{21}u_{12}-u_{12}u_{22}\otimes u_{21}u_{22}\\
&=&u_{11}u_{22}\otimes u_{11}u_{22}+u_{21}u_{12}\otimes u_{21}u_{12}
+u_{21}u_{12}\otimes u_{11}u_{22}+u_{11}u_{22}\otimes u_{21}u_{12}\\
&=&(u_{11}u_{22}+u_{21}u_{12})\otimes(u_{11}u_{22}+u_{21}u_{12})\\
&=&1\otimes1\\
&=&1
\end{eqnarray*}

Thus we can define indeed $\Delta$. Regarding $\varepsilon$, the verification here is trivial. Finally, regarding $S$, this follows from the fact that we imposed our relations on both $u,v$. 

\medskip

(2) At $N=3$ the verification is similar, but considerably more complicated, with the twisted determinant one condition for the variables $u_{ij}$ involving 6 terms, and with this same relation for the variables $U_{ij}$ involving $6\times 27=162$ terms. We will leave this as an instructive exercise, and we will come back to $N=3$ later in this chapter, with a simple, ad-hoc proof for the fact that $SO_3'=O_3^+\cap SL_3'$ is indeed a quantum group. 

\medskip

(3) In the general case now, where $N\in\mathbb N$ is arbitrary, the proof is similar. Consider as usual the following variables, that we must study in order to construct $\Delta$:
$$U_{ij}=\sum_ku_{ik}\otimes u_{kj}$$

We must prove that these variables $U_{ij}$ satisfy the twisted determinant one condition, provided that the original variables $u_{ij}$ do. For this purpose, observe that we have:
\begin{eqnarray*}
\sum_{\sigma\in S_N}U_{1\sigma(1)}\ldots U_{N\sigma(N)}
&=&\sum_{\sigma\in S_N}\sum_{k_1,\ldots,k_N}(u_{1k_1}\otimes u_{k_1\sigma(1)})\ldots(u_{Nk_N}\otimes u_{k_N\sigma(N)})\\
&=&\sum_{\sigma\in S_N}\sum_{k_1,\ldots,k_N}u_{1k_1}\ldots u_{Nk_N}\otimes u_{k_1\sigma(1)}\ldots u_{k_N\sigma(N)}\\
&=&\sum_{k_1,\ldots,k_N}u_{1k_1}\ldots u_{Nk_N}\otimes\sum_{\sigma\in S_N}u_{k_1\sigma(1)}\ldots u_{k_N\sigma(N)}
\end{eqnarray*}

The point now is that, due to the commutation and anticommutation relations for the variables $u_{ij}$, from Definition 7.1, the sum on the right vanishes, unless the indices $k_1,\ldots,k_N$ are all distinct. Thus, we can finish the computation in the following way, again by using the commutation and anticommutation relations for the variables $u_{ij}$, and then the twisted determinant one condition for these variables, at the very end:
\begin{eqnarray*}
\sum_{\sigma\in S_N}U_{1\sigma(1)}\ldots U_{N\sigma(N)}
&=&\sum_{k_1,\ldots,k_N\ distinct}u_{1k_1}\ldots u_{Nk_N}\otimes\sum_{\sigma\in S_N}u_{k_1\sigma(1)}\ldots u_{k_N\sigma(N)}\\
&=&\sum_{\tau\in S_N}u_{1\tau(1)}\ldots u_{N\tau(N)}\otimes\sum_{\sigma\in S_N}u_{\tau(1)\sigma(1)}\ldots u_{\tau(N)\sigma(N)}\\
&=&\sum_{\tau\in S_N}u_{1\tau(1)}\ldots u_{N\tau(N)}\otimes\sum_{\rho\in S_N}u_{1\rho(1)}\ldots u_{N\rho(N)}\\
&=&1\otimes1\\
&=&1
\end{eqnarray*}

Thus we can define indeed $\Delta$. Regarding $\varepsilon$, the verification here is trivial. Finally, regarding $S$, this follows from the fact that we imposed our relations on both $u,v$. 
\end{proof}

Summarizing, we have some basic twisting theory up and working.

\section*{7c. Twisted rotations}

Let us discuss now the twisted versions of the orthogonal groups $O_N$, and more generally of the groups $O_J$ with $J^t=\pm J$, considered in chapter 5. We first have:

\begin{theorem}
We can construct a twisted orthogonal group by setting
$$O_N'=O_N^+\cap GL_N'$$
with the intersection being computed inside $GL_N^+$.
\end{theorem}

\begin{proof}
The result as stated is of course something trivial, simply obtained by applying the intersection operation, that we constructed in chapter 4. However, at the level of the examples, at small values of $N$, things are quite tricky. We will be back to this. 
\end{proof}

We have as well a twisted special orthogonal group, constructed as follows:

\begin{theorem}
We can construct a twisted special orthogonal group by setting
$$SO_N'=O_N'\cap SL_N'$$
with the intersection being computed inside $GL_N'$.
\end{theorem}

\begin{proof}
As before, the result as stated is something trivial, obtained by applying the intersection operation, that we constructed in chapter 4. The interesting things happen at the level of examples, at small values of $N$, and we will be back to this. 
\end{proof}

As a more technical question now, we would like to perform the above two operations, $O\to O'\to SO'$, to the general quantum groups $O_J$ considered in chapter 5. So, let us recall from there that, barring some examples which are not very interesting, with the fundamental corepresentation $u$ being trivially reducible, the quantum groups $O_J^+\subset GL_N^+$, with $J^t=\pm J$, appear via the following relation, with $\bar{u}=(t\otimes id)u^{-1}$: 
$$u=J\bar{u}J^{-1}$$

Since we have by definition $O_J^+\subset GL_N^+$, we can use the above twisiting and intersection technology, as to successively get twisted orthogonal quantum groups, as follows:

\begin{theorem}
We can construct a twisted orthogonal group by setting
$$O_J'=O_J^+\cap GL_N'$$
inside $GL_N^+$, and then a twisted special orthogonal group by setting
$$SO_J'=O_J'\cap SL_N'$$
inside $GL_N'$, for any matrix $J\in GL_N$ satisfying $J^t=\pm J$.
\end{theorem}

\begin{proof}
As before, the result as stated is something trivial, obtained by applying the intersection operation, that we constructed in chapter 4. However, at the level of the examples, at small values of $N$, things are quite tricky. We will be back to this. 
\end{proof}

Going now towards some applications, the idea is that, with the above twisting technology, we can say many interesting things at $N=4$. Let us formulate:

\index{twisting}
\index{twisted determinant}
\index{anticommutation}

\begin{definition}
We let $SO_3'\subset O_3'$ be the subgroup coming from the relation
$$\sum_{\sigma\in S_3}u_{1\sigma(1)}u_{2\sigma(2)}u_{3\sigma(3)}=1$$
called twisted determinant one condition.
\end{definition}

Normally, when formulating such a definition, we should prove that $F(SO_3')$ is indeed an affine Hopf algebra. This is of course possible, directly, and actually we already did such verifications, in the beginning of this chapter, for some similar quantum groups.

\bigskip

However, the point is that this follows as well from the following result, from \cite{bb1}:

\index{Klein group}

\begin{theorem}
We have an isomorphism of compact quantum groups
$$S_4^+=SO_3'$$
given by the Fourier transform over the Klein group $K=\mathbb Z_2\times\mathbb Z_2$.
\end{theorem}

\begin{proof}
Consider the following matrix, coming from the action of $SO_3'$ on $F^4$:
$$u^+=\begin{pmatrix}1&0\\0&u\end{pmatrix}$$

We apply to this matrix the Fourier transform over the Klein group $K=\mathbb Z_2\times\mathbb Z_2$: 
$$v=
\frac{1}{4}
\begin{pmatrix}
1&1&1&1\\
1&-1&-1&1\\
1&-1&1&-1\\
1&1&-1&-1
\end{pmatrix}
\begin{pmatrix}
1&0&0&0\\
0&u_{11}&u_{12}&u_{13}\\
0&u_{21}&u_{22}&u_{23}\\
0&u_{31}&u_{32}&u_{33}
\end{pmatrix}
\begin{pmatrix}
1&1&1&1\\
1&-1&-1&1\\
1&-1&1&-1\\
1&1&-1&-1
\end{pmatrix}$$

By performing the matrix products, we are led to a formula of the following type, with each entry $v_{ij}$ being a certain linear combination of the entries $u_{ij}$:
$$v=\begin{pmatrix}
v_{11}&v_{12}&v_{13}&v_{14}\\
v_{21}&v_{22}&v_{23}&v_{24}\\
v_{31}&v_{32}&v_{33}&v_{34}\\
v_{41}&v_{42}&v_{43}&v_{44}
\end{pmatrix}$$

This matrix is then magic, and vice versa, so the Fourier transform over $K$ converts the relations in Definition 7.17 into the magic relations. But this gives the result.
\end{proof}

As a continuation of this, which is something more technical, and far-reaching, we have the following classification result, over $F=\mathbb C$, for the subgroups of $S_4^+=SO_3'$:

\index{ADE classification}

\begin{theorem}
The subgroups of $S_4^+=SO_3'$ are as follows:
\begin{enumerate}
\item Infinite quantum groups: $S_4^+$, $O_2'$, $\widehat{D}_\infty$.

\item Finite groups: $S_4$, and its subgroups.

\item Finite group twists: $S_4'$, $A_5'$.

\item Series of twists: $D_{2n}'$ $(n\geq 3)$, $DC_{2n}'$ $(n\geq 2)$.

\item A group dual series: $\widehat{D}_n$, with $n\geq 3$.
\end{enumerate}
Moreover, these quantum groups are subject to an ADE classification result.
\end{theorem}

\begin{proof}
This is something quite technical, from \cite{bb1}, and we will be back to this in the next chapter, with some further explanations on the various twists appearing above, by using some more advanced twisting technology, to be developed at that time.
\end{proof}

Let us explore now some consequences of the above result. By restricting the attention to the transitive case, which is the most important one, we obtain the following result:

\begin{theorem}
The transitive quantum groups $G\subset S_N^+$ are as follows:
\begin{enumerate}
\item At $N=1,2,3$ we have $\{1\}$, $\mathbb Z_2$, $\mathbb Z_3$, $S_3$.

\item At $N=4$ we have $\mathbb Z_2\times\mathbb Z_2,\mathbb Z_4,D_4,A_4,S_4,O_2',S_4^+$ and $S_4'$, $A_5'$.
\end{enumerate}
\end{theorem}

\begin{proof}
This follows from the above result, the idea being as follows:

\medskip

(1) This follows from the fact that we have $S_N=S_N^+$ at $N\leq3$.

\medskip

(2) This follows from the ADE classification of the subgroups $G\subset S_4^+$, with all the twists appearing in the statement being standard twists.
\end{proof}

Inspired by this, let us study now the transitive subgroups $G\subset S_5^+$. This is something substantially more complicated, and we first have here the following result:

\begin{proposition}
The transitive subgroups $G\subset S_5$ are
$$\mathbb Z_5,D_5,GA_1(5),A_5,S_5$$
with $GA_1(5)=\mathbb Z_5\rtimes\mathbb Z_4$ being the general affine group of $\mathbb F_5$.
\end{proposition}

\begin{proof}
This is something well-known, and elementary. Observe that the group $GA_1(5)=\mathbb Z_5\rtimes\mathbb Z_4$, appearing as the general affine group of $\mathbb F_5$, is indeed transitive.
\end{proof}

In the quantum case now, the study is something quite difficult, but we have here the following result, again over $F=\mathbb C$, based on advanced subfactor theory:

\begin{theorem}
The set of principal graphs of transitive subgroups $G\subset S_5^+$ coincide with the set of principal graphs of the following subgroups,
$$\mathbb Z_5,D_5,GA_1(5),A_5,S_5,S_5^+$$
and this gives the full list of transitive subgroups $G\subset S_5^+$.
\end{theorem}

\begin{proof}
This is something quite advanced, the idea being that subfactor theory tells us that the principal graphs of the irreducible index $5$ subfactors are as follows:

\medskip

-- $A_\infty$, and a non-extremal perturbation of $A_\infty^{(1)}$.

\medskip

-- The McKay graphs of $\mathbb Z_5,D_5,GA_1(5),A_5,S_5$.

\medskip

-- The twists of the McKay graphs of $A_5,S_5$.

\medskip

Thus, we must take the this list, and exclude the graphs which cannot be realized by a transitive subgroup $G\subset S_5^+$. And here, we have 3 cases, as follows:

\medskip

-- The graph $A_\infty$ corresponds to $S_5^+$ itself. As for the perturbation of $A_\infty^{(1)}$, this dissapears, because our notion of transitivity requires the subfactor extremality.

\medskip

-- For the McKay graphs of the groups $\mathbb Z_5,D_5,GA_1(5),A_5,S_5$ there is nothing to be done, all these graphs being solutions to our problem.

\medskip

-- The possible twists of $A_5,S_5$, coming from the subfactor graphs above, cannot contain $S_5$, because their cardinalities are smaller or equal than $|S_5|=120$.
\end{proof}

Still with me, I hope, and do not worry, we will come back later to subfactors, in this book. Now the point is that, as an interesting consequence of the above, we have:

\begin{theorem}
The inclusion of quantum groups
$$S_N\subset S_N^+$$
is maximal at $N\leq5$, in the sense that there is no quantum group in between.
\end{theorem}

\begin{proof} 
This follows indeed from the various classification results above:

\medskip

(1) At $N=2,3$ this simply follows from $S_N=S_N^+$.

\medskip

(2) At $N=4$ this follows from Theorem 7.20.

\medskip

(3) At $N=5$ this follows from Theorem 7.22.
\end{proof}

In general, it is conjectured that the quantum group inclusion $S_N\subset S_N^+$ is maximal at any $N\in\mathbb N$. We will be back to this, later in this book, with further results about it.

\section*{7d. Reflection groups}

Time now for some further applications of our twisting technology developed above, in higher dimensions. We have, following \cite{bbc}, the following result:

\index{twisted orthogonal group}

\begin{theorem}
The quantum symmetry group of the $N$-hypercube is
$$G^+(\square_N)=O_N'$$
with the corresponding coaction map on the vertex set being given by
$$\alpha:F[\mathbb Z_2^N]\to F[\mathbb Z_2^N]\otimes F(O_N')\quad,\quad 
g_i\to\sum_jg_j\otimes u_{ji}$$
via the standard identification $\square_N=\widehat{\mathbb Z_2^N}$. In particular we have $G^+(\square)=O_2'$.
\end{theorem} 

\begin{proof}
This result is from \cite{bbc}, with its $N=2$ particular case, corresponding to the last assertion, going back to \cite{bi1}, the idea being as follows:

\medskip

(1) To start with, the result in the statement is clear at $N=2$, coming as a consequence of the various results about $S_4^+$ and its subgroups, from the previous section.

\medskip

(2) In the general case now, where $N\in\mathbb N$ is arbitrary, our first claim is that the $N$-hypercube $\square_N$ is the Cayley graph of the following group:
$$\mathbb Z_2^N=<\tau_1,\ldots ,\tau_N>$$

Indeed, in order to prove this, let us recall that the vertices of this latter Cayley graph are by definition the products of group elements of the following form:
$$g=\tau_1^{i_1}\ldots\tau_N^{i_N}$$

Now observe that the sequence of 0-1 exponents defining such a group element uniquely determines a point of $F^N$, which is a vertex of the cube. We conclude that the vertices of the Cayley graph are indeed the vertices of the cube. 

\medskip

(3) Now regarding the edges, in the Cayley graph these are drawn between elements $g,h$ having the property $g=h\tau_i$ for some $i$. In terms of coordinates, the operation $h\to h\tau_i$ means to switch the sign of the $i$-th coordinate, and to keep the other coordinates fixed. In other words, we get in this way the edges of the cube, as desired.

\medskip

(4) Our second claim is that, when identifying the vector space spanned by the vertices of $\square_N$ with the algebra $F[\mathbb Z_2^N]$, the eigenvectors and eigenvalues of $\square_N$ are given by:
$$v_{i_1\ldots i_N}=\sum_{j_1\ldots j_N} (-1)^{i_1j_1
+\ldots+i_Nj_N}\tau_1^{j_1}\ldots\tau_N^{j_N}$$
$$\lambda_{i_1\ldots i_N}=(-1)^{i_1}+\ldots +(-1)^{i_N}$$

Indeed, let us recall that the action of $d$ on the functions on vertices is given by the following formula, with $q-p$ standing for the fact that $q,p$ are joined by an edge:
$$df(p)=\sum_{q-p}f(q)$$

Now by identifying the vertices with the the elements of $\mathbb Z_2^N$, hence the functions on the vertices with the elements of the algebra
$F[\mathbb Z_2^N]$, we get the following formula:
$$dv=\tau_1v+\ldots +\tau_Nv$$

With $v_{i_1\ldots i_N}$ as above, we have the following computation, which proves our claim:
\begin{eqnarray*}
dv_{i_1\ldots i_N}
&=&\sum_s\tau_s\sum_{j_1\ldots j_N}(-1)^{i_1j_1+\ldots+i_Nj_N}\tau_1^{j_1}\ldots\tau_N^{j_N}\\
&=&\sum_s\sum_{j_1\ldots j_N}(-1)^{i_1j_1+\ldots+i_Nj_N}\tau_1^{j_1}\ldots\tau_s^{j_s+1}\ldots\tau_N^{j_N}\\
&=&\sum_s\sum_{j_1\ldots j_N}(-1)^{i_s}(-1)^{i_1j_1+\ldots+i_Nj_N}\tau_1^{j_1}\ldots\tau_s^{j_s}\ldots\tau_N^{j_N}\\
&=&\sum_s(-1)^{i_s}\sum_{j_1\ldots j_N}(-1)^{i_s}(-1)^{i_1j_1+\ldots+i_Nj_N}\tau_1^{j_1}\ldots\tau_N^{j_N}\\
&=&\lambda_{i_1\ldots i_N}v_{i_1\ldots i_N}
\end{eqnarray*}

(5) We prove now that the quantum group $O_N'$ acts on the hypercube $\square_N$. For this purpose, observe first that we have a morphism of algebras as follows:
$$\alpha:F[\mathbb Z_2^N]\to F[\mathbb Z_2^N]\otimes F(O_N')
\quad,\quad\tau_i\to\sum_j\tau_j\otimes u_{ji}$$

It is routine to check that for $i_1\neq i_2\neq\ldots\neq i_l$ we have:
$$\alpha(\tau_{i_1}\ldots\tau_{i_l})=\sum_{j_1\neq\ldots\neq j_l}\tau_{j_1}\ldots\tau_{j_l} \otimes u_{j_1i_1}\ldots u_{j_li_l}$$

In terms of eigenspaces $E_s$ of the adjacency matrix, this gives:
$$\alpha(E_s)\subset E_s\otimes F(O_N')$$

Thus $\alpha$ preserves the adjacency matrix of $\square_N$, so is a coaction on $\square_N$, as claimed.

\medskip

(6) Conversely now, consider the universal coaction on the cube:
$$\beta:F[\mathbb Z_2^N]\to F[\mathbb Z_2^N]\otimes F(G)\quad,\quad\tau_i\to\sum_j\tau_j\otimes u_{ji}$$

By applying $\beta$ to the relation $\tau_i\tau_j=\tau_j\tau_i$ we get $u^tu=1$, so the matrix $u=(u_{ij})$ is orthogonal. By applying $\beta$ to the relation $\tau_i^2=1$ we get:
$$1\otimes\sum_ku_{ki}^2+\sum_{k<l}\tau_k\tau_l\otimes(u_{ki}u_{li}+u_{li}u_{ki})=1\otimes 1$$

This gives $u_{ki}u_{li}=-u_{li}u_{ki}$ for $i\neq j$, $k\neq l$, and by using the antipode we get $u_{ik}u_{il}=-u_{il}u_{ik}$ for $k\neq l$. Also, by applying $\beta$ to $\tau_i\tau_j=\tau_j\tau_i$ with $i\neq j$ we get:
$$\sum_{k<l}\tau_k\tau_l\otimes(u_{ki}u_{lj}+u_{li}u_{kj})=\sum_{k<l}\tau_k\tau_l\otimes
(u_{kj}u_{li}+u_{lj}u_{ki})$$

It follows that for $i\neq j$ and $k\neq l$, we have:
$$u_{ki}u_{lj}+u_{li}u_{kj}=u_{kj}u_{li}+u_{lj}u_{ki}$$

In other words, we have $[u_{ki},u_{lj}]=[u_{kj},u_{li}]$. By using the antipode we get:
$$[u_{jl},u_{ik}]=[u_{il},u_{jk}]$$

Now by combining these relations we get:
$$[u_{il},u_{jk}]=[u_{ik},u_{jl}]=[u_{jk},u_{il}]=-[u_{il},u_{jk}]$$

Thus $[u_{il},u_{jk}]=0$, so the elements $u_{ij}$ satisfy the relations for $F(O_N')$, as desired. Thus, we are led to the conclusions in the statement.
\end{proof}

With the above discussed, which is quite intriguing, our purpose now, assuming some basic representation theory knowledge, will be to understand which exact representation of $O_N$ produces by twisting the magic representation of $O_N'$, appearing above. 

\bigskip

In order to solve this question, we will need the following standard fact:

\begin{proposition}
The Fourier transform over $\mathbb Z_2^N$ is the map
$$\varphi:F(\mathbb Z_2^N)\to F[\mathbb Z_2^N]\quad,\quad 
\delta_{g_1^{i_1}\ldots g_N^{i_N}}\to\frac{1}{2^N}\sum_{j_1\ldots j_N}(-1)^{<i,j>}g_1^{j_1}\ldots g_N^{j_N}$$
with the usual convention $<i,j>=\sum_ki_kj_k$, and its inverse is the map
$$\psi:F[\mathbb Z_2^N]\to F(\mathbb Z_2^N)\quad,\quad 
g_1^{i_1}\ldots g_N^{i_N}\to\sum_{j_1\ldots j_N}(-1)^{<i,j>}\delta_{g_1^{j_1}\ldots g_N^{j_N}}$$
with all the exponents being binary, $i_1,\ldots,i_N,j_1,\ldots,j_N\in\{0,1\}$.
\end{proposition}

\begin{proof}
Observe first that the group $\mathbb Z_2^N$ can be written as follows:
$$\mathbb Z_2^N=\left\{g_1^{i_1}\ldots g_N^{i_N}\Big|i_1,\ldots,i_N\in\{0,1\}\right\}$$

Thus both $\varphi,\psi$ are well-defined, and it is elementary to check that both are morphisms of algebras. Also, we have $\varphi\psi=\psi\varphi=id$, coming from the following standard formula:
$$\frac{1}{2^N}\sum_{j_1\ldots j_N}(-1)^{<i,j>}
=\prod_{k=1}^N\left(\frac{1}{2}\sum_{j_r}(-1)^{i_rj_r}\right)
=\delta_{i0}$$

Thus we have indeed a pair of inverse Fourier morphisms, as claimed.
\end{proof}

By using now these Fourier transforms, we obtain following result:

\begin{proposition}
The magic matrix for the embedding $O_N'\subset S_{2^N}^+$ is given by
$$w_{i_1\ldots i_N,k_1\ldots k_N}=\frac{1}{2^N}\sum_{j_1\ldots j_N}\sum_{b_1\ldots b_N}(-1)^{<i+k_b,j>}\left(\frac{1}{N}\right)^{\#(0\in j)}u_{1b_1}^{j_1}\ldots u_{Nb_N}^{j_N}$$
where $k_b=(k_{b_1},\ldots,k_{b_N})$, with respect to multi-indices $i,k\in\{0,1\}^N$ as above.
\end{proposition}

\begin{proof}
By composing the coaction map $\alpha$ from Theorem 7.24 with the above Fourier transform isomorphisms $\varphi,\psi$, we have a diagram as follows:
$$\xymatrix@R=20mm@C=25mm{
F[\mathbb Z_2^N]\ar[r]^\alpha&F[\mathbb Z_2^N]\otimes F(O_N') \ar[d]^{\psi\otimes id}\\
F(\mathbb Z_2^N)\ar[u]^\varphi\ar@.[r]^\beta&F(\mathbb Z_2^N)\otimes F(O_N')}$$

In order to compute the composition on the bottom $\beta$, we first recall from Theorem 7.24 that the coaction map $\alpha$ is defined by the following formula:
$$\alpha(g_b)=\sum_ag_a\otimes u_{ab}$$

Now by making products of such quantities, we obtain the following global formula for $\alpha$, valid for any exponents $i_1,\ldots,i_N\in\{1,\ldots,N\}$:
$$\Phi(g_1^{i_1}\ldots g_N^{i_N})=\left(\frac{1}{N}\right)^{\#(0\in i)}\sum_{b_1\ldots b_N}g_{b_1}^{i_1}\ldots g_{b_N}^{i_N}\otimes u_{1b_1}^{i_1}\ldots u_{Nb_N}^{i_N}$$

The term on the right can be put in ``standard form'' as follows:
$$g_{b_1}^{i_1}\ldots g_{b_N}^{i_N}=g_1^{\sum_{b_x=1}i_x}\ldots g_N^{\sum_{b_x}i_x}$$

We therefore obtain the following formula for the coaction map $\alpha$:
$$\alpha(g_1^{i_1}\ldots g_N^{i_N})=\left(\frac{1}{N}\right)^{\#(0\in i)}
\sum_{b_1\ldots b_N}g_1^{\sum_{b_x=1}i_x}\ldots g_N^{\sum_{b_x=N}i_x}\otimes u_{1b_1}^{i_1}\ldots u_{Nb_N}^{i_N}$$

Now by applying the Fourier transforms, we obtain the following formula:
\begin{eqnarray*}
&&\alpha(\delta_{g_1^{i_1}\ldots g_N^{i_N}})\\
&=&(\psi\otimes id)\Phi\left(\frac{1}{2^N}\sum_{j_1\ldots j_N}(-1)^{<i,j>}g_1^{j_1}\ldots g_N^{j_N}\right)\\
&=&\frac{1}{2^N}\sum_{j_1\ldots j_N}\sum_{b_1\ldots b_N}(-1)^{<i,j>}
\left(\frac{1}{N}\right)^{\#(0\in j)}
\psi\left( g_1^{\sum_{b_x=1}j_x}\ldots g_N^{\sum_{b_x=N}j_x}\right)\otimes u_{1b_1}^{j_1}\ldots u_{Nb_N}^{j_N}
\end{eqnarray*}

By using now the formula of $\psi$ from Proposition 7.25, we obtain:
\begin{eqnarray*}
\beta(\delta_{g_1^{i_1}\ldots g_N^{i_N}})
&=&\frac{1}{2^N}\sum_{j_1\ldots j_N}\sum_{b_1\ldots b_N}\sum_{k_1\ldots k_N}\left(\frac{1}{N}\right)^{\#(0\in j)}\\
&&(-1)^{<i,j>}(-1)^{<(\sum_{b_x=1}j_x,\ldots,\sum_{b_x=N}j_x),(k_1,\ldots,k_N)>}\\
&&\delta_{g_1^{k_1}\ldots g_N^{k_N}}\otimes u_{1b_1}^{j_1}\ldots u_{Nb_N}^{j_N}
\end{eqnarray*}

Now observe that, with the notation $k_b=(k_{b_1},\ldots,k_{b_N})$, we have:
$$\left<\left(\sum_{b_x=1}j_x,\ldots,\sum_{b_x=N}j_x\right),(k_1,\ldots,k_N)\right>=<j,k_b>$$

Thus, we obtain the following formula for our map $\beta$:
$$\beta(\delta_{g_1^{i_1}\ldots g_N^{i_N}})
=\frac{1}{2^N}\sum_{j_1\ldots j_N}\sum_{b_1\ldots b_N}\sum_{k_1\ldots k_N}(-1)^{<i+k_b,j>}\left(\frac{1}{N}\right)^{\#(0\in j)}\delta_{g_1^{k_1}\ldots g_N^{k_N}}\otimes u_{1b_1}^{j_1}\ldots u_{Nb_N}^{j_N}$$

But this gives the formula in the statement for the corresponding magic matrix, with respect to the basis $\{\delta_{g_1^{i_1}\ldots g_N^{i_N}}\}$ of the algebra $F(\mathbb Z_2^N)$, and we are done.
\end{proof}

We can now solve our original question, namely understanding where the magic representation of $O_N'$ really comes from, with the following final answer to it:

\index{antisymmetric representation}
\index{twisting}

\begin{theorem}
The magic representation of $O_N'$, coming from its action on the $N$-cube, corresponds to the antisymmetric representation of $O_N$, via twisting.
\end{theorem}

\begin{proof}
This follows from the formula of $w$ in Proposition 7.26, by computing the character, and then interpreting the result via twisting, as follows:

\medskip

(1) By applying the trace to the formula of $w$, we obtain:
$$\chi
=\sum_{j_1\ldots j_N}\sum_{b_1\ldots b_N}\left(\frac{1}{2^N}\sum_{i_1\ldots i_N}(-1)^{<i+i_b,j>}\right)\left(\frac{1}{N}\right)^{\#(0\in j)}u_{1b_1}^{j_1}\ldots u_{Nb_N}^{j_N}$$

(2) By computing the Fourier sum in the middle, we are led to the following formula, with binary indices $j_1,\ldots,j_N\in\{0,1\}$, and plain indices $b_1,\ldots,b_N\in\{1,\ldots,N\}$:
$$\chi=\sum_{j_1\ldots j_N}\sum_{b_1\ldots b_N}\left(\frac{1}{N}\right)^{\#(0\in j)}\delta_{j_1,\sum_{b_x=1}j_x}\ldots\delta_{j_N,\sum_{b_x=N}j_x}u_{1b_1}^{j_1}\ldots u_{Nb_N}^{j_N}$$

(3) With the notation $r=\#(1\in j)$ we obtain a decomposition as follows:
$$\chi=\sum_{r=0}^N\chi_r$$

To be more precise, the variables $\chi_r$ are as follows:
$$\chi_r=\frac{1}{N^{N-r}}\sum_{\#(1\in j)=r}\sum_{b_1\ldots b_N}\delta_{j_1,\sum_{b_x=1}j_x}\ldots\delta_{j_N,\sum_{b_x=N}j_x}u_{1b_1}^{j_1}\ldots u_{Nb_N}^{j_N}$$

(4) Consider now the set $A\subset\{1,\ldots,N\}$ given by:
$$A=\left\{a\Big|j_a=1\right\}$$

The binary multi-indices $j\in\{0,1\}^N$ satisfying $\#(1\in j)=r$ being in bijection with such subsets $A$, satisfying $|A|=r$, we can replace the sum over $j$ with a sum over such subsets $A$. We obtain a formula as follows, where $j$ is the index corresponding to $A$:
$$\chi_r=\frac{1}{N^{N-r}}\sum_{|A|=r}\sum_{b_1\ldots b_N}\delta_{j_1,\sum_{b_x=1}j_x}\ldots\delta_{j_N,\sum_{b_x=N}j_x}\prod_{a\in A}u_{ab_a}$$

(5) Let us identify $b$ with the corresponding function $b:\{1,\ldots,N\}\to\{1,\ldots,N\}$, via $b(a)=b_a$. Then for any $p\in\{1,\ldots,N\}$ we have:
$$\delta_{j_p,\sum_{b_x=p}j_x}=1
\iff |b^{-1}(p)\cap A|=\chi_A(p)\ ({\rm mod}\ 2)$$

We conclude that the multi-indices $b\in\{1,\ldots,N\}^N$ which effectively contribute to the sum are those coming from the functions satisfying $b<A$. Thus, we have:
$$\chi_r=\frac{1}{N^{N-r}}\sum_{|A|=r}\sum_{b<A}\prod_{a\in A}u_{ab_a}$$

(6) We can further split each $\chi_r$ over the sets $A\subset\{1,\ldots,N\}$ satisfying $|A|=r$. The point is that for each of these sets we have:
$$\frac{1}{N^{N-r}}\sum_{b<A}\prod_{a\in A}u_{ab_a}=\sum_{\sigma\in S_N^A}\prod_{a\in A}u_{a\sigma(a)}$$

Thus, the magic character of $O_N'$ splits as $\chi=\sum_{r=0}^N\chi_r$, the components being:
$$\chi_r=\sum_{|A|=r}\sum_{\sigma\in S_N^A}\prod_{a\in A}u_{a\sigma(a)}$$

(7) The twisting operation $O_N\to O_N'$ makes correspond the following products:
$$\varepsilon(\sigma)\prod_{a\in A}u_{a\sigma(a)}\to\prod_{a\in A}u_{a\sigma(a)}$$

Now by summing over sets $A$ and permutations $\sigma$, we conclude that the twisting operation $O_N\to O_N'$ makes correspond the following quantities:
$$\sum_{|A|=r}\sum_{\sigma\in S_N^A}\varepsilon(\sigma)\prod_{a\in A}u_{a\sigma(a)}\to \sum_{|A|=r}\sum_{\sigma\in S_N^A}\prod_{a\in A}u_{a\sigma(a)}$$

Thus, we are led to the conclusion in the statement.
\end{proof}

The above result is not the end of the story with the cube, because we still have the question of understanding how the hyperoctahedral group $H_N=\mathbb Z_2\wr S_N$ correctly liberates. But, as explained in chapter 6, the correct liberation is $H_N^+=\mathbb Z_2\wr_*S_N^+$.

\section*{7e. Exercises}

We had an advanced algebraic chapter here, and as exercises on this, we have:

\begin{exercise}
Clarify our twisting results for the groups $GL_N$.
\end{exercise}

\begin{exercise}
Work out twisting results in the parametric context.
\end{exercise}

\begin{exercise}
Check all the details for the twisted half-liberation operation.
\end{exercise}

\begin{exercise}
Work out the anticommutation twisting question for $Sp_N$.
\end{exercise}

\begin{exercise}
Find a more conceptual approach to anticommutation twisting.
\end{exercise}

\begin{exercise}
Find other graphs whose quantum symmetry groups are twists.
\end{exercise}

As bonus exercise, learn about braided Hopf algebras, and $R$-matrix deformation.

\chapter{Cocycle twists}

\section*{8a. Cocycle twists}

We have seen in the previous chapter how to construct many interesting twists, by using a simple idea, namely suitably replacing commutation by anticommutation:
$$ab=ba\quad\to\quad ab=\pm ba$$

Moreover, we have seen as well some variations of this idea, such as twisting the half-liberations, by suitably replacing half-commutation by half-anticommutation:
$$abc=cba\quad\to\quad abc=\pm cba$$

However, this is not the end of the story, with twisting. There are several other methods for doing this, notably by using cocycles, and we will explain this here, with the material that follows standing as a useful complement to what we learned in chapter 7.

\bigskip 

There are several ways of doing the cocycle twisting, at various levels of generality. In what concerns us, let us start with something fairly general, from \cite{doi}, as follows:

\index{2-cocycle}
\index{Sweedler notation}

\begin{definition}
A $2$-cocycle on a Hopf algebra $A$ is a linear map 
$$\sigma:A\otimes A\to F$$
which is convolution invertible and satisfies the following conditions,
$$\sum\sigma(x_1\otimes y_1)\sigma(x_2y_2\otimes z)=\sum\sigma(y_1\otimes z_1)\sigma(x\otimes y_2z_2)$$
$$\sigma(x\otimes 1)=\sigma(1\otimes x)=\varepsilon(x)$$
in Sweedler notation $\Delta(x)=\sum x_1\otimes x_2$.
\end{definition}

To be more precise, the first cocycle condition in the statement, which must hold for any $x,y,z\in A$, is taken with respect to the following Sweedler writing formulae:
$$\Delta(x)=\sum x_1\otimes x_2\quad,\quad 
\Delta(y)=\sum y_1\otimes y_2\quad,\quad 
\Delta(z)=\sum z_1\otimes z_2$$

Thus, the first cocycle condition is in fact something quite complicated. As for the second cocycle condition, that simply uses the counit map $\varepsilon$, as indicated.

\bigskip

There are many examples of such cocycles, the general idea being that, for the Hopf algebras of type $A=F(G)$ or $A=F[H]$, these cocycles come from group-theoretical cocycles, on the underlying groups $G,H$. More on this later, when discussing examples.

\bigskip

Getting now to the twisting operation, once a cocycle $\sigma:A\otimes A\to F$ as above is given, this is something quite subtle, which is best done in 3 steps, as follows:

\bigskip

(1) First, by suitably perturbing the multiplication of $A$, by using the cocycle $\sigma$, we can construct a certain left twisted algebra $_{\sigma}\!\,A$, which is an associative algebra.

\bigskip

(2) By perturbing in a similar way the multiplication of $A$, by using $\sigma^{-1}$, we can construct a certain right twisted algebra $A_{\sigma^{-1}}$, which is an associative algebra too.

\bigskip

(3) And the point now is that, when jointly performing the above two operations, the associative algebra $A^\sigma={_{\sigma}\!\,A}_{\sigma^{-1}}$ that we obtain is a Hopf algebra.

\bigskip

So, this was for the idea, and in practice now, let us see how this works. As a first task, we must construct the left twisted algebra $_{\sigma}\!\,A$, and this can be done as follows:

\begin{proposition}
Given $\sigma:A\otimes A\to F$, we can construct the left twisted algebra $_{\sigma}\!\,A$, which as a vector space is $_{\sigma}\!\,A=A$, and with product being defined by 
$$\{x\}\{y\}=\sum\sigma(x_1\otimes y_1)\{x_2y_2\}$$ 
with the convention that any $x\in A$ is denoted $\{x\}$, when viewed as an element 
of $_{\sigma}\!\,A$.
\end{proposition}

\begin{proof}
This is something quite straightforward, the idea being as follows:

\medskip

(1) Observe first that $_{\sigma}\!\,A$ is indeed an associative algebra, as shown by the following computation, based on the associativity condition satisfied by $\sigma$:
\begin{eqnarray*}
(\{x\}\{y\})\{z\}
&=&\sum\sigma(x_1\otimes y_1)\{x_2y_2\}\{z\}\\
&=&\sum\sigma(x_1\otimes y_1)\sigma((x_2y_2)_1\otimes z_1)\{(x_2y_2)_2z_2\}\\
&=&\sum\sigma(x_1\otimes(y_2z_2)_1)\sigma(y_1\otimes z_1)\{x_2(y_2z_2)_2\}\\
&=&\sum\sigma(y_1\otimes z_1)\{x\}\{y_2z_2\}\\
&=&\{x\}(\{y\}\{z\})
\end{eqnarray*}

(2) The unit stays the same when performing the operation $A\to{\ }_{\sigma}\!\,A$. Indeed, we have the following computation, based on the first unitality condition satisfied by $\sigma$:
\begin{eqnarray*}
\{x\}\{1\}
&=&\sum\sigma(x_1\otimes 1)\{x_2\}\\
&=&\sum\varepsilon(x_1)\{x_2\}\\
&=&\sum(\varepsilon\otimes id)\sum x_1\otimes x_2\\
&=&(\varepsilon\otimes id)\Delta(x)\\
&=&x
\end{eqnarray*}

(3) Similarly, still regarding the unit, we have as well the following computation, based this time on the second unitality condition satisfied by $\sigma$:
\begin{eqnarray*}
\{1\}\{x\}
&=&\sum\sigma(1\otimes x_1)\{x_2\}\\
&=&\sum\varepsilon(x_1)\{x_2\}\\
&=&\sum(\varepsilon\otimes id)\sum x_1\otimes x_2\\
&=&(\varepsilon\otimes id)\Delta(x)\\
&=&x
\end{eqnarray*}

(4) Thus, we have an associative algebra, as claimed, and with the cocycle conditions in Definition 8.1 being in fact exactly those which are needed, for having this.
\end{proof}

Next, and still following our three-step twisting plan explained above, we must construct the right twisted algebra $A_{\sigma^{-1}}$. This can be done as follows:

\begin{proposition}
Given $\sigma:A\otimes A\to F$, we can construct the right twisted algebra $A_{\sigma^{-1}}$, which as a vector space is $A_{\sigma^{-1}}=A$, and with product being defined by 
$$<x ><y>=\sum\sigma^{-1}(x_2\otimes y_2)<x_1y_1>$$
where $\sigma^{-1}$ is the convolution inverse of $\sigma$, with the convention that any $x\in A$ is denoted $<x>$, when viewed as an element of $A_{\sigma^{-1}}$.
\end{proposition}

\begin{proof}
As before with the construction of the left twisted algebra $_{\sigma}\!\,A$, there are some elementary verifications to be performed here, the idea being as follows:

\medskip

(1) To start with, the present result can be proved exactly as we did before for Proposition 8.2, by performing some elementary computations, using the cocycle axioms and the Sweedler notation, that we will leave here as an instructive exercise. 

\medskip

(2) Alternatively, at a more conceptual level, let us first review Definition 8.1 and Proposition 8.2, with various right and left considerations in mind. We recall from Definition 8.1 that the axioms there for the $2$-cocycles on Hopf algebras were as follows:
$$\sum\sigma(x_1\otimes y_1)\sigma(x_2y_2\otimes z)=\sum\sigma(y_1\otimes z_1)\sigma(x\otimes y_2z_2)$$
$$\sigma(x\otimes 1)=\sigma(1\otimes x)=\varepsilon(x)$$

The point now is that, technically speaking, these are in fact the axioms for the left $2$-cocycles. We have as well the notion of right $2$-cocycle, which must satisfy:
$$\sum\tau(x_1y_1\otimes z)\tau(x_1\otimes y_2)=\sum\tau(x\otimes y_1z_1)\tau(y_2\otimes z_2)$$
$$\tau(x\otimes 1)=\tau(1\otimes x)=\varepsilon(x)$$

In what follows we will only need the notion of left 2-cocycle, that from Definition 8.1, and this is why we simply called there 2-cocycles, the objects axiomatized there.

\medskip

(3) However, and here comes the point, in the context of the present result, the distinction between left and right 2-cocycles can be something useful, avoiding us some of the computations evoked in (1). Indeed, it follows from definitions that if $\sigma:A\otimes A\to F$ is a left 2-cocycle, then its convolution inverse $\sigma^{-1}:A\otimes A\to F$ is a right 2-cocycle. Thus, this convolution inverse is subject to the following conditions:
$$\sum\sigma^{-1}(x_1y_1\otimes z)\sigma^{-1}(x_1\otimes y_2)=\sum\sigma^{-1}(x\otimes y_1z_1)\sigma^{-1}(y_2\otimes z_2)$$
$$\sigma^{-1}(x\otimes 1)=\sigma^{-1}(1\otimes x)=\varepsilon(x)$$

But with these latter conditions in hand, it is quite clear that the proof of Proposition 8.2 can be adapted, with minimal changes, to the present setting, and gives the result.

\medskip

(4) Finally, as yet another proof, we can see that $A_{\sigma^{-1}}$ is indeed an associative algebra, as a particular case of what we did before for $_{\sigma}\!\,A$, by suitably changing both $A$ and $\sigma$, and again we will leave the clarification of all this as an instructive exercise.
\end{proof}

And with this, good news, we can now perform the third and final step of our twisting program, with the final twisting statement, which is self-contained, formally making no reference to Proposition 8.2 and Proposition 8.3, being as follows:

\index{cocycle twist}
\index{Hopf algebra twist}
\index{category of comodules}
\index{twisted antipode}

\begin{theorem}
Given a Hopf algebra $A$, and a $2$-cocycle $\sigma:A\otimes A\to F$, we can construct a Hopf algebra $A^\sigma$, which as a coalgebra is $A^\sigma=A$, with product 
$$[x][y]=\sum\sigma(x_1\otimes y_1)\sigma^{-1}(x_3\otimes y_3)[x_2y_2]$$ 
and with the antipode being given by the following formula,
$$S^{\sigma}([x])=\sum\sigma(x_1\otimes S(x_2))\sigma^{-1}(S(x_4)\otimes x_5)[S(x_3)]$$
where an element $x\in A$ is denoted $[x]$, when viewed as an element of $A^{\sigma}$. \end{theorem}

\begin{proof}
This is a mixture of trivial and non-trivial facts, basically going back to the work of Doi \cite{doi}, the idea with this being as follows:

\medskip

(1) The first assertion is something standard, the point being that the Hopf algebra in the statement appears as follows, in terms of the algebra constructions before:
$$A^\sigma={_{\sigma}\!\,A}_{\sigma^{-1}}$$

(2) Indeed, recall from the above that the multiplication of $A_{\sigma^{-1}}$ is given by: 
$$<x ><y>=\sum\sigma^{-1}(x_2\otimes y_2)<x_1y_1>$$

On the other hand, the multiplication of an algebra of type $_{\sigma}\!\,B$ is given by:
$$\{x\}\{y\}=\sum\sigma(x_1\otimes y_1)\{x_2y_2\}$$ 

Now with $B=A_{\sigma^{-1}}$ as input, for this latter construction, we deduce that the multiplication of the algebra $_{\sigma}\!\,B={_{\sigma}\!\,A}_{\sigma^{-1}}$ is given by the following formula:
\begin{eqnarray*}
\{<x>\}\{<y>\}
&=&\sum\sigma(x_1\otimes y_1)\{<x>_2<y>_2\}\\
&=&\sum\sigma(x_1\otimes y_1)\{(<x><y>)_2\}\\
&=&\sum\sigma(x_1\otimes y_1)\sigma^{-1}(x_3\otimes y_3)\{<x_2y_2>\}
\end{eqnarray*}

But this is exactly the formula in the statement, with the convention $[x]=<\{x\}>$, which means that $x\in A$ is denoted $[x]$, when viewed as an element of $A^{\sigma}$.

\medskip

(3) Summarizing, we have proved that $A^\sigma$ with the multiplication in the statement is indeed an associative algebra, with this coming from $A^\sigma={_{\sigma}\!\,A}_{\sigma^{-1}}$. 

\medskip

(4) Next, as indicated in the statement, what we have is in fact a bialgebra, with the usual comultiplication map $\Delta$ and counit map $\varepsilon$ of the original Hopf algebra $A$. Indeed, the fact that $\Delta$ and $\varepsilon$ remain morphisms of algebras, when changing the multiplication as above, is something standard, that we will leave here as an exercise.

\medskip

(5) Regarding now the remaining bialgebra axioms to be satisfied, let us have a look at the Hopf algebra axioms for $A$, concerning $\Delta$ and $\varepsilon$, namely:
$$(\Delta\otimes id)\Delta=(id\otimes \Delta)\Delta$$
$$(\varepsilon\otimes id)\Delta=(id\otimes\varepsilon)\Delta=id$$

Now since these axioms do not make reference to the multiplicative structure, they will remain unchanged, when perturbing the multiplicative structure, as we did.

\medskip

(6) Getting now to the final, and most technical point, the existence of the antipode $S$ requires some work, and this follows from the explicit formula for the antipode in the statement. To be more precise, we obtain in this way an antimultiplicative map $S:A^\sigma\to A^\sigma$, which satisfies the usual Hopf algebra axiom for the antipode, namely:
$$m(S\otimes id)\Delta=m(id\otimes S)\Delta=\varepsilon(.)1$$

Thus, we obtain the result. For more on all this, details and examples, and some further generalizations too, we refer to Doi \cite{doi}, and related work.
\end{proof}

As a conclusion to all this, we have some general twisting theory up and working, and we refer to the specialized Hopf algebra literature, for more on all this. We will see later that the above formalism is quite general, in particular with the anticommutation twists constructed in chapter 7 appearing in this way, for some suitable cocycles.

\section*{8b. Basic examples}

As an application of the above cocycle twisting methods, we can go back now to the isomorphism $S_4^+=SO_3'$, briefly explained in chapter 7, and the subsequent ADE classification result for the subgroups of this quantum group, similar to the ADE classification of the subgroups of $SO_3$, briefly explained there too, this time with more details.

\bigskip

Let us begin with the following definition, that we already met in chapter 7:

\index{twisting}
\index{twisted determinant}
\index{anticommutation}

\begin{definition}
We let $SO_3'\subset O_3'$ be the subgroup coming from the relation
$$\sum_{\sigma\in S_3}u_{1\sigma(1)}u_{2\sigma(2)}u_{3\sigma(3)}=1$$
called twisted determinant one condition.
\end{definition}

Normally, when formulating such a definition, we should prove that $F(SO_3')$ is indeed an affine Hopf algebra. This is of course possible, directly, and actually we already did such verifications, in the beginning of the chapter 7, for some similar quantum groups.

\bigskip

However, the point is that this follows as well from the following result, from \cite{bb1}:

\index{Klein group}

\begin{theorem}
We have an isomorphism of compact quantum groups
$$S_4^+=SO_3'$$
given by the Fourier transform over the Klein group $K=\mathbb Z_2\times\mathbb Z_2$.
\end{theorem}

\begin{proof}
This is something which is quite routine, the idea being as follows:

\medskip

(1) Consider the following matrix, coming from the action of $SO_3'$ on $F^4$:
$$u^+=\begin{pmatrix}1&0\\0&u\end{pmatrix}$$

We apply to this matrix the Fourier transform over the Klein group $K=\mathbb Z_2\times\mathbb Z_2$: 
$$v=
\frac{1}{4}
\begin{pmatrix}
1&1&1&1\\
1&-1&-1&1\\
1&-1&1&-1\\
1&1&-1&-1
\end{pmatrix}
\begin{pmatrix}
1&0&0&0\\
0&u_{11}&u_{12}&u_{13}\\
0&u_{21}&u_{22}&u_{23}\\
0&u_{31}&u_{32}&u_{33}
\end{pmatrix}
\begin{pmatrix}
1&1&1&1\\
1&-1&-1&1\\
1&-1&1&-1\\
1&1&-1&-1
\end{pmatrix}$$

By performing the matrix products, we are led to a formula of the following type, with each entry $v_{ij}$ being a certain linear combination of the entries $u_{ij}$:
$$v=\begin{pmatrix}
v_{11}&v_{12}&v_{13}&v_{14}\\
v_{21}&v_{22}&v_{23}&v_{24}\\
v_{31}&v_{32}&v_{33}&v_{34}\\
v_{41}&v_{42}&v_{43}&v_{44}
\end{pmatrix}$$

Our claim now, which will prove the result, is that this matrix is magic, and vice versa, so that the above Fourier transform over the Klein group $K=\mathbb Z_2\times\mathbb Z_2$ converts the relations in Definition 8.5 into the magic relations.

\medskip

(2) In order to prove our claim, let us begin with some preliminaries. For $i,j\in\{1,2,3\}$
with $i\neq j$, we let $<i,j>$ be the unique  element in $\{1,2,3\}$ such that:
$$\{i,j,<i,j>\}=\{1,2,3\}$$

Now let $i,j,k,l \in \{1,2,3\}$ with $i\neq j$ and $k\neq l$. Our claim is that we have:
$$u_{<i,j><k,l>}=u_{ik}u_{jl}+u_{jk}u_{il}$$  

Indeed, for $i,j\in\{1,2,3\}$, let $i_1,i_2,j_1,j_2$ be the elements of $\{1,2,3\}$ satisfying:
$$\{i,i_1,i_2\}=\{j,j_1,j_2\}=\{1,2,3\} $$

Then, the well-known formula for the antipode of the quantum group $SU_3'$ gives:
$$S(u_{ji})=u_{i_1j_1}u_{i_2j_2}+u_{i_2j_1}u_{i_1j_2}$$

But since the matrix $u=(u_{ij})$ is orthogonal, we have the following formula:
$$u_{ij}=S(u_{ij})$$

Thus, we have proved the above claim, regarding the value of $u_{<i,j><k,l>}$.  

\medskip

(3) Let us prove now that $SO_3'$ is the universal quantum group acting on $F^4$. Let $e_1, e_2, e_3, e_4$ be the canonical basis of $F^4$, and consider as well the following basis:
\begin{eqnarray*}
1&=&e_1+e_2+e_3+e_4\\
\varepsilon_1&=&e_1-e_2-e_3+e_4\\
\varepsilon_2&=&e_1-e_2+e_3-e_4\\
\varepsilon_3&=&e_1+e_2-e_3-e_4
\end{eqnarray*}

Observe that this latter basis is precisely the one obtained by using the Fourier transform of the Klein group $K=\mathbb Z_2 \times\mathbb Z_2$, mentioned in (1). We have the following formulae, regarding this latter basis, which define a presentation of the algebra $F^4$:
$$\varepsilon_i^2=1\quad,\quad\forall i$$
$$\varepsilon_i\varepsilon_j=\varepsilon_j \varepsilon_i=\varepsilon_{<i, j>}\quad,\quad \forall i \neq j$$

Consider now the linear map $\alpha :F^4 \to F^4 \otimes F(SO_3')$ defined as follows:
$$\alpha(1)=1\otimes 1\quad,\quad \alpha(\varepsilon_i)=\sum_j\varepsilon_j\otimes u_{ji}$$

Our claim is that $\alpha$ is a coaction. Indeed, it is clear that $\alpha$ is coassociative, so it remains to check that $\alpha$ is a morphism of algebras. First, we have:
\begin{eqnarray*}
\alpha(\varepsilon_i)^2
&=&\sum_{kl}\varepsilon_k\varepsilon_l\otimes u_{ki}u_{li}\\
&=&\sum_k\varepsilon_k^2\otimes u_{ki}^2+\sum_{k\neq l}\varepsilon_k\varepsilon_l\otimes u_{ki}u_{li}\\
&=&1\otimes\left(\sum_ku_{ki}^2\right)+\sum_{k<l}\varepsilon_k\varepsilon_l\otimes(u_{ki}u_{li} +u_{li}u_{ki})\\
&=&1\otimes 1\\
&=&\alpha(\varepsilon_i^2)
\end{eqnarray*}

Also, by using the formula found in (2), we have, for any $i\neq j$:
\begin{eqnarray*}
\alpha(\varepsilon_i)\alpha(\varepsilon_j)
&=&\sum_{kl}\varepsilon_k\varepsilon_l\otimes u_{ki}u_{lj}\\
&=&\sum_k\varepsilon_k^2\otimes u_{ki}u_{kj}+\sum_{k\neq l}\varepsilon_k\varepsilon_l\otimes
 u_{ki}u_{lj}\\
&=&1\otimes\left(\sum_ku_{ki}u_{kj}\right)+\sum_{k<l}\varepsilon_k\varepsilon_l\otimes(u_{ki}u_{lj}+u_{li}u_{kj})\\
&=&\sum_{k<l}\varepsilon_{<k,l>}\otimes u_{<k,l\><i,j>}\\
&=&\sum_k\varepsilon_k\otimes u_{k< i,j>}\\
&=&\alpha(\varepsilon_{<i,j>})
\end{eqnarray*}

We conclude that our map $\alpha$ constructed above is indeed a coaction.

\medskip

(4) Our claim now is that $\alpha$ is the universal coaction on $F^4$. In order to prove this, consider a Hopf algebra $A$ coacting on $F^4$, with coaction denoted as follows:
$$\beta :F^4 \to F^4 \otimes A$$

We set $\varepsilon_0 =1$ and we write, with $x_{j0}=\delta_{j0}$: 
$$\beta(\varepsilon_i)=\sum_j\varepsilon_j\otimes x_{ji}$$

Let $\phi :F^4 \to F$ be the map given by $\phi(e_i) =1/4$, so that $\phi(\varepsilon_i)=0$ if $i>0$. This linear map $\phi$ is $A$-colinear, and so the linear map $F^4 \otimes F^4 \to F$ given by $x\otimes y\to\phi(xy)$ is $A$-colinear too. It follows that the matrix $x = (x_{ij})$ is orthogonal, and so:
$$x_{01}^2+x_{02}^2 + x_{03}^2=0$$

Therefore $x_{0i} = \delta_{0i}$, the matrix $x' =(x_{ij})_{1\leq i,j\leq3}$ is orthogonal and:
$$\beta(\varepsilon_i) = \sum_{j=1}^3 \varepsilon_j \otimes x_{ji}\quad,\quad\forall i\in\{1,2,3\}$$
 
By using these formulae, we have the following computation:
\begin{eqnarray*}
\beta(1)
&=&\beta(\varepsilon_i^2)\\
&=&\beta(\varepsilon_i)^2\\
&=&\sum_{kl}\varepsilon_k\varepsilon_l\otimes x_{ki}x_{li}\\
&=&\sum_k\varepsilon_k^2\otimes x_{ki}^2+\sum_{k\neq l}\varepsilon_k\varepsilon_l\otimes
 x_{ki}x_{li}\\
&=&1\otimes\left(\sum_kx_{ki}^2\right)+\sum_{k<l}\varepsilon_k\varepsilon_l\otimes (x_{ki}x_{li} +x_{li}x_{ki})\\
&=&1\otimes 1+\sum_{k<l}\varepsilon_k\varepsilon_l\otimes(x_{ki}x_{li}+x_{li}x_{ki})
\end{eqnarray*}

But $\beta(1)=1\otimes1$, so we deduce that for $k\neq l$, we have the following formula:
$$x_{ki}x_{li}=-x_{li}x_{ki}$$

By using now the antipode $S$, which satisfies $S(x_{ij})=x_{ji}$ since the matrix $x'$ is orthogonal, we have as well the following formula:
$$x_{ik}x_{il}=-x_{il}x_{ik}$$

(5) Next, for any $i\neq j$, we have the following computation:
\begin{eqnarray*}
\alpha(\varepsilon_{<i,j>})
&=&\alpha(\varepsilon_i \varepsilon_j)\\
&=&\alpha(\varepsilon_i)\alpha(\varepsilon_j)\\
&=&\sum_{kl}\varepsilon_k\varepsilon_l\otimes x_{ki}x_{lj}\\
&=&\sum_k\varepsilon_k^2\otimes x_{ki}x_{kj}+\sum_{k \neq l}\varepsilon_k\varepsilon_l\otimes
 x_{ki}x_{lj}\\
&=&1\otimes\left(\sum_kx_{ki}x_{kj}\right)+\sum_{k<l}\varepsilon_k\varepsilon_l\otimes(x_{ki}x_{lj}+x_{li}x_{kj})\\
&=&\sum_{k<l}\varepsilon_{<k,l>}\otimes(x_{ki}x_{lj}+x_{li}x_{kj})
\end{eqnarray*}

We conclude from this computation that we have the following formula:
$$x_{<k,l><i,j>}=x_{ki}x_{lj}+x_{li}x_{kj}$$

Similarly, since $\varepsilon_i\varepsilon_j=\varepsilon_j\varepsilon_i$, we have as well the following formula:
$$x_{<k,l><i,j>}=x_{kj}x_{li}+x_{lj}x_{ki}$$
 
By combining now these two relations that we found, we get, for $i\neq j$ and $k\neq l$:
$$[x_{ki},x_{lj}]=[x_{kj},x_{li}]$$

Now by using the antipode, we have as well the following formula:
$$[x_{ki},x_{lj}]=[x_{li},x_{kj}]$$

Summarizing, we have established the following equalities:
$$[x_{ki},x_{lj}]=[x_{li},x_{kj}]=[x_{kj},x_{li}]=-[x_{li},x_{kj}]$$

But this shows that we have the following equality:
$$x_{li}x_{kj} = x_{kj}x_{li}$$

(6) We can now finish. Indeed, we have the following computation, in relation with the quantum determinant one condition, from Definition 8.5:
\begin{eqnarray*}
&&x_{11}x_{22}x_{33}+x_{11}x_{23}x_{32}+x_{12}x_{21}x_{33}+x_{12}x_{23}x_{31}+x_{13}x_{22}x_{31}+x_{13}x_{23}x_{31}\\
&=&x_{11}(x_{22}x_{33}+x_{23}x_{32})+x_{12}(x_{21}x_{33}+x_{23}x_{31})+x_{13}(x_{22}x_{31} +x_{23}x_{31})\\
&=&x_{11}^2+x_{12}^2+x_{13}^2\\
&=&1
\end{eqnarray*}

Thus the quantum determinant one condition is satisfied, so we obtain a Hopf algebra morphism 
$F(SO_3') \to A$ commuting with the respective coactions, and we are done.
\end{proof}

We will prove now that the quantum group $SO_3'$ appears from the group $SO_3$, via cocycle twisting. Let $A$ be the algebra of representative functions on $SO_3$, with the canonical coordinate functions denoted  $x_{ij}$. Consider as well the Klein group:
$$K=\mathbb Z_2\times\mathbb Z_2=\left<t_1,t_2\ \Big|\ t_1^2=t_2^2 =1,\ t_1t_2=t_2t_1\right>$$

We set $t_3 = t_1t_2$. Now observe that the restriction of the functions on $SO_3$ to the diagonal subgroup of $SO_3$, which is $K$, gives a Hopf algebra surjection, as follows:
$$\pi_d:A\to F[K]\quad,\quad x_{ij}\to\delta_{ij}t_i$$

In view of this, let $\sigma:K\times K\to F^*$ be the unique bimultiplicative map such that:
$$\sigma(t_i,t_j)=\begin{cases}
-1&{\rm if} \ i\leq j\\
1&{\rm otherwise}
\end{cases}$$

Then $\sigma$ is a usual group 2-cocycle, and its unique linear extension to the group algebra
$F[K]\otimes F[K]$ is a 2-cocycle in the previous sense, still denoted $\sigma$. We get in this way a cocycle on $A$, given by the following formula, and satisfying $\sigma_d=\sigma_d^{-1}$:
$$\sigma_d=\sigma\circ(\pi_d \otimes \pi_d)$$

We can now formulate our twisting result regarding $SO_3'$, as follows:

\begin{theorem}
We have an isomorphism of Hopf algebras
$$F(SO_3')=A^{\sigma_d}$$
and so $SO_3'$ appears as a cocycle twist of $SO_3$.
\end{theorem}

\begin{proof} 
With the notations above, we have the following computation:
\begin{eqnarray*}
[x_{ij}][x_{kl}]
&=&\sum_{p,q,r,s} \sigma_d(x_{ip}\otimes x_{kq})\sigma_d^{-1}(x_{rj}\otimes x_{sl})[x_{pr}x_{qs}]\\
&=&\sigma(t_i,t_k)\sigma(t_j,t_l)[x_{ij}x_{kl}] 
\end{eqnarray*}
 
More generally, we have in fact the following formula, valid for any $r\in\mathbb N$:
 $$[x_{i_1j_1}] [x_{i_2j_2}] \ldots [x_{i_rj_r}] = 
 \left(\prod_{p<q}\sigma(t_{i_p},t_{i_q})\right)
 \left(\prod_{p<q}\sigma(t_{j_p},t_{j_q})\right)
 [x_{i_1j_1} x_{i_2j_2} \ldots x_{i_rj_r}]$$
 
Thus, we see that we have a Hopf algebra morphism, as follows:
 $$F(SO_3')\to A^{\sigma_d}\quad,\quad u_{ij}\to[x_{ij}]$$

This morphism is then clearly surjective, and a standard representation theory argument, for which we refer to \cite{bb1} and the related literature, including \cite{bi4}, shows that this morphism must be in fact an isomorphism. Thus, the result follows.
\end{proof}

As a continuation of this, which is something more technical, we have the following classification result, over $F=\mathbb C$, for the subgroups of $S_4^+=SO_3'$:

\index{ADE classification}

\begin{theorem}
The subgroups of $S_4^+=SO_3'$ are as follows:
\begin{enumerate}
\item Infinite quantum groups: $S_4^+$, $O_2'$, $\widehat{D}_\infty$.

\item Finite groups: $S_4$, and its subgroups.

\item Finite group twists: $S_4'$, $A_5'$.

\item Series of twists: $D_{2n}'$ $(n\geq 3)$, $DC_{2n}'$ $(n\geq 2)$.

\item A group dual series: $\widehat{D}_n$, with $n\geq 3$.
\end{enumerate}
Moreover, these quantum groups are subject to an ADE classification result.
\end{theorem}

\begin{proof}
This is something quite technical, based on the twisting technology that we developed above, the idea with all this being as follows:

\medskip

(1) Regarding the statement, the idea is that, with prime standing for twists, which all unique, and double prime denotes pseudo-twists, the classification is as follows:

\medskip

-- A case: $\mathbb Z_1$, $\mathbb Z_2$, $\mathbb Z_3$, $K$, $\widehat{D}_n$ ($n=2,3,\ldots,\infty$), $S_4^+$.

\medskip

-- D case: $\mathbb Z_4$, $D_{2n}'$, $D_{2n}''$ ($n=2,3,\ldots$), $H_2^+$, $D_1$, $S_3$.

\medskip

-- E case: $A_4$, $S_4$, $S_4'$, $A_5'$.

\medskip

(2) There are many comments to be made here, regarding conventions, as follows:

\medskip

-- To start with, the 2-element group $\mathbb Z_2=\{1,\tau\}$ can act in 2 ways on 4 points: either with the transposition $\tau$ acting without fixed point, and we use here the notation $\mathbb Z_2$, or with $\tau$ acting with 2 fixed points, and we use here the notation $D_1$.

\medskip

-- Similarly, the Klein group $K=\mathbb Z_2\times\mathbb Z_2$ can act in 2 ways on 4 points: either with 2 non-trivial elements having 2 fixed points each, and we use here the notation $K$, or with all non-trivial elements having no fixed points, and we use here the notation $D_2=\widehat{D}_2$.

\medskip

-- We have $D_4'=D_4$, and $D_4''=G_0$, the Kac-Paljutkin quantum group. Besides being a pseudo-twist of $D_{2n}$, the quantum group $D_{2n}''$ with $n\geq2$ is known to be as well a pseudo-twist of the dicyclic, or binary cyclic group $DC_{2n}$.

\medskip

-- Also, the definition of $D_{2n}'$, $D_{2n}''$ can be extended at $n=1,\infty$, and we formally have $D_2'=D_2''=K$, and $D_\infty'=D_\infty''=H_2^+$, but these conventions are not very useful. Finally, the groups $D_1,S_3$ are a bit special in the D case of the classification.

\medskip

(3) So, this was for the statement of the theorem, and with the various twists involved being constructed by using the general cocycle twisting procedure from the beginning of this chapter. In what regards now the proof, this is something quite technical, making a heavier use of the twisting technology, and we refer here to \cite{bb1} and related papers.
\end{proof}

We will be back to all this, which remains something non-trivial, later in this book.

\section*{8c. Group algebras} 

As a continuation of the above, let us develop some more applications. Following \cite{bbs}, we will prove in what follows that we have a general twisting result, as follows: 
$$PO_n^+=(S_{n^2}^+)'$$

In order to explain this material, which is quite technical, requiring good algebraic knowledge, let us begin with some generalities. We first have:

\begin{proposition}
Given a finite group $G$, the algebra $F(S_{\widehat{G}}^+)$ is isomorphic to the abstract algebra presented by generators $x_{gh}$ with $g,h\in G$, with the following relations:
$$x_{1g}=x_{g1}=\delta_{1g}\quad,\quad
x_{s,gh}=\sum_{t\in G}x_{st^{-1},g}x_{th}\quad,\quad 
x_{gh,s}=\sum_{t\in G}x_{gt^{-1}}x_{h,ts}$$
The comultiplication, counit and antipode are given by the formulae
$$\Delta(x_{gh})=\sum_{s\in G}x_{gs}\otimes x_{sh}\quad,\quad 
\varepsilon(x_{gh})=\delta_{gh}\quad,\quad
S(x_{gh})=x_{h^{-1}g^{-1}}$$
on the standard generators $x_{gh}$.
\end{proposition}

\begin{proof}
This follows indeed from a direct verification, based on the coaction axioms that must be satisfied by the algebra $F(S_{\widehat{G}}^+)$, that we will leave here as an exercise.
\end{proof}

Let us discuss now the twisted version of the above result. Consider a 2-cocycle on $G$, which is by definition a map $\sigma:G\times G\to F^*$ satisfying the following conditions:
$$\sigma_{gh,s}\sigma_{gh}=\sigma_{g,hs}\sigma_{hs}\quad,\quad 
\sigma_{g1}=\sigma_{1g}=1$$

Given such a cocycle, we can construct the associated twisted group algebra $F(\widehat{G}_\sigma)$, as being the vector space $F(\widehat{G})=F[G]$, with product $e_ge_h=\sigma_{gh}e_{gh}$. We have:

\begin{proposition}
The algebra $F(S_{\widehat{G}_\sigma}^+)$ is isomorphic to the abstract algebra presented by generators $x_{gh}$ with $g,h\in G$, with the relations $x_{1g}=x_{g1}=\delta_{1g}$ and:
$$\sigma_{gh}x_{s,gh}=\sum_{t\in G}\sigma_{st^{-1},t}x_{st^{-1},g}x_{th}\quad,\quad
\sigma_{gh}^{-1}x_{gh,s}=\sum_{t\in G}\sigma_{t^{-1},ts}^{-1}x_{gt^{-1}}x_{h,ts}$$
The comultiplication, counit and antipode are given by the formulae
$$\Delta(x_{gh})=\sum_{s\in G}x_{gs}\otimes x_{sh}\quad,\quad 
\varepsilon(x_{gh})=\delta_{gh}\quad,\quad
S(x_{gh})=\sigma_{h^{-1}h}\sigma_{g^{-1}g}^{-1}x_{h^{-1}g^{-1}}$$
on the standard generators $x_{gh}$.
\end{proposition}

\begin{proof}
Once again, this follows from a direct verification, which is quite similar to the verifications did before, and that we will leave here as an exercise.
\end{proof}

\index{twisting}
\index{cocycle twisting}

In what follows, we will prove that the quantum groups $S_{\widehat{G}}^+$ and $S_{\widehat{G}_\sigma}^+$ are related by a cocycle twisting operation. Let $A$ be a Hopf algebra. We recall that a left 2-cocycle is a convolution invertible linear map $\sigma:A\otimes A\to F$ satisfying the following conditions:
$$\sum\sigma(x_1\otimes y_1)\sigma(x_2y_2\otimes z)=\sum\sigma(y_1\otimes z_1)\sigma(x\otimes y_2z_2)$$ 
$$\sigma(x\otimes1)=\sigma(1\otimes x)=\varepsilon(x)$$

As explained before, $\sigma$ is a left 2-cocycle precisely when $\sigma^{-1}$, the convolution inverse of $\sigma$, is a right 2-cocycle, in the sense that we have:
$$\sum\sigma^{-1}(x_1y_1\otimes z)\sigma^{-1}(x_1\otimes y_2)=\sum\sigma^{-1}(x\otimes y_1z_1)\sigma^{-1}(y_2\otimes z_2)$$
$$\sigma^{-1}(x\otimes 1)=\sigma^{-1}(1\otimes x)=\varepsilon(x)$$

We recall that given a left 2-cocycle $\sigma$ on $A$, one can form the 2-cocycle twist $A^\sigma$ as follows. As a coalgebra, $A^\sigma=A$, and an element $x\in A$, when considered in $A^\sigma$, is denoted $[x]$. The product in $A^\sigma$ is then defined, in Sweedler notation, as follows: 
$$[x][y]=\sum\sigma(x_1\otimes y_1)\sigma^{-1}(x_3\otimes y_3)[x_2y_2]$$ 

As for the antipode of $A^\sigma$, this is given by the following formula:
$$S^{\sigma}([x])=\sum\sigma(x_1\otimes S(x_2))\sigma^{-1}(S(x_4)\otimes x_5)[S(x_3)]$$

We can now state and prove a main theorem, from \cite{bbs}, regarding the twisting of the quantum automorphism groups of the group algebras, as follows:

\begin{theorem}
If $G$ is a finite group and $\sigma$ is a $2$-cocycle on $G$, the Hopf algebras
$$F(S_{\widehat{G}}^+)\quad,\quad F(S_{\widehat{G}_\sigma}^+)$$
are $2$-cocycle twists of each other, in the above sense.
\end{theorem}

\begin{proof}
In order to prove this result, we use the following Hopf algebra map: 
$$\pi:F(S_{\widehat{G}}^+)\to F(\widehat{G})\quad,\quad
x_{gh}\to\delta_{gh}e_g$$

Our 2-cocycle $\sigma:G\times G\to F^*$ can be extended by linearity into a linear map as follows, which is a left and right 2-cocycle in the above sense:
$$\sigma:F(\widehat{G})\otimes F(\widehat{G})\to F$$

Consider now the following composition:
$$\alpha=\sigma(\pi \otimes \pi):F(S_{\widehat{G}}^+)\otimes F(S_{\widehat{G}}^+)\to F(\widehat{G})\otimes F(\widehat{G})\to F$$

Then $\alpha$ is a left and right 2-cocycle, because it is induced by a cocycle on a group algebra, and so is its convolution inverse $\alpha^{-1}$. Thus we can construct the twisted algebra $F(S_{\widehat{G}}^+)^{\alpha^{-1}}$, and inside this algebra we have the following computation:
\begin{eqnarray*}
[x_{gh}][x_{rs}]
&=&\alpha^{-1}(x_g\otimes x_r)\alpha(x_h\otimes x_s)[x_{gh}x_{rs}]\\
&=&\sigma^{-1}(g\otimes r)\sigma(h\otimes s)[x_{gh}x_{rs}]
\end{eqnarray*}

By using this identity, we obtain the following formula:
\begin{eqnarray*}
\sum_{t\in G}\sigma(st^{-1}\otimes t)[x_{st^{-1},g}][x_{th}]
&=&\sum_{t\in G}\sigma(st^{-1}\otimes t)\sigma^{-1}(st^{-1}\otimes t)\sigma_{gh}[x_{st^{-1},g}x_{th}]\\
&=&\sigma(g\otimes h)[x_{s,gh}]
\end{eqnarray*}

Similarly, we have the following formula:
$$\sum_{t\in G}\sigma^{-1}(t^{-1}\otimes ts)[x_{g,t^{-1}}][x_{h,ts}]=\sigma^{-1}(g\otimes h)[x_{gh,s}]$$

We deduce from this that there exists a Hopf algebra map, as follows:
$$\Phi:F(S_{\widehat{G}_\sigma}^+)\to F(S_{\widehat{G}}^+)^{\alpha^{-1}}\quad,\quad 
x_{gh}\to [x_{g,h}]$$

This map is clearly surjective, and is injective as well, by a standard fusion semiring argument, because both Hopf algebras have the same fusion semiring. See \cite{bbs}.
\end{proof}

Thus, we have proved our twisting result. As a first consequence, we have: 

\begin{proposition}
If $G$ is a finite group and $\sigma$ is a $2$-cocycle on $G$, then
$$\Phi(x_{g_1h_1}\ldots x_{g_mh_m})=\Omega(g_1,\ldots,g_m)^{-1}\Omega(h_1,\ldots,h_m)x_{g_1h_1}\ldots x_{g_mh_m}$$
with the coefficients on the right being given by the formula
$$\Omega(g_1,\ldots,g_m)=\prod_{k=1}^{m-1}\sigma_{g_1\ldots g_k,g_{k+1}}$$
is a coalgebra isomorphism $F(S_{\widehat{G}_\sigma}^+)\to F(S_{\widehat{G}}^+)$, commuting with the Haar integrals.
\end{proposition}

\begin{proof}
This is indeed just a technical reformulation of Theorem 8.11.
\end{proof}

Let us discuss now some concrete applications of these twisting results, over $F=\mathbb C$. Consider the group $G=\mathbb Z_n^2$, let $w=e^{2\pi i/n}$, and consider the following cocycle: 
$$\sigma:G\times G\to\mathbb T\quad,\quad 
\sigma_{(ij)(kl)}=w^{jk}$$ 

In order to understand what is the formula that we obtain, we must do some computations. Let $E_{ij}$ with $i,j \in\mathbb Z_n$ be the standard basis of $M_n(\mathbb C)$. We first have:

\begin{proposition}
The linear map given by
$$\psi(e_{(i,j)})=\sum_{k=0}^{n-1}{w}^{ki}E_{k,k+j}$$
defines an isomorphism of algebras $\psi:C(\widehat{G}_\sigma)\simeq M_n(\mathbb C)$. 
\end{proposition}

\begin{proof}
Consider indeed the following linear map:
$$\psi'(E_{ij})=\frac{1}{n}\sum_{k=0}^{n-1}{w}^{-ik}e_{(k,j-i)}$$
 
It is routine to check that both $\psi,\psi'$ are morphisms of algebras, and that these maps are inverse to each other. In particular, $\psi$ is an isomorphism of algebras, as stated.
\end{proof}

Next in line, we have the following result:

\begin{proposition}
The algebra map given by
$$\varphi(u_{ij}u_{kl}) = \frac{1}{n}\sum_{a,b=0}^{n-1}{w}^{ai-bj}x_{(a,k-i),(b,l-j)}$$
defines a Hopf algebra isomorphism $\varphi:C(S_{M_n}^+)\simeq C(S_{\widehat{G}_\sigma}^+)$.
\end{proposition}

\begin{proof}
Consider the universal coactions on the two algebras in the statement:
\begin{eqnarray*}
\alpha:M_n(\mathbb C)&\to&M_n({\mathbb C})\otimes C(S_{M_n}^+)\\
\beta:C(\widehat{G}_\sigma)&\to&C(\widehat{G}_\sigma)\otimes C(S_{\widehat{G}_\sigma}^+)
 \end{eqnarray*}
 
In terms of the standard bases, these coactions are given by:
\begin{eqnarray*}
\alpha(E_{ij})&=&\sum_{kl}E_{kl}\otimes u_{ki}u_{lj}\\
\beta(e_{(i,j)})&=&\sum_{kl} e_{(k,l)}\otimes x_{(k,l),(i,j)}
\end{eqnarray*}

We use now the identification $C(\widehat{G}_\sigma)\simeq M_n(\mathbb C)$ from Proposition 8.13. This identification produces a coaction map, as follows:
$$\gamma:M_n(\mathbb C)\to M_n(\mathbb C)\otimes C(S_{\widehat{G}_\sigma}^+)$$

Now observe that this map is given by the following formula:
 $$\gamma(E_{ij})=\frac{1}{n}\sum_{ab}E_{ab}\otimes\sum_{kr}w^{ar-ik} x_{(r,b-a),(k,j-i)}$$

By comparing with the formula of $\alpha$, we obtain the isomorphism in the statement.
\end{proof}

We will need one more result of this type, as follows:

\begin{proposition}
The algebra map given by
$$\rho(x_{(a,b),(i,j)})=\frac{1}{n^2}\sum_{klrs}w^{ki+lj-ra-sb}p_{(r,s),(k,l)}$$
defines a Hopf algebra isomorphism $\rho:C(S_{\widehat{G}}^+)\simeq C(S_G^+)$.
\end{proposition}

\begin{proof}
We have a Fourier transform isomorphism, as follows:
$$C(\widehat{G})\simeq C(G)$$

Thus the algebras in the statement are indeed isomorphic.
\end{proof}

As a conclusion to all this, we have the following result:

\begin{theorem}
Let $n\geq 2$ and $w=e^{2\pi i/n}$. Then
$$\Theta(u_{ij}u_{kl})=\frac{1}{n}\sum_{ab=0}^{n-1}w^{-a(k-i)+b(l-j)}p_{ia,jb}$$
defines a coalgebra isomorphism $C(PO_n^+)\to C(S_{n^2}^+)$ commuting with the Haar integrals.
\end{theorem}
 
\begin{proof}
We know that we have identifications as follows, where the projective version of $(A,u)$ is the pair $(PA,v)$, with $PA=<v_{ij}>$ and $v=u\otimes\bar{u}$:
$$PO_n^+=PU_n^+=S_{M_n}^+$$

With this in hand, the result follows from Theorem 8.11 and Proposition 8.12, by combining them with the various isomorphisms established above.
\end{proof}

Many other things can be said, as a continuation of the above. We will be back to this later in this book, when systematically discussing what happens over $F=\mathbb C$.

\section*{8d. Matrix models} 

Let us discuss now some further applications of the above methods and results, to the matrix modeling questions for our quantum permutation groups, again in the case $F=\mathbb C$. Let us first recall from chapter 4 that we have the following notion:

\begin{definition}
A matrix model for an affine Hopf algebra $(A,u)$ is a morphism of associative algebras as follows, with $T$ being a certain space:
$$\pi:A\to M_K(F(T))\quad,\quad u_{ij}\to U_{ij}$$
When this morphism $\pi$ is an inclusion, we say that our model is faithful.
\end{definition}

Obviously, we would like our models to be faithful, in order for our computations in the random matrix algebra $M_K(F(T))$ to be relevant, to our questions regarding $A$. 

\bigskip

But, the problem is that this situation is not always possible, due to a number of analytic reasons, the idea here being that the random matrix algebras $M_K(F(T))$ are quite ``thin'', from a certain functional analytic viewpoint, while the algebra $A$ to be modeled might be ``thick'', from the same functional analytic viewpoint.

\bigskip

However, and here comes our point, the notion of inner faithfulness introduced and studied in chapter 4 provide us with a potential solution to this, as follows:

\begin{definition}
We say that a matrix model as above,
$$\pi:A\to M_K(F(T))\quad,\quad u_{ij}\to U_{ij}$$
is inner faithful when there is no Hopf algebra factorization as follows:
$$\pi:A\to B\to M_K(F(T))\quad,\quad u_{ij}\to v_{ij}\to U_{ij}$$
That is, we can use our inner faithfulness notion, for the matrix models.
\end{definition} 

And the point is that, with this notion in hand, we can model far more affine Hopf algebras than before, with the above-mentioned analytic obstructions dissapearing. In fact, there is no known obstruction on the algebras $A$ than can be modeled as above.

\bigskip

Still in relation with the notions from chapter 4, let us record as well:

\begin{theorem}
Given an arbitrary matrix model, as before,
$$\pi:A\to M_K(F(T))\quad,\quad u_{ij}\to U_{ij}$$
we can always factorize it via a smallest Hopf algebra, as follows,
$$\pi:A\to B\to M_K(F(T))\quad,\quad u_{ij}\to v_{ij}\to U_{ij}$$
and the resulting factorized model $B\to M_K(F(T))$ is then inner faithful.
\end{theorem} 

\begin{proof}
This is indeed something self-explanatory, coming from definitions.
\end{proof}

According to our general matrix model philosophy, we must look for concrete realizations of the standard coordinates $u_{ij}\in C(G)$. In the quantum permutation group case, $G\subset S_N^+$, these standard coordinates are projections, so a natural idea is that of looking at the rank of these projections, in the model. We are led to the following notion:

\begin{definition}
Given a subgroup $G\subset S_N^+$, a random matrix model of type
$$\pi:C(G)\to M_K(C(T))$$
is called flat when the fibers $P_{ij}^x=\pi(u_{ij})(x)$ all have rank $1$. 
\end{definition}

We will see in what follows that there are many interesting examples of such models. Observe however, as an obstruction, that in order for such a model to exist, the subgroup $G\subset S_N^+$ must be transitive, due to the following trivial fact:
$$P_{ij}\neq0\implies u_{ij}\neq0$$

In view of this, it makes sense to talk, more generally, about matrix models which are quasi-flat, in the sense that the fibers $P_{ij}^x=\pi(u_{ij})(x)$ all have rank $\leq1$, and we refer here to the literature. In what follows we will only use flat models, as axiomatized above.

\bigskip

In order to construct now some explicit models, following \cite{bne}, which will lead to some interesting new twists too, we will need the following standard notion:

\index{Weyl matrix}
\index{Pauli matrix}

\begin{definition}
Given a finite abelian group $H$, the associated Weyl matrices are
$$W_{ia}:e_b\to<i,b>e_{a+b}$$
where $i\in H$, $a,b\in\widehat{H}$, and where $(i,b)\to<i,b>$ is the Fourier coupling $H\times\widehat{H}\to\mathbb T$.
\end{definition}

As a basic illustration here, consider the simplest cyclic group, namely:
$$H=\mathbb Z_2=\{0,1\}$$

In this case, the Fourier coupling is given by $<i,b>=(-1)^{ib}$, and so the Weyl matrices act in the following way, in the standard basis:
$$W_{00}:e_b\to e_b\quad,\quad 
W_{10}:e_b\to(-1)^be_b$$
$$W_{11}:e_b\to(-1)^be_{b+1}\quad,\quad 
W_{01}:e_b\to e_{b+1}$$

Thus, we have the following formulae for the Weyl matrices, in this case:
$$W_{00}=\begin{pmatrix}1&0\\0&1\end{pmatrix}\quad,\quad 
W_{10}=\begin{pmatrix}1&0\\0&-1\end{pmatrix}$$
$$W_{11}=\begin{pmatrix}0&-1\\1&0\end{pmatrix}\quad,\quad
W_{01}=\begin{pmatrix}0&1\\1&0\end{pmatrix}$$

We recognize here, up to some multiplicative factors, the four Pauli matrices. 

\bigskip

Now back to the general case, we will need the following well-known result:

\begin{proposition}
The Weyl matrices are unitaries, and satisfy:
\begin{enumerate}
\item $W_{ia}^*=<i,a>W_{-i,-a}$.

\item $W_{ia}W_{jb}=<i,b>W_{i+j,a+b}$.

\item $W_{ia}W_{jb}^*=<j-i,b>W_{i-j,a-b}$.

\item $W_{ia}^*W_{jb}=<i,a-b>W_{j-i,b-a}$.
\end{enumerate}
\end{proposition}

\begin{proof}
The unitary follows from (3,4), and the rest of the proof goes as follows:

\medskip

(1) We have indeed the following computation:
\begin{eqnarray*}
W_{ia}^*
&=&\left(\sum_b<i,b>E_{a+b,b}\right)^*\\
&=&\sum_b<-i,b>E_{b,a+b}\\
&=&\sum_b<-i,b-a>E_{b-a,b}\\
&=&<i,a>W_{-i,-a}
\end{eqnarray*}

(2) Here the verification goes as follows:
\begin{eqnarray*}
W_{ia}W_{jb}
&=&\left(\sum_d<i,b+d>E_{a+b+d,b+d}\right)\left(\sum_d<j,d>E_{b+d,d}\right)\\
&=&\sum_d<i,b><i+j,d>E_{a+b+d,d}\\
&=&<i,b>W_{i+j,a+b}
\end{eqnarray*}

(3,4) By combining the above two formulae, we obtain:
\begin{eqnarray*}
W_{ia}W_{jb}^*
&=&<j,b>W_{ia}W_{-j,-b}\\
&=&<j,b><i,-b>W_{i-j,a-b}
\end{eqnarray*}

We obtain as well the following formula:
\begin{eqnarray*}
W_{ia}^*W_{jb}
&=&<i,a>W_{-i,-a}W_{jb}\\
&=&<i,a><-i,b>W_{j-i,b-a}
\end{eqnarray*}

But this gives the formulae in the statement, and we are done.
\end{proof}

Observe that, with $n=|H|$, we can use an isomorphism $l^2(\widehat{H})\simeq\mathbb C^n$ as to view each $W_{ia}$ as a usual matrix, $W_{ia}\in M_n(\mathbb C)$, and hence as a usual unitary, $W_{ia}\in U_n$. 

\bigskip

Given a vector $\xi$, we denote by $Proj(\xi)$ the orthogonal projection onto $\mathbb C\xi$. With these conventions, we have the following result:

\begin{theorem}
Given a closed subgroup $E\subset U_n$, we have a representation
$$\pi_H:C(S_N^+)\to M_N(C(E))$$
$$w_{ia,jb}\to[U\to Proj(W_{ia}UW_{jb}^*)]$$
where $n=|H|,N=n^2$, and where $W_{ia}$ are the Weyl matrices associated to $H$.
\end{theorem}

\begin{proof}
The Weyl matrices being given by $W_{ia}:e_b\to<i,b>e_{a+b}$, we have:
$$tr(W_{ia})=\begin{cases}
1&{\rm if}\ (i,a)=(0,0)\\
0&{\rm if}\ (i,a)\neq(0,0)
\end{cases}$$

Together with the formulae in Proposition 8.22, this shows that the Weyl matrices are pairwise orthogonal with respect to the following scalar product on $M_n(\mathbb C)$:
$$<x,y>=tr(xy^*)$$

Thus, these matrices form an orthogonal basis of $M_n(\mathbb C)$, consisting of unitaries:
$$W=\left\{W_{ia}\Big|i\in H,a\in\widehat{H}\right\}$$

Thus, each row and each column of the matrix $\xi_{ia,jb}=W_{ia}UW_{jb}^*$ is an orthogonal basis of $M_n(\mathbb C)$, and so the corresponding projections form a magic unitary, as claimed.
\end{proof}

Still following \cite{bne}, let us discuss now some interesting generalizations of the Weyl matrix models. We will need the following standard definition:

\begin{definition}
A $2$-cocycle on a group $G$ is a function $\sigma:G\times G\to\mathbb T$ satisfying:
$$\sigma(gh,k)\sigma(g,h)=\sigma(g,hk)\sigma(h,k)\quad,\quad 
\sigma(g,1)=\sigma(1,g)=1$$
The algebra $C^*(G)$, with multiplication given by $g\cdot h=\sigma(g,h)gh$, and with the involution making the standard generators $g\in C^*_\sigma(G)$ unitaries, is denoted $C^*_\sigma(G)$.
\end{definition}

As explained in \cite{bne}, we have the following general construction:

\begin{proposition}
Given a finite group $G=\{g_1,\ldots,g_N\}$ and a $2$-cocycle on it, $\sigma:G\times G\to\mathbb T$, we have a matrix model as follows, 
$$\pi:C(S_N^+)\to M_N(C(E))\quad,\quad 
w_{ij}\to[x\to Proj(g_ixg_j^*)]$$
for any closed subgroup $E\subset U_A$, where $A=C^*_\sigma(G)$. 
\end{proposition}

\begin{proof}
This is clear from definitions, because the standard generators $\{g_1,\ldots,g_N\}$ are pairwise orthogonal with respect to the canonical trace of $A$. See \cite{bne}.
\end{proof}

Now with the above construction in hand, we have the following result:

\begin{theorem}
For any intermediate closed subgroup $G\subset E\subset U_A$, the model 
$$\pi:C(S_N^+)\to M_N(C(E))$$
constructed above is faithful.
\end{theorem}

\begin{proof}
This is something quite tricky, from \cite{bne}, the idea being as follows:

\medskip

(1) To start with, in the present analytic setting, over $F=\mathbb C$, the faithfulness of a matrix model $u_{ij}\to U_{ij}$ can be verified with respect to the trace, the criterion here being that we must have the formula $T_p^2=T_p$, for any $p\in\mathbb N$ and any $e\in\{1,*\}^p$, where:
$$(T_p)_{i_1\ldots i_p,j_1\ldots j_p}=\left(tr\otimes\int_T\right)(U_{i_1j_1}^{e_1}\ldots U_{i_pj_p}^{e_p})$$

We will be back to this faithfulness criterion, with full details, in Part III below, when discussing more in detail what happens in the analytic setting, over $F=\mathbb C$.

\medskip

(2) In our case, according to the definition of $T_p$, we have the following formula:
\begin{eqnarray*}
(T_p)_{i_1\ldots i_p,j_1\ldots j_p}
&=&\left(tr\otimes\int_E\right)\left(Proj(g_{i_1}xg_{j_1}^*)\ldots Proj(g_{i_p}xg_{j_p}^*)\right)dx\\
&=&\frac{1}{N}\int_E<g_{i_1}xg_{j_1}^*,g_{i_2}xg_{j_2}^*>\ldots<g_{i_p}xg_{j_p}^*,g_{i_1}xg_{j_1}^*>dx
\end{eqnarray*}

In order to compute the scalar products, we can use the following formula:
$$g_ig_{i^{-1}k}=\sigma(i,i^{-1}k)g_k$$

Indeed, we obtain from this the following formula:
$$g_i^*g_k=\overline{\sigma(i,i^{-1}k)}g_{i^{-1}k}$$

We therefore obtain the following formula, for the above scalar products:
\begin{eqnarray*}
<g_ixg_j^*,g_kxg_l^*>
&=&tr(g_jx^*g_i^*g_kxg_l^*)\\
&=&tr(g_i^*g_kxg_l^*g_jx^*)\\
&=&\overline{\sigma(i,i^{-1}k)}\cdot\overline{\sigma(l,l^{-1}j)}\cdot tr(g_{i^{-1}k}xg_{l^{-1}j}x^*)
\end{eqnarray*}

By plugging these quantities into the formula of $T_p$, we obtain the following formula:
\begin{eqnarray*}
(T_p)_{i_1\ldots i_p,j_1\ldots j_p}
&=&\overline{\sigma(i_1,i_1^{-1}i_2)}\ldots\overline{\sigma(i_p,i_p^{-1}i_1)}\cdot\overline{\sigma(j_2,j_2^{-1}j_1)}\ldots\overline{\sigma(j_1,j_1^{-1}j_p)}\\
&&\frac{1}{N}\int_Etr(g_{i_1^{-1}i_2}xg_{j_2^{-1}j_1}x^*)\ldots\ldots tr(g_{i_p^{-1}i_1}xg_{j_1^{-1}j_p}x^*)dx
\end{eqnarray*}

(3) Now let us write $(T_p)_{ij}=\rho(i,j)(T_p^\circ)_{ij}$, where $\rho(i,j)$ is the product of $\sigma$ terms appearing in the above formula. With this convention, we have:
\begin{eqnarray*}
(T_p^2)_{ij}
&=&\sum_k(T_p)_{ik}(T_p)_{kj}\\
&=&\sum_k\rho(i,k)\rho(k,j)(T_p^\circ)_{ik}(T_p^\circ)_{kj}
\end{eqnarray*}

Let us first compute the $\rho$ term. We have:
\begin{eqnarray*}
\rho(i,k)\rho(k,j)
&=&\overline{\sigma(i_1,i_1^{-1}i_2)}\ldots\overline{\sigma(i_p,i_p^{-1}i_1)}\cdot\overline{\sigma(k_2,k_2^{-1}k_1)}\ldots\overline{\sigma(k_1,k_1^{-1}k_p)}\\
&&\overline{\sigma(k_1,k_1^{-1}k_2)}\ldots\overline{\sigma(k_p,k_p^{-1}k_1)}\cdot\overline{\sigma(j_2,j_2^{-1}j_1)}\ldots\overline{\sigma(j_1,j_1^{-1}j_p)}\\
&=&\sigma(i,j)\cdot\overline{\sigma(k_2,k_2^{-1}k_1)}\cdot\overline{\sigma(k_1,k_1^{-1}k_2)}\ldots\ldots\overline{\sigma(k_1,k_1^{-1}k_p)}\cdot\overline{\sigma(k_p,k_p^{-1}k_1)}
\end{eqnarray*}

In order to compute now this quantity, observe that by multiplying the formulae $\sigma(i,i^{-1}k)g_i^*g_k=g_{i^{-1}k}$ and $\sigma(k,k^{-1}i)g_k^*g_i=g_{k^{-1}i}$ we obtain the following formula:
$$\sigma(i,i^{-1}k)\sigma(k,k^{-1}i)=\sigma(i^{-1}k,k^{-1}i)$$

Thus, our expression further simplifies, as follows:
$$\rho(i,k)\rho(k,j)=\sigma(i,j)\cdot\overline{\sigma(k_2^{-1}k_1,k_1^{-1}k_2)}\ldots\ldots\overline{\sigma(k_1^{-1}k_p,k_p^{-1}k_1)}$$

On the other hand, the $T^\circ$ term can be written as follows:
\begin{eqnarray*}
(T_p^\circ)_{ik}(T_p^\circ)_{kj}
&=&\frac{1}{N^2}\int_E\int_Etr(g_{i_1^{-1}i_2}xg_{k_2^{-1}k_1}x^*)tr(g_{k_1^{-1}k_2}yg_{j_2^{-1}j_1}y^*)\\
&&\hskip40mm\ldots\ldots\\
&&\hskip10mm tr(g_{i_p^{-1}i_1}xg_{k_1^{-1}k_1p}x^*)tr(g_{k_p^{-1}k_1}yg_{j_1^{-1}j_p}y^*)dxdy
\end{eqnarray*}

We therefore conclude that we have the following formula:
\begin{eqnarray*}
(T_p^2)_{ij}
&=&\frac{\sigma(i,j)}{N^2}\int_E\int_E\sum_{k_1\ldots k_p}
\overline{\sigma(k_2^{-1}k_1,k_1^{-1}k_2)}tr(g_{i_1^{-1}i_2}xg_{k_2^{-1}k_1}x^*)tr(g_{k_1^{-1}k_2}yg_{j_2^{-1}j_1}y^*)\\
&&\hskip50mm\vdots\\
&&\hskip10mm\overline{\sigma(k_1^{-1}k_p,k_p^{-1}k_1)}tr(g_{i_p^{-1}i_1}xg_{k_1^{-1}k_p}x^*)tr(g_{k_p^{-1}k_1}yg_{j_1^{-1}j_p}y^*)dxdy
\end{eqnarray*}

(4) By using now $g_i^*=\overline{\sigma(i,i^{-1})}g_{i^{-1}}$, and moving the $x^*$ variables at left, we obtain:
\begin{eqnarray*}
(T_p^2)_{ij}
&=&\frac{\sigma(i,j)}{N^2}\int_E\int_E\sum_{k_1\ldots k_p}
tr(x^*g_{i_1^{-1}i_2}xg_{k_2^{-1}k_1})tr(g_{k_2^{-1}k_1}^*yg_{j_2^{-1}j_1}y^*)\\
&&\hskip50mm\vdots\\
&&\hskip10mm tr(x^*g_{i_p^{-1}i_1}xg_{k_1^{-1}k_p})tr(g_{k_1^{-1}k_p}^*yg_{j_1^{-1}j_p}y^*)dxdy
\end{eqnarray*}

We can compute the products of traces by using the following formula:
\begin{eqnarray*}
tr(Ag_k)tr(g_k^*B)
&=&\sum_{qs}<g_q,Ag_k><g_s,g_k^*B>\\
&=&\sum_{qs}tr(g_q^*Ag_k)tr(g_s^*g_k^*B)
\end{eqnarray*}

Thus are left with an integral involving the variable $z=xy$, which gives $T_p^\circ$.
\end{proof}

The above construction is something quite general, and in relation with the Weyl matrix models discussed before, we have the following result:

\begin{theorem}
Given a finite abelian group $H$, consider the product $G=H\times\widehat{H}$, and endow it with its standard Fourier cocycle. 
\begin{enumerate}
\item With $E=U_n$, where $n=|H|$, the model $\pi:C(S_N^+)\to M_N(C(U_n))$ constructed above, where $N=n^2$, is the Weyl matrix model associated to $H$.

\item When assuming in addition that $H$ is cyclic, $H=\mathbb Z_n$, we obtain in this way the matrix model for $C(S_N^+)$ coming from the usual Weyl matrices.

\item In the particular case $H=\mathbb Z_2$, the model $\pi:C(S_4^+)\to M_N(C(U_2))$ constructed above is the matrix model for $C(S_4^+)$ coming from the Pauli matrices.
\end{enumerate}
\end{theorem}

\begin{proof}
All this is indeed clear from definitions, and from the above results.
\end{proof}

Many more things can be said about the above models, and all this suggests that the quantum group associated to a Weyl matrix model, as above, should appear as a suitable twist of $PE$. For more on all this, we refer to \cite{bb1}, \cite{bco}, \cite{bne} and related papers.

\section*{8e. Exercises}

We had an advanced algebraic chapter here, and as exercises on this, we have:

\begin{exercise}
Further build on our twisted action results above.
\end{exercise}

\begin{exercise}
Work out the general theory of pseudo-2-cocycle twists.
\end{exercise}

\begin{exercise}
Fill in all the details for the ADE classification result.
\end{exercise}

\begin{exercise}
Find matrix models for the group dual subgroups of $S_N^+$.
\end{exercise}

\begin{exercise}
Work out the Weyl twisting results in the Pauli matrix case.
\end{exercise}

\begin{exercise}
Discuss twisting questions for the general Weyl matrix models.
\end{exercise}

As bonus exercise, learn about braided Hopf algebras, and $R$-matrix deformation.

\part{Quantum groups}

\ \vskip50mm

\begin{center}
{\em Voi canta 

Pentru mileniul trei

Vreau, vreau sa-i las 

Un semn din anii mei}
\end{center}

\chapter{Complex algebras}

\section*{9a. Complex algebras}

Welcome to quantum groups. These are typically defined over the complex numbers, so let us first discuss what happens when the field of scalars is $F=\mathbb C$, and when our algebras are $*$-algebras. We will do this in this chapter, with a quick remake of Parts I-II, by insisting on various simplifications, and other new phenomena that might appear. 

\bigskip

Following the material in chapter 1, let us start with the following definition:

\begin{definition}
A Hopf $*$-algebra is a $*$-algebra $A$, with morphisms of $*$-algebras
$$\Delta:A\to A\otimes A\quad,\quad 
\varepsilon:A\to\mathbb C\quad,\quad 
S:A\to A^{opp}$$
called comultiplication, counit and antipode, satisfying the following conditions:
$$(\Delta\otimes id)\Delta=(id\otimes \Delta)\Delta$$
$$(\varepsilon\otimes id)\Delta=(id\otimes\varepsilon)\Delta=id$$
$$m(S\otimes id)\Delta=m(id\otimes S)\Delta=\varepsilon(.)1$$
If the square of the antipode is the identity, $S^2=id$, we say that $A$ is undeformed. Otherwise, in the case $S^2\neq id$, we say that $A$ is deformed.
\end{definition}

In other words, we are using here exactly the same formalism as in chapter 1, and in the whole first half of this book, save for our assumption that $A$ is a $*$-algebra, over $F=\mathbb C$. It is in fact possible to talk about $*$-algebras $A$ defined over an algebrically closed field $F$ of characteristic 0, and most of our results below will extend to this setting. 

\bigskip

However, the field of complex numbers $F=\mathbb C$ has its own magic, hope we agree on this, that will come into play later, so we will stay with it, as indicated in Definition 9.1. To be more precise, for putting things squarely, we are currently looking for the quantum groups, which are most likely quite complicated objects, so in our quest we will not bother much, and use of course the smartest field of scalars of them all, which is $F=\mathbb C$.

\bigskip

At the level of the examples now, we first have the following result:

\begin{theorem}
Given a finite group $G$, the $*$-algebra of complex functions on it, $\mathbb C(G)=\{\varphi:G\to\mathbb C\}$, with the usual pointwise product of functions,
$$(\varphi\psi)(g)=\varphi(g)\psi(g)$$
is a Hopf $*$-algebra, with comultiplication, counit and antipode as follows:
$$\Delta(\varphi)=[(g,h)\to \varphi(gh)]$$
$$\varepsilon(\varphi)=\varphi(1)$$
$$S(\varphi)=[g\to\varphi(g^{-1})]$$
This Hopf algebra is finite dimensional, commutative, and undeformed.
\end{theorem}

\begin{proof}
This is something that we know from chapter 1, with the only new input coming from the fact that $\mathbb C(G)$ is a $*$-algebra, with involution as follows:
$$\varphi^*(g)=\overline{\varphi(g)}$$

Indeed, with this formula in hand, we can see that the morphisms in the statement $\Delta,\varepsilon,S$ are morphisms of $*$-algebras, so we have indeed a Hopf $*$-algebra, as stated.
\end{proof}

Next, still following the material in chapter 1, we have as well the following result:

\begin{theorem}
Given a group $H$, which can be finite or not, its group algebra
$$\mathbb C[H]=span(H)$$
is a Hopf $*$-algebra, with structural maps given on group elements as follows:
$$\Delta(g)=g\otimes g\quad,\quad
\varepsilon(g)=1\quad,\quad
S(g)=g^{-1}$$
This Hopf algebra is cocommutative, and undeformed.
\end{theorem}

\begin{proof}
As before, this is something that we know from chapter 1, with the only new input coming from the fact that $\mathbb C[H]$ is a $*$-algebra, with involution as follows:
$$\left(\sum_g\lambda_gg\right)^*=\sum_g\bar{\lambda}_gg^{-1}$$

Indeed, with this formula in hand, we can see that the morphisms in the statement $\Delta,\varepsilon,S$ are morphisms of $*$-algebras, so we have indeed a Hopf $*$-algebra, as stated.
\end{proof}

We refer to chapter 1 for various technical reformulations of the above result, say by using convolution products, as well as for all sorts of comments that can be made.

\bigskip

In order to discuss now the precise relation between Theorem 9.2 and Theorem 9.3, we will need the following fact, that we also basically know from chapter 1:

\begin{theorem}
If $G,H$ are finite abelian groups, dual to each other via Pontrjagin duality, in the sense that each of them is the character group of the other,
$$G=\big\{\chi:H\to\mathbb T\big\}\quad,\quad H=\big\{\rho:G\to\mathbb T\big\}$$
we have an identification of Hopf $*$-algebras as follows:
$$\mathbb C(G)=\mathbb C[H]$$
In the case $G=H=\mathbb Z_N$, this identification is the usual discrete Fourier transform isomorphism. In general, we obtain a tensor product of such Fourier transforms.
\end{theorem}

\begin{proof}
As before, this is something that we know from chapter 1, with the only new input coming from the fact that we are now dealing with Hopf $*$-algebras. Indeed, the isomorphism $\mathbb C(G)\simeq\mathbb C[H]$ that we constructed in chapter 1 is easily seen to be a $*$-algebra morphism, so we have an isomorphism of Hopf $*$-algebras, as stated.
\end{proof}

Moving ahead, in what regards now the main theorem established in chapter 1, which was the duality one, in the finite dimensional case, things again extend well, as follows:

\begin{theorem}
Given a finite dimensional Hopf $*$-algebra $A$, its dual space
$$A^*=\Big\{\varphi:A\to\mathbb C\ {\rm linear}\Big\}$$
is also a finite dimensional Hopf $*$-algebra, with multiplication and unit as follows,
$$\Delta^t:A^*\otimes A^*\to A^*\quad,\quad \varepsilon^t:\mathbb C\to A^*$$
and with comultiplication, counit and antipode as follows:
$$m^t:A^*\to A^*\otimes A^*\quad,\quad
u^t:A^*\to\mathbb C\quad,\quad 
S^t:A^*\to A^*$$
This duality makes correspond the commutative algebras to the cocommutative algebras. Also, this duality makes correspond $\mathbb C(G)$ to $\mathbb C[G]$, for any finite group $G$.
\end{theorem}

\begin{proof}
Again, this is something that we know from chapter 1, with the few verifications in relation with the $*$-structure being all trivial, basically coming from the fact that the transpose of a $*$-algebra morphism is again a $*$-algebra morphism.
\end{proof}

As a conclusion to all this, again following chapter 1, let us formulate:

\begin{speculation}
We can think of any finite dimensional Hopf $*$-algebra $A$, not necessarily commutative or cocommutative, as being of the form
$$A=\mathbb C(G)=\mathbb C[H]$$
with $G,H$ being finite quantum groups, related by a generalized Pontrjagin duality. And with this generalizing what we know about the abelian groups.
\end{speculation}

And good news, this is all, in what concerns the quick remake of chapters 1-2. There are of course still a number of things there, waiting to be discussed in the present $*$-algebra setting, especially in what regards the Lie algebra material in chapter 2. But, we will keep this for later, in chapter 12 below, where we will discuss Lie algebra aspects.

\bigskip

Getting now into the material from chapters 3-4, and afterwards, we have split the remake of that, in the present $*$-algebra setting, into 3 parts: we will first discuss corepresentations, then various product operations, and then symmetry groups and twists. 

\bigskip

So, let us first talk about corepresentations of Hopf $*$-algebras. We will choose these corepresentations to be unitary, and we have here the following definition, that we know perfectly well, since chapter 3, save for the added unitarity assumption:

\begin{definition}
A unitary corepresentation of a Hopf $*$-algebra $A$ is a unitary matrix $u\in M_N(A)$ satisfying the following relations, with the first one implying the other two:
$$(id\otimes\Delta)u=u_{12}u_{13}\quad,\quad (id\otimes \varepsilon)u=1\quad,\quad (id\otimes S)u=u^{-1}$$
The contragradient, or conjugate corepresentation to $u$, is constructed as follows:
$$\bar{u}=(t\otimes id)u^{-1}=(t\otimes S)u=(u_{ij}^*)$$
We say that our algebra $A$ is affine when it is of the form $A=<u_{ij}>$, with $u\in M_N(A)$ being a certain corepresentation, called fundamental corepresentation.
\end{definition}

We refer to chapters 3-4 for some basic theory and observations, regarding the corepresentations, and their relation with the equivalent notion of coaction. Among others, let us recall from there that when $F=\mathbb C$, and assuming that $A$ is a $*$-algebra, and that $u$ is a unitary corepresentation, $u^*=u^{-1}$, as above, the conjugate is simply given by:
$$\bar{u}_{ij}=u_{ij}^*$$

The corepresentations are subject to a number of operations, as follows:

\index{sum of corepresentations}
\index{product of corepresentations}
\index{tensor product of corepresentations}
\index{conjugate of corepresentation}
\index{spinned corepresentation}
\index{spinning by unitaries}

\begin{theorem}
The unitary corepresentations are subject to:
\begin{enumerate}
\item Making sums, $u+v=diag(u,v)$.

\item Making tensor products, $(u\otimes v)_{ia,jb}=u_{ij}v_{ab}$.

\item Taking conjugates, $\bar{u}=(u_{ij}^*)$.

\item Spinning by unitary scalar matrices, $u\to VuV^*$.
\end{enumerate}
\end{theorem}

\begin{proof}
We already know this from chapter 3, the only new input coming from the fact that the above operations preserve the unitarity, which can be checked as follows:

\medskip

(1) The unitarity of $u+v=diag(u,v)$ is indeed clear.

\medskip

(2) In order to check that $(u\otimes v)_{ia,jb}=u_{ij}v_{ab}$ is unitary, we first have:
\begin{eqnarray*}
\sum_{jb}(u\otimes v)_{ia,jb}(u\otimes v)_{kc,jb}^*
&=&\sum_{jb}u_{ij}v_{ab}v_{cb}^*u_{kj}^*\\
&=&\delta_{ac}\sum_ju_{ij}u_{kj}^*\\
&=&\delta_{ik}\delta_{ac}
\end{eqnarray*}

In the other sense, the computation is similar, as follows, which gives the result:
\begin{eqnarray*}
\sum_{jb}(u\otimes v)^*_{jb,ia}(u\otimes v)_{jb,kc}
&=&\sum_{jb}v_{ba}^*u_{ji}^*u_{jk}v_{bc}\\
&=&\delta_{ik}\sum_bv_{ba}^*v_{bc}\\
&=&\delta_{ik}\delta_{ac}
\end{eqnarray*}

(3) Our more precise claim here is that a corepresentation $u=(u_{ij})$ is a biunitary, meaning unitary, with unitary transpose. Indeed, the idea is that $u^t=\bar{u}^{-1}$ comes from $u^*=u^{-1}$, by applying the antipode. To be more precise, by denoting $(a,b)\to a\cdot b$ the multiplication of the opposite algebra $A^{opp}$, we have the following computation:
\begin{eqnarray*}
(uu^*)_{ij}=\delta_{ij}
&\implies&\sum_ku_{ik}u_{jk}^*=\delta_{ij}\\
&\implies&\sum_kS(u_{ik})\cdot S(u_{jk}^*)=\delta_{ij}\\
&\implies&\sum_ku_{kj}u_{ki}^*=\delta_{ij}\\
&\implies&(u^t\bar{u})_{ji}=\delta_{ij}
\end{eqnarray*}

Similarly, we have the following computation, which gives the result:
\begin{eqnarray*}
(u^*u)_{ij}=\delta_{ij}
&\implies&\sum_ku_{ki}^*u_{kj}=\delta_{ij}\\
&\implies&\sum_kS(u_{ki}^*)\cdot S(u_{kj})=\delta_{ij}\\
&\implies&\sum_ku_{jk}^*u_{ik}=\delta_{ij}\\
&\implies&(\bar{u}u^t)_{ji}=\delta_{ij}
\end{eqnarray*}

Thus, we are led to the conclusion in the statement.

\medskip

(4) Finally, the unitarity of $VuV^*$, when $V$ is taken unitary, is clear.
\end{proof}

In order to reach now to a more advanced theory, let us introduce as well:

\index{character}

\begin{definition}
Given a Hopf $*$-algebra $A$, and a corepresentation $u\in M_N(A)$, we define the character of $u$ to be its matrix trace:
$$\chi_u=\sum_iu_{ii}\in A$$
In the case where $A$ is an affine Hopf algebra, with $u\in M_N(A)$ being its fundamental corepresentation, we call $\chi_u\in A$ the main character of $A$. 
\end{definition}

Generally speaking, the idea in what follows will be that of proving that, under suitable assumptions on $A$, the corepresentations are uniquely determined by their characters. 

\bigskip

But more on this later, for the moment let us get more familiar with the operation $u\to\chi_u$. In relation with the various operations in Theorem 9.8, we have:

\index{character}
\index{equivalent corepresentations}

\begin{theorem}
The characters of corepresentations satisfy the formulae
$$\chi_{u+v}=\chi_u+\chi_v\quad,\quad
\chi_{u\otimes v}=\chi_u\chi_v\quad,\quad
\chi_{\bar{u}}=\chi_u^*\quad,\quad 
\chi_{VuV^*}=\chi_u$$
valid for any corepresentations $u,v$, and any unitary scalar matrix $V$. 
\end{theorem}

\begin{proof}
This is something elementary, the idea being as follows:

\medskip

(1) The first formula is clear from the definition of the sum, as follows:
$$\chi_{u+v}
=Tr(diag(u,v))
=Tr(u)+Tr(v)
=\chi_u+\chi_v$$

(2) The second formula comes from the following computation:
$$\chi_{u\otimes v}
=\sum_{ia}u_{ii}v_{aa}
=\sum_iu_{ii}\sum_av_{aa}
=\chi_u\chi_v$$

(3) The third formula comes from the following computation:
$$\chi_{\bar{u}}
=\sum_iu_{ii}^*
=\left(\sum_iu_{ii}\right)^*
=\chi_u^*$$

(4) Finally, the last formula comes from the following computation:
$$\chi_{VuV^*}
=Tr(VuV^*)
=Tr(u)
=\chi_u$$

Thus, we are led to the conclusions in the statement.
\end{proof}

As already mentioned, we will be back to such characters later, under some suitable extra assumptions on our Hopf $*$-algebras $A$, of compactness type, allowing us to develop an analogue of the Peter-Weyl theory for them. More on this later.

\bigskip

Next in line, also at the general level, we have the following key definition:

\index{intertwiner}
\index{Hom space}
\index{End space}
\index{Fix space}
\index{irreducible corepresentation}
\index{fixed points}
\index{self-intertwiner}
\index{equivalent corepresentations}

\begin{definition}
Given two corepresentations $u\in M_N(A),v\in M_M(A)$, we set 
$$Hom(u,v)=\left\{T\in M_{M\times N}(\mathbb C)\Big|Tu=vT\right\}$$
and we use the following conventions:
\begin{enumerate}
\item We use the notations $Fix(u)=Hom(1,u)$, and $End(u)=Hom(u,u)$.

\item We write $u\sim v$ when $m=n$, and $Hom(u,v)$ contains an invertible element.

\item We say that $u$ is irreducible, and write $u\in Irr(A)$, when $End(u)=\mathbb C1$.
\end{enumerate}
\end{definition}

In the function algebra case $A=\mathbb C(G)$ we obtain the usual notions concerning the representations of $G$. Observe also that in the group dual case, $A=\mathbb C[H]$, we have:
$$g\sim h\iff g=h$$

Finally, observe that $u\sim v$ means that $u,v$ are conjugated by an invertible matrix. Here are now a few basic results, regarding the above Hom spaces:

\index{tensor category}
\index{tensor $*$-category}

\begin{theorem}
We have the following results:
\begin{enumerate}
\item $T\in Hom(u,v),S\in Hom(v,w)\implies ST\in Hom(u,w)$.

\item $S\in Hom(p,q),T\in Hom(v,w)\implies S\otimes T\in Hom(p\otimes v,q\otimes w)$.

\item $T\in Hom(v,w)\implies T^*\in Hom(w,v)$.
\end{enumerate}
In other words, the Hom spaces between corepresentations form a tensor $*$-category.
\end{theorem}

\begin{proof}
These assertions are all elementary, as follows:

\medskip

(1) By using our assumptions $Tu=vT$ and $Sv=Ws$ we obtain, as desired:
$$STu
=SvT
=wST$$

(2) Assume indeed that we have $Sp=qS$ and $Tv=wT$. With standard tensor product notations, we have the following computation:
$$(S\otimes T)(p\otimes v)
=S_1T_2p_{13}v_{23}
=(Sp)_{13}(Tv)_{23}$$

On the other hand, we have as well the following computation:
$$(q\otimes w)(S\otimes T)
=q_{13}w_{23}S_1T_2
=(qS)_{13}(wT)_{23}$$

(3) This is something new, with respect to what we previously knew from Part I. Indeed, by conjugating, and then using the unitarity of $v,w$, we obtain, as desired:
\begin{eqnarray*}
Tv=wT
&\implies&v^*T^*=T^*w^*\\
&\implies&vv^*T^*w=vT^*w^*w\\
&\implies&T^*w=vT^*
\end{eqnarray*}

Finally, the last assertion follows from definitions, and from the obvious fact that, in addition to (1,2,3) above, the Hom spaces are linear spaces, and contain the units. In short, this is just a theoretical remark, that will be used only later on.
\end{proof}

We will be back to tensor $*$-categories, and their uses in the above context, later on. The general idea will be that all this leads us into linear algebra, as follows:

\bigskip

(1) To start with, Theorem 9.12 shows that the intertwiner algebras $End(u)\subset M_N(\mathbb C)$ are $*$-algebras. Thus, we can decompose such algebras $End(u)$ into components.

\bigskip

(2) But this ultimately this leads to Peter-Weyl theory, namely a decomposition of the corepresentations $u$ themselves into irreducible components, $u=u_1+\ldots+u_k$. 

\bigskip

(3) Moreover, some further work shows that the irreducible corepresentations are fully described by their characters, and with this completing the Peter-Weyl theory.

\bigskip

So, this was for the general idea, and in practice, we will discuss all this in detail in chapter 10 below, following Woronowicz \cite{wo1}, under some suitable extra assumptions on the Hopf $*$-algebras $A$ that we are considering, of compactness type.

\section*{9b. Basic operations}

Let us discuss now the operations that can be performed on our Hopf $*$-algebras. Following the material from chapter 3, we first have the following construction: 

\begin{theorem}
Given two Hopf $*$-algebras $A,B$, so is their tensor product 
$$C=A\otimes B$$
and as main illustrations for this operation, we have the following formulae:
\begin{enumerate}
\item $\mathbb C(G\times H)=\mathbb C(G)\otimes\mathbb C(H)$.

\item $\mathbb C[G\times H]=\mathbb C[G]\otimes \mathbb C[H]$.
\end{enumerate}
\end{theorem}

\begin{proof}
As before with many other things, this is something that we know from chapter 3, with the few checks in regards with the $*$-structure being all trivial.
\end{proof}

Next, and still following the material from chapter 3, we have:

\begin{theorem}
Given two Hopf $*$-algebras $A,B$, so is their free product 
$$C=A*B$$
and as main illustrations for this operation, we have the following formulae:
\begin{enumerate}
\item $\mathbb C(G\,\hat{*}\,H)=\mathbb C(G)*\mathbb C(H)$, standing as definition for $G\,\hat{*}\,H$, as quantum group.

\item $\mathbb C[G*H]=\mathbb C[G]*\mathbb C[H]$.
\end{enumerate}
\end{theorem}

\begin{proof}
Again, this is something that we know well from chapter 3, with the few needed checks in regards with the $*$-structure being all trivial.
\end{proof}

Still following the material from chapter 3, coming next, we have:

\begin{theorem}
Given a Hopf $*$-algebra $A$, so is its quotient $B=A/I$, provided that $I\subset A$ is a $*$-ideal satisfying the following conditions, called Hopf ideal conditions,
$$\Delta(I)\subset A\otimes I+I\otimes A\quad,\quad\varepsilon(I)=0\quad,\quad S(I)\subset I$$
and as main illustrations for this operation, we have the following formulae:
\begin{enumerate}
\item $\mathbb C(G)/I=\mathbb C(H)$, with $H\subset G$ being a certain subgroup.

\item $\mathbb C[G]/I=\mathbb C[H]$, with $G\to H$ being a certain quotient.
\end{enumerate}
\end{theorem}

\begin{proof}
As before, this is something that we know from chapter 3, with the few checks in regards with the $*$-structure being all trivial. Be said in passing, sure I know that all this is a bit boring for you, but please believe me that it is boring for me too, dumb copy/paste from chapter 3, by replacing $F\to\mathbb C$. By the way, as a funny story here, chapter 3, and in fact the whole Part I of this book, was basically made by starting with material from my quantum group books \cite{ba1}, \cite{ba2}, \cite{ba3}, copy/paste by replacing $\mathbb C\to F$. So now I have to switch back to $\mathbb C$, life as a math book author is not funny every day.
\end{proof}

Still following the material from chapter 3, and as something more interesting now, substantially clarifying our considerations there, on this matter, we have:

\begin{theorem}
Given a Hopf $*$-algebra $A$, and a finite dimensional unitary corepresentation $u=(u_{ij})$, the subalgebra generated by the coefficients of $u$,
$$B=<u_{ij}>\subset A$$
is a Hopf algebra. As main illustrations for this operation, we obtain subalgebras:
\begin{enumerate}
\item $\mathbb C(H)\subset\mathbb C(G)$, with $G\to H$ being certain quotients.

\item $\mathbb C[H]\subset \mathbb C[G]$, with $H\subset G$ being certain subgroups.
\end{enumerate}
\end{theorem}

\begin{proof}
This is something that we know from chapter 3, with the few checks in regards with the $*$-structure being trivial, and clarifying our previous considerations there. 
\end{proof}

Getting now to the material from chapter 4, and again as something more interesting now, clarifying our considerations there, on this matter, we have:

\begin{theorem}
Given an affine Hopf $*$-algebra $(A,u)$, we can construct its free complexification $(\widetilde{A},\tilde{u})$ as follows,
$$\widetilde{A}=<\tilde{u}>\subset\mathbb C[\mathbb Z]*A\quad,\quad\tilde{u}=zu$$
with $z=1\in\mathbb C[\mathbb Z]$. As main illustrations for this operation, we have:
\begin{enumerate}
\item $\widetilde{\mathbb C}(G)=\mathbb C(\widetilde{G})$, standing as definition for $\widetilde{G}$, as quantum group.

\item $\widetilde{\mathbb C}[H]=\mathbb C[\widetilde{H}]$, with $\widetilde{H}\subset\mathbb Z*H$ being constructed similarly.
\end{enumerate}
\end{theorem}

\begin{proof}
This is something that we know from chapter 4, with the few checks in regards with the $*$-structure being trivial, and clarifying our previous considerations there. 
\end{proof}

Still following the material from chapter 4, coming next, we have:

\begin{theorem}
Given an affine Hopf $*$-algebra $(A,u)$, we can construct its projective version $(PA,v)$ by setting
$$PA=<v_{ia,jb}>\subset A\quad,\quad v_{ia,jb}=u_{ij}\bar{u}_{ab}$$
and as main illustrations for this construction, we have the following formulae:
\begin{enumerate}
\item $P\mathbb C(G)=\mathbb C(PG)$, with $PG=G/(G\cap\mathbb T^N)$.

\item $P\mathbb C[H]=\mathbb C[PH]$, with $PH=<g_ig_j^{-1}>$, assuming $H=<g_i>$.
\end{enumerate}
\end{theorem}

\begin{proof}
Again, this is something from chapter 4, with the few checks in regards with the $*$-structure being trivial, and clarifying our previous considerations there. 
\end{proof}

As in chapter 4, we can formulate a result about these latter operations, as follows:

\begin{theorem}
Given an affine Hopf $*$-algebra $(A,u)$, construct its free complexification $(\widetilde{A},\tilde{u})$, with $\tilde{u}=zu$. We have then an identification
$$P\widetilde{A}=PA$$
and the same happens at the level of right projective versions.
\end{theorem}

\begin{proof}
As before, this is something from chapter 4, with the few checks in regards with the $*$-structure being trivial, and clarifying our previous considerations there. 
\end{proof}

Getting now to the more specialized operations on the Hopf algebras, we have:

\begin{theorem}
The quantum subgroups of a given quantum group $G$, in the general Hopf algebra sense, are subject to operations as follows:
\begin{enumerate}
\item Intersection: $H\cap K$ is the biggest quantum subgroup of $H,K$.

\item Generation: $<H,K>$ is the smallest quantum group containing $H,K$.
\end{enumerate}
\end{theorem}

\begin{proof}
As before, this is something from chapter 4, with the few checks in regards with the $*$-structure being trivial, and clarifying our previous considerations there. 
\end{proof}

Next, we have the following result, which is straightforward too:

\begin{theorem}
Given a representation $\pi:\mathbb C(G)\to A$, there exists a smallest Hopf $*$-algebra quotient $\mathbb C(L)=\mathbb C(G)/I$ producing a factorization as follows,
$$\pi:\mathbb C(G)\to\mathbb C(L)\to A$$
called Hopf image of $\pi$. More generally, we can talk about the Hopf image of a family
$$\pi_i:\mathbb C(G)\to\mathbb C(L)\to A_i$$
constructed in a similar way, by dividing by a suitable Hopf ideal.
\end{theorem}

\begin{proof}
Again, this is something from chapter 4, with the few checks in regards with the $*$-structure being trivial, and clarifying our previous considerations there. 
\end{proof}

Finally, still following chapter 4, we have the following key notion:

\begin{definition}
We say that a representation $\pi:\mathbb C(G)\to A$ is inner faithful when there is no proper factorization of type
$$\pi:\mathbb C(G)\to\mathbb C(L)\to A$$
that is, when its Hopf image is $\mathbb C(G)$ itself. Alternatively, we say that a representation $\pi:\mathbb C[H]\to A$ is inner faithful when there is no proper factorization of type
$$\pi:\mathbb C[H]\to\mathbb C[K]\to A$$
that is, when its Hopf image is $\mathbb C[H]$ itself. 
\end{definition}

As in chapter 4, as a basic result about Hopf images, we have:

\begin{theorem}
Assuming $H,K\subset G$, the quantum group $<H,K>$ is such that
$$\mathbb C(G)\to\mathbb C(H\cap K)\to\mathbb C(H),\mathbb C(K)$$
is the joint Hopf image of the following quotient maps:
$$\mathbb C(G)\to\mathbb C(H),\mathbb C(K)$$
A similar result holds for a family of quantum subgroups $H_i\subset G$.
\end{theorem}

\begin{proof}
Again, this is something from chapter 4, with the few checks in regards with the $*$-structure being trivial, and clarifying our previous considerations there. 
\end{proof}

And good news, that is all. Job well done, everything in chapter 4 extends well.

\section*{9c. Normed algebras}

Getting now to the material from Part II, in order to discuss what happens to the constructions there, in the case of the Hopf $*$-algebras $A$ over $F=\mathbb C$, it is actually convenient to further fine-tune our formalism, by adding a norm to the picture.

\bigskip

So, let us formulate the following definition, which is standard functional analysis: 

\begin{definition}
A $C^*$-algebra is a complex algebra $A$, having a norm $||.||$ making it a Banach algebra, and an involution $*$, related to the norm by the formula 
$$||aa^*||=||a||^2$$
which must hold for any $a\in A$.
\end{definition}

As a basic example, the algebra $M_N(\mathbb C)$ of the complex $N\times N$ matrices is a $C^*$-algebra, with the usual matrix norm and involution of matrices, namely:
$$||M||=\sup_{||x||=1}||Mx||\quad,\quad 
(M^*)_{ij}=\bar{M}_{ji}$$

More generally, any $*$-subalgebra $A\subset M_N(\mathbb C)$ is automatically closed, and so is a $C^*$-algebra. In fact, in finite dimensions, the situation is as follows:

\begin{theorem}
The finite dimensional $C^*$-algebras are exactly the algebras
$$A=M_{n_1}(\mathbb C)\oplus\ldots\oplus M_{n_k}(\mathbb C)$$
with norm $||(a_1,\ldots,a_k)||=\sup_i||a_i||$, and involution $(a_1,\ldots,a_k)^*=(a_1^*,\ldots,a_k^*)$.
\end{theorem}

\begin{proof}
In one sense this is clear. In the other sense, this comes by splitting the unit of our algebra $A$ as a sum of central minimal projections, $1=p_1+\ldots+p_k$. Indeed, when doing so, each of the $*$-algebras $A_i=p_iAp_i$ follows to be a matrix algebra, $A_i\simeq M_{n_i}(\mathbb C)$, and this gives the direct sum decomposition in the statement.
\end{proof}

As a generalization now of the basic $A=M_N(\mathbb C)$ example, the algebra $A=B(H)$ formed by the bounded linear operators $T:H\to H$ over a complex Hilbert space $H$ is a $C^*$-algebra, with the usual norm and involution of linear operators, given by:
$$||T||=\sup_{||x||=1}||Tx||\quad,\quad 
<T^*x,y>=<x,Ty>$$

More generally, any normed closed $*$-algebra $A\subset B(H)$ is a $C^*$-algebra. Now inspired by this, we can develop some spectral theory for the $C^*$-algebras, as follows:

\begin{proposition}
Given an element $a\in A$ of a $C^*$-algebra, define its spectrum as:
$$\sigma(a)=\left\{\lambda\in\mathbb C\Big|a-\lambda\notin A^{-1}\right\}$$
The following spectral theory results hold, exactly as in the $A=B(H)$ case:
\begin{enumerate}
\item We have $\sigma(ab)\cup\{0\}=\sigma(ba)\cup\{0\}$.

\item We have $\sigma(f(a))=f(\sigma(a))$, for any $f\in\mathbb C(X)$ having poles outside $\sigma(a)$.

\item The spectrum $\sigma(a)$ is compact, non-empty, and contained in $D_0(||a||)$.

\item The spectra of unitaries $(u^*=u^{-1})$ and self-adjoints $(a=a^*)$ are on $\mathbb T,\mathbb R$.

\item The spectral radius of normal elements $(aa^*=a^*a)$ is given by $\rho(a)=||a||$.
\end{enumerate}
In addition, assuming $a\in A\subset B$, the spectra of $a$ with respect to $A$ and to $B$ coincide.
\end{proposition}

\begin{proof}
Here the assertions (1-5), which are formulated a bit informally, are well-known for the full operator algebra $A=B(H)$, and the proof in general is similar:

\medskip

(1) Assuming that $1-ab$ is invertible, with inverse $c$, we have $abc=cab=c-1$, and it follows that $1-ba$ is invertible too, with inverse $1+bca$. Thus $\sigma(ab),\sigma(ba)$ agree on $1\in\mathbb C$, and by linearity, it follows that $\sigma(ab),\sigma(ba)$ agree on any point $\lambda\in\mathbb C^*$.

\medskip

(2) The formula $\sigma(f(a))=f(\sigma(a))$ is clear for polynomials, $f\in\mathbb C[X]$, by factorizing $f-\lambda$, with $\lambda\in\mathbb C$. Then, the extension to the rational functions is straightforward, because $P(a)/Q(a)-\lambda$ is invertible precisely when $P(a)-\lambda Q(a)$ is.

\medskip

(3) By using $1/(1-b)=1+b+b^2+\ldots$ for $||b||<1$ we obtain that $a-\lambda$ is invertible for $|\lambda|>||a||$, and so $\sigma(a)\subset D_0(||a||)$. It is also clear that $\sigma(a)$ is closed, so what we have is a compact set. Finally, assuming $\sigma(a)=\emptyset$ the function $f(\lambda)=\varphi((a-\lambda)^{-1})$ is well-defined, for any $\varphi\in A^*$, and by Liouville we get $f=0$, contradiction.

\medskip

(4) Assuming $u^*=u^{-1}$ we have $||u||=1$, and so $\sigma(u)\subset D_0(1)$. But with $f(z)=z^{-1}$ we obtain via (2) that we have as well $\sigma(u)\subset f(D_0(1))$, and this gives $\sigma(u)\subset\mathbb T$. As for the result regarding the self-adjoints, this can be obtained from the result for the unitaries, by using (2) with functions of type $f(z)=(z+it)/(z-it)$, with $t\in\mathbb R$.

\medskip

(5) It is routine to check, by integrating quantities of type $z^n/(z-a)$ over circles centered at the origin, and estimating, that the spectral radius is given by $\rho(a)=\lim||a^n||^{1/n}$. But in the self-adjoint case, $a=a^*$, this gives $\rho(a)=||a||$, by using exponents of type $n=2^k$, and then the extension to the general normal case is straightforward.

\medskip 

(6) Regarding now the last assertion, the inclusion $\sigma_B(a)\subset\sigma_A(a)$ is clear. For the converse, assume $a-\lambda\in B^{-1}$, and set $b=(a-\lambda )^*(a-\lambda )$. We have then:
$$\sigma_A(b)-\sigma_B(b)=\left\{\mu\in\mathbb C-\sigma_B(b)\Big|(b-\mu)^{-1}\in B-A\right\}$$

Thus this difference in an open subset of $\mathbb C$. On the other hand $b$ being self-adjoint, its two spectra are both real, and so is their difference. Thus the two spectra of $b$ are equal, and in particular $b$ is invertible in $A$, and so $a-\lambda\in A^{-1}$, as desired.
\end{proof}

We can now a prove a key result regarding the $C^*$-algebras, as follows:

\begin{theorem}[Gelfand]
If $X$ is a compact space,  the algebra $C(X)$ of continuous functions on it $f:X\to\mathbb C$ is a $C^*$-algebra, with usual norm and involution, namely:
$$||f||=\sup_{x\in X}|f(x)|\quad,\quad 
f^*(x)=\overline{f(x)}$$
Conversely, any commutative $C^*$-algebra is of this form, $A=C(X)$, with 
$$X=\Big\{\chi:A\to\mathbb C\ ,\ {\rm normed\ algebra\ character}\Big\}$$
with topology making continuous the evaluation maps $ev_a:\chi\to\chi(a)$.
\end{theorem}

\begin{proof}
There are several things going on here, the idea being as follows:

\medskip

(1) The first assertion is clear from definitions. Observe that we have indeed:
$$||ff^*||
=\sup_{x\in X}|f(x)|^2
=||f||^2$$

Observe also that the algebra $C(X)$ is commutative, because $fg=gf$.

\medskip

(2) Conversely, given a commutative $C^*$-algebra $A$, let us define $X$ as in the statement. Then $X$ is compact, and $a\to ev_a$ is a morphism of algebras, as follows:
$$ev:A\to C(X)$$

(3) We first prove that $ev$ is involutive. We use the following formula, which is similar to the $z=Re(z)+iIm(z)$ decomposition formula for usual complex numbers:
$$a=\frac{a+a^*}{2}+i\cdot\frac{a-a^*}{2i}$$

Thus it is enough to prove $ev_{a^*}=ev_a^*$ for the self-adjoint elements $a$. But this is the same as proving that $a=a^*$ implies that $ev_a$ is a real function, which is in turn true, by Proposition 9.26, because $ev_a(\chi)=\chi(a)$ is an element of $\sigma(a)$, contained in $\mathbb R$.

\medskip

(4) Since $A$ is commutative, each element is normal, so $ev$ is isometric:
$$||ev_a||
=\rho(a)
=||a||$$

It remains to prove that $ev$ is surjective. But this follows from the Stone-Weierstrass theorem, because $ev(A)$ is a closed subalgebra of $C(X)$, which separates the points.
\end{proof}

As a second main result now about the $C^*$-algebras, we have:

\begin{theorem}
Any $C^*$-algebra appears as an operator algebra:
$$A\subset B(H)$$
Moreover, when $A$ is separable, which is usually the case, $H$ can be taken separable.
\end{theorem}

\begin{proof}
This result, called GNS representation theorem after Gelfand-Naimark-Segal, comes as a continuation of the Gelfand theorem, the idea being as follows:

\medskip

(1) Let us first prove that the result holds in the commutative case, $A=C(X)$. Here, we can pick a positive measure on $X$, and construct our embedding as follows:
$$C(X)\subset B(L^2(X))\quad,\quad f\to[g\to fg]$$

(2) In general the proof is similar, the idea being that given a $C^*$-algebra $A$ we can construct a Hilbert space $H=L^2(A)$, and then an embedding as above:
$$A\subset B(L^2(A))\quad,\quad a\to[b\to ab]$$

(3) Finally, the last assertion is clear, because when $A$ is separable, meaning that it has a countable algebraic basis, so does the associated Hilbert space $H=L^2(A)$.
\end{proof}

Many more things can be said about the $C^*$-algebras, and for more on this, we refer to any operator algebra book. Also, and in answer to a question that you surely might have, right now, in what regards the Hopf algebra axioms, these can be indeed upgraded, from $*$-algebras to $C^*$-algebras. However, all this is a bit technical, especially in what regards the tensor products, and we will discuss this later, in chapter 10 below.

\section*{9d. Symmetry groups}

Getting back now to the material from chapter 5, at the level of basic examples of quantum groups, we have the orthogonal and unitary groups, and their free versions, and more generally, their liberations and twists. The study here is quite standard, and for more on this, we refer to the Hopf algebra literature, including Bichon's paper \cite{bi3}.

\bigskip

In what follows, we would like to discuss a key related question, namely the relation between the quantum orthogonal and symplectic groups. As a starting point here, and by slightly upgrading our formalism, by using $C^*$-algebras as above, we have:

\begin{theorem}
The algebra of continuous functions on $SU_2$ appears as
$$C(SU_2)=C^*\left((u_{ij})_{i,j=1,2}\Big|u=J\bar{u}J^{-1}={\rm unitary}\right)$$
where $J$ is the following matrix,
$$J=\begin{pmatrix}0&1\\ -1&0\end{pmatrix}$$
called super-identity matrix. 
\end{theorem}

\begin{proof}
This can be done in several steps, as follows:

\medskip

(1) Let us first compute $SU_2$. Consider an arbitrary $2\times2$ complex matrix:
$$U=\begin{pmatrix}a&b\\c&d\end{pmatrix}$$

Assuming $\det U=1$, the unitarity condition $U^{-1}=U^*$ reads:
$$\begin{pmatrix}d&-b\\-c&a\end{pmatrix}
=\begin{pmatrix}\bar{a}&\bar{c}\\\bar{b}&\bar{d}\end{pmatrix}$$

Thus we must have $d=\bar{a}$, $c=-\bar{b}$, and we obtain the following formula:
$$SU_2=\left\{\begin{pmatrix}a&b\\-\bar{b}&\bar{a}\end{pmatrix}\Big|\ |a|^2+|b|^2=1\right\}$$

(2) Now observe that we have the following equality:
$$\begin{pmatrix}a&b\\ -\bar{b}&\bar{a}\end{pmatrix}
\begin{pmatrix}0&1\\ -1&0\end{pmatrix}
=\begin{pmatrix}-b&a\\-\bar{a}&-\bar{b}\end{pmatrix}
=\begin{pmatrix}0&1\\ -1&0\end{pmatrix}
\begin{pmatrix}\bar{a}&\bar{b}\\ -b&a\end{pmatrix}$$

Thus, with $J$ being as in the statement, we have $uJ=J\bar{u}$, and so:
$$u=J\bar{u}J^{-1}$$

We conclude that, if $A$ is the universal algebra in the statement, we have:
$$A\to C(SU_2)$$

(3) Conversely now, let us compute the universal algebra $A$ in the statement. For this purpose, let us write its fundamental corepresentation as follows:
$$u=\begin{pmatrix}a&b\\c&d\end{pmatrix}$$

We have $uJ=J\bar{u}$, with these quantities being respectively given by:
$$uJ=\begin{pmatrix}a&b\\c&d\end{pmatrix}
\begin{pmatrix}0&1\\ -1&0\end{pmatrix}
=\begin{pmatrix}-b&a\\-d&c\end{pmatrix}$$
$$J\bar{u}=\begin{pmatrix}0&1\\ -1&0\end{pmatrix}
\begin{pmatrix}a^*&b^*\\c^*&d^*\end{pmatrix}
=\begin{pmatrix}c^*&d^*\\-a^*&-b^*\end{pmatrix}$$

Thus we must have $d=a^*$, $c=-b^*$, and we obtain the following formula:
$$u=\begin{pmatrix}a&b\\-b^*&a^*\end{pmatrix}$$

We also know that this matrix must be unitary, and we have:
$$uu^*=\begin{pmatrix}a&b\\-b^*&a^*\end{pmatrix}
\begin{pmatrix}a^*&-b\\b^*&a\end{pmatrix}
=\begin{pmatrix}aa^*+bb^*&ba-ab\\a^*b^*-b^*a^*&a^*a+b^*b\end{pmatrix}$$
$$u^*u=\begin{pmatrix}a^*&-b\\b^*&a\end{pmatrix}
\begin{pmatrix}a&b\\-b^*&a^*\end{pmatrix}
=\begin{pmatrix}a^*a+bb^*&a^*b-ba^*\\b^*a-ab^*&aa^*+b^*b\end{pmatrix}$$

Thus, the unitarity equations for $u$ are as follows:
$$aa^*=a^*a=1-bb^*=1-b^*b$$
$$ab=ba,a^*b=ba^*,ab^*=a^*b,a^*b^*=b^*a^*$$

It follows that $a,b,a^*,b^*$ commute, so our algebra is commutative. Now since this algebra is commutative, the involution $*$ becomes the usual conjugation $-$, and so:
$$u=\begin{pmatrix}a&b\\-\bar{b}&\bar{a}\end{pmatrix}$$

But this tells us that we have $A=C(X)$ with $X\subset SU_2$, and so we have a quotient map $C(SU_2)\to A$, which is inverse to the map constructed in (2), as desired.
\end{proof}

Now with the above result in hand, we can see right away the relation with $O_N^+$, and more specifically with $O_2^+$. Indeed, this latter quantum group appears as follows:
$$C(O_2^+)=C^*\left((u_{ij})_{i,j=1,2}\Big|u=\bar{u}={\rm unitary}\right)$$

Thus, $SU_2$ appears from $O_2^+$ by replacing the identity with the super-identity, or perhaps vice versa. Which looks like something interesting, worth some more study.

\bigskip

In order to better understand now what happens in the orthogonal case, as a continuation of the material from chapter 5, following \cite{ba1}, \cite{bdv}, we first have:

\begin{proposition}
Given a closed subgroup $G\subset U_N^+$, with irreducible fundamental corepresentation $u=(u_{ij})$, this corepresentation is self-adjoint, $u\sim\bar{u}$, precisely when 
$$u=J\bar{u}J^{-1}$$
for some unitary matrix $J\in U_N$, satisfying the following condition:
$$J\bar{J}=\pm 1$$
Moreover, when $N$ is odd we must have $J\bar{J}=1$. 
\end{proposition}

\begin{proof}
Since $u$ is self-adjoint, $u\sim\bar{u}$, we must have $u=J\bar{u}J^{-1}$, for a certain matrix $J\in GL_N(\mathbb C)$. We obtain from this, by using our assumption that $u$ is irreducible:
\begin{eqnarray*}
u=J\bar{u}J^{-1}
&\implies&\bar{u}=\bar{J}u\bar{J}^{-1}\\
&\implies&u=(J\bar{J})u(J\bar{J})^{-1}\\
&\implies&J\bar{J}=c1\\
&\implies&\bar{J}J=\bar{c}1\\
&\implies&c\in\mathbb R
\end{eqnarray*}

Now by rescaling we can assume $c=\pm1$, so we have proved so far that:
$$J\bar{J}=\pm 1$$

In order to establish now the formula $JJ^*=1$, we can proceed as follows:
\begin{eqnarray*}
(id\otimes S)u=u^*
&\implies&(id\otimes S)\bar{u}=u^t\\
&\implies&(id\otimes S)(J\bar{u}J^{-1})=Ju^tJ^{-1}\\
&\implies&u^*=Ju^tJ^{-1}\\
&\implies&u=(J^*)^{-1}\bar{u}J^*\\
&\implies&\bar{u}=J^*u(J^*)^{-1}\\
&\implies&\bar{u}=J^*J\bar{u}J^{-1}(J^*)^{-1}\\
&\implies&JJ^*=d1
\end{eqnarray*}

We have $JJ^*>0$, so $d>0$. On the other hand, from $J\bar{J}=\pm 1$, $JJ^*=d1$ we get:
$$|\det J|^2=\det(J\bar{J})=(\pm1)^N$$
$$|\det J|^2=\det(JJ^*)=d^N$$

Since $d>0$ we obtain from this $d=1$, and so $JJ^*=1$ as claimed. We obtain as well that when $N$ is odd the sign must be 1, and so $J\bar{J}=1$, as claimed.
\end{proof}

It is convenient to diagonalize $J$. Once again following Bichon-De Rijdt-Vaes \cite{bdv}, up to an orthogonal base change, we can assume that our matrix is as follows, where $N=2p+q$ and $\varepsilon=\pm 1$, with the $1_q$ block at right disappearing if $\varepsilon=-1$:
$$J=\begin{pmatrix}
0&1\ \ \ \\
\varepsilon 1&0_{(0)}\\
&&\ddots\\
&&&0&1\ \ \ \\
&&&\varepsilon 1&0_{(p)}\\
&&&&&1_{(1)}\\
&&&&&&\ddots\\
&&&&&&&1_{(q)}
\end{pmatrix}$$

To be more precise, in the case $\varepsilon=1$, the super-identity is the following matrix:
$$J=\begin{pmatrix}
0&1\ \ \ \\
1&0_{(1)}\\
&&\ddots\\
&&&0&1\ \ \ \\
&&&1&0_{(p)}\\
&&&&&1_{(1)}\\
&&&&&&\ddots\\
&&&&&&&1_{(q)}
\end{pmatrix}$$

In the case $\varepsilon=-1$ now, the diagonal terms vanish, and the super-identity is:
$$J=\begin{pmatrix}
0&1\ \ \ \\
-1&0_{(1)}\\
&&\ddots\\
&&&0&1\ \ \ \\
&&&-1&0_{(p)}
\end{pmatrix}$$

We are therefore led into the following definition:

\index{super-space}
\index{super-identity}

\begin{definition}
The ``super-space'' $\mathbb C^N_J$ is the usual space $\mathbb C^N$, with its standard basis $\{e_1,\ldots,e_N\}$, with a chosen sign $\varepsilon=\pm 1$, and a chosen involution on the set of indices,
$$i\to\bar{i}$$
with $J$ being the ``super-identity'' matrix, $J_{ij}=\delta_{i\bar{j}}$ for $i\leq j$ and $J_{ij}=\varepsilon\delta_{i\bar{j}}$ for $i\geq j$.
\end{definition}

In what follows we will usually assume that $J$ is the explicit matrix appearing in the above. Indeed, up to a permutation of the indices, we have a decomposition $n=2p+q$ such that the involution is, in standard permutation notation, given by:
$$(12)\ldots (2p-1,2p)(2p+1)\ldots (q)$$

Now let us construct some basic compact quantum groups, orthogonal and symplectic, in our ``super'' setting. We can indeed formulate the following definition:

\index{super-orthogonal group}
\index{super-orthogonal quantum group}
\index{symplectic group}
\index{free symplectic group}

\begin{definition}
Associated to the super-space $\mathbb C^N_J$ are the following objects,
\begin{enumerate}
\item The super-orthogonal group, given by:
$$O_J=\left\{U\in U_N\Big|U=J\bar{U}J^{-1}\right\}$$

\item The super-orthogonal quantum group, given by:
$$C(O_J^+)=C\left((u_{ij})_{i,j=1,\ldots,n}\Big|u=J\bar{u}J^{-1}={\rm unitary}\right)$$
\end{enumerate}
with $J$ standing as usual for the super-identity matrix.
\end{definition}

It it possible to considerably extend this list, but for our purposes here, this is what we need for the moment. We have indeed the following result:

\begin{theorem}
The basic orthogonal groups and quantum groups are as follows:
\begin{enumerate}
\item At $\varepsilon=-1$ we have $O_J=Sp_N$ and $O_J^+=Sp_N^+$.

\item At $\varepsilon=-1$ and $N=2$ we have $O_J=O_J^+=SU_2$.

\item At $\varepsilon=1$ we have $O_J=O_N$ and $O_J^+=O_N^+$.
\end{enumerate}
\end{theorem}

\begin{proof}
These results are all elementary, as follows:

\medskip

(1) At $\varepsilon=-1$ this follows from definitions, because $Sp_N\subset U_N$ is given by:
$$Sp_N=\left\{U\in U_N\Big|U=J\bar{U}J^{-1}\right\}$$

(2) Still at $\varepsilon=-1$, the equation $U=J\bar{U}J^{-1}$ tells us that the symplectic matrices $U\in Sp_N$ are exactly the unitaries $U\in U_N$ which are patterned as follows:
$$U=\begin{pmatrix}
a&b&\ldots\\
-\bar{b}&\bar{a}\\
\vdots&&\ddots
\end{pmatrix}$$

In particular, the symplectic matrices at $N=2$ are as follows:
$$U=\begin{pmatrix}
a&b\\
-\bar{b}&\bar{a}
\end{pmatrix}$$

Thus we have $Sp_2=U_2$, and the formula $Sp_2^+=Sp_2$ is elementary as well.

\medskip

(3) At $\varepsilon=1$ now, consider the root of unity $\rho=e^{\pi i/4}$, and set:
$$K=\frac{1}{\sqrt{2}}\begin{pmatrix}\rho&\rho^7\\ \rho^3&\rho^5\end{pmatrix}$$

This matrix $K$ is then unitary, and we have:
$$K\begin{pmatrix}0&1\\1&0\end{pmatrix}K^t=1$$

Thus the following matrix is unitary as well, and satisfies $LJL^t=1$:
$$L=\begin{pmatrix}K^{(1)}\\&\ddots\\&&K^{(p)}\\&&&1_q\end{pmatrix}$$

Thus in terms of the matrix $V=LUL^*$ we have:
$$U=J\bar{U}J^{-1}={\rm unitary}
\quad\iff\quad V=\bar{V}={\rm unitary}$$

We obtain in this way an isomorphism $O_J^+=O_N^+$ as in the statement, and by passing to classical versions, we obtain as well $O_J=O_N$, as desired.
\end{proof}

As a last topic for this chapter, let us discuss now the quantum permutation and quantum reflection groups. For dealing with quantum permutations, we will need:

\begin{proposition}
If $p_1,\ldots,p_N$ are projections, $p^2=p^*=p$, satisfying
$$p_1+\ldots+p_N=1$$
then these projections are pairwise orthogonal, $p_ip_j=0$ for $i\neq j$.
\end{proposition}

\begin{proof}
We have indeed the following computation, using our assumptions:
\begin{eqnarray*}
p_i
&=&p_i\cdot1\cdot p_i\\
&=&p_i\left(\sum_jp_j\right)p_i\\
&=&p_i+\sum_{i\neq j}p_ip_jp_i\\
&=&p_i+\sum_{i\neq j}(p_ip_j)(p_ip_j)^*
\end{eqnarray*}

Thus by positivity we get $p_ip_j=0$ for any $i\neq j$, as desired.
\end{proof}

Now by using the above result, we can simplify the definition of $S_N^+$, as follows:

\begin{theorem}
The quantum permutation group $S_N^+$ is given by the formula
$$C(S_N^+)=C^*\left((u_{ij})_{i,j=1,\ldots,N}\Big|u={\rm magic}\right)$$
with magic meaning formed of projections, summing up to $1$ on each row and column.
\end{theorem}

\begin{proof}
This is something that we know well since chapter 6, with the only change coming from the exact magic conditions which are used, which are precisely those from chapter 6, altered by our findings from Proposition 9.34, which simplify them a bit.
\end{proof}

Getting now to the quantum reflection groups, here we have some simplifications too. It is actually convenient to rewrite here the whole theory, by starting with:

\begin{definition}
The algebra $C(H_N^{s+})$ is the universal $C^*$-algebra generated by $N^2$ normal elements $u_{ij}$, subject to the following relations,
\begin{enumerate}
\item $u=(u_{ij})$ is unitary,

\item $u^t=(u_{ji})$ is unitary,

\item $p_{ij}=u_{ij}u_{ij}^*$ is a projection,

\item $u_{ij}^s=p_{ij}$,
\end{enumerate}
with Hopf algebra maps $\Delta,\varepsilon,S$ constructed by universality.
\end{definition}

Here we allow the value $s=\infty$, with the convention that the last axiom simply disappears in this case. Observe that at $s<\infty$ the normality condition is actually redundant, because a partial isometry $a$ subject to the relation $aa^*=a^s$ is normal.

\bigskip

Let us prove now that $H_N^{s+}$ with $s<\infty$ is a quantum permutation group. For this purpose, we must change the fundamental representation. Let us start with:

\index{sudoku matrix}

\begin{definition}
A $(s,N)$-sudoku matrix is a magic unitary of size $sN$, of the form
$$m=\begin{pmatrix}
a^0&a^1&\ldots&a^{s-1}\\
a^{s-1}&a^0&\ldots&a^{s-2}\\
\vdots&\vdots&&\vdots\\
a^1&a^2&\ldots&a^0
\end{pmatrix}$$
where $a^0,\ldots,a^{s-1}$ are $N\times N$ matrices.
\end{definition}

Now with this notion in hand, we have the following result, which shows that the quantum groups constructed in Definition 9.36 coincide with those from chapter 6:

\begin{theorem}
The following happen:
\begin{enumerate}
\item The algebra $C(H_N^s)$ is isomorphic to the universal commutative $C^*$-algebra generated by the entries of a $(s,N)$-sudoku matrix.

\item The algebra $C(H_N^{s+})$ is isomorphic to the universal $C^*$-algebra generated by the entries of a $(s,N)$-sudoku matrix.
\end{enumerate}
\end{theorem}

\begin{proof}
The first assertion follows from the second one. In order to prove now the second assertion, consider the universal algebra in the statement, namely:
$$A=C^*\left(a_{ij}^p\ \Big\vert \left(a^{q-p}_{ij}\right)_{pi,qj}=(s,N)-\mbox{sudoku }\right)$$

Consider also the algebra $C(H_N^{s+})$. According to Definition 9.36, this is presented by certain relations $R$, that we will call here level $s$ cubic conditions:
$$C(H_N^{s+})=C^*\left(u_{ij}\ \Big\vert\  u=N\times N\mbox{ level $s$ cubic }\right)$$

We can construct a pair of inverse morphisms between these algebras, as follows:

\medskip

(1) Our first claim is that $U_{ij}=\sum_pw^{-p}a^p_{ij}$ is a level $s$ cubic unitary. Indeed, by using the sudoku condition, the verification of (1-4) in Definition 9.36 is routine.

\medskip

(2) Our second claim is that the elements $A^p_{ij}=\frac{1}{s}\sum_rw^{rp}u^r_{ij}$, with the convention $u_{ij}^0=p_{ij}$, form a level $s$ sudoku unitary. Once again, the proof here is routine.

\medskip

(3) According to the above, we can define a morphism $\Phi:C(H_N^{s+})\to A$ by the formula $\Phi(u_{ij})=U_{ij}$, and a morphism $\Psi:A\to C(H_N^{s+})$ by the formula $\Psi(a^p_{ij})=A^p_{ij}$.

\medskip

(4) We check now the fact that $\Phi,\Psi$ are indeed inverse morphisms:
\begin{eqnarray*}
\Psi\Phi(u_{ij})
&=&\sum_pw^{-p}A^p_{ij}\\
&=&\frac{1}{s}\sum_pw^{-p}\sum_rw^{rp}u_{ij}^r\\
&=&\frac{1}{s}\sum_{pr}w^{(r-1)p}u_{ij}^r\\
&=&u_{ij}
\end{eqnarray*}

(5) As for the other composition, we have the following computation:
\begin{eqnarray*}
\Phi\Psi(a^p_{ij})
&=&\frac{1}{s}\sum_rw^{rp}U_{ij}^r\\
&=&\frac{1}{s}\sum_rw^{rp}\sum_qw^{-rq}a_{ij}^q\\
&=&\frac{1}{s}\sum_qa_{ij}^q\sum_rw^{r(p-q)}\\
&=&a^p_{ij}
\end{eqnarray*}

Thus we have an isomorphism $C(H_N^{s+})=A$, as claimed.
\end{proof}

Many other things can be said, as a continuation of the above, and we refer here to the quantum permutation group literature. We will be back to this.

\section*{9e. Exercises}

This was a quite elementary chapter, and as exercises here, we have:

\begin{exercise}
What happens to the above, over a field $F$ with $char(F)\neq0$?
\end{exercise}

\begin{exercise}
What happens to the above, over a field $F$ not algebrically closed?
\end{exercise}

\begin{exercise}
What happens to the above, over a field $F$ having no involution?
\end{exercise}

\begin{exercise}
All in all, can the above be extended, to other fields $F\neq\mathbb C$?
\end{exercise}

\begin{exercise}
Learn more about tensor categories, over $\mathbb C$ and other fields.
\end{exercise}

\begin{exercise}
Learn more about symplectic groups, and their various properties.
\end{exercise}

As bonus exercise, in relation with this, learn about the field of $p$-adic numbers.

\chapter{Quantum groups}

\section*{10a. Quantum groups}

We have been talking about quantum groups since the beginning of this book, but these were quite often speculatory objects, rather mentioned in order to get a better understanding of the Hopf algebra theory that we were developing.

\bigskip

Time now to come up with a more precise definition for them, or at least to come up with a precise definition for a large class of quantum groups, which covers many interesting examples. The idea is very simple, based on Gelfand duality, as follows:

\begin{principle}
According to Gelfand duality, any $C^*$-algebra can be thought of as being of the form $A=C(X)$, with $X$ being a compact quantum space. Thus, a Hopf $C^*$-algebra should be of the form $A=C(G)$, with $G$ being a compact quantum group.
\end{principle}

In practice now, formulating the Hopf algebra axioms in the $C^*$-algebra context is not exactly an easy task, due to various topological issues. So, we will choose here to trick, by restricting the attention to the compact Lie group case, and by following:

\begin{principle}[update]
The compact Lie groups, $G\subset U_N$, have their multiplication, unit and inverse operations given by the formulae
$$(UV)_{ij}=\sum_kU_{ik}V_{kj}\quad,\quad 
(1_N)_{ij}=\delta_{ij}\quad,\quad 
(U^{-1})_{ij}=U_{ji}^*$$
and based on this, we should axiomatize the algebras $A=C(G)$, and then lift the commutativity assumption on $A$, as to have axioms for the compact quantum Lie groups.
\end{principle}

In short, we are proposing here to forget about the Hopf algebra axioms, or rather leave these axioms for later, and follow instead the above straightforward approach. As we will see in a moment, this is something that will majestically work, but in case you still have doubts about our method, here is for you, some more philosophy:

\begin{comment}
The groups appeared in mathematics as transformation groups, 
$$G\subset U_N$$
and abstract groups came afterwards. The same should happen for quantum groups.
\end{comment}

But probably too much philosophy, let us get to work, following Principle 10.2. Getting directly to the answer, and with the Gelfand duality details to be explained in a moment, we are led in this way to the following key definition, due to Woronowicz \cite{wo1}:

\index{Woronowicz algebra}
\index{comultiplication}
\index{counit}
\index{antipode}

\begin{definition}
A Woronowicz algebra is a $C^*$-algebra $A$, given with a unitary matrix $u\in M_N(A)$ whose coefficients generate $A$, such that the formulae
$$\Delta(u_{ij})=\sum_ku_{ik}\otimes u_{kj}\quad,\quad
\varepsilon(u_{ij})=\delta_{ij}\quad,\quad 
S(u_{ij})=u_{ji}^*$$
define morphisms of $C^*$-algebras as follows,
$$\Delta:A\to A\otimes A\quad,\quad 
\varepsilon:A\to\mathbb C\quad,\quad 
S:A\to A^{opp}$$
called comultiplication, counit and antipode.
\end{definition}

Obviously, this is something tricky, that will take us some time, to understand. Here are some quick comments in relation with this, and more on everything later:

\bigskip

(1) To start with, the above formulae of $\Delta,\varepsilon,S$ are straightforward, obtained from those of $m,u,i$ from Principle 10.2, by transposition. So, at least we know one thing.

\bigskip

(2) Regarding the tensor product $\otimes$ appearing above, that can be any topological tensor product, and with the choice here being irrelevant. More on this later. 

\bigskip

(3) Finally, as we will see in a moment, the above axioms force $S^2=id$, in analogy with $i^2=id$, happening in the context of Principle 10.2. More on this later too.

\bigskip

We say that $A$ is cocommutative when $\Sigma\Delta=\Delta$, where $\Sigma(a\otimes b)=b\otimes a$ is the flip. With this convention, we have the following key result, from Woronowicz \cite{wo1}:

\index{cocommutative algebra}
\index{commutative algebra}

\begin{theorem}
The following are Woronowicz algebras:
\begin{enumerate}
\item $C(G)$, with $G\subset U_N$ compact Lie group. Here the structural maps are:
$$\Delta(\varphi)=(g,h)\to \varphi(gh)$$
$$\varepsilon(\varphi)=\varphi(1)$$
$$S(\varphi)=g\to\varphi(g^{-1})$$

\item $C^*(\Gamma)$, with $F_N\to\Gamma$ finitely generated group. Here the structural maps are:
$$\Delta(g)=g\otimes g$$
$$\varepsilon(g)=1$$
$$S(g)=g^{-1}$$

\end{enumerate}
Moreover, we obtain in this way all the commutative/cocommutative algebras.
\end{theorem}

\begin{proof}
Everything here, or almost, is elementary, the idea being as follows:

\medskip

(1) Given a closed subgroup $G\subset U_N$, consider the matrix $u=(u_{ij})$ formed by the standard coordinates, $u_{ij}(g)=g_{ij}$. That is, define $u=(u_{ij})$ by the following formula:
$$g=\begin{pmatrix}
u_{11}(g)&\ldots&u_{1N}(g)\\
\vdots&&\vdots\\
u_{N1}(g)&\ldots&u_{NN}(g)
\end{pmatrix}$$

Let us check now the conditions in Definition 10.4, with the comultiplication, counit and antipode being the maps in the statement, appearing as transposes, as follows:
$$\Delta=m^t\quad,\quad\varepsilon=u^t\quad,\quad S=i^t$$

To strat with, the fact that the transpose of the multiplication map $m^t$ satisfies the condition in Definition 10.4 is clear, coming from the following computation:
\begin{eqnarray*}
m^t(u_{ij})(U\otimes V)
&=&(UV)_{ij}\\
&=&\sum_kU_{ik}V_{kj}\\
&=&\sum_k(u_{ik}\otimes u_{kj})(U\otimes V)
\end{eqnarray*}

Regarding now the transpose of the unit map $u^t$, the verification of the condition in Definition 10.4 is trivial, coming from the following equalities:
$$u^t(u_{ij})
=1_{ij}
=\delta_{ij}$$

Finally, the transpose of the inversion map $i^t$ verifies the condition in Definition 10.4, because we have the following computation, valid for any $U\in G$:
$$i^t(u_{ij})(U)
=(U^{-1})_{ij}
=\bar{U}_{ji}
=u_{ji}^*(U)$$

Summarizing, main assertion proved. Regarding now the converse, assume that a Woronowicz algebra $A$ is commutative. By using the Gelfand theorem, we can write $A=C(G)$, with $G$ being a certain compact space. By using now the coordinates $u_{ij}$, we obtain an embedding $G\subset U_N$. Finally, by using $\Delta,\varepsilon,S$, it follows that the subspace $G\subset U_N$ that we have obtained is in fact a closed subgroup, and we are done.

\medskip

(2) Given a finitely generated group $\Gamma=<g_1,\ldots,g_N>$, as in the statement, consider the diagonal matrix $u=(u_{ij})$ formed by the generators of $\Gamma$:
$$u=\begin{pmatrix}
g_1&&0\\
&\ddots&\\
0&&g_N
\end{pmatrix}$$

Let us verify now the axioms in Definition 10.4. In order to deal with the first axiom, consider the following map, which is a unitary representation:
$$\Gamma\to C^*(\Gamma)\otimes C^*(\Gamma)\quad,\quad 
g\to g\otimes g$$

This representation extends, as desired, into a morphism of algebras, as follows:
$$\Delta:C^*(\Gamma)\to C^*(\Gamma)\otimes C^*(\Gamma)\quad,\quad 
\Delta(g)=g\otimes g$$

The situation for $\varepsilon$ is similar, because this comes from the trivial representation:
$$\Gamma\to\{1\}\quad,\quad 
g\to1$$

Finally, the antipode $S$ comes from the following unitary representation:
$$\Gamma\to C^*(\Gamma)^{opp}\quad,\quad 
g\to g^{-1}$$

Summarizing, we have shown that we have a Woronowicz algebra, with $\Delta,\varepsilon,S$ being as in the statement. Regarding now the last assertion, observe that we have:
$$\Sigma\Delta(g)
=\Sigma(g\otimes g)
=g\otimes g
=\Delta(g)$$

Thus the formula $\Sigma\Delta=\Delta$ holds on the group elements $g\in\Gamma$, and by linearity and continuity, this formula must hold on the whole algebra $C^*(\Gamma)$, as desired. As for the converse, stating that any cocommutative Woronowicz algebra must be of this form, this is something more technical, that we will discuss later in this chapter.
\end{proof}

In general now, the structural maps $\Delta,\varepsilon,S$ have the following properties:

\index{square of antipode}

\begin{theorem}
Let $(A,u)$ be a Woronowicz algebra.
\begin{enumerate} 
\item $\Delta,\varepsilon$ satisfy the usual axioms for a comultiplication and a counit, namely:
$$(\Delta\otimes id)\Delta=(id\otimes \Delta)\Delta$$
$$(\varepsilon\otimes id)\Delta=(id\otimes\varepsilon)\Delta=id$$

\item $S$ satisfies the antipode axiom, on the $*$-subalgebra generated by entries of $u$: 
$$m(S\otimes id)\Delta=m(id\otimes S)\Delta=\varepsilon(.)1$$

\item In addition, the square of the antipode is the identity, $S^2=id$.
\end{enumerate}
\end{theorem}

\begin{proof}
Observe first that the result holds in the case where the algebra $A$ is commutative. Indeed, by using Theorem 10.5 (1) we can write:
$$\Delta=m^t\quad,\quad 
\varepsilon=u^t\quad,\quad
S=i^t$$

And with this done, the 3 conditions in the statement simply come then by transposition from the basic 3 group theory conditions satisfied by $m,u,i$, namely:
$$m(m\times id)=m(id\times m)$$
$$m(id\times u)=m(u\times id)=id$$
$$m(id\times i)\delta=m(i\times id)\delta=1$$

As for the last condition, $S^2=id$, this is satisfied too, coming from $i^2=id$. Now observe that the result holds as well in the case where $A$ is cocommutative, by using Theorem 10.5 (2). Indeed, the 3 formulae in the statement are all trivial, and the condition $S^2=id$ follows once again from $i^2=id$. In general now, the proof goes as follows:

\medskip

(1) We have the following computation:
$$(\Delta\otimes id)\Delta(u_{ij})
=\sum_l\Delta(u_{il})\otimes u_{lj}
=\sum_{kl}u_{ik}\otimes u_{kl}\otimes u_{lj}$$

We have as well the following computation, which gives the first formula:
$$(id\otimes\Delta)\Delta(u_{ij})
=\sum_ku_{ik}\otimes\Delta(u_{kj})
=\sum_{kl}u_{ik}\otimes u_{kl}\otimes u_{lj}$$

On the other hand, we have the following computation:
$$(id\otimes\varepsilon)\Delta(u_{ij})
=\sum_ku_{ik}\otimes\varepsilon(u_{kj})
=u_{ij}$$

We have as well the following computation, which gives the second formula:
$$(\varepsilon\otimes id)\Delta(u_{ij})
=\sum_k\varepsilon(u_{ik})\otimes u_{kj}
=u_{ij}$$

(2) By using the fact that the matrix $u=(u_{ij})$ is unitary, we obtain:
\begin{eqnarray*}
m(id\otimes S)\Delta(u_{ij})
&=&\sum_ku_{ik}S(u_{kj})\\
&=&\sum_ku_{ik}u_{jk}^*\\
&=&(uu^*)_{ij}\\
&=&\delta_{ij}
\end{eqnarray*}

We have as well the following computation, which gives the result:
\begin{eqnarray*}
m(S\otimes id)\Delta(u_{ij})
&=&\sum_kS(u_{ik})u_{kj}\\
&=&\sum_ku_{ki}^*u_{kj}\\
&=&(u^*u)_{ij}\\
&=&\delta_{ij}
\end{eqnarray*}

(3) Finally, the formula $S^2=id$ holds as well on the generators, and we are done.
\end{proof}

Still talking basics, let us record as well the following useful technical result:

\index{biunitary}

\begin{proposition}
Given a Woronowicz algebra $(A,u)$, we have 
$$u^t=\bar{u}^{-1}$$
so $u$ is biunitary, in the sense that it is unitary, with unitary transpose.
\end{proposition}

\begin{proof}
We have the following computation, based on the fact that $u$ is unitary:
\begin{eqnarray*}
(uu^*)_{ij}=\delta_{ij}
&\implies&\sum_kS(u_{ik}u_{jk}^*)=\delta_{ij}\\
&\implies&\sum_ku_{kj}u_{ki}^*=\delta_{ij}\\
&\implies&(u^t\bar{u})_{ji}=\delta_{ij}
\end{eqnarray*}

Similarly, we have the following computation, once again using the unitarity of $u$:
\begin{eqnarray*}
(u^*u)_{ij}=\delta_{ij}
&\implies&\sum_kS(u_{ki}^*u_{kj})=\delta_{ij}\\
&\implies&\sum_ku_{jk}^*u_{ik}=\delta_{ij}\\
&\implies&(\bar{u}u^t)_{ji}=\delta_{ij}
\end{eqnarray*}

Thus, we are led to the conclusion in the statement.
\end{proof}

Summarizing, the Woronowicz algebras appear to have nice properties. In view of Theorem 10.5 and of Theorem 10.6, we can formulate the following definition:

\index{quantum group}
\index{compact quantum group}
\index{discrete quantum group}
\index{Pontrjagin dual}
\index{group algebra}

\begin{definition}
Given a Woronowicz algebra $A$, we formally write
$$A=C(G)=C^*(\Gamma)$$
and call $G$ compact quantum group, and $\Gamma$ discrete quantum group.
\end{definition}

When $A$ is commutative and cocommutative, $G$ and $\Gamma$ are usual abelian groups, dual to each other. In general, we still agree to write $G=\widehat{\Gamma},\Gamma=\widehat{G}$, but in a formal sense. As a final piece of general theory now, let us complement Definitions 10.4 and 10.8 with:

\begin{definition}
Given two Woronowicz algebras $(A,u)$ and $(B,v)$, we write 
$$A\simeq B$$
and identify the corresponding quantum groups, when we have an isomorphism 
$$<u_{ij}>\simeq<v_{ij}>$$
of $*$-algebras, mapping standard coordinates to standard coordinates.
\end{definition}

With this convention, any compact or discrete quantum group corresponds to a unique Woronowicz algebra, up to equivalence. Also, we can see now why in Definition 10.4 the choice of the exact topological tensor product $\otimes$ was irrelevant. Indeed, no matter what tensor product $\otimes$ we choose to use there, we end up with the same Woronowicz algebra, and the same compact and discrete quantum groups, up to equivalence. 

\bigskip

Moving ahead now, let us call corepresentation of $A$ any unitary matrix $v\in M_n(\mathcal A)$, where $\mathcal A=<u_{ij}>$, satisfying the same conditions are those satisfied by $u$, namely:
$$\Delta(v_{ij})=\sum_kv_{ik}\otimes v_{kj}\quad,\quad 
\varepsilon(v_{ij})=\delta_{ij}\quad,\quad 
S(v_{ij})=v_{ji}^*$$

\index{representation}
\index{corepresentation}

These corepresentations can be thought of as corresponding to the finite dimensional unitary smooth representations of the underlying compact quantum group $G$. Following Woronowicz \cite{wo1}, we can now formulate a first Peter-Weyl theorem, as follows:

\index{Peter-Weyl theorem}
\index{complete reducibility}
\index{cosemisimplicity}
\index{multimatrix algebra}
\index{sum of matrix algebras}
\index{finite dimensional algebra}

\begin{theorem}[PW1]
Let $u\in M_N(A)$ be a corepresentation, consider the $*$-algebra $B=End(u)$, and write its unit as $1=p_1+\ldots+p_k$. We have then
$$u=u_1+\ldots+u_k$$
with each $u_i$ being an irreducible corepresentation, obtained by restricting $u$ to $Im(p_i)$.
\end{theorem}

\begin{proof}
This is something very classical, well-known to hold for the compact groups, and the proof in general can be established in a similar way, as follows:

\medskip

(1) We first associate to our corepresentation $u\in M_N(A)$ the corresponding coaction map $\alpha:\mathbb C^N\to\mathbb C^N\otimes A$, given by the following formula:
$$\alpha(e_i)=\sum_je_j\otimes u_{ji}$$

We say that a linear subspace $V\subset\mathbb C^N$ is invariant under $u$ if we have:
$$\alpha(V)\subset V\otimes A$$

In this case, we can consider the following restriction map:
$$\alpha_{|V}:V\to V\otimes A$$

This is a coaction map too, which must come from a subcorepresentation $v\subset u$.

\medskip

(2) Consider now a projection $p\in End(u)$. From $pu=up$ we obtain that the linear space $V=Im(p)$ is invariant under $u$, and so this space must come from a subcorepresentation $v\subset u$. It is routine to check that the operation $p\to v$ maps subprojections to subcorepresentations, and minimal projections to irreducible corepresentations.

\medskip

(3) With these preliminaries in hand, let us decompose now the $*$-algebra $End(u)$ as a multimatrix algebra, by using the decomposition of 1 into minimal projections:
$$1=p_1+\ldots+p_k$$

Consider now the following vector spaces, obtained as images of these projections:
$$V_i=Im(p_i)$$

If we denote by $u_i\subset u$ the subcorepresentations coming from these vector spaces, then we obtain in this way a decomposition $u=u_1+\ldots+u_k$, as in the statement.
\end{proof}

In order to formulate our second Peter-Weyl type theorem, we will need:

\index{Peter-Weyl representation}
\index{Peter-Weyl corepresentation}
\index{colored integer}

\begin{definition}
Given a Woronowicz algebra $(A,u)$, we denote by $u^{\otimes k}$, with $k=\circ\bullet\bullet\circ\ldots$ being a colored integer, the tensor products between $u,\bar{u}$, with the rules 
$$u^{\otimes\emptyset}=1\quad,\quad 
u^{\otimes\circ}=u\quad,\quad 
u^{\otimes\bullet}=\bar{u}$$
and multiplicativity, $u^{\otimes kl}=u^{\otimes k}\otimes u^{\otimes l}$, and call them Peter-Weyl corepresentations. 
\end{definition}

Here are a few examples of such corepresentations, namely those coming from the colored integers of length 2, to be often used in what follows:
$$u^{\otimes\circ\circ}=u\otimes u\quad,\quad 
u^{\otimes\circ\bullet}=u\otimes\bar{u}$$
$$u^{\otimes\bullet\circ}=\bar{u}\otimes u\quad,\quad 
u^{\otimes\bullet\bullet}=\bar{u}\otimes\bar{u}$$

Here is now our second Peter-Weyl theorem, adding to Theorem 10.10:

\index{Peter-Weyl theorem}

\begin{theorem}[PW2]
Each irreducible corepresentation of $A$ appears as:
$$v\subset u^{\otimes k}$$
That is, $v$ appears inside a certain Peter-Weyl corepresentation.
\end{theorem}

\begin{proof}
Given an arbitrary corepresentation $v\in M_n(A)$, consider its space of coefficients, $C(v)=span(v_{ij})$. It is routine to check that the construction $v\to C(v)$ is functorial, in the sense that it maps subcorepresentations into subspaces. By definition of the Peter-Weyl corepresentations, we have an equality as follows:
$$\mathcal A=\sum_{k\in\mathbb N*\mathbb N}C(u^{\otimes k})$$

Now given a corepresentation $v\in M_n(A)$, the corresponding coefficient space is a finite dimensional subspace $C(v)\subset\mathcal A$, and so we must have, for certain $k_1,\ldots,k_p$:
$$C(v)\subset C(u^{\otimes k_1}\oplus\ldots\oplus u^{\otimes k_p})$$

We deduce from this that we have an inclusion of corepresentations, as follows:
$$v\subset u^{\otimes k_1}\oplus\ldots\oplus u^{\otimes k_p}$$

Together with Theorem 10.10, this leads to the conclusion in the statement.
\end{proof}

\section*{10b. Haar integration}

Our goal now will be to prove that our Woronowicz algebras $A$ have a Haar integration functional, meaning a functional satisfying the following invariance condition:
$$\left(\int_A\otimes id\right)\Delta=\left(id\otimes\int_A\right)\Delta=\int_A(.)1$$

Our source of inspiration will be the classical proof of the existence of the Haar measure on the compact Lie groups $G\subset U_N$. Following Woronowicz \cite{wo1}, we first have:

\index{Ces\`aro limit}
\index{convolution}

\begin{proposition}
Given an arbitrary unital linear form $\varphi\in A^*$, the limit
$$\int_\varphi a=\lim_{n\to\infty}\frac{1}{n}\sum_{k=1}^n\varphi^{*k}(a)$$
exists, and for a coefficient of a corepresentation $a=(\tau\otimes id)v$, we have
$$\int_\varphi a=\tau(P)$$
where $P$ is the orthogonal projection onto the $1$-eigenspace of $(id\otimes\varphi)v$.
\end{proposition}

\begin{proof}
By linearity, it is enough to prove the first assertion for elements of the following type, where $v$ is a Peter-Weyl corepresentation, and $\tau$ is a linear form:
$$a=(\tau\otimes id)v$$

Thus we are led into the second assertion, and more precisely we can have the whole result proved if we can establish the following formula, with $a=(\tau\otimes id)v$:
$$\lim_{n\to\infty}\frac{1}{n}\sum_{k=1}^n\varphi^{*k}(a)=\tau(P)$$

In order to prove this latter formula, observe that we have:
$$\varphi^{*k}(a)
=(\tau\otimes\varphi^{*k})v
=\tau((id\otimes\varphi^{*k})v)$$

Consider now the following matrix, which is a usual complex matrix:
$$M=(id\otimes\varphi)v$$

In terms of this matrix, we have the following formula:
$$((id\otimes\varphi^{*k})v)_{i_0i_{k+1}}
=\sum_{i_1\ldots i_k}M_{i_0i_1}\ldots M_{i_ki_{k+1}}
=(M^k)_{i_0i_{k+1}}$$

Thus for any $k\in\mathbb N$ we have the following formula:
$$(id\otimes\varphi^{*k})v=M^k$$

It follows that our Ces\`aro limit is given by the following formula:
\begin{eqnarray*}
\lim_{n\to\infty}\frac{1}{n}\sum_{k=1}^n\varphi^{*k}(a)
&=&\lim_{n\to\infty}\frac{1}{n}\sum_{k=1}^n\tau(M^k)\\
&=&\tau\left(\lim_{n\to\infty}\frac{1}{n}\sum_{k=1}^nM^k\right)
\end{eqnarray*}

Now since the matrix $v$ is unitary we have $||v||=1$, and we conclude from this that we have $||M||\leq1$. Thus, by standard calculus, the above Ces\`aro limit on the right exists, and equals the orthogonal projection onto the $1$-eigenspace of $M$:
$$\lim_{n\to\infty}\frac{1}{n}\sum_{k=1}^nM^k=P$$

Thus our initial Ces\`aro limit converges as well, to $\tau(P)$, as desired. 
\end{proof}

When $\varphi$ is faithful, we have the following finer result, also from Woronowicz \cite{wo1}:

\begin{proposition}
Given a faithful unital linear form $\varphi\in A^*$, the limit
$$\int_\varphi a=\lim_{n\to\infty}\frac{1}{n}\sum_{k=1}^n\varphi^{*k}(a)$$
exists, and is independent of $\varphi$, given on coefficients of corepresentations by
$$\left(id\otimes\int_\varphi\right)v=P$$
where $P$ is the orthogonal projection onto $Fix(v)=\{\xi\in\mathbb C^n|v\xi=\xi\}$.
\end{proposition}

\begin{proof}
In view of Proposition 10.13, it remains to prove that when $\varphi$ is faithful, the $1$-eigenspace of $M=(id\otimes\varphi)v$ equals $Fix(v)$. But we can do this as follows:

\medskip

``$\supset$'' This is clear, and for any $\varphi$, because we have:
$$v\xi=\xi\implies M\xi=\xi$$

``$\subset$'' Here we must prove that, when $\varphi$ is faithful, we have:
$$M\xi=\xi\implies v\xi=\xi$$

For this purpose, we use a positivity trick. Consider the following element:
$$a=\sum_i\left(\sum_jv_{ij}\xi_j-\xi_i\right)\left(\sum_kv_{ik}\xi_k-\xi_i\right)^*$$

We want to prove that we have $a=0$. Since $v$ is biunitary, we have:
\begin{eqnarray*}
a
&=&\sum_i\left(\sum_j\left(v_{ij}\xi_j-\frac{1}{N}\xi_i\right)\right)\left(\sum_k\left(v_{ik}^*\bar{\xi}_k-\frac{1}{N}\bar{\xi}_i\right)\right)\\
&=&\sum_{ijk}v_{ij}v_{ik}^*\xi_j\bar{\xi}_k-\frac{1}{N}v_{ij}\xi_j\bar{\xi}_i-\frac{1}{N}v_{ik}^*\xi_i\bar{\xi}_k+\frac{1}{N^2}\xi_i\bar{\xi}_i\\
&=&\sum_j|\xi_j|^2-\sum_{ij}v_{ij}\xi_j\bar{\xi}_i-\sum_{ik}v_{ik}^*\xi_i\bar{\xi}_k+\sum_i|\xi_i|^2\\
&=&||\xi||^2-<v\xi,\xi>-\overline{<v\xi,\xi>}+||\xi||^2\\
&=&2(||\xi||^2-Re(<v\xi,\xi>))
\end{eqnarray*}

By using now our assumption $M\xi=\xi$, we obtain from this:
\begin{eqnarray*}
\varphi(a)
&=&2\varphi(||\xi||^2-Re(<v\xi,\xi>))\\
&=&2(||\xi||^2-Re(<M\xi,\xi>))\\
&=&2(||\xi||^2-||\xi||^2)\\
&=&0
\end{eqnarray*}

Thus $a=0$, and by positivity we obtain $v\xi=\xi$, as desired.
\end{proof}

We can now formulate the general Haar measure result, due to Woronowicz \cite{wo1}:

\index{Haar measure}
\index{Haar integration}
\index{Ces\`aro limit}
\index{uniform measure}
\index{uniform integration}

\begin{theorem}
Any Woronowicz algebra has a unique Haar integration, which can be constructed by starting with any faithful positive unital state $\varphi\in A^*$, and setting
$$\int_G=\lim_{n\to\infty}\frac{1}{n}\sum_{k=1}^n\varphi^{*k}$$
where $\phi*\psi=(\phi\otimes\psi)\Delta$. Moreover, for any corepresentation $v$ we have
$$\left(id\otimes\int_G\right)v=P$$
where $P$ is the orthogonal projection onto $Fix(v)=\{\xi\in\mathbb C^n|v\xi=\xi\}$.
\end{theorem}

\begin{proof}
Let us first go back to the general context of Proposition 10.13. Since convolving one more time with $\varphi$ will not change the Ces\`aro limit appearing there, the functional $\int_\varphi\in A^*$ constructed there has the following invariance property:
$$\int_\varphi*\,\varphi=\varphi*\int_\varphi=\int_\varphi$$

In the case where $\varphi$ is assumed to be faithful, as in Proposition 10.14, our claim is that we have the following formula, valid this time for any $\psi\in A^*$:
$$\int_\varphi*\,\psi=\psi*\int_\varphi=\psi(1)\int_\varphi$$

By linearity, it is enough to prove this formula on a coefficient of a corepresentation, $a=(\tau\otimes id)v$. In order to do so, consider the following matrices:
$$P=\left(id\otimes\int_\varphi\right)v\quad,\quad 
Q=(id\otimes\psi)v$$

In terms of these matrices, we have the following formula:
$$\left(\int_\varphi*\,\psi\right)a
=\left(\tau\otimes\int_\varphi\otimes\,\psi\right)(v_{12}v_{13})
=\tau(PQ)$$

Similarly, we have the following computation, for the term in the middle:
$$\left(\psi*\int_\varphi\right)a
=\left(\tau\otimes\psi\otimes\int_\varphi\right)(v_{12}v_{13})
=\tau(QP)$$

Finally, regarding the term on the right, this is given by:
$$\psi(1)\int_\varphi a=\psi(1)\tau(P)$$

Thus, our claim above is equivalent to the following equality:
$$PQ=QP=\psi(1)P$$

But this latter equality follows from the fact, coming from Proposition 10.14, that $P=(id\otimes\int_\varphi)v$ equals the orthogonal projection onto $Fix(v)$. Thus, we have proved our claim. Now observe that our formula can be written as:
$$\psi\left(\int_\varphi\otimes id\right)\Delta=\psi\left(id\otimes\int_\varphi\right)\Delta=\psi\int_\varphi(.)1$$

This formula being true for any $\psi\in A^*$, we can simply delete $\psi$, and we conclude that the invariance formula holds indeed, with $\int_G=\int_\varphi$. Finally, assuming that we have two invariant integrals $\int_G,\int_G'$, we have the following computation:
\begin{eqnarray*}
\left(\int_G\otimes\int_G'\right)\Delta
&=&\left(\int_G'\otimes\int_G\right)\Delta\\
&=&\int_G(.)1\\
&=&\int_G'(.)1
\end{eqnarray*}

Thus we have $\int_G=\int_G'$, and this finishes the proof.
\end{proof}

Let us record as well the following result, which is very useful in practice:

\index{Weingarten formula}
\index{Tannakian category}
\index{Gram matrix}

\begin{theorem}
The integration over $G\subset_uU_N^+$ is given by the Weingarten formula
$$\int_Gu_{i_1j_1}^{e_1}\ldots u_{i_kj_k}^{e_k}=\sum_{\pi,\sigma\in D_k}\delta_\pi(i)\delta_\sigma(j)W_k(\pi,\sigma)$$
for any colored integer $k=e_1\ldots e_k$ and indices $i,j$, where $D_k$ is a linear basis of $Fix(u^{\otimes k})$, 
$$\delta_\pi(i)=<\pi,e_{i_1}\otimes\ldots\otimes e_{i_k}>$$
and $W_k=G_k^{-1}$, with $G_k(\pi,\sigma)=<\pi,\sigma>$.
\end{theorem}

\begin{proof}
We know from Theorem 10.15 that the integrals in the statement form altogether the orthogonal projection $P^k$ onto the following space:
$$Fix(u^{\otimes k})=span(D_k)$$

Consider now the following linear map, with $D_k=\{\xi_k\}$ being as in the statement:
$$E(x)=\sum_{\pi\in D_k}<x,\xi_\pi>\xi_\pi$$

By a standard linear algebra computation, it follows that we have $P=WE$, where $W$ is the inverse on $span(T_\pi|\pi\in D_k)$ of the restriction of $E$. But this restriction is the linear map given by $G_k$, and so $W$ is the linear map given by $W_k$, and this gives the result.
\end{proof}

The above result is something quite far-reaching, allowing many explicit computations, and we will be back to it in chapter 11, with examples and applications.

\section*{10c. Peter-Weyl theory}

Let us get back now to algebra, and establish two more Peter-Weyl theorems. We will need the following result, which is something very useful, of independent interest:

\index{Frobenius isomorphism}

\begin{proposition}
We have a Frobenius type isomorphism
$$Hom(v,w)\simeq Fix(\bar{v}\otimes w)$$
valid for any two corepresentations $v,w$.
\end{proposition}

\begin{proof}
According to the definitions, we have the following equivalence:
\begin{eqnarray*}
T\in Hom(v,w)
&\iff&Tv=wT\\
&\iff&\sum_jT_{aj}v_{ji}=\sum_bw_{ab}T_{bi}
\end{eqnarray*}

On the other hand, we have as well the following equivalence:
\begin{eqnarray*}
T\in Fix(\bar{v}\otimes w)
&\iff&(\bar{v}\otimes w)T=T\\
&\iff&\sum_{kb}v_{ik}^*w_{ab}T_{bk}=T_{ai}
\end{eqnarray*}

With these formulae in hand, we must prove that we have:
$$\sum_jT_{aj}v_{ji}=\sum_bw_{ab}T_{bi}\iff \sum_{kb}v_{ik}^*w_{ab}T_{bk}=T_{ai}$$

(1) In one sense, the computation is as follows, using the unitarity of $v^t$:
\begin{eqnarray*}
\sum_{kb}v_{ik}^*w_{ab}T_{bk}
&=&\sum_kv_{ik}^*\sum_bw_{ab}T_{bk}\\
&=&\sum_kv_{ik}^*\sum_jT_{aj}v_{jk}\\
&=&\sum_j(\bar{v}v^t)_{ij}T_{aj}\\
&=&T_{ai}
\end{eqnarray*}

(2) In the other sense we have, once again by using the unitarity of $v^t$:
\begin{eqnarray*}
\sum_jT_{aj}v_{ji}
&=&\sum_jv_{ji}\sum_{kb}v_{jk}^*w_{ab}T_{bk}\\
&=&\sum_{kb}(v^t\bar{v})_{ik}w_{ab}T_{bk}\\
&=&\sum_bw_{ab}T_{bi}
\end{eqnarray*}

Thus, we are led to the conclusion in the statement.
\end{proof}

We can now establish a third Peter-Weyl theorem, as follows:

\index{Peter-Weyl theorem}
\index{direct sum decomposition}
\index{Peter-Weyl algebra}
\index{scalar product of functions}

\begin{theorem}[PW3]
The dense subalgebra $\mathcal A\subset A$ decomposes as a direct sum 
$$\mathcal A=\bigoplus_{v\in Irr(A)}M_{\dim(v)}(\mathbb C)$$
with this being an isomorphism of $*$-coalgebras, and with the summands being pairwise orthogonal with respect to the scalar product given by
$$<a,b>=\int_Gab^*$$
where $\int_G$ is the Haar integration over $G$.
\end{theorem}

\begin{proof}
By combining the previous Peter-Weyl results, from Theorem 10.10 and Theorem 10.12, we deduce that we have a linear space decomposition as follows:
$$\mathcal A
=\sum_{v\in Irr(A)}C(v)
=\sum_{v\in Irr(A)}M_{\dim(v)}(\mathbb C)$$

Thus, in order to conclude, it is enough to prove that for any two irreducible corepresentations $v,w\in Irr(A)$, the corresponding spaces of coefficients are orthogonal:
$$v\not\sim w\implies C(v)\perp C(w)$$ 

But this follows from Theorem 10.15, via Proposition 10.17. Let us set indeed:
$$P_{ia,jb}=\int_Gv_{ij}w_{ab}^*$$

Then $P$ is the orthogonal projection onto the following vector space:
$$Fix(v\otimes\bar{w})
\simeq Hom(\bar{v},\bar{w})
=\{0\}$$

Thus we have $P=0$, and this gives the result.
\end{proof}

Finally, we have the following result, completing the Peter-Weyl theory:

\index{Peter-Weyl theorem}
\index{central function}
\index{character}
\index{character of corepresentation}
\index{Peter-Weyl basis}
\index{orthonormal basis}

\begin{theorem}[PW4]
The characters of irreducible corepresentations belong to
$$\mathcal A_{central}=\left\{a\in\mathcal A\Big|\Sigma\Delta(a)=\Delta(a)\right\}$$
called ``algebra of smooth central functions'', and form an orthonormal basis of it.
\end{theorem}

\begin{proof}
As a first remark, the linear space $\mathcal A_{central}$ defined above is indeed an algebra. In the classical case, we obtain in this way the usual algebra of smooth central functions. Also, in the group dual case, where we have $\Sigma\Delta=\Delta$, we obtain the whole convolution algebra. Regarding now the proof, in general, this goes as follows:

\medskip

(1) The algebra $\mathcal A_{central}$ contains indeed all the characters, because we have:
\begin{eqnarray*}
\Sigma\Delta(\chi_v)
&=&\Sigma\left(\sum_{ij}v_{ij}\otimes v_{ji}\right)\\
&=&\sum_{ij}v_{ji}\otimes v_{ij}\\
&=&\Delta(\chi_v)
\end{eqnarray*}

(2) Conversely, consider an element $a\in\mathcal A$, written as follows:
$$a=\sum_{v\in Irr(A)}a_v$$

The condition $a\in\mathcal A_{central}$ is then equivalent to the following conditions:
$$a_v\in\mathcal A_{central}\quad,\forall v\in Irr(A)$$

But each condition $a_v\in\mathcal A_{central}$ means that $a_v$ must be a scalar multiple of the corresponding character $\chi_v$, and so the characters form a basis of $\mathcal A_{central}$, as stated.

\medskip

(3) The fact that we have an orthogonal basis follows from Theorem 10.18. 

\medskip

(4) Finally, regarding the norm 1 assertion, consider the following integrals:
$$P_{ik,jl}=\int_Gv_{ij}v_{kl}^*$$

We know from Theorem 10.15 that these integrals form the orthogonal projection onto the following vector space, computed via Proposition 10.17:
$$Fix(v\otimes\bar{v})
\simeq End(\bar{v})
=\mathbb C1$$

By using this fact, we obtain the following formula:
\begin{eqnarray*}
\int_G\chi_v\chi_v^*
&=&\sum_{ij}\int_Gv_{ii}v_{jj}^*\\
&=&\sum_i\frac{1}{N}\\
&=&1
\end{eqnarray*}

Thus the characters have indeed norm 1, and we are done.
\end{proof}

We can now solve a problem that we left open before, namely:

\index{cocommutative algebra}

\begin{proposition}
The cocommutative Woronowicz algebras appear as the quotients
$$C^*(\Gamma)\to A\to C^*_{red}(\Gamma)$$
given by $A=C^*_\pi(\Gamma)$ with $\pi\otimes\pi\subset\pi$, with $\Gamma$ being a discrete group.
\end{proposition}

\begin{proof}
This follows from the Peter-Weyl theory, and clarifies a number of things said before, in Theorem 10.5. Indeed, for a cocommutative Woronowicz algebra the irreducible corepresentations are all 1-dimensional, and this gives the result.
\end{proof}

As another consequence of the above results, once again by following Woronowicz \cite{wo1}, we have the following statement, dealing with functional analysis aspects, and extending some standard facts regarding the $C^*$-algebras of the usual discrete groups:

\index{full algebra}
\index{reduced algebra}
\index{amenable quantum group}
\index{coamenable quantum group}
\index{Kesten amenability}

\begin{theorem}
Let $A_{full}$ be the enveloping $C^*$-algebra of $\mathcal A$, and $A_{red}$ be the quotient of $A$ by the null ideal of the Haar integration. The following are then equivalent:
\begin{enumerate}
\item The Haar functional of $A_{full}$ is faithful.

\item The projection map $A_{full}\to A_{red}$ is an isomorphism.

\item The counit map $\varepsilon:A_{full}\to\mathbb C$ factorizes through $A_{red}$.

\item We have $N\in\sigma(Re(\chi_v))$, the spectrum being taken inside $A_{red}$.
\end{enumerate}
If this is the case, we say that the underlying discrete quantum group $\Gamma$ is amenable.
\end{theorem}

\begin{proof}
This is well-known in the group dual case, $A=C^*(\Gamma)$, with $\Gamma$ being a usual discrete group. In general, the result follows by adapting the group dual case proof:

\medskip

$(1)\iff(2)$ This simply follows from the fact that the GNS construction for the algebra $A_{full}$ with respect to the Haar functional produces the algebra $A_{red}$.

\medskip

$(2)\iff(3)$ Here $\implies$ is trivial, and conversely, a counit $\varepsilon:A_{red}\to\mathbb C$ produces an isomorphism $\Phi:A_{red}\to A_{full}$, by slicing the map $\widetilde{\Delta}:A_{red}\to A_{red}\otimes A_{full}$.

\medskip

$(3)\iff(4)$ Here $\implies$ is clear, coming from $\varepsilon(N-Re(\chi (u)))=0$, and the converse can be proved by doing some functional analysis. See \cite{wo1}.
\end{proof}

The above was quite short, but we will be back to this, later, with further details. Also, many other things can be said about amenability. We will be back to amenability on several occasions, in what follows, and notably in chapter 12 below.

\section*{10d. Tannakian duality}

In this section we discuss Tannakian duality, in its various forms, and there are many of them, together with its various abstract applications. This will further clarify some of the theory that we learned both in Parts I-II, and in the present chapter. 

\bigskip

Let us first discuss Tannakian duality for the usual compact Lie groups. We use the following notion, that already we met in the above, in more general settings:

\index{Peter-Weyl representation}
\index{Peter-Weyl corepresentation}

\begin{definition}
Given a closed subgroup $G\subset U_N$, its Peter-Weyl representations are the tensor products between the fundamental representation and its conjugate,
$$u:G\subset U_N\quad,\quad \bar{u}:G\subset U_N$$
denoted $u^{\otimes k}$, and indexed by colored integers $k=\circ\bullet\bullet\circ\ldots$ according to the formulae
$$u^{\otimes\emptyset}=1\quad,\quad 
u^{\otimes\circ}=u\quad,\quad 
u^{\otimes\bullet}=\bar{u}$$
and multiplicativity, $u^{\otimes kl}=u^{\otimes k}\otimes u^{\otimes l}$.
\end{definition}

Here are a few examples of such corepresentations, namely those coming from the colored integers of length 2, to be often used in what follows:
$$u^{\otimes\circ\circ}=u\otimes u\quad,\quad 
u^{\otimes\circ\bullet}=u\otimes\bar{u}$$
$$u^{\otimes\bullet\circ}=\bar{u}\otimes u\quad,\quad 
u^{\otimes\bullet\bullet}=\bar{u}\otimes\bar{u}$$

As something new appearing in the classical setting, with respect to what we already knew about the representations of compact quantum groups, observe that we have the following formula, due to the commutativity of the algebra $C(G)$:
$$u^{\otimes\circ\bullet}\sim u^{\otimes\bullet\circ}$$

In view of this, we can index if we want the Peter-Weyl representations by indices $l\in\mathbb N\times\mathbb N$, by moving all $\circ$ symbols to the left, and all $\bullet$ symbols to the right. However, since we want to generalize this afterwards to the case of quantum groups, which are not necessarily braided, it is better to use colored exponents $k=\circ\bullet\bullet\circ\ldots\,$, as above.

\bigskip

We can now formulate a key result, namely Tannakian duality for the compact groups, in its various abstract and concrete formulations, as follows:

\index{Tannakian duality}
\index{tensor category}
\index{symmetries of functor}
\index{Doplicher-Roberts theorem}
\index{Deligne theorem}

\begin{theorem}
Given a compact group $G$, the following happen:
\begin{enumerate}
\item $G$ can be reconstructed from the knowledge of the associated tensor category of representations $\mathcal C$, as group of symmetries of the functor $\mathcal C\to Vect(\mathbb C)$. 

\item When $G\subset U_N$, the group $G$ can be reconstructed from the knowledge of the Peter-Weyl subcategory $C\subset\mathcal C$, in a similar way. 

\item Moreover, there is in fact no need for the functor $\mathcal C\to Vect(\mathbb C)$, with the group $G$ appearing as the plain group of symmetries of $\mathcal C$, or of $C$.
\end{enumerate}
\end{theorem}

\begin{proof}
This is something quite technical, the idea being as follows:

\medskip

(1) This is something very classical, due to Tannaka and Krein, with the statement being self-explanatory category theory, and with the proof being long, but straightforward too. In what follows we will not really need this result, so we refer here to the classical Tannakian literature, or to the paper of Woronowicz \cite{wo2}, for generalizations.

\medskip

(2) This is something sharper, that will be our central formulation of Tannakian duality, for our purposes here, which follows from (1), or from the following formula:
$$G=\left\{g\in U_N\Big| Tg^{\otimes k}=g^{\otimes l}T,\forall k,l,\forall T\in Hom(u^{\otimes k},u^{\otimes l})\right\}$$

To be more precise, our first claim is that in this formula the set on the right is an intermediate subgroup $G\subset\widetilde{G}\subset U_N$. Indeed, consider the set on the right:
$$\widetilde{G}=\left\{g\in U_N\Big| Tg^{\otimes k}=g^{\otimes l}T,\forall k,l,\forall T\in Hom(u^{\otimes k},u^{\otimes l})\right\}$$

Assuming now that we have $g,h\in\widetilde{G}$, it follows that we have $gh\in\widetilde{G}$, due to the following computation, valid for any $k,l$ and any $T\in C_{kl}$:
\begin{eqnarray*}
T(gh)^{\otimes k}
&=&Tg^{\otimes k}h^{\otimes k}\\
&=&g^{\otimes l}Th^{\otimes k}\\
&=&g^{\otimes l}h^{\otimes l}T\\
&=&(gh)^{\otimes l}T
\end{eqnarray*}

Also, we have $1\in\widetilde{G}$, trivially. Finally, assuming $g\in\widetilde{G}$, we have:
\begin{eqnarray*}
T(g^{-1})^{\otimes k}
&=&(g^{-1})^{\otimes l}[g^{\otimes l}T](g^{-1})^{\otimes k}\\
&=&(g^{-1})^{\otimes l}[Tg^{\otimes k}](g^{-1})^{\otimes k}\\
&=&(g^{-1})^{\otimes l}T
\end{eqnarray*}

Thus we have $g^{-1}\in\widetilde{G}$, and so $\widetilde{G}$ is a group, as claimed. Finally, the fact that we have an inclusion $G\subset\widetilde{G}$, and that $\widetilde{G}\subset U_N$ is closed, are both clear. Thus, we have an intermediate subgroup as follows, that we want to prove to be equal to $G$ itself:
$$G\subset\widetilde{G}\subset U_N$$

In order to prove this, consider the Tannakian category of $\widetilde{G}$, namely:
$$\widetilde{C}_{kl}=\left\{T\in\mathcal L(H^{\otimes k},H^{\otimes l})\Big|Tg^{\otimes k}=g^{\otimes l}T,\forall g\in\widetilde{G}\right\}$$

By functoriality, from $G\subset\widetilde{G}$ we obtain $\widetilde{C}\subset C$. On the other hand, according to the definition of $\widetilde{G}$, we have $C\subset\widetilde{C}$. Thus, we have the following equality:
$$C=\widetilde{C}$$

Assume now by contradiction that $G\subset\widetilde{G}$ is not an equality. Then, at the level of algebras of functions, the following quotient map is not an isomorphism either:
$$C(\widetilde{G})\to C(G)$$ 

On the other hand, we know from Peter-Weyl that we have decompositions as follows, with the sums being over all the irreducible unitary representations:
$$C(\widetilde{G})=\overline{\bigoplus}_{\pi\in Irr(\widetilde{G})}M_{\dim\pi}(\mathbb C)
\quad,\quad 
C(G)=\overline{\bigoplus}_{\nu\in Irr(G)}M_{\dim\nu}(\mathbb C)$$

Now observe that each unitary representation $\pi:\widetilde{G}\to U_K$ restricts into a certain representation $\pi':G\to U_K$. Since the quotient map $C(\widetilde{G})\to C(G)$ is not an isomorphism, we conclude that there is at least one representation $\pi$ satisfying:
$$\pi\in Irr(\widetilde{G})\quad,\quad \pi'\notin Irr(G)$$

We are now in position to conclude. By using Peter-Weyl theory again, the above representation $\pi\in Irr(\widetilde{G})$ appears in a certain tensor power of the fundamental representation $u:\widetilde{G}\subset U_N$. Thus, we have inclusions of representations, as follows:
$$\pi\in u^{\otimes k}\quad,\quad\pi'\in u'^{\otimes k}$$

Now since we know that $\pi$ is irreducible, and that $\pi'$ is not, by using one more time Peter-Weyl theory, we conclude that we have a strict inequality, as follows:
\begin{eqnarray*}
\dim(\widetilde{C}_{kk})
&=&dim(End(u^{\otimes k}))\\
&<&dim(End(u'^{\otimes k}))\\
&=&\dim(C_{kk})
\end{eqnarray*}

But this contradicts the equality $C=\widetilde{C}$ found above, which finishes the proof.

\medskip

(3) This is something far more tricky, due to Doplicher-Roberts \cite{dro} and Deligne \cite{del}, and whose proof is non-trivial. In what follows, we will not really need this result, so we will simply refer here to the above-mentioned papers, and the related literature.
\end{proof}

There are many more things that can be said here, and in particular, we have:

\index{planar algebra}

\begin{fact}
The diagonal part of $C=(C_{kl})$, formed by the algebras
$$C_{kk}=\left\{T\in\mathcal L(H^{\otimes k})\Big|Tg^{\otimes k}=g^{\otimes k}T,\forall g\in G\right\}$$
does not determine $G$. For instance $G=\{1\},\mathbb Z_2$ are not distinguished by it.
\end{fact}

Obviously, this is something quite annoying, because there are countless temptations to use $\Delta C=(C_{kk})$ instead of $C$, for instance because the spaces $C_{kk}$ are algebras, and also, at a more advanced level, because $\Delta C$ is a planar algebra in the sense of Jones \cite{jo3}. But, we are not allowed to do this, at least in general. More on this later.

\bigskip

Our purpose now will be to establish Tannakian duality results which are analogous to those above, for the compact quantum groups. We first have the following result:

\begin{theorem}
Any compact quantum group $G$ can be reconstructed from the knowledge of the associated tensor category of representations $\mathcal C$, as being, in a suitable sense, the quantum group of symmetries of the functor $\mathcal C\to Vect(\mathbb C)$.
\end{theorem}

\begin{proof}
This is something a bit technical, although following the same ideas as for the usual compact groups, and we refer here to Woronowicz \cite{wo2} and related papers. In what follows, technically speaking, we will not really need this result.
\end{proof}

Importantly now, and contrary to what happensfor the usual groups, the need for the functor $\mathcal C\to Vect(\mathbb C)$ in the above is essential, at least at our level of generality in this book. For more on all this, need for this functor, we refer to \cite{egn} and related papers.

\bigskip

Summarizing, we are left with extending Theorem 10.23 (2). In order to discuss this, let us start with the following definition, inspired from the Peter-Weyl categories:

\begin{definition}
Let $H$ be a finite dimensional Hilbert space. A tensor category over $H$ is a collection $C=(C_{kl})$ of subspaces 
$$C_{kl}\subset\mathcal L(H^{\otimes k},H^{\otimes l})$$
satisfying the following conditions:
\begin{enumerate}
\item $S,T\in C$ implies $S\otimes T\in C$.

\item If $S,T\in C$ are composable, then $ST\in C$.

\item $T\in C$ implies $T^*\in C$.

\item Each $C_{kk}$ contains the identity operator.

\item $C_{\emptyset,\circ\bullet}$ and $C_{\emptyset,\bullet\circ}$ contain the operator $R:1\to\sum_ie_i\otimes e_i$.
\end{enumerate}
\end{definition}

To be more precise, the point with this is that, given a Woronowicz algebra $(A,u)$ with $u\in\mathcal L(H)\otimes A$, the conditions (1-5) are those satisfied by the following spaces:
$$C_{kl}=Hom(u^{\otimes k},u^{\otimes l})$$

We have the following key result, coming from the work of Woronowicz in \cite{wo2}:

\index{Peter-Weyl category}
\index{Peter-Weyl subcategory}
\index{restricted duality}
\index{restricted Tannakian duality}

\begin{theorem}
The following operations are inverse to each other:
\begin{enumerate}
\item The construction $G\to C$, which associates to a compact quantum group $G$ the tensor category formed by the intertwiner spaces $C_{kl}=Hom(u^{\otimes k},u^{\otimes l})$.

\item The construction $C\to G$, associating to a tensor category $C$ the compact quantum group $G$ coming from the relations $T\in Hom(u^{\otimes k},u^{\otimes l})$, with $T\in C_{kl}$.
\end{enumerate}
\end{theorem}

\begin{proof}
This is something quite deep, going back to Woronowicz \cite{wo2} in a slightly different form, and to Malacarne \cite{mal} in the simplified form above:

\medskip

(1) Regarding the precise statement of the result, the idea is that, fixing $N\in\mathbb N$ as the matrix size of our quantum groups, we have indeed a construction $G\to C_G$, whose output is a tensor category over the Hilbert space $H=\mathbb C^N$, as follows:
$$(C_G)_{kl}=Hom(u^{\otimes k},u^{\otimes l})$$

We have as well a construction $C\to G_C$, producing out of a Tannakian category $C$ a certain compact quantum group $G_C$, according to the following formula:
$$C(G_C)=C^*\Big((u_{ij})_{i,j=1,\ldots,N}\Big)\Big/\Big<Tu^{\otimes k}=u^{\otimes l}T,\forall k,l,\forall T\in C_{kl}\Big>$$

Thus, we have correspondences as in the statement, and the result itself states that these correspondences are inverse to each other, in the sense that we have:
$$C_{G_C}=C\quad,\quad G_{C_G}=G$$

(2) Regarding now the proof, this can be deduced from Theorem 10.25, but following Malacarne \cite{mal}, we can prove this as well directly. Indeed, up to some standard study, we must establish an inclusion as follows, for any Tannakian category $C$:
$$C_{G_C}\subset C$$

Now for this purpose, given a tensor category $C=(C_{kl})$ over a Hilbert space $H=\mathbb C^N$, consider the following direct sum, which is easily seen to be a $*$-algebra:
$$E_C
=\bigoplus_{k,l}C_{kl}
\subset\bigoplus_{k,l}B(H^{\otimes k},H^{\otimes l})
\subset B\left(\bigoplus_kH^{\otimes k}\right)$$

Consider also, inside this $*$-algebra, the following $*$-subalgebra:
$$E_C^{(s)}
=\bigoplus_{|k|,|l|\leq s}C_{kl}
\subset\bigoplus_{|k|,|l|\leq s}B(H^{\otimes k},H^{\otimes l})
=B\left(\bigoplus_{|k|\leq s}H^{\otimes k}\right)$$

It is then routine to check that we have equivalences as follows:
\begin{eqnarray*}
C_{G_C}\subset C
&\iff&E_{C_{G_C}}\subset E_C\\
&\iff&E_{C_{G_C}}^{(s)}\subset E_C^{(s)},\forall s\\
&\iff&E_{C_{G_C}}^{(s)'}\supset E_C^{(s)'},\forall s
\end{eqnarray*}

Summarizing, we would like to prove that we have the following inclusion:
$$E_C^{(s)'}\subset E_{C_{G_C}}^{(s)'}$$

(3) In order to prove this, let us first study the commutant on the right. We have the following equality, between subalgebras of the algebra $B\left(\bigoplus_{|k|\leq s}H^{\otimes k}\right)$:
$$E_{C_G}^{(s)}=End\left(\bigoplus_{|k|\leq s}u^{\otimes k}\right)$$

We have to compute the commutant of this algebra. But for this purpose, we can use the standard fact that, given a unitary corepresentation $v\in M_n(C(G))$, we have an algebra representation as follows, whose image is given by $Im(\pi_v)=End(v)'$:
$$\pi_v:C(G)^*\to M_n(\mathbb C)\quad,\quad 
\varphi\to (\varphi(v_{ij}))_{ij}$$

(4) Generally speaking, in order to prove anything about $G_C$, we are in need of an explicit model for this quantum group. In order to construct such a model, let $<u_{ij}>$ be the free $*$-algebra over $\dim(H)^2$ variables, with bialgebra structure as follows:
$$\Delta(u_{ij})=\sum_ku_{ik}\otimes u_{kj}\quad,\quad 
\varepsilon(u_{ij})=\delta_{ij}$$

Consider also the following pair of dual vector spaces:
$$V=\bigoplus_kB\left(H^{\otimes k}\right)\quad,\quad 
V^*=\bigoplus_kB\left(H^{\otimes k}\right)^*$$

Finally, let $f_{ij},f_{ij}^*\in V^*$ be the standard generators of $B(H)^*,B(\bar{H})^*$. Then:

\medskip

-- $V^*$ is a $*$-algebra, with multiplication $\otimes$ and involution as follows:
$$f_{ij}\leftrightarrow f_{ij}^*$$

-- $V^*$ is in fact a $*$-bialgebra, with $*$-bialgebra operations as follows:
$$\Delta(f_{ij})=\sum_kf_{ik}\otimes f_{kj}\quad,\quad 
\varepsilon(f_{ij})=\delta_{ij}$$

-- We have a $*$-bialgebra isomorphism $<u_{ij}>\simeq V^*$, given by $u_{ij}\to f_{ij}$.

\medskip

(5) The point now is that the smooth part of the algebra $A_C=C(G_C)$ is given by $\mathcal A_C\simeq V^*/J$, where $J\subset V^*$ is the ideal coming from the following relations, for any $i,j$, one for each pair of colored integers $k,l$, and each $T\in C(k,l)$:
$$\sum_{p_1,\ldots,p_k}T_{i_1\ldots i_l,p_1\ldots p_k}f_{p_1j_1}\otimes\ldots\otimes  f_{p_kj_k}
=\sum_{q_1,\ldots,q_l}T_{q_1\ldots q_l,j_1\ldots j_k}f_{i_1q_1}\otimes\ldots\otimes f_{i_lq_l}$$

Moreover, the linear space $\mathcal A_C^*$ is given by the following formula:
$$\mathcal A_C^*=\left\{a\in V\Big|Ta_k=a_lT,\forall T\in C(k,l)\right\}$$

(6) Now consider the representation that we are interested in, namely:
$$\pi_v:\mathcal A_C^*\to B\left(\bigoplus_{|k|\leq s}H^{\otimes k}\right)$$

This representation appears then diagonally, by truncating, as follows:
$$\pi_v:a\to (a_k)_{kk}$$

(7) In order to further advance, consider the following vector spaces:
$$V_s=\bigoplus_{|k|\leq s}B\left(H^{\otimes k}\right)\quad,\quad 
V^*_s=\bigoplus_{|k|\leq s}B\left(H^{\otimes k}\right)^*$$

If we denote by $a\to a_s$ the truncation operation $V\to V_s$, it is routine to check, by using the categorical axioms for $C$, that the following inclusions and equalities hold:

\medskip

-- $E_C^{(s)'}\subset V_s$.

\medskip

-- $E_C'\subset V$.

\medskip

-- $\mathcal A_C^*=E_C'$.

\medskip

-- $Im(\pi_v)=(E_C')_s$.

\medskip

(8) We can now prove the duality. According to our results above, we have to prove that, for any Tannakian category $C$, and any $s\in\mathbb N$, we have an inclusion as follows:
$$E_C^{(s)'}\subset(E_C')_s$$

By taking duals, this is the same as proving that we have:
$$\left\{f\in V_s^*\Big|f_{|(E_C')_s}=0\right\}\subset\left\{f\in V_s^*\Big|f_{|E_C^{(s)'}}=0\right\}$$

(9) In order to do so, we use the following formula, established above: 
$$\mathcal A_C^*=E_C'$$

We also know from the above that we have an identification as follows:
$$\mathcal A_C=V^*/J$$

We conclude that the ideal $J$ is given by the following formula:
$$J=\left\{f\in V^*\Big|f_{|E_C'}=0\right\}$$

(10) Our claim now is that we have the following formula, for any $s\in\mathbb N$:
$$J\cap V_s^*=\left\{f\in V_s^*\Big|f_{|E_C^{(s)'}}=0\right\}$$

Indeed, let us denote by $X_s$ the spaces on the right. The axioms for $C$ show that these spaces are increasing, that their union $X=\cup_sX_s$ is an ideal, and that:
$$X_s=X\cap V_s^*$$

(11) We must prove that we have $J=X$, and this can be done as follows:

\medskip

``$\subset$'' This follows from the following fact, for any $T\in C(k,l)$ with $|k|,|l|\leq s$:
\begin{eqnarray*}
(f_T)_{|\{T\}'}=0
&\implies&(f_T)_{|E_C^{(s)'}}=0\\
&\implies&f_T\in X_s
\end{eqnarray*}

``$\supset$'' This follows from our description of $J$, because from $E_C^{(s)}\subset E_C$ we obtain:
$$f_{|E_C^{(s)'}}=0\implies f_{|E_C'}=0$$

(12) Summarizing, we have proved our claim. On the other hand, we have:
\begin{eqnarray*}
J\cap V_s^*
&=&\left\{f\in V^*\Big|f_{|E_C'}=0\right\}\cap V_s^*\\
&=&\left\{f\in V_s^*\Big|f_{|E_C'}=0\right\}\\
&=&\left\{f\in V_s^*\Big|f_{|(E_C')_s}=0\right\}
\end{eqnarray*}

Thus, our claim is exactly the inclusion that we wanted to prove, and we are done.
\end{proof}

\section*{10e. Exercises}

We had a lot of advanced algebra in this chapter, and as exercises, we have:

\begin{exercise}
Learn about further examples of semisimple quantum groups.
\end{exercise}

\begin{exercise}
Discuss as well the Haar functionals of these quantum groups.
\end{exercise}

\begin{exercise}
Can we have some Lie theory going on, in our setting?
\end{exercise}

\begin{exercise}
Learn with full details the classical Tannakian duality.
\end{exercise}

\begin{exercise}
Learn as well the Doplicher-Roberts and Deligne theorems.
\end{exercise}

\begin{exercise}
Clarify all the details in relation with Tannakian duality.
\end{exercise}

As bonus exercise, have a look at the Drinfeld-Jimbo theory of deformations.

\chapter{Brauer theorems}

\section*{11a. Diagrams, easiness}

Good news, with the theory that we developed so far in this book, we can now talk about easiness. As an answer to a question that you might have right away, namely what exactly is easiness, well, this is something quite philosophical, as follows:

\begin{answer}
The real, true quantum groups, which do not depend on the ground field $F$, and which in addition are something simple and easy, are the easy ones.
\end{answer}

In practice now, yes I know, this rather sounds like an enigmatic statement coming from my cat, so more details are certainly needed. Here they are, in brief:

\bigskip

(1) We have learned in the previous chapter, under the various technical assumptions there, the Tannakian duality correspondence, $G\leftrightarrow C$. In this correspondence both the quantum group $G$ and the tensor category $C$ were of course taken over $F=\mathbb C$.

\bigskip

(2) However, with a bit more work, and some extra assumptions, we can get rid of $F=\mathbb C$. Indeed, we can call a tensor category $C$ easy when it is of the form $C=span(D)$, with $D$ being a certain purely combinatorial structure, called category of partitions.

\bigskip

(3) But then, by using the Tannakian correspondence $G\leftrightarrow C$, we can call a quantum group $G$ easy when its corresponding Tannakian category $C$ is easy. That is, we can call $G$ easy when we have $C=span(D)$, for a certain category of partitions $D$.

\bigskip

(4) Thus, we have our basically scalarless, easy quantum groups $G$. And with this being something powerful, and useful, because the easiness condition $C=span(D)$ corresponds every time to a subtle and useful theorem, called Brauer theorem for $G$.

\bigskip

So, this was for the general idea, and we will see the details in a moment. Before starting, however, you might still wonder why using $F=\mathbb C$ in our plan above, instead of simply developing easiness over an arbitrary field $F$. Good point, and in answer:

\begin{comment}
Normally easiness can be developed over an arbitrary field $F$, provided that we have a good knowledge of Tannakian duality there. In our situation, however, which is rather beginner level, we will use $F=\mathbb C$, as indicated above.
\end{comment}

And for more details here, developing easiness over $F=\mathbb C$, which is what we will do, is something quite standard, following \cite{bsp} and related papers, which themselves are basically a slight reformulation and generalization of some old work of Weyl, Brauer and others. As for easiness over an arbitrary field $F$, this is something more subtle, requiring far more algebra know-how, and we refer here to the papers of Bichon and his collaborators \cite{bi1}, \cite{bi2}, \cite{bi3}, \cite{bi4}, \cite{bi5}, \cite{bdv}, \cite{bd1}, \cite{bd2}, \cite{bna}, \cite{bny}, which basically do this.

\bigskip

In short, modesty, what we will be doing here will be rather an introduction to easiness, instead of the presentation of the final, and probably trivial easiness theory, which undoubtedly exists somewhere, on the heights of abstract algebra, say next to schemes.

\bigskip

And with this discussion made, time I guess to get to work. But not however before asking the cat about what he thinks, about our plan. And cat declares:

\begin{cat}
Human beginners should stay indeed with $F=\mathbb C$, but you forgot to mention two important points, over your $F=\mathbb C$, in your easiness plan:
\begin{enumerate}
\item Easier than easiness should be super-easiness, covering all Lie types.

\item Easiness should apply too to spheres, and other types of manifolds.
\end{enumerate}
\end{cat}

Okay, thanks cat, not that I really understand what you are saying, but we will certainly add these two topics, which look quite interesting, on our to-do list. And in any case, this is good to hear, so with $F=\mathbb C$ we are definitely safe, no need for more.

\bigskip

Getting started now, as indicated above, we first need to make sense of the formula $C=span(D)$, for the tensor categories $C$ axiomatized in chapter 10. But here, the answer is quite straightforward, obtained by delinearizing our axioms for $C$, as follows:

\index{category of partitions}
\index{colored integer}

\begin{definition}
Let $P(k,l)$ be the set of partitions between an upper colored integer $k$, and a lower colored integer $l$. A collection of subsets 
$$D=\bigsqcup_{k,l}D(k,l)$$
with $D(k,l)\subset P(k,l)$ is called a category of partitions when it has the following properties:
\begin{enumerate}
\item Stability under the horizontal concatenation, $(\pi,\sigma)\to[\pi\sigma]$.

\item Stability under vertical concatenation $(\pi,\sigma)\to[^\sigma_\pi]$, with matching middle symbols.

\item Stability under the upside-down turning $*$, with switching of colors, $\circ\leftrightarrow\bullet$.

\item Each set $P(k,k)$ contains the identity partition $||\ldots||$.

\item The sets $P(\emptyset,\circ\bullet)$ and $P(\emptyset,\bullet\circ)$ both contain the semicircle $\cap$.
\end{enumerate}
\end{definition} 

As basic examples of such categories, with the verification of the axioms being trivial, in each case, we have the category $P$ of all partitions, the category $NC$ of all noncrossing partitions, the category $P_2$ of all pairings, and the category $NC_2$ of all noncrossing pairings. Observe that these four basic categories form a diagram, as follows:
$$
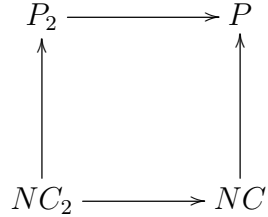
\xymatrix@R=18.5mm@C=16.2mm{
P_2\ar[r]&P\\
NC_2\ar[u]\ar[r]&NC\ar[u]}$$

There are many other examples of such categories, and more on this later. Getting now to the point, the relation with the Tannakian categories of linear maps comes from the fact that we can associate linear maps to the partitions, as follows:

\begin{definition}
Associated to any $\pi\in P(k,l)$ and any $N\in\mathbb N$ is the linear map 
$$T_\pi:(\mathbb C^N)^{\otimes k}\to(\mathbb C^N)^{\otimes l}$$
given by the following formula, with $\{e_1,\ldots,e_N\}$ being the standard basis of $\mathbb C^N$,
$$T_\pi(e_{i_1}\otimes\ldots\otimes e_{i_k})=\sum_{j_1\ldots j_l}\delta_\pi\begin{pmatrix}i_1&\ldots&i_k\\ j_1&\ldots&j_l\end{pmatrix}e_{j_1}\otimes\ldots\otimes e_{j_l}$$
and with the Kronecker symbols $\delta_\pi\in\{0,1\}$ depending on whether the indices fit or not. 
\end{definition}

\index{Kronecker symbol}

To be more precise here, in the definition of the Kronecker symbols, we agree to put the two multi-indices on the two rows of points of the partition, in the obvious way. The Kronecker symbols are then defined by $\delta_\pi=1$ when all the strings of $\pi$ join equal indices, and by $\delta_\pi=0$ otherwise. Observe that all this is independent of the coloring.

\bigskip

In practice now, here are a few basic examples of such linear maps:

\index{semicircle partition}
\index{basic crossing}

\begin{proposition}
The correspondence $\pi\to T_\pi$ has the following properties, where $R$ is the operator mapping $1\to\sum_ie_i\otimes e_i$, and $\Sigma(a\otimes b)=b\otimes a$ is the flip operator:
\begin{enumerate}
\item $T_\cap=R$.

\item $T_\cup=R^*$.

\item $T_{||\ldots||}=id$.

\item $T_{\slash\hskip-1.5mm\backslash}=\Sigma$.
\end{enumerate}
\end{proposition}

\begin{proof}
We can assume if we want that all the upper and lower legs of $\pi$ are colored $\circ$. With this assumption made, the proof goes as follows:

\medskip

(1) We have $\cap\in P_2(\emptyset,\circ\circ)$, and so the corresponding operator is a certain linear map $T_\cap:\mathbb C\to\mathbb C^N\otimes\mathbb C^N$. The formula of this map is as follows:
\begin{eqnarray*}
T_\cap(1)
&=&\sum_{ij}\delta_\cap(i\ j)e_i\otimes e_j\\
&=&\sum_{ij}\delta_{ij}e_i\otimes e_j\\
&=&\sum_ie_i\otimes e_i
\end{eqnarray*}

We recognize here the formula of $R(1)$, and so we have $T_\cap=R$, as claimed.

\medskip

(2) Here we have $\cup\in P_2(\circ\circ,\emptyset)$, and so the corresponding operator is a certain linear form $T_\cap:\mathbb C^N\otimes\mathbb C^N\to\mathbb C$. The formula of this linear form is as follows:
\begin{eqnarray*}
T_\cap(e_i\otimes e_j)
&=&\delta_\cap(i\ j)\\
&=&\delta_{ij}
\end{eqnarray*}

Since this is the same as $R^*(e_i\otimes e_j)$, we have $T_\cup=R^*$, as claimed.

\medskip

(3) Consider indeed the ``identity'' pairing $||\ldots||\in P_2(k,k)$, with $k=\circ\circ\ldots\circ\circ$. The corresponding linear map is then the identity, because we have:
\begin{eqnarray*}
T_{||\ldots||}(e_{i_1}\otimes\ldots\otimes e_{i_k})
&=&\sum_{j_1\ldots j_k}\delta_{||\ldots||}\begin{pmatrix}i_1&\ldots&i_k\\ j_1&\ldots&j_k\end{pmatrix}e_{j_1}\otimes\ldots\otimes e_{j_k}\\
&=&\sum_{j_1\ldots j_k}\delta_{i_1j_1}\ldots\delta_{i_kj_k}e_{j_1}\otimes\ldots\otimes e_{j_k}\\
&=&e_{i_1}\otimes\ldots\otimes e_{i_k}
\end{eqnarray*}

(4) In the case of the basic crossing $\slash\hskip-2.0mm\backslash\in P_2(\circ\circ,\circ\circ)$, the corresponding linear map $T_{\slash\hskip-1.5mm\backslash}:\mathbb C^N\otimes\mathbb C^N\to\mathbb C^N\otimes\mathbb C^N$ can be computed as follows:
\begin{eqnarray*}
T_{\slash\hskip-1.5mm\backslash}(e_i\otimes e_j)
&=&\sum_{kl}\delta_{\slash\hskip-1.5mm\backslash}\begin{pmatrix}i&j\\ k&l\end{pmatrix}e_k\otimes e_l\\
&=&\sum_{kl}\delta_{il}\delta_{jk}e_k\otimes e_l\\
&=&e_j\otimes e_i
\end{eqnarray*}

Thus we obtain the flip operator $\Sigma(a\otimes b)=b\otimes a$, as claimed.
\end{proof}

Summarizing, the correspondence $\pi\to T_\pi$ provides us with some simple formulae for the duality operators $R,R^*$ from chapter 10, and for other important operators, such as the flip map $\Sigma(a\otimes b)=b\otimes a$, and has as well some interesting categorical properties. 

\bigskip

Let us further explore now these properties, and make the link with the Tannakian categories. We have the following result, from \cite{bsp}:

\begin{proposition}
The assignement $\pi\to T_\pi$ is categorical, in the sense that we have
$$T_\pi\otimes T_\sigma=T_{[\pi\sigma]}\quad,\quad 
T_\pi T_\sigma=N^{c(\pi,\sigma)}T_{[^\sigma_\pi]}\quad,\quad 
T_\pi^*=T_{\pi^*}$$
where $c(\pi,\sigma)$ are certain integers, coming from the erased components in the middle.
\end{proposition}

\begin{proof}
The formulae in the statement are all elementary, as follows:

\medskip

(1) The concatenation axiom follows from the following computation:
\begin{eqnarray*}
&&(T_\pi\otimes T_\sigma)(e_{i_1}\otimes\ldots\otimes e_{i_p}\otimes e_{k_1}\otimes\ldots\otimes e_{k_r})\\
&=&\sum_{j_1\ldots j_q}\sum_{l_1\ldots l_s}\delta_\pi\begin{pmatrix}i_1&\ldots&i_p\\j_1&\ldots&j_q\end{pmatrix}\delta_\sigma\begin{pmatrix}k_1&\ldots&k_r\\l_1&\ldots&l_s\end{pmatrix}e_{j_1}\otimes\ldots\otimes e_{j_q}\otimes e_{l_1}\otimes\ldots\otimes e_{l_s}\\
&=&\sum_{j_1\ldots j_q}\sum_{l_1\ldots l_s}\delta_{[\pi\sigma]}\begin{pmatrix}i_1&\ldots&i_p&k_1&\ldots&k_r\\j_1&\ldots&j_q&l_1&\ldots&l_s\end{pmatrix}e_{j_1}\otimes\ldots\otimes e_{j_q}\otimes e_{l_1}\otimes\ldots\otimes e_{l_s}\\
&=&T_{[\pi\sigma]}(e_{i_1}\otimes\ldots\otimes e_{i_p}\otimes e_{k_1}\otimes\ldots\otimes e_{k_r})
\end{eqnarray*}

(2) The composition axiom follows from the following computation:
\begin{eqnarray*}
&&T_\pi T_\sigma(e_{i_1}\otimes\ldots\otimes e_{i_p})\\
&=&\sum_{j_1\ldots j_q}\delta_\sigma\begin{pmatrix}i_1&\ldots&i_p\\j_1&\ldots&j_q\end{pmatrix}
\sum_{k_1\ldots k_r}\delta_\pi\begin{pmatrix}j_1&\ldots&j_q\\k_1&\ldots&k_r\end{pmatrix}e_{k_1}\otimes\ldots\otimes e_{k_r}\\
&=&\sum_{k_1\ldots k_r}N^{c(\pi,\sigma)}\delta_{[^\sigma_\pi]}\begin{pmatrix}i_1&\ldots&i_p\\k_1&\ldots&k_r\end{pmatrix}e_{k_1}\otimes\ldots\otimes e_{k_r}\\
&=&N^{c(\pi,\sigma)}T_{[^\sigma_\pi]}(e_{i_1}\otimes\ldots\otimes e_{i_p})
\end{eqnarray*}

(3) Finally, the involution axiom follows from the following computation:
\begin{eqnarray*}
&&T_\pi^*(e_{j_1}\otimes\ldots\otimes e_{j_q})\\
&=&\sum_{i_1\ldots i_p}<T_\pi^*(e_{j_1}\otimes\ldots\otimes e_{j_q}),e_{i_1}\otimes\ldots\otimes e_{i_p}>e_{i_1}\otimes\ldots\otimes e_{i_p}\\
&=&\sum_{i_1\ldots i_p}\delta_\pi\begin{pmatrix}i_1&\ldots&i_p\\ j_1&\ldots& j_q\end{pmatrix}e_{i_1}\otimes\ldots\otimes e_{i_p}\\
&=&T_{\pi^*}(e_{j_1}\otimes\ldots\otimes e_{j_q})
\end{eqnarray*}

Summarizing, our correspondence is indeed categorical.
\end{proof}

In relation now with quantum groups, we have the following result, also from \cite{bsp}:

\index{Tannakian duality}

\begin{theorem}
Each category of partitions $D=(D(k,l))$ produces a family of quantum groups $G=(G_N)$, one for each $N\in\mathbb N$, via the formula
$$Hom(u^{\otimes k},u^{\otimes l})=span\left(T_\pi\Big|\pi\in D(k,l)\right)$$
which produces a Tannakian category, and the Tannakian duality correspondence. We call such quantum groups, coming from a category of partitions $D=(D(k,l))$, easy.
\end{theorem}

\begin{proof}
This follows indeed from the Tannakian duality established in the previous chapter. Indeed, consider the following collection of linear spaces:
$$C_{kl}=span\left(T_\pi\Big|\pi\in D(k,l)\right)$$

By using the axioms in Definition 11.4, and the categorical properties of the operation $\pi\to T_\pi$, from Propositions 11.6 and 11.7, we deduce that $C=(C_{kl})$ is a Tannakian category. Thus the Tannakian duality applies, and gives the result.
\end{proof}

As a first comment, the terminology in the above result, which is from \cite{bsp}, and has endured so far several unfortunate attempts of change, by younger mathematicians who hate God, life,  beauty and easiness, is motivated by the fact that, from the point of view of Tannakian duality, the above quantum groups are indeed the ``easiest'' ones. 

\bigskip

Observe also that the category of pairings producing an easy quantum group is not unique, for instance because at $N=1$ all the possible categories of pairings produce the same easy group, namely the trivial group $G=\{1\}$. We will be back to this.

\bigskip

And with the above theory in hand, we are ready to go. On the menu, proving that the various groups and quantum groups that we know well are indeed easy. Followed by some further investigations, suggested by Cat 11.3, Comment 11.2 and Answer 11.1.

\section*{11b. Rotations, reflections}

Getting now to applications, let us first talk about rotation groups. We already know a bit about them, from chapters 5 and 9, but it is better to start all over again, with fresh, self-contained definitions. Following Wang \cite{wan}, we have the following result:

\index{unitary quantum group}
\index{free unitary group}
\index{quantum unitary group}
\index{orthogonal quantum group}
\index{quantum orthogonal group}
\index{free orthogonal group}

\begin{theorem}
The following universal algebras are Woronowicz algebras,
$$C(O_N^+)=C^*\left((u_{ij})_{i,j=1,\ldots,N}\Big|u=\bar{u},u^t=u^{-1}\right)$$
$$C(U_N^+)=C^*\left((u_{ij})_{i,j=1,\ldots,N}\Big|u^*=u^{-1},u^t=\bar{u}^{-1}\right)$$
so the underlying quantum spaces $O_N^+,U_N^+$ are compact quantum groups.
\end{theorem}

\begin{proof}
This follows from the elementary fact that if a matrix $u=(u_{ij})$ is orthogonal or biunitary, then so must be the following matrices:
$$(u^\Delta)_{ij}=\sum_ku_{ik}\otimes u_{kj}\quad,\quad 
(u^\varepsilon)_{ij}=\delta_{ij}\quad,\quad
(u^S)_{ij}=u_{ji}^*$$

Thus, we can indeed define Woronowicz algebra morphisms $\Delta,\varepsilon,S$ as in chapter 10, by using the universal properties of $C(O_N^+)$, $C(U_N^+)$, and this gives the result.
\end{proof}

Getting now to Brauer theorems, we have the following result about the above quantum groups, with the convention that $\mathcal P_2$ stands for the matching pairings, meaning the pairings which pair $\circ-\circ$ or $\bullet-\bullet$ on the vertical, and $\circ-\bullet$ on the horizontal:

\index{Brauer theorem}
\index{noncrossing partition}
\index{noncrossing pairing}
\index{partition with even blocks}

\begin{theorem}
The basic classical and quantum rotation groups are all easy,
$$\xymatrix@R=51pt@C=50pt{
O_N^+\ar[r]&U_N^+\\
O_N\ar[u]\ar[r]&U_N\ar[u]}
\qquad\xymatrix@R=25pt@C=50pt{\\ :\\}\qquad
\xymatrix@R=18.5mm@C=16.2mm{
NC_2\ar[d]&\mathcal{NC}_2\ar[l]\ar[d]\\
P_2&\mathcal P_2\ar[l]}$$
with the quantum groups on the left corresponding to the categories on the right.
\end{theorem}

\begin{proof}
This is quite standard, by using Tannakian duality, as follows:

\medskip

(1) The quantum group $U_N^+$ is defined via the following relations:
$$u^*=u^{-1}\quad,\quad 
u^t=\bar{u}^{-1}$$ 

But, via our correspondence between partitions and maps, these relations tell us that the following two operators must be in the associated Tannakian category $C$:
$$T_\pi\quad,\quad \pi={\ }^{\,\cap}_{\circ\bullet}\ ,\,{\ }^{\,\cap}_{\bullet\circ}$$

Thus the associated Tannakian category is $C=span(T_\pi|\pi\in D)$, with:
$$D
=<{\ }^{\,\cap}_{\circ\bullet}\,\,,{\ }^{\,\cap}_{\bullet\circ}>
={\mathcal NC}_2$$

(2) The quantum group $O_N^+\subset U_N^+$ is defined by imposing the following relations:
$$u_{ij}=\bar{u}_{ij}$$

Thus, the following operators must be in the associated Tannakian category $C$:
$$T_\pi\quad,\quad\pi=|^{\hskip-1.32mm\circ}_{\hskip-1.32mm\bullet}\ ,\,|_{\hskip-1.32mm\circ}^{\hskip-1.32mm\bullet}$$

Thus the associated Tannakian category is $C=span(T_\pi|\pi\in D)$, with:
$$D
=<\mathcal{NC}_2,|^{\hskip-1.32mm\circ}_{\hskip-1.32mm\bullet},|_{\hskip-1.32mm\circ}^{\hskip-1.32mm\bullet}>
=NC_2$$

(3) The group $U_N\subset U_N^+$ is defined via the following relations:
$$[u_{ij},u_{kl}]=0\quad,\quad
[u_{ij},\bar{u}_{kl}]=0$$

Thus, the following operators must be in the associated Tannakian category $C$:
$$T_\pi\quad,\quad \pi={\slash\hskip-2.1mm\backslash}^{\hskip-2.5mm\circ\circ}_{\hskip-2.5mm\circ\circ}\ ,\,{\slash\hskip-2.1mm\backslash}^{\hskip-2.5mm\circ\bullet}_{\hskip-2.5mm\bullet\circ}$$

Thus the associated Tannakian category is $C=span(T_\pi|\pi\in D)$, with:
$$D
=<\mathcal{NC}_2,{\slash\hskip-2.1mm\backslash}^{\hskip-2.5mm\circ\circ}_{\hskip-2.5mm\circ\circ},{\slash\hskip-2.1mm\backslash}^{\hskip-2.5mm\circ\bullet}_{\hskip-2.5mm\bullet\circ}>
=\mathcal P_2$$

(4) In order to deal now with $O_N$, we can simply use the following formula: 
$$O_N=O_N^+\cap U_N$$

Indeed, at the categorical level, this formula tells us that the associated Tannakian category is given by $C=span(T_\pi|\pi\in D)$, with:
$$D
=<NC_2,\mathcal P_2>
=P_2$$

Thus, we are led to the conclusions in the statement.
\end{proof}

Getting now towards quantum reflections, we first have the following result, for the quantum permutations, discussed in chapter 6, and then in chapter 9, which provides among others a more reasonable explanation for the liberation operation $S_N\to S_N^+$:

\index{quantum permutation group}
\index{quantum symmetric group}
\index{free permutation group}
\index{free symmetric group}
\index{noncrossing partition}
\index{Brauer theorem}
\index{liberation}
\index{fork partition}

\begin{theorem}
The following hold, for the quantum permutation groups:
\begin{enumerate}
\item The quantum groups $S_N,S_N^+$ are both easy, coming respectively from the categories $P,NC$ of partitions, and noncrossing partitions.  

\item Thus, $S_N\to S_N^+$ is just a regular easy quantum group liberation, coming from $D\to D\cap NC$ at the level of the associated categories of partitions.
\end{enumerate}
\end{theorem}

\begin{proof}
We proceed a bit as for the rotation groups. Let us first prove the result for $S_N^+$. In order to do so, recall that the subgroup $S_N^+\subset O_N^+$ appears as follows:
$$C(S_N^+)=C(O_N^+)\Big\slash\Big<u={\rm magic}\Big>$$

In order to interpret the magic condition, consider the fork partition:
$$Y\in P(2,1)$$

Given a corepresentation $u$, we have the following formulae:
$$(T_Yu^{\otimes 2})_{i,jk}
=\sum_{lm}(T_Y)_{i,lm}(u^{\otimes 2})_{lm,jk}
=u_{ij}u_{ik}$$
$$(uT_Y)_{i,jk}
=\sum_lu_{il}(T_Y)_{l,jk}
=\delta_{jk}u_{ij}$$

We conclude from this that we have the following equivalence:
$$T_Y\in Hom(u^{\otimes 2},u)\iff u_{ij}u_{ik}=\delta_{jk}u_{ij},\forall i,j,k$$

The condition on the right being equivalent to the magic condition, we obtain:
$$C(S_N^+)=C(O_N^+)\Big\slash\Big<T_Y\in Hom(u^{\otimes 2},u)\Big>$$

Thus $S_N^+$ is indeed easy, the corresponding category of partitions being:
$$D=<Y>=NC$$

Finally, observe that this proves the result for $S_N$ too, because from the formula $S_N=S_N^+\cap O_N$ we obtain that the group $S_N$ is easy, coming from the category of partitions $D=<NC,P_2>=P$. Thus, we are led to the conclusions in the statement.
\end{proof}

With this understood, we can get now into quantum reflection groups. We recall from chapter 9 that we have quantum groups $H_N^+,K_N^+$, constructed as follows:
$$C(H_N^+)=C^*\left((u_{ij})_{i,j=1,\ldots,N}\Big|u_{ij}=u_{ij}^*,\,(u_{ij}^2)={\rm magic}\right)$$
$$C(K_N^+)=C^*\left((u_{ij})_{i,j=1,\ldots,N}\Big|[u_{ij},u_{ij}^*]=0,\,(u_{ij}u_{ij}^*)={\rm magic}\right)$$

These are liberations of $H_N,K_N$, and we have $H_N^+=\mathbb Z_2\wr_*S_N^+$ and $K_N^+=\mathbb T\wr_*S_N^+$, in analogy with the formulae $H_N=\mathbb Z_2\wr S_N$ and  $K_N=\mathbb T\wr S_N$ from the classical case.

\bigskip

The point now is that, by working out Brauer theorems for these quantum groups, we reach to something quite conceptual, namely having a nice 3D cube, as follows:

\index{standard cube}
\index{quantum rotation}
\index{quantum reflection}
\index{Brauer theorem}
\index{quantum reflection group}
\index{free reflection group}
\index{complex reflection group}
\index{free complex reflection group}

\begin{theorem}
The main quantum rotation and reflection groups,
$$
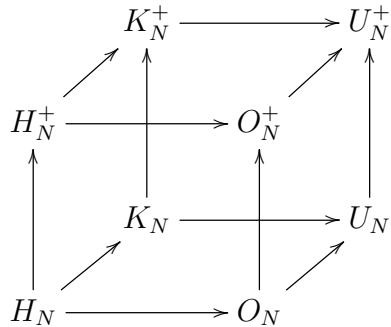
\xymatrix@R=18pt@C=18pt{
&K_N^+\ar[rr]&&U_N^+\\
H_N^+\ar[rr]\ar[ur]&&O_N^+\ar[ur]\\
&K_N\ar[rr]\ar[uu]&&U_N\ar[uu]\\
H_N\ar[uu]\ar[ur]\ar[rr]&&O_N\ar[uu]\ar[ur]
}$$
are all easy, the corresponding categories of partitions being as follows,
$$
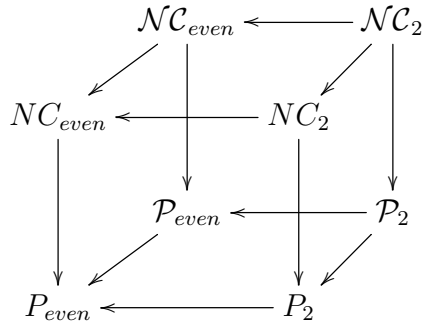
\xymatrix@R=19pt@C5pt{
&\mathcal{NC}_{even}\ar[dl]\ar[dd]&&\mathcal {NC}_2\ar[dl]\ar[ll]\ar[dd]\\
NC_{even}\ar[dd]&&NC_2\ar[dd]\ar[ll]\\
&\mathcal P_{even}\ar[dl]&&\mathcal P_2\ar[dl]\ar[ll]\\
P_{even}&&P_2\ar[ll]
}$$
with on top, the symbol $NC$ standing everywhere for noncrossing partitions.
\end{theorem}

\begin{proof}
This is something that we know for the objects on the right, and for those on the left, this follows by using the same arguments as for $S_N,S_N^+$, as follows:

\medskip

(1) We know that $H_N^+\subset O_N^+$ appears via the cubic relations, namely:
$$u_{ij}u_{ik}=u_{ji}u_{ki}=0\quad,\quad\forall j\neq k$$

Our claim now is that, in Tannakian terms, these relations reformulate as follows, with $H\in P(2,2)$ being the 1-block partition, joining all 4 points:
$$T_H\in End(u^{\otimes 2})$$

In order to prove our claim, observe first that we have, by definition of $T_H$:
$$T_H(e_i\otimes e_j)=\delta_{ij}e_i\otimes e_i$$

With this formula in hand, we have the following computation:
\begin{eqnarray*}
T_Hu^{\otimes 2}(e_i\otimes e_j\otimes1)
&=&T_H\left(\sum_{abij}e_{ai}\otimes e_{bj}\otimes u_{ai}u_{bj}\right)(e_i\otimes e_j\otimes1)\\
&=&T_H\sum_{ab}e_a\otimes e_b\otimes u_{ai}u_{bj}\\
&=&\sum_ae_a\otimes e_a\otimes u_{ai}u_{aj}
\end{eqnarray*}

On the other hand, we have as well the following computation:
\begin{eqnarray*}
u^{\otimes 2}T_H(e_i\otimes e_j\otimes1)
&=&\delta_{ij}u^{\otimes 2}(e_i\otimes e_j\otimes1)\\
&=&\delta_{ij}\left(\sum_{abij}e_{ai}\otimes e_{bj}\otimes u_{ai}u_{bj}\right)(e_i\otimes e_j\otimes1)\\
&=&\delta_{ij}\sum_{ab}e_a\otimes e_b\otimes u_{ai}u_{bi}
\end{eqnarray*}

We conclude that $T_Hu^{\otimes 2}=u^{\otimes 2}T_H$ means that $u$ is cubic, as desired. Thus, our claim is proved. But this shows that $H_N^+$ is easy, coming from the following category:
$$D=<H>=NC_{even}$$

(2) Regarding now the full quantum reflection group $K_N^+$, the proof here is similar, leading this time to the category $\mathcal{NC}_{even}$ of noncrossing matching partitions. 

\medskip

(3) Finally, both the passages $H_N^+\to H_N$ and $K_N^+\to K_N$ come by adding the crossing, at the level of partitions, and this leads to the conclusions in the statement.
\end{proof}

So long for basic Brauer theorems. We should mention that much more is true, with the general reflection groups $H_N^s,H_N^{s+}$ being easy too, and with the same going for the bistochastic groups $B_N,C_N$ and $B_N^+,C_N^+$. We will be back to this, later in this chapter.

\section*{11c. Liberation theory}

With the above results in hand, we can now review the half-liberation and twisting operations, that we met in Part II, via some ad-hoc constructions. Let us start with the half-liberation operation. We have here the following key result, to start with: 

\begin{theorem}
We have an intermediate easy quantum group $U_N\subset U_N^*\subset U_N^+$, given by the following formula, and called half-classical unitary group,
$$C(U_N^*)=C(U_N^+)\Big/\left<abc=cba\Big|\forall a,b,c\in\{u_{ij},u_{ij}^*\}\right>$$
corresponding to the category of matching pairings $\mathcal P_2^*$ having the property that when relabelling clockwise the legs $\circ\bullet\circ\bullet\ldots$, the formula $\#\circ=\#\bullet$ must hold in each block.
\end{theorem}

\begin{proof}
Here the fact that $U_N^*$ as constructed above is indeed a quantum group, lying as a proper intermediate subgroup $U_N\subset U_N^*\subset U_N^+$, can be checked via a routine computation, but the best is to view this via Tannakian duality. Indeed, the  half-commutation relations $abc=cba$ come from the map $T_{\slash\hskip-1.6mm\backslash\hskip-1.1mm|\hskip0.5mm}$ associated to the half-classical crossing:
$$\slash\hskip-2.0mm\backslash\hskip-1.7mm|\hskip0.5mm\in P(3,3)$$

Thus, by Tannakian duality, we are led to the first conclusions in the statement, and with the category of partitions associated to $U_N^*$ being as follows, with the convention that the symbol $\slash\hskip-2.0mm\backslash\hskip-1.7mm|\hskip0.5mm$ stands here for all 8 possible colorings of the diagram $\slash\hskip-2.0mm\backslash\hskip-1.7mm|\hskip0.5mm$:
$$D=<\slash\hskip-2.0mm\backslash\hskip-1.7mm|\hskip0.5mm>$$

Regarding now the explicit computation of $D$, observe that no matter how we color the legs of $\slash\hskip-2.0mm\backslash\hskip-1.7mm|\hskip0.5mm$, of course as for strings to join $\circ-\circ$ or $\bullet-\bullet$, we have a matching pairing, having in addition the property that when relabelling clockwise the legs $\circ\bullet\circ\bullet\ldots$, the formula $\#\circ=\#\bullet$ must hold in each block. Thus, we have an inclusion as follows:
$$D\subset\mathcal P_2^*$$

On the other hand, by doing some standard combinatorics, we see that any element of $\mathcal P_2^*$ can be decomposed into a composition of diagrams of type $\slash\hskip-2.0mm\backslash\hskip-1.7mm|\hskip0.5mm$, appearing with all its 8 possible colorings, as above. Thus, we have as well an inclusion as follows:
$$\mathcal P_2^*\subset D$$

But this shows that we have $D=\mathcal P_2^*$, which proves the last assertion.
\end{proof}

Now with the above construction in hand, we can perform a ``cutting the cube'' operation, and we are led to the following statement, improving what we have so far:

\begin{theorem}
We have easy quantum groups as follows, obtained via the commutation relations $abc=cba$, applied to the standard coordinates and their adjoints,
$$
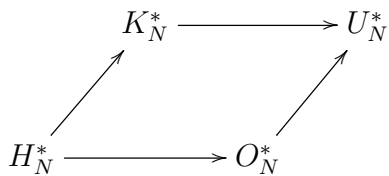
\xymatrix@R=32pt@C=18pt{
&K_N^*\ar[rr]&&U_N^*\\
H_N^*\ar[rr]\ar[ur]&&O_N^*\ar[ur]
}$$
which fit horizontally, in the middle, into the diagram of basic easy quantum groups
$$
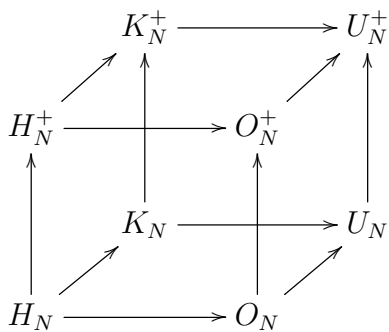
\xymatrix@R=18pt@C=18pt{
&K_N^+\ar[rr]&&U_N^+\\
H_N^+\ar[rr]\ar[ur]&&O_N^+\ar[ur]\\
&K_N\ar[rr]\ar[uu]&&U_N\ar[uu]\\
H_N\ar[uu]\ar[ur]\ar[rr]&&O_N\ar[uu]\ar[ur]
}$$
with the enlarged diagram being an intersection/easy generation diagram.
\end{theorem}

\begin{proof}
There are several things going on here, and we will be quite brief:

\medskip

(1) First, the fact that we have indeed quantum groups as in the statement, which are all easy, follows from our previous results, and from Theorem 11.13. 

\medskip

(2) Once again by using our previous results and Theorem 11.13, we conclude as well that the categories of partitions for our new quantum groups are as follows:
$$
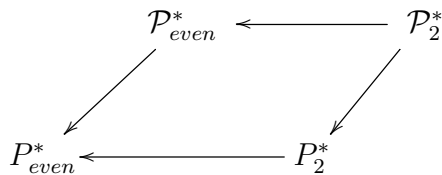
\xymatrix@R=32pt@C=18pt{
&\ \mathcal P_{even}^*\ \ar[dl]&&\ \mathcal P_2^*\ar[ll]\ar[dl]\\
P_{even}^*\!&&P_2^*\ar[ll]
}$$

(3) Now let us insert this latter square diagram into the standard cube of categories of partitions, on the horizontal, in the middle:
$$
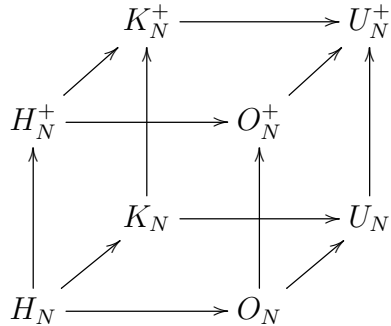
\xymatrix@R=18pt@C=18pt{
&K_N^+\ar[rr]&&U_N^+\\
H_N^+\ar[rr]\ar[ur]&&O_N^+\ar[ur]\\
&K_N\ar[rr]\ar[uu]&&U_N\ar[uu]\\
H_N\ar[uu]\ar[ur]\ar[rr]&&O_N\ar[uu]\ar[ur]
}$$

(4) We obtain in this way an intersection and generation diagram. Thus, the diagram formed by quantum groups is an intersection/easy generation diagram, as stated.
\end{proof}

Many other things can be said here, for instance with the construction of ``hybrid'' quantum groups between real and complex too, again cutting the standard cube. 

\bigskip

Also, in the non-easy case, we have available the modern approach to half-liberation, by Bichon and Dubois-Violette \cite{bd1}, based on crossed products and $2\times2$ matrix models. We will be back to this later, when discussing more in detail the matrix models.

\bigskip

We would like to end this section with the discussion of the bistochastic groups and quantum groups $B_N,C_N,B_N^+,C_N^+$, which are quite interesting objects. 

\bigskip

These groups and quantum groups are all easy, and are not half-liberable, so our discussion here will be quite simple and pleasant, basically consisting in explaining what $B_N,C_N,B_N^+,C_N^+$ are, and then working out their easiness property. Let us start with:

\begin{proposition}
We have groups as follows:
\begin{enumerate}
\item $B_N\subset O_N$, consisting of the orthogonal matrices which are bistochastic.

\item $C_N\subset U_N$, consisting of the unitary matrices which are bistochastic.
\end{enumerate}
\end{proposition}

\begin{proof}
We know that the sets of bistochastic matrices $B_N,C_N$ in the statement appear as follows, with $\xi$ being the all-one vector:
$$B_N=\left\{U\in O_N\Big|U\xi=\xi\right\}\quad,\quad 
C_N=\left\{U\in U_N\Big|U\xi=\xi\right\}$$

It is then clear that both $B_N,C_N$ are stable under the multiplication, contain the unit, and are stable by inversion. Thus, we have indeed groups, as stated.
\end{proof}

Getting back now to quantum groups, we would like to talk about the free analogues $B_N^+,C_N^+$ of the orthogonal and unitary bistochastic groups. Following \cite{bsp}, we have the following result, including as well the above construction of $B_N,C_N$:

\begin{theorem}
We have the following groups and quantum groups:
\begin{enumerate}
\item $B_N\subset O_N$, consisting of the orthogonal matrices which are bistochastic.

\item $C_N\subset U_N$, consisting of the unitary matrices which are bistochastic.

\item $B_N^+\subset O_N^+$, coming via $u\xi=\xi$, where $\xi$ is the all-one vector.

\item $C_N^+\subset U_N^+$, coming via $u\xi=\xi$, where $\xi$ is the all-one vector.
\end{enumerate}
Also, we have inclusions $B_N\subset B_N^+$ and $C_N\subset C_N^+$, which are both liberations.
\end{theorem}

\begin{proof}
There are several things to be proved, the idea being as follows:

\medskip

(1) We already know from Proposition 11.15 that $B_N,C_N$ are indeed groups. In what regards now $B_N^+,C_N^+$, these appear by definition as follows:
$$C(B_N^+)=C(O_N^+)\Big/\Big<\xi\in Fix(u)\Big>$$
$$C(C_N^+)=C(U_N^+)\Big/\Big<\xi\in Fix(u)\Big>$$

But since the relation $\xi\in Fix(u)$ is categorical, we have indeed quantum groups.

\medskip

(2) In what regards the last assertion, since we already know that $O_N\subset O_N^+$ and $U_N\subset U_N^+$ are liberations, we must prove that we have isomorphisms as follows:
$$C(B_N)=C(O_N)\Big/\Big<\xi\in Fix(u)\Big>$$
$$C(C_N)=C(U_N)\Big/\Big<\xi\in Fix(u)\Big>$$

But these isomorphisms are both clear from the formulae of $B_N,C_N$ in (1).
\end{proof}

The above might seem a bit puzzling, because the ``basic'' compact groups can be of type ABCD, and there was no mention there, in that Lie algebra theory and classification, of bistochastic groups as above. Good point, and in answer, we have indeed:

\begin{theorem}
We have isomorphisms as follows:
\begin{enumerate}
\item $B_N\simeq O_{N-1}$.

\item $B_N^+\simeq O_{N-1}^+$.

\item $C_N\simeq U_{N-1}$.

\item $C_N^+\simeq U_{N-1}^+$.
\end{enumerate}
\end{theorem}

\begin{proof}
Let us pick indeed a matrix $F\in U_N$ satisfying the following condition, where $\xi$ is the all-one vector, and $e_0=(1,0,\ldots,0)$ is the first vector of the standard basis of $\mathbb C^N$, written with indices $0,1,\ldots,N-1$, as usual in discrete Fourier analysis:
$$Fe_0=\frac{1}{\sqrt{N}}\xi$$

Such matrices exist of course, the basic example being the Fourier matrix:
$$F_N=\frac{1}{\sqrt{N}}(w^{ij})_{ij}\quad,\quad w=e^{2\pi i/N}$$

We have then the following computation, for any corepresentation $u$:
\begin{eqnarray*}
u\xi=\xi
&\iff&uFe_0=Fe_0\\
&\iff&F^*uFe_0=e_0\\
&\iff&F^*uF=diag(1,w)
\end{eqnarray*}

Thus we have an isomorphism given by $w_{ij}\to(F^*uF)_{ij}$, as desired.
\end{proof}

Getting now to the real thing, the point is that our newly constructed bistochastic groups and quantum groups are all easy, and following \cite{bsp}, \cite{twe}, we have:

\begin{theorem}
The classical and quantum bistochastic groups are all easy, with the quantum groups on the left corresponding to the categories on the right,
$$\xymatrix@R=51pt@C=50pt{
B_N^+\ar[r]&C_N^+\\
B_N\ar[u]\ar[r]&C_N\ar[u]}
\qquad\xymatrix@R=25pt@C=50pt{\\ :\\}\qquad
\xymatrix@R=18.5mm@C=16.2mm{
NC_{12}\ar[d]&\mathcal{NC}_{12}\ar[l]\ar[d]\\
P_{12}&\mathcal P_{12}\ar[l]}$$
where the symbol $12$ stands for ``category of singletons and pairings''.
\end{theorem}

\begin{proof}
This comes from the fact that the all-one vector $\xi$ used in the constructions in Theorem 11.16 is the vector associated to the singleton partition:
$$\xi=T_|$$

Indeed, we obtain that $B_N,C_N,B_N^+,C_N^+$ are inded easy, appearing from the categories of partitions for $O_N,U_N,O_N^+,U_N^+$, by adding singletons. Thus, we get the result.
\end{proof}

\section*{11d. Character laws}

This book was supposed to be about algebra, but we would like to end this chapter with some analysis. Here is something analytic, which is truly beautiful:

\index{main character}
\index{derangement}
\index{Poisson law}
\index{random permutation}

\begin{theorem}
Consider the symmetric group $S_N$, regarded as a compact group of matrices, $S_N\subset O_N$, via the standard permutation matrices.
\begin{enumerate}
\item The main character $\chi\in C(S_N)$, defined as usual as $\chi=\sum_iu_{ii}$, counts the number of fixed points, $\chi(\sigma)=\#\{i|\sigma(i)=i\}$.

\item The probability for a permutation $\sigma\in S_N$ to be a derangement, meaning to have no fixed points at all, becomes, with $N\to\infty$, equal to $1/e$.

\item The law of the main character $\chi\in C(S_N)$ becomes with $N\to\infty$ the Poisson law $p_1=\frac{1}{e}\sum_k\delta_k/k!$, with respect to the counting measure.
\end{enumerate}
\end{theorem}

\begin{proof}
This is something very classical, the proof being as follows:

\medskip

(1) We have indeed the following computation, which gives the result:
$$\chi(\sigma)
=\sum_iu_{ii}(\sigma)
=\sum_i\delta_{\sigma(i)i}
=\#\left\{i\Big|\sigma(i)=i\right\}$$

(2) We use the inclusion-exclusion principle. Consider the following sets:
$$S_N^i=\left\{\sigma\in S_N\Big|\sigma(i)=i\right\}$$

The probability that we are interested in is then given by:
\begin{eqnarray*}
P(\chi=0)
&=&\frac{1}{N!}\left(|S_N|-\sum_i|S_N^i|+\sum_{i<j}|S_N^i\cap S_N^j|-\sum_{i<j<k}|S_N^i\cap S_N^j\cap S_N^k|+\ldots\right)\\
&=&\frac{1}{N!}\sum_{r=0}^N(-1)^r\sum_{i_1<\ldots<i_r}(N-r)!\\
&=&\frac{1}{N!}\sum_{r=0}^N(-1)^r\binom{N}{r}(N-r)!\\
&=&\sum_{r=0}^N\frac{(-1)^r}{r!}
\end{eqnarray*}

Since we have here the expansion of $1/e$, this gives the result.

\medskip

(3) This follows by generalizing the computation in (2). Indeed, we get:
$$\lim_{N\to\infty}P(\chi=k)=\frac{1}{k!e}$$

Thus, we obtain in the limit a Poisson law of parameter 1, as stated.
\end{proof}

The above result is very beautiful, but we can in fact do even better, as follows:

\index{truncated character}
\index{derangement}
\index{Poisson law}
\index{random permutation}

\begin{theorem}
For the symmetric group $S_N\subset O_N$, the truncated character
$$\chi_t=\sum_{i=1}^{[tN]}u_{ii}$$
with $t\in[0,1]$ follows with $N\to\infty$ the Poisson law $p_t=e^{-t}\sum_k\delta_kt^k/k!$.
\end{theorem}

\begin{proof}
This follows by suitably modifying the proof of Theorem 11.19. Indeed, according to our definition of $\chi_t$, we first have the following formula:
$$\chi(\sigma)=\#\left\{i\in\{1,\ldots,[tN]\}\Big|\sigma(i)=i\right\}$$

Then, the inclusion-exclusion principle gives the following formula:
$$\lim_{N\to\infty}P(\chi_t=k)=\frac{t^k}{k!e^t}$$

Thus, we obtain in the limit a Poisson law of parameter $t$, as stated.
\end{proof}

Summarizing, we have some interesting theory going on here, and this is probably just the tip of the iceberg. So, let us formulate the following question: 

\begin{question}
Given $G\subset U_N^+$, what is the law of the truncated character
$$\chi_t=\sum_{i=1}^{[tN]}u_{ii}$$
with respect to Haar integration? Also, in the easy case, what is the asymptotic law
$$\lim_{N\to\infty}law(\chi_t)$$
say regarded as a formal measure, determined by its moments?
\end{question}

Our goal in what follows will be that of answering this question, by extending what we have about $S_N$ to all the easy groups and easy quantum groups that we know.

\bigskip

In order to reach to our goals, we first need to know how to integrate on the easy quantum groups. The formula here, coming from Peter-Weyl theory, is as follows:

\index{Gram matrix}
\index{Weingarten matrix}
\index{Weingarten formula}
\index{Haar functional}

\begin{theorem}
Assuming that a closed subgroup $G\subset U_N^+$ is easy, coming from a category of partitions $D\subset P$, we have the Weingarten formula
$$\int_Gu_{i_1j_1}^{e_1}\ldots u_{i_kj_k}^{e_k}=\sum_{\pi,\sigma\in D(k)}\delta_\pi(i)\delta_\sigma(j)W_{kN}(\pi,\sigma)$$
where $\delta\in\{0,1\}$ are the usual Kronecker type symbols, and where the Weingarten matrix $W_{kN}=G_{kN}^{-1}$ is the inverse of the Gram matrix $G_{kN}(\pi,\sigma)=N^{|\pi\vee\sigma|}$. 
\end{theorem}

\begin{proof}
We know from chapter 10 that the integrals in the statement form altogether the orthogonal projection $P^k$ onto the following space:
$$Fix(u^{\otimes k})=span\left(\xi_\pi\Big|\pi\in D(k)\right)$$

In order to prove the result, consider the following linear map:
$$E(x)=\sum_{\pi\in D(k)}<x,\xi_\pi>\xi_\pi$$

By a standard linear algebra computation, it follows that we have $P=WE$, where $W$ is the inverse on $Fix(u^{\otimes k})$ of the restriction of $E$. But this restriction is the linear map given by $G_{kN}$, and so $W$ is the linear map given by $W_{kN}$, and this gives the result.
\end{proof}

In relation now with truncated characters, we have the following formula:

\index{truncated character}

\begin{proposition}
The moments of truncated characters are given by the formula
$$\int_G(u_{11}+\ldots +u_{ss})^k=Tr(W_{kN}G_{ks})$$
where $G_{kN}$ and $W_{kN}=G_{kN}^{-1}$ are the associated Gram and Weingarten matrices.
\end{proposition}

\begin{proof}
We have indeed the following computation:
\begin{eqnarray*}
\int_G(u_{11}+\ldots +u_{ss})^k
&=&\sum_{i_1=1}^{s}\ldots\sum_{i_k=1}^s\int u_{i_1i_1}\ldots u_{i_ki_k}\\
&=&\sum_{\pi,\sigma\in D(k)}W_{kN}(\pi,\sigma)\sum_{i_1=1}^{s}\ldots\sum_{i_k=1}^s\delta_\pi(i)\delta_\sigma(i)\\
&=&\sum_{\pi,\sigma\in D(k)}W_{kN}(\pi,\sigma)G_{ks}(\sigma,\pi)\\
&=&Tr(W_{kN}G_{ks})
\end{eqnarray*}

Thus, we have obtained the formula in the statement.
\end{proof}

All this is very good, and normally we have here what is needed in order to answer Question 11.21, in the easy case. In practice, however, before doing so, we will need some training in probability theory.  Let us start with something very standard, as follows:

\begin{definition}
Let $A$ be a $C^*$-algebra, given with a trace $tr:A\to\mathbb C$.
\begin{enumerate}
\item The elements $a\in A$ are called random variables.

\item The moments of such a variable are the numbers $M_k(a)=tr(a^k)$.

\item The law of such a variable is the functional $\mu:P\to tr(P(a))$.
\end{enumerate}
\end{definition}

Here $k=\circ\bullet\bullet\circ\ldots$ is by definition a colored integer, and the corresponding powers $a^k$ are defined by the following formulae, and multiplicativity: 
$$a^\emptyset=1\quad,\quad
a^\circ=a\quad,\quad
a^\bullet=a^*$$

As for the polynomial $P$, this is a noncommuting $*$-polynomial in one variable:
$$P\in\mathbb C<X,X^*>$$

Observe that the law is uniquely determined by the moments, because we have:
$$P(X)=\sum_k\lambda_kX^k\implies\mu(P)=\sum_k\lambda_kM_k(a)$$

Generally speaking, the above definition is something quite abstract, but there is no other way of doing things, at least at this level of generality. However, in certain special cases, the formalism simplifies, and we recover more familiar objects, as follows:

\index{spectral theorem}
\index{spectral measure}

\begin{proposition}
Assuming that $a\in A$ is normal, $aa^*=a^*a$, its law corresponds to a probability measure on its spectrum $\sigma(a)\subset\mathbb C$, according to the following formula:
$$tr(P(a))=\int_{\sigma(a)}P(x)d\mu(x)$$
When the trace is faithful we have $supp(\mu)=\sigma(a)$. Also, in the particular case where the variable is self-adjoint, $a=a^*$, this law is a real probability measure.
\end{proposition}

\begin{proof}
Since the $C^*$-algebra $<a>$ generated by $a$ is commutative, the Gelfand theorem applies to it, and gives an identification of $C^*$-algebras, as follows:
$$<a>=C(\sigma(a))$$

Now by using the Riesz theorem, the restriction of $tr$ to this algebra must come from a probability measure $\mu$ as in the statement, and this gives all the assertions.
\end{proof}

Getting now where we wanted to get, we want our ``quantum probability'' theory to apply to the two main cases that we have in mind, namely classical and free. So, following Voiculescu \cite{vdn}, let us introduce the following two notions of independence:

\index{independence}
\index{freeness}

\begin{definition}
Two subalgebras $A,B\subset C$ are called independent when
$$tr(a)=tr(b)=0\implies tr(ab)=0$$
holds for any $a\in A$ and $b\in B$, and free when
$$tr(a_i)=tr(b_i)=0\implies tr(a_1b_1a_2b_2\ldots)=0$$
holds for any $a_i\in A$ and $b_i\in B$.
\end{definition}
 
In other words, what we have here is a straightforward noncommutative extension of the usual notion of independence, along with a natural free analogue of it. 

\bigskip

In order to understand now what is going on, let us discuss some basic models for independence and freeness. We have the following result, which clarifies things:

\begin{proposition}
Given two algebras $(A,tr)$ and $(B,tr)$, the following hold:
\begin{enumerate}
\item $A,B$ are independent inside their tensor product $A\otimes B$.

\item $A,B$ are free inside their free product $A*B$.
\end{enumerate}
\end{proposition}

\begin{proof}
Both the assertions are clear from definitions, after some standard discussion regarding the tensor product and free product trace. See Voiculescu \cite{vdn}.
\end{proof}

In relation with groups and algebra, we have the following result:

\begin{proposition}
We have the following results, valid for group algebras:
\begin{enumerate}
\item $C^*(\Gamma),C^*(\Lambda)$ are independent inside $C^*(\Gamma\times\Lambda)$.

\item $C^*(\Gamma),C^*(\Lambda)$ are free inside $C^*(\Gamma*\Lambda)$.
\end{enumerate}
\end{proposition}

\begin{proof}
This follows from the general results in Proposition 11.27, along with the following two isomorphisms, which are both standard:
$$C^*(\Gamma\times\Lambda)=C^*(\Lambda)\otimes C^*(\Gamma)\quad,\quad 
C^*(\Gamma*\Lambda)=C^*(\Lambda)*C^*(\Gamma)$$

Alternatively, we can prove this directly, by using the fact that each algebra is spanned by the corresponding group elements, and checking the result on group elements.
\end{proof}

In order to study independence and freeness, our main tool will be: 

\begin{theorem}
The convolution is linearized by the log of the Fourier transform,
$$F_f(x)=E(e^{ixf})$$
and the free convolution is linearized by the $R$-transform, which is given by
$$G_\mu(\xi)=\int_\mathbb R\frac{d\mu(t)}{\xi-t}\implies G_\mu\left(R_\mu(\xi)+\frac{1}{\xi}\right)=\xi$$
and so is the inverse of the Cauchy transform, up to a $\xi^{-1}$ factor.
\end{theorem}

\begin{proof}
For the first assertion, if $f,g$ are independent, we have indeed:
\begin{eqnarray*}
F_{f+g}(x)
&=&\int_\mathbb Re^{ixz}d(\mu_f*\mu_g)(z)\\
&=&\int_{\mathbb R\times\mathbb R}e^{ix(z+t)}d\mu_f(z)d\mu_g(t)\\
&=&F_f(x)F_g(x)
\end{eqnarray*}

For the second assertion, we need a good model for free convolution, and the best is to use the semigroup algebra of the free semigroup on two generators:
$$A=C^*(\mathbb N*\mathbb N)$$

Indeed, we have some freeness in the semigroup setting, a bit in the same way as for the group algebras $C^*(\Gamma*\Lambda)$, from Proposition 11.28. In addition to this fact, and to what happens in the group algebra case, the following things happen:

\medskip

(1) The variables of type $S^*+f(S)$, with $S\in C^*(\mathbb N)$ being the shift, and with $f\in\mathbb C[X]$ being a polynomial, model in moments all the distributions $\mu:\mathbb C[X]\to\mathbb C$. This is indeed something elementary, which can be checked via a direct algebraic computation.

\medskip

(2) Given $f,g\in\mathbb C[X]$, the variables $S^*+f(S)$ and $T^*+g(T)$, where $S,T\in C^*(\mathbb N*\mathbb N)$ are the shifts corresponding to the generators of $\mathbb N*\mathbb N$, are free, and their sum has the same law as $S^*+(f+g)(S)$. This follows indeed by using a $45^\circ$ argument.

\medskip

(3) But with this in hand, we can see that $\mu\to f$ linearizes the free convolution. We are therefore left with a computation inside $C^*(\mathbb N)$, whose conclusion is that $R_\mu=f$ can be recaptured from $\mu$ via the Cauchy transform $G_\mu$, as stated. See \cite{vdn}.
\end{proof}

As a first main result now, we have the following statement, called Central Limit Theorem (CLT), which in the free case is due to Voiculescu \cite{vdn}:

\begin{theorem}[CLT]
Given self-adjoint random variables $x_1,x_2,x_3,\ldots$ which are i.i.d./f.i.d., centered, with variance $t>0$, we have, with $n\to\infty$, in moments,
$$\frac{1}{\sqrt{n}}\sum_{i=1}^nx_i\sim g_t/\gamma_t$$
where the limiting laws $g_t/\gamma_t$ are the following measures,
$$g_t=\frac{1}{\sqrt{2\pi t}}e^{-x^2/2t}dx\quad,\quad 
\gamma_t=\frac{1}{2\pi t}\sqrt{4t^2-x^2}dx$$
called normal, or Gaussian, and Wigner semicircle law of parameter $t$.
\end{theorem}

\begin{proof}
This is routine, by using the linearization properties of the Fourier transform and the $R$-transform from Theorem 11.29, and for details here, we refer to any classical probability book for the classical result, and to \cite{vdn} for the free result.
\end{proof}

Next, we have the following complex version of the CLT:

\begin{theorem}[CCLT]
Given variables $x_1,x_2,x_3,\ldots$ which are i.i.d./f.i.d., centered, with variance $t>0$, we have, with $n\to\infty$, in moments,
$$\frac{1}{\sqrt{n}}\sum_{i=1}^nx_i\sim G_t/\Gamma_t$$
where $G_t/\Gamma_t$ are the complex normal and Voiculescu circular law of parameter $t$, given by:
$$G_t=law\left(\frac{1}{\sqrt{2}}(a+ib)\right)\quad,\quad 
\Gamma_t=law\left(\frac{1}{\sqrt{2}}(\alpha+i\beta)\right)$$
where $a,b/\alpha,\beta$ are independent/free, each following the law $g_t/\gamma_t$.
\end{theorem}

\begin{proof}
This follows indeed from the CLT, by taking real and imaginary parts of all the variables involved, and for details and more here, including the combinatorics of the Voiculescu circular law $\Gamma_t$, which is quite subtle, we refer again to \cite{vdn}.
\end{proof}

We denote by $\boxplus$ the free convolution operation for the real probability measures, given by the fact that we have the following formula, when $a,b$ are free:
$$\mu_{a+b}=\mu_a\boxplus\mu_b$$

With this convention, and adding to our collection of limiting results, we have the following discrete version of the CLT, called Poisson Limit Theorem (PLT):

\begin{theorem}[PLT]
The following Poisson limits converge, for any $t>0$,
$$p_t=\lim_{n\to\infty}\left(\left(1-\frac{t}{n}\right)\delta_0+\frac{t}{n}\delta_1\right)^{*n}\quad,\quad 
\pi_t=\lim_{n\to\infty}\left(\left(1-\frac{t}{n}\right)\delta_0+\frac{t}{n}\delta_1\right)^{\boxplus n}$$
the limiting measures being the Poisson law $p_t$, and the Marchenko-Pastur law $\pi_t$,
$$p_t=\frac{1}{e^t}\sum_{k=0}^\infty\frac{t^k\delta_k}{k!}\quad,\quad 
\pi_t=\max(1-t,0)\delta_0+\frac{\sqrt{4t-(x-1-t)^2}}{2\pi x}\,dx$$
with at $t=1$, the Marchenko-Pastur law being $\pi_1=\frac{1}{2\pi}\sqrt{4x^{-1}-1}\,dx$. 
\end{theorem}

\begin{proof}
This is again routine, by using the Fourier and $R$-transform, and as before we refer here to any classical probability book, and to \cite{vdn}.
\end{proof}

Finally, we have the following ``compound'' generalization of the PLT:

\begin{theorem}[CPLT]
Given a compactly supported positive measure $\rho$, of mass $c=mass(\rho)$, the following compound Poisson limits converge,
$$p_\rho=\lim_{n\to\infty}\left(\left(1-\frac{c}{n}\right)\delta_0+\frac{1}{n}\rho\right)^{*n}
\quad,\quad 
\pi_\rho=\lim_{n\to\infty}\left(\left(1-\frac{c}{n}\right)\delta_0+\frac{1}{n}\rho\right)^{\boxplus n}$$
and if we write $\rho=\sum_{i=1}^sc_i\delta_{z_i}$ with $c_i>0$ and $z_i\in\mathbb R$ we have the formula
$$p_\rho/\pi_\rho={\rm law}\left(\sum_{i=1}^sz_i\alpha_i\right)$$
where the variables $\alpha_i$ are Poisson/free Poisson$(c_i)$, independent/free.
\end{theorem}

\begin{proof}
As before, this follows by using the Fourier and the $R$-transform, and details can be found in any probability book, and in \cite{vdn}.
\end{proof}

So long for limiting results in classical and free probability. To finish with, for our purposes here, we will need the following notions, coming from Theorem 11.33:

\begin{definition}
The Bessel and free Bessel laws, depending on parameters $s\in\mathbb N\cup\{\infty\}$ and $t>0$, are the following compound Poisson and free Poisson laws,
$$b^s_t=p_{t\varepsilon_s}\quad,\quad 
\beta^s_t=\pi_{t\varepsilon_s}$$
with $\varepsilon_s$ being the uniform measure on the $s$-th roots of unity. In particular:
\begin{enumerate}
\item At $s=1$ we recover the Poisson laws $p_t,\pi_t$.

\item At $s=2$ we have the real Bessel laws $b_t,\beta_t$.

\item At $s=\infty$ we have the complex Bessel laws $B_t,\mathfrak B_t$.
\end{enumerate}
\end{definition}

Here the terminology comes from the fact that the density of the measure $b_t$ from (2) is a Bessel function of the first kind, the formula, from \cite{bbc}, being as follows:
$$b_t=e^{-t}\sum_{r=-\infty}^\infty\delta_r\sum_{p=0}^\infty \frac{(t/2)^{|r|+2p}}{(|r|+p)!p!}$$

Good news, with the above general theory in hand, we can now formulate our truncated character results for the main examples of easy quantum groups, as follows:

\index{standard cube}
\index{truncated character}

\begin{theorem}
For the main quantum rotation and reflection groups,
$$
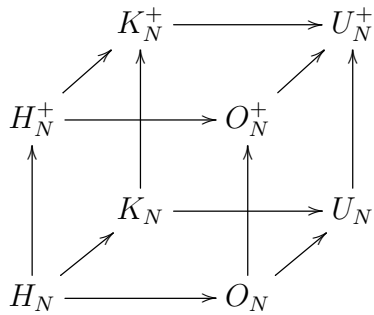
\xymatrix@R=16pt@C=16pt{
&K_N^+\ar[rr]&&U_N^+\\
H_N^+\ar[rr]\ar[ur]&&O_N^+\ar[ur]\\
&K_N\ar[rr]\ar[uu]&&U_N\ar[uu]\\
H_N\ar[uu]\ar[ur]\ar[rr]&&O_N\ar[uu]\ar[ur]
}$$
the corresponding truncated characters follow with $N\to\infty$ the laws
$$
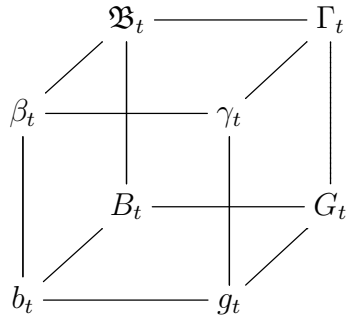
\xymatrix@R=18pt@C=20pt{
&\mathfrak B_t\ar@{-}[rr]\ar@{-}[dd]&&\Gamma_t\ar@{-}[dd]\\
\beta_t\ar@{-}[rr]\ar@{-}[dd]\ar@{-}[ur]&&\gamma_t\ar@{-}[dd]\ar@{-}[ur]\\
&B_t\ar@{-}[rr]\ar@{-}[uu]&&G_t\ar@{.}[uu]\\
b_t\ar@{-}[uu]\ar@{-}[ur]\ar@{-}[rr]&&g_t\ar@{-}[uu]\ar@{-}[ur]
}$$
which are the main limiting laws in classical and free probability.
\end{theorem}

\begin{proof}
We know from earlier in this chapter that the above quantum groups are all easy, coming from the following categories of partitions:
$$
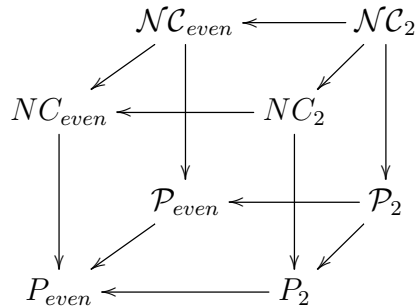
\xymatrix@R=17pt@C4pt{
&\mathcal{NC}_{even}\ar[dl]\ar[dd]&&\mathcal {NC}_2\ar[dl]\ar[ll]\ar[dd]\\
NC_{even}\ar[dd]&&NC_2\ar[dd]\ar[ll]\\
&\mathcal P_{even}\ar[dl]&&\mathcal P_2\ar[dl]\ar[ll]\\
P_{even}&&P_2\ar[ll]
}$$

(1) At $t=1$, we can use the following formula, coming from easiness:
$$\lim_{N\to\infty}\int_{G_N}\chi^k=|D(k)|$$

But this gives the laws in the statement, via some standard calculus. 

\medskip

(2) In order to compute now the asymptotic laws of truncated characters, at any $t>0$, we can use the general moment formula in Proposition 11.23, namely:
$$\int_G(u_{11}+\ldots +u_{ss})^k=Tr(W_{kN}G_{ks})$$

To be more precise, what happens is that in each of the cases under consideration, the Gram matrix is asymptotically diagonal, and so the Weingarten matrix is asymptotically diagonal too. Thus, in the limit we obtain the following moment formula:
$$\lim_{N\to\infty}\int_{G_N}\chi_t^k=\sum_{\pi\in D(k)}t^{|\pi|}$$

But this gives the laws in the statement, via some standard calculus.
\end{proof}

\section*{11e. Exercises}

We had a lot of combinatorics in this chapter, and as exercises, we have:

\begin{exercise}
Fill in the missing details for our various Brauer theorems.
\end{exercise}

\begin{exercise}
Find easiness results for other groups and quantum groups.
\end{exercise}

\begin{exercise}
Discuss the half-liberation operation, in the unitary setting.
\end{exercise}

\begin{exercise}
Learn about the Idel-Wolf theorem, involving the group $C_N$.
\end{exercise}

\begin{exercise}
Learn more about free probability, and random matrices too.
\end{exercise}

\begin{exercise}
Do all computations for the asymptotic laws of characters.
\end{exercise}

As bonus exercise, try developing a general theory of Lie-Brauer algebras.

\chapter{Advanced results}

\section*{12a. Deformation, twists}

We have already seen some non-trivial quantum group theory in this book, with this certainly making it for an introduction to the subject. The present chapter will be our last, regarding quantum groups, with Part IV below rather dealing with quantum manifolds. We have many things that we can talk about, here, and we have chosen 4, as follows:

\bigskip

(1) Deformation, twists. We will talk here about Drinfeld-Jimbo deformations, involving a parameter $q\in\mathbb C-\{0\}$, then about Woronowicz deformations, where $q>0$, and finally we will discuss a certain Schur-Weyl twisting procedure, where $q=-1$.

\bigskip

(2) Amenability, growth. We have seen some analysis, of probabilistic flavor, at the end of chapter 11, but that is not the end of the story. The bulk of analysis on quantum groups regards the amenability and growth of $\Gamma=\widehat{G}$, and we will discuss this here.

\bigskip

(3) Planar algebras. Besides the above-mentioned deformations, there are many versions and generalizations of the compact quantum groups discussed before, and notably the tensor categories, and the planar algebras. We will have a look into these.

\bigskip

(4) Matrix models. Finally, with planar algebras and the related subfactors being something more general, and a bit more concrete too, than our quantum groups, what to do. Well, we will strike back with some sharp results, regarding the matrix models.

\bigskip

So, this was for the general idea with this chapter, our final one on quantum groups, and in practice, many things to be discussed, and do not expect of course many proofs, or so. However, we will provide links to the relevant literature, on each subject.

\bigskip

Getting started now, in relation with deformations and twists, the first question that we would like to understand is something very simple and natural, as follows:

\begin{question}
Can a compact Lie group $G$ be deformed or not, within the framework of Hopf algebras, or so to say, within the framework of quantum groups?
\end{question}

As a first observation, this is certainly not possible within the framework of classical Lie groups, because these are subject to a discrete classification, ABCDEFG. In fact, this follows as well from an elementary result, stating that any Lie algebra is rigid, in the sense that when attempting to change a bit the structure constants, this is impossible.

\bigskip

In the Hopf algebra setting, however, everything changes. The discovery here, coming from the work of Faddeev and the Leningrad School of Physics from the late 70s \cite{fad}, then that of Drinfeld \cite{dri} and Jimbo \cite{jim} from the mid 80s, is as follows:

\index{Drinfeld-Jimbo deformation}
\index{Leningrad School}
\index{Faddeev}
\index{generic parameter}
\index{square of antipode}

\begin{theorem}
Given a classical Lie group $G$, with Lie algebra $\mathfrak g$, the enveloping Lie algebra $U\mathfrak g$, regarded as Hopf algebra, admits a non-trivial deformation $U_q\mathfrak g$, with $q\in F-\{0\}$. This latter deformation can be in fact recovered as being something canonical, and can be thought of as being the enveloping Lie algebra of a quantum group $G_q$.
\end{theorem}

\begin{proof}
This is something coming from the above-mentioned papers of Drinfeld \cite{dri} and Jimbo \cite{jim}, and which is best understood for $G=SU_2$, the idea being as follows:

\medskip

(1) Deforming $SU_2$ with the help of a generic parameter $q\in F$ is something not very complicated, that you can actually recover by yourself, by doing some computations, once of course someone told you in advance that this is possible. 

\medskip

(2) For more complicated Lie groups, however, some Lie algebra theory is needed. Also, switching from the generic parameter $q\in F$ to the concrete parameter $q\in F-\{0\}$ is something quite delicate too, again requiring some good Lie algebra knowledge. 

\medskip

(3) Finally, in what regards the interpretation of the output, navigating between Lie groups, Lie algebras and enveloping Lie algebras, at the quantum level, requires a good knowledge of both Hopf algebras and Lie groups and Lie algebra theory.

\medskip

(4) In short, a bit complicated all this, with lots of algebra involved, but again, the case of $SU_2$, which is the most important one, is quite easy to be fully understood, and for more, you have several books dedicated to this, such as Lusztig \cite{lus}.
\end{proof}

At the representation theory level now, which is the one that matters, as we previously know from the present book, things are quite tricky, as follows:

\index{roots of unity}
\index{not semisimple}
\index{semisimple}
\index{root of unity}

\begin{theorem}
When the deformation parameter $q\in F-\{0\}$ is not a root of unity, $G_q$ is semisimple, with representations similar to those of the original Lie group $G=G_1$. However, when $q\in F-\{0,1\}$ is a root of unity, $G_q$ is not semisimple.
\end{theorem}

\begin{proof}
Again, there is a long story here, the idea being that semisimplicity is easy to establish in the framework of the deformations with generic parameter $q\in F$, but then the study of the specializations, at the various values of $q\in F-\{0\}$, is something quite delicate. As before, we refer here to Lusztig \cite{lus}, and other books on the subject.
\end{proof}

Many other things can be said about the above quantum groups, notably with the following statement, clarifying the relation with the Woronowicz quantum groups:

\index{subfactor}
\index{planar algebra}
\index{deformation}

\begin{theorem}
The quantum groups $G_q$ with $q\in\mathbb C-\{0\}$ are as follows,
\begin{enumerate}
\item At $q=1$, these are the original Lie groups $G$.

\item At $q>0$, these are semisimple, but with $S^2\neq id$ at $q\neq1$.

\item Still at $q>0$, these are covered by an extension of the Woronowicz formalism,

\item With this being best seen via Tannakian duality, as explained by Rosso.

\item At $q\not>0$, no Woronowicz type algebra formalism is available,

\item However, subfactors and planar algebras do cover these $q\not>0$ beasts,

\item And in particular, they do cover the case where $q\neq 1$ is a root of unity.
\end{enumerate}
and with all this coming, as we will soon discover, with an uncertainty at $q=-1$.
\end{theorem}

\begin{proof}
Crowded statement that we have here, and for more on all this, we refer to the papers of Woronowicz \cite{wo1}, \cite{wo2} and of Rosso \cite{ro1}, \cite{ro2}, which are a must-read, if you are interested in $q$-deformations, and to Jones \cite{jo1}, \cite{jo3} and Kirillov Jr. \cite{kir} and Wenzl \cite{wen} for the last part, in relation with subfactors, which are again a must-read.
\end{proof}

With this discussed, which is more than borderline with respect to what a mathematics book is supposed to contain, Theorems coming with Proofs, let us focus now on a modest question, with the aim of discussing it with full details, as mathematicians do:

\begin{question}
What exactly happens at $q=-1$?
\end{question}

To be more precise, we have seen in chapters 7-8 that anticommutation twists, which normally should correspond to the parameter value $q=-1$, do exist, and are very interesting and natural objects, in the general affine Hopf algebra setting. And the point now is that, if you get back to that chapters 7-8 material, you will certainly discover that, in most of the cases, what we have there, over $F=\mathbb C$, are Woronowicz algebras.

\bigskip

Summarizing, we have here a bit of contradiction with what Drinfeld-Jimbo and the others are saying, coming from the fact that their $q=-1$ twists are not semisimple, due to the fact that $q=-1$ is a root of unity, while our $q=-1$ twists are semisimple.

\bigskip

Not an easy situation that we have here, with each camp, they and us, claiming that their theory is the correct one at $q=-1$. So, which theory is correct? God only knows, but if that can help, here is what cat declares, about all this:

\index{anticommutation}
\index{quantum mechanics}

\begin{cat}
Anticommutation comes from quantum mechanics, that you humans would like to understand, with your quantum groups. So go with it, and by the way, you should update your anticommutation twisting methods, now that you know about easiness. 
\end{cat}

Thanks cat, not only I am pleased to learn this, but I eventually have a plan, for the remainder of this section. So, let us review now the anticommutation twisting operation $G\to G'$, that we know about since chapters 7-8, by using Tannakian duality.

\bigskip

Given a partition $\tau\in P(k,l)$, let us call ``switch'' the operation which consists in switching two neighbors, belonging to different blocks, in the upper row, or in the lower row. With this convention, we have the following key result:

\index{signature}
\index{number of crossings}

\begin{proposition}
There is a signature map $\varepsilon:P_{even}\to\{-1,1\}$, given by 
$$\varepsilon(\tau)=(-1)^c$$
where $c$ is the number of switches needed to make $\tau$ noncrossing. In addition:
\begin{enumerate}
\item For $\tau\in S_k$, this is the usual signature.

\item For $\tau\in P_2$ we have $(-1)^c$, where $c$ is the number of crossings.

\item For $\tau\leq\pi\in NC_{even}$, the signature is $1$.
\end{enumerate}
\end{proposition}

\begin{proof}
The fact that $\varepsilon$ is indeed well-defined comes from the fact that the number $c$ in the statement is well-defined modulo 2, which is standard combinatorics. In order to prove now the remaining assertion, observe that any partition $\tau\in P(k,l)$ can be put in ``standard form'', by ordering its blocks according to the appearence of the first leg in each block, counting clockwise from top left, and then by performing the switches as for block 1 to be at left, then for block 2 to be at left, and so on:
$$
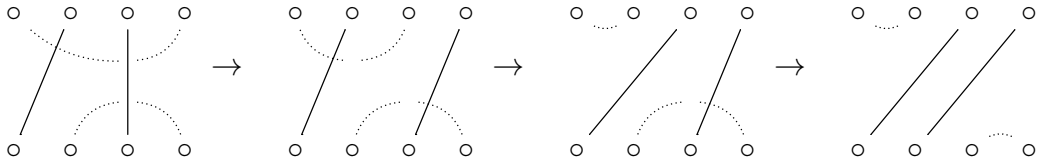
\xymatrix@R=3mm@C=3mm{\circ\ar@/_/@{.}[drr]&\circ\ar@{-}[dddl]&\circ\ar@{-}[ddd]&\circ\\
&&\ar@/_/@{.}[ur]&\\
&&\ar@/^/@{.}[dr]&\\
\circ&\circ\ar@/^/@{.}[ur]&\circ&\circ}
\xymatrix@R=4mm@C=1mm{&\\\to\\&\\& }
\xymatrix@R=3mm@C=3mm{\circ\ar@/_/@{.}[dr]&\circ\ar@{-}[dddl]&\circ&\circ\ar@{-}[dddl]\\
&\ar@/_/@{.}[ur]&&\\
&&\ar@/^/@{.}[dr]&\\
\circ&\circ\ar@/^/@{.}[ur]&\circ&\circ}
\xymatrix@R=4mm@C=1mm{&\\\to\\&\\&}
\xymatrix@R=3mm@C=3mm{\circ\ar@/_/@{.}[r]&\circ&\circ\ar@{-}[dddll]&\circ\ar@{-}[dddl]\\
&&&\\
&&\ar@/^/@{.}[dr]&\\
\circ&\circ\ar@/^/@{.}[ur]&\circ&\circ}
\xymatrix@R=4mm@C=1mm{&\\\to\\&\\& }
\xymatrix@R=3mm@C=3mm{\circ\ar@/_/@{.}[r]&\circ&\circ\ar@{-}[dddll]&\circ\ar@{-}[dddll]\\
&&&\\
&&&\\
\circ&\circ&\circ\ar@/^/@{.}[r]&\circ}$$

With this convention, the proof of the remaining assertions is as follows:

\medskip

(1) For $\tau\in S_k$ the standard form is $\tau'=id$, and the passage $\tau\to id$ comes by composing with a number of transpositions, which gives the signature. 

\medskip

(2) For a general $\tau\in P_2$, the standard form is of type $\tau'=|\ldots|^{\cup\ldots\cup}_{\cap\ldots\cap}$, and the passage $\tau\to\tau'$ requires $c$ mod 2 switches, where $c$ is the number of crossings. 

\medskip

(3) Assuming that $\tau\in P_{even}$ comes from $\pi\in NC_{even}$ by merging a certain number of blocks, we can prove that the signature is 1 by proceeding by recurrence.
\end{proof}

With the above result in hand, we can now formulate:

\index{twisted linear map}
\index{twisted Kronecker symbol}

\begin{definition}
Associated to any partition $\pi\in P_{even}(k,l)$ is the linear map
$$T_\pi':(\mathbb C^N)^{\otimes k}\to(\mathbb C^N)^{\otimes l}$$
given by the following formula, with $e_1,\ldots,e_N$ being the standard basis of $\mathbb C^N$,
$$T_\pi'(e_{i_1}\otimes\ldots\otimes e_{i_k})=\sum_{j_1\ldots j_l}\delta_\pi'\begin{pmatrix}i_1&\ldots&i_k\\ j_1&\ldots&j_l\end{pmatrix}e_{j_1}\otimes\ldots\otimes e_{j_l}$$
and where $\delta_\pi'\in\{-1,0,1\}$ is $\delta_\pi'=\varepsilon(\tau)$ if $\tau\geq\pi$, and $\delta_\pi'=0$ otherwise, with $\tau=\ker(^i_j)$.
\end{definition}

In other words, what we are doing here is to add signatures to the usual formula of $T_\pi$. Indeed, observe that the usual formula for $T_\pi$ can be written as folllows:
$$T_\pi(e_{i_1}\otimes\ldots\otimes e_{i_k})=\sum_{j:\ker(^i_j)\geq\pi}e_{j_1}\otimes\ldots\otimes e_{j_l}$$

Now by inserting signs, coming from the signature map $\varepsilon:P_{even}\to\{\pm1\}$, we are led to the following formula, which coincides with the one given above:
$$T_\pi'(e_{i_1}\otimes\ldots\otimes e_{i_k})=\sum_{\tau\geq\pi}\varepsilon(\tau)\sum_{j:\ker(^i_j)=\tau}e_{j_1}\otimes\ldots\otimes e_{j_l}$$

We will be back later to this analogy, with more details on what can be done with it. For the moment, we must first prove a key categorical result, as follows:

\begin{proposition}
The assignement $\pi\to T_\pi'$ is categorical, in the sense that
$$T_\pi'\otimes T_\sigma'=T_{[\pi\sigma]}'\quad,\quad 
T_\pi'T_\sigma'=N^{c(\pi,\sigma)}T_{[^\sigma_\pi]}'\quad,\quad
(T_\pi')^*=T_{\pi^*}'$$
where $c(\pi,\sigma)$ are certain positive integers.
\end{proposition}

\begin{proof}
We have to go back to the proof from the untwisted case, from chapter 11, and insert signs. We have to check three conditions, as follows:

\medskip

\underline{1. Concatenation}. In the untwisted case, this was based on the following formula:
$$\delta_\pi\begin{pmatrix}i_1\ldots i_p\\ j_1\ldots j_q\end{pmatrix}
\delta_\sigma\begin{pmatrix}k_1\ldots k_r\\ l_1\ldots l_s\end{pmatrix}
=\delta_{[\pi\sigma]}\begin{pmatrix}i_1\ldots i_p&k_1\ldots k_r\\ j_1\ldots j_q&l_1\ldots l_s\end{pmatrix}$$

In the twisted case, it is enough to check the following formula:
$$\varepsilon\left(\ker\begin{pmatrix}i_1\ldots i_p\\ j_1\ldots j_q\end{pmatrix}\right)
\varepsilon\left(\ker\begin{pmatrix}k_1\ldots k_r\\ l_1\ldots l_s\end{pmatrix}\right)=
\varepsilon\left(\ker\begin{pmatrix}i_1\ldots i_p&k_1\ldots k_r\\ j_1\ldots j_q&l_1\ldots l_s\end{pmatrix}\right)$$

Let us denote by $\tau,\nu$ the partitions on the left, so that the partition on the right is of the form $\rho\leq[\tau\nu]$. Now by switching to the noncrossing form, $\tau\to\tau'$ and $\nu\to\nu'$, the partition on the right transforms into $\rho\to\rho'\leq[\tau'\nu']$. Now since the partition $[\tau'\nu']$ is noncrossing, we can use Proposition 12.7 (3), and we obtain the result.

\medskip

\underline{2. Composition}. In the untwisted case, this was based on the following formula:
$$\sum_{j_1\ldots j_q}\delta_\pi\begin{pmatrix}i_1\ldots i_p\\ j_1\ldots j_q\end{pmatrix}
\delta_\sigma\begin{pmatrix}j_1\ldots j_q\\ k_1\ldots k_r\end{pmatrix}
=N^{c(\pi,\sigma)}\delta_{[^\pi_\sigma]}\begin{pmatrix}i_1\ldots i_p\\ k_1\ldots k_r\end{pmatrix}$$

In order to prove now the result in the twisted case, it is enough to check that the signs match. More precisely, we must establish the following formula:
$$\varepsilon\left(\ker\begin{pmatrix}i_1\ldots i_p\\ j_1\ldots j_q\end{pmatrix}\right)
\varepsilon\left(\ker\begin{pmatrix}j_1\ldots j_q\\ k_1\ldots k_r\end{pmatrix}\right)
=\varepsilon\left(\ker\begin{pmatrix}i_1\ldots i_p\\ k_1\ldots k_r\end{pmatrix}\right)$$

Let $\tau,\nu$ be the partitions on the left, so that the partition on the right is of the form $\rho\leq[^\tau_\nu]$. Our claim is that we can jointly switch $\tau,\nu$ to the noncrossing form. Indeed, we can first switch as for $\ker(j_1\ldots j_q)$ to become noncrossing, and then switch the upper legs of $\tau$, and the lower legs of $\nu$, as for both these partitions to become noncrossing. Now observe that when switching in this way to the noncrossing form, $\tau\to\tau'$ and $\nu\to\nu'$, the partition on the right transforms into $\rho\to\rho'\leq[^{\tau'}_{\nu'}]$. Now since the partition $[^{\tau'}_{\nu'}]$ is noncrossing, we can apply Proposition 12.7 (3), and we obtain the result.

\medskip

\underline{3. Involution}. Here we must prove the following formula:
$$\delta_\pi'\begin{pmatrix}i_1\ldots i_p\\ j_1\ldots j_q\end{pmatrix}=\delta_{\pi^*}'\begin{pmatrix}j_1\ldots j_q\\ i_1\ldots i_p\end{pmatrix}$$

But this is clear from the definition of $\delta_\pi'$, and we are done.
\end{proof}

As a conclusion, our twisted construction $\pi\to T_\pi'$ has all the needed properties for producing quantum groups, via Tannakian duality, and we can now formulate:

\begin{theorem}
Given a category of partitions $D\subset P_{even}$, the construction
$$Hom(u^{\otimes k},u^{\otimes l})=span\left(T_\pi'\Big|\pi\in D(k,l)\right)$$
produces via Tannakian duality a quantum group $G_N'\subset U_N^+$, for any $N\in\mathbb N$.
\end{theorem}

\begin{proof}
This follows indeed from the Tannakian results from chapter 10, exactly as in the easy case, by using this time Proposition 12.9 as technical ingredient.
\end{proof}

We can unify the easy quantum groups, or at least the examples coming from categories $D\subset P_{even}$, with the quantum groups constructed above, as follows:

\begin{definition}
A closed subgroup $G\subset U_N^+$ is called $q$-easy, or quizzy, with deformation parameter $q=\pm1$, when its tensor category appears as follows,
$$Hom(u^{\otimes k},u^{\otimes l})=span\left(\dot{T}_\pi\Big|\pi\in D(k,l)\right)$$
for a certain category of partitions $D\subset P_{even}$, where, for $q=-1,1$:
$$\dot{T}=T',T$$
The Schur-Weyl twist of $G$ is the quizzy quantum group $G'\subset U_N^+$ obtained via $q\to-q$.
\end{definition}

Let us clarify now the relation between the maps $T_\pi,T_\pi'$. We recall that the M\"obius function of any lattice, and in particular of $P_{even}$, is given by:
$$\mu(\sigma,\pi)=\begin{cases}
1&{\rm if}\ \sigma=\pi\\
-\sum_{\sigma\leq\tau<\pi}\mu(\sigma,\tau)&{\rm if}\ \sigma<\pi\\
0&{\rm if}\ \sigma\not\leq\pi
\end{cases}$$

With this notation, we have the following useful result:

\index{M\"obius formula}

\begin{theorem}
For any partition $\pi\in P_{even}$ we have the formula
$$T_\pi'=\sum_{\tau\leq\pi}\alpha_\tau T_\tau$$
where $\alpha_\sigma=\sum_{\sigma\leq\tau\leq\pi}\varepsilon(\tau)\mu(\sigma,\tau)$, with $\mu$ being the M\"obius function of $P_{even}$.
\end{theorem}

\begin{proof}
The linear combinations $T=\sum_{\tau\leq\pi}\alpha_\tau T_\tau$ acts on tensors as follows:
\begin{eqnarray*}
T(e_{i_1}\otimes\ldots\otimes e_{i_k})
&=&\sum_{\tau\leq\pi}\alpha_\tau T_\tau(e_{i_1}\otimes\ldots\otimes e_{i_k})\\
&=&\sum_{\tau\leq\pi}\alpha_\tau\sum_{\sigma\leq\tau}\sum_{j:\ker(^i_j)=\sigma}e_{j_1}\otimes\ldots\otimes e_{j_l}\\
&=&\sum_{\sigma\leq\pi}\left(\sum_{\sigma\leq\tau\leq\pi}\alpha_\tau\right)\sum_{j:\ker(^i_j)=\sigma}e_{j_1}\otimes\ldots\otimes e_{j_l}
\end{eqnarray*}

Thus, in order to have $T_\pi'=\sum_{\tau\leq\pi}\alpha_\tau T_\tau$, we must have $\varepsilon(\sigma)=\sum_{\sigma\leq\tau\leq\pi}\alpha_\tau$, for any $\sigma\leq\pi$. But this problem can be solved by using the M\"obius inversion formula, and we obtain the numbers $\alpha_\sigma=\sum_{\sigma\leq\tau\leq\pi}\varepsilon(\tau)\mu(\sigma,\tau)$ in the statement.
\end{proof}

We can now twist $O_N,U_N$, and $O_N^*,U_N^*$ too, the result here being as follows:

\index{twisted orthogonal group}
\index{twisted unitary group}

\begin{theorem}
The twists of $O_N,U_N$ are obtained by replacing the commutation relations $ab=ba$ between the coordinates $u_{ij}$ and their adjoints $u_{ij}^*$ with the relations
$$ab=\pm ba$$
with anticommutation on rows and columns, and commutation otherwise, so we obtain the previous twists $O_N',U_N'$. Similarly, for $O_N^*,U_N^*$ we obtain $O_N^{*\prime},U_N^{*\prime}$.
\end{theorem}

\begin{proof}
The basic crossing, $\ker\binom{ij}{ji}$ with $i\neq j$, comes from the transposition $\tau\in S_2$, so its signature is $-1$. As for its degenerated version $\ker\binom{ii}{ii}$, this is noncrossing, so here the signature is $1$. We conclude that the linear map associated to the basic crossing is:
$$\bar{T}_{\slash\!\!\!\backslash}(e_i\otimes e_j)
=\begin{cases}
-e_j\otimes e_i&{\rm for}\ i\neq j\\
e_j\otimes e_i&{\rm otherwise}
\end{cases}$$

We can proceed now as in the untwisted case, and since the intertwining relations coming from $\bar{T}_{\slash\!\!\!\backslash}$ correspond to the relations defining $O_N',U_N'$, we obtain the result. The study is similar for $O_N^*,U_N^*$, where we obtain the previous twists $O_N^{*\prime},U_N^{*\prime}$.
\end{proof}

There are many more things that can be said about Schur-Weyl twists, notably with the fact that the quantum reflection groups equal their own twists. Also, at the analytic level, the Gram and Weingarten matrices stay the same, under twisting, with the signs appearing as coefficients. For more on all this, we refer to the twisting literature.

\section*{12b. Amenability, growth}

Switching topics, as a second theme of discussion for the present chapter, let us discuss analysis on discrete quantum groups, notably with the notions of amenability, and growth. We have already met amenability in chapter 10, with our discussion there being quite brief. Here is the main result, again, coming this time with a complete proof, or almost:

\index{amenability}
\index{growth}
\index{semisimplicity}
\index{cosemisimplicity}

\begin{theorem}
Let $A_{full}$ be the enveloping $C^*$-algebra of $\mathcal A$, and let $A_{red}$ be the quotient of $A$ by the null ideal of the Haar integration. The following are equivalent:
\begin{enumerate}
\item The Haar functional of $A_{full}$ is faithful.

\item The projection map $A_{full}\to A_{red}$ is an isomorphism.

\item The counit map $\varepsilon:A\to\mathbb C$ factorizes through $A_{red}$.

\item We have $N\in\sigma(Re(\chi_u))$, the spectrum being taken inside $A_{red}$.

\item $||ax_k-\varepsilon(a)x_k||\to0$ for any $a\in\mathcal A$, for certain norm $1$ vectors $x_k\in L^2(A)$.
\end{enumerate}
If this is the case, we say that the underlying discrete quantum group $\Gamma$ is amenable.
\end{theorem}

\begin{proof}
We will work out explicitly all the possible implications between (1-5), whenever possible, adding to the global formal proof, which will be linear, as follows:
$$(1)\implies(2)\implies(3)\implies(4)\implies(5)\implies(1)$$

In order to prove these implications, and the other ones too, the general idea is that this is is well-known in the group dual case, $A=C^*(\Gamma)$, with $\Gamma$ being a usual discrete group, and in general, the result follows by adapting the group dual case proof. 

\medskip

$(1)\iff(2)$ This follows from the fact that the GNS construction for the algebra $A_{full}$ with respect to the Haar functional produces the algebra $A_{red}$.

\medskip

$(2)\implies(3)$ This is trivial, because we have quotient maps $A_{full}\to A\to A_{red}$, and so our assumption $A_{full}=A_{red}$ implies that we have $A=A_{red}$. 

\medskip 

$(3)\implies(2)$ Assume indeed that we have a counit map $\varepsilon:A_{red}\to\mathbb C$. In order to prove $A_{full}=A_{red}$, we can use the right regular corepresentation. Indeed, we can define such a corepresentation by the following formula:
$$W(a\otimes x)=\Delta(a)(1\otimes x)$$

This corepresentation is unitary, so we can define a morphism as follows: 
$$\Delta':A_{red}\to A_{red}\otimes A_{full}\quad,\quad 
a\to W(a\otimes1)W^*$$

Now by composing with $\varepsilon\otimes id$, we obtain a morphism as follows:
$$(\varepsilon\otimes id)\Delta':A_{red}\to A_{full}\quad,\quad 
u_{ij}\to u_{ij}$$

Thus, we have our inverse for the canonical projection $A_{full}\to A_{red}$, as desired.

\medskip

$(3)\implies(4)$ This implication is clear, because we have:
\begin{eqnarray*}
\varepsilon(Re(\chi_u))
&=&\frac{1}{2}\left(\sum_{i=1}^N\varepsilon(u_{ii})+\sum_{i=1}^N\varepsilon(u_{ii}^*)\right)\\
&=&\frac{1}{2}(N+N)\\
&=&N
\end{eqnarray*}

Thus the element $N-Re(\chi_u)$ is not invertible in $A_{red}$, as claimed.

\medskip

$(4)\implies(3)$ In terms of the corepresentation $v=u+\bar{u}$, whose dimension is $2N$ and whose character is $2Re(\chi_u)$, our assumption $N\in\sigma(Re(\chi_u))$ reads:
$$\dim v\in\sigma(\chi_v)$$

By functional calculus the same must hold for $w=v+1$, and then once again by functional calculus, the same must hold for any tensor power of $w$:
$$w_k=w^{\otimes k}$$ 

Now choose for each $k\in\mathbb N$ a state $\varepsilon_k\in A_{red}^*$ having the following property:
$$\varepsilon_k(w_k)=\dim w_k$$

By Peter-Weyl we must have $\varepsilon_k(r)=\dim r$ for any $r\leq w_k$, and since any irreducible corepresentation appears in this way, the sequence $\varepsilon_k$ converges to a counit map: 
$$\varepsilon:A_{red}\to\mathbb C$$

$(4)\implies(5)$ Consider the following elements of $A_{red}$, which are positive:
$$a_i=1-Re(u_{ii})$$
 
Our assumption $N\in\sigma(Re(\chi_u))$ tells us that $a=\sum a_i$ is not invertible, and so there exists a sequence $x_k$ of norm one vectors in $L^2(A)$ such that: 
$$<ax_k,x_k>\to 0$$

Since the summands $<a_ix_k,x_k>$ are all positive, we must have, for any $i$:
$$<a_ix_k,x_k>\to0$$

We can go back to the variables $u_{ii}$ by using the following general formula:
$$||vx-x||^2=||vx||^2 +2<(1-Re(v))x,x>-1$$

Indeed, with $v=u_{ii}$ and $x=x_k$ the middle term on the right goes to 0, and so the whole term on the right becomes asymptotically negative, and so we must have:
$$||u_{ii}x_k-x_k||\to0$$

Now let $M_n(A_{red})$ act on $\mathbb C^n\otimes L^2(A)$. Since $u$ is unitary we have:
$$\sum_i||u_{ij}x_k||^2
=||u(e_j\otimes x_k)||
=1$$

From $||u_{ii}x_k||\to1$ we obtain $||u_{ij}x_k||\to0$ for $i\neq j$. Thus we have, for any $i,j$:
$$||u_{ij}x_k-\delta_{ij}x_k||\to0$$

Now by remembering that we have $\varepsilon(u_{ij})=\delta_{ij}$, this formula reads:
$$||u_{ij}x_k-\varepsilon(u_{ij})x_k||\to0$$

By linearity, multiplicativity and continuity, we must have, for any $a\in\mathcal  A$, as desired:
$$||ax_k-\varepsilon(a)x_k||\to0$$

$(5)\implies(1)$ This is well-known, coming from standard functional analysis.

\medskip

$(1)\implies(5)$ Once again this is something well-known, and standard. See \cite{wo1}.
\end{proof}

Let us focus now on the Kesten amenability criterion, the one from Theorem 12.14 (4), which will be our main amenability criterion, in what follows. We will need:

\index{main character}

\begin{proposition}
Given a Woronowicz algebra $(A,u)$, with $u\in M_N(A)$, the moments of the main character $\chi=\sum_iu_{ii}$ are given by:
$$\int_G\chi^k=\dim\left(Fix(u^{\otimes k})\right)$$
In the case $u\sim\bar{u}$ the law of $\chi$ is a usual probability measure, supported on $[-N,N]$.
\end{proposition}

\begin{proof}
There are two assertions here, the proof being as follows:

\medskip

(1) The first assertion follows from the Peter-Weyl theory, which tells us that we have the following formula, valid for any corepresentation $v\in M_n(A)$:
$$\int_G\chi_v=\dim(Fix(v))$$

Indeed, for $v=u^{\otimes k}$ we have $\chi_v=\chi^k$, so we obtain the formula in the statement.

\medskip

(2) As for the second assertion, if we assume $u\sim\bar{u}$ then we have $\chi=\chi^*$, and so the general theory, explained above, tells us that $law(\chi)$ is in this case a real probability measure, supported by the spectrum of $\chi$. But, since $u\in M_N(A)$ is unitary, we have:
$$uu^*=1
\implies||u_{ij}||\leq 1,\forall i,j
\implies||\chi||\leq N$$

Thus the spectrum of the character satisfies $\sigma(\chi)\subset [-N,N]$, as claimed. 
\end{proof}

In relation now with the notion of amenability, we have the following result:

\index{Kesten amenability}
\index{amenability}

\begin{theorem}
A Woronowicz algebra $(A,u)$, with $u\in M_N(A)$, is amenable when:
$$N\in supp\Big(law\left[Re(\chi)\right]\Big)$$
Also, the support on the right depends only on $law(\chi)$.
\end{theorem}

\begin{proof}
There are two assertions here, the proof being as follows:

\medskip

(1) According to the Kesten amenability criterion, from Theorem 12.14 (4), the algebra $A$ is amenable precisely when the following condition is satisfied:
$$N\in\sigma(Re(\chi))$$

Now since $Re(\chi)$ is self-adjoint, we know from spectral theory that the support of its spectral measure $law(Re(\chi))$ is precisely its spectrum $\sigma(Re(\chi))$, as desired:
$$supp(law(Re(\chi)))=\sigma(Re(\chi))$$

(2) Regarding the second assertion, once again the variable $Re(\chi)$ being self-adjoint, its law depends only on the moments $\int_GRe(\chi)^p$, with $p\in\mathbb N$. But, we have:
$$\int_GRe(\chi)^p
=\int_G\left(\frac{\chi+\chi^*}{2}\right)^p
=\frac{1}{2^p}\sum_{|k|=p}\int_G\chi^k$$

Thus $law(Re(\chi))$ depends only on $law(\chi)$, and this gives the result.
\end{proof}

As a first illustration, let us work out now in detail the group dual case. Here we obtain a very interesting measure, called Kesten measure of the group, as follows:

\begin{theorem}
In the case $A=C^*(\Gamma)$ and $u=diag(g_1,\ldots,g_N)$, and by enlarging if necessary the generating set $g_1,\ldots,g_N$, as for the following to happen,
$$1\in u=\bar{u}$$
the moments of the main character are given by the formula
$$\int_{\widehat{\Gamma}}\chi^p=\#\left\{i_1,\ldots,i_p\Big|g_{i_1}\ldots g_{i_p}=1\right\}$$
counting the loops based at $1$, having lenght $p$, on the corresponding Cayley graph.
\end{theorem}

\begin{proof}
Consider indeed a discrete group $\Gamma=<g_1,\ldots,g_N>$. The main character of $A=C^*(\Gamma)$, with fundamental corepresentation $u=diag(g_1,\ldots,g_N)$, is then:
$$\chi=g_1+\ldots+g_N$$

Given a colored integer $k=e_1\ldots e_p$, the corresponding moment is given by:
$$\int_{\widehat{\Gamma}}\chi^k
=\int_{\widehat{\Gamma}}(g_1+\ldots+g_N)^k
=\#\left\{i_1,\ldots,i_p\Big|g_{i_1}^{e_1}\ldots g_{i_p}^{e_p}=1\right\}$$

In the self-adjoint case, $u\sim\bar{u}$, we are only interested in the moments with respect to usual integers, $p\in\mathbb N$, and the above formula becomes:
$$\int_{\widehat{\Gamma}}\chi^p=\#\left\{i_1,\ldots,i_p\Big|g_{i_1}\ldots g_{i_p}=1\right\}$$

Assume now that we have in addition $1\in u$, so that the condition $1\in u=\bar{u}$ in the statement is satisfied. At the level of the generating set $S=\{g_1,\ldots,g_N\}$ this means:
$$1\in S=S^{-1}$$

Thus the corresponding Cayley graph is well-defined, with the elements of $\Gamma$ as vertices, and with the edges $g-h$ appearing when the condition $gh^{-1}\in S$ is satisfied. A loop on this graph based at 1, having lenght $p$, is then a sequence as follows:
$$(1)-(g_{i_1})-(g_{i_1}g_{i_2})-\ldots-(g_{i_1}\ldots g_{i_{p-1}})-(g_{i_1}\ldots g_{i_p}=1)$$

Thus the moments of $\chi$ count indeed such loops, as claimed.
\end{proof}

In order to generalize the above result, we will need the following standard fact:

\begin{theorem}
Let $(A,u)$ be a Woronowicz algebra, and assume, by enlarging if necessary $u$, that we have $1\in u=\bar{u}$. The following formula
$$d(v,w)=\min\left\{k\in\mathbb N\Big|1\subset\bar{v}\otimes w\otimes u^{\otimes k}\right\}$$
defines then a distance on $Irr(A)$, which coincides with the geodesic distance on the associated Cayley graph. In the group dual case we obtain the usual distance.
\end{theorem}

\begin{proof}
There are several assertions here, the idea being as follows:

\medskip

(1) The fact that the lengths are finite follows from Peter-Weyl theory, and the other verifications are standard too, using Frobenius duality, as a main ingredient. 

\medskip

(2) In the group dual case now, where our algebra is of the form $A=C^*(\Gamma)$, with $\Gamma=<S>$ being a finitely generated discrete group, our normalization condition $1\in u=\bar{u}$ from the statement means that the generating set $S\subset\Gamma$ must satisfy:
$$1\in S=S^{-1}$$

But this is precisely the normalization condition made before for the discrete groups, and the fact that we obtain the same metric space is clear.
\end{proof}

We can now formulate a generalization of Theorem 12.17, as follows:

\begin{theorem}
Let $(A,u)$ be a Woronowicz algebra, with the normalization assumption $1\in u=\bar{u}$ made. The moments of the main character, 
$$\int_G\chi^p=\dim\left(Fix(u^{\otimes p})\right)$$
count then the loops based at $1$, having lenght $p$, on the corresponding Cayley graph.
\end{theorem}

\begin{proof}
Here the formula of the moments, with $p\in\mathbb N$, is the one coming from Peter-Weyl, and the Cayley graph interpretation comes from Theorem 12.18.
\end{proof}

As an application of all this, corepresentation theory used for ``discrete'' questions, we can introduce the notion of growth for discrete quantum groups, as follows:

\begin{definition}
Given a closed subgroup $G\subset U_N^+$, with $1\in u=\bar{u}$, consider the series  whose coefficients are the ball volumes on the corresponding Cayley graph,
$$f(z)=\sum_kb_kz^k\quad,\quad 
b_k=\sum_{l(v)\leq k}\dim(v)^2$$
and call it growth series of the discrete quantum group $\widehat{G}$. In the group dual case, $G=\widehat{\Gamma}$, we obtain in this way the usual growth series of $\Gamma$. 
\end{definition}

As a first result about this, in relation with the notion of amenability, we have:

\begin{theorem}
Polynomial growth implies amenability.
\end{theorem}

\begin{proof}
We recall from Theorem 12.18 that the Cayley graph of $\widehat{G}$ has by definition the elements of $Irr(G)$ as vertices, and the distance is as follows:
$$d(v,w)=\min\left\{k\in\mathbb N\Big|1\subset\bar{v}\otimes w\otimes u^{\otimes k}\right\}$$

By taking $w=1$ and by using Frobenius reciprocity, the lenghts are given by:
$$l(v)=\min\left\{k\in\mathbb N\Big|v\subset u^{\otimes k}\right\}$$

By Peter-Weyl we have a decomposition as follows, where $B_k$ is the ball of radius $k$, and $m_k(v)\in\mathbb N$ are certain multiplicities:
$$u^{\otimes k}=\sum_{v\in B_k}m_k(v)\cdot v$$

By using now Cauchy-Schwarz, we obtain the following inequality:
\begin{eqnarray*}
m_{2k}(1)b_k
&=&\sum_{v\in B_k}m_k(v)^2\sum_{v\in B_k}\dim(v)^2\\
&\geq&\left(\sum_{v\in B_k}m_k(v)\dim(v)\right)^2\\
&=&N^{2k}
\end{eqnarray*}

But shows that if $b_k$ has polynomial growth, then the following happens:
$$\limsup_{k\to\infty}\, m_{2k}(1)^{1/2k}\geq N$$

Thus, the Kesten type criterion applies, and gives the result.
\end{proof}

For more on these topics, we refer to \cite{cdp}, \cite{dpr} and related papers.

\section*{12c. Planar algebras}

Switching again topics, let us discuss now the Jones theory of subfactors and planar algebras \cite{jo1}, \cite{jo2}, \cite{jo3}, \cite{jo4}, which can be regarded as being a far-reaching generalization of the usual quantum group theory, as developed in the previous chapters.

\bigskip

For this purpose, we will need some basic operator algebra theory. The foundational result here is the bicommutant theorem of von Neumann, which is as follows:

\index{operator algebra}
\index{Hilbert space}
\index{bounded operator}
\index{adjoint operator}
\index{von Neumann algebra}
\index{weak topology}
\index{algebra of operators}
\index{bicommutant theorem}
\index{Spectral theorem}
\index{normal operator}

\begin{theorem}
For a $*$-algebra of operators $A\subset B(H)$ the following conditions are equivalent, and if satisfied, we say that $A$ is a von Neumann algebra:
\begin{enumerate}
\item $A$ is closed under the weak topology, making each $T\to Tx$ continuous.

\item $A$ equals its bicommutant, $A=A''$, computed inside $B(H)$.
\end{enumerate}
\end{theorem}

\begin{proof}
The above equivalence is indeed something very standard, by using an amplification trick, and for details here, we refer to any operator algebra book. 
\end{proof}

Moving ahead, the continuation of the story involves an accumulation of non-trivial results, due to Murray and von Neumann, from the 1930s and 1940s, and then due to Connes, much later, in the 1970s, the main conclusions being as follows:

\begin{theorem}
The von Neumann algebras are as follows:
\begin{enumerate}
\item In the commutative case, these are the algebras $A=L^\infty(X)$, with $X$ measured space, represented on $H=L^2(X)$, up to a multiplicity.

\item If we write the center as $Z(A)=L^\infty(X)$, then we have a decomposition of type $A=\int_XA_x\,dx$, with the fibers $A_x$ having trivial center,  $Z(A_x)=\mathbb C$.

\item The factors, $Z(A)=\mathbb C$, can be fully classified in terms of ${\rm II}_1$ factors, which are those satisfying $\dim A=\infty$, and having a faithful trace $tr:A\to\mathbb C$.
\end{enumerate}
\end{theorem}

\begin{proof}
This is something quite heavy, the idea being as follows:

\medskip

(1) To start with, it is clear that $L^\infty(X)$ is indeed a von Neumann algebra on $H=L^2(X)$. The converse can be proved as well, by using spectral theory, one way of viewing this being by saying that, given a commutative von Neumann algebra $A\subset B(H)$, its elements $T\in A$ are commuting normal operators, so the Spectral Theorem for such operators applies, and gives $A=L^\infty(X)$, for some measured space $X$.

\medskip

(2) This is von Neumann's reduction theory main result, whose statement is already quite hard to understand, and whose proof uses advanced functional analysis. To be more precise, in finite dimensions this is something that we know well, with the formula $A=\int_XA_x\,dx$ corresponding to our usual direct sum decomposition, namely:
$$A=M_{n_1}(\mathbb C)\oplus\ldots\oplus M_{n_k}(\mathbb C)$$

In infinite dimensions, things are more complicated, but the idea remains the same, namely using (1) for the commutative von Neumann algebra $Z(A)$, as to get a measured space $X$, and then making your way towards a decomposition of type $A=\int_XA_x\,dx$.

\medskip

(3) This is something fairly heavy, due to Murray-von Neumann and Connes, the idea being that the other factors can be basically obtained via crossed product constructions. To be more precise, the various type of factors can be classified as follows:

\medskip

-- Type I. These are the matrix algebras $M_N(\mathbb C)$, called of type ${\rm I}_N$, and their infinite generalization, $B(H)$ with $H$ infinite dimensional, called of type ${\rm I}_\infty$. Although these factors are very interesting and difficult mathematical objects, from the perspective of the general von Neumann algebra classification work, they are dismissed as ``trivial''.

\medskip

-- Type II. These are the infinite dimensional factors having a trace, which is a usual trace $tr:A\to\mathbb C$ in the type ${\rm II}_1$ case, and is something more technical, possibly infinite, in the remaining case, the type ${\rm II}_\infty$ one, with these latter factors being of the form $B(H)\otimes A$, with $A$ being a ${\rm II}_1$ factor, and with $H$ being an infinite dimensional Hilbert space.

\medskip

-- Type III. These are the factors which are infinite dimensional, and do not have a trace $tr:A\to\mathbb C$. Murray and von Neumann struggled a lot with such beasts, with even giving an example being a non-trivial task, but later Connes came and classified them, basically showing that they appear from ${\rm II}_1$ factors, via crossed product constructions.
\end{proof}

So long for basic, or rather advanced but foundational, von Neumann algebra theory. In what follows we will focus on the ${\rm II}_1$ factors, according to the following principle:

\begin{principle}
The building blocks of the von Neumann algebra theory are the ${\rm II}_1$ factors, which are the von Neumann algebras having the following properties:
\begin{enumerate}
\item They are infinite dimensional, $\dim A=\infty$.

\item They are factors, their center being $Z(A)=\mathbb C$.

\item They have a faithful trace $tr:A\to\mathbb C$.
\end{enumerate}
\end{principle}

But you might perhaps ask, is it even clear that such beasts exist? Good point, and in answer, given a discrete group $\Gamma$, you can talk about its von Neumann algebra, obtained by talking the weak closure of the usual group algebra, or group $C^*$-algebra:
$$L(\Gamma)\subset B(l^2(\Gamma))$$

This algebra is then infinite dimensional when $\Gamma$ is infinite, and also has a trace, given on group elements by $tr(g)=\delta_{g1}$. As for the center, this consists of the functions on $\Gamma$ which are constant on the conjugacy classes, so when $\Gamma$ has infinite conjugacy classes, called ICC property, what we have is a factor. Thus, as a conclusion, when $\Gamma$ is infinite and has the ICC property, its von Neumann algebra $L(\Gamma)$ is a ${\rm II}_1$ factor.

\bigskip

Summarizing, we have our objects, the ${\rm II}_1$ factors, but what about morphisms. And here, a natural idea is that of looking at the inclusions of such factors:

\begin{definition}
A subfactor is an inclusion of ${\rm II}_1$ factors $A_0\subset A_1$.
\end{definition}

So, these will be the objects that we will be interested in, in what follows. Now given a subfactor $A_0\subset A_1$, a first question is that of defining its index, measuring how big $A_1$ is, when compared to $A_0$. But this can be done as follows:

\index{subfactor}
\index{index of subfactor}
\index{coupling constant}
\index{index of subfactor}

\begin{theorem}
Given an inclusion of ${\rm II}_1$ factors $A_0\subset A_1$, the number
$$N=\frac{\dim_{A_0}H}{\dim_{A_1}H}$$
is independent of the ambient Hilbert space $H$, and is called index.
\end{theorem}

\begin{proof}
This is standard, with the fact that the index as defined by the above formula is independent of the ambient Hilbert space $H$ coming from the various properties of the coupling constant $\dim_AH$, coming from the work of Murray and von Neumann.
\end{proof}

With this discussed, time now for a systematic study of subfactors. Following Jones \cite{jo1}, let us start with the following standard result:

\index{conditional expectation}
\index{Jones projection}

\begin{proposition}
Associated to any subfactor $A_0\subset A_1$ is the orthogonal projection
$$e:L^2(A_1)\to L^2(A_0)$$
producing the conditional expectation $E:A_1\to A_0$ via the following formula:
$$exe=E(x)e$$
This projection is called Jones projection for the subfactor $A_0\subset A_1$.
\end{proposition}

\begin{proof}
This is indeed somehing quite routine. See \cite{jo1}.
\end{proof}

Quite remarkably, the subfactor $A_0\subset A_1$, as well as its commutant, can be recovered from the knowledge of this projection, in the following way:

\begin{proposition}
Given a subfactor $A_0\subset A_1$, with Jones projection $e$, we have
$$A_0=A_1\cap\{e\}'\quad,\quad 
A_0'=(A_1'\cap\{e\})''$$
as equalities of von Neumann algebras, acting on the space $L^2(A_1)$.
\end{proposition}

\begin{proof}
The above two formulae both follow from $exe=E(x)e$, via some elementary computations, and for details here, we refer to Jones' paper \cite{jo1}.
\end{proof}

We are now ready to formulate a key definition, as follows:

\index{basic construction}
\index{Jones basic construction}

\begin{definition}
Associated to any subfactor $A_0\subset A_1$ is the basic construction
$$A_0\subset_eA_1\subset A_2$$
with $A_2=<A_1,e>$ being the algebra generated by $A_1$ and by the Jones projection
$$e:L^2(A_1)\to L^2(A_0)$$
acting on the Hilbert space $L^2(A_1)$.
\end{definition}

The idea now, following \cite{jo1}, will be that $A_1\subset A_2$ appears as a kind of ``reflection'' of $A_0\subset A_1$, and also that the basic construction can be iterated, and with all this leading to non-trivial results. To be more precise, following \cite{jo1}, we have:

\index{Jones tower}

\begin{theorem}
Associated to any subfactor $A_0\subset A_1$ is the Jones tower
$$A_0\subset_{e_1}A_1\subset_{e_2}A_2\subset_{e_3}A_3\subset\ldots\ldots$$
with the Jones projections having the following properties:
\begin{enumerate}
\item $e_i^2=e_i=e_i^*$.

\item $e_ie_j=e_je_i$ for $|i-j|\geq2$.

\item $e_ie_{i\pm1}e_i=[A_1:A_0]^{-1}e_i$.

\item $tr(we_{n+1})=[A_1:A_0]^{-1}tr(w)$, for any word $w\in<e_1,\ldots,e_n>$.
\end{enumerate}
\end{theorem}

\begin{proof}
This follows indeed by doing some routine computations. See \cite{jo1}.
\end{proof}

The relations found in Theorem 12.30 are in fact well-known, from the standard theory of the Temperley-Lieb algebra. This algebra, discovered by Temperley and Lieb in the context of statistical mechanics, has a very simple definition, as follows:

\index{Temperley-Lieb algebra}
\index{noncrossing pairings}
\index{closed circle}
\index{value of circle}

\begin{definition}
The Temperley-Lieb algebra of index $N\in[1,\infty)$ is defined as
$$TL_N(k)=span(NC_2(k,k))$$
with product given by vertical concatenation, with the rule
$$\bigcirc=N$$
for the closed circles that might appear when concatenating.
\end{definition}

As already mentioned, this algebra was discovered by Temperley and Lieb in the context of general statistical mechanics, and we refer here to the physics literature. In what concerns us, still following Jones' paper \cite{jo1}, we have the following result:

\index{subfactor}
\index{Temperley-Lieb algebra}
\index{Jones tower}

\begin{theorem}
Given a subfactor $A_0\subset A_1$, construct its the Jones tower:
$$A_0\subset_{e_1}A_1\subset_{e_2}A_2\subset_{e_3}A_3\subset\ldots\ldots$$
The rescaled sequence of projections $e_1,e_2,e_3,\ldots\in B(H)$ produces then a representation 
$$TL_N\subset B(H)$$
of the Temperley-Lieb algebra of index $N=[A_1:A_0]$.
\end{theorem}

\begin{proof}
We know from Theorem 12.30 that the rescaled sequence of Jones projections $e_1,e_2,e_3,\ldots\in B(H)$ behaves algebrically exactly as the following $TL_N$ diagrams:
$$\varepsilon_1={\ }^\cup_\cap\quad,\quad  
\varepsilon_2=|\!{\ }^\cup_\cap\quad,\quad 
\varepsilon_3=||\!{\ }^\cup_\cap\quad,\quad 
\ldots$$

But these diagrams generate $TL_N$, and so we have an embedding $TL_N\subset B(H)$, where $H$ is the Hilbert space where our subfactor $A_0\subset A_1$ lives, as claimed.
\end{proof}

Quite remarkably, the planar algebra structure of $TL_N$, taken in an intuitive sense, of composing diagrams, extends to a planar algebra structure on $P$. In order to discuss this, let us start with the axioms for planar algebras. Following Jones \cite{jo3}, we have:

\index{planar tangle}
\index{planar algebra}

\begin{definition}
The planar algebras are defined as follows:
\begin{enumerate}
\item We consider rectangles in the plane, with the sides parallel to the coordinate axes, and taken up to planar isotopy, and we call such rectangles boxes.

\item A labeled box is a box with $2n$ marked points on its boundary, $n$ on its upper side, and $n$ on its lower side, for some integer $n\in\mathbb N$.

\item A tangle is labeled box, containing a number of labeled boxes, with all marked points, on the big and small boxes, being connected by noncrossing strings.

\item A planar algebra is a sequence of finite dimensional vector spaces $P=(P_n)$, together with linear maps $P_{n_1}\otimes\ldots\otimes P_{n_k}\to P_n$, one for each tangle, such that the gluing of tangles corresponds to the composition of linear maps.
\end{enumerate}
\end{definition}

In fact things are a bit more complicated than this, but for getting started, this will do. As a basic example of a planar algebra, we have the Temperley-Lieb algebra:

\index{Temperley-Lieb planar algebra}

\begin{theorem}
The Temperley-Lieb algebra $TL_N$, viewed as graded algebra
$$TL_N=(TL_N(n))_{n\in\mathbb N}$$
is a planar algebra, with the corresponding linear maps associated to the planar tangles
$$TL_N(n_1)\otimes\ldots\otimes TL_N(n_k)\to TL_N(n)$$
appearing by putting the various $TL_N(n_i)$ diagrams into the small boxes of the given tangle, which produces a $TL_N(n)$ diagram.
\end{theorem}

\begin{proof}
This is something trivial, which follows from definitions:

\medskip

(1) Assume indeed that we are given a planar tangle $\pi$, as in Definition 12.33, consisting of a box having $2n$ marked points on its boundary, and containing $k$ small boxes, having respectively $2n_1,\ldots,2n_k$ marked points on their boundaries, and then a total of $n+\Sigma n_i$ noncrossing strings, connecting the various $2n+\Sigma 2n_i$ marked points.

\medskip

(2) We want to associate to this tangle $\pi$ a linear map as follows:
$$T_\pi:TL_N(n_1)\otimes\ldots\otimes TL_N(n_k)\to TL_N(n)$$

For this purpose, by linearity, it is enough to construct elements as follows, for any choice of Temperley-Lieb diagrams $\sigma_i\in TL_N(n_i)$, with $i=1,\ldots,k$:
$$T_\pi(\sigma_1\otimes\ldots\otimes\sigma_k)\in TL_N(n)$$

(3) But constructing such an element is obvious, just by putting the various diagrams $\sigma_i\in TL_N(n_i)$ into the small boxes the given tangle $\pi$. Indeed, this procedure produces a certain diagram in $TL_N(n)$, that we can call $T_\pi(\sigma_1\otimes\ldots\otimes\sigma_k)$, as above.

\medskip

(4) Finally, we have to check that everything is well-defined up to planar isotopy, and that the gluing of tangles corresponds to the composition of linear maps. But both these checks are trivial, coming from the definition of $TL_N$, and we are done.
\end{proof}

As a conclusion to all this, $P=TL_N$ is indeed a planar algebra, and of somewhat ``trivial'' type, with the triviality coming from the fact that, in this case, the elements of $P$ are planar diagrams themselves, so the planar structure appears trivially.

\bigskip

In relation now with subfactors, the result, which extends Theorem 12.32, and which was found by Jones in \cite{jo3}, almost 20 years after \cite{jo1}, is as follows:

\index{higher commutant}
\index{planar algebra}

\begin{theorem} 
Given a subfactor $A_0\subset A_1$, the collection $P=(P_n)$ of linear spaces 
$$P_n=A_0'\cap A_n$$
has a planar algebra structure, extending the planar algebra structure of $TL_N$.
\end{theorem}

\begin{proof}
We know from Theorem 12.32 that we have an inclusion as follows, coming from the basic construction, and with $TL_N$ itself being a planar algebra:
$$TL_N\subset P$$

Thus, the whole point is that of proving that the trivial planar algebra structure of $TL_N$ extends into a planar algebra structure of $P$. But this can be done via a long algebraic study, and for the full computation here, we refer to Jones' paper \cite{jo3}.
\end{proof}

Getting now to quantum groups, the above machinery is interesting for us. We will need the construction of the tensor planar algebra $\mathcal T_N$, which is as follows:

\index{tensor planar algebra}
\index{Kronecker symbol}

\begin{definition}
The tensor planar algebra $\mathcal T_N$ is the sequence of vector spaces 
$$P_k=M_N(\mathbb C)^{\otimes k}$$
with the multilinear maps $T_\pi:P_{k_1}\otimes\ldots\otimes P_{k_r}\to P_k$
being given by the formula
$$T_\pi(e_{i_1}\otimes\ldots\otimes e_{i_r})=\sum_j\delta_\pi(i_1,\ldots,i_r:j)e_j$$
with the Kronecker symbols $\delta_\pi$ being $1$ if the indices fit, and being $0$ otherwise.
\end{definition}

In other words, we are using here a construction which is very similar to the construction $\pi\to T_\pi$ that we used for easy quantum groups. We put the indices of the basic tensors on the marked points of the small boxes, in the obvious way, and the coefficients of the output tensor are then given by Kronecker symbols, exactly as in the easy case.

\bigskip

The fact that we have indeed a planar algebra, in the sense that the gluing of tangles corresponds to the composition of linear maps, as required by Definition 12.33, is something elementary, in the same spirit as the verification of the functoriality properties of the correspondence $\pi\to T_\pi$, discussed in chapter 11, and we refer here to Jones \cite{jo3}. 

\bigskip

We have the following result, making the link with quantum groups:

\begin{theorem}
The following happen, at the algebraic level:
\begin{enumerate}
\item The closed subgroups $G\subset O_N^+$ produce planar algebras $P\subset\mathcal T_N$, via the following formula, and any subalgebra $P\subset\mathcal T_N$ appears in this way:
$$P_k=End(u^{\otimes k})$$

\item Also, the closed subgroups $G\subset U_N^+$ produce planar algebras $P\subset\mathcal T_N$, via the following formula, and any subalgebra $P\subset\mathcal T_N$ appears in this way:
$$P_k=End(\underbrace{u\otimes\bar{u}\otimes u\otimes\ldots}_{k\ terms})$$

\item In fact, the closed subgroups $G\subset PO_N^+\simeq PU_N^+$ are in correspondence with the subalgebras $P\subset\mathcal T_N$, with $G\to P$ being given by $P_k=Fix(u^{\otimes k})$.
\end{enumerate}
\end{theorem}

\begin{proof}
There is a long story with this result, the idea being as follows:

\medskip

(1) This is something quite routine, ultimately appearing as a suitable modification of Woronowicz's Tannakian duality in \cite{wo2}. Note that the correspondence is not bijective, because the spaces $P_k$ determine $PG\subset PO_N^+$, but not $G\subset O_N^+$ itself.

\medskip

(2) This is an extension of (1), and the same comments apply. With the extra comment that the fact that the subgroups $PG\subset PO_N^+$ produce the same planar algebras as the subgroups $PG\subset PU_N^+$ should not be surprising, due to $PO_N^+=PU_N^+$.

\medskip

(3) This is an extension of (2), and a further extension of (1), and is in fact the best result on the subject, due to the fact that we have there a true, bijective correspondence. As before, this ultimately comes from Woronowicz's Tannakian duality in \cite{wo2}.
\end{proof}

Finally, in relation with subfactors, the result here is as follows:

\index{fixed point subfactor}
\index{minimal action}
\index{hyperfinite factor}
\index{R}

\begin{theorem}
The planar algebras coming from the subgroups $G\subset PO_N^+=PU_N^+$ appear from fixed point subfactors, of the following type,
$$A^G\subset(M_N(\mathbb C)\otimes A)^G$$
with the action $G\curvearrowright A$ being assumed to be minimal, $(A^G)'\cap A=\mathbb C$.
\end{theorem}

\begin{proof}
Again, there is a long story with this result, and for details about this, and for various generalizations as well, we refer here to \cite{ba1}, \cite{twa} and related papers.
\end{proof}

Finally, let us mention that an important question, which is still open, is that of understanding whether the above subfactors can be taken to be hyperfinite, $A^G\simeq R$. This is related to the axiomatization of hyperfinite subfactors, another open question, which is of central importance in von Neumann algebras. See Jones \cite{jo1}.

\section*{12d. Matrix models}

With the subfactors and planar algebras discussed, which are at the same time more general than the quantum groups, and a bit more concrete as well, with everything happening on a Hilbert space $H$, which is something very nice and concrete, you might wonder, does it really make sense now to go back to algebra, and to quantum groups.

\bigskip

Good point, and in answer, here is a good and healthy quantum group principle:

\begin{principle}
We have concrete Hilbert spaces and concrete operator algebras in quantum groups too, via the notion of matrix model,
$$\pi:C(G)\to M_K(C(T))$$
which involves the random matrix algebra $M_K(C(T))$. So, for anything advanced, truly competing with subfactors and planar algebras, we should work more on such models.
\end{principle}

Getting now to work, we already know a bit about such models, from the previous chapters, but time now to have a more systematic look at all this, using our recent quantum group knowledge. As a first task, let us review the notion of Hopf image:
$$\pi:C(G)\to C(H)\to M_K(C(T))$$

We recall that the existence and uniqueness of the Hopf image come by dividing the algebra $C(G)$ by a suitable ideal. Alternatively, in Tannakian terms, we have:

\index{Tannakian duality}

\begin{theorem}
Assuming $G\subset U_N^+$, with fundamental corepresentation $u=(u_{ij})$, the Hopf image of a model $\pi:C(G)\to M_K(C(T))$ comes from the Tannakian category
$$C_{kl}=Hom(U^{\otimes k},U^{\otimes l})$$
where $U_{ij}=\pi(u_{ij})$, and where the spaces on the right are taken in a formal sense.
\end{theorem}

\begin{proof}
Since the morphisms increase the intertwining spaces, when defined either in a representation theory sense, or just formally, we have inclusions as follows:
$$Hom(u^{\otimes k},u^{\otimes l})\subset Hom(U^{\otimes k},U^{\otimes l})$$

More generally, we have such inclusions when replacing $(G,u)$ with any pair producing a factorization of $\pi$. Thus, by Tannakian duality, the Hopf image must be given by the fact that the intertwining spaces must be the biggest, subject to the above inclusions. On the other hand, since $u$ is biunitary, so is $U$, and it follows that the spaces on the right form a Tannakian category. Thus, we have a quantum group $(H,v)$ given by:
$$Hom(v^{\otimes k},v^{\otimes l})=Hom(U^{\otimes k},U^{\otimes l})$$

By the above discussion, $C(H)$ follows to be the Hopf image of $\pi$, as claimed.
\end{proof}

Regarding now the study of the inner faithful models, a key problem is that of computing the Haar integration functional, and we have here the following result:

\index{truncated integrals}

\begin{theorem}
Given an inner faithful model $\pi:C(G)\to M_K(C(T))$, we have
$$\int_G=\lim_{k\to\infty}\frac{1}{k}\sum_{r=1}^k\int_G^r$$
with the truncations of the integration on the right being given by
$$\int_G^r=(\varphi\circ\pi)^{*r}$$
with $\phi*\psi=(\phi\otimes\psi)\Delta$, and with $\varphi=tr\otimes\int_T$ being the random matrix trace.
\end{theorem}

\begin{proof}
This is something quite tricky, the idea being as follows:

\medskip

(1) As a first observation, there is an obvious similarity here with the Woronowicz construction of the Haar measure, explained in chapter 10. In fact, the above result holds more generally for any model $\pi:C(G)\to B$, with $\varphi\in B^*$ being a faithful trace.

\medskip

(2) In order to prove now the result, we can proceed as in chapter 10. If we denote by $\int_G'$ the limit in the statement, we must prove that this limit converges, and that:
$$\int_G'=\int_G$$

It is enough to check this on the coefficients of the Peter-Weyl corepresentations, and if we let $v=u^{\otimes k}$ be one of these corepresentations, we must prove that we have:
$$\left(id\otimes\int_G'\right)v=\left(id\otimes\int_G\right)v$$

(3) In order to prove this, we already know, from the Haar measure theory from chapter 1, that the matrix on the right is the orthogonal projection onto $Fix(v)$:
$$\left(id\otimes\int_G\right)v=Proj\Big[Fix(v)\Big]$$

Regarding now the matrix on the left, the trick in \cite{wo1} applied to the linear form $\varphi\pi$ tells us that this is the orthogonal projection onto the $1$-eigenspace of $(id\otimes\varphi\pi)v$:
$$\left(id\otimes\int_G'\right)v=Proj\Big[1\in (id\otimes\varphi\pi)v\Big]$$

(4) Now observe that, if we set $V_{ij}=\pi(v_{ij})$, we have the following formula:
$$(id\otimes\varphi\pi)v=(id\otimes\varphi)V$$

Thus, we can apply the trick in \cite{wo1}, and we conclude that the $1$-eigenspace that we are interested in equals $Fix(V)$. But, according to Theorem 12.40, we have:
$$Fix(V)=Fix(v)$$

Thus, we have proved that we have $\int_G'=\int_G$, as desired.
\end{proof}

As a further piece of general theory, regarding models, let us formulate:

\begin{definition}
A model $\pi:C(G)\to M_K(C(T))$ is called stationary when
$$\int_G=\left(tr\otimes\int_T\right)\pi$$
where $\int_T$ is the integration with respect to a given probability measure on $T$.
\end{definition}

This is of course a quite specialized notion, because by basic functional analysis, the model must be faithful, and the algebra $C(G)$ must be of type I. However, there are many interesting examples of such models, as we have already seen in chapter 8.

\bigskip

In order to detect the stationary models, we have the following useful criterion, mixing linear algebra and analysis, based on the Ces\`aro formula from Theorem 12.41:

\index{idempotent state}
\index{stationary on its image}

\begin{theorem}
For a model $\pi:C(G)\to M_K(C(T))$, the following are equivalent:
\begin{enumerate}
\item $Im(\pi)$ is a Hopf algebra, and the Haar integration on it is:
$$\psi=\left(tr\otimes\int_T\right)\pi$$

\item The linear form $\psi=(tr\otimes\int_T)\pi$ satisfies the idempotent state property:
$$\psi*\psi=\psi$$

\item We have $T_e^2=T_e$, $\forall p\in\mathbb N$, $\forall e\in\{1,*\}^p$, where:
$$(T_e)_{i_1\ldots i_p,j_1\ldots j_p}=\left(tr\otimes\int_T\right)(U_{i_1j_1}^{e_1}\ldots U_{i_pj_p}^{e_p})$$
\end{enumerate}
If these conditions are satisfied, we say that $\pi$ is stationary on its image.
\end{theorem}

\begin{proof}
Given a matrix model $\pi:C(G)\to M_K(C(T))$ as in the statement, we can factorize it via its Hopf image, as in chapter 4, or in Theorem 12.40:
$$\pi:C(G)\to C(H)\to M_K(C(T))$$

Now observe that (1,2,3) above depend only on the factorized representation:
$$\nu:C(H)\to M_K(C(T))$$

Thus, we can assume in practice that we have $G=H$, which means that we can assume that $\pi$ is inner faithful. With this assumption made, the formula in Theorem 12.41 applies to our situation, and the proof of the equivalences goes as follows:

\medskip

$(1)\implies(2)$ This is clear from definitions, because the Haar integration on any compact quantum group satisfies the following idempotent state equation:
$$\psi*\psi=\psi$$

$(2)\implies(1)$ Assuming $\psi*\psi=\psi$, we have $\psi^{*r}=\psi$ for any $r\in\mathbb N$, and Theorem 12.41 gives $\int_G=\psi$. Thus, via some standard functional analysis, we obtain the result.

\medskip

In order to establish now $(2)\Longleftrightarrow(3)$, we use the following elementary formula, which comes from the definition of the convolution operation:
$$\psi^{*r}(u_{i_1j_1}^{e_1}\ldots u_{i_pj_p}^{e_p})=(T_e^r)_{i_1\ldots i_p,j_1\ldots j_p}$$

\medskip

$(2)\implies(3)$ Assuming $\psi*\psi=\psi$, by using the above formula at $r=1,2$ we obtain that the matrices $T_e$ and $T_e^2$ have the same coefficients, and so they are equal.

\medskip

$(3)\implies(2)$ Assuming $T_e^2=T_e$, by using the above formula at $r=1,2$ we obtain that the linear forms $\psi$ and $\psi*\psi$ coincide on any product of coefficients $u_{i_1j_1}^{e_1}\ldots u_{i_pj_p}^{e_p}$. Now since these coefficients span a dense subalgebra of $C(G)$, this gives the result.
\end{proof}

So long for the basic theory of the random matrix models. At the level of concrete applications now, we have already seen some in chapter 8, and with the above theorems providing the key for fully understanding that material, say at a second reading. 

\bigskip

Also, following Bichon-Dubois-Violette \cite{bd1} and related work, there are some interesting applications of this to the notion of half-liberation, the idea being that the half-liberations have stationary models, obtained by using antidiagonal $2\times2$ matrices.

\bigskip

Finally, in relation with Principle 12.39 that we started with, there are some deep truths hidden there, the idea being that our quantum groups and matrix models can potentially help in relation with certain questions from statistical mechanics. And here, again, we will leave some exploration of what is known and what is not, to you reader.

\section*{12e. Exercises}

We had a fairly advanced chapter here, and as exercises, we have:

\begin{exercise}
Learn, with details, the Drinfeld-Jimbo deformation theory.
\end{exercise}

\begin{exercise}
Prove that the quantum reflection groups equal their own twists.
\end{exercise}

\begin{exercise}
Learn more about discrete quantum groups, and their growth.
\end{exercise}

\begin{exercise}
Learn von Neumann algebra theory, as much as you can.
\end{exercise}

\begin{exercise}
Learn about the classification of subfactors of index $N\leq4$.
\end{exercise}

\begin{exercise}
Learn about the $2\times2$ matrix models for half-liberations.
\end{exercise}

As bonus exercise, learn about the hyperfinite factor $R$, and its subfactors.

\part{Quantum spaces}

\ \vskip50mm

\begin{center}
{\em Sa te gandesti din cand in cand la mine

Si-al stelelor stralucitor sirag

Sa se aprinda noaptea pentru tine

Ca niste ochi ce te asteapta-n prag}
\end{center}

\chapter{Finite spaces}

\section*{13a. Twisted spaces}

Welcome to noncommutative geometry. In this last part of the present book we discuss how what we learned to far, concerning the quantum groups, can be extended to more general classes of ``quantum manifolds''. That is, we would like to reach to a full theory of noncommutative geometry, and with the hope that this can be useful for physics.

\bigskip

We will be quite modest in our purposes and goals, which will be basically introductory. Among others, with the art of geometry, even in the classical case, being quite hard to develop over arbitrary fields $F$, we will develop here our theory over $F=\mathbb C$, as a continuation of the material from Part III, which was over $F=\mathbb C$ too. 

\bigskip

However, and here comes the point, abstract algebraic aspects will not be forgotten, and there will be a philosophical relation too with the material from Parts I-II, for the most dealing with arbitrary fields $F$. We will talk about this in the final chapter, 16.

\bigskip

In order to get started now, let us first talk about the finite quantum spaces. In view of the general $C^*$-algebra theory explained before, we have the following definition:

\begin{definition}
A finite quantum space $Z$ is the abstract dual of a finite dimensional $C^*$-algebra $B$, according to the following formula:
$$C(Z)=B$$
The formal number of elements of such a space is $|Z|=\dim B$. By decomposing the algebra $B$, we have a formula of the following type:
$$C(Z)=M_{n_1}(\mathbb C)\oplus\ldots\oplus M_{n_k}(\mathbb C)$$
With $n_1=\ldots=n_k=1$ we obtain in this way the space $Z=\{1,\ldots,k\}$. Also, when $k=1$ the equation is $C(Z)=M_n(\mathbb C)$, and the solution will be denoted $Z=M_n$.
\end{definition}

In order to do some mathematics on such spaces, the very first observation is that we can talk about the formal number of points of such a space, as follows:
$$|Z|=\dim B$$

Alternatively, by decomposing the algebra $B$ as a sum of matrix algebras, as in Definition 13.1, we have the following formula for the formal number of points:
$$|Z|=n_1^2+\ldots+n_k^2$$

Pictorially, this suggests representing $Z$ as a set of $|Z|$ points in the plane, arranged in squares having sides $n_1,\ldots,n_k$, coming from the matrix blocks of $B$, as follows:
$$\begin{matrix}
\circ&\circ&\circ\\
\circ&\circ&\circ&&\ldots&&\circ&\circ\\
\circ&\circ&\circ&&&&\circ&\circ
\end{matrix}$$

As a second piece of mathematics, we can talk about counting measures, as follows:

\begin{definition}
Given a finite quantum space $Z$, we construct the functional
$$tr:C(Z)\to B(l^2(Z))\to\mathbb C$$
obtained by applying the regular representation, and the normalized matrix trace, and we call it integration with respect to the normalized counting measure on $Z$.
\end{definition}

To be more precise, consider the algebra $B=C(Z)$, which is by definition finite dimensional. We can make act $B$ on itself, by left multiplication:
$$\pi:B\to\mathcal L(B)\quad,\quad 
a\to(b\to ab)$$

The target of $\pi$ being a matrix algebra, $\mathcal L(B)\simeq M_N(\mathbb C)$ with $N=\dim B$, we can further compose with the normalized matrix trace, and we obtain $tr$:
$$tr=\frac{1}{N}\,Tr\circ\pi$$

As basic examples, for both $Z=\{1,\ldots,N\}$ and $Z=M_N$ we obtain the usual trace. In general, with $C(Z)=M_{n_1}(\mathbb C)\oplus\ldots\oplus M_{n_k}(\mathbb C)$, the weights of $tr$ are:
$$c_i=\frac{n_i^2}{\sum_in_i^2}$$

Pictorially, this suggests fine-tuning our previous picture of $Z$, by adding to each point the unnormalized trace of the corresponding element of $B$, as follows:
$$\begin{matrix}
\bullet_{n_1}&\circ_0&\circ_0\\
\circ_0&\bullet_{n_1}&\circ_0&&\ldots&&\bullet_{n_k}&\circ_0\\
\circ_0&\circ_0&\bullet_{n_1}&&&&\circ_0&\bullet_{n_k}
\end{matrix}$$

Here we have represented the points on the diagonals with solid circles, since they are of different nature from the off-diagonal ones, the attached numbers being nonzero. However, this picture is not complete either, and we can do better, as follows:

\begin{definition}
Given a finite quantum space $Z$, coming via a formula of type
$$C(Z)=M_{n_1}(\mathbb C)\oplus\ldots\oplus M_{n_k}(\mathbb C)$$
we use the following equivalent conventions for drawing $Z$:
\begin{enumerate}
\item Triple indices. We represent $Z$ as a set of $N=|Z|$ points, with each point being decorated with a triple index $ija$, coming from the standard basis $\{e_{ij}^a\}\subset B$. 

\item Double indices. As before, but by ignoring the index $a$, with the convention that $i,j$ belong to various indexing sets, one for each of the matrix blocks of $B$.

\item Single indices. As before, but with each point being now decorated with a single index, playing the role of the previous triple indices $ija$, or double indices $ij$.
\end{enumerate}
\end{definition}

As an illustration, consider the space $Z=\{1,\ldots,k\}$. According to our single index convention, we can represent this space as a set of $k$ points, decorated by some indices, which must be chosen different. Thus, we are led to the following picture:
$$\bullet_1\quad\bullet_2\quad\ldots\quad\bullet_k$$

As another illustration, consider the space $Z=M_n$. Here the picture is as follows, using double indices, which can be regarded as well as being single indices:
$$\begin{matrix}
\bullet_{11}&\circ_{12}&\circ_{13}\\
\circ_{21}&\bullet_{22}&\circ_{23}\\
\circ_{31}&\circ_{32}&\bullet_{33}
\end{matrix}$$

\smallskip

As yet another illustration, for the space $Z=M_3\sqcup M_2$, which appears by definition from the algebra $B=M_3(\mathbb C)\oplus M_2(\mathbb C)$, we are in need of triple indices, which can be of course regarded as single indices, in order to label all the points, and the picture is:
$$\begin{matrix}
\bullet_{111}&\circ_{121}&\circ_{131}\\
\circ_{211}&\bullet_{221}&\circ_{231}&&\ &&\bullet_{112}&\circ_{122}\\
\circ_{311}&\circ_{321}&\bullet_{331}&&&&\circ_{212}&\bullet_{222}
\end{matrix}$$

\smallskip

Let us study now the quantum group actions $G\curvearrowright Z$. If we denote by $\mu,\eta$ the multiplication and unit map of the algebra $C(Z)$, we have the following result:

\begin{proposition}
Consider a linear map $\Phi:C(Z)\to C(Z)\otimes C(G)$, written as
$$\Phi(e_i)=\sum_je_j\otimes u_{ji}$$
with $\{e_i\}$ being a linear space basis of $C(Z)$, chosen orthonormal with respect to $tr$.
\begin{enumerate}
\item $\Phi$ is a linear space coaction $\iff$ $u$ is a corepresentation.

\item $\Phi$ is multiplicative $\iff$ $\mu\in Hom(u^{\otimes 2},u)$.

\item $\Phi$ is unital $\iff$ $\eta\in Hom(1,u)$.

\item $\Phi$ leaves invariant $tr$ $\iff$ $\eta\in Hom(1,u^*)$.

\item If these conditions hold, $\Phi$ is involutive $\iff$ $u$ is unitary.
\end{enumerate}
\end{proposition}

\begin{proof}
This is similar to the proof for $S_N^+$ from chapter 6, as follows:

\medskip

(1) There are two axioms to be processed here, and we have indeed:
$$(id\otimes\Delta)\Phi=(\Phi\otimes id)\Phi
\iff\Delta(u_{ji})=\sum_ku_{jk}\otimes u_{ki}$$
$$(id\otimes\varepsilon)\Phi=id
\iff\varepsilon(u_{ji})=\delta_{ji}$$

(2) By using $\Phi(e_i)=u(e_i\otimes 1)$ we have the following identities, which give the result:
$$\Phi(e_ie_k)
=u(\mu\otimes id)(e_i\otimes e_k\otimes 1)$$
$$\Phi(e_i)\Phi(e_k)
=(\mu\otimes id)u^{\otimes 2}(e_i\otimes e_k\otimes 1)$$

(3) From $\Phi(e_i)=u(e_i\otimes1)$ we obtain by linearity, as desired:
$$\Phi(1)=u(1\otimes1)$$

(4) This follows from the following computation, by applying the involution:
\begin{eqnarray*}
(tr\otimes id)\Phi(e_i)=tr(e_i)1
&\iff&\sum_jtr(e_j)u_{ji}=tr(e_i)1\\
&\iff&\sum_ju_{ji}^*1_j=1_i\\
&\iff&(u^*1)_i=1_i\\
&\iff&u^*1=1
\end{eqnarray*}

(5) Assuming that (1-4) are satisfied, and that $\Phi$ is involutive, we have:
\begin{eqnarray*}
(u^*u)_{ik}
&=&\sum_lu_{li}^*u_{lk}\\
&=&\sum_{jl}tr(e_j^*e_l)u_{ji}^*u_{lk}\\
&=&(tr\otimes id)\sum_{jl}e_j^*e_l\otimes u_{ji}^*u_{lk}\\
&=&(tr\otimes id)(\Phi(e_i)^*\Phi(e_k))\\
&=&(tr\otimes id)\Phi(e_i^*e_k)\\
&=&tr(e_i^*e_k)1\\
&=&\delta_{ik}
\end{eqnarray*}

Thus $u^*u=1$, and since we know from (1) that $u$ is a corepresentation, it follows that $u$ is unitary. The proof of the converse is standard too, by using a similar computation.
\end{proof}

Following now \cite{ba2}, \cite{wan}, we have the following result, extending the basic theory of $S_N^+$ from chapter 6 to the present finite quantum space setting:

\begin{theorem}
Given a finite quantum space $Z$, there is a universal compact quantum group $S_Z^+$ acting on $Z$, and leaving the counting measure invariant. We have
$$C(S_Z^+)=C(U_N^+)\Big/\Big<\mu\in Hom(u^{\otimes2},u),\eta\in Fix(u)\Big>$$
where $N=|Z|$, and where $\mu,\eta$ are the multiplication and unit maps of the algebra $C(Z)$. For the classical space $Z=\{1,\ldots,N\}$ we have $S_Z^+=S_N^+$. 
\end{theorem}

\begin{proof}
Here the first two assertions follow from Proposition 13.4, by using the standard fact that the complex conjugate of a corepresentation is a corepresentation too. As for the last assertion, regarding $S_N^+$, this follows from the results in chapter 6.
\end{proof}

The above result is quite conceptual, and we will see some applications in a moment. However, for many concrete questions, nothing beats multimatrix bases and indices. So, following the original paper of Wang \cite{wan}, let us discuss this. We first have:

\begin{definition}
Given a finite quantum space $Z$, we let $\{e_i\}$ be the standard basis of $B=C(Z)$, so that the multiplication, involution and unit of $B$ are given by
$$e_ie_j=e_{ij}\quad,\quad 
e_i^*=e_{\bar{i}}\quad,\quad
1=\sum_{i=\bar{i}}e_i$$
where $(i,j)\to ij$ is the standard partially defined multiplication on the indices, with the convention $e_\emptyset=0$, and where $i\to\bar{i}$ is the standard involution on the indices.
\end{definition}

To be more precise, let $\{e_{ab}^r\}\subset B$ be the multimatrix basis. We set $i=(abr)$, and with this convention, the multiplication, coming from $e_{ab}^re_{cd}^p=\delta_{rp}\delta_{bc}e_{ad}^r$, is given by:
$$(abr)(cdp)=\begin{cases}
(adr)&{\rm if}\ b=c,\ r=p\\
\emptyset&{\rm otherwise}
\end{cases}$$

Regarding now the generalized quantum permutation groups $S_Z^+$, the construction in Theorem 13.5 reformulates as follows, by using the above formalism:

\begin{proposition}
Given a finite quantum space $Z$, with basis $\{e_i\}\subset C(Z)$ as above, the algebra $C(S_Z^+)$ is generated by variables $u_{ij}$ with the following relations,
$$\sum_{ij=p}u_{ik}u_{jl}=u_{p,kl}\quad,\quad\sum_{kl=p}u_{ik}u_{jl}=u_{ij,p}$$
$$\sum_{i=\bar{i}}u_{ij}=\delta_{j\bar{j}}\quad,\quad\sum_{j=\bar{j}}u_{ij}=\delta_{i\bar{i}}$$
$$u_{ij}^*=u_{\bar{i}\hskip0.3mm\bar{j}}$$
with the fundamental corepresentation being the matrix $u=(u_{ij})$. We call a matrix $u=(u_{ij})$ satisfying the above relations ``generalized magic''.
\end{proposition}

\begin{proof}
We recall from Theorem 13.5 that the algebra $C(S_Z^+)$ appears as follows, where $N=|Z|$, and where $\mu,\eta$ are the multiplication and unit maps of $C(Z)$:
$$C(S_Z^+)=C(U_N^+)\Big/\Big<\mu\in Hom(u^{\otimes2},u),\eta\in Fix(u)\Big>$$

But the relations $\mu\in Hom(u^{\otimes2},u)$ and $\eta\in Fix(u)$ produce the 1st and 4th relations in the statement, then the biunitarity of $u$ gives the 5th relation, and finally the 2nd and 3rd relations follow from the 1st and 4th relations, by using the antipode.
\end{proof}

As an illustration, consider the case $Z=\{1,\ldots,N\}$. Here our relations are as follows, corresponding to the standard magic conditions on a matrix $u=(u_{ij})$:
$$u_{ik}u_{il}=\delta_{kl}u_{ik}\quad,\quad u_{ik}u_{jk}=\delta_{ij}u_{ik}$$
$$\sum_iu_{ij}=1\quad,\quad\sum_ju_{ij}=1$$
$$u_{ij}^*=u_{ij}$$

As a second illustration now, which is something new, we have:

\begin{theorem}
For the space $Z=M_2$, coming via $C(Z)=M_2(\mathbb C)$, we have 
$$S_Z^+=SO_3$$
with the action $SO_3\curvearrowright M_2(\mathbb C)$ being the standard one, coming from $SU_2\to SO_3$.
\end{theorem}

\begin{proof}
This is something quite tricky, the idea being as follows:

\medskip

(1) First, we have an action by conjugation $SU_2\curvearrowright M_2(\mathbb C)$, and this action produces, via the canonical quotient map $SU_2\to SO_3$, an action as follows:
$$SO_3\curvearrowright M_2(\mathbb C)$$

(2) Then, it is routine to check, by using computations like those from the proof of $S_N^+=S_N$ at $N\leq3$, from chapter 6, that any action $G\curvearrowright M_2(\mathbb C)$ must come from a classical group. Thus the action $SO_3\curvearrowright M_2(\mathbb C)$ is universal, as claimed.

\medskip

(3) This was for the idea, and we will actually come back to this in a moment, in a more general setting, and with a new proof, complete this time. 
\end{proof}

Let us develop now some basic theory for the quantum symmetry groups $S_Z^+$, and their closed subgroups $G\subset S_Z^+$. We have here the following key result, from \cite{ba1}:

\begin{theorem}
The quantum groups $S_Z^+$ have the following properties: 
\begin{enumerate}
\item The associated Tannakian categories are $TL_N$, with $N=|F|$.

\item The main character follows the Marchenko-Pastur law $\pi_1$, when $|Z|\geq4$.

\item The fusion rules for $S_Z^+$ with $|Z|\geq4$ are the same as for $SO_3$.
\end{enumerate}
\end{theorem}

\begin{proof}
This result is from \cite{ba1}, the idea being as follows:

\medskip

(1) Let us pick our orthogonal basis $\{e_i\}$ as in Definition 13.6, so that we have, for a certain involution $i\to\bar{i}$ on the index set, the following formula:
$$e_i^*=e_{\bar{i}}$$

With this convention, we have the following computation:
\begin{eqnarray*}
\Phi(e_i)=\sum_je_j\otimes u_{ji}
&\implies&\Phi(e_i)^*=\sum_je_j^*\otimes u_{ji}^*\\
&\implies&\Phi(e_{\bar{i}})=\sum_je_{\bar{j}}\otimes u_{ji}^*\\
&\implies&\Phi(e_i)=\sum_je_j\otimes u_{\bar{i}\hskip0.3mm\bar{j}}^*
\end{eqnarray*}

Thus $u_{ji}^*=u_{\bar{i}\hskip0.3mm\bar{j}}$, so $u\sim\bar{u}$. Now with this result in hand, the proof goes as for the proof for $S_N^+$, from chapter 6. To be more precise, the result follows from the fact that the multiplication and unit of any complex algebra, and in particular of the algebra $C(Z)$ that we are interested in here, can be modelled by the following two diagrams:
$$m=|\cup|\qquad,\qquad u=\cap$$

Indeed, this is certainly true algebrically, and well-known, with as an illustration here, the associativity formula $m(m\otimes id)=(id\otimes m)m$ being checked as follows:
$$\begin{matrix}
|&\cup&|\ |&|\\
|&&\cup&|
\end{matrix}
\ \ =\ \ \begin{matrix}
|&|\ |&\cup&|\\
|&\cup&&|
\end{matrix}$$

As in what regards the $*$-structure, things here are fine too, because our choice for the trace from Definition 13.2 leads to the following formula regarding the adjoints, corresponding to $mm^*=N$, and so to the basic Temperley-Lieb calculus rule $\bigcirc=N$:
$$\mu\mu^*=N\cdot id$$

We conclude that the Tannakian category associated to $S_Z^+$ is, as claimed:
\begin{eqnarray*}
C
&=&<\mu,\eta>\\
&=&<m,u>\\
&=&<|\cup|,\cap>\\
&=&TL_N
\end{eqnarray*}

(2) The proof here is exactly as for $S_N^+$, by using moments. To be more precise, according to (1) these moments are the Catalan numbers, which are the moments of $\pi_1$.

\medskip

(3) Once again same proof as for $S_N^+$, by using the fact that the moments of $\chi$ are the Catalan numbers, which naturally leads to the Clebsch-Gordan rules.
\end{proof}

We can merge and reformulate our main results so far in the following way:

\begin{theorem}
The quantun groups $S_Z^+$ have the following properties: 
\begin{enumerate}
\item For $Z=\{1,\ldots,N\}$ we have $S_Z^+=S_N^+$.

\item For the space $Z=M_N$ we have $S_Z^+=PO_N^+=PU_N^+$.

\item In particular, for the space $Z=M_2$ we have $S_Z^+=SO_3$.

\item The fusion rules for $S_Z^+$ with $|Z|\geq4$ are independent of $Z$.

\item Thus, the fusion rules for $S_Z^+$ with $|Z|\geq4$ are the same as for $SO_3$.
\end{enumerate}
\end{theorem}

\begin{proof}
This is basically a compact form of what has been said above, with a new result added, and with some technicalities left aside, the idea being as follows:

\medskip

(1) This is something that we know from Theorem 13.5.

\medskip

(2) We recall from chapter 4 that we have $PO_N^+=PU_N^+$. Consider the standard vector space action of the free unitary group $U_N^+$, and its adjoint action:
$$U_N^+\curvearrowright\mathbb C^N\quad,\quad 
PU_N^+\curvearrowright M_N(\mathbb C)$$

By universality of $S_{M_N}^+$, we must have inclusions as follows:
$$PO_N^+\subset PU_N^+\subset S_{M_N}^+$$

On the other hand, the main character of $O_N^+$ with $N\geq2$ being semicircular, the main character of $PO_N^+$ must be Marchenko-Pastur. Thus the inclusion $PO_N^+\subset S_{M_N}^+$ has the property that it keeps fixed the law of main character, and by Peter-Weyl theory we conclude that this inclusion must be an isomorphism, as desired.

\medskip

(3) This is something that we know from Theorem 13.9, and that can be deduced as well from (2), by using the formula $PO_2^+=SO_3$, which is something elementary. Alternatively, this follows without computations from (4) below, because the inclusion of quantum groups $SO_3\subset S_{M_2}^+$ has the property that it preserves the fusion rules.

\medskip

(4) This is something that we know from Theorem 13.9.

\medskip

(5) This follows from (3,4), as already pointed out in Theorem 13.9.
\end{proof}

As an application of our extended formalism, the Cayley theorem for the finite quantum groups holds in the $S_Z^+$ setting. We have indeed the following result, from \cite{bbn}:

\begin{theorem}
Any finite quantum group $G$ has a Cayley embedding, as follows:
$$G\subset S_G^+$$
However, there are finite quantum groups which are not quantum permutation groups.
\end{theorem}

\begin{proof}
There are two statements here, the idea being as follows:

\medskip

(1) We have an action $G\curvearrowright G$, which leaves invariant the Haar measure. Now since the counting measure is left and right invariant, so is the Haar measure. We conclude that $G\curvearrowright G$ leaves invariant the counting measure, and so $G\subset S_G^+$, as claimed.

\medskip

(2) Regarding the second assertion, this is something non-trivial, the simplest counterexample being a certain quantum group $G$ appearing as a split abelian extension associated to the factorization $S_4=\mathbb Z_4S_3$, having cardinality $|G|=24$. 
\end{proof}

Finally, some interesting phenomena appear in the ``homogeneous'' case, where our quantum space is of the form $Z=M_K\times \{1\ldots,L\}$. Here we first have:

\begin{proposition}
The classical symmetry group of $Z=M_K\times \{1\ldots,L\}$ is
$$S_Z=PU_K\wr S_L$$
with on the right a wreath product, equal by definition to $PU_K^L\rtimes S_L$.
\end{proposition}

\begin{proof}
The fact that we have an inclusion $PU_K\wr S_L\subset S_Z$ is standard, and this follows as well by taking the classical version of the inclusion $PU_K^+\wr_*S_L^+\subset S_Z^+$, established below. As for the fact that this inclusion $PU_K\wr S_L\subset S_Z$ is an isomorphism, this can be proved by picking an arbitrary element $g\in S_Z$, and decomposing it.
\end{proof}

Quite surprisingly, the quantum analogue of the above result fails:

\begin{theorem}
The quantum symmetry group of $Z=M_K\times \{1\ldots,L\}$ satisfies:
$$PU_K^+\wr_*S_L^+\subset S_Z^+$$
However, this inclusion is not an isomorphism at $K,L\geq2$.
\end{theorem}

\begin{proof}
We have several assertions to be proved, the idea being as follows:

\medskip

(1) The fact that we have $PU_K^+\wr_*S_L^+\subset S_Z^+$ is well-known and routine, by checking the fact that the matrix $w_{ija,klb}=u_{ij,kl}^{(a)}v_{ab}$ is a generalized magic unitary.

\medskip

(2) The inclusion $PU_K^+\wr_*S_L^+\subset S_Z^+$ is not an isomorphism, by using standard results, along with the fact that $\pi_1\boxtimes\pi_1\neq\pi_1$ where $\pi_1$ is the Marchenko-Pastur law.
\end{proof}

\section*{13b. Quantum graphs}

Let us start with the following straightforward extension of the usual notion of finite graph, obtained by using a finite quantum space as set of vertices:

\begin{definition}
We call ``finite quantum graph'' a pair of type
$$X=(Z,d)$$
with $Z$ being a finite quantum space, and with $d\in M_N(\mathbb C)$ being a matrix.
\end{definition}

This is of course something quite general. In the case $Z=\{1,\ldots,N\}$ for instance, what we have here is a directed graph, with the edges $i\to j$ colored by complex numbers $d_{ij}\in\mathbb C$, and with self-edges $i\to i$ allowed too, again colored by numbers $d_{ii}\in\mathbb C$. In the general case, however, where $Z$ is arbitrary, the need for extra conditions of type $d=d^*$, or $d_{ii}=0$, or $d\in M_N(\mathbb R)$, or $d\in M_N(0,1)$ and so on, is not very natural, as we will soon discover, and it is best to use Definition 13.14 as such, with no restrictions on $d$.

\bigskip

In general, a quantum graph can be represented as a colored oriented graph on $\{1,\ldots,N\}$, where $N=|Z|$, with the vertices being decorated by single indices $i$, and with the colors being complex numbers, namely the entries of $d$. This is similar to the formalism from before, but there is a discussion here in what regards the exact choice of the colors, which are usually irrelevant in connection with our symmetry problematics, and so can be true colors instead of complex numbers. More on this later.

\bigskip

With the above notion in hand, we have the following definition:

\begin{definition}
The quantum automorphism group of $X=(Z,d)$ is the subgroup 
$$G^+(X)\subset S_Z^+$$
obtained via the relation $du=ud$, where $u=(u_{ij})$ is the fundamental corepresentation.
\end{definition}

Again, this is something very natural, coming as a continuation of the constructions for usual graphs. We refer to the literature for more on this notion, and for a number of advanced computations, in relation with free wreath products. At an elementary level, a first problem is that of working out the basics of the correspondence $X\to G^+(X)$. There are several things to be done here, namely simplices, complementation, color independence, multi-simplices, and with a few twists, all this basically extends well. 

\bigskip

Let us start with the simplices. As we will soon discover, things are quite tricky here, leading us in particular to the conclusion that the simplex based on an arbitrary finite quantum space $F$ is not a usual graph, with $d\in M_N(0,1)$ where $N=|Z|$, but rather a sort of ``signed graph'', with $d\in M_N(-1,0,1)$. Let us start our study with:

\begin{theorem}
Given a finite quantum space $Z$, we have
$$G^+(Z_{empty})=G^+(Z_{full})=S_Z^+$$
where $Z_{empty}$ is the empty graph on the vertex set $Z$, coming from the matrix $d=0$, and where $Z_{full}$ is the simplex on the vertex set $Z$, coming from the matrix 
$$d=NP_1-1_N$$
where $N=|Z|$, and where $P_1$ is the orthogonal projection on the unit $1\in C(Z)$.
\end{theorem}

\begin{proof}
This is something quite tricky, the idea being as follows:

\medskip

(1) First of all, the formula $G^+(Z_{empty})=S_Z^+$ is clear from definitions, because the commutation of $u$ with the matrix $d=0$ is automatic.

\medskip

(2) Regarding $G^+(Z_{full})=S_Z^+$, let us first discuss the classical case, $Z=\{1,\ldots,N\}$. Here the simplex $Z_{full}$ is the graph having having edges between any two vertices, whose adjacency matrix is $d=\mathbb I_N-1_N$, where $\mathbb I_N$ is the all-1 matrix. The commutation of $u$ with $1_N$ being automatic, and the commutation with $\mathbb I_N$ being automatic too, $u$ being bistochastic, we have $[u,d]=0$, and so $G^+(Z_{full})=S_Z^+$ in this case, as stated.

\medskip

(3) In the general case, we know from Theorem 13.5 that we have $\eta\in Fix(u)$, with $\eta:\mathbb C\to C(Z)$ being the unit map. Thus we have $P_1\in End(u)$, and so the condition $[u,P_1]=0$ is automatic. Together with the fact that in the classical case we have $\mathbb I_N=NP_1$, this suggests to define the adjacency matrix of the simplex as being $d=NP_1-1_N$, and with this definition, we have indeed $G^+(Z_{full})=S_Z^+$, as claimed.
\end{proof}

Let us study now the simplices $F_{full}$ found in Theorem 13.16. In the classical case, $Z=\{1,\ldots,N\}$, what we have is of course the usual simplex. However, in the general case things are more mysterious, the first result here being as follows:

\begin{proposition}
The adjacency matrix of the simplex $Z_{full}$, given by definition by $d=NP_1-1_N$, is a matrix $d\in M_N(-1,0,1)$, which can be computed as follows:
\begin{enumerate}
\item In single index notation, $d_{ij}=\delta_{i\bar{i}}\delta_{j\bar{j}}-\delta_{ij}$.

\item In double index notation, $d_{ab,cd}=\delta_{ab}\delta_{cd}-\delta_{ac}\delta_{bd}$.

\item In triple index notation, $d_{abp,cdq}=\delta_{ab}\delta_{cd}-\delta_{ac}\delta_{bd}\delta_{pq}$.
\end{enumerate}
\end{proposition}

\begin{proof}
According to our single index conventions, from Definition 13.3, the adjacency matrix of the simplex is the one in the statement, namely:
\begin{eqnarray*}
d_{ij}
&=&(NP_1-1_N)_{ij}\\
&=&\bar{1}_i1_j-\delta_{ij}\\
&=&\delta_{i\bar{i}}\delta_{j\bar{j}}-\delta_{ij}
\end{eqnarray*}

In double index notation now, with $i=(ab)$ and $j=(cd)$, and $a,b,c,d$ being usual matrix indices, each thought to be attached to the corresponding matrix block of $C(Z)$, the formula that we obtain in the second one in the statement, namely:
\begin{eqnarray*}
d_{ab,cd}
&=&\delta_{ab,ba}\delta_{cd,dc}-\delta_{ab,cd}\\
&=&\delta_{ab}\delta_{cd}-\delta_{ac}\delta_{bd}
\end{eqnarray*}

Finally, in standard triple index notation, $i=(abp)$ and $j=(cdq)$, with $a,b,c,d$ being now usual numeric matrix indices, ranging in $1,2,3,\ldots\,$, and with $p,q$ standing for corresponding blocks of the algebra $C(Z)$, the formula that we obtain is:
\begin{eqnarray*}
d_{abp,cdq}
&=&\delta_{abp,bap}\delta_{cdq,dcq}-\delta_{abp,cdq}\\
&=&\delta_{ab}\delta_{cd}-\delta_{ac}\delta_{bd}\delta_{pq}
\end{eqnarray*}

Thus, we are led to the conclusions in the statement.
\end{proof}

At the level of examples, for $Z=\{1,\ldots,N\}$ the best is to use the above formula (1). The involution on the index set is $\bar{i}=i$, and we obtain, as we should:
$$d_{ij}=1-\delta_{ij}$$

As a more interesting example now, for the quantum space $Z=M_n$, coming by definition via the formula $C(Z)=M_n(\mathbb C)$, the situation is as follows:

\begin{proposition}
The simplex $Z_{full}$ with $Z=M_n$ is as follows: 
\begin{enumerate}
\item The vertices are $n^2$ points in the plane, arranged in square form.

\item Usual edges, worth $1$, are drawn between distinct points on the diagonal.

\item In addition, each off-diagonal point comes with a self-edge, worth $-1$.
\end{enumerate}
\end{proposition}

\begin{proof}
Here the most convenient is to use the double index formula from Proposition 13.17 (2), which tells us that $d$ is as follows, with indices $a,b,c,d\in\{1,\ldots,n\}$:
$$d_{ab,cd}=\delta_{ab}\delta_{cd}-\delta_{ac}\delta_{bd}$$

This quantity can be $-1,0,1$, and the study goes as follows:

\medskip

-- Case $d_{ab,cd}=1$. This can only happen when $\delta_{ab}\delta_{cd}=1$ and $\delta_{ac}\delta_{bd}=0$, corresponding to a formula of type $d_{aa,cc}=0$, with $a\neq c$, and so to the edges in (2).

\medskip

-- Case $d_{ab,cd}=-1$. This can only happen when $\delta_{ab}\delta_{cd}=0$ and $\delta_{ac}\delta_{bd}=1$, corresponding to a formula of type $d_{ab,ab}=0$, with $a\neq b$, and so to the self-edges in (3).
\end{proof}

The above result is quite interesting, and as an illustration, here is the pictorial representation of the simplex $Z_{full}$ on the vertex set $Z=M_3$, with the convention that the solid arrows are worth $-1$, and the dashed arrows are worth 1:
$$
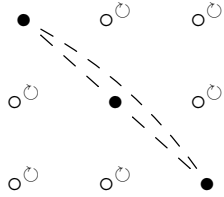
\xymatrix@R=15pt@C=15pt
{\bullet\ar@{--}[dr]\ar@/^/@{--}[ddrr]&\circ^\circlearrowright&\circ^\circlearrowright\\
\circ^\circlearrowright&\bullet\ar@{--}[dr]&\circ^\circlearrowright\\
\circ^\circlearrowright&\circ^\circlearrowright&\bullet}$$

More generally, we can in fact compute $Z_{full}$ for any finite quantum space $Z$, with the result here, which will be our final saying on the subject, being as follows:

\begin{theorem}
Consider a finite quantum space $Z$, and write it as follows, according to the decomposition formula $C(Z)=M_{n_1}(\mathbb C)\oplus\ldots\oplus M_{n_k}(\mathbb C)$ for its function algebra:
$$Z=M_{n_1}\sqcup\ldots\sqcup M_{n_k}$$
The simplex $Z_{full}$ is then the classical simplex formed by the points lying on the diagonals of $M_{n_1},\ldots,M_{n_k}$, with self-edges added, each worth $-1$, at the non-diagonal points.
\end{theorem}

\begin{proof}
The study here is quite similar to the one from the proof of Proposition 13.18, but by using this time the triple index formula from Proposition 13.17 (3), namely:
$$d_{abp,cdq}=\delta_{ab}\delta_{cd}-\delta_{ac}\delta_{bd}\delta_{pq}$$

Indeed, this quantity can be $-1,0,1$, and the $1$ case appears precisely as follows, leading to the classical simplex mentioned in the statement:
$$d_{aap,ccq}=1\quad,\quad\forall ap\neq cq$$

As for the remaining $-1$ case, this appears precisely as follows, leading this time to the self-edges worth $-1$, also mentioned in the statement:
$$d_{abp,abp}=1\quad,\quad\forall a\neq b$$

Thus, we are led to the conclusion in the statement.
\end{proof}

As an illustration, here is the simplex on the vertex set $Z=M_3\sqcup M_2$, with again the convention that the solid arrows are worth $-1$, and the dashed arrows are worth 1:
$$
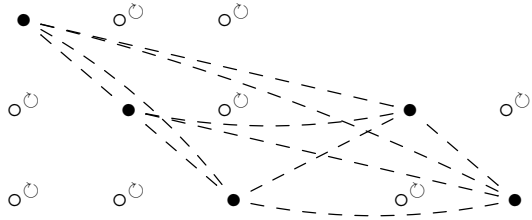
\xymatrix@R=18pt@C=20pt
{\bullet\ar@{--}[dr]\ar@/^/@{--}[ddrr]
\ar@{--}[drrrr]\ar@/^/@{--}[ddrrrrr]
&\circ^\circlearrowright&\circ^\circlearrowright\\
\circ^\circlearrowright&\bullet\ar@{--}[dr]\ar@/_/@{--}[rrr]\ar@{--}[drrrr]
&\circ^\circlearrowright&&\bullet\ar@{--}[dr]&\circ^\circlearrowright\\
\circ^\circlearrowright&\circ^\circlearrowright&\bullet\ar@{--}[urr]\ar@/_/@{--}[rrr]
&&\circ^\circlearrowright&\bullet}$$

\medskip

Long story short, we know what the simplex $Z_{full}$ is, and we have the formula $G^+(Z_{empty})=G^+(Z_{full})=S_Z^+$, exactly as in the $Z=\{1,\ldots,N\}$ case. Now with the above results in hand, we can talk as well about complementation, as follows:

\begin{theorem}
For any finite quantum graph $X$ we have the formula
$$G^+(X)=G^+(X^c)$$
where $X\to X^c$ is the complementation operation, given by $d_X+d_{X^c}=d_{Z_{full}}$.
\end{theorem}

\begin{proof}
This follows from Theorem 13.16, and more specifically from the following commutation relation, which is automatic, as explained there:
$$[u,d_{Z_{full}}]=0$$

Let us mention too that, in what concerns the pictorial representation of $X^c$, this can be deduced from what we have Theorem 13.19, in the obvious way.
\end{proof}

\section*{13c. Color independence} 

Following now \cite{ba2}, let us discuss an important point, namely the ``independence on the colors'' question. The idea indeed is that given a classical graph $X$ with edges colored by complex numbers, or other types of colors, $G(X)$ does not change when changing the colors. This is obvious, and a quantum analogue of this fact, involving $G^+(X)$, can be shown to hold as well, as explained in \cite{ba2}, and in chapter 6. 

\bigskip

In the quantum graph setting things are more complicated, as we will soon discover. Let us start with the following technical definition:

\begin{definition}
We call a quantum graph $X=(Z,d)$ washable if, with
$$d=\sum_ccd_c$$
being the color decomposition of $d$, we have the equivalence
$$[u,d]=0\iff[u,d_c]=0,\forall c$$
valid for any magic unitary matrix $u$, having size $|Z|$.
\end{definition}

Obviously, this is something which is not very beautiful, but the point is that some quantum graphs are washable, and some other are not, and so we have to deal with the above definition, as stated. As a first observation, we have the following result:

\begin{proposition}
Assuming that $X=(Z,d)$ is washable, its quantum symmetry group $G^+(X)$ does not depend on the precise colors of $X$. That is, whenever we have another quantum graph $X'=(Z,d')$ with same color scheme, in the sense that 
$$d_{ij}=d_{kl}\iff d'_{ij}=d'_{kl}$$
we have $G^+(X)=G^+(X')$.
\end{proposition}

\begin{proof}
This is something which is clear from the definition of $G^+(X)$, namely:
$$C(G^+(X))=C(S_Z^+)\Big/\Big<[u,d]=0\Big>$$

Indeed, assuming that our graph is washable in the above sense, we have:
$$C(G^+(X))=C(S_Z^+)\Big/\Big<[u,d_c]=0,\forall c\Big>$$

Thus, we are led to the conclusion in the statement.
\end{proof}

As already mentioned, it was proved in \cite{ba2} that in the classical case, $Z=\{1,\ldots,N\}$, all graphs are washable. This is a key result, and this for several reasons: 

\medskip

(1) First, this gives some intuition on what is going on with respect to colors, in analogy with what happens for $G(X)$. Also, it allows the use of true colors, like black, blue, red and so on, when drawing colored graphs, instead of complex numbers.

\medskip

(2) Also, this can be combined with the fact that $G^+(X)$ is invariant as well via some similar changes in the spectral decomposition of $d$, at the level of eigenvalues, with all this leading to some powerful combinatorial methods for the computation of $G^+(X)$.

\medskip

All these things do not necessarily hold in general, and to start with, we have:

\begin{theorem}
There are quantum graphs, such as the simplex in the homogeneous quantum space case, where
$$Z=M_K\times\{1,\ldots,L\}$$
with $K,L\geq2$, which are not washable.
\end{theorem}

\begin{proof}
We know that the simplex, in the case $Z=M_K\times\{1,\ldots,L\}$, has as adjacency matrix a certain matrix $d\in M_N(-1,0,1)$, with $N=K^2L$. Moreover, assuming $K,L\geq2$ as above, entries of all types, $-1,0,1$, are possible. Thus, the color decomposition of the adjacency matrix is as follows, with all 3 components being nonzero:
$$d=-1\cdot d_{-1}+0\cdot d_0+1\cdot d_1$$

Now assume that our simplex $X=Z_{full}$ is washable, and let $u$ be the fundamental corepresentation of $G^+(X)$. We have then the following commutation relations:
$$d_{-1},d_0,d_1\in End(u)$$

Now since the matrices $d_{-1},d_0,d_1$ are all nonzero, we deduce from this that:
$$\dim(End(u))\geq3$$

On the other hand, we know that we have $G^+(X)=S_Z^+$. Also, we know from the above that the Tannakian category of $S_Z^+$ is the Temperley-Lieb category $TL_N$, with the index being $N=K^2L$ as above. By putting these two results together, we obtain:
$$\dim(End(u))=\dim\Big(span\Big(\ |\ |\ ,{\ }^\cup_\cap\ \Big)\,\Big)\leq2$$

Thus, we have a contradiction, and so our simplex is not washable, as claimed.
\end{proof}

In order to come up with some positive results as well, the general idea will be that of using the method in \cite{ba2}. We have the following statement, coming from there:

\begin{theorem}
The following matrix belongs to $End(u)$, for any $n\in\mathbb N$:
$$d^{\times n}_{ij}=\sum_{i=k_1\ldots k_n}\sum_{j=l_1\ldots l_n}d_{k_1l_1}\ldots d_{k_nl_n}$$
In particular, in the classical case, $Z=\{1,\ldots,N\}$, all graphs are washable.
\end{theorem}

\begin{proof}
We have two assertions here, the idea being as follows:

\medskip

(1) Consider the multiplication and comultiplication maps of the algebra $C(Z)$, which in single index notation are given by:
$$\mu(e_i\otimes e_j)=e_{ij}\quad,\quad 
\gamma(e_i)=\sum_{i=jk}e_j\otimes e_k$$

Observe that we have $\mu^*=\gamma$, with the adjoint taken with respect to the scalar product coming from the canonical trace. We conclude that we have:
$$\mu\in Hom(u^{\otimes 2},u)\quad,\quad 
\gamma\in Hom(u,u^{\otimes 2})$$

The point now is that we can consider the iterations $\mu^{(n)},\gamma^{(n)}$ of $\mu,\gamma$, constructed in the obvious way, and we have then, for any $n\in\mathbb N$:
$$\mu^{(n)}\in Hom(u^{\otimes n},u)\quad,\quad 
\gamma^{(n)}\in Hom(u,u^{\otimes n})$$

Now if we assume that we have $d\in End(u)$, we have $d^{\otimes n}\in End(u^{\otimes n})$ for any $n\in\mathbb N$, and we conclude that we have the following formula:
$$\mu^{(n)}d^{\otimes n}\gamma^{(n)}\in End(u)$$

But, in single index notation, we have the following formula:
$$(\mu^{(n)}d^{\otimes n}\gamma^{(n)})_{ij}=\sum_{i=k_1\ldots k_n}\sum_{j=l_1\ldots l_n}d_{k_1l_1}\ldots d_{k_nl_n}$$

Thus, we are led to the conclusion in the statement.

\medskip

(2) Assuming that we are in the case $Z=\{1,\ldots,N\}$, the matrix $d^{\times n}$ in the statement is simply the componentwise $n$-th power of $d$, given by:
$$d^{\times n}_{ij}=d_{ij}^n$$ 

As explained in \cite{ba2} or in chapter 6, a simple analytic argument based on this, using $n\to\infty$ and then a recurrence on the number of colors, shows that we have washability indeed. Thus, we are led to the conclusions in the statement.
\end{proof}

In order now to further exploit the findings in Theorem 13.24, an idea would be that of assuming that we are in the homogeneous case, $Z=M_K\times\{1,\ldots,L\}$, and that the adjacency matrix is split, in the sense that one of the following happens:
$$d_{ab,cd}=e_{ab}f_{cd}\quad,\quad
d_{ab,cd}=e_{ac}f_{bd}\quad,\quad 
d_{ab,cd}=e_{ad}f_{bc}$$

Normally the graph should be washable in this case, but the computations are quite complex, and there is no clear result known in this sense. Thus, as a conclusion to all this, the basic theory of the quantum groups $G^+(X)$ from chapter 6 extends well to the present quantum graph setting, modulo some subtleties in connection with the colors. 

\section*{13d. Twisted reflections}

With the above technology in hand, we can talk about twisted quantum reflections. The idea will be that the twisted analogues of the quantum reflection groups $H_N^{s+}\subset S_{sN}^+$ will be the quantum automorphism groups $S_{Z\to Y}^+$ of the fibrations of finite quantum spaces $Z\to Y$, which correspond by definition to the Markov inclusions of finite dimensional $C^*$-algebras $C(Y)\subset C(Z)$. In order to discuss this, let us start with:

\begin{definition}
A fibration of finite quantum spaces $Z\to Y$ corresponds to an inclusion of finite dimensional $C^*$-algebras
$$C(Y)\subset C(Z)$$
which is Markov, in the sense that it commutes with the canonical traces.
\end{definition}  

Here the commutation condition with the canonical traces means that the composition $C(Y)\subset C(Z)\to\mathbb C$ should equal the canonical trace $C(Y)\to\mathbb C$. At the level of the corresponding quantum spaces, this means that the quotient map $Z\to Y$ must commute with the corresponding counting measures, and this is where our term ``fibration'' comes from. In order to talk now about the quantum symmetry groups $S_{Z\to Y}^+$, we will need:

\begin{proposition}
Given a fibration $Z\to Y$, a closed subgroup $G\subset S_Z^+$ leaves invariant $Y$ precisely when its magic unitary $u=(u_{ij})$ satisfies the condition
$$e\in End(u)$$
where $e:C(Z)\to C(Z)$ is the Jones projection, onto the subalgebra $C(Y)\subset C(Z)$.
\end{proposition}

\begin{proof}
This is something that we know well, in the commutative case, where $Z$ is a usual finite set, and the proof in general is similar.
\end{proof}

We can now talk about twisted quantum reflection groups, as follows:

\begin{theorem}
Any fibration of finite quantum spaces $Z\to Y$ has a quantum symmetry group, which is the biggest acting on $Z$ by leaving $Y$ invariant:
$$S_{Z\to Y}^+\subset S_Z^+$$
At the level of algebras of functions, this quantum group $S_{Z\to Y}^+$ is obtained as follows, with $e:C(Z)\to C(Y)$ being the Jones projection:
$$C(S_{Z\to Y}^+)=C(S_Z^+)\Big/\Big<e\in End(u)\Big>$$
We call these quantum groups $S_{Z\to Y}^+$ twisted quantum reflection groups. 
\end{theorem}

\begin{proof}
This follows indeed from Proposition 13.26.
\end{proof}

As a basic example, let us discuss the commutative case. Here we have:

\begin{proposition}
In the commutative case, the fibration $Z\to Y$ must be of the following special form, with $N,s$ being certain integers,
$$\{1,\ldots,N\}\times\{1,\ldots,s\}\to\{1,\ldots,N\}\quad,\quad (i,a)\to i$$
and we obtain the usual quantum reflection groups,
$$(S_{Z\to Y}^+\subset S_Z^+)\ = \ (H_N^{s+}\subset S_{sN}^+)$$
via some standard identifications.
\end{proposition}

\begin{proof}
In the commutative case our fibration must be a usual fibration of finite spaces, $\{1,\ldots,M\}\to\{1,\ldots,N\}$, commuting with the counting measures. But this shows that our fibration must be of the following special form, with $N,s\in\mathbb N$:
$$\{1,\ldots,N\}\times\{1,\ldots,s\}\to\{1,\ldots,N\}\quad,\quad (i,a)\to i$$

Regarding now the quantum symmetry group, we have the following formula for it, with $e:\mathbb C^N\otimes\mathbb C^s\to\mathbb C^N$ being the Jones projection for the inclusion $\mathbb C^N\subset\mathbb C^N\otimes\mathbb C^s$:
$$C(S_{Z\to Y}^+)=C(S_{sN}^+)\Big/\Big<e\in End(u)\Big>$$

On the other hand, recall that the quantum reflection group $H_N^{s+}\subset S_{sN}^+$ appears via the condition that the corresponding magic matrix must be sudoku:
$$u=\begin{pmatrix}
a^0&a^1&\ldots&a^{s-1}\\
a^{s-1}&a^0&\ldots&a^{s-2}\\
\vdots&\vdots&&\vdots\\
a^1&a^2&\ldots&a^0
\end{pmatrix}$$

But, as explained in \cite{ba2}, this is the same as saying that the quantum group $H_N^{s+}\subset S_{sN}^+$ appears as the symmetry group of the multi-simplex associated to the fibration $\{1,\ldots,N\}\times\{1,\ldots,s\}\to\{1,\ldots,N\}$, so we have an identification as follows:
$$(S_{Z\to Y}^+\subset S_Z^+)\ = \ (H_N^{s+}\subset S_{sN}^+)$$

Thus, we are led to the conclusions in the statement.
\end{proof}

Observe that in Proposition 13.28 the fibration $Z\to Y$ is ``trivial'', in the sense that it is of the following special form:
$$Y\times T\to Y\quad,\quad (i,a)\to i$$

However, in the general quantum case, there are many interesting fibrations $Z\to Y$ which are not trivial, and in what follows we will not make any assumption on our fibrations, and use Definition 13.26 and Theorem 13.27 as stated.

\bigskip

Following \cite{ba2}, we will prove now that the Tannakian category of $S_{Z\to Y}^+$, which is by definition a generalization of $S_Z^+$, is the Fuss-Catalan category, which is a generalization of the Temperley-Lieb category, introduced by Bisch and Jones. Let us start with: 

\begin{theorem}
Any Markov inclusion of finite dimensional algebras $D\subset B$ has a quantum symmetry group $S_{D\subset B}^+$. The corresponding Woronowicz algebra is generated by the coefficients of a biunitary matrix $v=(v_{ij})$ subject to the conditions
$$m\in Hom(v^{\otimes 2},v)\quad,\quad 
u\in Hom(1,v)\quad,\quad 
e\in End(v)$$
where $m:B\otimes B\to B$ is the multiplication, $u:\mathbb C\to B$ is the unit and $e:B\to B$ is the projection onto $D$, with respect to the scalar product $<x,y>=tr(xy^*)$.
\end{theorem}

\begin{proof}
This is a reformulation of Theorem 13.27, with several modifications made. Indeed, by using the algebras $D=C(Y)$, $B=C(Z)$ instead of the quantum spaces $Y,Z$ used there, and also by calling the fundamental corepresentation $v=(v_{ij})$, in order to avoid confusion with the unit $u:\mathbb C\to B$, the formula in Theorem 13.27 reads:
$$C(S_{D\subset B}^+)=C(S_B^+)\Big/\Big<e\in End(v)\Big>$$

Also, we know from Theorem 13.5 that we have the following formula, again by using $B$ instead of $F$, and by calling the fundamental corepresentation $v=(v_{ij})$:
$$C(S_B^+)=C(U_N^+)\Big/\Big<m\in Hom(v^{\otimes2},v),u\in Fix(v)\Big>$$

Thus, we are led to the conclusion in the statement.
\end{proof}

Let us first discuss in detail the Temperley-Lieb algebra, as a continuation of the material above. In the present context, we have the following definition:

\begin{definition}
The $\mathbb N$-algebra $TL^2$ of index $\delta >0$ is defined as follows:
\begin{enumerate}
\item The space $TL^2(m,n)$ consists of linear combinations
of noncrossing pairings between $2m$ points and $2n$ points:
$${{TL^2}}(m,n)=
\left\{\sum\,\alpha\,\,
\begin{matrix}\cdots\cdots&\leftarrow&2m\,\,\ {\rm points}\\
\mathfrak{W}&\leftarrow&m+n\ \ {\rm strings}\\
\cdot\,\cdot&\leftarrow&2n\,\,{\rm points}
\end{matrix}\right\}$$

\item The operations $\circ$, $\otimes$, $*$ are induced by the vertical and horizontal concatenation and the upside-down turning of diagrams:
$$A\circ B=\binom{B}{A}\quad,\quad A\otimes B=AB\quad,\quad A^*=\forall$$

\item With the rule $\bigcirc=\delta$, erasing a circle is the same as multiplying by $\delta$.
\end{enumerate}
\end{definition}

Our first task will be that of finding a suitable presentation for this algebra. Consider the following two elements $u\in TL^2(0,1)$ and $m\in TL^2(2,1)$:
$$u=\delta^{-\frac{1}{2}}\,\cap
\quad,\quad
m=\delta^{\frac{1}{2}}\,|\cup|$$ 

With this convention, we have the following result:

\begin{theorem}
The following relations are a presentation of $TL^2$ by the above rescaled diagrams $u\in TL^2(0,1)$ and $m\in TL^2(2,1)$:
\begin{enumerate}
\item $mm^*=\delta^2$.

\item $u^*u=1$.

\item $m(m\otimes1)=m(1\otimes m)$.

\item $m(1\otimes u)=m(u\otimes1)=1$.

\item $(m\otimes1)(1\otimes m^*)=(1\otimes m)(m^*\otimes1)=m^*m$.
\end{enumerate}
\end{theorem}

\begin{proof}
This is something well-known, and elementary, obtained by drawing diagrams, and for details here, we refer to Bisch and Jones \cite{bjo}.
\end{proof}

In more concrete terms, the above result says that $u,m$ satisfy the above relations, which is something clear, and that if $C$ is a $\mathbb N$-algebra and $v\in C(0,1)$ and $n\in C(2,1)$ satisfy the same relations then there exists a $\mathbb N$-algebra morphism as follows:
$$TL^2\to C\quad,\quad 
u\to v\quad,\quad 
m\to n$$

Now let $B$ be a finite dimensional $C^*$-algebra, with its canonical trace. We have a scalar product $<x,y>=tr(xy^*)$ on $B$, so $B$ is an object in the category of finite dimensional Hilbert spaces. Consider the unit $u$ and the multiplication $m$ of $B$:
$$u\in\mathbb NB(0,1)\quad,\quad 
m\in\mathbb NB(2,1)$$

The relations in Theorem 13.31 are then satisfied, and one can deduce from this that in this case, the category of representations of $S_B^+$ is the completion of $TL^2$, as we already know. Getting now to Fuss-Catalan algebras, we have here:

\begin{definition}
A Fuss-Catalan diagram is a planar diagram formed by an upper row of $4m$ points, a lower row of $4n$ points, both colored
$$\circ\bullet\bullet\circ\circ\bullet\bullet\ldots$$
and by $2m+2n$ noncrossing strings joining these $4m+4n$ points, with the rule that the points which are joined must have the same color.
\end{definition}

Fix $\beta>0$ and $\omega>0$. The $\mathbb N$-algebra $FC$ is defined as follows. The spaces $FC(m,n)$ consist of linear combinations of Fuss-Catalan diagrams:
$$FC(m,n)=\left\{\sum\,\alpha\,\,
\begin{matrix}\circ\bullet\bullet\circ\circ\bullet\bullet\circ\ldots\ldots&\leftarrow&4m\,\,\ {\rm colored\ points}\\
\ &\ &m+n {\rm\ black\ strings}\\
\mathfrak{W}\mbox{\ \ }&\leftarrow&{\rm and}\\
\ & \ & m+n{\rm\ white\ strings}\\
\circ\bullet\bullet\circ\circ\bullet\bullet\circ\ldots\ldots&\leftarrow&4n\,\,{\rm colored\ points}\end{matrix}\right\}$$

As before with the Temperley-Lieb algebra, the operations $\circ$, $\otimes$, $*$ are induced by vertical and horizontal concatenation and upside-down turning of diagrams, but this time with the rule that erasing a black/white circle is the same as multiplying by $\beta$/$\omega$:
$$A\circ B=\binom{B}{A}\quad,\quad 
A\otimes B=AB\quad,\quad 
A^*=\forall$$
$$\mbox{\tiny{black}}\rightarrow\bigcirc=\beta
\qquad,\qquad
\mbox{\tiny{white}}\rightarrow\bigcirc=\omega$$

Let $\delta=\beta\omega$. We have the following bicolored analogues of the elements $u,m$:
$$u=\delta^{-\frac{1}{2}}\,\,\bigcap\;\!\!\!\!\!\!\!\!\cap
\quad,\quad
m=\delta^{\frac{1}{2}}\,\,||\bigcup\;\!\!\!\!\!\!\!\!\cup\ ||$$

Consider also the black and white Jones projections, namely:
$$e=\omega^{-1}\,\,|\ \begin{matrix}\cup\cr\cap\end{matrix}\ |
\quad,\quad 
f=\beta^{-1}\,\,|||\ \begin{matrix}\cup\cr\cap\end{matrix}\ |||$$

For simplifying writing we identify $x$ and $x\otimes 1$. We have the following result:

\begin{theorem}
The following relations, with $f=\beta^{-2}(1\otimes me)m^*$, are a presentation of $FC$ by $m\in FC(2,1)$, $u\in FC(0,1)$ and $e\in FC(1)$:
\begin{enumerate}
\item The relations in Theorem 13.31, with $\delta =\beta\omega$.

\item $e=e^2=e^*$, $f=f^*$ and $(1\otimes f)f=f(1\otimes f)$.

\item $eu=u$.

\item $mem^*=m(1\otimes e)m^*=\beta^2$.

\item $mm(e\otimes e\otimes e)=emm(e\otimes1\otimes e)$.
\end{enumerate}
\end{theorem}

\begin{proof}
As for any presentation result, we have to prove two assertions:

\medskip

(I) The elements $m,u,e$ satisfy the relations (1-5) and generate the
$\mathbb N$-algebra $FC$.

\medskip

(II) If  $M$, $U$ and $E$ in a $\mathbb N$-algebra $C$ satisfy the relations
(1-5), then there exists a morphism of $\mathbb N$-algebras $FC\to C$
sending $m\to M$, $u\to U$, $e\to E$.

\medskip

(I) First, the relations (1-5) are easily verified by drawing pictures. Let us show now that the $\mathbb N$-subalgebra $C=<m,u,e>$ of $FC$ is equal to $FC$. First, $C$ contains the infinite sequence of black and white Jones projections:
$$p_1=e=\omega^{-1}\,\,|\ \begin{matrix}\cup\\\cap\end{matrix}\ |$$
$$p_2=f=\beta^{-1}\,\,|||\ \begin{matrix}\cup\\\cap\end{matrix}\ |||$$
$$p_3=1\otimes e=\omega^{-1}\,\,|||||\ \begin{matrix}\cup\\\cap\end{matrix}\ |$$
$$p_4=1\otimes f=\beta^{-1}\,\,|||||||\ \begin{matrix}\cup\cr\cap\end{matrix}\ |||$$
$$\ldots$$

The algebra $C$ contains as well the infinite sequence of bicolored Jones projections:
$$e_1=uu^*=\delta^{-1}\,\,\begin{matrix}\bigcup\hskip -3.9mm\cup\\ \bigcap\hskip -3.9mm\cap\end{matrix}$$
$$e_2=\delta^{-2}m^*m=\delta^{-1}\,\,||\ \begin{matrix}\bigcup\hskip -3.86mm\cup\\\bigcap\hskip -3.86mm\cap\end{matrix}\ ||$$
$$e_3=1\otimes uu^*=\delta^{-1}\,\,||||\ \begin{matrix}\bigcup\hskip -3.86mm\cup\\ \bigcap\hskip -3.86mm\cap\end{matrix}$$
$$e_4=\delta^{-2}(1\otimes m^*m)=\delta^{-1}\,\,||||||\ \begin{matrix}\bigcup\hskip -3.86mm\cup\\\bigcap\hskip -3.86mm\cap\end{matrix}\ ||$$
$$\ldots$$

By standard results of Bisch and Jones, these latter projections generate the diagonal $\mathbb N$-algebra $\Delta FC$. But this gives the result, via a standard study.
\end{proof}

Getting back now to the inclusions $D\subset B$, we have the following result:

\begin{theorem}
Given a Markov inclusion $D\subset B$, we have
$$<m,u,e>=FC$$
as an equality of $\mathbb N$-algebras.
\end{theorem}

\begin{proof}
It is routine to check that the linear maps $m,u,e$ associated to an inclusion $D\subset B$ as in the statement satisfy the relations (1-5) in Theorem 13.33. Thus, we obtain a certain $\mathbb N$-algebra surjective morphism, as follows:
$$J:FC\to<m,u,e>$$

But it is routine to prove that this morphism $J$ is faithful on $\Delta FC$, and then by Frobenius reciprocity faithfulness has to hold on the whole $FC$.
\end{proof}

Getting back now to quantum groups, we have:

\begin{theorem}
Given a Markov inclusion of finite dimensional algebras $D\subset B$, the category of representations of its quantum symmetry group 
$$S_{D\subset B}^+\subset S_B^+$$
is the completion of the Fuss-Catalan category $FC$.
\end{theorem}

\begin{proof}
Since $S_{D\subset B}^+$ comes by definition from the relations corresponding to $m,u,e$, its tensor category of corepresentations is the completion of the tensor category $<m,u,e>$. Thus Theorem 13.34 applies, and gives an isomorphism $<m,u,e>\simeq FC$.
\end{proof}

In terms of finite quantum spaces and quantum graphs, the conclusion is that the quantum automorphism groups $S_{Z\to Y}^+$ of the Markov fibrations $Z\to Y$, which can be thought of as being the ``twisted versions'' of the quantum reflection groups $H_N^{s+}$, correspond to the Fuss-Catalan algebras. We refer here to \cite{ba2} and related papers.

\section*{13e. Exercises}

This was a quite elementary chapter, and as exercises on this, we have:

\begin{exercise}
Come up with further speculations on the finite quantum spaces.
\end{exercise}

\begin{exercise}
Further study the homogeneous case, $Z=M_K\times\{1,\ldots,L\}$.
\end{exercise}

\begin{exercise}
Find some further results on the color independence question.
\end{exercise}

\begin{exercise}
Learn about the Frucht theorem, and its quantum extensions.
\end{exercise}

\begin{exercise}
Fill in all the technical details at our main Fuss-Catalan proof.
\end{exercise}

\begin{exercise}
Clarify the axiomatization of general quantum reflection groups.
\end{exercise}

As bonus exercise, learn more graph theory, and generalize that, to our setting.

\chapter{Quantum spheres}

\section*{14a. Spheres and tori}

We have seen in the previous chapter how discrete quantum geometry looks like, and many other things can be said, along that lines. For instance, there are many further things that can be said, in relation with the subfactor considerations in chapter 12. But, time now to get into the real thing, namely continuous quantum geometry. 

\bigskip

As a first observation, we are already familiar with continuous quantum geometry since Part III, because our considerations there regarding $O_N,U_N,O_N^+,U_N^+$ and related quantum groups, such as $B_N,C_N,B_N^+,C_N^+$ are, obviously, of ``continuous'' nature. However, such things remain quite specialized, and we will be here rather interested in foundations.

\bigskip

So, foundations. To start with, we can talk about quantum tori, as follows:

\begin{definition}
The classical and quantum, real and complex tori are
$$\xymatrix@R=17mm@C=19mm{
T_N^+\ar[r]&\mathbb T_N^+\\
T_N\ar[r]\ar[u]&\mathbb T_N\ar[u]
}$$
with these standing by definition for the abstract duals of the groups
$$\xymatrix@R=16mm@C=18mm{
\mathbb Z_2^{*N}\ar[d]&\mathbb Z^{*N}\ar[d]\ar[l]\\
\mathbb Z_2^N&\mathbb Z^N\ar[l]
}$$
in the sense of the abstract duality between compact and discrete quantum groups.
\end{definition}

In other words, we have indeed quantum tori, as a corollary, or rather as the simplest particular case, of the quantum group theory developed in Part III. Getting now to noncommutative geometry, over these tori, things here rather belong to algebra and probability, instead of differential geometry, as one might naively expect, and this because $\mathbb T_N^+$, which is the dual of the free group $F_N=\mathbb Z^{*N}$, has of course the same mathematics as this free group $F_N$, which is notoriously a business of algebra and analysis.

\bigskip 

Nevermind. So, as a first conclusion, the noncommutative geometry that we would like to develop will be rather algebraic geometry, with a touch of probability, and with this being not an issue, as kindly explained by our quantum physics collaborator:

\begin{cat}
No wonder true noncommutative geometry is about algebra and probability, quantum mechanics being the same, something about algebra and probability.
\end{cat}

With this discussed, instead of further insisting on the tori, whose mathematics is widely well-known, and is also one level below what we did in Part III, with far more general compact quantum groups, let us get now into the real thing, quantum spheres.

\bigskip

According to our findings above, and to Cat 14.2 too, we will take an algebraic geometry approach, with a touch of probability, to the quantum spheres. With this idea in mind, and following \cite{bgo}, we are led to the following notions:

\begin{definition}
We have free real and complex spheres, defined via
$$C(S^{N-1}_{\mathbb R,+})=C^*\left(x_1,\ldots,x_N\Big|x_i=x_i^*,\sum_ix_i^2=1\right)$$
$$C(S^{N-1}_{\mathbb C,+})=C^*\left(x_1,\ldots,x_N\Big|\sum_ix_ix_i^*=\sum_ix_i^*x_i=1\right)$$
where the symbol $C^*$ stands for universal enveloping $C^*$-algebra.
\end{definition}

Here the fact that these algebras are indeed well-defined comes from the following estimate, which shows that the biggest $C^*$-norms on these $*$-algebras are bounded:
$$||x_i||^2
=||x_ix_i^*||
\leq\left|\left|\sum_ix_ix_i^*\right|\right|
=1$$

As a first result now, regarding the above free spheres, we have:

\begin{theorem}
We have embeddings of compact quantum spaces, as follows,
$$\xymatrix@R=15mm@C=15mm{
S^{N-1}_{\mathbb R,+}\ar[r]&S^{N-1}_{\mathbb C,+}\\
S^{N-1}_\mathbb R\ar[r]\ar[u]&S^{N-1}_\mathbb C\ar[u]
}$$
and the spaces on top appear as liberations of the spaces on the bottom.
\end{theorem}

\begin{proof}
The first assertion, regarding the inclusions, comes from the fact that at the level of the associated $C^*$-algebras, we have surjective maps, as follows:
$$\xymatrix@R=15mm@C=15mm{
C(S^{N-1}_{\mathbb R,+})\ar[d]&C(S^{N-1}_{\mathbb C,+})\ar[d]\ar[l]\\
C(S^{N-1}_\mathbb R)&C(S^{N-1}_\mathbb C)\ar[l]
}$$

For the second assertion, we must establish the following isomorphisms, where the symbol $C^*_{comm}$ stands for ``universal commutative $C^*$-algebra generated by'':
$$C(S^{N-1}_\mathbb R)=C^*_{comm}\left(x_1,\ldots,x_N\Big|x_i=x_i^*,\sum_ix_i^2=1\right)$$
$$C(S^{N-1}_\mathbb C)=C^*_{comm}\left(x_1,\ldots,x_N\Big|\sum_ix_ix_i^*=\sum_ix_i^*x_i=1\right)$$

It is enough to establish the second isomorphism. So, consider the second universal commutative $C^*$-algebra $A$ constructed above. Since the standard coordinates on $S^{N-1}_\mathbb C$ satisfy the defining relations for $A$, we have a quotient map of as follows:
$$A\to C(S^{N-1}_\mathbb C)$$

Conversely, let us write $A=C(S)$, by using the Gelfand theorem. The variables $x_1,\ldots,x_N$ become in this way true coordinates, providing us with an embedding $S\subset\mathbb C^N$. Also, the quadratic relations become $\sum_i|x_i|^2=1$, so we have $S\subset S^{N-1}_\mathbb C$. Thus, we have a quotient map $C(S^{N-1}_\mathbb C)\to A$, as desired, and this gives all the results.
\end{proof}

By using the free spheres constructed above, we can now formulate:

\begin{definition}
A real algebraic manifold $X\subset S^{N-1}_{\mathbb C,+}$ is a closed quantum subspace defined, at the level of the corresponding $C^*$-algebra, by a formula of type
$$C(X)=C(S^{N-1}_{\mathbb C,+})\Big/\Big<f_i(x_1,\ldots,x_N)=0\Big>$$
for certain family of noncommutative polynomials, as follows:
$$f_i\in\mathbb C<x_1,\ldots,x_N>$$
We denote by $\mathcal C(X)$ the $*$-subalgebra of $C(X)$ generated by the coordinates $x_1,\ldots,x_N$. 
\end{definition}

As a basic example here, we have the free real sphere $S^{N-1}_{\mathbb R,+}$. The classical spheres $S^{N-1}_\mathbb C,S^{N-1}_\mathbb R$, and their real submanifolds, are covered as well by this formalism. At the level of the general theory, we have the following version of the Gelfand theorem:

\begin{theorem}
If $X\subset S^{N-1}_{\mathbb C,+}$ is an algebraic manifold, as above, we have
$$X_{class}=\left\{x\in S^{N-1}_\mathbb C\Big|f_i(x_1,\ldots,x_N)=0\right\}$$
and $X$ appears as a liberation of $X_{class}$.
\end{theorem}

\begin{proof}
This is something that we already met, in the context of the free spheres. In general, the proof is similar, by using the Gelfand theorem. Indeed, if we denote by $X_{class}'$ the manifold constructed in the statement, then we have a quotient map of $C^*$-algebras as follows, mapping standard coordinates to standard coordinates:
$$C(X_{class})\to C(X_{class}')$$

Conversely now, from $X\subset S^{N-1}_{\mathbb C,+}$ we obtain $X_{class}\subset S^{N-1}_\mathbb C$. Now since the relations defining $X_{class}'$ are satisfied by $X_{class}$, we obtain an inclusion $X_{class}\subset X_{class}'$. Thus, at the level of algebras of continuous functions, we have a quotient map of $C^*$-algebras as follows, mapping standard coordinates to standard coordinates:
$$C(X_{class}')\to C(X_{class})$$

Thus, we have constructed a pair of inverse morphisms, and we are done.
\end{proof}

Finally, once again at the level of the general theory, we have:

\begin{definition}
We agree to identify two real algebraic submanifolds $X,Y\subset S^{N-1}_{\mathbb C,+}$ when we have a $*$-algebra isomorphism between $*$-algebras of coordinates
$$f:\mathcal C(Y)\to\mathcal C(X)$$
mapping standard coordinates to standard coordinates.
\end{definition}

We will see later the reasons for making this convention, coming from amenability. Now back to the tori, as constructed before, we can see that these are examples of algebraic manifolds, in the sense of Definition 14.5. In fact, we have the following result:

\begin{theorem}
The four main quantum spheres produce the main quantum tori
$$\xymatrix@R=15mm@C=15mm{
S^{N-1}_{\mathbb R,+}\ar[r]&S^{N-1}_{\mathbb C,+}\\
S^{N-1}_\mathbb R\ar[r]\ar[u]&S^{N-1}_\mathbb C\ar[u]
}\qquad
\xymatrix@R=8mm@C=15mm{\\ \to}
\qquad
\xymatrix@R=16mm@C=18mm{
T_N^+\ar[r]&\mathbb T_N^+\\
T_N\ar[r]\ar[u]&\mathbb T_N\ar[u]
}$$
via the formula $T=S\cap\mathbb T_N^+$, with the intersection being taken inside $S^{N-1}_{\mathbb C,+}$.
\end{theorem}

\begin{proof}
This comes from the above results, the situation being as follows:

\medskip

(1) Free complex case. Here the formula in the statement reads $\mathbb T_N^+=S^{N-1}_{\mathbb C,+}\cap\mathbb T_N^+$. But this is something trivial, because we have $\mathbb T_N^+\subset S^{N-1}_{\mathbb C,+}$.

\medskip

(2) Free real case. Here the formula in the statement reads $T_N^+=S^{N-1}_{\mathbb R,+}\cap\mathbb T_N^+$. But this is clear as well, the real version of $\mathbb T_N^+$ being $T_N^+$.

\medskip

(3) Classical complex case. Here the formula in the statement reads $\mathbb T_N=S^{N-1}_\mathbb C\cap\mathbb T_N^+$. But this is clear as well, the classical version of $\mathbb T_N^+$ being $\mathbb T_N$.

\medskip

(4) Classical real case. Here the formula in the statement reads $T_N=S^{N-1}_\mathbb R\cap\mathbb T_N^+$. But this follows by intersecting the formulae from the proof of (2) and (3).
\end{proof}

In order to discuss now the relation between the quantum rotations and reflections, constructed and studied in Parts II-III, with the spheres and tori, we will need:

\begin{proposition}
Given an algebraic manifold $X\subset S^{N-1}_{\mathbb C,+}$, the category of closed subgroups $G\subset U_N^+$ acting affinely on $X$, in the sense that the formula
$$\alpha(x_i)=\sum_jx_j\otimes u_{ji}$$ 
defines a morphism of $C^*$-algebras $\alpha:C(X)\to C(X)\otimes C(G)$, has a universal object, denoted $G^+(X)$, and called affine quantum isometry group of $X$.
\end{proposition}

\begin{proof}
Assume indeed that our manifold $X\subset S^{N-1}_{\mathbb C,+}$ comes as follows:
$$C(X)=C(S^{N-1}_{\mathbb C,+})\Big/\Big<f_\iota(x_1,\ldots,x_N)=0\Big>$$

In order to prove the result, consider the following variables:
$$X_i=\sum_jx_j\otimes u_{ji}\in C(X)\otimes C(U_N^+)$$

Our claim is that the quantum group in the statement $G=G^+(X)$ appears as:
$$C(G)=C(U_N^+)\Big/\Big<f_\iota(X_1,\ldots,X_N)=0\Big>$$

In order to prove this, pick one of the defining polynomials, and write it as follows:
$$f_\iota(x_1,\ldots,x_N)=\sum_r\sum_{i_1^r\ldots i_{s_r}^r}\lambda_r\cdot x_{i_1^r}\ldots x_{i_{s_r}^r}$$

With $X_i=\sum_jx_j\otimes u_{ji}$ as above, we have the following formula:
$$f_\iota(X_1,\ldots,X_N)=\sum_r\sum_{i_1^r\ldots i_{s_r}^r}\lambda_r\sum_{j_1^r\ldots j_{s_r}^r}x_{j_1^r}\ldots x_{j_{s_r}^r}\otimes u_{j_1^ri_1^r}\ldots u_{j_{s_r}^ri_{s_r}^r}$$

Since the variables on the right span a certain finite dimensional space, the relations $f_\iota(X_1,\ldots,X_N)=0$ correspond to certain relations between the variables $u_{ij}$. Thus, we have indeed a closed subspace $G\subset U_N^+$, with a universal map, as follows:
$$\alpha:C(X)\to C(X)\otimes C(G)$$

In order to show now that $G$ is a quantum group, consider the following elements:
$$u_{ij}^\Delta=\sum_ku_{ik}\otimes u_{kj}\quad,\quad u_{ij}^\varepsilon=\delta_{ij}\quad,\quad u_{ij}^S=u_{ji}^*$$

Consider as well the following elements, with $\gamma\in\{\Delta,\varepsilon,S\}$:
$$X_i^\gamma=\sum_jx_j\otimes u_{ji}^\gamma$$

From the relations $f_\iota(X_1,\ldots,X_N)=0$ we deduce that we have:
$$f_\iota(X_1^\gamma,\ldots,X_N^\gamma)
=(id\otimes\gamma)f_\iota(X_1,\ldots,X_N)
=0$$

Thus we can map $u_{ij}\to u_{ij}^\gamma$ for any $\gamma\in\{\Delta,\varepsilon,S\}$, and we are done.
\end{proof}

We can now formulate a result about spheres and rotations, as follows:

\begin{theorem}
The quantum isometry groups of the basic spheres are
$$\xymatrix@R=15mm@C=14mm{
S^{N-1}_{\mathbb R,+}\ar[r]&S^{N-1}_{\mathbb C,+}\\
S^{N-1}_\mathbb R\ar[r]\ar[u]&S^{N-1}_\mathbb C\ar[u]
}
\qquad
\xymatrix@R=8mm@C=15mm{\\ \to}
\qquad
\xymatrix@R=16mm@C=18mm{
O_N^+\ar[r]&U_N^+\\
O_N\ar[r]\ar[u]&U_N\ar[u]}$$
modulo identifying, as usual, the various $C^*$-algebraic completions.
\end{theorem}

\begin{proof}
We have 4 results to be proved, the idea being as follows:

\medskip

\underline{$S^{N-1}_{\mathbb C,+}$}. Let us first construct an action $U_N^+\curvearrowright S^{N-1}_{\mathbb C,+}$. We must prove here that the variables $X_i=\sum_jx_j\otimes u_{ji}$ satisfy the defining relations for $S^{N-1}_{\mathbb C,+}$, namely:
$$\sum_ix_ix_i^*=\sum_ix_i^*x_i=1$$

By using the biunitarity of $u$, we have the following computation:
$$\sum_iX_iX_i^*
=\sum_{ijk}x_jx_k^*\otimes u_{ji}u_{ki}^*
=\sum_jx_jx_j^*\otimes1
=1\otimes1$$

Once again by using the biunitarity of $u$, we have as well:
$$\sum_iX_i^*X_i
=\sum_{ijk}x_j^*x_k\otimes u_{ji}^*u_{ki}
=\sum_jx_j^*x_j\otimes1
=1\otimes1$$

Thus we have an action $U_N^+\curvearrowright S^{N-1}_{\mathbb C,+}$, which gives $G^+(S^{N-1}_{\mathbb C,+})=U_N^+$, as desired. 

\medskip

\underline{$S^{N-1}_{\mathbb R,+}$}. Let us first construct an action $O_N^+\curvearrowright S^{N-1}_{\mathbb R,+}$. We already know that the variables $X_i=\sum_jx_j\otimes u_{ji}$ satisfy the defining relations for $S^{N-1}_{\mathbb C,+}$, so we just have to check that these variables are self-adjoint. But this is clear from $u=\bar{u}$, as follows:
$$X_i^*
=\sum_jx_j^*\otimes u_{ji}^*
=\sum_jx_j\otimes u_{ji}
=X_i$$

Conversely, assume that we have an action $G\curvearrowright S^{N-1}_{\mathbb R,+}$, with $G\subset U_N^+$. The variables $X_i=\sum_jx_j\otimes u_{ji}$ must be then self-adjoint, and the above computation shows that we must have $u=\bar{u}$. Thus our quantum group must satisfy $G\subset O_N^+$, as desired.

\medskip

\underline{$S^{N-1}_\mathbb C$}. The fact that we have an action $U_N\curvearrowright S^{N-1}_\mathbb C$ is clear. Conversely, assume that we have an action $G\curvearrowright S^{N-1}_\mathbb C$, with $G\subset U_N^+$. We must prove that this implies $G\subset U_N$, and we will use a standard trick of Bhowmick-Goswami \cite{bgo}. We have:
$$\alpha(x_i)=\sum_jx_j\otimes u_{ji}$$

By multiplying this formula with itself we obtain the following formulae:
$$\alpha(x_ix_k)=\sum_{jl}x_jx_l\otimes u_{ji}u_{lk}$$
$$\alpha(x_kx_i)=\sum_{jl}x_lx_j\otimes u_{lk}u_{ji}$$

Now since the variables $x_i$ commute, these formulae can be written as:
$$\alpha(x_ix_k)=\sum_{j<l}x_jx_l\otimes(u_{ji}u_{lk}+u_{li}u_{jk})+\sum_jx_j^2\otimes u_{ji}u_{jk}$$
$$\alpha(x_ix_k)=\sum_{j<l}x_jx_l\otimes(u_{lk}u_{ji}+u_{jk}u_{li})+\sum_jx_j^2\otimes u_{jk}u_{ji}$$

Since the tensors at left are linearly independent, we must have:
$$u_{ji}u_{lk}+u_{li}u_{jk}=u_{lk}u_{ji}+u_{jk}u_{li}$$

By applying now the antipode to this formula, then applying the involution, and then relabeling the indices, we succesively obtain:
$$u_{kl}^*u_{ij}^*+u_{kj}^*u_{il}^*=u_{ij}^*u_{kl}^*+u_{il}^*u_{kj}^*$$
$$u_{ij}u_{kl}+u_{il}u_{kj}=u_{kl}u_{ij}+u_{kj}u_{il}$$
$$u_{ji}u_{lk}+u_{jk}u_{li}=u_{lk}u_{ji}+u_{li}u_{jk}$$

Now by comparing with the original formula, we obtain from this:
$$u_{li}u_{jk}=u_{jk}u_{li}$$

In order to finish, it remains to prove that the coordinates $u_{ij}$ commute as well with their adjoints. For this purpose, we use a similar method. We have:
$$\alpha(x_ix_k^*)=\sum_{jl}x_jx_l^*\otimes u_{ji}u_{lk}^*$$
$$\alpha(x_k^*x_i)=\sum_{jl}x_l^*x_j\otimes u_{lk}^*u_{ji}$$

Since the variables on the left are equal, we deduce from this that we have:
$$\sum_{jl}x_jx_l^*\otimes u_{ji}u_{lk}^*=\sum_{jl}x_jx_l^*\otimes u_{lk}^*u_{ji}$$

Thus we have $u_{ji}u_{lk}^*=u_{lk}^*u_{ji}$, and so $G\subset U_N$, as claimed.

\medskip

\underline{$S^{N-1}_\mathbb R$}. The fact that we have an action $O_N\curvearrowright S^{N-1}_\mathbb R$ is clear. In what regards the converse, this follows by combining the results that we already have, as follows:
\begin{eqnarray*}
G\curvearrowright S^{N-1}_\mathbb R
&\implies&G\curvearrowright S^{N-1}_{\mathbb R,+},S^{N-1}_\mathbb C\\
&\implies&G\subset O_N^+,U_N\\
&\implies&G\subset O_N^+\cap U_N=O_N 
\end{eqnarray*}

Thus, we conclude that we have $G^+(S^{N-1}_\mathbb R)=O_N$, as desired.
\end{proof}

\section*{14b. Integration theory} 

Let us discuss now the correspondence $U\to S$. In the classical case the situation is very simple, because the sphere $S=S^{N-1}$ appears by rotating the point $x=(1,0,\ldots,0)$ by the isometries in $U=U_N$. Moreover, the stabilizer of this action is the subgroup $U_{N-1}\subset U_N$ acting on the last $N-1$ coordinates, and so the sphere $S=S^{N-1}$ appears from the corresponding rotation group $U=U_N$ as an homogeneous space, as follows:
$$S^{N-1}=U_N/U_{N-1}$$

In functional analytic terms, all this becomes even simpler, the correspondence $U\to S$ being obtained, at the level of algebras of functions, as follows:
$$C(S^{N-1})\subset C(U_N)\quad,\quad 
x_i\to u_{1i}$$

In general now, the straightforward homogeneous space interpretation of $S$ as above fails. However, we can have some theory going by using the functional analytic viewpoint, with an embedding $x_i\to u_{1i}$ as above. Let us start with the following result:

\begin{proposition}
For the basic spheres, we have a diagram as follows,
$$\xymatrix@R=50pt@C=50pt{
C(S)\ar[r]^\alpha\ar[d]^\gamma&C(S)\otimes C(U)\ar[d]^{\gamma\otimes id}\\
C(U)\ar[r]^\Delta&C(U)\otimes C(U)
}$$
where on top $\alpha(x_i)=\sum_jx_j\otimes u_{ji}$, and on the left $\gamma(x_i)=u_{1i}$.
\end{proposition}

\begin{proof}
The diagram in the statement commutes indeed on the standard coordinates, the corresponding arrows being as follows, on these coordinates:
$$\xymatrix@R=50pt@C=50pt{
x_i\ar[r]\ar[d]&\sum_jx_j\otimes u_{ji}\ar[d]\\
u_{1i}\ar[r]&\sum_ju_{1j}\otimes u_{ji}
}$$

Thus by linearity and multiplicativity, the whole the diagram commutes.
\end{proof}

As a consequence, we have the following result:

\begin{proposition}
We have a quotient map and an inclusion $U\to S_U\subset S$, with 
$$C(S_U)=<u_{1i}>\subset C(U)$$
at the level of the corresponding algebras of functions.
\end{proposition}

\begin{proof}
At the algebra level, we have an inclusion and a quotient map as follows:
$$C(S)\to C(S_U)\subset C(U)$$

Thus, we obtain the result, by transposing.
\end{proof}

We will prove in what follows that the inclusion $S_U\subset S$ constructed above is an isomorphism. This will produce the correspondence $U\to S$ that we are currently looking for. In order to do so, we will use the uniform integration over $S$, which can be introduced, in analogy with what happens in the classical case, in the folowing way:

\index{integration over spheres}

\begin{definition}
We endow each of the algebras $C(S)$ with the functional
$$\int_S:C(S)\to C(U)\to\mathbb C$$
obtained by composing the morphism $x_i\to u_{1i}$ with the Haar integration of $C(U)$.
\end{definition}

In order to efficiently integrate over the sphere $S$, and in the lack of some trick like spherical coordinates, we need to know how to efficiently integrate over the corresponding quantum isometry group $U$. There is a long story here, and following \cite{bgo}, we have:

\index{integration over spheres}
 
\begin{theorem}
The integration over the basic spheres is given by
$$\int_Sx_{i_1}^{e_1}\ldots x_{i_k}^{e_k}=\sum_\pi\sum_{\sigma\leq\ker i}W_{kN}(\pi,\sigma)$$
with $\pi,\sigma\in D(k)$, where $W_{kN}=G_{kN}^{-1}$ is the inverse of $G_{kN}(\pi,\sigma)=N^{|\pi\vee\sigma|}$. 
\end{theorem}

\begin{proof}
According to our conventions, the integration over $S$ is a particular case of the integration over $U$, via $x_i=u_{1i}$. By using the Weingarten formula, we obtain:
\begin{eqnarray*}
\int_Sx_{i_1}^{e_1}\ldots x_{i_k}^{e_k}
&=&\int_Uu_{1i_1}^{e_1}\ldots u_{1i_k}^{e_k}\\
&=&\sum_{\pi,\sigma\in D(k)}\delta_\pi(1)\delta_\sigma(i)W_{kN}(\pi,\sigma)\\
&=&\sum_{\pi,\sigma\in D(k)}\delta_\sigma(i)W_{kN}(\pi,\sigma)
\end{eqnarray*}

Thus, we are led to the formula in the statement.
\end{proof}

\index{ergodicity}

Again following \cite{bgo}, we have the following key result:

\begin{theorem}
The integration functional of $S$ has the ergodicity property
$$\left(id\otimes\int_U\right)\alpha(x)=\int_Sx$$
where $\alpha:C(S)\to C(S)\otimes C(U)$ is the universal affine coaction map.
\end{theorem}

\begin{proof}
In the real case, $x_i=x_i^*$, it is enough to check the equality in the statement on an arbitrary product of coordinates, $x_{i_1}\ldots x_{i_k}$. The left term is as follows:
\begin{eqnarray*}
\left(id\otimes\int_U\right)\alpha(x_{i_1}\ldots x_{i_k})
&=&\sum_{j_1\ldots j_k}x_{j_1}\ldots x_{j_k}\int_Uu_{j_1i_1}\ldots u_{j_ki_k}\\
&=&\sum_{j_1\ldots j_k}\ \sum_{\pi,\sigma\in D(k)}\delta_\pi(j)\delta_\sigma(i)W_{kN}(\pi,\sigma)x_{j_1}\ldots x_{j_k}\\
&=&\sum_{\pi,\sigma\in D(k)}\delta_\sigma(i)W_{kN}(\pi,\sigma)\sum_{j_1\ldots j_k}\delta_\pi(j)x_{j_1}\ldots x_{j_k}
\end{eqnarray*}

Let us look now at the last sum on the right. The situation is as follows:

\medskip

(1) In the free case we have to sum quantities of type $x_{j_1}\ldots x_{j_k}$, over all choices of multi-indices $j=(j_1,\ldots,j_k)$ which fit into our given noncrossing pairing $\pi$, and just by using the condition $\sum_ix_i^2=1$, we conclude that the sum is 1. 

\medskip

(2) The same happens in the classical case. Indeed, our pairing $\pi$ can now be crossing, but we can use the commutation relations $x_ix_j=x_jx_i$, and the sum is again 1.

\medskip

Thus the sum on the right is 1, in all cases, and we obtain:
$$\left(id\otimes\int_U\right)\alpha(x_{i_1}\ldots x_{i_k})
=\sum_{\pi,\sigma\in D(k)}\delta_\sigma(i)W_{kN}(\pi,\sigma)$$

But, another application of the Weingarten formula gives, as desired:
\begin{eqnarray*}
\int_Sx_{i_1}\ldots x_{i_k}
&=&\int_Uu_{1i_1}\ldots u_{1i_k}\\
&=&\sum_{\pi,\sigma\in D(k)}\delta_\pi(1)\delta_\sigma(i)W_{kN}(\pi,\sigma)\\
&=&\sum_{\pi,\sigma\in D(k)}\delta_\sigma(i)W_{kN}(\pi,\sigma)
\end{eqnarray*}

In the complex case the proof is similar, by adding exponents everywhere.
\end{proof} 

Still following \cite{bgo}, we have as well the following useful abstract result:

\index{integration over spheres}

\begin{theorem}
There is a unique positive unital trace $tr:C(S)\to\mathbb C$ satisfying
$$(tr\otimes id)\alpha(x)=tr(x)1$$
where $\alpha$ is the coaction map of the corresponding quantum isometry group,
$$\alpha:C(S)\to C(S)\otimes C(U)$$
and this is the canonical integration, as constructed in Definition 14.13.
\end{theorem}

\begin{proof}
First of all, it follows from the Haar integral invariance condition for $U$ that the canonical integration has indeed the invariance property in the statement, namely:
$$(tr\otimes id)\alpha(x)=tr(x)1$$

In order to prove now the uniqueness, let $tr$ be as in the statement. We have:
\begin{eqnarray*}
tr\left(id\otimes\int_U\right)\alpha(x)
&=&\int_U(tr\otimes id)\alpha(x)\\
&=&\int_U(tr(x)1)\\
&=&tr(x)
\end{eqnarray*}

On the other hand, according to Theorem 14.15, we have as well:
$$
tr\left(id\otimes\int_U\right)\alpha(x)
=tr\left(\int_Sx\right)
=\int_Sx$$

We therefore conclude that $tr$ equals the standard integration, as claimed.
\end{proof}

Getting back now to our axiomatization questions, we have:

\begin{theorem}
The operation $S\to S_U$ produces a correspondence as follows,
$$\xymatrix@R=15mm@C=15mm{
S^{N-1}_{\mathbb R,+}\ar[r]&S^{N-1}_{\mathbb C,+}\\
S^{N-1}_\mathbb R\ar[r]\ar[u]&S^{N-1}_\mathbb C\ar[u]}
\qquad
\xymatrix@R=8mm@C=15mm{\\ \to}
\qquad
\xymatrix@R=17mm@C=16mm{
O_N^+\ar[r]&U_N^+\\
O_N\ar[r]\ar[u]&U_N\ar[u]
}$$
between basic unitary groups and the basic noncommutative spheres.
\end{theorem}

\begin{proof}
We use the ergodicity formula from Theorem 14.15, namely:
$$\left(id\otimes\int_U\right)\alpha=\int_S$$

We know that $\int_U$ is faithful on $\mathcal C(U)$, and that we have:
$$(id\otimes\varepsilon)\alpha=id$$

The coaction map $\alpha$ follows to be faithful as well. Thus for any $x\in\mathcal C(S)$ we have:
$$\int_Sxx^*=0\implies x=0$$

Thus $\int_S$ is faithful on $\mathcal C(S)$. But this shows that we have:
$$S=S_U$$

Thus, we are led to the conclusion in the statement.
\end{proof}

Getting now to probabilistic aspects, it is well-known that for the real sphere $S^{N-1}_\mathbb R$ the laws of the standard coordinates, which are called hyperspherical laws, become normal and independent, with $N\to\infty$. Here is the free analogue of this result:

\begin{theorem}
For the free sphere $S^{N-1}_{\mathbb R,+}$, the rescaled coordinates 
$$y_i=\sqrt{N}x_i$$
become semicircular and free, in the $N\to\infty$ limit.
\end{theorem}

\begin{proof}
The Weingarten formula for the free sphere, discussed above, together with the standard fact that the Gram matrix, and hence the Weingarten matrix too, is asymptotically diagonal, gives the following estimate, for the joint moments of coordinates:
$$\int_{S^{N-1}_{\mathbb R,+}}x_{i_1}\ldots x_{i_k}\,dx\simeq N^{-k/2}\sum_{\sigma\in NC_2(k)}\delta_\sigma(i_1,\ldots,i_k)$$

With this formula in hand, we can compute the asymptotic moments of each coordinate $x_i$. Indeed, by setting $i_1=\ldots=i_k=i$, all Kronecker symbols are 1, and we obtain:
$$\int_{S^{N-1}_{\mathbb R,+}}x_i^k\,dx\simeq N^{-k/2}|NC_2(k)|$$

Thus the rescaled coordinates $y_i=\sqrt{N}x_i$ become semicircular in the $N\to\infty$ limit, as claimed. As for the asymptotic freeness result, this follows as well from the above general joint moment estimate, via standard free probability theory. See \cite{bco}, \cite{bgo}.
\end{proof}

For the complex spheres the results are similar, with the standard coordinates on $S^{N-1}_\mathbb C$ becoming complex normal and independent, with $N\to\infty$, and with those on $S^{N-1}_{\mathbb C,+}$ becoming circular and free, with $N\to\infty$. As before, we refer here to \cite{bco}, \cite{bgo}.

\bigskip

Summarizing, we have good results for the free spheres, with $N\to\infty$. The problem now is that of computing the moments of the coordinates of the free spheres at fixed values of $N\in\mathbb N$. In the real case, the answer here, which is non-trivial, is as follows:

\index{free hyperspherical law}
\index{special functions}
\index{twisting}

\begin{theorem}
The moments of the free hyperspherical law are given by
$$\int_{S^{N-1}_{\mathbb R,+}}x_1^{2l}\,dx=\frac{1}{(N+1)^l}\cdot\frac{q+1}{q-1}\cdot\frac{1}{l+1}\sum_{r=-l-1}^{l+1}(-1)^r\begin{pmatrix}2l+2\cr l+r+1\end{pmatrix}\frac{r}{1+q^r}$$
where $q\in [-1,0)$ is such that $q+q^{-1}=-N$.
\end{theorem}

\begin{proof}
This is something quite tricky. To be more precise, this follow in 4 steps, none of which is something trivial, which are as follows:

\medskip

(1) $x_1\in C(S^{N-1}_{\mathbb R,+})$ has the same law as $u_{11}\in C(O_N^+)$.

\medskip

(2) $u_{11}\in C(O_N^+)$ has the same law as a certain variable $w\in C(SU^q_2)$.

\medskip

(3) $w\in C(SU^q_2)$ can be in turn modeled by an explicit operator $T\in B(l^2(\mathbb N))$.

\medskip

(4) The law of $T\in B(l^2(\mathbb N))$ can be computed by using advanced calculus.

\medskip

Let us first explain the relation between $O_N^+$ and $SU^q_2$. To any matrix $F\in GL_N(\mathbb R)$ satisfying $F^2=1$ we associate the following universal algebra:
$$C(O_F^+)=C^*\left((u_{ij})_{i,j=1,\ldots,N}\Big|u=F\bar{u}F={\rm unitary}\right)$$

Observe that $O_{I_N}^+=O_N^+$. In general, the above algebra satisfies Woronowicz's generalized axioms in \cite{wo1}, which do not include the strong antipode axiom $S^2=id$.

\medskip

At $N=2$, up to a trivial equivalence relation on the matrices $F$, and on the quantum groups $O_F^+$, we can assume that $F$ is as follows, with $q\in [-1,0)$:
$$F=\begin{pmatrix}0&\sqrt{-q}\\
1/\sqrt{-q}&0\end{pmatrix}$$

Our claim is that for this matrix we have an isomorphism as follows:
$$O_F^+=SU^q_2$$

Indeed, the relations $u=F\bar{u}F$ tell us that $u$ must be of the following form:
$$u=\begin{pmatrix}\alpha&-q\gamma^*\cr \gamma&\alpha^*\end{pmatrix}$$

Thus $C(O_F^+)$ is the universal algebra generated by two elements $\alpha,\gamma$, with the relations making the above matrix $u$ unitary. But these unitarity conditions are:
$$\alpha\gamma=q\gamma\alpha\quad,\quad
\alpha\gamma^*=q\gamma^*\alpha\quad,\quad 
\gamma\gamma^*=\gamma^*\gamma$$
$$\alpha^*\alpha+\gamma^*\gamma=1\quad,\quad 
\alpha\alpha^*+q^2\gamma\gamma^*=1$$

We recognize here the relations in \cite{wo1} defining the algebra $C(SU^q_2)$, and it follows that we have an isomorphism of Hopf $C^*$-algebras, as follows:
$$C(O_F^+)\simeq C(SU^q_2)$$

Now back to the general case, let us try to understand the integration over $O_F^+$. Given $\pi\in NC_2(2k)$ and $i=(i_1,\ldots,i_{2k})$, we set:
$$\delta_\pi^F(i)=\prod_{s\in\pi}F_{i_{s_l}i_{s_r}}$$

Here the product is over all strings of $\pi$, denoted as follows:
$$s=\{s_l\curvearrowright s_r\}$$

Our claim now is that the following family of vectors, with $\pi\in NC_2(2k)$, spans the space of fixed vectors of $u^{\otimes 2k}$:
$$\xi_\pi=\sum_i\delta_\pi^F(i)e_{i_1}\otimes\ldots\otimes e_{i_{2k}}$$ 

Indeed, having $\xi_\cap$ fixed by $u^{\otimes 2}$ is equivalent to assuming that $u=F\bar{u}F$ is unitary. By using now the above vectors, we obtain the following Weingarten formula:
$$\int_{O_F^+}u_{i_1j_1}\ldots u_{i_{2k}j_{2k}}=\sum_{\pi\sigma}\delta_\pi^F(i)\delta_\sigma^F(j)W_{kN}(\pi,\sigma)$$

With these preliminaries in hand, let us start the computation. Let $N\in\mathbb N$, and consider the number $q\in [-1,0)$ satisfying:
$$q+q^{-1}=-N$$

Our claim is that we have the following formula:
$$\int_{O_N^+}\varphi(\sqrt{N+2}\,u_{ij})=\int_{SU^q_2}\varphi(\alpha+\alpha^*+\gamma-q\gamma^*)$$

Indeed, the moments of the variable on the left are given by:
$$\int_{O_N^+}u_{ij}^{2k}=\sum_{\pi\sigma}W_{kN}(\pi,\sigma)$$

On the other hand, the moments of the variable on the right, which in terms of the fundamental corepresentation $v=(v_{ij})$ is given by $w=\sum_{ij}v_{ij}$, are given by:
$$\int_{SU^q_2}w^{2k}=\sum_{ij}\sum_{\pi\sigma}\delta_\pi^F(i)\delta_\sigma^F(j)W_{kN}(\pi,\sigma)$$

We deduce that $w/\sqrt{N+2}$ has the same moments as $u_{ij}$, which proves our claim. In order to do now the computation over $SU^q_2$, we can use a matrix model due to Woronowicz \cite{wo1}, where the standard generators $\alpha,\gamma$ are mapped as follows:
\begin{eqnarray*}
\pi_u(\alpha)e_k&=&\sqrt{1-q^{2k}}e_{k-1}\\
\pi_u(\gamma)e_k&=&uq^k e_k
\end{eqnarray*}

Here $u\in\mathbb T$ is a parameter, and $(e_k)$ is the standard basis of $l^2(\mathbb N)$. The point with this representation is that it allows the computation of the Haar functional. Indeed, if $D$ is the diagonal operator given by $D(e_k)=q^{2k}e_k$, then the formula is as follows:
$$\int _{SU^q_2}x=(1-q^2)\int_{\mathbb T}tr(D\pi_u(x))\frac{du}{2\pi iu}$$

With the above model in hand, the law of the variable that we are interested in is of the following form:
$$\int_{SU^q_2}\varphi(\alpha+\alpha^*+\gamma-q\gamma^*)=(1-q^2)\int_{\mathbb T}tr(D\varphi(M))\frac{du}{2\pi iu}$$

To be more precise, this formula holds indeed, with:
$$M(e_k)=e_{k+1}+q^k(u-qu^{-1})e_k+(1-q^{2k})e_{k-1}$$

The point now is that the integral on the right can be computed, by using advanced calculus methods, and this gives the result.
\end{proof}

The computation of the joint free hyperspherical laws remains an open problem. Open as well is the question of finding a more conceptual proof for the above formula. For more on all this, and related topics, we refer to \cite{bbs}, \cite{bco}, \cite{bgo} and related papers.

\section*{14c. Further spheres}

Getting back now to algebra and philosophy, based on our quantum sphere work above, we have evidence for the existence of four ``main geometries'', namely real and complex, classical and free. But the question is, can we axiomatize such geometries?

\bigskip

In order to discuss this question, the idea will be that of starting with what we have, as geometric objects so far, namely the spheres $S$, the tori $T$, the unitary groups $U$, and the reflection groups $K$. Not much, but as we will soon discover, these will do.

\bigskip

To be more precise, we first have the following result:

\begin{theorem}
We have basic quadruplets $(S,T,U,K)$, as follows:
\begin{enumerate}
\item A classical real and a classical complex quadruplet, as follows:
$$
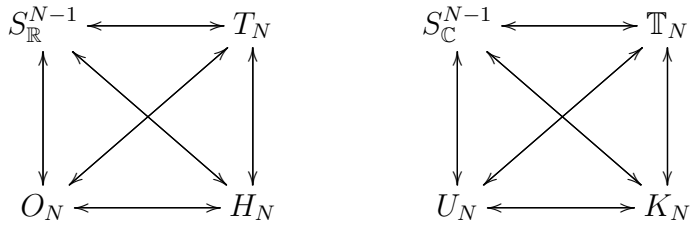
\xymatrix@R=50pt@C=50pt{
S^{N-1}_\mathbb R\ar[r]\ar[d]\ar[dr]&T_N\ar[l]\ar[d]\ar[dl]\\
O_N\ar[u]\ar[ur]\ar[r]&H_N\ar[l]\ar[ul]\ar[u]}
\qquad\qquad 
\xymatrix@R=50pt@C=50pt{
S^{N-1}_\mathbb C\ar[r]\ar[d]\ar[dr]&\mathbb T_N\ar[l]\ar[d]\ar[dl]\\
U_N\ar[u]\ar[ur]\ar[r]&K_N\ar[l]\ar[ul]\ar[u]}$$

\item A free real and a free complex quadruplet, as follows:
$$
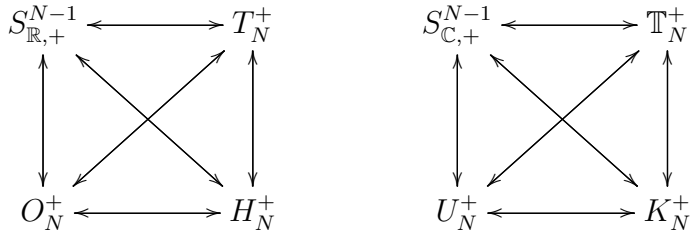
\xymatrix@R=50pt@C=50pt{
S^{N-1}_{\mathbb R,+}\ar[r]\ar[d]\ar[dr]&T_N^+\ar[l]\ar[d]\ar[dl]\\
O_N^+\ar[u]\ar[ur]\ar[r]&H_N^+\ar[l]\ar[ul]\ar[u]}
\qquad\qquad
\xymatrix@R=50pt@C=50pt{
S^{N-1}_{\mathbb C,+}\ar[r]\ar[d]\ar[dr]&\mathbb T_N^+\ar[l]\ar[d]\ar[dl]\\
U_N^+\ar[u]\ar[ur]\ar[r]&K_N^+\ar[l]\ar[ul]\ar[u]}$$
\end{enumerate}
Moreover, for all these quadruplets, the unitary quantum group $U$ is easy.
\end{theorem}

\begin{proof}
Here the various objects appearing in the above diagrams are objects that we know well, constructed at various places, in this book, and the last assertion, regarding easiness, is something that we know well too. As for the fact that, in each of the 4 cases under investigation, we have indeed a full set of 12 correspondences between these objects, this is something quite routine, and for details here, we refer to \cite{ba3}.
\end{proof}

The point now is that we can reformulate Theorem 14.20 as follows:

\begin{theorem}
We have $4$ basic noncommutative geometries,
$$
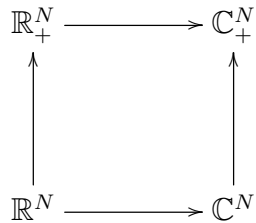
\xymatrix@R=52pt@C=52pt{
\mathbb R^N_+\ar[r]&\mathbb C^N_+\\
\mathbb R^N\ar[u]\ar[r]&\mathbb C^N\ar[u]
}$$
called classical real and complex, and free real and complex.
\end{theorem}

\begin{proof}
This is indeed a reformulation of Theorem 14.20, with the convention that a quadruplet as there corresponds to a noncommutative geometry, that we can denote and call as we like, and with the best notations and terminology being as above.
\end{proof}

Getting now to more technical aspects, when trying to construct noncommutative geometries by starting with more complicated unitary quantum groups $U$, things can be quite tricky. Clarifying all this was in fact an open problem, all over the 10s, with several papers written on the subject. In order to discuss the solution, which came in the late 10s, let us start with more details regarding our axioms for noncommutative geometries. As explained in \cite{ba3}, these axioms can be taken as follows:

\begin{definition}
A quadruplet $(S,T,U,K)$ is said to produce a noncommutative geometry when one can pass from each object to all the other objects, as follows,
$$\begin{matrix}
S&=&S_{<O_N,T>}&=&S_U&=&S_{<O_N,K>}\\
\\
S\cap\mathbb T_N^+&=&T&=&U\cap\mathbb T_N^+&=&K\cap\mathbb T_N^+\\
\\
G^+(S)&=&<O_N,T>&=&U&=&<O_N,K>\\
\\
K^+(S)&=&K^+(T)&=&U\cap K_N^+&=&K
\end{matrix}$$
with the usual convention that all this is up to the equivalence relation.
\end{definition}

There axioms can look a bit complicated, at a first glance, but they are in fact very simple and natural, inspired from what happens in the classical case, and with a look at the free case too. To be more precise, what we have above are all sorts of objects and operations that we know well, along with the operation $U\to S_U$, which consists in taking the first row space, the operation $X\to G^+(X)$, which consists in taking the quantum isometry group, and finally the operation $X\to K^+(X)$, which consists in taking the reflection isometry group, $K^+(X)=G^+(X)\cap K_N^+$. As for the fact that these axioms are indeed satisfied in the 4 main cases of interest, as claimed in Theorem 14.21, this is something well-known in the classical case, and is routine to check in the free case.

\bigskip

In the easy case now, in order to reach to concrete results, it is convenient to formulate an independent definition, using the easy generation operation $\{\,,\}$ instead of the usual generation operation $<\,,>$. This definition, which is more or less a particular case of Definition 14.22, modulo the usual issues with the differences between $\{\,,\}$ and $<\,,>$, that we have met several times, and are quite familiar with, is as follows:

\begin{definition}
A quadruplet $(S,T,U,K)$ is said to produce an easy geometry when $U,K$ are easy, and one can pass from each object to all the other objects, as follows,
$$\begin{matrix}
S&=&S_{\{O_N,K^+(T)\}}&=&S_U&=&S_{\{O_N,K\}}\\
\\
S\cap\mathbb T_N^+&=&T&=&U\cap\mathbb T_N^+&=&K\cap\mathbb T_N^+\\
\\
G^+(S)&=&\{O_N,K^+(T)\}&=&U&=&\{O_N,K\}\\
\\
K^+(S)&=&K^+(T)&=&U\cap K_N^+&=&K
\end{matrix}$$
with the usual convention that all this is up to the equivalence relation.
\end{definition}

Getting now into classification results, the idea is to focus on the quantum group content of the above definition. Indeed, we know that both the quantum groups $U,K$ are easy, and that the following easy generation formula must be satisfied:
$$U=\{O_N,K\}$$

Combinatorially, this leads to the following statement:

\index{easy geometry}
\index{category of partitions}
\index{category of pairings}
\index{easy generation}

\begin{proposition}
An easy geometry is uniquely determined by a pair $(D,E)$ of categories of partitions, which must be as follows,
$$\mathcal{NC}_2\subset D\subset P_2\quad,\quad 
\mathcal{NC}_{even}\subset E\subset P_{even}$$
and which are subject to the following intersection and generation conditions,
$$D=E\cap P_2\quad,\quad 
E=<D,\mathcal{NC}_{even}>$$
and to the usual axioms for the associated quadruplet $(S,T,U,K)$, where $U,K$ are respectively the easy quantum groups associated to the categories $D,E$.
\end{proposition}

\begin{proof}
This comes from the following conditions, with the first one being the one mentioned above, and with the second one being part of our general axioms: 
$$U=\{O_N,K\}\quad,\quad 
K=U\cap K_N^+$$

Indeed, $U,K$ must be easy, coming from certain categories of partitions $D,E$. It is clear that $D,E$ must appear as intermediate categories, as in the statement, and the fact that the intersection and generation conditions must be satisfied follows from:
\begin{eqnarray*}
U=\{O_N,K\}&\iff&D=E\cap P_2\\
K=U\cap K_N^+&\iff&E=<D,\mathcal{NC}_{even}>
\end{eqnarray*}

Thus, we are led to the conclusion in the statement.
\end{proof}

In order to discuss now classification results, we would need some technical results regarding the intermediate easy quantum groups as follows:
$$O_N\subset U\subset U_N^+\quad,\quad 
H_N\subset K\subset K_N^+$$

But we do have such results, as explained in the quantum group literature, and by using this, after some work, we are led to the following classification result:

\begin{theorem}
Under strong combinatorial axioms, of easiness and uniformity type, we have only $9$ noncommutative geometries, namely:
$$\xymatrix@R=40pt@C=40pt{
\mathbb R^N_+\ar[r]&\mathbb T\mathbb R^N_+\ar[r]&\mathbb C^N_+\\
\mathbb R^N_*\ar[u]\ar[r]&\mathbb T\mathbb R^N_*\ar[u]\ar[r]&\mathbb C^N_*\ar[u]\\
\mathbb R^N\ar[u]\ar[r]&\mathbb T\mathbb R^N\ar[u]\ar[r]&\mathbb C^N\ar[u]
}$$
Moreover, under even stronger combinatorial axioms, including a slicing condition, the $4$ basic geometries, those at the corners, are the only ones.
\end{theorem}

\begin{proof}
This is something quite technical, for which we refer to \cite{ba3}, but that we can basically understand with our quantum group knowledge, the idea being as follows:

\medskip

(1) To start with, the geometries in the statement correspond to the main examples of intermediate quantum groups $O_N\subset U\subset U_N^+$, which are as follows:
$$\xymatrix@R=13mm@C=13mm{
O_N^+\ar[r]&\mathbb TO_N^+\ar[r]&U_N^+\\
O_N^*\ar[r]\ar[u]&\mathbb TO_N^*\ar[r]\ar[u]&U_N^*\ar[u]\\
O_N\ar[r]\ar[u]&\mathbb TO_N\ar[r]\ar[u]&U_N\ar[u]}$$

(2) As for the corresponding reflection groups, these correspond to the main examples of intermediate quantum groups $H_N\subset U\subset K_N^+$, which are as follows:
$$\xymatrix@R=13mm@C=13mm{
H_N^+\ar[r]&\mathbb TH_N^+\ar[r]&K_N^+\\
H_N^*\ar[r]\ar[u]&\mathbb TH_N^*\ar[r]\ar[u]&K_N^*\ar[u]\\
H_N\ar[r]\ar[u]&\mathbb TH_N\ar[r]\ar[u]&K_N\ar[u]}$$

(3) With these conventions made, telling us who the quantum groups $U,K$ are, in each of the 9 cases under investigation, we can complete our quadruplets with objects $S,T$, by using either of the formulae involving them from Definition 14.23, and the verification of the axioms from Definition 14.23 is straightforward, in each of these 9 cases.

\medskip

(4) Finally, the classification assertions are more technical, whose proofs are basically based on the study of the correspondence $U\leftrightarrow K$ from Definition 14.23. To be more precise, we know from Proposition 14.24 that the unitary group $U$ of our easy geometry must come from a category of pairings $D\subset P_2$ satisfying the following condition:
$$D=<D,\mathcal{NC}_{even}>\cap  P_2$$

(5) But this equation can be solved by using the classification results available from the quantum group literature, and we are led to the conclusions in the statement.

\medskip

(6) So, this was for the idea, and in practice now, all this needs a massive amount of routine verifications, at each single step, and all this is explained in \cite{ba3}. 
\end{proof}

Very nice all this, and as a question that you might have, right now, sure that we love our 4 main geometries, those at the corners of the square in Theorem 14.25, but can we get rid somehow of the other ones, the intermediate ones, and to have something really beautiful, as final philosophical conclusion? Good point, and in answer:

\bigskip

(1) As already mentioned at the end of Theorem 14.25, it is possible to formally do this, with the help of some extra combinatorical axioms.

\bigskip

(2) More drastically, when passing from affine to projective, all intermediate objects dissapear, and in addition, the real and complex geometries get identified too.

\bigskip

In short, the story if not over with the above results, quite the opposite, and we will discuss the continuation of all this later, in the projective setting, in chapter 16 below.

\section*{14d. Free quadrics}

Getting back now to the free spheres, which remain our main quantum manifolds constructed so far, let us discuss, more generally, the compact hypersurfaces $X\subset\mathbb R^N_+$. These hypersurfaces fit into the $C^*$-algebra formalism, as follows:

\index{free hypersurface}
\index{noncommutative hypersurface}

\begin{definition}
A real compact hypersurface in $N$ variables, denoted $X_f\subset\mathbb R^N_+$, is the abstract spectrum of a universal $C^*$-algebra of the following type,
$$C(X_f)=C^*\left(x_1,\ldots,x_N\Big|x_i=x_i^*,f(x_1,\ldots,x_N)=0\right)$$
with the noncommutative polynomial $f\in\mathbb R<x_1,\ldots,x_N>$ being such the maximal $C^*$-norm on the complex $*$-algebra $\mathbb C<x_1,\ldots,x_N>/(f)$ is bounded.
\end{definition}

Summarizing, good news, and job for us, to study such beasts. As a first result now, regarding such hypersurfaces, coming from the Gelfand theorem, we have:

\begin{theorem}
In order for $X_f$ to exist, the real algebraic manifold
$$X_f^\times=\left\{x\in\mathbb R^N\Big|f(x_1,\ldots,x_N)=0\right\}$$
must be compact. In addition, in this case we have $||x_i||_\times\leq||x_i||$, for any $i$.
\end{theorem}

\begin{proof}
Assuming that $X_f$ exists, the Gelfand theorem shows that the algebra of continuous functions on the manifold $X_f^\times$ in the statement appears as:
$$C(X_f^\times)=C(X_f)\Big/\Big<[x_i,x_j]=0\Big>$$

Thus we have an embedding of compact quantum spaces $X_f^\times\subset X_f$, and the norm estimate is clear as well, because such embeddings increase the norms.
\end{proof}

Let us first discuss the quadratic case. We have here:

\index{free quadric}

\begin{theorem}
Given a quadratic polynomial $f\in\mathbb R<x_1,\ldots,x_N>$, written as
$$f=\sum_{ij}A_{ij}x_ix_j+\sum_iB_ix_i+C$$
the following conditions are equivalent:
\begin{enumerate}
\item $X_f$ exists.

\item $X_f^\times$ is compact.

\item The symmetric matrix $Q=\frac{A+A^t}{2}$ is positive or negative.
\end{enumerate}
\end{theorem}

\begin{proof}
The implication $(1)\implies(2)$ being known from Theorem 14.27, and the implication $(2)\iff(3)$ being well-known, we are left with proving $(3)\implies(1)$.  As a first remark here, by applying the adjoint, our manifold $X_f$ is defined by:
$$\begin{cases}
\sum_{ij}A_{ij}x_ix_j+\sum_iB_ix_i+C=0\\
\sum_{ij}A_{ij}x_jx_i+\sum_iB_ix_i+C=0
\end{cases}$$

In terms of $P=\frac{A-A^t}{2}$ and $Q=\frac{A+A^t}{2}$, these equations can be written as:
$$\begin{cases}
\sum_{ij}P_{ij}x_ix_j=0\\
\sum_{ij}Q_{ij}x_ix_j+\sum_iB_ix_i+C=0
\end{cases}$$

Let us first examine the second equation. When regarding $x$ as a column vector, and $B$ as a row vector, this equation becomes an equality of $1\times1$ matrices, as follows:
$$x^tQx+Bx+C=0$$

Now let us assume that $Q$ is positive or negative. Up to a sign change, we can assume $Q>0$. We can write $Q=UDU^t$, with $D=diag(d_i)$ and $d_i>0$, and with $U\in O_N$. In terms of the vector $y=U^tx$, and with $E=BU$, our equation becomes:
$$y^tDy+Ey+C=0$$

By reverting back to sums and indices, this equation reads:
$$\sum_id_iy_i^2+\sum_ie_iy_i+C=0$$

Now by making squares, this equation takes the following form:
$$\sum_id_i\left(y_i+\frac{e_i}{2d_i}\right)^2=c$$

By positivity, we deduce that we have the following estimate:
$$\left\|y_i+\frac{e_i}{2d_i}\right\|^2\leq\frac{|c|}{d_i}$$

Thus our hypersurface $X_f$ is well-defined, and we are done. 
\end{proof}

We have in fact the following result, which is a bit better:

\begin{theorem}
Up to linear changes of coordinates, the free compact quadrics in $\mathbb R^N_+$ are the empty set, the point, the standard free sphere $S^{N-1}_{\mathbb R,+}$, defined by
$$\sum_ix_i^2=1$$
and some intermediate spheres $S^{N-1}_\mathbb R\subset S\subset S^{N-1}_{\mathbb R,+}$, which can be explicitly characterized. Moreover, for all these free quadrics, we have $||x_i||=||x_i||_\times$, for any $i$.
\end{theorem}

\begin{proof}
We use the computations from the proof of Theorem 14.28. The first equation there, making appear the matrix $P=\frac{A-A^t}{2}$, is as follows:
$$\sum_{ij}P_{ij}x_ix_j=0$$ 

As for the second equation, up to a linear change of the coordinates, this reads:
$$\sum_iz_i^2=c$$

At $c<0$ we obtain the empty set. At $c=0$ we must have $z=0$, and depending on whether the first equation is satisfied or not, we obtain either a point, or the empty set. At $c>0$ now, we can assume by rescaling $c=1$, and our second equation reads:
$$X_f\subset S^{N-1}_{\mathbb R,+}$$

As a conclusion, the solutions here are certain subspaces $S\subset S^{N-1}_{\mathbb R,+}$ which appear via equations of type $\sum_{ij}P_{ij}x_ix_j=0$, with $P\in M_N(\mathbb R)$ being antisymmetric, and with $x_1,\ldots,x_N$ appearing via $z_1,\ldots,z_N$ via a linear change of variables. Now observe that when redoing the above computation with $X_f^\times$ at the place of $X_f$, we obtain $X_f=S^{N-1}_\mathbb R$, and this, because the equations $\sum_{ij}P_{ij}x_ix_j=0$ are trivial for commuting variables. We conclude that our subspaces $S\subset S^{N-1}_{\mathbb R,+}$ must satisfy the following condition:
$$S^{N-1}_\mathbb R\subset S\subset S^{N-1}_{\mathbb R,+}$$

Thus, we are left with investigating which such subspaces can indeed be solutions. Observe that both the extreme cases can appear as solutions, as shown by:
\begin{eqnarray*}
X_{2x^2+y^2+\frac{3}{2}xy+\frac{1}{2}yx}&=&S^1_\mathbb R\\
X_{2x^2+y^2+xy+yx}&=&S^1_{\mathbb R,+}
\end{eqnarray*}

Finally, the last assertion is clear for the empty set and for the point, and for the remaining hypersurfaces, this follows from $S^{N-1}_\mathbb R\subset S\subset S^{N-1}_{\mathbb R,+}$.
\end{proof}

Here is yet another version of Theorem 14.28, which is interesting as well:

\index{sum of squares}
\index{SOS}

\begin{theorem}
Given $M$ real linear functions $L_1,\ldots,L_M$ in $N$ noncommuting variables $x_1,\ldots,x_N$, the following are equivalent:
\begin{enumerate}
\item $\sum_kL_k(x_1,\ldots,x_N)^2=1$ defines a compact hypersurface in $\mathbb R^N$.

\item $\sum_kL_k(x_1,\ldots,x_N)^2=1$ defines a compact quantum hypersurface.

\item The matrix formed by the coefficients of $L_1,\ldots,L_M$ has rank $N$.
\end{enumerate}
\end{theorem}

\begin{proof}
The equivalence $(1)\iff(2)$ follows from the equivalence $(1)\iff(2)$ in Theorem 14.28, because the surfaces under investigation are quadrics. As for the equivalence $(2)\iff(3)$, this is well-known too.
\end{proof}

Many other things can be said, along the above lines, the idea being that all this eventually brings us into advanced operator theory, and notably into SOS results.

\section*{14e. Exercises}

We had a quite elementary chapter here, and as exercises on this, we have:

\begin{exercise}
Try defining free balls, in the $C^*$-algebra framework.
\end{exercise}

\begin{exercise}
Work out all the details, for the free hyperspherical laws.
\end{exercise}

\begin{exercise}
Try half-liberating and twisting our various spheres.
\end{exercise}

\begin{exercise}
Can you talk about orientation, on our free spheres?
\end{exercise}

\begin{exercise}
What about a Dirac operator, in the sense of Connes?
\end{exercise}

\begin{exercise}
Learn SOS tricks, and apply them to our hypersurface questions.
\end{exercise}

As bonus exercise, learn some algebraic geometry. We will need this, soon.

\chapter{Free manifolds}

\section*{15a. Quotient spaces}

We have seen that we have four main abstract geometries, namely real and complex, classical and free. In this chapter we develop the free geometries, real and complex.

\bigskip

Let us begin with some generalities regarding the quotient spaces, and more general homogeneous spaces. Regarding the quotients, we have the following construction:

\index{quantum subgroup}
\index{quotient space}
\index{fixed points}

\begin{proposition}
Given a quantum subgroup $H\subset G$, with associated quotient map $\rho:C(G)\to C(H)$, if we define the quotient space $X=G/H$ by setting
$$C(X)=\left\{f\in C(G)\Big|(\rho\otimes id)\Delta f=1\otimes f\right\}$$
then we have a coaction map as follows,
$$\alpha:C(X)\to C(X)\otimes C(G)$$
obtained as the restriction of the comultiplication of $C(G)$. In the classical case, we obtain in this way the usual quotient space $X=G/H$.
\end{proposition}

\begin{proof}
Observe that the linear subspace $C(X)\subset C(G)$ defined in the statement is indeed a subalgebra, because it is defined via a relation of type $\varphi(f)=\psi(f)$, with both $\varphi,\psi$ being morphisms of algebras. Observe also that in the classical case we obtain the algebra of continuous functions on the quotient space $X=G/H$, because:
\begin{eqnarray*}
(\rho\otimes id)\Delta f=1\otimes f
&\iff&(\rho\otimes id)\Delta f(h,g)=(1\otimes f)(h,g),\forall h\in H,\forall g\in G\\
&\iff&f(hg)=f(g),\forall h\in H,\forall g\in G\\
&\iff&f(hg)=f(kg),\forall h,k\in H,\forall g\in G
\end{eqnarray*}

Regarding now the construction of $\alpha$, observe that for $f\in C(X)$ we have: 
\begin{eqnarray*}
(\rho\otimes id\otimes id)(\Delta\otimes id)\Delta f
&=&(\rho\otimes id\otimes id)(id\otimes\Delta)\Delta f\\
&=&(id\otimes\Delta)(\rho\otimes id)\Delta f\\
&=&(id\otimes\Delta)(1\otimes f)\\
&=&1\otimes\Delta f
\end{eqnarray*}

Thus the condition $f\in C(X)$ implies $\Delta f\in C(X)\otimes C(G)$, and this gives the existence of the coaction map $\alpha$. Finally, the other assertions are all clear.
\end{proof}

As an illustration, in the group dual case we have the following result:

\begin{proposition}
Assume that $G=\widehat{\Gamma}$ is a discrete group dual.
\begin{enumerate}
\item The quantum subgroups of $G$ are $H=\widehat{\Lambda}$, with $\Gamma\to\Lambda$ being a quotient group.

\item For such a quantum subgroup $\widehat{\Lambda}\subset\widehat{\Gamma}$, we have $\widehat{\Gamma}/\widehat{\Lambda}=\widehat{\Theta}$, where:
$$\Theta=\ker(\Gamma\to\Lambda)$$
\end{enumerate}
\end{proposition}

\begin{proof}
This is well-known, the idea being as follows:

\medskip

(1) In one sense, this is clear. Conversely, since the algebra $C(G)=C^*(\Gamma)$ is cocommutative, so are all its quotients, and this gives the result.

\medskip

(2) Consider a quotient map $r:\Gamma\to\Lambda$, and denote by $\rho:C^*(\Gamma)\to C^*(\Lambda)$ its extension. Consider a group algebra element, written as follows:
$$f=\sum_{g\in\Gamma}\lambda_g\cdot g\in C^*(\Gamma)$$

We have then the following computation:
\begin{eqnarray*}
f\in C(\widehat{\Gamma}/\widehat{\Lambda})
&\iff&(\rho\otimes id)\Delta(f)=1\otimes f\\
&\iff&\sum_{g\in\Gamma}\lambda_g\cdot r(g)\otimes g=\sum_{g\in\Gamma}\lambda_g\cdot 1\otimes g\\
&\iff&\lambda_g\cdot r(g)=\lambda_g\cdot 1,\forall g\in\Gamma\\
&\iff&supp(f)\subset\ker(r)
\end{eqnarray*}

But this means that we have $\widehat{\Gamma}/\widehat{\Lambda}=\widehat{\Theta}$, with $\Theta=\ker(\Gamma\to\Lambda)$, as claimed.
\end{proof}

Given two compact quantum spaces $X,Y$, we say that $X$ is a quotient space of $Y$ when we have an embedding of $C^*$-algebras $\gamma:C(X)\subset C(Y)$. We have:

\index{homogeneous space}

\begin{definition}
We call a quotient space $G\to X$ homogeneous when
$$\Delta(C(X))\subset C(X)\otimes C(G)$$
where $\Delta:C(G)\to C(G)\otimes C(G)$ is the comultiplication map.
\end{definition}

In other words, an homogeneous quotient space $G\to X$ is a quantum space coming from a subalgebra $C(X)\subset C(G)$, which is stable under the comultiplication. The relation with the quotient spaces from Proposition 15.1 is as follows:

\begin{theorem}
The following results hold:
\begin{enumerate}
\item The quotient spaces $X=G/H$ are homogeneous.

\item In the classical case, any homogeneous space is of type $G/H$.

\item In general, there are homogeneous spaces which are not of type $G/H$.
\end{enumerate}
\end{theorem}

\begin{proof}
Once again these results are well-known, the proof being as follows:

\medskip

(1) This is clear from what we have in Proposition 15.1.

\medskip

(2) Consider a quotient map $p:G\to X$. The invariance condition in the statement tells us that we must have an action $G\curvearrowright X$, given by:
$$g(p(g'))=p(gg')$$

Thus, we have the following implication:
$$p(g')=p(g'')\implies p(gg')=p(gg''),\ \forall g\in G$$

Now observe that the following subset $H\subset G$ is a subgroup:
$$H=\left\{g\in G\Big|p(g)=p(1)\right\}$$

Indeed, $g,h\in H$ implies that we have:
$$p(gh)=p(g)=p(1)$$

Thus we have $gh\in H$, and the other axioms are satisfied as well. Our claim now is that we have an identification $X=G/H$, obtained as follows:
$$p(g)\to Hg$$

Indeed, the map $p(g)\to Hg$ is well-defined and bijective, because $p(g)=p(g')$ is equivalent to $p(g^{-1}g')=p(1)$, and so to $Hg=Hg'$, as desired. 

\medskip

(3) Given a discrete group $\Gamma$ and an arbitrary subgroup $\Theta\subset\Gamma$, the quotient space $\widehat{\Gamma}\to\widehat{\Theta}$ is homogeneous. Now by using Proposition 15.2, we can see that if $\Theta\subset\Gamma$ is not normal, the quotient space $\widehat{\Gamma}\to\widehat{\Theta}$ is not of the form $G/H$.
\end{proof}

With the above formalism in hand, let us try now to understand the general properties of the homogeneous spaces $G\to X$, in the sense of Theorem 15.4. We have:

\index{coaction map}

\begin{proposition}
Assume that a quotient space $G\to X$ is homogeneous.
\begin{enumerate}
\item We have a coaction map as follows, obtained as restriction of $\Delta$:
$$\alpha:C(X)\to C(X)\otimes C(G)$$

\item We have the following implication:
$$\alpha(f)=f\otimes 1\implies f\in\mathbb C1$$

\item We have as well the following formula:
$$\left(id\otimes\int_G\right)\alpha f=\int_Gf$$

\item The restriction of $\int_G$ is the unique form satisfying $(\tau\otimes id)\alpha=\tau(.)1$.
\end{enumerate} 
\end{proposition}

\begin{proof}
These results are all elementary, the proof being as follows:

\medskip

(1) This is clear from definitions, because $\Delta$ itself is a coaction.

\medskip

(2) Assume that $f\in C(G)$ satisfies $\Delta(f)=f\otimes 1$. By applying the counit we obtain:
$$(\varepsilon\otimes id)\Delta f=(\varepsilon\otimes id)(f\otimes 1)$$

We conclude from this that we have $f=\varepsilon(f)1$, as desired.

\medskip

(3) The formula in the statement, $(id\otimes\int_G)\alpha f=\int_Gf$, follows indeed from the left invariance property of the Haar functional of $C(G)$, namely:
$$\left(id\otimes\int_G\right)\Delta f=\int_Gf$$

(4) We use here the right invariance of the Haar functional of $C(G)$, namely:
$$\left(\int_G\otimes id\right)\Delta f=\int_Gf$$

Indeed, we obtain from this that $tr=(\int_G)_{|C(X)}$ is $G$-invariant, in the sense that:
$$(tr\otimes id)\alpha f=tr(f)1$$

Conversely, assuming that $\tau:C(X)\to\mathbb C$ satisfies $(\tau\otimes id)\alpha f=\tau(f)1$, we have:
\begin{eqnarray*}
\left(\tau\otimes\int_G\right)\alpha(f)
&=&\int_G(\tau\otimes id)\alpha(f)\\
&=&\int_G(\tau(f)1)\\
&=&\tau(f)
\end{eqnarray*}

On the other hand, we can compute the same quantity as follows:
\begin{eqnarray*}
\left(\tau\otimes\int_G\right)\alpha(f)
&=&\tau\left(id\otimes\int_G\right)\alpha(f)\\
&=&\tau(tr(f)1)\\
&=&tr(f)
\end{eqnarray*}

Thus we have $\tau(f)=tr(f)$ for any $f\in C(X)$, and this finishes the proof.
\end{proof}

Summarizing, we have a notion of noncommutative homogeneous space, which perfectly covers the classical case. In general, however, the group dual case shows that our formalism is more general than that of the quotient spaces $G/H$.

\bigskip

We discuss now an extra issue, of analytic nature. The point indeed is that for one of the most basic examples of actions, namely $O_N^+\curvearrowright S^{N-1}_{\mathbb R,+}$, the associated morphism $\gamma:C(X)\to C(G)$ is not injective. The same is true for other basic actions, in the free setting. In order to include such examples, we must relax our axioms, as follows:

\index{extended homogeneous space}

\begin{definition}
An extended homogeneous space over a compact quantum group $G$ consists of a morphism of $C^*$-algebras, and a coaction map, as follows,
$$\gamma:C(X)\to C(G)\quad,\quad \alpha:C(X)\to C(X)\otimes C(G)$$
such that the following diagram commutes
$$\xymatrix@R=16mm@C=20mm{
C(X)\ar[r]^\alpha\ar[d]_\gamma&C(X)\otimes C(G)\ar[d]^{\gamma\otimes id}\\
C(G)\ar[r]^\Delta&C(G)\otimes C(G)
}$$
and such that the following diagram commutes as well,
$$\xymatrix@R=16mm@C=20mm{
C(X)\ar[r]^\alpha\ar[d]_\gamma&C(X)\otimes C(G)\ar[d]^{id\otimes\int}\\
C(G)\ar[r]^{\int(.)1}&C(X)
}$$
where $\int$ is the Haar integration over $G$. We write then $G\to X$.
\end{definition}

As a first observation, when the above morphism $\gamma$ is injective we obtain an homogeneous space in the previous sense. The examples with $\gamma$ not injective, which motivate the above formalism, include the standard action $O_N^+\curvearrowright S^{N-1}_{\mathbb R,+}$, and the standard action $U_N^+\curvearrowright S^{N-1}_{\mathbb C,+}$. Here are now a few general remarks, on the above axioms:

\index{coassociativity}
\index{ergodicity}

\begin{proposition}
Assume that we have morphisms of $C^*$-algebras 
$$\gamma:C(X)\to C(G)\quad,\quad 
\alpha:C(X)\to C(X)\otimes C(G)$$ 
satisfying the coassociativity condition $(\gamma\otimes id)\alpha=\Delta\gamma$.
\begin{enumerate}
\item If $\gamma$ is injective on a dense $*$-subalgebra $A\subset C(X)$, and $\alpha(A)\subset A\otimes C(G)$, then $\alpha$ is automatically a coaction map, and is unique.

\item The ergodicity type condition $(id\otimes\int)\alpha=\int\gamma(.)1$ is equivalent to the existence of a linear form $\lambda:C(X)\to\mathbb C$ such that $(id\otimes\int)\alpha=\lambda(.)1$.
\end{enumerate}
\end{proposition}

\begin{proof}
This is something elementary, the idea being as follows:

\medskip

(1) Assuming that we have a dense $*$-subalgebra $A\subset C(X)$ as in the statement, satisying $\alpha(A)\subset A\otimes C(G)$, the restriction $\alpha_{|A}$ is given by:
$$\alpha_{|A}=(\gamma_{|A}\otimes id)^{-1}\Delta\gamma_{|A}$$

This restriction and is therefore coassociative, and unique. By continuity, the morphism $\alpha$ itself follows to be coassociative and unique, as desired.

\medskip

(2) Assuming $(id\otimes\int)\alpha=\lambda(.)1$, we have:
$$\left(\gamma\otimes\int\right)\alpha=\lambda(.)1$$

On the other hand, we have as well the following formula:
$$\left(\gamma\otimes\int\right)\alpha=\left(id\otimes\int\right)\Delta\gamma=\int\gamma(.)1$$

Thus we obtain $\lambda=\int\gamma$, as claimed.
\end{proof}

Given an extended homogeneous space $G\to X$ in our sense, with associated map $\gamma:C(X)\to C(G)$, we can consider the image of this latter map:
$$\gamma:C(X)\to C(Y)\subset C(G)$$

Equivalently, at the level of the associated noncommutative spaces, we can factorize the corresponding quotient map $G\to Y\subset X$. With these conventions, we have:

\begin{proposition}
Consider an extended homogeneous space $G\to X$.
\begin{enumerate}
\item $\alpha(f)=f\otimes 1\implies f\in\mathbb C1$.

\item $tr=\int\gamma$ is the unique unital $G$-invariant form on $C(X)$.

\item The image space obtained by factorizing, $G\to Y$, is homogeneous.
\end{enumerate}
\end{proposition}

\begin{proof}
We have several assertions to be proved, the idea being as follows:

\medskip

(1) This follows indeed from $(id\otimes\int)\alpha(f)=\int\gamma(f)1$, which gives $f=\int\gamma(f)1$.

\medskip

(2) The fact that $tr=\int\gamma$ is indeed $G$-invariant can be checked as follows:
\begin{eqnarray*}
(tr\otimes id)\alpha f
&=&\left(\int\gamma\otimes id\right)\alpha f\\
&=&\left(\int\otimes id\right)\Delta\gamma f\\
&=&\int\gamma(f)1\\
&=&tr(f)1
\end{eqnarray*}

As for the uniqueness assertion, this follows as before.

\medskip

(3) The condition $(\gamma\otimes id)\alpha=\Delta\gamma$, together with the fact that $i$ is injective, allows us to factorize $\Delta$ into a morphism $\Psi$, as follows:
$$\xymatrix@R=12mm@C=30mm{
C(X)\ar[r]^\alpha\ar[d]_\gamma&C(X)\otimes C(G)\ar[d]^{\gamma\otimes id}\\
C(Y)\ar@.[r]^\Psi\ar[d]_i&C(Y)\otimes C(G)\ar[d]^{i\otimes id}\\
C(G)\ar[r]^\Delta&C(G)\otimes C(G)
}$$

Thus the image space $G\to Y$ is indeed homogeneous, and we are done.
\end{proof}

Finally, we have the following result:

\index{GNS construction}

\begin{theorem}
Let $G\to X$ be an extended homogeneous space, and construct quotients $X\to X'$, $G\to G'$ by performing the GNS construction with respect to $\int\gamma,\int$. Then $\gamma$ factorizes into an inclusion $\gamma':C(X')\to C(G')$, and we have an homogeneous space.
\end{theorem}

\begin{proof}
We factorize $G\to Y\subset X$ as above. By performing the GNS construction with respect to $\int i\gamma,\int i,\int$, we obtain a diagram as follows:
$$\xymatrix@R=12mm@C=30mm{
C(X)\ar[r]^p\ar[d]_\gamma&C(X')\ar[d]^{\gamma'}\ar[dr]^{tr'}\\
C(Y)\ar[r]^q\ar[d]_i&C(Y')\ar[d]^{i'}&\mathbb C\\
C(G)\ar[r]^r&C(G')\ar[ur]_{\int'}
}$$

Indeed, with $tr=\int\gamma$, the GNS quotient maps $p,q,r$ are defined respectively by:
\begin{eqnarray*}
\ker p&=&\left\{f\in C(X)\Big|tr(f^*f)=0\right\}\\
\ker q&=&\left\{f\in C(Y)\Big|\smallint(f^*f)=0\right\}\\
\ker r&=&\left\{f\in C(G)\Big|\smallint(f^*f)=0\right\}
\end{eqnarray*}

Next, we can define factorizations $i',\gamma'$ as above. Observe that $i'$ is injective, and that $\gamma'$ is surjective. Our claim now is that $\gamma'$ is injective as well. Indeed:
\begin{eqnarray*}
\gamma'p(f)=0
&\implies&q\gamma(f)=0\\
&\implies&\int\gamma(f^*f)=0\\
&\implies&tr(f^*f)=0\\
&\implies&p(f)=0
\end{eqnarray*}

We conclude that we have $X'=Y'$, and this gives the result.
\end{proof}

\section*{15b. Partial isometries}

Getting now to what we wanted to do in this chapter, namely develop the free geometries, our task will be that of finding a suitable collection of ``free homogeneous spaces'', generalizing at the same time the free spheres $S$, and the free unitary groups $U$. 

\bigskip

This can be done at several levels of generality, and central here is the construction of the free spaces of partial isometries, which can be done in fact for any easy quantum group. In order to explain this, let us start with the classical case. We have here:

\index{partial isometry}

\begin{definition}
Associated to any integers $L\leq M,N$ are the spaces
$$O_{MN}^L=\left\{T:E\to F\ {\rm isometry}\Big|E\subset\mathbb R^N,F\subset\mathbb R^M,\dim_\mathbb RE=L\right\}$$
$$U_{MN}^L=\left\{T:E\to F\ {\rm isometry}\Big|E\subset\mathbb C^N,F\subset\mathbb C^M,\dim_\mathbb CE=L\right\}$$
where the notion of isometry is with respect to the usual real/complex scalar products.
\end{definition}

As a first observation, at $L=M=N$ we obtain the groups $O_N,U_N$:
$$O_{NN}^N=O_N\quad,\quad 
U_{NN}^N=U_N$$ 

Another interesting specialization is $L=M=1$. Here the elements of $O_{1N}^1$ are the isometries $T:E\to\mathbb R$, with $E\subset\mathbb R^N$ one-dimensional. But such an isometry is uniquely determined by $T^{-1}(1)\in\mathbb R^N$, which must belong to $S^{N-1}_\mathbb R$. Thus, we have $O_{1N}^1=S^{N-1}_\mathbb R$. Similarly, in the complex case we have $U_{1N}^1=S^{N-1}_\mathbb C$, and so our results here are:
$$O_{1N}^1=S^{N-1}_\mathbb R\quad,\quad 
U_{1N}^1=S^{N-1}_\mathbb C$$

Yet another interesting specialization is $L=N=1$. Here the elements of $O_{1N}^1$ are the isometries $T:\mathbb R\to F$, with $F\subset\mathbb R^M$ one-dimensional. But such an isometry is uniquely determined by $T(1)\in\mathbb R^M$, which must belong to $S^{M-1}_\mathbb R$. Thus, we have $O_{M1}^1=S^{M-1}_\mathbb R$. Similarly, in the complex case we have $U_{M1}^1=S^{M-1}_\mathbb C$, and so our results here are:
$$O_{M1}^1=S^{M-1}_\mathbb R\quad,\quad
U_{M1}^1=S^{M-1}_\mathbb C$$

In general, the most convenient is to view the elements of $O_{MN}^L,U_{MN}^L$ as rectangular matrices, and to use matrix calculus for their study. We have indeed:

\begin{proposition}
We have identifications of compact spaces
$$O_{MN}^L\simeq\left\{U\in M_{M\times N}(\mathbb R)\Big|UU^t={\rm projection\ of\ trace}\ L\right\}$$
$$U_{MN}^L\simeq\left\{U\in M_{M\times N}(\mathbb C)\Big|UU^*={\rm projection\ of\ trace}\ L\right\}$$
with each partial isometry being identified with the corresponding rectangular matrix.
\end{proposition}

\begin{proof}
We can indeed identify the partial isometries $T:E\to F$ with their corresponding extensions $U:\mathbb R^N\to\mathbb R^M$, $U:\mathbb C^N\to\mathbb C^M$, obtained by setting $U_{E^\perp}=0$. Then, we can identify these latter maps $U$ with the corresponding rectangular matrices.
\end{proof}

As an illustration, at $L=M=N$ we recover in this way the usual matrix description of $O_N,U_N$. Also, at $L=M=1$ we obtain the usual description of $S^{N-1}_\mathbb R,S^{N-1}_\mathbb C$, as row spaces over the corresponding groups $O_N,U_N$. Finally, at $L=N=1$ we obtain the usual description of $S^{N-1}_\mathbb R,S^{N-1}_\mathbb C$, as column spaces over the corresponding groups $O_N,U_N$. 

\bigskip

Now back to the general case, we have the following result:

\begin{proposition}
We have action maps as follows, which are both transitive,
$$O_M\times O_N\curvearrowright O_{MN}^L\quad,\quad 
(A,B)U=AUB^t$$
$$U_M\times U_N\curvearrowright U_{MN}^L\quad,\quad 
(A,B)U=AUB^*$$
whose stabilizers are respectively $O_L\times O_{M-L}\times O_{N-L}$ and $U_L\times U_{M-L}\times U_{N-L}$.
\end{proposition}

\begin{proof}
We have indeed action maps as in the statement, which are transitive. Let us compute now the stabilizer $G$ of the following point:
$$U=\begin{pmatrix}1&0\\0&0\end{pmatrix}$$

Since $(A,B)\in G$ satisfy $AU=UB$, their components must be of the following form:
$$A=\begin{pmatrix}x&*\\0&a\end{pmatrix}\quad,\quad 
B=\begin{pmatrix}x&0\\ *&b\end{pmatrix}$$

Now since $A,B$ are unitaries, these matrices follow to be block-diagonal, and so:
$$G=\left\{(A,B)\Big|A=\begin{pmatrix}x&0\\0&a\end{pmatrix},B=\begin{pmatrix}x&0\\ 0&b\end{pmatrix}\right\}$$

The stabilizer of $U$ is parametrized by triples $(x,a,b)$ belonging to $O_L\times O_{M-L}\times O_{N-L}$ and $U_L\times U_{M-L}\times U_{N-L}$, and we are led to the conclusion in the statement.
\end{proof}

Finally, let us work out the quotient space description of $O_{MN}^L,U_{MN}^L$. We have here:

\begin{theorem}
We have isomorphisms of homogeneous spaces as follows,
\begin{eqnarray*}
O_{MN}^L&=&(O_M\times O_N)/(O_L\times O_{M-L}\times O_{N-L})\\
U_{MN}^L&=&(U_M\times U_N)/(U_L\times U_{M-L}\times U_{N-L})
\end{eqnarray*}
with the quotient maps being given by $(A,B)\to AUB^*$, where $U=(^1_0{\ }^0_0)$.
\end{theorem}

\begin{proof}
This is just a reformulation of Proposition 15.12, by taking into account the fact that the fixed point used in the proof there was $U=(^1_0{\ }^0_0)$.
\end{proof}

Once again, the basic examples here come from the cases $L=M=N$ and $L=M=1$. At $L=M=N$ the quotient spaces at right are respectively:
$$O_N\quad,\quad U_N$$

At $L=M=1$  the quotient spaces at right are respectively:
$$O_N/O_{N-1}\quad,\quad U_N/U_{N-1}$$

In fact, in the general $L=M$ case we obtain the following spaces:
$$O_{MN}^M=O_N/O_{N-M}
\quad,\quad 
U_{MN}^M=U_N/U_{N-M}$$

Similarly, the examples coming from the cases $L=M=N$ and $L=N=1$ are particular cases of the general $L=N$ case, where we obtain the following spaces:
$$O_{MN}^N=O_N/O_{M-N}
\quad,\quad 
U_{MN}^N=U_N/U_{M-N}$$

Summarizing, we have here some basic homogeneous spaces, unifying the spheres with the rotation groups. The point now is that we can liberate these spaces, as follows:

\index{free partial isometry}

\begin{definition}
Associated to any integers $L\leq M,N$ are the algebras
\begin{eqnarray*}
C(O_{MN}^{L+})&=&C^*\left((u_{ij})_{i=1,\ldots,M,j=1,\ldots,N}\Big|u=\bar{u},uu^t={\rm projection\ of\ trace}\ L\right)\\
C(U_{MN}^{L+})&=&C^*\left((u_{ij})_{i=1,\ldots,M,j=1,\ldots,N}\Big|uu^*,\bar{u}u^t={\rm projections\ of\ trace}\ L\right)
\end{eqnarray*}
with the trace being by definition the sum of the diagonal entries.
\end{definition}

Observe that the above universal algebras are indeed well-defined, as it was previously  the case for the free spheres, and this due to the trace conditions, which read: 
$$\sum_{ij}u_{ij}u_{ij}^*
=\sum_{ij}u_{ij}^*u_{ij}
=L$$

We have inclusions between the various spaces constructed so far, as follows:
$$\xymatrix@R=15mm@C=15mm{
O_{MN}^{L+}\ar[r]&U_{MN}^{L+}\\
O_{MN}^L\ar[r]\ar[u]&U_{MN}^L\ar[u]}$$

At the level of basic examples now, at $L=M=1$ and at $L=N=1$ we obtain the following diagrams, showing that our formalism covers indeed the free spheres:
$$\xymatrix@R=15mm@C=15mm{
S^{N-1}_{\mathbb R,+}\ar[r]&S^{N-1}_{\mathbb C,+}\\
S^{N-1}_\mathbb R\ar[r]\ar[u]&S^{N-1}_\mathbb C\ar[u]}
\qquad\qquad 
\xymatrix@R=15mm@C=15mm{
S^{M-1}_{\mathbb R,+}\ar[r]&S^{M-1}_{\mathbb C,+}\\
S^{M-1}_\mathbb R\ar[r]\ar[u]&S^{M-1}_\mathbb C\ar[u]}$$

We have as well the following result, in relation with the free rotation groups:

\begin{proposition}
At $L=M=N$ we obtain the diagram
$$\xymatrix@R=15mm@C=15mm{
O_N^+\ar[r]&U_N^+\\
O_N\ar[r]\ar[u]&U_N\ar[u]}$$
consisting of the groups $O_N,U_N$, and their liberations.
\end{proposition}

\begin{proof}
We recall that the various quantum groups in the statement are constructed as follows, with the symbol $\times$ standing once again for ``commutative'' and ``free'':
\begin{eqnarray*}
C(O_N^\times)&=&C^*_\times\left((u_{ij})_{i,j=1,\ldots,N}\Big|u=\bar{u},uu^t=u^tu=1\right)\\
C(U_N^\times)&=&C^*_\times\left((u_{ij})_{i,j=1,\ldots,N}\Big|uu^*=u^*u=1,\bar{u}u^t=u^t\bar{u}=1\right)
\end{eqnarray*}

On the other hand, according to Proposition 15.11 and to Definition 15.14, we have the following presentation results:
\begin{eqnarray*}
C(O_{NN}^{N\times})&=&C^*_\times\left((u_{ij})_{i,j=1,\ldots,N}\Big|u=\bar{u},uu^t={\rm projection\ of\ trace}\ N\right)\\
C(U_{NN}^{N\times})&=&C^*_\times\left((u_{ij})_{i,j=1,\ldots,N}\Big|uu^*,\bar{u}u^t={\rm projections\ of\ trace}\ N\right)
\end{eqnarray*}

We use now the standard fact that if $p=aa^*$ is a projection then $q=a^*a$ is a projection too. We use as well the following formulae:
$$Tr(uu^*)=Tr(u^t\bar{u})\quad,\quad 
Tr(\bar{u}u^t)=Tr(u^*u)$$

We therefore obtain the following formulae:
\begin{eqnarray*}
C(O_{NN}^{N\times})&=&C^*_\times\left((u_{ij})_{i,j=1,\ldots,N}\Big|u=\bar{u},\ uu^t,u^tu={\rm projections\ of\ trace}\ N\right)\\
C(U_{NN}^{N\times})&=&C^*_\times\left((u_{ij})_{i,j=1,\ldots,N}\Big|uu^*,u^*u,\bar{u}u^t,u^t\bar{u}={\rm projections\ of\ trace}\ N\right)
\end{eqnarray*}

Now observe that, in tensor product notation, the conditions at right are all of the form $(tr\otimes id)p=1$. Thus, $p$ must be follows, for the above conditions:
$$p=uu^*,u^*u,\bar{u}u^t,u^t\bar{u}$$

We therefore obtain that, for any faithful state $\varphi$, we have $(tr\otimes\varphi)(1-p)=0$. It follows from this that the following projections must be all equal to the identity:
$$p=uu^*,u^*u,\bar{u}u^t,u^t\bar{u}$$

But this leads to the conclusion in the statement.
\end{proof}

Regarding now the homogeneous space structure of $O_{MN}^{L\times},U_{MN}^{L\times}$, the situation here is a bit more complicated in the free case than in the classical case, due to a number of algebraic and analytic issues. We first have the following result:

\begin{proposition}
The spaces $U_{MN}^{L\times}$ have the following properties:
\begin{enumerate}
\item We have an action $U_M^\times\times U_N^\times\curvearrowright U_{MN}^{L\times}$, given by $u_{ij}\to\sum_{kl}u_{kl}\otimes a_{ki}\otimes b_{lj}^*$.

\item We have a map $U_M^\times\times U_N^\times\to U_{MN}^{L\times}$, given by $u_{ij}\to\sum_{r\leq L}a_{ri}\otimes b_{rj}^*$.
\end{enumerate}
Similar results hold for the spaces $O_{MN}^{L\times}$, with all the $*$ exponents removed.
\end{proposition}

\begin{proof}
In the classical case, consider the following action and quotient maps:
$$U_M\times U_N\curvearrowright U_{MN}^L\quad,\quad 
U_M\times U_N\to U_{MN}^L$$

The transposes of these two maps are as follows, where $J=(^1_0{\ }^0_0)$:
\begin{eqnarray*}
\varphi&\to&((U,A,B)\to\varphi(AUB^*))\\
\varphi&\to&((A,B)\to\varphi(AJB^*))
\end{eqnarray*}

But with $\varphi=u_{ij}$ we obtain precisely the formulae in the statement. The proof in the orthogonal case is similar. Regarding now the free case, the proof goes as follows:

\medskip

(1) Assuming $uu^*u=u$, consider the following variables:
$$U_{ij}=\sum_{kl}u_{kl}\otimes a_{ki}\otimes b_{lj}^*$$

We have then the following computation:
\begin{eqnarray*}
(UU^*U)_{ij}
&=&\sum_{pq}\sum_{klmnst}u_{kl}u_{mn}^*u_{st}\otimes a_{ki}a_{mq}^*a_{sq}\otimes b_{lp}^*b_{np}b_{tj}^*\\
&=&\sum_{klmt}u_{kl}u_{ml}^*u_{mt}\otimes a_{ki}\otimes b_{tj}^*\\
&=&\sum_{kt}u_{kt}\otimes a_{ki}\otimes b_{tj}^*\\
&=&U_{ij}
\end{eqnarray*}

Also, assuming that we have $\sum_{ij}u_{ij}u_{ij}^*=L$, we obtain:
\begin{eqnarray*}
\sum_{ij}U_{ij}U_{ij}^*
&=&\sum_{ij}\sum_{klst}u_{kl}u_{st}^*\otimes a_{ki}a_{si}^*\otimes b_{lj}^*b_{tj}\\
&=&\sum_{kl}u_{kl}u_{kl}^*\otimes1\otimes1\\
&=&L
\end{eqnarray*}

(2) Assuming $uu^*u=u$, consider the following variables:
$$V_{ij}=\sum_{r\leq L}a_{ri}\otimes b_{rj}^*$$

We have then the following computation:
\begin{eqnarray*}
(VV^*V)_{ij}
&=&\sum_{pq}\sum_{x,y,z\leq L}a_{xi}a_{yq}^*a_{zq}\otimes b_{xp}^*b_{yp}b_{zj}^*\\
&=&\sum_{x\leq L}a_{xi}\otimes b_{xj}^*\\
&=&V_{ij}
\end{eqnarray*}

Also, assuming that we have $\sum_{ij}u_{ij}u_{ij}^*=L$, we obtain:
\begin{eqnarray*}
\sum_{ij}V_{ij}V_{ij}^*
&=&\sum_{ij}\sum_{r,s\leq L}a_{ri}a_{si}^*\otimes b_{rj}^*b_{sj}\\
&=&\sum_{l\leq L}1\\
&=&L
\end{eqnarray*}

By removing all the $*$ exponents, we obtain as well the orthogonal results.
\end{proof}

Let us examine now the relation between the above maps. In the classical case, given a quotient space $X=G/H$, the associated action and quotient maps are given by:
$$\begin{cases}
a:X\times G\to X&:\quad (Hg,h)\to Hgh\\
p:G\to X&:\quad g\to Hg
\end{cases}$$

Thus we have $a(p(g),h)=p(gh)$. In our context, a similar result holds: 

\begin{theorem}
With $G=G_M\times G_N$ and $X=G_{MN}^L$, where $G_N=O_N^\times,U_N^\times$, we have
$$\xymatrix@R=15mm@C=30mm{
G\times G\ar[r]^m\ar[d]_{p\times id}&G\ar[d]^p\\
X\times G\ar[r]^a&X
}$$
where $a,p$ are the action map and the map constructed in Proposition 15.16.
\end{theorem}

\begin{proof}
At the level of the associated algebras of functions, we must prove that the following diagram commutes, where $\alpha,\gamma$ are morphisms of algebras induced by $a,p$:
$$\xymatrix@R=15mm@C=25mm{
C(X)\ar[r]^\alpha\ar[d]_\gamma&C(X\times G)\ar[d]^{\gamma\otimes id}\\
C(G)\ar[r]^\Delta&C(G\times G)
}$$

When going right, and then down, the composition is as follows:
\begin{eqnarray*}
(\gamma\otimes id)\alpha(u_{ij})
&=&(\gamma\otimes id)\sum_{kl}u_{kl}\otimes a_{ki}\otimes b_{lj}^*\\
&=&\sum_{kl}\sum_{r\leq L}a_{rk}\otimes b_{rl}^*\otimes a_{ki}\otimes b_{lj}^*
\end{eqnarray*}

On the other hand, when going down, and then right, the composition is as follows, where $F_{23}$ is the flip between the second and the third components:
\begin{eqnarray*}
\Delta\gamma(u_{ij})
&=&F_{23}(\Delta\otimes\Delta)\sum_{r\leq L}a_{ri}\otimes b_{rj}^*\\
&=&F_{23}\left(\sum_{r\leq L}\sum_{kl}a_{rk}\otimes a_{ki}\otimes b_{rl}^*\otimes b_{lj}^*\right)
\end{eqnarray*}

Thus the above diagram commutes indeed, and this gives the result.
\end{proof}

\section*{15c. Partial permutations}

As already mentioned in the beginning of the previous section, our above constructions using the quantum groups $O_N^\times,U_N^\times$ work in fact for any easy quantum group. In order to discuss this, we first need to work out the discrete extensions of our constructions. 

\bigskip

The starting point for all this is the semigroup $\widetilde{S}_N$ of partial permutations. This is a quite familiar object in combinatorics, defined as follows:

\index{partial permutation}
\index{semigroup of partial permutations}

\begin{definition}
$\widetilde{S}_N$ is the semigroup of partial permutations of $\{1\,\ldots,N\}$,
$$\widetilde{S}_N=\left\{\sigma:X\simeq Y\Big|X,Y\subset\{1,\ldots,N\}\right\}$$
with the usual composition operation, $\sigma'\sigma:\sigma^{-1}(X'\cap Y)\to\sigma'(X'\cap Y)$.
\end{definition}

Observe that $\widetilde{S}_N$ is not simplifiable, because the null permutation $\emptyset\in\widetilde{S}_N$, having the empty set as domain/range, satisfies $\emptyset\sigma=\sigma\emptyset=\emptyset$, for any $\sigma\in\widetilde{S}_N$. Observe also that $\widetilde{S}_N$ has a ``subinverse'' map, sending $\sigma:X\to Y$ to its usual inverse $\sigma^{-1}:Y\simeq X$.

\bigskip

A first interesting result about this semigroup $\widetilde{S}_N$, which shows that we are dealing here with some non-trivial combinatorics, is as follows:

\begin{proposition}
The number of partial permutations is given by
$$|\widetilde{S}_N|=\sum_{k=0}^Nk!\binom{N}{k}^2$$
that is, $1,2,7,34,209,\ldots\,$, and we have the cardinality estimate
$$|\widetilde{S}_N|\simeq N!\sqrt{\frac{\exp(4\sqrt{N}-1)}{4\pi\sqrt{N}}}$$
in the $N\to\infty$ limit.
\end{proposition}

\begin{proof}
The first assertion is clear, because in order to construct a partial permutation $\sigma:X\to Y$ we must choose an integer $k=|X|=|Y|$, then we must pick two subsets $X,Y\subset\{1,\ldots,N\}$ having cardinality $k$, and there are $\binom{N}{k}$ choices for each, and finally we must construct a bijection $\sigma:X\to Y$, and there are $k!$ choices here. As for the estimate, which is non-trivial, this is however something standard, and well-known.
\end{proof}

Another result, which is trivial, but quite fundamental, is as follows:

\begin{proposition}
We have a semigroup embedding $u:\widetilde{S}_N\subset M_N(0,1)$, defined by 
$$u_{ij}(\sigma)=
\begin{cases}
1&{\rm if}\ \sigma(j)=i\\
0&{\rm otherwise}
\end{cases}$$
whose image are the matrices having at most one nonzero entry, on each row and column.
\end{proposition}

\begin{proof}
This is trivial from definitions, with $u:\widetilde{S}_N\subset M_N(0,1)$ extending the standard embedding $u:S_N\subset M_N(0,1)$, that we have been heavily using, so far.
\end{proof}

Many other things can be said, about the partial permutations. Getting now to what we wanted to do, partial isometry spaces, let us formulate the following definition:

\index{partial permutation}

\begin{definition}
Associated to a partial permutation, $\sigma:I\simeq J$ with $I\subset\{1,\ldots,N\}$ and $J\subset\{1,\ldots,M\}$, is the real/complex partial isometry
$$T_\sigma:span\left(e_i\Big|i\in I\right)\to span\left(e_j\Big|j\in J\right)$$
given on the standard basis elements by $T_\sigma(e_i)=e_{\sigma(i)}$.
\end{definition}

Let $S_{MN}^L$ be the set of partial permutations $\sigma:I\simeq J$ as above, with range $I\subset\{1,\ldots,N\}$ and target $J\subset\{1,\ldots,M\}$, and with $L=|I|=|J|$. We have:

\begin{proposition}
The space of partial permutations signed by elements of $\mathbb Z_s$,
$$H_{MN}^{sL}=\left\{T(e_i)=w_ie_{\sigma(i)}\Big|\sigma\in S_{MN}^L,w_i\in\mathbb Z_s\right\}$$
is isomorphic to the quotient space 
$$(H_M^s\times H_N^s)/(H_L^s\times H_{M-L}^s\times H_{N-L}^s)$$
via a standard isomorphism.
\end{proposition}

\begin{proof}
This follows by adapting the computations in the proof of Proposition 15.12 and Theorem 15.13. Indeed, we have an action map as follows, which is transitive:
$$H_M^s\times H_N^s\to H_{MN}^{sL}\quad,\quad 
(A,B)U=AUB^*$$

Consider now the following point, as we did before in the continuous case:
$$U=\begin{pmatrix}1&0\\0&0\end{pmatrix}$$

The stabilizer of this point follows to be the following group:
$$H_L^s\times H_{M-L}^s\times H_{N-L}^s$$

To be more precise, this group is embedded via:
$$(x,a,b)\to\left[\begin{pmatrix}x&0\\0&a\end{pmatrix},\begin{pmatrix}x&0\\0&b\end{pmatrix}\right]$$

But this gives the result.
\end{proof}

In the free case now, the idea is similar, by using inspiration from the construction of the quantum group $H_N^{s+}=\mathbb Z_s\wr_*S_N^+$. The result here is as follows:

\index{free partial permutation}

\begin{proposition}
The compact quantum space $H_{MN}^{sL+}$ associated to the algebra
$$C(H_{MN}^{sL+})=C(U_{MN}^{L+})\Big/\left<u_{ij}u_{ij}^*=u_{ij}^*u_{ij}=p_{ij}={\rm projections},u_{ij}^s=p_{ij}\right>$$
has an action map, and is the target of a quotient map, as in Theorem 15.17.
\end{proposition}

\begin{proof}
We must show that if the variables $u_{ij}$ satisfy the relations in the statement, then these relations are satisfied as well for the following variables: 
$$U_{ij}=\sum_{kl}u_{kl}\otimes a_{ki}\otimes b_{lj}^*\quad,\quad 
V_{ij}=\sum_{r\leq L}a_{ri}\otimes b_{rj}^*$$

We use the fact that the standard coordinates $a_{ij},b_{ij}$ on the quantum groups $H_M^{s+},H_N^{s+}$ satisfy the following relations, for any $x\neq y$ on the same row or column of $a,b$:
$$xy=xy^*=0$$
 
We obtain, by using these relations, the following formula:
$$U_{ij}U_{ij}^*
=\sum_{klmn}u_{kl}u_{mn}^*\otimes a_{ki}a_{mi}^*\otimes b_{lj}^*b_{mj}
=\sum_{kl}u_{kl}u_{kl}^*\otimes a_{ki}a_{ki}^*\otimes b_{lj}^*b_{lj}$$

On the other hand, we have as well the following formula:
$$V_{ij}V_{ij}^*
=\sum_{r,t\leq L}a_{ri}a_{ti}^*\otimes b_{rj}^*b_{tj}
=\sum_{r\leq L}a_{ri}a_{ri}^*\otimes b_{rj}^*b_{rj}$$

In terms of the projections $x_{ij}=a_{ij}a_{ij}^*$, $y_{ij}=b_{ij}b_{ij}^*$, $p_{ij}=u_{ij}u_{ij}^*$, we have:
$$U_{ij}U_{ij}^*=\sum_{kl}p_{kl}\otimes x_{ki}\otimes y_{lj}\quad,\quad 
V_{ij}V_{ij}^*=\sum_{r\leq L}x_{ri}\otimes y_{rj}$$

By repeating the computation, we conclude that these elements are projections. Also, a similar computation shows that $U_{ij}^*U_{ij},V_{ij}^*V_{ij}$ are given by the same formulae. Finally, once again by using the relations of type $xy=xy^*=0$, we have:
$$U_{ij}^s
=\sum_{k_rl_r}u_{k_1l_1}\ldots u_{k_sl_s}\otimes a_{k_1i}\ldots a_{k_si}\otimes b_{l_1j}^*\ldots b_{l_sj}^*
=\sum_{kl}u_{kl}^s\otimes a_{ki}^s\otimes(b_{lj}^*)^s$$

On the other hand, we have as well the following formula:
$$V_{ij}^s
=\sum_{r_l\leq L}a_{r_1i}\ldots a_{r_si}\otimes b_{r_1j}^*\ldots b_{r_sj}^*
=\sum_{r\leq L}a_{ri}^s\otimes(b_{rj}^*)^s$$

Thus the conditions of type $u_{ij}^s=p_{ij}$ are satisfied as well, and we are done.
\end{proof}

Let us discuss now the general case. We have the following result:

\begin{proposition}
The various spaces $G_{MN}^L$ constructed so far appear by imposing to the standard coordinates of $U_{MN}^{L+}$ the relations
$$\sum_{i_1\ldots i_s}\sum_{j_1\ldots j_s}\delta_\pi(i)\delta_\sigma(j)u_{i_1j_1}^{e_1}\ldots u_{i_sj_s}^{e_s}=L^{|\pi\vee\sigma|}$$
with $s=(e_1,\ldots,e_s)$ ranging over all the colored integers, and with $\pi,\sigma\in D(0,s)$.
\end{proposition}

\begin{proof}
According to the various constructions above, the relations defining the quantum space $G_{MN}^L$ can be written as follows, with $\sigma$ ranging over a family of generators, with no upper legs, of the corresponding category of partitions $D$:
$$\sum_{j_1\ldots j_s}\delta_\sigma(j)u_{i_1j_1}^{e_1}\ldots u_{i_sj_s}^{e_s}=\delta_\sigma(i)$$

We therefore obtain the relations in the statement, as follows:
\begin{eqnarray*}
\sum_{i_1\ldots i_s}\sum_{j_1\ldots j_s}\delta_\pi(i)\delta_\sigma(j)u_{i_1j_1}^{e_1}\ldots u_{i_sj_s}^{e_s}
&=&\sum_{i_1\ldots i_s}\delta_\pi(i)\sum_{j_1\ldots j_s}\delta_\sigma(j)u_{i_1j_1}^{e_1}\ldots u_{i_sj_s}^{e_s}\\
&=&\sum_{i_1\ldots i_s}\delta_\pi(i)\delta_\sigma(i)\\
&=&L^{|\pi\vee\sigma|}
\end{eqnarray*}

As for the converse, this follows by using the relations in the statement, by keeping $\pi$ fixed, and by making $\sigma$ vary over all the partitions in the category.
\end{proof}

In the general case now, where $G=(G_N)$ is an arbitary uniform easy quantum group, we can construct spaces $G_{MN}^L$ by using the above relations, and we have:

\begin{theorem}
The spaces $G_{MN}^L\subset U_{MN}^{L+}$ constructed by imposing the relations 
$$\sum_{i_1\ldots i_s}\sum_{j_1\ldots j_s}\delta_\pi(i)\delta_\sigma(j)u_{i_1j_1}^{e_1}\ldots u_{i_sj_s}^{e_s}=L^{|\pi\vee\sigma|}$$
with $\pi,\sigma$ ranging over all the partitions in the associated category, having no upper legs, are subject to an action map/quotient map diagram, as in Theorem 15.17.
\end{theorem}

\begin{proof}
We proceed as in the proof of Proposition 15.23. We must prove that, if the variables $u_{ij}$ satisfy the relations in the statement, then so do the following variables:
$$U_{ij}=\sum_{kl}u_{kl}\otimes a_{ki}\otimes b_{lj}^*\quad,\quad 
V_{ij}=\sum_{r\leq L}a_{ri}\otimes b_{rj}^*$$

Regarding the variables $U_{ij}$, the computation here goes as follows:
\begin{eqnarray*}
&&\sum_{i_1\ldots i_s}\sum_{j_1\ldots j_s}\delta_\pi(i)\delta_\sigma(j)U_{i_1j_1}^{e_1}\ldots U_{i_sj_s}^{e_s}\\
&=&\sum_{i_1\ldots i_s}\sum_{j_1\ldots j_s}\sum_{k_1\ldots k_s}\sum_{l_1\ldots l_s}u_{k_1l_1}^{e_1}\ldots u_{k_sl_s}^{e_s}\otimes \delta_\pi(i)\delta_\sigma(j)a_{k_1i_1}^{e_1}\ldots a_{k_si_s}^{e_s}\otimes(b_{l_sj_s}^{e_s}\ldots b_{l_1j_1}^{e_1})^*\\
&=&\sum_{k_1\ldots k_s}\sum_{l_1\ldots l_s}\delta_\pi(k)\delta_\sigma(l)u_{k_1l_1}^{e_1}\ldots u_{k_sl_s}^{e_s}\\
&=&L^{|\pi\vee\sigma|}
\end{eqnarray*}

For the variables $V_{ij}$ the proof is similar, as follows:
\begin{eqnarray*}
&&\sum_{i_1\ldots i_s}\sum_{j_1\ldots j_s}\delta_\pi(i)\delta_\sigma(j)V_{i_1j_1}^{e_1}\ldots V_{i_sj_s}^{e_s}\\
&=&\sum_{i_1\ldots i_s}\sum_{j_1\ldots j_s}\sum_{l_1,\ldots,l_s\leq L}\delta_\pi(i)\delta_\sigma(j)a_{l_1i_1}^{e_1}\ldots a_{l_si_s}^{e_s}\otimes(b_{l_sj_s}^{e_s}\ldots b_{l_1j_1}^{e_1})^*\\
&=&\sum_{l_1,\ldots,l_s\leq L}\delta_\pi(l)\delta_\sigma(l)\\
&=&L^{|\pi\vee\sigma|}
\end{eqnarray*}

Thus we have constructed an action map, and a quotient map, as in Proposition 15.23, and the commutation of the diagram in Theorem 15.17 is then trivial.
\end{proof}

Summarizing, our partial isometry space theory works well for any easy quantum group, and we have full extensions of our algebraic results regarding the spheres.

\section*{15d. Integration results}

Getting now ot analytic matters, again in analogy with what we did before for the spheres, let us discuss now the integration over the spaces $G_{MN}^L$. We first have:

\index{Haar integration}

\begin{definition}
The integration functional of $G_{MN}^L$ is the composition
$$\int_{G_{MN}^L}:C(G_{MN}^L)\to C(G_M\times G_N)\to\mathbb C$$
of the representation $u_{ij}\to\sum_{r\leq L}a_{ri}\otimes b_{rj}^*$ with the Haar functional of $G_M\times G_N$.
\end{definition}

Observe that in the case $L=M=N$ we obtain the integration over $G_N$. Also, at $L=M=1$, or at $L=N=1$, we obtain the integration over the sphere. 

\bigskip

In the general case now, we first have the following result:

\index{left invariance}
\index{right invariance}

\begin{proposition}
The integration functional of $G_{MN}^L$ has the invariance property 
$$\left(\int_{G_{MN}^L}\!\otimes\ id\right)\alpha(x)=\int_{G_{MN}^L}x$$
with respect to the coaction map $\alpha(u_{ij})=\sum_{kl}u_{kl}\otimes a_{ki}\otimes b_{lj}^*$.
\end{proposition}

\begin{proof}
We restrict the attention to the orthogonal case, the proof in the unitary case being similar. We must check the following formula:
$$\left(\int_{G_{MN}^L}\!\otimes\ id\right)\alpha(u_{i_1j_1}\ldots u_{i_sj_s})=\int_{G_{MN}^L}u_{i_1j_1}\ldots u_{i_sj_s}$$

Let us compute the left term. This is given by:
\begin{eqnarray*}
X
&=&\left(\int_{G_{MN}^L}\!\otimes\ id\right)\sum_{k_xl_x}u_{k_1l_1}\ldots u_{k_sl_s}\otimes a_{k_1i_1}\ldots a_{k_si_s}\otimes b_{l_1j_1}^*\ldots b_{l_sj_s}^*\\
&=&\sum_{k_xl_x}\sum_{r_x\leq L}a_{k_1i_1}\ldots a_{k_si_s}\otimes b_{l_1j_1}^*\ldots b_{l_sj_s}^*\int_{G_M}a_{r_1k_1}\ldots a_{r_sk_s}\int_{G_N}b_{r_1l_1}^*\ldots b_{r_sl_s}^*\\
&=&\sum_{r_x\leq L}\sum_{k_x}a_{k_1i_1}\ldots a_{k_si_s}\int_{G_M}a_{r_1k_1}\ldots a_{r_sk_s}
\otimes\sum_{l_x}b_{l_1j_1}^*\ldots b_{l_sj_s}^*\int_{G_N}b_{r_1l_1}^*\ldots b_{r_sl_s}^*
\end{eqnarray*}

By using now the invariance property of the Haar functionals of $G_M,G_N$, we obtain:
\begin{eqnarray*}
X
&=&\sum_{r_x\leq L}\left(\int_{G_M}\!\otimes\ id\right)\Delta(a_{r_1i_1}\ldots a_{r_si_s})
\otimes\left(\int_{G_N}\!\otimes\ id\right)\Delta(b_{r_1j_1}^*\ldots b_{r_sj_s}^*)\\
&=&\sum_{r_x\leq L}\int_{G_M}a_{r_1i_1}\ldots a_{r_si_s}\int_{G_N}b_{r_1j_1}^*\ldots b_{r_sj_s}^*\\
&=&\left(\int_{G_M}\otimes\int_{G_N}\right)\sum_{r_x\leq L}a_{r_1i_1}\ldots a_{r_si_s}\otimes b_{r_1j_1}^*\ldots b_{r_sj_s}^*
\end{eqnarray*}

But this gives the formula in the statement, and we are done.
\end{proof}

We will prove now that the above functional is in fact the unique positive unital invariant trace on $C(G_{MN}^L)$. For this purpose, we will need the Weingarten formula:

\index{Weingarten formula}

\begin{theorem}
We have the Weingarten type formula
$$\int_{G_{MN}^L}u_{i_1j_1}\ldots u_{i_sj_s}=\sum_{\pi\sigma\tau\nu}L^{|\pi\vee\tau|}\delta_\sigma(i)\delta_\nu(j)W_{sM}(\pi,\sigma)W_{sN}(\tau,\nu)$$
where the matrices on the right are given by $W_{sM}=G_{sM}^{-1}$, with $G_{sM}(\pi,\sigma)=M^{|\pi\vee\sigma|}$.
\end{theorem}

\begin{proof}
We make use of the usual quantum group Weingarten formula, that we know from chapter 11. By using this formula for $G_M,G_N$, we obtain:
\begin{eqnarray*}
\int_{G_{MN}^L}u_{i_1j_1}\ldots u_{i_sj_s}
&=&\sum_{l_1\ldots l_s\leq L}\int_{G_M}a_{l_1i_1}\ldots a_{l_si_s}\int_{G_N}b_{l_1j_1}^*\ldots b_{l_sj_s}^*\\
&=&\sum_{l_1\ldots l_s\leq L}\sum_{\pi\sigma}\delta_\pi(l)\delta_\sigma(i)W_{sM}(\pi,\sigma)\sum_{\tau\nu}\delta_\tau(l)\delta_\nu(j)W_{sN}(\tau,\nu)\\
&=&\sum_{\pi\sigma\tau\nu}\left(\sum_{l_1\ldots l_s\leq L}\delta_\pi(l)\delta_\tau(l)\right)\delta_\sigma(i)\delta_\nu(j)W_{sM}(\pi,\sigma)W_{sN}(\tau,\nu)
\end{eqnarray*}

The coefficient being $L^{|\pi\vee\tau|}$, we obtain the formula in the statement.
\end{proof}

We can now derive an abstract characterization of the integration, as follows:

\begin{theorem}
The integration of $G_{MN}^L$ is the unique positive unital trace 
$$C(G_{MN}^L)\to\mathbb C$$
which is invariant under the action of the quantum group $G_M\times G_N$.
\end{theorem}

\begin{proof}
We use a standard method, from \cite{bgo}, that we already met in the context of the spheres, the point being to show that we have the following ergodicity formula: 
$$\left(id\otimes\int_{G_M}\otimes\int_{G_N}\right)\alpha(x)=\int_{G_{MN}^L}x$$

We restrict the attention to the orthogonal case, the proof in the unitary case being similar. We must verify that the following holds:
$$\left(id\otimes\int_{G_M}\otimes\int_{G_N}\right)\alpha(u_{i_1j_1}\ldots u_{i_sj_s})=\int_{G_{MN}^L}u_{i_1j_1}\ldots u_{i_sj_s}$$

By using the Weingarten formula, the left term can be written as follows:
\begin{eqnarray*}
X
&=&\sum_{k_1\ldots k_s}\sum_{l_1\ldots l_s}u_{k_1l_1}\ldots u_{k_sl_s}\int_{G_M}a_{k_1i_1}\ldots a_{k_si_s}\int_{G_N}b_{l_1j_1}^*\ldots b_{l_sj_s}^*\\
&=&\sum_{k_1\ldots k_s}\sum_{l_1\ldots l_s}u_{k_1l_1}\ldots u_{k_sl_s}\sum_{\pi\sigma}\delta_\pi(k)\delta_\sigma(i)W_{sM}(\pi,\sigma)\sum_{\tau\nu}\delta_\tau(l)\delta_\nu(j)W_{sN}(\tau,\nu)\\
&=&\sum_{\pi\sigma\tau\nu}\delta_\sigma(i)\delta_\nu(j)W_{sM}(\pi,\sigma)W_{sN}(\tau,\nu)\sum_{k_1\ldots k_s}\sum_{l_1\ldots l_s}\delta_\pi(k)\delta_\tau(l)u_{k_1l_1}\ldots u_{k_sl_s}
\end{eqnarray*}

By using now the summation formula in Theorem 15.28, we obtain:
$$X=\sum_{\pi\sigma\tau\nu}L^{|\pi\vee\tau|}\delta_\sigma(i)\delta_\nu(j)W_{sM}(\pi,\sigma)W_{sN}(\tau,\nu)$$

Now by comparing with the Weingarten formula for $G_{MN}^L$, this proves our claim. Assume now that $\tau:C(G_{MN}^L)\to\mathbb C$ satisfies the invariance condition. We have:
\begin{eqnarray*}
\tau\left(id\otimes\int_{G_M}\otimes\int_{G_N}\right)\alpha(x)
&=&\left(\tau\otimes\int_{G_M}\otimes\int_{G_N}\right)\alpha(x)\\
&=&\left(\int_{G_M}\otimes\int_{G_N}\right)(\tau\otimes id)\alpha(x)\\
&=&\left(\int_{G_M}\otimes\int_{G_N}\right)(\tau(x)1)\\
&=&\tau(x)
\end{eqnarray*}

On the other hand, according to the formula established above, we have as well:
\begin{eqnarray*}
\tau\left(id\otimes\int_{G_M}\otimes\int_{G_N}\right)\alpha(x)
&=&\tau(tr(x)1)\\
&=&tr(x)
\end{eqnarray*}

Thus we obtain $\tau=tr$, and this finishes the proof.
\end{proof}

As a main application of the above results, we have:

\index{overlapping coordinates}
\index{sum of coordinates}

\begin{proposition}
For a sum of coordinates of the following type,
$$\chi_E=\sum_{(ij)\in E}u_{ij}$$
with the coordinates not overlapping on rows and columns, we have
$$\int_{G_{MN}^L}\chi_E^s=\sum_{\pi\sigma\tau\nu}K^{|\pi\vee\tau|}L^{|\sigma\vee\nu|}W_{sM}(\pi,\sigma)W_{sN}(\tau,\nu)$$
where $K=|E|$ is the cardinality of the indexing set.
\end{proposition}

\begin{proof}
With $K=|E|$, we can write $E=\{(\alpha(i),\beta(i))\}$, for certain embeddings:
$$\alpha:\{1,\ldots,K\}\subset\{1,\ldots,M\}$$
$$\beta:\{1,\ldots,K\}\subset\{1,\ldots,N\}$$

In terms of these maps $\alpha,\beta$, the moment in the statement is given by:
$$M_s=\int_{G_{MN}^L}\left(\sum_{i\leq K}u_{\alpha(i)\beta(i)}\right)^s$$

By using the Weingarten formula, we can write this quantity as follows:
\begin{eqnarray*}
&&M_s\\
&=&\int_{G_{MN}^L}\sum_{i_1\ldots i_s\leq K}u_{\alpha(i_1)\beta(i_1)}\ldots u_{\alpha(i_s)\beta(i_s)}\\
&=&\sum_{i_1\ldots i_s\leq K}\sum_{\pi\sigma\tau\nu}L^{|\sigma\vee\nu|}\delta_\pi(\alpha(i_1),\ldots,\alpha(i_s))\delta_\tau(\beta(i_1),\ldots,\beta(i_s))W_{sM}(\pi,\sigma)W_{sN}(\tau,\nu)\\
&=&\sum_{\pi\sigma\tau\nu}\left(\sum_{i_1\ldots i_s\leq K}\delta_\pi(i)\delta_\tau(i)\right)L^{|\sigma\vee\nu|}W_{sM}(\pi,\sigma)W_{sN}(\tau,\nu)
\end{eqnarray*}

But, as explained before, in the proof of Theorem 15.28, the coefficient on the left in the last formula is $C=K^{|\pi\vee\tau|}$. We therefore obtain the formula in the statement.
\end{proof}

At a more concrete level now, we have the following conceptual result, making the link with the Bercovici-Pata bijection \cite{bpa}:

\index{Bercovici-Pata bijection}
\index{non-overlapping coordinates}

\begin{theorem}
In the context of the liberation operations 
$$O_{MN}^L\to O_{MN}^{L+}\quad,\quad 
U_{MN}^L\to U_{MN}^{L+}\quad,\quad 
H_{MN}^{sL}\to H_{MN}^{sL+}$$ 
the laws of the sums of non-overlapping coordinates,
$$\chi_E=\sum_{(ij)\in E}u_{ij}$$
are in Bercovici-Pata bijection, in the 
$$|E|=\kappa N,L=\lambda N,M=\mu N$$
regime and $N\to\infty$ limit.
\end{theorem}

\begin{proof}
We use various formulae from \cite{bbc}, \cite{bsp}. According to Proposition 15.30, in terms of $K=|E|$, the moments of the variables in the statement are given by:
$$M_s=\sum_{\pi\sigma\tau\nu}K^{|\pi\vee\tau|}L^{|\sigma\vee\nu|}W_{sM}(\pi,\sigma)W_{sN}(\tau,\nu)$$

We use now two standard facts, from \cite{bbc} and related papers, namely the fact that in the $N\to\infty$ limit the Weingarten matrix $W_{sN}$ is concentrated on the diagonal, and the fact that we have an inequality as follows, with equality precisely when $\pi=\sigma$:
$$|\pi\vee\sigma|\leq\frac{|\pi|+|\sigma|}{2}$$

Indeed, with these two ingredients in hand, we can now look in detail at what happens to our moment $M_s$ in the regime from the statement, namely:
$$K=\kappa N,L=\lambda N,M=\mu N,N\to\infty$$

In this regime, we obtain the following estimate:
\begin{eqnarray*}
M_s
&\simeq&\sum_{\pi\tau}K^{|\pi\vee\tau|}L^{|\pi\vee\tau|}M^{-|\pi|}N^{-|\tau|}\\
&\simeq&\sum_\pi K^{|\pi|}L^{|\pi|}M^{-|\pi|}N^{-|\pi|}\\
&=&\sum_\pi\left(\frac{\kappa\lambda}{\mu}\right)^{|\pi|}
\end{eqnarray*}

In order to interpret this formula, we use general theory from \cite{bbc}, \cite{bsp}:

\medskip

(1) For $G_N=O_N,O_N^+$, the above variables $\chi_E$ follow to be asymptotically Gaussian/semicircular, of parameter $\frac{\kappa\lambda}{\mu}$, and hence in Bercovici-Pata bijection.

\medskip

(2) For $G_N=U_N,U_N^+$ the situation is similar, with $\chi_E$ being asymptotically complex Gaussian/circular, of parameter $\frac{\kappa\lambda}{\mu}$, and in Bercovici-Pata bijection. 

\medskip

(3) Finally, for $G_N=H_N^s,H_N^{s+}$, the variables $\chi_E$ are asymptotically Bessel/free Bessel of parameter $\frac{\kappa\lambda}{\mu}$, and once again in Bercovici-Pata bijection.  
\end{proof}

There are several possible extensions of the above results, for instance by using twisting operations as well. We refer here to \cite{bgo} and related papers.

\section*{15e. Exercises}

We had a quite technical chapter here, and as exercises on this, we have:

\begin{exercise}
Learn more about the various quantum homogeneous spaces.
\end{exercise}

\begin{exercise}
Work out the construction of twisted partial isometry spaces.
\end{exercise}

\begin{exercise}
Learn more about the semigroup $\widetilde{S}_N$, and its properties.
\end{exercise}

\begin{exercise}
Construct and study the free quantum semigroup $\widetilde{S_N^+}$.
\end{exercise}

\begin{exercise}
Learn more about the Bercovici-Pata bijection, in general.
\end{exercise}

\begin{exercise}
Learn about affine homogeneous spaces, generalizing the above.
\end{exercise}

As bonus exercise, try to develop some abstract free algebraic geometry.

\chapter{Threefold way}

\section*{16a. Rotation groups}

Time to end this book, and we have chosen to talk about something quite mysterious, in connection with virtually everything done since we started in chapter 1.

\bigskip

Hang on, this is what we would like to talk about:

\begin{principle}[Threefold way]
There are three main geometries, or master ground fields $F$ if you prefer, namely:
\begin{enumerate}
\item Real.

\item Complex.

\item Free.
\end{enumerate}
\end{principle}

Obviously, this is not the sort of statement that you will meet every day, so this is definitely worth a detailed look, from all the possible perspectives. Our purpose in this chapter will be to do so, thoroughly examine this threefold way principle.

\bigskip

As a first disclaimer, Principle 16.1 is not something of arithmetic nature, but rather something quantum, coming from our quantum algebra experience.

\bigskip

Getting started now, we will start with what we best know to do, namely quantum groups. Our purpose will be to discuss the isomorphism $PO_N^+=PU_N^+$, which is something that we have already met, and which brings heavy evidence for Principle 16.1.

\bigskip

For this purpose, let us first discuss the passage from real to complex, in the group setting. The passage $O_N\to U_N$ cannot be understood directly. In order to understand this, we must pass through the corresponding Lie algebras, as follows:

\begin{theorem}
The passage $O_N\to U_N$ appears via a Lie algebra complexification,
$$O_N\to\mathfrak o_N\to\mathfrak u_n\to U_N$$
with the Lie algebra $\mathfrak u_N$ being a complexification of the Lie algebra $\mathfrak o_N$.
\end{theorem}

\begin{proof}
This is something rather philosophical, the idea being as follows:

\medskip

(1) The orthogonal and unitary groups $O_N,N_N$ are both Lie groups, and the corresponding Lie algebras $\mathfrak o_N,\mathfrak u_N$ can be computed by differentiating the equations defining $O_N,U_N$, with the conclusion being as follows:
$$\mathfrak o_N=\left\{ A\in M_N(\mathbb R)\Big|A^t=-A\right\}$$
$$\mathfrak u_N=\left\{ B\in M_N(\mathbb C)\Big|B^*=-B\right\}$$

(2) This was for the correspondences $O_N\to\mathfrak o_N$ and $U_N\to\mathfrak u_N$. In the other sense, the correspondences $\mathfrak o_N\to O_N$ and $\mathfrak u_N\to U_N$ appear by exponentiation, the result here stating that, around 1, the orthogonal matrices can be written as $U=e^A$, with $A\in\mathfrak o_N$, and the unitary matrices can be written as $U=e^B$, with $B\in\mathfrak u_N$. 

\bigskip

(3) In view of all this, in order to understand the passage $O_N\to U_N$ it is enough to understand the passage $\mathfrak o_N\to\mathfrak u_N$. But, in view of the above explicit formulae for $\mathfrak o_N,\mathfrak u_N$, this is basically an elementary linear algebra problem. Indeed, let us pick an arbitrary matrix $B\in M_N(\mathbb C)$, and write it as follows, with $A,C\in M_N(\mathbb R)$:
$$B=A+iC$$

In terms of $A,C$, the equation $B^*=-B$ defining the Lie algebra $\mathfrak u_N$ reads:
$$A^t=-A$$
$$C^t=C$$

(4) As a first observation, we must have $A\in\mathfrak o_N$. Regarding now $C$, let us decompose it as follows, with $D$ being its diagonal, and $C'$ being the remainder:
$$C=D+C'$$

The remainder $C'$ being symmetric with 0 on the diagonal, by switching all the signs below the main diagonal we obtain a certain matrix $C'_-\in\mathfrak o_N$. Thus, we have decomposed $B\in\mathfrak u_N$ as follows, with $A,C'\in\mathfrak o_N$, and with $D\in M_N(\mathbb R)$ being diagonal:
$$B=A+iD+iC'_-$$

(5) As a conclusion now, we have shown that we have a direct sum decomposition of real linear spaces as follows, with $\Delta\subset M_N(\mathbb R)$ being the diagonal matrices:
$$\mathfrak u_N\simeq\mathfrak o_N\oplus\Delta\oplus\mathfrak o_N$$

Thus, we can stop our study here, and say that we have reached the conclusion in the statement, namely that $\mathfrak u_N$ appears as a ``complexification'' of $\mathfrak o_N$.
\end{proof}

In the free case now, the situation is much simpler, and we have:

\begin{theorem}
The passage $O_N^+\to U_N^+$ appears via free complexification,
$$U_N^+=\widetilde{O_N^+}$$
where the free complexification of a pair $(G,u)$ is by definition the pair $(\widetilde{G},\widetilde{u})$ with
$$C(\widetilde{G})=<zu_{ij}>\subset C(\mathbb T)*C(G)\quad,\quad 
\widetilde{u}=zu$$
where $z\in C(\mathbb T)$ is the standard generator, given by $x\to x$ for any $x\in\mathbb T$.
\end{theorem}

\begin{proof}
We have embeddings as follows, with the first one coming by using the counit, and with the second one coming from the universality property of $U_N^+$:
$$O_N^+
\subset\widetilde{O_N^+}
\subset U_N^+$$

We must prove that the embedding on the right is an isomorphism. Let us recall that if we denote by $v,u$ the fundamental corepresentations of $O_N^+,U_N^+$, then we have:
$$Fix(v^{\otimes k})=span\left(\xi_\pi\Big|\pi\in NC_2(k)\right)$$
$$Fix(u^{\otimes k})=span\left(\xi_\pi\Big|\pi\in\mathcal{NC}_2(k)\right)$$

Moreover, the above vectors $\xi_\pi$ are known to be linearly independent at $N\geq2$, and so the above results provide us with bases, and we obtain:
$$\dim(Fix(v^{\otimes k}))=|NC_2(k)|\quad,\quad 
\dim(Fix(u^{\otimes k}))=|\mathcal{NC}_2(k)|$$ 

Now since integrating the character of a corepresentation amounts in counting the fixed points, the above two formulae can be rewritten as follows:
$$\int_{O_N^+}\chi_v^k=|NC_2(k)|\quad,\quad 
\int_{U_N^+}\chi_u^k=|\mathcal{NC}_2(k)|$$ 

But this shows, via standard free probability theory, that $\chi_v$ must follow the Winger semicircle law $\gamma_1$, and that $\chi_u$ must follow the Voiculescu circular law $\Gamma_1$:
$$\chi_v\sim\gamma_1\quad,\quad 
\chi_u\sim\Gamma_1$$

On the other hand, when freely multiplying a semicircular variable by a Haar unitary we obtain a circular variable. Thus, the main character of $\widetilde{O_N^+}$ is circular:
$$\chi_{zv}\sim\Gamma_1$$

Now by forgetting about circular variables and free probability, the conclusion is that the inclusion $\widetilde{O_N^+}\subset U_N^+$ preserves the law of the main character:
$$law(\chi_{zv})=law(u)$$

Thus by Peter-Weyl we obtain that the inclusion $\widetilde{O_N^+}\subset U_N^+$ must be an isomorphism, modulo the usual equivalence relation for quantum groups.
\end{proof}

As a main consequence of the above result, which is of interest for us, we have:

\index{projective unitary group}
\index{projective quantum group}

\begin{theorem}
We have an identification as follows,
$$PO_N^+=PU_N^+$$
modulo the usual equivalence relation for compact quantum groups.
\end{theorem}

\begin{proof}
As before, we have several proofs for this result, as follows:

\medskip

(1) This follows from Theorem 16.3, because we have:
$$PU_N^+=P\widetilde{O_N^+}=PO_N^+$$

(2) We can deduce this as well directly. With notations as before, we have:
$$Hom\left((v\otimes v)^k,(v\otimes v)^l\right)=span\left(T_\pi\Big|\pi\in NC_2((\circ\bullet)^k,(\circ\bullet)^l)\right)$$
$$Hom\left((u\otimes\bar{u})^k,(u\otimes\bar{u})^l\right)=span\left(T_\pi\Big|\pi\in \mathcal{NC}_2((\circ\bullet)^k,(\circ\bullet)^l)\right)$$

The sets on the right being equal, we conclude that the inclusion $PO_N^+\subset PU_N^+$ preserves the corresponding Tannakian categories, and so must be an isomorphism.
\end{proof}

Many other things can be said, along these lines, notably with an axiomatization of the notion of projective easiness, and some related classification results, for instance for the projective free quantum groups. We will be back to this, later in this chapter.

\section*{16b. Monomial spheres}

Getting now into noncommutative geometry, and again with Principle 16.1 in mind, looking back at the definition of the spheres that we have, and at the precise relations between the coordinates, we are led into the following notion:

\index{monomial sphere}
\index{free sphere}
\index{noncommutative sphere}

\begin{definition}
A monomial sphere is a subset $S\subset S^{N-1}_{\mathbb C,+}$ obtained via relations
$$x_{i_1}^{e_1}\ldots x_{i_k}^{e_k}=x_{i_{\sigma(1)}}^{f_1}\ldots x_{i_{\sigma(k)}}^{f_k}\quad,\quad\forall (i_1,\ldots,i_k)\in\{1,\ldots,N\}^k$$
with $\sigma\in S_k$ being certain permutations, and with $e_r,f_r\in\{1,*\}$ being certain exponents.
\end{definition}

This definition is quite broad, and we have for instance as example the sphere $S^{N-1}_{\mathbb C,\times}$ coming from the relations $ab^*c=cb^*a$, corresponding to the following diagram: 
$$
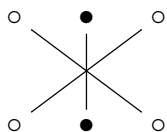
\xymatrix@R=10mm@C=5mm{\circ\ar@{-}[drr]&\bullet\ar@{-}[d]&\circ\ar@{-}[dll]\\\circ&\bullet&\circ}$$ 

In view of these difficulties, we will restrict now the attention to the real case. Let us first recall that we have the following fundamental result, dealing with the real case:

\index{real geometry}
\index{easy geometry}

\begin{theorem}
There are exactly $3$ real easy geometries, namely
$$\mathbb R^N\subset\mathbb R^N_*\subset\mathbb R^N_+$$
coming from $P_2\supset P_2^*\supset NC_2$, whose associated spheres are
$$S^{N-1}_\mathbb R\subset S^{N-1}_{\mathbb R,*}\subset S^{N-1}_{\mathbb R,+}$$
and whose tori, unitary and reflection groups are given by similar formulae.
\end{theorem}

\begin{proof}
This is something that we know well, coming from the fact that $G=O_N^*$ is the unique intermediate easy quantum group $O_N\subset G\subset O_N^+$.
\end{proof}

Let us focus now on the spheres, and try to better understand their ``easiness'' property, with results getting beyond what has been done above, in the general easy context. That is, our objects of interest in what follows will be the 3 real spheres, namely:
$$S^{N-1}_\mathbb R\subset S^{N-1}_{\mathbb R,*}\subset S^{N-1}_{\mathbb R,+}$$

Our purpose in what follows we will be that of proving that these spheres are the only monomial ones. In order to best talk about monomiality, in the present real case, it is convenient to introduce the following infinite group:
$$S_\infty=\bigcup_{k\geq0}S_k$$

To be more precise, this group appears by definition as an inductive limit, with the inclusions $S_k\subset S_{k+1}$ that we use being given by: 
$$\sigma\in S_k\implies\sigma(k+1)=k+1$$

In terms of $S_\infty$, the definition of the monomial spheres reformulates as follows:

\index{monomial sphere}

\begin{proposition}
The monomial spheres are the algebraic manifolds $S\subset S^{N-1}_{\mathbb R,+}$ obtained via relations of type
$$x_{i_1}\ldots x_{i_k}=x_{i_{\sigma(1)}}\ldots x_{i_{\sigma(k)}},\ \forall (i_1,\ldots,i_k)\in\{1,\ldots,N\}^k$$
associated to certain elements $\sigma\in S_\infty$, where $k\in\mathbb N$ is such that $\sigma\in S_k$. 
\end{proposition}

\begin{proof}
We must prove that the relations $x_{i_1}\ldots x_{i_k}=x_{i_{\sigma(1)}}\ldots x_{i_{\sigma(k)}}$ are left unchanged when replacing $k\to k+1$. But this follows from $\sum_ix_i^2=1$, because:
\begin{eqnarray*}
&&x_{i_1}\ldots x_{i_k}x_{i_{k+1}}=x_{i_{\sigma(1)}}\ldots x_{i_{\sigma(k)}}x_{i_{k+1}}\\
&\implies&x_{i_1}\ldots x_{i_k}x_{i_{k+1}}^2=x_{i_{\sigma(1)}}\ldots x_{i_{\sigma(k)}}x_{i_{k+1}}^2\\
&\implies&\sum_{i_{k+1}}x_{i_1}\ldots x_{i_k}x_{i_{k+1}}^2=\sum_{i_{k+1}}x_{i_{\sigma(1)}}\ldots x_{i_{\sigma(k)}}x_{i_{k+1}}^2\\
&\implies&x_{i_1}\ldots x_{i_k}=x_{i_{\sigma(1)}}\ldots x_{i_{\sigma(k)}}
\end{eqnarray*}

Thus we can indeed ``simplify at right'', and this gives the result.
\end{proof}

As already mentioned, our goal in what follows will be that of proving that the 3 main spheres are the only monomial ones. In order to prove this result, we will use group theory methods. We call a subgroup $G\subset S_\infty$ filtered when it is stable under concatenation, in the sense that when writing $G=(G_k)$ with $G_k\subset S_k$, we have:
$$\sigma\in G_k,\pi\in G_l\implies \sigma\pi\in G_{k+l}$$

With this convention, we have the following result:

\index{filtered group}
\index{infinite permutation group}

\begin{theorem}
The monomial spheres are the subsets $S_G\subset S^{N-1}_{\mathbb R,+}$ given by
$$C(S_G)=C(S^{N-1}_{\mathbb R,+})\Big/\Big<x_{i_1}\ldots x_{i_k}=x_{i_{\sigma(1)}}\ldots x_{i_{\sigma(k)}},\forall (i_1,\ldots,i_k)\in\{1,\ldots,N\}^k,\forall\sigma\in G_k\Big>$$
where $G=(G_k)$ is a filtered subgroup of $S_\infty=(S_k)$.
\end{theorem}

\begin{proof}
We know from Proposition 16.7 that the construction in the statement produces a monomial sphere. Conversely, given a monomial sphere $S\subset S^{N-1}_{\mathbb R,+}$, let us set:
$$G_k=\left\{\sigma\in S_k\Big|x_{i_1}\ldots x_{i_k}=x_{i_{\sigma(1)}}\ldots x_{i_{\sigma(k)}},\forall (i_1,\ldots,i_k)\in\{1,\ldots,N\}^k\right\}$$

With $G=(G_k)$ we have then $S=S_G$. Thus, it remains to prove that $G$ is a filtered group. But since the relations $x_{i_1}\ldots x_{i_k}=x_{i_{\sigma(1)}}\ldots x_{i_{\sigma(k)}}$ can be composed and reversed, each $G_k$ follows to be stable under composition and inversion, and is therefore a group. Also, since the relations $x_{i_1}\ldots x_{i_k}=x_{i_{\sigma(1)}}\ldots x_{i_{\sigma(k)}}$ can be concatenated as well, our group $G=(G_k)$ is stable under concatenation, and we are done. 
\end{proof}

At the level of examples, according to our definitions, the simplest filtered groups, namely $\{1\}\subset S_\infty$, produce the simplest real spheres, namely:
$$S^{N-1}_{\mathbb R,+}\supset S^{N-1}_\mathbb R$$

In order to discuss now the half-classical case, we need to introduce and study a certain privileged intermediate filtered group $\{1\}\subset S_\infty^*\subset S_\infty$, which will eventually produce the intermediate sphere $S^{N-1}_{\mathbb R,+}\supset S^{N-1}_{\mathbb R,*}\supset S^{N-1}_\mathbb R$. This can be done as follows:

\begin{proposition}
Let $S_\infty^*\subset S_\infty$ be the set of permutations having the property that when labelling cyclically the legs as follows
$$\bullet\circ\bullet\circ\ldots$$
each string joins a black leg to a white leg.
\begin{enumerate}
\item $S_\infty^*$ is a filtered subgroup of $S_\infty$, generated by the half-classical crossing.

\item We have $S_{2k}^*\simeq S_k\times S_k$, and $S^*_{2k+1}\simeq S_k\times S_{k+1}$, for any $k\in\mathbb N$.
\end{enumerate}
\end{proposition}

\begin{proof}
The fact that $S_\infty^*$ is indeed a subgroup of $S_\infty$, which is filtered, is clear. Observe now that the half-classical crossing has the ``black-to-white'' joining property:
$$
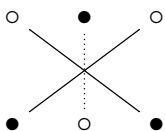
\xymatrix@R=10mm@C=5mm{\circ\ar@{-}[drr]&\bullet\ar@{.}[d]&\circ\ar@{-}[dll]\\\bullet&\circ&\bullet}$$ 

Thus this crossing belongs to $S_3^*$, and it is routine to check that the filtered subgroup of $S_\infty$ generated by it is the whole $S_\infty^*$. Regarding now the last assertion, observe first that the filtered subgroups $S_3^*,S_4^*$ consist of the following permutations:
$$
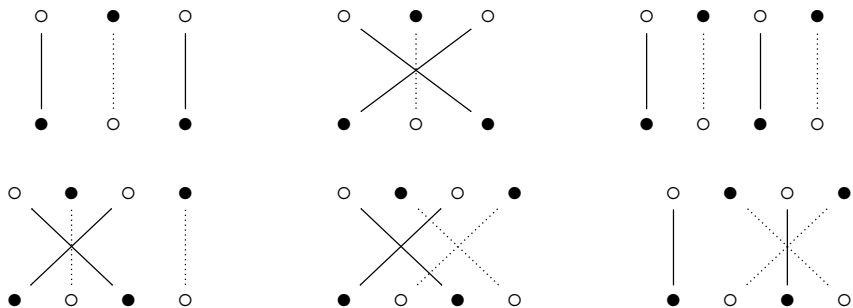
\xymatrix@R=10mm@C=5mm{\circ\ar@{-}[d]&\bullet\ar@{.}[d]&\circ\ar@{-}[d]\\\bullet&\circ&\bullet}\qquad\qquad
\xymatrix@R=10mm@C=5mm{\circ\ar@{-}[drr]&\bullet\ar@{.}[d]&\circ\ar@{-}[dll]\\\bullet&\circ&\bullet}\qquad\qquad 
\xymatrix@R=10mm@C=3mm{\circ\ar@{-}[d]&\bullet\ar@{.}[d]&\circ\ar@{-}[d]&\bullet\ar@{.}[d]\\\bullet&\circ&\bullet&\circ}$$
$$\xymatrix@R=10mm@C=3mm{\circ\ar@{-}[drr]&\bullet\ar@{.}[d]&\circ\ar@{-}[dll]&\bullet\ar@{.}[d]\\\bullet&\circ&\bullet&\circ}\qquad\qquad 
\xymatrix@R=10mm@C=3mm{\circ\ar@{-}[drr]&\bullet\ar@{.}[drr]&\circ\ar@{-}[dll]&\bullet\ar@{.}[dll]\\\bullet&\circ&\bullet&\circ}\qquad\qquad
\xymatrix@R=10mm@C=3mm{\circ\ar@{-}[d]&\bullet\ar@{.}[drr]&\circ\ar@{-}[d]&\bullet\ar@{.}[dll]\\\bullet&\circ&\bullet&\circ}$$

Thus we have $S_3^*=S_1\times S_2$ and $S_4^*=S_2\times S_2$, with the first component coming from dotted permutations, and with the second component coming from the solid line permutations. The same argument works in general, and gives the last assertion.
\end{proof}

Now back to the main 3 real spheres, the result is as follows:

\begin{proposition}
The basic monomial real spheres, namely 
$$S^{N-1}_\mathbb R\subset S^{N-1}_{\mathbb R,*}\subset S^{N-1}_{\mathbb R,+}$$
come respectively from the filtered groups $S_\infty\supset S_\infty^*\supset\{1\}$.
\end{proposition}

\begin{proof}
This is clear by definition in the classical and in the free cases. In the half-liberated case, the result follows from Proposition 16.9 (1).
\end{proof}

Now back to the general case, with the idea in mind of proving the uniqueness of the above spheres, consider a monomial sphere $S_G\subset S^{N-1}_{\mathbb R,+}$, with the filtered group $G\subset S_\infty$ taken to be maximal, as in the proof of Theorem 16.8. We have the following result:

\begin{proposition}
The filtered group $G\subset S_\infty$ associated to a monomial sphere $S\subset S^{N-1}_{\mathbb R,+}$ is stable under the following operations, on the corresponding diagrams:
\begin{enumerate}
\item Removing outer strings.

\item Removing neighboring strings.
\end{enumerate}
\end{proposition}

\begin{proof}
Both these results follow by using the quadratic condition:

\medskip

(1) Regarding the outer strings, by summing over $a$, we have:
\begin{eqnarray*}
Xa=Ya
&\implies&Xa^2=Ya^2\\
&\implies&X=Y
\end{eqnarray*}

We have as well the following computation:
\begin{eqnarray*}
aX=aY
&\implies&a^2X=a^2Y\\
&\implies&X=Y
\end{eqnarray*}

(2) Regarding the neighboring strings, once again by summing over $a$, we have:
\begin{eqnarray*}
XabY=ZabT
&\implies&Xa^2Y=Za^2T\\
&\implies&XY=ZT
\end{eqnarray*}

We have as well the following computation:
\begin{eqnarray*}
XabY=ZbaT
&\implies&Xa^2Y=Za^2T\\
&\implies&XY=ZT
\end{eqnarray*}

Thus $G=(G_k)$ has both the properties in the statement.
\end{proof}

We can now state and prove a main result, as follows:

\begin{theorem}
There is only one intermediate monomial sphere
$$S^{N-1}_\mathbb R\subset S\subset S^{N-1}_{\mathbb R,+}$$
namely the half-classical real sphere $S^{N-1}_{\mathbb R,*}$.
\end{theorem}

\begin{proof}
We will prove that the only filtered groups $G\subset S_\infty$ satisfying the conditions in Proposition 16.11 are those correspoding to our 3 spheres, namely:
$$\{1\}\subset S_\infty^*\subset S_\infty$$

In order to do so, consider such a filtered group $G\subset S_\infty$. We assume this group to be non-trivial, $G\neq\{1\}$, and we want to prove that we have $G=S_\infty^*$ or $G=S_\infty$.

\medskip

\underline{Step 1}. Our first claim is that $G$ contains a 3-cycle. Assume indeed that two permutations $\pi,\sigma\in S_\infty$ have support overlapping on exactly one point, say:
$$supp(\pi)\cap supp(\sigma)=\{i\}$$

The point is then that the commutator $\sigma^{-1}\pi^{-1}\sigma\pi$ is a 3-cycle, namely:
$$(i,\sigma^{-1}(i),\pi^{-1}(i))$$

Indeed the computation of the commutator goes as follows:
$$
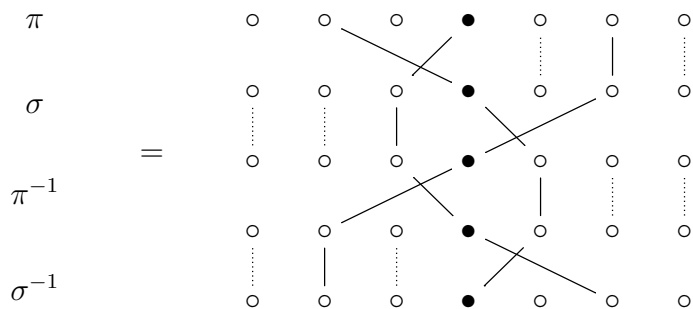
\xymatrix@R=7mm@C=5mm{\pi\\ \sigma\\ \pi^{-1}\\ \sigma^{-1}}\qquad
\xymatrix@R=6mm@C=5mm{\\ \\ =}\qquad
\xymatrix@R=5mm@C=5mm{
\circ&\circ\ar@{-}[drr]&\circ&\bullet\ar@{-}[dl]&\circ\ar@{.}[d]&\circ\ar@{-}[d]&\circ\ar@{.}[d]\\
\circ\ar@{.}[d]&\circ\ar@{.}[d]&\circ\ar@{-}[d]&\bullet\ar@{-}[dr]&\circ&\circ\ar@{-}[dll]&\circ\\
\circ&\circ&\circ\ar@{-}[dr]&\bullet\ar@{-}[dll]&\circ\ar@{-}[d]&\circ\ar@{.}[d]&\circ\ar@{.}[d]\\
\circ\ar@{.}[d]&\circ\ar@{-}[d]&\circ\ar@{.}[d]&\bullet\ar@{-}[drr]&\circ\ar@{-}[dl]&\circ&\circ\\
\circ&\circ&\circ&\bullet&\circ&\circ&\circ
}$$

Now let us pick a non-trivial element $\tau\in G$. By removing outer strings at right and at left we obtain permutations $\tau'\in G_k,\tau''\in G_s$ having a non-trivial action on their right/left leg, and the trick applies, with:
$$\pi=\tau'\otimes id_{s-1}\quad,\quad 
\sigma=id_{k-1}\otimes\tau''$$

Thus, $G$ contains a 3-cycle, as claimed.

\medskip

\underline{Step 2}. Our second claim is $G$ must contain one of the following permutations:
$$
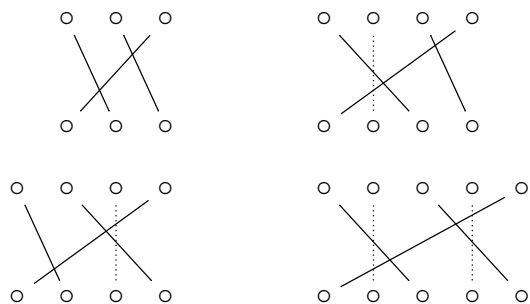
\xymatrix@R=10mm@C=2mm{
\circ\ar@{-}[dr]&\circ\ar@{-}[dr]&\circ\ar@{-}[dll]\\
\circ&\circ&\circ}\qquad\qquad
\xymatrix@R=10mm@C=2mm{
\circ\ar@{-}[drr]&\circ\ar@{.}[d]&\circ\ar@{-}[dr]&\circ\ar@{-}[dlll]\\
\circ&\circ&\circ&\circ}$$
$$\xymatrix@R=10mm@C=2mm{
\circ\ar@{-}[dr]&\circ\ar@{-}[drr]&\circ\ar@{.}[d]&\circ\ar@{-}[dlll]\\
\circ&\circ&\circ&\circ}\qquad\qquad
\xymatrix@R=10mm@C=2mm{
\circ\ar@{-}[drr]&\circ\ar@{.}[d]&\circ\ar@{-}[drr]&\circ\ar@{.}[d]&\circ\ar@{-}[dllll]\\
\circ&\circ&\circ&\circ&\circ}$$

Indeed, consider the 3-cycle that we just constructed. By removing all outer strings, and then all pairs of adjacent vertical strings, we are left with these permutations.

\medskip

\underline{Step 3}. Our claim now is that we must have $S_\infty^*\subset G$. Indeed, let us pick one of the permutations that we just constructed, and apply to it our various diagrammatic rules. From the first permutation we can obtain the basic crossing, as follows:
$$
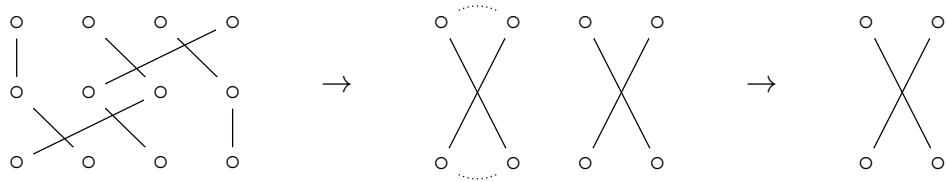
\xymatrix@R=5mm@C=5mm{
\circ\ar@{-}[d]&\circ\ar@{-}[dr]&\circ\ar@{-}[dr]&\circ\ar@{-}[dll]\\
\circ\ar@{-}[dr]&\circ\ar@{-}[dr]&\circ\ar@{-}[dll]&\circ\ar@{-}[d]\\
\circ&\circ&\circ&\circ}
\qquad
\xymatrix@R=5mm@C=5mm{
\\ \to\\}\qquad
\xymatrix@R=6mm@C=5mm{
\circ\ar@{-}[ddr]\ar@/^/@{.}[r]&\circ\ar@{-}[ddl]&\circ\ar@{-}[ddr]&\circ\ar@{-}[ddl]\\
\\
\circ\ar@/_/@{.}[r]&\circ&\circ&\circ}
\qquad
\xymatrix@R=5mm@C=5mm{
\\ \to\\}\qquad
\xymatrix@R=6mm@C=5mm{
\circ\ar@{-}[ddr]&\circ\ar@{-}[ddl]\\
\\
\circ&\circ}$$

Also, by removing a suitable $\slash\hskip-2.1mm\backslash$ shaped configuration, which is represented by dotted lines in the diagrams below, we can obtain the basic crossing from the second and third permutation, and the half-liberated crossing from the fourth permutation:
$$
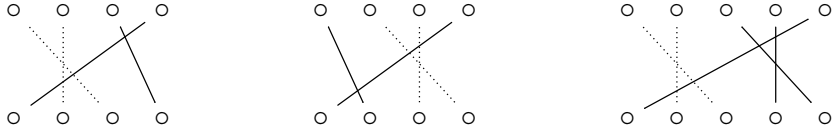
\xymatrix@R=10mm@C=2mm{
\circ\ar@{.}[drr]&\circ\ar@{.}[d]&\circ\ar@{-}[dr]&\circ\ar@{-}[dlll]\\
\circ&\circ&\circ&\circ}\qquad\qquad
\xymatrix@R=10mm@C=2mm{
\circ\ar@{-}[dr]&\circ\ar@{.}[drr]&\circ\ar@{.}[d]&\circ\ar@{-}[dlll]\\
\circ&\circ&\circ&\circ}\qquad\qquad
\xymatrix@R=10mm@C=2mm{
\circ\ar@{.}[drr]&\circ\ar@{.}[d]&\circ\ar@{-}[drr]&\circ\ar@{-}[d]&\circ\ar@{-}[dllll]\\
\circ&\circ&\circ&\circ&\circ}$$

Thus, in all cases we have a basic or half-liberated crossing, and so, as desired: 
$$S_\infty^*\subset G$$

\underline{Step 4}. Our last claim, which will finish the proof, is that there is no proper intermediate subgroup as follows:
$$S_\infty^*\subset G\subset S_\infty$$

In order to prove this, observe that $S_\infty^*\subset S_\infty$ is the subgroup of parity-preserving permutations, in the sense that ``$i$ even $\implies$ $\sigma(i)$ even''. 

\medskip

Now let us pick an element $\sigma\in S_k-S_k^*$, with $k\in\mathbb N$. We must prove that the group $G=<S_\infty^*,\sigma>$ equals the whole $S_\infty$. In order to do so, we use the fact that $\sigma$ is not parity preserving. Thus, we can find $i$ even such that $\sigma(i)$ is odd. In addition, up to passing to $\sigma|$, we can assume that $\sigma(k)=k$, and then, up to passing one more time to $\sigma|$, we can further assume that $k$ is even. Since both the indices $i,k$ are even we have:
$$(i,k)\in S_k^*$$

We conclude that the following element belongs to $G$:
$$\sigma(i,k)\sigma^{-1}=(\sigma(i),k)$$

But, since $\sigma(i)$ is odd, by deleting an appropriate number of vertical strings, $(\sigma(i),k)$ reduces to the basic crossing $(1,2)$. Thus $G=S_\infty$, and we are done.
\end{proof}

As already mentioned in the above, the story is not over with this kind of result, because the complex case still remains to be worked out.

\bigskip

Moving on, some even better results can be obtained by looking at the various intersections between spheres, twisted and untwisted. These intersections cannot be thought of in general as being smooth, let us introduce the following objects:

\begin{definition}
The polygonal spheres are real algebraic manifolds, defined as
$$S^{N-1,d-1}_\mathbb R=\left\{x\in S^{N-1}_\mathbb R\Big|x_{i_0}\ldots x_{i_d}=0,\forall i_0,\ldots,i_d\ {\rm distinct}\right\}$$
depending on integers $1\leq d\leq N$.
\end{definition}

More generally, we have the following construction of ``generalized polygonal spheres'', which applies to the half-classical and twisted cases too:
$$C(\dot{S}^{N-1,d-1}_{\mathbb R,\times})=C\big(\dot{S}^{N-1}_{\mathbb R,\times}\big)\Big/\Big<x_{i_0}\ldots x_{i_d}=0,\forall i_0,\ldots,i_d\ {\rm distinct}\Big>$$

With these conventions, we have the following result, dealing with all the spheres that we have so far in real case, namely twisted, untwisted and intersections:

\begin{theorem}
The diagram obtained by intersecting the real spheres is
$$\xymatrix@R=13mm@C=16mm{
S^{N-1}_\mathbb R\ar[r]&S^{N-1}_{\mathbb R,*}\ar[r]&S^{N-1}_{\mathbb R,+}\\
S^{N-1,1}_\mathbb R\ar[r]\ar[u]&S^{N-1,1}_{\mathbb R,*}\ar[r]\ar[u]&\bar{S}^{N-1}_{\mathbb R,*}\ar[u]\\
S^{N-1,0}_\mathbb R\ar[r]\ar[u]&\bar{S}^{N-1,1}_\mathbb R\ar[r]\ar[u]&\bar{S}^{N-1}_\mathbb R\ar[u]}$$
and so all these spheres are generalized polygonal spheres.
\end{theorem}

\begin{proof}
Consider the 4-diagram obtained by intersecting the 5 main spheres:
$$\xymatrix@R=13mm@C=13mm{
S^{N-1}_\mathbb R\cap\bar{S}^{N-1}_{\mathbb R,*}\ar[r]&S^{N-1}_{\mathbb R,*}\cap\bar{S}^{N-1}_{\mathbb R,*}\\
S^{N-1}_\mathbb R\cap\bar{S}^{N-1}_\mathbb R\ar[r]\ar[u]&S^{N-1}_{\mathbb R,*}\cap\bar{S}^{N-1}_\mathbb R\ar[u]}$$

We want to prove that this diagram coincides with the 4-diagram at bottom left in the statement. But this is clear, because combining commutation and anticommutation leads to the vanishing relations defining the spheres of type $\dot{S}^{N-1,d-1}_{\mathbb R,\times}$. More precisely:

\medskip

(1) $S^{N-1}_\mathbb R\cap\bar{S}^{N-1}_\mathbb R$ consists of the points $x\in S^{N-1}_\mathbb R$ such that, for any $i\neq j$:
$$x_ix_j=-x_jx_i$$

Now since we have as well $x_ix_j=x_jx_i$, for any $i,j$, this relation reads $x_ix_j=0$ for $i\neq j$, which means that we have $x\in S^{N-1,0}_\mathbb R$, as desired.

\medskip

(2) $S^{N-1}_\mathbb R\cap\bar{S}^{N-1}_{\mathbb R,*}$ consists of the points $x\in S^{N-1}_\mathbb R$ such that, for $i,j,k$ distinct:
$$x_ix_jx_k=-x_kx_jx_i$$

Once again by commutativity, this relation is equivalent to $x\in S^{N-1,1}_\mathbb R$, as desired.

\medskip

(3) $S^{N-1}_{\mathbb R,*}\cap\bar{S}^{N-1}_\mathbb R$ is obtained from $\bar{S}^{N-1}_\mathbb R$ by imposing to the standard coordinates the half-commutation relations $abc=cba$. On the other hand, we know from $\bar{S}^{N-1}_\mathbb R\subset \bar{S}^{N-1}_{\mathbb R,*}$ that the standard coordinates on $\bar{S}^{N-1}_\mathbb R$ satisfy $abc=-cba$ for $a,b,c$ distinct, and $abc=cba$ otherwise. Thus, the relations brought by intersecting with $S^{N-1}_{\mathbb R,*}$ reduce to the relations $abc=0$ for $a,b,c$ distinct, and so we are led to the sphere $\bar{S}^{N-1,1}_\mathbb R$.

\medskip

(4) $S^{N-1}_{\mathbb R,*}\cap\bar{S}^{N-1}_{\mathbb R,*}$ is obtained from $\bar{S}^{N-1}_{\mathbb R,*}$ by imposing the relations $abc=-cba$ for $a,b,c$ distinct, and $abc=cba$ otherwise. Since we know that $abc=cba$ for any $a,b,c$, the extra relations reduce to $abc=0$ for $a,b,c$ distinct, and so we are led to $S^{N-1,1}_{\mathbb R,*}$.
\end{proof} 

We have now all the needed ingredients for axiomatizing the various spheres appearing so far, namely the twisted and untwisted ones, and their intersections:

\begin{definition}
We have $3$ types of quantum spheres $S\subset S^{N-1}_{\mathbb R,+}$, as follows:
\begin{enumerate}
\item Monomial, namely $\dot{S}^{N-1}_{\mathbb R,E}$, with $E\subset S_\infty$, obtained via the following relations:
$$\left\{\dot{\mathcal R}_\sigma\Big|\sigma\in E\right\}$$

\item Mixed monomial, which appear as intersections as follows, with $E,F\subset S_\infty$:
$$S^{N-1}_{\mathbb R,E,F}=S^{N-1}_{\mathbb R,E}\cap\bar{S}^{N-1}_{\mathbb R,F}$$

\item Polygonal, which are again intersections, with $E,F\subset S_\infty$, and $d\in\{1,\ldots,N\}$:
$$S^{N-1,d-1}_{\mathbb R,E,F}=S^{N-1}_{\mathbb R,E,F}\cap S^{N-1,d-1}_{\mathbb R,+}$$
\end{enumerate}
\end{definition}

With the above notions, we cover all spheres appearing so far. Observe also that the set of mixed monomial spheres is closed under intersections. The same holds for the set of polygonal spheres, because we have the following formula:
$$S^{N-1,d-1}_{\mathbb R,E,F}\cap S^{N-1,d'-1}_{\mathbb R,E',F'}=S^{N-1,min(d,d')-1}_{\mathbb R,E\cup E',F\cup F'}$$

Let us try now to understand the structure of the various types of spheres, by using the real sphere technology developed before. We call a group of permutations $G\subset S_\infty$ filtered if, with $G_k=G\cap S_k$, we have $G_k\times G_l\subset G_{k+l}$, for any $k,l$. We have:

\begin{proposition}
The various spheres can be parametrized by groups, as follows:
\begin{enumerate}
\item Monomial case: $\dot{S}^{N-1}_{\mathbb R,G}$, with $G\subset S_\infty$ filtered group.

\item Mixed monomial case: $S^{N-1}_{\mathbb R,G,H}$, with $G,H\subset S_\infty$ filtered groups.

\item Polygonal case: $S^{N-1,d-1}_{\mathbb R,G,H}$, with $G,H\subset S_\infty$ filtered groups, and $d\in\{1,\ldots,N\}$.
\end{enumerate}
\end{proposition}

\begin{proof}
This basically follows from the theory developed before, as follows:

\medskip

(1) As explained before, in order to prove this assertion, for a monomial sphere $S=\dot{S}_{\mathbb R,E}$, we can take $G\subset S_\infty$ to be the set of permutations $\sigma\in S_\infty$ having the property that the relations $\dot{\mathcal R}_\sigma$ hold for the standard coordinates of $S$. We have then $E\subset G$, we have as well $S=\dot{S}^{N-1}_{\mathbb R,G}$, and the fact that $G$ is a filtered group is clear as well.

\medskip

(2) This follows from (1), by taking intersections.

\medskip

(3) Once again this follows from (1), by taking intersections.
\end{proof}

The idea in what follows will be that of writing the 9 main polygonal spheres as in Proposition 16.16 (2), as to reach to a ``standard parametrization'' for our spheres. 

\bigskip

We recall that the permutations $\sigma\in S_\infty$ having the property that when labelling clockwise their legs $\circ\bullet\circ\bullet\ldots$, and string joins a white leg to a black leg, form a filtered group, denoted $S_\infty^*\subset S_\infty$. This group comes from the general half-liberation considerations from above, and its algebraic structure is very simple, as follows:
$$S_{2n}^*\simeq S_n\times S_n\quad,\quad 
S_{2n+1}^*\simeq S_n\times S_{n+1}$$

Let us formulate as well the following definition:

\index{standard parametrization}

\begin{definition}
We call a mixed monomial sphere parametrization 
$$S=S^{N-1}_{\mathbb R,G,H}$$
standard when both filtered groups $G,H\subset S_\infty$ are chosen to be maximal.
\end{definition}

In this case, Proposition 16.16 and its proof tell us that $G,H$ encode all the monomial relations which hold in $S$. With these conventions, we have the following result, extending some previous findings from above, regarding the untwisted spheres:

\begin{theorem}
The standard parametrization of the $9$ main spheres is
$$
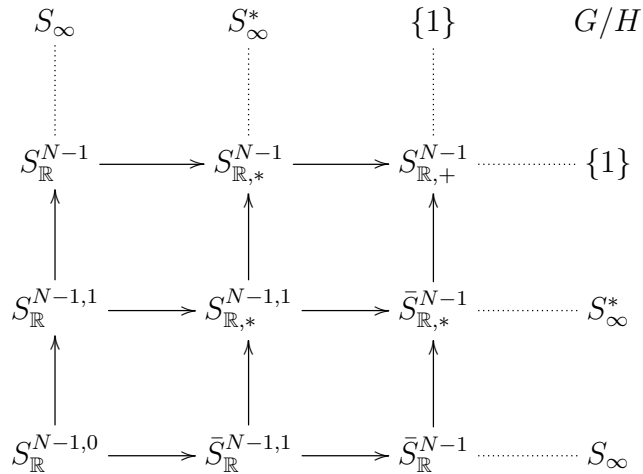
\xymatrix@R=11.5mm@C=11.5mm{
S_\infty\ar@{.}[d]&S_\infty^*\ar@{.}[d]&\{1\}\ar@{.}[d]&G/H\\
S^{N-1}_\mathbb R\ar[r]&S^{N-1}_{\mathbb R,*}\ar[r]&S^{N-1}_{\mathbb R,+}&\{1\}\ar@{.}[l]\\
S^{N-1,1}_\mathbb R\ar[r]\ar[u]&S^{N-1,1}_{\mathbb R,*}\ar[r]\ar[u]&\bar{S}^{N-1}_{\mathbb R,*}\ar[u]&S_\infty^*\ar@{.}[l]\\
S^{N-1,0}_\mathbb R\ar[r]\ar[u]&\bar{S}^{N-1,1}_\mathbb R\ar[r]\ar[u]&\bar{S}^{N-1}_\mathbb R\ar[u]&S_\infty\ar@{.}[l]}$$
so these spheres come from the $3\times 3=9$ pairs of groups among $\{1\}\subset S_\infty^*\subset S_\infty$.
\end{theorem}

\begin{proof}
The fact that we have parametrizations as above is known to hold for the 5 untwisted and twisted spheres. For the remaining 4 spheres the result follows by intersecting, by using the following formula, valid for any $E,F\subset S_\infty$:
$$S^{N-1}_{\mathbb R,E,F}\cap S^{N-1}_{\mathbb R,E',F'}=S^{N-1}_{\mathbb R,E\cup E',F\cup F'}$$

In order to prove now that the parametrizations are standard, we must compute the following two filtered groups, and show that we get the groups in the statement:
$$G=\left\{\sigma\in S_\infty\Big|{\rm the\ relations\ }\mathcal R_\sigma\ {\rm hold\ over\ }S\right\}$$ 
$$H=\left\{\sigma\in S_\infty\Big|{\rm the\ relations\ }\bar{\mathcal R}_\sigma\ {\rm hold\ over\ }S\right\}$$ 

As a first observation, by using the various inclusions between spheres, we just have to compute $G$ for the spheres on the bottom, and $H$ for the spheres on the left:
$$X=S^{N-1,0}_\mathbb R,\bar{S}^{N-1,1}_\mathbb R,\bar{S}^{N-1}_\mathbb R\implies G=S_\infty,S_\infty^*,\{1\}$$
$$X=S^{N-1,0}_\mathbb R,S^{N-1,1}_\mathbb R,S^{N-1}_\mathbb R\implies H=S_\infty,S_\infty^*,\{1\}$$

The results for $S^{N-1,0}_\mathbb R$ being clear, we are left with computing the remaining 4 groups, for the spheres $S^{N-1}_\mathbb R,\bar{S}^{N-1}_\mathbb R,S^{N-1,1}_\mathbb R,\bar{S}^{N-1,1}_\mathbb R$. The proof here goes as follows:

\medskip

(1) $S^{N-1}_\mathbb R$. According to the definition of $H=(H_k)$, we have:
\begin{eqnarray*}
H_k
&=&\left\{\sigma\in S_k\Big|x_{i_1}\ldots x_{i_k}=\varepsilon\left(
\ker\begin{pmatrix}i_1&\ldots&i_k\\
i_{\sigma(1)}&\ldots&i_{\sigma(k)}
\end{pmatrix}\right)
x_{i_{\sigma(1)}}\ldots x_{i_{\sigma(k)}},\forall i_1,\ldots,i_k\right\}\\
&=&\left\{\sigma\in S_k\Big|
\varepsilon\left(
\ker\begin{pmatrix}i_1&\ldots&i_k\\
i_{\sigma(1)}&\ldots&i_{\sigma(k)}
\end{pmatrix}\right)
=1,\forall i_1,\ldots,i_k\right\}\\
&=&\left\{\sigma\in S_k\Big|\varepsilon(\tau)=1,\forall\tau\leq\sigma\right\}\end{eqnarray*}

Now observe that for any permutation $\sigma\in S_k,\sigma\neq1_k$, we can always find a partition $\tau\leq\sigma$ satisfying the following condition:
$$\varepsilon(\tau)=-1$$

We deduce that we have $H_k=\{1_k\}$, and so $H=\{1\}$, as desired.

\medskip

(2) $\bar{S}^{N-1}_\mathbb R$. The proof of $G=\{1\}$ here is similar to the proof of $H=\{1\}$ in (1) above, by using the same combinatorial ingredient at the end.

\medskip

(3) $S^{N-1,1}_\mathbb R$. By definition of $H=(H_k)$, a permutation $\sigma\in S_k$ belongs to $H_k$ when the following condition is satisfied, for any choice of the indices $i_1,\ldots,i_k$:
$$x_{i_1}\ldots x_{i_k}=\varepsilon\left(\ker\begin{pmatrix}i_1&\ldots&i_k\\ i_{\sigma(1)}&\ldots&i_{\sigma(k)}\end{pmatrix}\right)x_{i_{\sigma(1)}}\ldots x_{i_{\sigma(k)}}$$

We have three cases here, as follows:

\medskip

-- When $|\ker i|=1$ this formula reads $x_r^k=x_r^k$, which is true. 

\medskip

-- When $|\ker i|\geq3$ this formula is automatically satisfied as well, because by using the relations $ab=ba$, and $abc=0$ for $a,b,c$ distinct, which both hold over $S^{N-1,1}_\mathbb R$, this formula reduces to $0=0$. 

\medskip

-- Thus, we are left with studying the case $|\ker i|=2$. Here the quantities on the left $x_{i_1}\ldots x_{i_k}$ will not vanish, so the sign on the right must be 1, and we therefore have:
$$H_k=\left\{\sigma\in S_k\Big|\varepsilon(\tau)=1,\forall\tau\leq\sigma,|\tau|=2\right\}$$

Now by coloring the legs of $\sigma$ clockwise $\circ\bullet\circ\bullet\ldots$, the above condition is satisfied when each string of $\sigma$ joins a white leg to a black leg. Thus $H_k=S_k^*$, as desired.

\medskip

(4) $\bar{S}^{N-1,1}_\mathbb R$. The proof of $G=S_\infty^*$ here is similar to the proof of $H=S_\infty^*$ in (3) above, by using the same combinatorial ingredient at the end.
\end{proof}

\section*{16c. Projective spaces}

Getting now to projective geometry, which is where our threefold principle really shines, let us start with something that you surely know, namely:

\begin{definition}
A projective space is a space consisting of points and lines, subject to the following conditions:
\begin{enumerate}
\item Each $2$ points determine a line.

\item Each $2$ lines cross, on a point.
\end{enumerate}
\end{definition}

As a basic example we have the usual projective plane $P^2_\mathbb R$, which is best seen as being the space of lines in $\mathbb R^3$ passing through the origin. To be more precise, let us call each of these lines in $\mathbb R^3$ passing through the origin a ``point'' of $P^2_\mathbb R$, and let us also call each plane in $\mathbb R^3$ passing through the origin a ``line'' of $P^2_\mathbb R$. Now observe the following:

\bigskip

(1) Each $2$ points determine a line. Indeed, 2 points in our sense means 2 lines in $\mathbb R^3$ passing through the origin, and these 2 lines obviously determine a plane in $\mathbb R^3$ passing through the origin, namely the plane they belong to, which is a line in our sense.

\bigskip

(2) Each $2$ lines cross, on a point. Indeed, 2 lines in our sense means 2 planes in $\mathbb R^3$ passing through the origin, and these 2 planes obviously determine a line in $\mathbb R^3$ passing through the origin, namely their intersection, which is a point in our sense.

\bigskip

Thus, what we have is a projective space in the sense of Definition 16.19. More generally, we have the following construction, in arbitrary dimensions:

\begin{theorem}
We can define the projective space $P^{N-1}_\mathbb R$ as being the space of lines in $\mathbb R^N$ passing through the origin, and in small dimensions:
\begin{enumerate}
\item $P^1_\mathbb R$ is the usual circle.

\item $P^2_\mathbb R$ is some sort of twisted sphere.
\end{enumerate}
\end{theorem}

\begin{proof}
We have several assertions here, the idea being as follows:

\medskip

(1) To start with, the fact that the space $P^{N-1}_\mathbb R$ constructed in the statement is indeed a projective space in the sense of Definition 16.19 follows from definitions.

\medskip

(2) At $N=2$ now, a line in $\mathbb R^2$ passing through the origin corresponds to 2 opposite points on the unit circle $\mathbb T\subset\mathbb R^2$, according to the following scheme:
$$
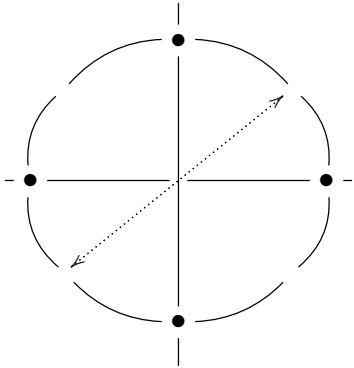
\xymatrix@R=9pt@C=3pt{
&&&&&&\ar@{-}[d]\\
&&&&&&\bullet\ar@{-}[dddddd]\ar@{-}@/_/[dllll]&\\
&&\ar@{-}@/_/[ddl]&&&&&&&&\ar@.[ddddllllllll]\ar@{-}@/_/[ullll]\\
&&&&&&&&&&\\
\ar@{-}[r]&\bullet\ar@{-}@/_/[ddr]\ar@{-}[rrrrr]&&&&&\ar@{-}[rrrrr]&&&&&\bullet\ar@{-}@/_/[uul]\ar@{-}[r]&\\
&&&&&&&&&&\\
&&\ar@{-}@/_/[drrrr]\ar@.[uuuurrrrrrrr]&&&&&&&&\ar@{-}@/_/[uur]\\
&&&&&&\bullet\ar@{-}@/_/[urrrr]\ar@{-}[d]\\
&&&&&&}$$

Thus, $P^1_\mathbb R$ corresponds to the upper semicircle of $\mathbb T$, with the endpoints identified, and so we obtain a circle, $P^1_\mathbb R=\mathbb T$, according to the following scheme: 
$$
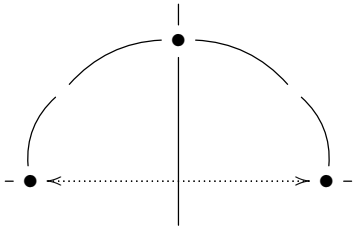
\xymatrix@R=9pt@C=3pt{
&&&&&&\ar@{-}[d]\\
&&&&&&\bullet\ar@{-}[dddd]\ar@{-}@/_/[dllll]&\\
&&\ar@{-}@/_/[ddl]&&&&&&&&\ar@{-}@/_/[ullll]\\
&&&&&&&&&&\\
\ar@{-}[r]&\bullet\ar@.[rrrrrrrrrr]&&&&&&&&&&\bullet\ar@{-}@/_/[uul]\ar@{-}[r]\ar@.[llllllllll]&\\
&&&&&&&&
}$$

(3) At $N=3$, the space $P^2_\mathbb R$ corresponds to the upper hemisphere of the sphere $S^2_\mathbb R\subset\mathbb R^3$, with the points on the equator identified via $x=-x$. Topologically speaking, we can deform if we want the hemisphere into a square, with the equator becoming the boundary of this square, and in this picture, the $x=-x$ identification corresponds to a ``identify opposite edges, with opposite orientations'' folding method for the square:
$$
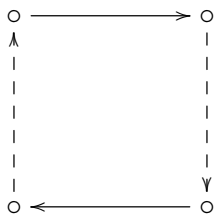
\xymatrix@R=60pt@C=60pt{
\circ\ar[r]&\circ\ar@{-->}[d]\\
\circ\ar@{-->}[u]&\circ\ar[l]}$$

(4) Thus, we have our space. In order to understand now what this beast is, let us look first at the other 3 possible methods of folding the square, which are as follows:
$$
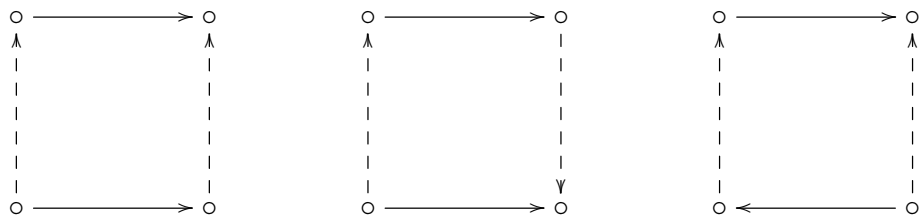
\xymatrix@R=60pt@C=60pt{
\circ\ar[r]&\circ\\
\circ\ar@{-->}[u]\ar[r]&\circ\ar@{-->}[u]}\qquad\qquad 
\xymatrix@R=60pt@C=60pt{
\circ\ar[r]&\circ\ar@{-->}[d]\\
\circ\ar[r]\ar@{-->}[u]&\circ}\qquad\qquad 
\xymatrix@R=60pt@C=60pt{
\circ\ar[r]&\circ\\
\circ\ar@{-->}[u]&\circ\ar[l]\ar@{-->}[u]}$$

Regarding the first space, the one on the left, things here are quite simple. Indeed, when identifying the solid edges we get a cylinder, and then when further identifying the dotted edges, what we get is some sort of closed cylinder, which is a torus.

\medskip

(5) Regarding the second space, the one in the middle, things here are more tricky. Indeed, when identifying the solid edges we get again a cylinder, but then when further identifying the dotted edges, we obtain some sort of ``impossible'' closed cylinder, called Klein bottle. This Klein bottle obviously cannot be drawn in 3 dimensions, but with a bit of imagination, you can see it, in its full splendor, in 4 dimensions.

\medskip

(6) Finally, regarding the third space, the one on the right, we know by symmetry that this must be the Klein bottle too. But we can see this as well via our standard folding method, namely identifying solid edges first, and dotted edges afterwards. Indeed, we first obtain in this way a M\"obius strip, and then, well, the Klein bottle.

\medskip

(7) With these preliminaries made, and getting back now to the projective space $P^2_\mathbb R$, we can see that this is something more complicated, of the same type, reminding the torus and the Klein bottle. So, we will call it ``sort of twisted sphere'', as in the statement, and exercise for you to imagine how this beast looks like, in 4 dimensions.
\end{proof}

In order to reach to free projective spaces, our starting point will be the following functional analytic description of the real and complex projective spaces $P^{N-1}_\mathbb R,P^{N-1}_\mathbb C$:

\index{real projective space}
\index{complex projective space}

\begin{theorem}
We have presentation results as follows,
\begin{eqnarray*}
C(P^{N-1}_\mathbb R)&=&C^*_{comm}\left((p_{ij})_{i,j=1,\ldots,N}\Big|p=\bar{p}=p^t=p^2,Tr(p)=1\right)\\
C(P^{N-1}_\mathbb C)&=&C^*_{comm}\left((p_{ij})_{i,j=1,\ldots,N}\Big|p=p^*=p^2,Tr(p)=1\right)
\end{eqnarray*}
for the algebras of continuous functions on the real and complex projective spaces.
\end{theorem}

\begin{proof}
We use the fact that the projective spaces $P^{N-1}_\mathbb R,P^{N-1}_\mathbb C$ can be respectively identified with the spaces of rank one projections in $M_N(\mathbb R),M_N(\mathbb C)$. With this picture in mind, it is clear that we have arrows $\leftarrow$. In order to construct now arrows $\to$, consider the universal algebras on the right, $A_R,A_C$. These algebras being both commutative, by the Gelfand theorem we can write, with $X_R,X_C$ being certain compact spaces:
$$A_R=C(X_R)\quad,\quad 
A_C=C(X_C)$$

Now by using the coordinate functions $p_{ij}$, we conclude that $X_R,X_C$ are certain spaces of rank one projections in $M_N(\mathbb R),M_N(\mathbb C)$. In other words, we have embeddings:
$$X_R\subset P^{N-1}_\mathbb R\quad,\quad 
X_C\subset P^{N-1}_\mathbb C$$

By transposing we obtain arrows $\to$, as desired.
\end{proof}

The point now is that the above result suggests the following definition:

\begin{definition}
Associated to any $N\in\mathbb N$ is the following universal algebra,
$$C(P^{N-1}_+)=C^*\left((p_{ij})_{i,j=1,\ldots,N}\Big|p=p^*=p^2,Tr(p)=1\right)$$
whose abstract spectrum is called ``free projective space''.
\end{definition}

Observe that, according to our presentation results for the real and complex projective spaces $P^{N-1}_\mathbb R$ and $P^{N-1}_\mathbb C$, we have embeddings of compact quantum spaces, as follows:
$$P^{N-1}_\mathbb R\subset P^{N-1}_\mathbb C\subset P^{N-1}_+$$

Let us first discuss the relation with the spheres. Given a closed subset $X\subset S^{N-1}_{\mathbb R,+}$, its projective version is by definition the quotient space $X\to PX$ determined by the fact that $C(PX)\subset C(X)$ is the subalgebra generated by the following variables:
$$p_{ij}=x_ix_j$$

In order to discuss the relation with the spheres, let us start with:

\begin{theorem}
The projective versions of the $3$ real spheres are as follows,
$$\xymatrix@R=15mm@C=15mm{
S^{N-1}_\mathbb R\ar[r]\ar[d]&S^{N-1}_{\mathbb R,*}\ar[r]\ar[d]&S^{N-1}_{\mathbb R,+}\ar[d]\\
P^{N-1}_\mathbb R\ar[r]&P^{N-1}_\mathbb C\ar[r]&P^{N-1}_+}$$
modulo the standard equivalence relation for the quantum algebraic manifolds.
\end{theorem}

\begin{proof}
The assertion at left is true by definition. For the assertion at right, we have to prove that the variables $p_{ij}=z_iz_j$ over the free sphere $S^{N-1}_{\mathbb R,+}$ satisfy the defining relations for $C(P^{N-1}_+)$, from Definition 16.22, namely:
$$p=p^*=p^2\quad,\quad 
Tr(p)=1$$

We first have the following computation:
$$(p^*)_{ij}
=p_{ji}^*
=(z_jz_i)^*
=z_iz_j
=p_{ij}$$

We have as well the following computation:
$$(p^2)_{ij}
=\sum_kp_{ik}p_{kj}
=\sum_kz_iz_k^2z_j
=z_iz_j\\
=p_{ij}$$

Finally, we have as well the following computation:
$$Tr(p)
=\sum_kp_{kk}
=\sum_kz_k^2
=1$$

Regarding now $PS^{N-1}_{\mathbb R,*}=P^{N-1}_\mathbb C$, the inclusion ``$\subset$'' follows from $abcd=cbad=cbda$. In the other sense now, the point is that we have a matrix model, as follows:
$$\pi:C(S^{N-1}_{\mathbb R,*})\to M_2(C(S^{N-1}_\mathbb C))\quad,\quad 
x_i\to\begin{pmatrix}0&z_i\\ \bar{z}_i&0\end{pmatrix}$$ 

But this gives the missing inclusion ``$\supset$'', and we are done. See \cite{bgo}.
\end{proof}

In addition to the above result, let us mention that passing to the complex case brings nothing new. This is because the projective version of the free complex sphere is equal to the free projective space constructed above. The same goes for the ``hybrid'' spheres.

\bigskip

In fact, things become considerably simpler in the projective geometry setting. Consider indeed the diagram of 9 main affine geometries, that we found above:
$$\xymatrix@R=40pt@C=42pt{
\mathbb R^N_+\ar[r]&\mathbb T\mathbb R^N_+\ar[r]&\mathbb C^N_+\\
\mathbb R^N_*\ar[u]\ar[r]&\mathbb T\mathbb R^N_*\ar[u]\ar[r]&\mathbb C^N_*\ar[u]\\
\mathbb R^N\ar[u]\ar[r]&\mathbb T\mathbb R^N\ar[u]\ar[r]&\mathbb C^N\ar[u]
}$$

The point now is that when looking at the projective versions of these geometries,  the diagram drastically simplifies. To be more precise, the diagram of projective versions of the corresponding spheres are as follows, consisting of 3 objects only:
$$\xymatrix@R=36pt@C=33pt{
P^{N-1}_+\ar[r]&P^{N-1}_+\ar[r]&P^{N-1}_+\\
P^{N-1}_\mathbb C\ar[r]\ar[u]&P^{N-1}_\mathbb C\ar[r]\ar[u]&P^{N-1}_\mathbb C\ar[u]\\
P^{N-1}_\mathbb R\ar[r]\ar[u]&P^{N-1}_\mathbb R\ar[r]\ar[u]&P^{N-1}_\mathbb C\ar[u]}$$

Thus, we are led to the conclusion that, under suitable combinatorial axioms, there should be only 3 projective geometries, namely the real, complex and free ones:
$$P^{N-1}_\mathbb R\subset P^{N-1}_\mathbb C\subset P^{N-1}_+$$ 

We will further explore this phenomenon, with a number of more precise results on the subject, featuring some non-trivial combinatorics, in the next section.

\section*{16d. Threefold way}

We would like to discuss now the axiomatization of monomial projective spaces, projective easiness, and projective geometries. In order to get started, let us formulate:

\begin{definition}
A monomial projective space is a closed subset $P\subset P^{N-1}_+$ obtained via relations of type
$$p_{i_1i_2}\ldots p_{i_{k-1}i_k}=p_{i_{\sigma(1)}i_{\sigma(2)}}\ldots p_{i_{\sigma(k-1)}i_{\sigma(k)}},\ \forall (i_1,\ldots,i_k)\in\{1,\ldots,N\}^k$$
with $\sigma$ ranging over a certain subset of the infinite permutation group
$$\bigcup_{k\in2\mathbb N}S_k$$
which is stable under the operation $\sigma\to|\sigma|$.
\end{definition}

Observe the similarity with the corresponding notion for the spheres, axiomatized before. The only subtlety in the projective case is the stability under the ``sandwiching'' operation $\sigma\to|\sigma|$, which in practice means that if the above relation associated to $\sigma$ holds, then the following relation, associated to $|\sigma|$, must hold as well:
$$p_{i_0i_1}\ldots p_{i_ki_{k+1}}=p_{i_0i_{\sigma(1)}}p_{i_{\sigma(2)}i_{\sigma(3)}}\ldots p_{i_{\sigma(k-2)}i_{\sigma(k-1)}}p_{i_{\sigma(k)}i_{k+1}}$$

As an illustration, the basic projective spaces are all monomial:

\begin{proposition}
The $3$ projective spaces are all monomial, with the permutations
$$
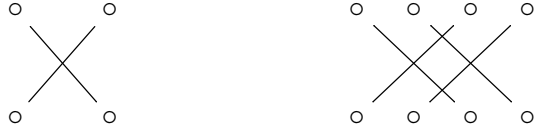
\xymatrix@R=10mm@C=8mm{\circ\ar@{-}[dr]&\circ\ar@{-}[dl]\\\circ&\circ}\qquad\ \qquad\ \qquad 
\xymatrix@R=10mm@C=3mm{\circ\ar@{-}[drr]&\circ\ar@{-}[drr]&\circ\ar@{-}[dll]&\circ\ar@{-}[dll]\\\circ&\circ&\circ&\circ}$$
producing respectively the spaces $P^{N-1}_\mathbb R,P^{N-1}_\mathbb C$.
\end{proposition}

\begin{proof}
We must divide the algebra $C(P^{N-1}_+)$ by the relations associated to the diagrams in the statement, as well as those associated to their shifted versions, given by:
$$
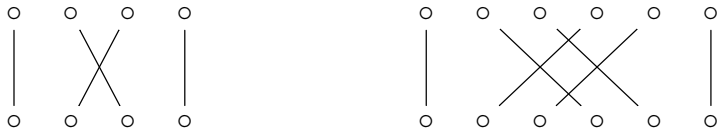
\xymatrix@R=10mm@C=3mm{\circ\ar@{-}[d]&\circ\ar@{-}[dr]&\circ\ar@{-}[dl]&\circ\ar@{-}[d]\\\circ&\circ&\circ&\circ}
\qquad\ \qquad\ \qquad 
\xymatrix@R=10mm@C=3mm{\circ\ar@{-}[d]&\circ\ar@{-}[drr]&\circ\ar@{-}[drr]&\circ\ar@{-}[dll]&\circ\ar@{-}[dll]&\circ\ar@{-}[d]\\\circ&\circ&\circ&\circ&\circ&\circ}$$ 

(1) The basic crossing, and its shifted version, produce the following relations: 
$$p_{ab}=p_{ba}\quad,\quad 
p_{ab}p_{cd}=p_{ac}p_{bd}$$

Now by using these relations several times, we obtain the following formula:
\begin{eqnarray*}
p_{ab}p_{cd}
&=&p_{ac}p_{bd}\\
&=&p_{ca}p_{db}\\
&=&p_{cd}p_{ab}
\end{eqnarray*}

Thus, the space produced by the basic crossing is classical, $P\subset P^{N-1}_\mathbb C$. By using one more time the relations $p_{ab}=p_{ba}$ we conclude that we have $P=P^{N-1}_\mathbb R$, as claimed.

\medskip

(2) The fattened crossing, and its shifted version, produce the following relations:
$$p_{ab}p_{cd}=p_{cd}p_{ab}$$
$$p_{ab}p_{cd}p_{ef}=p_{ad}p_{eb}p_{cf}$$

The first relations tell us that the projective space must be classical, $P\subset P^{N-1}_\mathbb C$. Now observe that with $p_{ij}=z_i\bar{z}_j$, the second relations read:
$$z_a\bar{z}_bz_c\bar{z}_dz_e\bar{z}_f=z_a\bar{z}_dz_e\bar{z}_bz_c\bar{z}_f$$

Since these relations are automatic, we have $P=P^{N-1}_\mathbb C$, and we are done.
\end{proof}

We can now formulate a main classification result, as follows:

\begin{theorem}
The basic projective spaces, namely 
$$P^{N-1}_\mathbb R\subset P^{N-1}_\mathbb C\subset P^{N-1}_+$$
are the only monomial ones.
\end{theorem}

\begin{proof}
We follow the proof from the affine case. Let $\mathcal R_\sigma$ be the collection of relations associated to a permutation $\sigma\in S_k$ with $k\in 2\mathbb N$, as in Definition 16.24. We fix a monomial projective space $P\subset P^{N-1}_+$, and we associate to it subsets $G_k\subset S_k$, as follows:
$$G_k=\begin{cases}
\{\sigma\in S_k|\mathcal R_\sigma\ {\rm hold\ over\ }P\}&(k\ {\rm even})\\
\{\sigma\in S_k|\mathcal R_{|\sigma}\ {\rm hold\ over\ }P\}&(k\ {\rm odd})
\end{cases}$$

As in the affine case, we obtain in this way a filtered group $G=(G_k)$, which is stable under removing outer strings, and under removing neighboring strings. Thus the computations from the affine case apply, and show that we have only 3 possible situations, corresponding to the 3 projective spaces in Proposition 16.25.
\end{proof}

Let us discuss now similar results for the projective quantum groups. We recall that given a closed subgroup $G\subset O_N^+$, its projective version $G\to PG$ is by definition given by the fact that $C(PG)\subset C(G)$ is the subalgebra generated by the following variables:
$$w_{ij,ab}=u_{ia}u_{jb}$$

In relation now with easiness, we can formulate the following definition:

\begin{definition}
A projective category of pairings is a collection of subsets 
$$NC_2(2k,2l)\subset E(k,l)\subset P_2(2k,2l)$$
which is stable under the usual categorical operations, and satisfying as well
$$\sigma\in E\implies |\sigma|\in E$$
called projective categorical condition.
\end{definition}

As basic illustrations for this notion, we have the following projective categories of pairings, where $P_2^*$ is the category of matching pairings:
$$NC_2\subset P_2^*\subset P_2$$

This follows indeed from definitions. Now with the above notion in hand, we can formulate the following projective analogue of our usual notion of easiness:

\index{projective easiness}

\begin{definition}
An intermediate quantum group $PO_N\subset H\subset PO_N^+$ is called projectively easy when its Tannakian category
$$span(NC_2(2k,2l))\subset Hom(v^{\otimes k},v^{\otimes l})\subset span(P_2(2k,2l))$$
comes via via the following formula, using the standard $\pi\to T_\pi$ construction,
$$Hom(v^{\otimes k},v^{\otimes l})=span(E(k,l))$$
for a certain projective category $E=(E(k,l))$.
\end{definition}

Observe that, given any easy quantum group $O_N\subset G\subset O_N^+$, its projective version $PO_N\subset PG\subset PO_N^+$ is projectively easy in our sense. In particular the quantum groups $PO_N\subset PU_N\subset PO_N^+$ are all projectively easy. We have in fact the following result:

\begin{theorem}
We have a bijective correspondence between the affine and projective categories of partitions, given by $G\to PG$ at the quantum group level.
\end{theorem}

\begin{proof}
The construction of correspondence $D\to E$ is clear, simply by setting:
$$E(k,l)=D(2k,2l)$$

Indeed, due to the axioms for the categories of partitions, the conditions in Definition 16.27 are satisfied. Conversely, given $E=(E(k,l))$ as in Definition 16.27, we can set:
$$D(k,l)=\begin{cases}
E(k,l)&(k,l\ {\rm even})\\
\{\sigma:|\sigma\in E(k+1,l+1)\}&(k,l\ {\rm odd})
\end{cases}$$

Our claim is that $D=(D(k,l))$ is a category of partitions. Indeed:

\medskip

(1) The composition action is clear. Indeed, when looking at the numbers of legs involved, in the even case this is clear, and in the odd case, this follows from:
$$|\sigma,|\sigma'\in E
\implies|^\sigma_\tau\in E
\implies{\ }^\sigma_\tau\in D$$

(2) For the tensor product axiom, we have 4 cases to be investigated, depending on the parity of the number of legs of $\sigma,\tau$, as follows:

\medskip

-- The even/even case is clear. 

\medskip

-- The odd/even case follows from the following computation:
$$|\sigma,\tau\in E
\implies|\sigma\tau\in E
\implies\sigma\tau\in D$$

-- Regarding now the even/odd case, this can be solved as follows:
\begin{eqnarray*}
\sigma,|\tau\in E
&\implies&|\sigma|,|\tau\in E\\
&\implies&|\sigma||\tau\in E\\
&\implies&|\sigma\tau\in E\\
&\implies&\sigma\tau\in D
\end{eqnarray*}

-- As for the remaining odd/odd case, here the computation is as follows:
\begin{eqnarray*}
|\sigma,|\tau\in E
&\implies&||\sigma|,|\tau\in E\\
&\implies&||\sigma||\tau\in E\\
&\implies&\sigma\tau\in E\\
&\implies&\sigma\tau\in D
\end{eqnarray*}

(3) Finally, the conjugation axiom is clear from definitions.

\medskip

It is clear that both compositions $D\to E\to D$ and $E\to D\to E$ are the identities, as claimed. As for the quantum group assertion, this is clear as well.
\end{proof}

We refer to \cite{ba2} and related papers for further details, and comments on the above correspondence, and for some classification results as well.

\bigskip

Finally, let us discuss the axiomatization question for the projective quadruplets of type $(P,PT,PU,PK)$. We first have a classical real quadruplet, as follows:
$$\xymatrix@R=50pt@C=50pt{
P^{N-1}_\mathbb R\ar@{-}[r]\ar@{-}[d]\ar@{-}[dr]&PT_N\ar@{-}[l]\ar@{-}[d]\ar@{-}[dl]\\
PO_N\ar@{-}[u]\ar@{-}[ur]\ar@{-}[r]&PH_N\ar@{-}[l]\ar@{-}[ul]\ar@{-}[u]}$$

We have then a classical complex quadruplet, which can be thought of as well as being a real half-classical quadruplet, which is as follows:
$$\xymatrix@R=50pt@C=50pt{
P^{N-1}_\mathbb C\ar@{-}[r]\ar@{-}[d]\ar@{-}[dr]&P\mathbb T_N\ar@{-}[l]\ar@{-}[d]\ar@{-}[dl]\\
PU_N\ar@{-}[u]\ar@{-}[ur]\ar@{-}[r]&PK_N\ar@{-}[l]\ar@{-}[ul]\ar@{-}[u]}$$

Finally, we have a free quadruplet, which can be thought of as being the same time real and complex, which is as follows:
$$\xymatrix@R=50pt@C=50pt{
P^{N-1}_+\ar@{-}[r]\ar@{-}[d]\ar@{-}[dr]&PT_N^+\ar@{-}[l]\ar@{-}[d]\ar@{-}[dl]\\
PO_N^+\ar@{-}[u]\ar@{-}[ur]\ar@{-}[r]&PH_N^+\ar@{-}[l]\ar@{-}[ul]\ar@{-}[u]}$$

Thus, modulo some axioms, which remain to be worked out, things are quite nice, because we have here only 3 geometries, namely real, complex and free.

\section*{16e. Exercises}

Congratulations for having read this book, and no exercises for this final chapter. However, if looking for some difficult exercises, in relation with this, you can check the literature referenced below, where many interesting questions are left open.

\baselineskip=14pt

\printindex

\end{document}